\documentclass[a4paper,11pt]{article}

\usepackage[T1]{fontenc}
\usepackage{theorem,amsmath,amssymb,euscript}
\usepackage[applemac]{inputenc}
\usepackage{stmaryrd}
 
\usepackage[french]{babel}
\usepackage{a4wide}
\usepackage{xspace}

\def \goth{\mathfrak}
\newcommand{\D}{\mathbb{D}}
\newcommand{\E}{\mathbb{E}}
\newcommand{\C}{\mathbb{C}}
\newcommand{\F}{\mathbb{F}}
\newcommand{\f}{\F^{*2}}
 
\newcommand{\Z}{\mathbb{Z}}
 \newcommand{\R}{\mathbb{R}}
 \newcommand{\Q}{\mathbb{Q}} 
 \newcommand{\N}{\mathbb{N}}
\newcommand{\∆}{\Delta}
\newcommand{\∑}{\Sigma}
 \newcommand{\g}{\goth g}
\newcommand{\noi}{ \noindent }
\newcommand{\ds}{ \displaystyle}
\newcommand{\Frac}{ \displaystyle\frac}
\newcommand{\fdem}{\hfill$\Box$ }
\newcommand{\é}{\' e}

\newcommand{\à}{\`a}
\newcommand{\ù}{\`u}

\newcommand{\û}{\^u}

\newcommand{\ssi}{$\Leftrightarrow\ $}
\newcommand{\sg}{\xspace sous-groupe\xspace}
\newcommand{\sgs}{sous-groupes\xspace}
\newcommand{\≥}{\geq}
\newcommand{\≤}{\leq}
\newcommand{\elt}{\'el\'ement\xspace}
\newcommand{\elts}{\'el\'ements\xspace }
\newcommand{\PV}{pr\éhomog\ène\xspace}
\newcommand{\PVs}{pr\éhomog\ènes\xspace}
\newcommand{\irf}{invariant relatif fondamental\xspace}
\newcommand{\irfs}{invariants relatifs fondamentaux\xspace}
 \newcommand{\sgp}{sous-groupe parabolique\xspace}
  \newcommand{\sgps}{sous-groupes paraboliques\xspace}

\newcommand{\Sl}{$sl_2$-triplet\xspace}
\newcommand{\Sls}{$sl_2$-triplets\xspace}
\newcommand{\dg}{degr\'e\xspace}
\newcommand{\demo}{d\émonstration }
\newcommand{\gog}{\xspace$(\goth g_0,\goth g_1)$\xspace}
\newcommand{\eq}{\'equation fonctionnelle\xspace}
\newcommand{\eqs}{ \'equations fonctionnelles\xspace}
\newcommand{\propo}{proposition \xspace}
\newcommand{\air}{absolument irr\'eductible r\'egulier\xspace}
\newcommand{\dpl}{d\'eploy\'e\xspace}
\newcommand{\dple}{d\'eploy\'ee\xspace}
\newcommand{\four}{\EuScript  F}
 \newcommand{\theor}{\xspace th\' eor\`eme\xspace}
 \newcommand{\FF}{ \F^*/\F^{*2}}
   \newcommand{\ES}{\EuScript S}
   
\newtheorem{theo}{Th\'eor\`eme}[subsection]
 
\newtheorem{prop}[theo]{Proposition} 
\newtheorem{lem}[theo]{Lemme} 
\newtheorem{cor}[theo]{Corollaire} 
 
\newtheorem{defi}[theo]{D\' efinition} 
\newtheorem{rien}[theo]{}
\newtheorem{rema}[theo]{Remarque} 
\newtheorem{rappel}[theo]{Rappel}

\def \dem {\paragraph{D\'emonstration: }}

   \date{ 24  janvier 2007}
   
\begin{document} 
 
  \centerline     {\bf  Une classe d'espaces pr\éhomog\ènes de type parabolique faiblement sph\ériques }
  \vskip 5mm   
   \centerline {Iris MULLER }
 \vskip 5mm 
  
\noi {Institut de Recherche Math\'ematique  
 Avanc\'ee (IRMA)}
  
\noi {Universit\'e Louis Pasteur et CNRS (UMR  7501)}

 \vskip 2cm
 \begin{abstract}
 Pour les alg\èbres de Lie absolument simples, de dimension finies, de rang au moins $2,$ d\éfinies sur un corps local de caract\éristique $0$ et admettant une graduation: $\g=\g_{-2}\oplus \g_{-1}\oplus \g_{0}\oplus \g_{1}\oplus \g_2,$ donn\ée par un \elt $H_0$ tel que $2H_0$ soit simple, on d\éfinit (\à $2$ exceptions pr\ès) des \sgs paraboliques $P,$ inclus dans le groupe des automorphismes de $\g$ et centralisant $H_0,$ dont l'action sur $\g_1$ et $\g_{-1}$ est g\éom\étriquement \PV.  On \étudie la structure de ces espaces \PVs. On montre que les fonctions Z\étas associ\ées aux \irfs de $P$ d\éfinis sur $\g_1$ et $\g_{-1}$ admettent des prolongements m\éromorphes qui v\érifient des \eqs abstraites et on donne le calcul explicite des coefficients des \eqs, des polyn\ômes de Bernstein associ\és aux \irfs de $P$ dans le cas archim\édien, par une m\éthode de descente  \à des centralisateurs de  paires d'alg\èbres de Lie de type $sl_2$ qui commutent.
 
 \noi Ceci g\én\éralise des r\ésultats bien connus lorsque $\g_2=\{0\}.$

 \end{abstract}

 \vskip 2cm
 \centerline{\bf Introduction}
\vskip 5mm
Comme il est dit dans \cite{sato6}, un probl\ème int\éressant en analyse et en th\éorie des nombres est le suivant: \étant donn\és $2$ polynomes homog\ènes $P$ et $P^*$ en $n$ variables de degr\é $d$ \à coefficients r\éels, trouver des conditions afin que la transform\ée de Fourier de $|P(x)|^s,$ $s$ \étant un nombre complexe,  v\érifie l'\égalit\é $\widehat{|P(x)|^s}=$  facteur gamma.$|P^*(x)|^{-s-\frac{n}{d}}$ au sens des distributions.\\

Il est bien connu qu'une r\éponse \à cette question est donn\ée par la th\éorie des espaces vectoriels \PVs d\ûe \à Mikio Sato et introduite dans les ann\ées 1960 (\cite{msato}, premi\ère version publi\ée en 1970 en Japonais) et dont une classification a \ét\é obtenue par M.Sato et T.Kimura (\cite{satokimura}).

Par la suite cette th\éorie a connu de grands d\éveloppements et dans de nombreuses directions.\\ \\
Dans ce travail, on s'int\éresse aux fonctions z\étas locales associ\ées \à l'action \PV de certains \sgps ceci dans le cadre des \PVs de type parabolique.\\ \\

 {\bf 1.} Les \PVs consid\ér\és sont de ``type parabolique" (\cite{rubthese},\cite{rublivre}),  plus pr\écis\ément soit $\goth g$ une alg\èbre de Lie  absolument simple de
dimension finie, d\éfinie sur un corps local  
 de caract\éristique z\éro que l'on note
$\F,$ on suppose que $\goth g$ est munie
d'une graduation :
$$\goth g=\oplus_{i\in
\Z}\goth g_i$$
o\ù les $\goth g_i$ sont des sous-espaces vectoriels de $\goth
g $
v\érifiant la condition :
$$[\goth g_i,\goth g_j]\subset \goth g_{i+j}\quad \forall
i,j\ \in \ \Z.$$ 
 
Les d\érivations de $\goth g$ \étant int\érieures, il existe un
unique \él\ément appartenant \à $\goth g_0$ qui d\éfinit la
graduation, on le note $H_0$ et on a:
$$\forall i\in \Z \quad \goth g_i= 
 \{x\in \goth g\ |\ [H_0,x]=ix\}.$$ 
Soit $G$ le centralisateur de $H_0$ dans le sous-groupe
$Aut_0(\goth g)$ (\cite{bourbakigal8}) des automorphismes de $\goth g.$ Ce
groupe op\ère sur chaque espace vectoriel $\goth g_i, i\in \Z,$ mais on consid\ère l'action de $G$ sur $\g_1$ et la
repr\ésentation correspondante est not\ée de mani\ère
infinit\ésimale : $(\goth g_0,\goth g_1,H_0)$ ou $(\goth
g_0,\goth g_1).$

Les sous-espaces vectoriels $\g_1$ et $\g_{-1}$ sont mis en dualit\é \à l'aide de la forme de Killing de $\g$ d'o\ù $(G,\g_{-1})$ est la repr\ésentation contragr\édiente de $(G,\g_1).$\\

 D'apr\ès un r\ésultat de E.B.Vinberg
(\cite{vinberg}), ces repr\ésentations sont g\éom\étriquement pr\éhomog\ènes
ce qui signifie que $G$ poss\ède une orbite ouverte dans
$\goth g_1$ muni de la topologie de Zariski, ceci lorsqu'on
se place sur une cl\ôture alg\ébrique de $\F,$ et le nombre d'orbites de $G$ dans $\g_{\pm1}$ est fini.

Ainsi
$(\goth g_0,\goth g_1)$ est une $\F$-forme d'un espace
pr\éhomog\ène, appel\é de type {\it parabolique},
en effet la sous-alg\èbre
$\oplus_{i\≥0}\goth g_i$ est une sous-alg\èbre parabolique de
$\goth g$ dont le radical nilpotent est
 la sous-alg\èbre : $\oplus_{i\≥1}\goth g_i$ et la
sous-alg\èbre r\éductive  $\goth g_0$ est l'alg\èbre de Lie du
groupe $G$ (\cite{rubthese}). Lorsque la
graduation est courte, c'est \à dire telle que $\goth g=\goth
g_{-1}\oplus \goth g_0\oplus\goth g_1,$  le \PV est dit de``type
commutatif" (car $\goth g_1$ est une alg\èbre commutative).\\

L'int\ér\êt des \PVs de
type paraboliques r\éside dans la pr\ésence de l'alg\èbre
simple $\goth g$ qui contient \à la fois l'alg\èbre de
Lie du groupe $G$ et l'espace $\goth g_1$ de la
repr\ésentation, permettant ainsi l'\étude et l'expression des
r\ésultats en termes de la structure  de $\goth g$ (cf. par
exemple les travaux de N. Bopp et H. Rubenthaler dans \cite{boppruben}). \\

Ainsi une description possible des
orbites de $G$ dans $\goth g_1$ est r\éalis\ée \à l'aide d'une
version 
$\ll$gradu\ée$\gg$ des
$sl_2-$triplets. On rappelle qu'un $sl_2-$triplet est un
triplet non nul $(x,h,y)$ v\érifiant les relations de
commutation suivantes :
$$[x,y]=-h,\ \ [h,x]=2x,\ \ [h,y]=-2y\quad (\cite{bourbakigal8}).$$
Par une g\én\éralisation du th\éor\ème de
Jacobson-Morozov, il est bien connu que pour tout $x$ non
nul de $\goth g_1$, il existe $h$ et $y$ appartenant
respectivement \à $\goth g_0$ et \à $\goth g_{-1}$ tels que
$(x,h,y)$ soit un $sl_2-$triplet.
 \vskip 3mm  
 \noi{\bf D\éfinition} {\it
 Un \'el\'ement $h$ de $\g$ est dit {\it 1-simple} si il existe un
$sl_2$-triplet
$(x,h,y)$ tel que $x$ (resp. $y)$ soit dans $\goth g_1$ (resp.
$\goth g_{-1}$) . Un tel $sl_2-$triplet est appel\'e {\it
1-adapt\'e}.}
\vskip 2mm
 Soit $\goth a$ une sous-alg\èbre ab\'elienne
d\'eploy\'ee maximale de $\goth g$ contenant $H_0,$ on note
$\Delta$ le syst\ème de racines correspondant. Ce dernier est
\'egalement gradu\'e par $H_0:$ 
$$\Delta_i=\{\lambda\in \Delta\ |\   \lambda(H_0)=i\}.$$ On choisit un ordre
sur $\Delta$ pour lequel les racines appartenant \à  
$\cup_{i>0}\Delta_i $ soient positives, soit $\Sigma$  une base de $\Delta$ et $\Sigma_i=\Sigma\cap \Delta_i ,$ $i\in \Z,$ on adopte la notation classique $(\Delta,\Sigma-\Sigma_0)$ (ou bien $(\Delta,\lambda_0)$ lorsque $\Sigma-\Sigma_0=\{\lambda_0\}$) pour d\ésigner le \PV $(\g_0,\g_1)$ (\cite{rubthese}), le \PV $(\g_0,\g_1)$ sera dit {\it d\éploy\é} lorsque $\g$ est d\éploy\ée.\\

Les \PVs de type parabolique absolument
irr\éductibles, c'est \à dire pour lesquels $\goth g_1$ {\it
est un}
$\goth g_0-${\it module absolument irr\éductible}, sont
r\éguliers (\cite{satokimura}) si et seulement si $2H_0$ est $1$-simple
(\cite{rubthese}). Dans ces conditions, on peut supposer que $\g_1$ et $\g_{-1}$ engendrent $\g$, ce que l'on fera syst\ématiquement ainsi le \PV sera de type $(\Delta,\lambda_0),$  il existe (\cite{rubthese}) un polynome
d\éfini sur $\goth g_1,$ relativement invariant par $G$ et de
degr\é minimal, que l'on note $F,$ tel que :
$$ \{x\in \goth g_1\ \hbox{admettant
}\ 2H_0\
\hbox{pour \él\ément $1$-simple}\} 
=\{x\in \goth g_1\  |\ F(x)\not=0\}. $$
Il existe un caract\ère $\chi$ de $G$ tel que:
$$\forall g\in G,\ \forall x\in \goth g_1: \
F(gx)=\chi(g)F(x)\quad (1)$$
et toute fraction rationnelle d\éfinie sur $\goth g_1$
v\érifiant une relation analogue \à (1) est (\à une constante
multiplicative pr\ès) une puissance enti\ère de $F.$
L'invariant relatif $F,$ qui lui aussi est d\éfini \à une
constante pr\ès, est dit {\it fondamental} et l'ensemble de
ses  z\éros est le lieu {\it singulier} not\é $S.$  \\ 

On a exactement la m\ême situation sur le \PV ``dual'', $(G,\g_{-1}),$ l'\irf sera not\é $F^*$ et son lieu singulier $S^*.$ \\ \\

Soient:

$O_1,...,O_\ell$ les orbites de $G$ dans $\g_1-S$ et $O^*_1,...,O^*_\ell$ les orbites de $G$ dans $\g_{-1}-S^*,$\\

$\omega$ un caract\ère unitaire de $\F^*,$ $\pi=\omega|\ |^s,$ avec  $s\in \C,$ un caract\ère de $\F^*,$\\

Pour $f$  appartenant \à $ \EuScript S(\g_1),$ l'espace de Schwartz de $\g_1$ et $g$ appartenant \à $ \EuScript S(\g_{-1}),$ on appelle {\it fonctions z\étas locales} les expressions:
$$Z_{O_i}(f;\pi)=\int_{O_i}f(x)\pi(F(x))dx \ ,\ Z_{O_i^*}(g;\pi)=\int_{O_i^*}g(x)\pi(F^*(x))dx $$
qui convergent pour partie r\éelle de $s$ assez grand, elles admettent un prolongement m\éromorphe sur $\C$ et il existe des fonctions m\éromorphes (en $s$), not\ées $a_{O_i^*,O_j}(\pi)$, telles que pour $i=1,...,\ell$ et $f$ appartenant \à $ \EuScript S(\g_1)$ on ait:
$$Z_{O_i^*}(\widehat f;\pi)=\sum_{j=1,...,\ell}a_{O_i^*,O_j}(\pi)Z_{O_i}(f;\pi^{-1}|\ |^{-N})    \ 
 \text{avec }N=\ds\frac{\text{dim}(\g_1)}{\text{degr\é de }F}\quad (\cite{sato-shintani}, \cite{formescubiques}, \cite{igusa5}).$$
 
 Dans le cas r\éel, les fonctions $a_{O_i^*,O_j}$ ont une forme particuli\ère (\cite{formescubiques}). Elles ont \ét\é calcul\ées dans un grand nombre de cas (dans le cas commutatif et dans le cadre des alg\èbres de Jordan \cite{safa}, cf. les travaux de M.Muro, par exemple \cite{muro1} , \cite{muro2} , \cite{kimuramuro}, par des techniques de micro-analyse ).  Dans le cas $\goth p$-adique, l'\étude fondamentale est d\ûe \à I.J.Igusa (cf. les travaux cit\és et en particulier \cite{igusa11} et \cite{igusa12} pour un panorama) et \à ses \él\èves (par ex. \cite{robinson}).\\ \\
 
Rappelons tr\ès bri\èvement que l'\étude des fonctions z\étas globales qui sont d\éfinies:
 
 - soit  \à partir de la situation d\écrite ci-dessus mais en prenant comme corps $\F$ un corps de nombres, la fonction z\éta est alors celle obtenue par extension sur les ad\èles de $\F,$  son prolongement m\éromorphe  v\érifie \également une \eq (cf. par exemple: appendice de \cite{mars}, \cite{rallisschiffmann}, \cite{wright}, \cite{igusa6}, \cite{kimurakogiso}, \cite{mullerRims}, \cite{saito2}, \cite{saito3}),
 
 - soit \à partir de la situation r\éelle et de s\éries de Dirichlet associ\ées \à des r\éseaux $G_\Z$-stable dans $\g_1$ (cf. par exemple \cite{sato-shintani} , \cite{formescubiques}, \cite{shintani}, \cite{sato1}, \cite{sato2}),\\
est nettement plus difficile et qu'il est regrettable que toute  la bibliographie que nous donnons \à titre indicatif soit non exhaustive.\\ \\

 Dans le cas particulier des \PVs commutatifs pour lesquels l'alg\èbre $\g$ n'est pas de rang $1,$ il existe (au moins) un \sgp $P$ de $G$ dont l'action sur $\g_{\pm 1}$ est encore g\éom\étriquement pr\éhomog\ène  (\cite{M-R-S}). Lorsque le \sgp est minimal parmi ceux-ci, les  \irfs correspondants, not\és $F_1,...,F_p$ sur $\g_1$ et $F^*_1,...,F^*_p$ sur $\g_{-1}$  donnent \à nouveau des fonctions z\étas \à multi-indice qui v\érifient \également des\eqs analogue aux pr\éc\édentes, ces r\ésultats, qui se trouvent d\éj\à dans \cite{shintani}, sont d\ûs dans le cas r\éel \à N.Bopp et H.Rubenthaler (\cite{boppruben}) ainsi qu'\à J.Faraut et A.Koranyi lorsque le \sgp est minimal (\cite{farautkoranyi}), et ce dans le cadre des alg\èbres de Jordan, et \à Y.Hironaka dans certains cas $\goth p$-adiques  (\cite{hironaka1}).\\
Dans cette situation commutative {\it uniquement}, chaque orbite ouverte, $O_i,i=1,...,\ell,$ est la r\éalisation d'un espace sym\étrique $G/H_i.$ Dans \cite{boppruben}, N.Bopp et H.Rubenthaler ont g\én\éralis\és les travaux de R.Godement et H.Jacquet (\cite{godementjacquet}) en associant des fonctions z\étas locales \à des vecteurs distributions  $H_i-$invariants associ\és \à la m\ême repr\ésentation sph\érique minimale de $G$ et ont \établis les \équations fonctionnelles correspondantes.

Lorsque le parabolique est minimal, les fonctions sph\ériques et les probl\èmes associ\és ont \ét\é largement \étudi\és par J.Faraut et A.Koranyi (\cite{farautkoranyi}) dans le cadre des alg\èbres de Jordan et Y.Hironaka dans certains cas $\goth p$-adiques.\\

Dans le cas r\éel, J.L. Clerc a r\éussi  \à g\én\éraliser certains  r\ésultats pr\éc\édents  \à une classe plus large de \PVs non commutatifs \à l'aide des repr\ésentations d'alg\èbres de Jordan (\cite{clerc}).
\vskip 3mm
 
{\bf 2.}  Dans ce travail, on se propose de donner une situation analogue au cas commutatif c'est \à dire d\éfinir pour chaque \PV de type parabolique absolument irr\éductible r\égulier (cf.1.) de graduation au plus $5$ c'est \à dire
$$\g= \g_{-2}\oplus  \g_{-1}\oplus  \g_{0}\oplus  \g_{1}\oplus  \g_{2}$$
de rang au moins $2$ dont l'\irf n'est pas une forme quadratique (i.e. degr\é($F)>2$) et qui n'est pas de type $\bf G_2,$ (au moins) un \sgp standard, $P,$ dont l'action est \PV au sens pr\éc\édent. 

D\ésignant par $F_1,...,F_p=F,F^*_1,...,F^*_p=F^*$ les \irfs de $(P,\g_1)$ et $(P,\g_{-1})$ rang\és dans un ``ordre croissant" (essentiellement chaque $F_i$ apparait comme la restriction de $F_{i+1}$), on montre classiquement que la fonction z\éta locale \à multi-indice associ\ée v\érifie une \eq . On d\éfinit des normalisations coh\érentes avec cette situation et on montre que les coefficients locaux (ainsi que les polyn\ômes de Bernstein associ\és dans le cas archim\édien) s'obtiennent par descente sur des \PVs construits \à partir de centralisateurs de paires de \Sls $1$-adapt\és qui commutent.

La m\éthode utilis\ée convient sur $\F$ sans distinction (archim\édien ou $\goth p$-adique),  elle est \él\émentaire (d\écomposition des \irfs et des mesures) et g\én\éralise celle du cas commutatif (\cite{muller1})\footnote{Ces r\ésultats  ont \ét\é partiellement expos\és aux Journ\ées Pr\éhomog\ènes Franco-Japonaises organis\ées par F.Sato-P.Kaplan-A.L.Mortajine  \à Tokyo en 1999}.

Il se peut qu'il y ait un lien avec certaines d\écompositions des \irfs consid\ér\ées par A.L.Mortajine (\cite{mortajine}) ainsi que par F.Sato (\cite{sato6}). \\

$(\Delta,\lambda_0)$ \étant le diagramme de Dynkin gradu\é associ\é au \PV $(\g_0,\g_1),$ on note  $({\overline \Delta}, \overline{\lambda_0})$ le diagramme de Dynkin gradu\é associ\é au pr\éhomo\-g\ène  d\éploy\é sur une extension convenable de $\F$ et admettant $(\g_0,\g_1)$ comme $\F-$forme.\\
Dans le cas $\goth p$-adique, les r\ésultats obtenus montrent le lien entre les p\ôles des coefficients $a_{O^*,O}(|\ |^{s_1},...,|\ |^{s_p})$ et certaines racines du polynome de Bernstein associ\é \à $F$ dans le \PV d\éploy\é r\éel de m\ême type $({\overline \Delta},\overline{\lambda_0}).$

Ainsi dans le cas complexe, on obtient le r\ésultat attendu, c'est \à dire que si on note:
$$F(\partial)(\prod_{i=1}^pF_i^{*s_i})=b_{\g,P}(s_1,...,s_p)(\prod_{i=1}^{p-1}F_i^{*s_i} )F_p^{*{s_p-1}}\quad (s_1,...,s_p)\in \C^p,$$
alors (\à $2$ exceptions pr\ès) on a 
 $$b_{\goth g,P_{\goth t}}(s_1,...,s_p)= \prod_{\ell=1}^p\biggl(\prod_{j=1}^{d_{p-\ell+1}-d_{p-\ell}}(s_{\ell}+...+s_p+
\lambda_{\ell,j})\biggr)$$
et si on d\éfinit sur $\F,$ archim\édien ou $\goth p$-adique, pour toute application polynomiale $B$ de la forme: $$B(s_1,...,s_p)=\prod_{1\≤j\≤r}(a_{1,j}s_1+...+a_{p,j}s_p+q_j) \text{ avec }a_{i,j}\in \Z, $$ 
et   pour tout caract\ère $\pi$ de ${(\F^*)}^p$, $\pi=(\pi_1 ,...,\pi_p ):$ $$\rho'_B(\pi)=\rho'_B(\pi_{ 1},...,\pi_{ p} )=
\prod_{j=1}^r \rho'(\pi_ 1^{a_{1,j}}....\pi_p^{a_{p,j}}|\ |^{  q_j+1})
 $$ avec $ \rho'(\pi_1)=\pi_1(-1)\rho(\pi_1),$ $\rho$ \étant le coefficient obtenu par J.Tate (\cite{tate}) dans l'\eq associ\ée \à la fonction z\éta sur $\F$ (dans ce cas $F(x)=x$ et les valeurs explicites de $\rho$ sont rappel\ées dans le $\S 3.6.1$)
  alors pour $f$ dans $\EuScript S(\g_1)$ on a:
$$(A)\ Z^*(\hat f;\pi )=\rho'_{b_{\goth g,P_{\goth t}}}(\pi)Z(f;\pi^*|\ |^{-N1_p})  \ \ \text{ avec }\ 
\pi^*=(\pi_{p-1} ,...,\pi_1,(\prod_{1\≤i\≤p}\pi_i)^{-1})\text{ et }1_p=(0,...,0,1).$$  

On applique la m\éthode de descente pour le  \sgp $P_0,$ minimal parmi les \sgps consid\ér\és (\à l'exception d'un cas qui ne semble pas significatif) ceci  nous conduit \à consid\érer successivement le cas commutatif (cf.tableaux 1 et 2), les cas classiques symplectiques ($n^\circ 13$ de la classification de \cite{satokimura} avec par exemple $Trig(2k)\times Sp(n-2k),$ avec $6k\≤2n,$ op\érant sur les matrices \à $2k$ lignes et $2n-4k$ colonnes ainsi que les diff\érentes $\F$ formes) et orthogonaux ($n^\circ 15$ de la classification de \cite{satokimura}, avec par exemple $Trig(k)\times SO(m),$ op\érant sur les matrices \à $k$ lignes et $m$ colonnes ainsi que les diff\érentes $\F$ formes avec quelques hypoth\èses techniques suppl\émentaires) (cf.tableau 2) puis les cas ``exceptionnels" ayant $\g_2$ de dimension $1$ 

\noi ($(\overline \Delta,\alpha_0)$ de type  $(F_4,\alpha_1)$ ou $(E_6,\alpha_2)$ ou $(E_7,\alpha_1)$ ou $(E_8,\alpha_8)$ avec les formes EIII, EVI et EIX lorsqu'elles existent, cf.les $n^\circ 14,5,23,29$ de la classification de \cite{satokimura}),

\noi puis les $\F-$formes de $(E_7,\alpha_6)$  ($n^\circ 20$ de la classification de \cite{satokimura}), le type d\éploy\é  $(E_7,\alpha_2)$  ($n^\circ 6$ de la classification de  \cite{satokimura}) et pour finir les $\F-$formes de $(E_8,\alpha_1)$ ($n^\circ 24$ de la classification de \cite{satokimura}) (cf.tableau 3).\\

\noi Pour \éviter d'avoir \à d\éterminer les orbites de $P_0$ dans le lieu non singulier de $\g_1$ et $\g_{-1},$ on a pr\éf\ér\é introduire pour $u=(u_1,...,u_p)\in (\FF)^p$ les ouverts (\éventuellement vides):
$$O_u=\{x\in \g_1\ |\ F_1(x)...F_i(x)\F^{*2}=u_1...u_i\ ,i=1,...,p\}$$
$$ O^*_u=\{x\in \g_{-1}\ |\ F^*_1(x)...F^*_i(x)\F^{*2}=u_p...u_{p-i+1}\ ,i=1,...,p\}$$
ainsi que les fonctions z\étas correspondantes, not\ées simplement $Z_u$ et $Z^*_u,$ alors il existe des fonctions m\éromorphes (au sens pr\éc\édent), $a_{v,u}(\pi)$ telles que
$$(B)\quad Z^*_u(\hat f;\pi)=\sum_{v\in (\FF)^p} a_{v,u}(\pi)Z_v(f;\pi^*|\ |^{-N1_p}),$$
elles sont d\étermin\ées explicitement puisque la m\éthode de descente ram\ène les calculs \à des cas bien connus (rang $1$ et on applique \cite{tate} , \cite{godementjacquet}, ou bien au cas d'une forme quadratique \cite{rallisschiffmann}).

\noi Cependant la structure particuli\ère de $\FF$ dans le cas $\goth p-$adique nous a conduit \à introduire certaines sommes qui ont \ét\é calcul\ées uniquement lorsque la caract\éristique r\ésiduelle est diff\érente de $2$ ce qui explique cette restriction chaque fois que l'on utilise les r\ésultats de \cite{rallisschiffmann} (cf.\S 3.6.2).\\

\noi Ce r\ésultat subsiste dans les cas exceptionnels pour le \PV  $(G,\g_{\pm 1})$ (ce sont les d\éfinitions ci-dessus avec $p=1$), et dans le cas $(E_7,\alpha_2)$ r\éel ou $\goth p$-adique  de caract\éristique r\ésiduelle diff\érente de $2,$ on obtient l'\équation:
$$ Z^*( \four(
f); \tilde\omega_b|\ |^s)=C(s)B(\tilde{\omega}_b,s)Z (f;  \tilde\omega_b|\ |^{-s - 5})\ \quad\text{
avec}:$$
$$ C(s)= |2|_\F^{-4s-10}\rho ( |\ |^{2s+4})\rho ( |\ |^{2s+7})\ ,  
B(\tilde\omega_b,s)=  \rho' ( \tilde\omega_b |\ |^{s+ 1})\rho' (\tilde\omega_b |\ |^{s+ 3})\rho' (\tilde\omega_b |\ |^{s+ 5})$$
$\tilde\omega_x$ d\ésignant, pour $x\in\FF $,   le caract\ère de $\FF$ d\éfini par  le symbole de Hilbert ($\tilde\omega_x(y)=(x,y)$).\\ \\

De plus, lorsque les racines de $b_{\g,P_0}(0,..,s)$ sont enti\ères, on a \également l'unique \eq:
$$\pi=(\pi_1,...,\pi_p)\ ,\ Z^*(\hat f;\pi )= C\rho'_{b_{\goth g,P_{\goth t}}}(\pi'_0\pi)Z(f;\pi^*\pi_0|\ |^{-N1_p}) $$
avec $\pi_0=\pi'_0=Id$ \à l'exception de l'unique cas commutatif $(\overline\Delta,\overline{\lambda_0})=(A_{2p-1},\alpha_p) $ et $(\Delta,\lambda_0)=(C_p,\alpha_p)$ dans les notations des planches de \cite{bourbakigal6} (cas des matrices hermitiennes) pour lequel $\begin{cases}\pi_0=(1,...,1,\tilde\omega_{\delta^{p-1} }),\\
\pi'_0=( \tilde\omega_\delta,..., \tilde\omega_\delta,1),\end{cases}$,  
$\F[\sqrt \delta]$ \étant l'extension sur laquelle $\g$ se d\éploie,  \\

\noi et $ C= 1$ lorsque $\g$ est d\éploy\ée ou de type exceptionnel, sinon \\

$C=\begin{cases}(\alpha(1)\alpha(-\delta))^{\frac{p(p-1)}{2}} \ \text{ lorsque }\pi'_0\not=1,\\
\\
(1)^{np}\ \text{ dans le cas CII }\  ((\overline\Delta,\overline{\lambda_0})=(C_n,\alpha_{2p}) \text{ et }  (\Delta,\lambda_0)=(C_\ell,\alpha_p))\\
\\
 (-1)^{\frac{p(p-1)}{2}} \ \text{ dans le cas DIII } ((\overline\Delta,\overline{\lambda_0})=(D_{2p},\alpha_{2p}), (\Delta,\lambda_0)=(C_p,\alpha_p))\\
\\
(-1)^{p(k-1)}\ \text{  dans le cas  AII }((\overline\Delta,\overline{\lambda_0})=(A_{2kp-1},\alpha_{pk}) \text{  et }(\Delta,\lambda_0)=(A_{2p-1},\alpha_p)).
\end{cases}$ 

\noi (r\ésultat connu partiellement).\\ \\
{\bf 3.} Indiquons bri\èvement le contenu de chaque partie.\\

$\bullet$ \underline{Dans la section 1} on motive le choix du parabolique particulier que l'on consid\ère et que l'on introduit de mani\ère inhabituelle, \à partir de certains \elts $1$-simples ``compatibles " avec la graduation de d\épart.

Dans la situation g\én\érale gradu\ée, cf. $\S 1.$ avec $\g$ semi-simple, o\ù l'on suppose $2H_0$ {\it $1$-simple}, \à tout \elt $1$-simple, $h,$ on associe la sous-alg\èbre parabolique : $p(h)=\oplus_{i\≥0}E_i(h)$ de partie nilpotente $n(h)=\oplus_{i>0}E_i(h),$ avec $E_i(h)=\{x\in\g_1|[h,x]=ix\}$ pour $i$ entier. Alors lorsque $h\not=2H_0$ et $n(2H_0)\subset n(h),$ $2H_0-h$ est \également $1$-simple et $exp(ad(n(2H_0-h)\cap\g_0)(E'_2(h)\cap \g_1\oplus E_2(2H_0-h)\cap \g_1)$ est un ouvert de Zariski de $\g_1,$ avec $E'_2(h)=\{x\in \g|(x,h,\ )$ se compl\ète en un \Sl$\}$ (lemme 1.1.1).

Un tel \elt $h$ est dit $1$-simple sp\écial. Lorsque $h$ et $2H_0-h$ sont $1$-simples sp\éciaux, on dit que $h$ est $1$-simple tr\ès sp\écial. Ces \elts tr\ès sp\éciaux sont associ\és \à des graduations courtes  et aux situations irr\éductibles, en effet si le \PV de d\épart est absolument irr\éductible, alors $\g_i=\{0\}$ pour $|i|\≥3$ (lemme 1.2.2) et le \PV associ\é $(E_0(h)\cap\g_0),E_2(h)\cap\g_1,\frac{h}{2})$ est encore absolument irr\éductible r\égulier (lemme 1.2.1) ainsi que le \PV obtenu en prenant le centralisateur d'un \Sl 1-adapt\é construit avec $h$ et muni de la graduation induite par $ad(H_0)$ \à quelques exceptions pr\ès (prop.1.2.4).

On v\érifie ensuite que sous les hypoth\èses \énonc\ées dans 2., un tel \elt $1$-simple tr\ès sp\écial existe toujours ($\S 1.3$).

Dans le $\S 1.4$ on introduit la classe des \PVs \étudi\és, $H_1,...H_p$ sont $p$ \elts $1$-simples tr\ès sp\éciaux de $\goth a$ de somme $2H_0,$ $\goth t=\oplus_{i=1}^p\F H_i,$ $\Delta_R(\goth t)$ les restrictions non nulles des racines de $\Delta$ \à $\goth t$ muni de l'ordre: $$ \lambda\succ 0 \ \Leftrightarrow\
\lambda(H_k)>0\quad \hbox{avec}\quad k=sup\{j\ ,\
\lambda(H_j)\not=0\} \ ,$$
alors $p(H_1,...,H_p)= E_0(\goth  t)\oplus n _{\goth t},$ avec $n _{\goth t}= \oplus_{\lambda\succ 0,\lambda\in
\Delta_0}\goth g^{\lambda} ,$ est une sous-alg\èbre parabolique de \sgp associ\é  : $P_{\goth t}=G_{\goth t}.N _{\goth
t},$ avec $N_{\goth t}=exp(ad(n _{\goth t})), $  
   $G_{\goth t}$ \étant le centralisateur de $\goth t$ dans $G;$ un tel sous-groupe parabolique (que l'on prendra standard) est appel\é ``tr\és sp\écial''. \\
   
 \noi La situation est analogue au cas commutatif  puisqu'on a les propri\ét\és suivantes.\\ 
    
   \noi Pour $k=1,...,p-1,$ soient $h_k=\sum_{1\≤j\≤k}H_j$ et 
 $F_k$  l'\irf du \PV \air $(E_0( h_k)\cap \g_0,E_2( h_k)\cap \g_1),\frac{1}{2}h_k)$  que l'on \étend naturellement \à $\g_1;$  on note $F_p$ l'\irf du \PV \gog. 
 
\noi  On a la situation duale avec pour $k=1,...,p-1,$
   $F^*_{p-k}$   l'\irf du \PV \air $(E_0( h_k)\cap \g_0,E_0(  h_k)\cap \g_{-1}),H_0-\frac{1}{2} h_k)$  que l'on \étend naturellement \à $\g_{-1};$ on note $F^*_p$ l'\irf du \PV $(\g_0,\g_{-1}).$ \\
      
\noi    Alors $F_1,...,F_p$ (resp. $F^*_1,...,F^*_p$) sont relativement invariants par $P(H_1,...,H_p)$ de caract\ères associ\és $\chi_k,k=1,...,p$ (resp. $\chi_{p-k}\chi_p^{-1},k=1,...,p$ en posant $\chi_0=1$) et sont les \irfs de l'action de $P(H_1,...,H_p)$ sur $\g_1$ (resp.$\g_{-1})$ (prop.1.4.5, lemme 1.4.7).\\
\noi Les orbites de $P(H_1,...,H_p)$ dans  $\g"_1=\{x\in \g_1\ |\ \prod_{i=1}^pF_i(x)\not=0\}$
 sont les orbites de $G_\goth t$ dans $W_{\goth t}=\{\sum_{i=1}^px_i\ |\ [x_i,x_j]=0$ pour $i,j=1,...,p$ et $x_i\in E'_2(H_i)\cap \g_1\}$, la situation est similaire dans $\g_{-1}$ (lemme 1.4.4).\\

 \noi Le \sgp tr\és sp\écial $P(H_p,...,H_1)(=\theta P(H_1,...,H_p) \theta^{-1}$ avec 
 $\theta=$ $ \prod_{i=1}^p\theta_{H_i}(-1)$ d\éfini dans l'introduction du $\S 1$) est associ\é \à l'ordre inverse et donne lieu \à la m\ême situation, pour $k=1,...,p-1:$\\

 $\bullet$ $P_{p-k}$  est l'\irf du \PV \air $(E_0( h_k)\cap \g_0,E_0(  h_k)\cap \g_1,H_0-\frac{1}{2} h_k)),$ que l'on \étend naturellement \à $\g_1,$\\
 
$\bullet$ $P^*_{k}$ celui du \PV  $(E_0( h_k)\cap \g_0,E_{-2}( h_k)\cap \g_{-1},\frac{1}{2}h_k),$  que l'on \étend naturellement \à   $\g_{-1}$.  \\
 
\noi Alors pour $1\≤k\≤p-1,$  $x\in E_2(  h_k)\cap \g_1$ et  $y\in E_0( h_k)\cap \g_1$ qui {\it commutent} on a $F_p(x+y)=F_k(x)P_{p-k}(y)$ (\à 2 exceptions pr\ès pour lesquelles il faut rajouter une puissance) (lemme 1.4.7).\\

 $\bullet$ \underline{Dans la section 2}, \étant donn\é un \PV\gog \   absolument irr\éductible et r\égulier  et $P_{\goth t}$ un \sgp standard tr\ès sp\écial, situation du $\S 1.4$ dont on reprend les notations, on \établit essentiellement le : \\
 
{\it \noi Th\'eor\`eme  2.1.1

\begin{enumerate}
 \item $P_{\goth t}$ a un nombre fini d'orbites dans
$\goth g_1$ et $\goth
g_{-1}.$

 \item Soient $S_{P_{\goth t}}=   \{ x\in {\goth g_{1}} \
 | \ \prod_{1\≤i\≤p}F_i(x) =0\}$ 
 et $S^*_{P_{\goth t}}
  = \{ x\in {\goth g_{-1}} \
 | \ \prod_{1\≤i\≤p}F^*_i(x) =0\}.$

\noi Pour tout \elt x de
$S_{P_{\goth t}}
  $  (resp. $S^*_{P_{\goth t}}$),
il existe $i\in \{1,...,p\}$  
tel que
$\chi_i /( (P_{\goth t})_x)^0\not=1$
 (resp.
$(\chi^*_i(/(P_{\goth t})_x)^0\not=1).$
\end{enumerate}}

\noi Ce th\'eor\`eme  \établit les conditions suffisantes d'existence des\eqs lorsque $\F$ est un corps $\goth p$-adique (th\éor\ème $k_{\goth p}$ de \cite{sato3}).\\

La \demo de ce th\'eor\`eme est faite par descente, en se ramenant \à des centralisateurs d'alg\èbres de type $sl_2$ et est commenc\ée dans le $\S2.1.$

\noi Elle n\écessite une forme relativement simple des repr\ésentants des orbites de $P_{\goth t}$ dans $S_{P_{\goth t}}$  obtenue par des consid\érations cas par cas. Ainsi on termine la \demo du th\'eor\`eme 2.1.1 cas par cas mais ceci nous am\ène \également \à obtenir une classification des \sgps standards tr\ès sp\éciaux et l'on retrouve ainsi certains exemples bien connus.\\

 \noi Le $\S 2.2$ contient quelques lemmes techniques g\én\éraux utiles pour les simplifications finales des repr\ésentants des orbites de $P_{\goth t}$ dans $S_{P_{\goth t}}$.\\
 
 \noi Le cas commutatif ($\g_2=\{0\}$) est trait\é dans le $\S 2.3, $ on y \établit \également que:  \\
 
{\noi \it Proposition 2.3.1
\noi Soit $N=exp(ad(\oplus_{\alpha\in
\∆_0^+}\ \goth g^{\alpha})).$ 

\noi Pour tout \elt x de $\goth
g_1$ (resp. $\goth g_{-1}$), il existe un ensemble S de
racines fortement orthogonales de
$\∆_1$  (resp. $\∆_{-1}$)  telles que
$N.x\cap\ (\oplus_{ \mu\in S
 }\goth g^{\mu})\ \not= \emptyset.$ }\\
 
 {\noi \it Lemme 2.3.2
 \begin{enumerate}
\item  Il existe un unique ensemble maximal ordonn\é de racines fortement orthogonales (longues) de
$\∆_1,$  $\lambda_1,...,\lambda_n,$ telles que le sous-groupe parabolique 

\noi
 $P_0=P( h_{\lambda_1},...,h_{\lambda_n})$ soit standard.

\item  Soit
$P_{\goth t}$ un sous-groupe parabolique  standard
 alors $P_{\goth t}\supset P_0$ et
  il existe  p, $ 2\≤p\≤n,$ et des
entiers:
$ l_0=0<1\≤l_1< ...<l_p=n$ tels que $P_{\goth
t}=P(H_1,...,H_p)$ avec  
$H_i=\sum_{l_{i-1}+1\≤j\≤l_i}h_{\lambda_j}$ pour $i=1,...,p.$
\end{enumerate}}
\vskip 2mm
\noi Les cas classiques sont trait\és dans le $\S 2.4.$

\noi Introduisons les notations des planches de II,III et IV de \cite{bourbakigal6}, $\∆$ est de type $ B_n,BC_n,C_n$ ou $D_n$ et $\lambda_0=\alpha_k= \epsilon_k- \epsilon_{k+1}$ alors
$\∆_1=\{\epsilon_i\pm \epsilon_j\ |\ 1\≤i\≤k<j\≤n\ ,\ \epsilon_i,1\≤i\≤k\}\cap \∆.$\\

\noi Soient: $p_0= \begin{cases}    \frac{k}{2} \  \textrm{ lorsque }  \goth g\textrm{ est d\éploy\ée de type  }C_n ,\\    k\    \textrm{ sinon}\    \end{cases}$
   et pour $i=1,...,p_0:$
$$H_i=\left\{ \begin{array}{ll}   2(h_{2\epsilon_{k-2i+1}}+h_{2\epsilon_{k-2i+2}})     \textrm{ lorsque }\  \goth g\textrm{ est d\éploy\ée de type  }C_n ,\\ 
   h_{\epsilon_{k-i+1}}\    \textrm{ sinon.}\    \end{array}\right.$$
  \noi Lorsque $2k\≤n,$  on introduit \également 
$H'_1=\sum_{i=1} ^kh_{\epsilon_i-\epsilon_{2k-i+1}}$  et
$H'_2=\sum_{i=1}^kh_{\epsilon_i+\epsilon_{2k-i+1}}.$\\

\noi  Les \sgs paraboliques $P_0=P(H_1,...,H_{p_0})$ et $P_0' =P(H'_1,H'_2)$  sont des \sgs 
 paraboliques standards tr\ès sp\éciaux et si 
  $P=P(H_1,...,H_p)$ est un \sg  parabolique standard tr\ès sp\écial alors 

$\bullet$  soit  $P=P'_0$ et $2k\≤n,$

$\bullet$  soit $P\supset P_0,$   $k\≥2$  et  il existe des
entiers:
$ l_{p+1}=0<l_p< ...<l_2<l_1=k$ tels que  $H_i=\sum_{l_{i+1}+1\≤j\≤l_i}h_{\epsilon_j}$ pour $i=1,...,p$ (prop.2.4.3).\\

\noi Les cas exceptionnels non commutatifs sont trait\és dans le $\S 2.5$ et dans tous ces cas il y a un unique sous-groupe parabolique standard tr\ès sp\écial,   $P(H_1,H_2)$ (cf.tableau 3).

\noi Lorsque $\g_2$ n'est pas de dimension $1,$   $H_2=2h_{\tilde \alpha},$   $\tilde \alpha$ \étant la plus grande racine de $\Delta.$
  
   \noi $P(H_1,H_2)$ est un \sg parabolique maximal de G \à l'exception du cas $(E_6,\alpha_2)$ pour lequel $\∑_0-\∑_{\goth t}=\{\alpha_1,\alpha_6\}$ et dans les autres cas, on a:
 
\noi  $\∑-\∑_{\goth t}= \left \{
    \begin{array}{l} \{\alpha_1,\alpha_4\}\ \hbox{lorsque }\ \∆=F_4 ,\\ 
  \{\alpha_1,\alpha_6\}\ \hbox{lorsque }\ \∆=E_7\ \hbox{et}\ \∑_1\subset   \{\alpha_1,\alpha_6\},\\ 
   \{\alpha_1,\alpha_8\}\ \hbox{lorsque }\ \∆=E_8 ,\\ 
   \{\alpha_1,\alpha_2\}\ \hbox {dans le cas  }\ (E_7,  \alpha_2).\\  \end{array} \right.    
$\\
\noi De plus, lorsque $(\∆,\lambda_0)$ est de type   $(E_6,\alpha_2),(E_7,\alpha_1),(E_8,\alpha_8),$  pour tout \elt $x$ de $\goth g_1$  il existe  $y$ dans $ P(H_1,H_2).x$ et  4 racines fortement orthogonales $\beta_1,...,\beta_4,$ tels que le support de $y$ comprenne au plus 4 racines appartenant \à $\{\tilde \alpha-\beta_1,\beta_i,1=1,...,4\}$ et dans le cas $(F_4,\alpha_1)$ il faut y adjoindre l'ensemble $\{\tilde \alpha-\lambda_1,\frac{1}{2}(\lambda_1+\lambda_3),\frac{1}{2}(\lambda_2+\lambda_4)\}.$  Ces $4$ cas correspondent \à $\g_2$ de dimension $1.$

\bigskip
 $\bullet$ \underline{La section 3} est consacr\ée aux r\ésultats classiques sur les fonctions Z\éta associ\ées aux \PVs\   absolument irr\éductibles et r\éguliers  $(P_{\goth t}=P(H_1,...,H_p),\g_{\pm 1})$ ce qui impose des normalisations: 
 
 $\bullet$ de la forme de Killing $B$ (\S 3.1, $\widetilde B=-\ds\frac{\text{degr\é de}F_p}{2B(H_0,H_0)}.B$), 
 
 $\bullet$ des mesures de Haar de $\g_1$ et $\g_{-1}$(\S 3.2), 
 
 $\bullet$ des \irfs $F_1,$...,$F_p,$ $F^*_1,$ ...,$F^*_p$ (\S 3.3).\\
 
Dans le cas archim\édien , on introduit les op\érateurs diff\érentiels $F_k(\partial),F^*_k(\partial),k=1,...,p$ ainsi que les polyn\ômes de Bernstein associ\és $b_k$ et $b^*_k,k=1,...,p$:

\noi pour $s=(s_1,...,s_p)\in \C^p$ et $F^s=\prod_{1\≤i\≤p}F_i^{s_i},$ 
$F^{*s}=\prod_{1\≤i\≤p}F_i^{*s_i},$ on a:
$$F_k(\partial)(F^{*s})=b_k(s){F^*}^{s-\bf 1_p+\bf 1_{p-k}}\ ,\ F^*_k(\partial)(F^{s})=b^*_k(s){F}^{s-\bf 1_p+\bf 1_{p-k}},$$
avec $\bf 1_0=(0,...,0)$ et $\bf 1_j$ \étant le vecteur de $\R^p$ dont toutes les coordonn\ées sont nulles sauf la j\ème qui vaut $1,$ avec un am\énagement pour $2$ exceptions (lemme 3.4.2). On pose $b_{\goth g,P_{\goth t}}:=b_p.$

\noi Lorsque $k$ est diff\érent de $p,$ comme dans le cas commutatif (\cite{rubschiff}), on \établit que $b_k(s_1,...,s_p)$ est proportionnel au polyn\ôme de Bernstein, $b_{\goth U,P(H_1,...,H_k)}(s_{p-k+1},...,s_p),$  $\goth U$ \étant le centralisateur d'une alg\èbre de type $sl_2$ telle que $\goth U_1\subset E_2(h_k)\cap \g_1$ et $(\goth U_0,\goth U_1)$ admet $F_k/\goth U_1$ comme invariant relatif fondamental (ainsi que la propri\ét\é analogue pour $b^*_k$) (prop.3.4.4).\\

\noi Le th\éor\ème 2.1.1 ainsi que les propri\ét\és de $F_k(\partial)$ et 
 $F^*_k(\partial),$ pour $k=1,...,p,$ dans le cas archim\édien, permettent d'obtenir le prolongement analytique et l'existence des\eqs des fonctions z\étas associ\ées \à $(P(H_1,...,H_p),\g_{\pm 1})$ (th.3.5.2 qui s'appuie sur \cite{sato3},\cite{boppruben1},\cite{boppruben}, cf. \également le lemme 3.5.3 pour l'ind\épendance des coefficients).
 
 \noi En application imm\édiate, pour $k$ strictement inf\érieur \à $p,$ on obtient   $b_{\g,P_{\goth t}}$  comme  produit des 2 polyn\ômes de bernstein, 
 
\noi $b_{\goth U,P(H_1,...,H_k)}(s_{p-k+1},...,s_p)$ et  $b_{\goth U',P(H_{k+1},...,H_p)}(s_1,...,s_{p-k-1}, \sum_{p-k\≤i\≤p}s_i+r_k),$ 

\noi$\goth U'$ \étant \à nouveau le centralisateur d'une alg\èbre de type $sl_2$ telle que $\goth U'_1\subset E_0(h_k)\cap \g_1$ et  $r_k$ \étant  un nombre rationnel explicitement d\éfini \à partir de dimensions (prop.3.7.3).\\
 
 \noi On rappelle dans le $\S 3.6$ les $2$ exemples fondamentaux utilis\és ult\érieurement:
 
$\bullet$ lorsque $\g_1$ est de dimension $1,$ r\ésultat d\û \à J.Tate(\cite{tate})  faisant intervenir les facteurs $\rho'(\pi) $ et $\rho'(\pi;x),$ $\pi\in \widehat{\F^*}$ et $x\in \FF$,\\
    
$\bullet$ lorsque l'\irf de $(\g_0,\g_1)$ est une forme quadratique $F$,  les r\ésultats de \cite{rallisschiffmann}) ont \ét\é reexprim\és (th.3.6.5, notamment le 5)) et font intervenir les facteurs $A^a_{\pi_1,\pi_2}(v,u,\delta),$ $\pi_1$ et $\pi_2\in \widehat{\F^*},$  $(v,u)\in (\FF)^2,$ $\delta=(-1)^{[N]}\text{disc}(F)$ (B du lemme 3.6.4,cf.\également lemmes 3.6.7 et 3.6.8).
    
    \bigskip
 $\bullet$ Dans \underline{ la section 4}  on \établit une d\écomposition des mesures sur $\g_1$ et $\g_{-1}$ ``adapt\ée"  \à l'action du \sgp $P(h_k,2H_0-h_k)$  et par cons\équent aux invariants relatifs  $F_k$ et $F^*_{p-k},$ $k$ \étant un entier fix\é tel que $1\≤k\≤p-1,$ (th.4.3.3) et qui s'exprime \à l'aide des mesures $G_{h_k}$-invariante \à gauche sur 
 
 \noi $W_k=\{x+y\ |\ x\in E_{2}(h_k)\cap \g_1\ ,\ y\in E_{0}(h_k)\cap \g_1\ ,\ [x,y]=0,F(x+y)\not=0\} $ et 
 
 \noi $W^*_k=\{x'+y'\ |\ x'\in E_{-2}(h_k)\cap \g_{-1}\ ,\ y'\in E_{0}(h_k)\cap \g_{-1}\ ,\ [x',y']=0,F^*(x'+y')\not=0\} $
 
 \noi (corollaires 4.4.2 et 4.4.3).
 
 \noi Avec cette d\écomposition, la transformation de Fourier se d\écompose en 2 transformations de Fourier partielles (th.4.3.5 dans lequel intervient le r\ésultat de A.Weil, \cite{weil}, sur la transformation de Fourier d'un caract\ère quadratique, cf.lemme 4.2.3).
 
 \noi Ces r\ésultats imposent des normalisations pr\écises des mesures sur les sous-espaces vectoriels $E_i(h_k)\cap E_j(2H_0-h_k)$ ($\S 4.2$).
   \bigskip
 
 $\bullet$ Dans \underline{la section 5}  on applique les r\ésultats de la section 4 au calcul des coefficients de l'\eq v\érifi\ée par les fonctions Z\étas. On obtient ainsi une relation explicite les exprimant \à partir de ceux associ\és \à des \PVs construits  avec des centralisateurs d'alg\èbres de Lie de type $sl_2$ (prop.5.1.1). 
 
 \noi  Ceci donne l'\équation (A) du cas complexe (th.5.2.2 et remarque 5.2.3).

\noi On \établit \également, sous certaines conditions v\érifi\ées dans la plupart des cas, une relation permettant le calcul explicite des coefficients $a_{v,u}$ (cf.relation (B) de l'introduction) (prop.5.3.2).

\noi Cette section se termine par l'\étude de  2 exemples r\éels particuli\èrement simples puisque la descente fait apparaitre $2$ formes quadratiques anisotropes dont les r\ésultats sont connus ($(E_6,\alpha_2),$  $E_6$ \étant de type III et $(E_7,\alpha_1)$ avec $E_7$ de type VII, prop.5.3.3).\\

 Les sections suivantes sont les applications des r\ésultats de la section 5.  Dans chaque cas on d\étermine les polyn\ômes de Bernstein, les coefficients de l'\eq \à multi-indice (c'est \à dire associ\ée aux   \irfs $F_1,...,F_p,F^*_1,...F^*_p$) et de l'\eq simple (c'est \à dire associ\ée aux aux $2$ \irfs $F,F^*$), la m\éthode employ\ée fait apparaitre ces coefficients sous forme de somme de produits de coefficients d'\eqs associ\ées \à des \PVs de rang plus petit ce qui impose:
 
 - la connaissance assez pr\écise de la structure et des r\ésultats pour ces derniers,
 
 - la simplification de ces sommes lorsqu'on passe \à l'\eq simple et qui est obtenue en utilisant les r\ésultats du $\S 3.6.2.$\\
 
 \noi Ceci est illustr\é dans les:\\
 
 $\bullet$ {\underline {section 6}} qui traite du cas commutatif (cf.tableaux 1 et 2) avec le \sgp $P_0,$  ainsi qu'un  approfondissement du cas ``symplectique"   commutatif (i.e. $\overline\Delta=C_n$) dans le cas $\goth p$-adique de caract\éristique r\ésiduelle diff\érente de $2,$\\
 
$\bullet$ {\underline {section 7}} qui traite des cas classiques (cf.tableau 2) avec l'action du \sgp $P_0.$ Une attention particuli\ère a \ét\é donn\ée au cas orthogonal BDI, malheureusement les r\ésultats sont incomplets dans le cas DIII $\goth p$-adique,\\
 
$\bullet$ {\underline {section 8}} qui traite des cas exceptionnels (cf.tableau 3) dont l'\étude est emboit\ée.\\

$\bullet$  {\underline {L'appendice 1}} contient un r\ésultat \él\émentaire  sur le rang du centralisateur d'un \Sl. Ce  r\ésultat est utilis\é dans l'\étude de $(E_7,\alpha_6).$\\

 $\bullet$  {\underline {Appendice 2:}}  Dans le cas commutatif ou bien, {\it lorsque $g_2$ est de dimension $1$ et $\g$ de rang au moins $4,$} on exprime les mesures sur $S$ et $S^*,$ invariantes par le noyau de $\chi$   \à l'aide de  la fonction Z\éta provenant d'un \PV inclus dans celui de d\épart (prop.1). Elles sont relativement invariantes par un sous-groupe d'indice fini de $G$ et temp\ér\ées.
 
 \noi Cette description permet dans le cas commutatif, lorsque $\overline d$ est pair, le calcul explicite de la transform\ée de Fourier de ces mesures pour les fonctions de $\EuScript S(\g_1-S)$  en appliquant les $\S 4.3$ et $6.2.$

    \newpage
    
     \centerline {\large \bf PLAN }

\bigskip

\noi {\bf {\large 1 }   Une classe d'espaces pr\éhomog\ènes de type parabolique faiblement sph\ériques} \\
 
\noi 1.1 \'El\ément 1-simple (tr\ès) sp\écial \\
  
\noi 1.2 Quelques propri\ét\és  \\
 
\noi 1.3 Existence  \\

\noi 1.4 Une classe de \PVs faiblement sph\ériques  \\ \\
 
 \noi {\bf {\large 2 }   Un r\ésultat sur les orbites des \sgs paraboliques tr\ès sp\éciaux  }\\

\noi 2.1 Le r\ésultat \\

\noi 2.2 Lemmes techniques ($\F$ est de caract\éristique $0$)\\

\noi 2.3 Le cas commutatif ($\goth g_2=\{0\}$)\\

\noi 2.4 Les cas classiques\\

\noi 2.5 Les cas exceptionnels non commutatifs

\noi  2.5.1  G\én\éralit\és

\noi  2.5.2  Les cas exceptionnels  pour lesquels $\g_2$ est de dimension 1

\noi 2.5.3 Les cas exceptionnels restants\\ \\

 \noi {\bf {\large 3 }  Fonctions Z\'etas associ\ées   }\\

 \noi  3.1 Normalisation de  la forme de Killing  \\

 \noi 3.2 Normalisation des mesures de $\goth g_1$ et
$\goth g_{-1}$  \\

 \noi 3.3 Normalisation des invariants relatifs fondamentaux  \\

 \noi 3.4 Le cas  archim\édien  \\
  
 \noi 3.5 Fonctions Z\étas: d\éfinition et \équations
fonctionnelles abstraites  \\

 \noi 3.6 Deux exemples fondamentaux 
 
 \noi 3.6.1 $\g_1$ est de dimension 1
 
 \noi 3.6.2  L'\irf est une forme quadratique\\
 
 \noi 3.7 Application au cas archim\édien\\ \\
 
  \noi {\bf {\large 4 } D\écomposition des mesures sur $\goth g_1$ et
$\goth g_{-1}$ et application }  \\
  
 \noi 4.1 Notations\\
   
  \noi 4.2 Mesures sur les sous espaces $E_{i,j}$ et $E_{-i,-j}$ pour $(i,j)\not=(0,0)$  \\
    
  \noi 4.3 D\écomposition des mesures sur $\goth g_1$ et
$\goth g_{-1}$  \\
 
  \noi 4.4 Expression des  mesures $\bf G_{h_k}$
invariantes sur $W_k$ et $W^*_k$  \\   \\ 
     
   \noi {\bf {\large 5 } Application aux fonctions Z\étas}\\

   \noi 5.1 Les coefficients de l'\équation fonctionnelle  \\

   \noi 5.2 Le cas complexe \\
   
 \noi 5.3 Un cas particulier  \\ \\
 
 \noi {\bf {\large 6 } Le cas commutatif}\\

  \noi 6.1 Structure
  
   \noi 6.1.1 Rappels
   
 \noi 6.1.2 Le type III

  \noi 6.1.3 Le cas $(C_{2n},\alpha_{2n})$ de type III

  \noi 6.1.4 Calcul de $\gamma_k$\\
 
  \noi 6.2 \'Equations fonctionnelles

  \noi 6.2.1 Le cas transitif
  
   \noi 6.2.2  Le cas non transitif\\ \\

  \noi {\bf {\large 7 } Les cas classiques}\\

  \noi  7.1 G\én\éralit\és

  \noi  7.1.1 Description des cas classiques consid\ér\és
 
  \noi  7.1.2 Sous-groupe parabolique standard tr\ès sp\écial et \PV associ\é\\
 
  \noi  7.2 Le cas symplectique\\
 
  \noi  7.3 Les cas
orthogonaux de type I  

   \noi  7.3.1 Normalisation des invariants relatifs

   \noi  7.3.2 Le r\ésultat ($\bf 3k\≤2(n-1)+\delta$)

   \noi  7.3.3 Le cas r\éel  

   \noi  7.3.4 Le cas $\goth p$-adique\\

   \noi  7.4 Le cas DIII

  \noi  7.4.1 Le cas $p=1$

  \noi  7.4.2 Le cas $p\≥2$\\ 

  \noi {\bf {\large 8 } Les cas exceptionnels}\\

 \noi  8.1 Le cas $(E_7,\alpha_6)$
 
 \noi  8.1.1 G\én\éralit\és
 
  \noi  8.1.2 Pr\éliminaires
  
 \noi  8.1.3 Le r\ésultat\\

  \noi  8.2 $\goth g_2$ de dimension $1$
  
   \noi  8.2.1 Structure

  \noi  8.2.2 Une premi\ère \équation fonctionnelle

  \noi  8.2.3 Orbites

   \noi  8.2.4 Le cas r\éel de rang $4$ $(\∆=F_4$ et $f$ est anisotrope)\\
 
  \noi  8.3 $ (E_7,\alpha_2)$\\
 
  \noi  8.4 $ (E_8,\alpha_1)$ 
 
   \noi  8.4.1 Le cas d\éploy\é

  \noi  8.4.2 Le cas r\éel non d\éploy\é\\
  
  \noi {\bf Appendice 1}\\
  
   \noi {\bf Appendice 2: Remarques sur les mesures relativement invariantes sur les orbites singuli\ères}\\
   
 \noi  {\bf Tableau 1 ($\bf \S 6$)}\\
 
  \noi  {\bf Tableau 2 ($\bf \S 7$)}\\
  
   \noi  {\bf Tableau 3 ($\bf \S 8$)}\\
  
\newpage

\section {\bf  Une classe d'espaces pr\éhomog\ènes de type parabolique faiblement sph\érique }

\bigskip
On introduit les \PVs de type  paraboliques consid\ér\és dans ce travail de mani\ère inhabituelle, \à partir de l'action de \sgps associ\és \à des \elts $1$-simples ayant de bonnes propri\ét\és relativement \à la graduation de $\g$ et qui nous semblent \être une g\én\éralisation du cas commutatif.

Dans  ce chapitre, sauf mention du contraire, on suppose que  $\g$ est une
alg\`ebre de Lie d\'efinie sur $\F$, semi-simple, de
dimension finie et gradu\ée:  
 $$\goth {g}=\oplus _{i\in {\Z}} \goth {g}_i\quad ([\goth
{g}_i,
\goth {g}_j]\subset \goth {g}_{i+j}),$$ 
on note $H_0,$ l'\elt appartenant \`a
$\goth {g}_0,$ qui d\'efinit la graduation et on suppose que $2H_0$ est {\it $1$-simple}, par cons\équent $\g_1$ et $\g_{-1}$ engendrent $\g.$

\noi Certains r\ésultats sont vrais sur un corps  de  caract\'eristique $0$.  \\

Soit $\goth a$ une sous-alg\èbre ab\'elienne
d\'eploy\'ee maximale de $\goth g$ contenant $H_0,$ on note
$\Delta$ le syst\ème de racines correspondant qui est
\'egalement gradu\'e par $H_0:$ 
$$\Delta_i=\{\lambda\in \Delta\ |\   \lambda(H_0)=i\}.$$ On choisit un ordre
sur $\Delta$ compatible avec la graduation, c'est \à dire que  $\cup_{i>0}\Delta_i \subset \Delta^+,$   $\Sigma$ \étant un ensemble de racines simples,  le \PV est de type $(\Delta,\Sigma_1=\Sigma\cap \Delta_1 ).$  \\

Introduisons quelques notations suppl\émentaires:\\
\begin{enumerate}

\item On note $B$ la forme de Killing de $\g.$\\

\item $\overline \F$ est une cl\ôture alg\ébrique de $\F$ et les \él\éments correspondants sont not\és avec une barre, par exemple $\overline{\goth g}=\goth g\otimes_{\overline \F}\F.$

\noi $\overline{\goth h}$  d\ésigne une sous-alg\èbre de Cartan de $\overline{\goth g}$ contenant $\goth a$, $\overline \∆$ le syst\ème de racines de $(\overline{\goth g},\overline{\goth h})$ muni de l'ordre habituel et,  lorsque $\overline \∆$ est irr\éductible, on note $\tilde \omega$ la plus grande racine .\\
\noi $(\overline \∆,\overline \∑_1)$ sera le diagramme de Dynkin gradu\é  du \PV $(\overline{\goth g}_0,\overline{\goth g}_1)$ et par abus \également du \PV  $((\goth g\otimes_{\E}\F)_0,(\goth g\otimes_{\E}\F)_1)$, $\E$ \étant une extension galoisienne de $\F$ pour laquelle l'alg\èbre $\goth g\otimes_{\E}\F$ est d\éployable.\\ 

$G$ est le centralisateur de $H_0$ dans le groupe $Aut_0(\goth
g)$ des automorphismes de $\goth g$ qui sont \'el\'ementaires sur  
 $\overline{\F}.$  \\
 
 \item Lorsque $\Delta$ est irr\éductible, 
la plus grande racine du syst\ème de racines r\éduit associ\é \à
$\Delta$ est not\ée  $\tilde \alpha.$\\ 

 \item Si $\goth c$ est une sous-alg\èbre  de $\goth g,$ on
d\éfinit
$\goth c_i=\goth c\cap \goth g_i.$\\
 
 \item Soit  $\goth c$ une sous-alg\èbre de $\goth g,$
r\éductive dans
$\goth g$ alors le centralisateur de $\goth c$  dans $\goth g,$
not\é $E_0(\goth c),$ est une alg\èbre r\éductive dans $\goth g$ (\cite{bourbakigal8}) dont
on note  $[E_0(\goth c),E_0(\goth c)]$ sa partie semi-simple.

\noi Lorsque $H_0$ est un \él\ément de la sous-alg\èbre $\goth c+E_0(\goth
c),$ $\goth c$ est gradu\ée : $\goth c=\oplus_{i\in \Z}\goth
c_i,$ $E_0(\goth c)$ l'est \également par $ad(H_0)$ et son centre est inclus dans $\goth g_0$ d'o\ù l'alg\èbre d\ériv\ée $ [E_0(\goth c),E_0(\goth c)]$ est \également gradu\ée : $\ 
 [E_0(\goth c),E_0(\goth c)]=\oplus_{i\in \Z}  [E_0(\goth c),E_0(\goth c)]_i$ et de plus pour $i$
non nul on a $   [E_0(\goth c),E_0(\goth c)]_i\subset E_0(\goth c)_i;$ soit $\goth U(\goth c)$ l'alg\èbre engendr\ée par $E_0(\goth c)_{\pm i}$ pour $i\geq 1,$ on a  $\goth U(\goth c)_{\pm i}=E_0(\goth c)_{\pm i}$ pour $i\geq 1.$  $\goth U(\goth c)$ \étant un $ [E_0(\goth c),E_0(\goth c)]-$module, $\goth U(\goth c)$ est un id\éal de $ [E_0(\goth c),E_0(\goth c)]$ donc $\goth U(\goth c)$ est une alg\èbre de Lie semi-simple gradu\ée : $\goth U(\goth c)=\oplus_{i\in \Z}  \goth U(\goth c)_i$  de \PV associ\é:  $(\goth
U(\goth c)_0,\goth
U(\goth c)_1)$ avec la convention habituelle lorsque $(\goth U(\goth c)_{\pm 1})$ engendre une alg\èbre semi-simple que $\goth U(\goth c)$ soit cette alg\èbre engendr\ée.\\  
Lorsque $\goth c=\F H,$  on le note \également $(E_0(H)_0,E_0( H)_1).$\\

\item    $P(\check {\goth g})$ d\ésigne l'ensemble $\{u\in
\goth g_0\ |\ ad(u)$ est diagonalisable \à valeurs propres
enti\ères$\},$  pour $i\in \Z$ on note $E_i(u)=\{x\in \g\ |[u,x]=ix\}$  et pour $t\in \F^*,$  $h_u(t)$ est l'\él\ément de
$G$ correspondant  (c'est \à dire d\éfini pour  $i\in \Z$ par $h_u(t)/_{
E_i(u)}=t^iid_{  E_i(u)}$).\\

 \item  Lorsque $(x,h,y)$ est un $sl_2-$triplet  on rappelle que  $H_0-\frac{h}{2}\in E_0( \F x+\F h+\F y),$ que $E_0( \F x+\F h+\F y)=E_0( \F x+\F h)$ et que l'\él\ément $y$ est
determin\é de mani\ère unique par $x$ et $h,$  $y$ est not\é parfois $x^{-1}$. \\ 

\noi  $\theta_{x,h}(t)$   (ou
bien simplement
$\theta_x(t)$   lorsqu'il n'y a pas
d'ambiguit\é ou $\theta_h(t)$) d\ésigne l'automorphisme \él\émentaire :
$$t\in \F^*\ exp(ad(tx))exp(ad(t^{-1}y)exp(ad(tx))\ (\text{alors }h_h(t)=\theta_{x,h}(t)\theta_{x,h}(-1)) .$$

\noi Notons que $B(2H_0-h,2H_0)=B(2H_0-h,2H_0-h)$ car
$B(2H_0-h, h)=$ $B(\theta_{x,h}(t)(2H_0-h),\theta_{x,h}(t)(
h))=$ $ -B(2H_0-h, h).$\\

\item  Lorsque $\alpha\in \Delta,$  $\g^\alpha$ est le sous-espace radiciel correspondant et $h_\alpha$ la co-racine ;   on  
note  simplement $h_{\alpha}(t)$ (resp.$\theta_{ \alpha} $)  l'\él\ément
$h_{h_{\alpha}}(t)$ (resp.$\theta_{h_{\alpha}}(-1)$) et on rappelle que $\theta_{ \alpha} (\goth g^{\beta})=\goth g^{s_{\alpha} (\beta)}$ pour toute racine $\beta $ de $\Delta.$
  
\end{enumerate}
\bigskip

 \subsection{\bf  \'El\ément  1-simple (tr\ès) sp\'ecial}

\bigskip
On commence par motiver les d\éfinitions sur les \elts $1$-simples par un lemme.\\

\noi  A un \él\ément $h$ $1$-simple, on peut associer  les
sous-alg\èbres paraboliques de $\goth g$ et $\g_0$:  
$$p(h)=\oplus_{i\geq 0}E_i(h)=E_0(h)\oplus n(h)\ ,\
   p_0(h)=p(h)\cap \goth g_0$$
$ P(h)=G_hexp(ad(\oplus_{i\≥1}E_i(h)\cap \g_0)$ le sous-groupe parabolique
de $G$  d'alg\èbre de Lie $p_0(h),$ ainsi que le \PV $(E_0(h)\cap
\goth g_0,E_2(h)\cap \goth g_1,{h\over 2})$  construit avec l'alg\èbre r\éductive   
 $\oplus_{i\in \Z}E_{2i}(h)$   
 dont le centre est inclus dans  $E_0(h)\cap \goth
g_0.$    \\
 
\noi On note $\pi_h$  la projection de
$\goth g_1$ sur $E_2(h)\cap \goth g_1$ 
 parall\èlement \à $ \oplus_{i\not=2} E_i(h)\cap
\goth g_1.$\\

\begin{lem}  
\begin{enumerate}
 \item Soit h un \elt 1-simple distinct de $2H_0$ tel que $ n(2H_0)\subset  p(h) $  alors la sous-vari\ét\é
  $$\tilde {W_h}=\{x+y\in E'_2 (h)\cap \goth g_1 \oplus   E_0
(h)\cap \goth g_1 \ |\  [x,y]=0\}$$   rencontre  $\goth
g'_1$ c'est \à dire que  $2H_0-h$  est aussi $1$-simple et
$\tilde{W}_h\cap
\goth g'_1=\tilde{W}_h\cap\tilde{W}_{2H_0-h}.$\\ 
 \item  Si $ n(2H_0)\subset p(h)\cap p(2H_0-h)$ alors  $P(h)$ et $P(2H_0-h)$ ont une orbite dense dans $\goth
g'_1$ et $\goth g'_{-1}$.\\ 
\item L'application $\Psi$ d\éfinie sur
$n_0(2H_0-h) \times \tilde{W}_h$
  par  $\Psi(A, w)=exp(ad(A)w)$ est 
 un hom\éomorphisme (diff\éomorphisme  dans le cas r\éel)  sur
un ouvert de Zariski   de $\goth g_1$  et
    $\Psi(n_0(2H_0-h) \times  W_h)=\Psi(n_0(2H_0-h)\times \tilde{W}_h)\cap
\goth g'_1$ avec $W_h=\tilde{W}_h\cap
 \tilde{W}_{2H_0-h}.$\\
\item Soit $x\in E'_2(h)\cap \goth g_1$   alors $\pi_h^{-1}(x)$
rencontre  
$\Psi(n_0(2H_0-h)
\times
  (\pi_h^{-1}(x)\cap \tilde{W}_h)).$
 \end{enumerate}
 \end{lem}
 
 \dem

\begin{enumerate}

 \item Soit $x+y$ un \él\ément de $ \tilde W_h,$ $\goth s$ la sous-alg\èbre
engendr\ée par $x,h,x^{-1},$ $z$ un \él\ément non singulier
du \PV $(\goth U(\goth s)_0 , \goth U(\goth s)_1).$ Il
suffit de v\érifier que $ad(x+z)$ est une surjection de
$\goth g_0$ sur $\goth g_1.$\\

\noi Soit $p_0$ le plus grand entier tel que le sous-espace
$E_{p_0}(h)\cap \goth g_1$ ne soit pas r\éduit \à $\{0\}$, on
v\érifie alors ais\ément que pour tout $y $ dans
$\goth g_1,$ l'\équation
$[X,x+z]=y$ a au moins une solution dans $E=\oplus_{
i=-2}^{i=p_0-2}E_{-i}(h)\cap
\goth g_0,$ en utilisant les propri\ét\és usuelles des
$sl_2-$triplets ainsi que la d\écomposition :
$E_0(h)=[x,E_{-2}(h)]\oplus E_0(\goth s).$\\

\noi Compl\étons $z$ en un $sl_2-$triplet $1-$adapt\é,   $(z,u,v),$ au sens du
\PV  $ (\goth U(\goth s)_0 , \goth U(\goth s)_1),$
  en raison des relations de commutation le triplet 
$(z+x,h+u,v+x^{-1})$ est encore un $sl_2-$triplet $1-$adapt\é
d'o\ù $h+u=2H_0.$\\
\item Remarquons que les conditions impos\ées donnent  l'\égalit\é: 
$$\goth g_1=\oplus_{i\geq 0,2-i\geq 0}E_i(h)\cap
E_{2-i}(2H_0-h) \quad \hbox{d'o\ù}\quad E\subset n_0(2H_0-h).$$
Il suffit d'appliquer le r\ésultat d\émontr\é en 1) pour
montrer que  $P(2H_0-h)$ a une orbite dense dans $\goth
g'_1$. Comme $h$ et $2H_0-h$ ont le m\ême r\ôle,  $P(h)$  a
\également une orbite dense dans $\goth g'_1.$\\
\item Pour le dernier point, r\ésolvons l'\équation :
$$\Psi(A_1+A_2,x+y)=x_2+x_0+x_1\quad \hbox{avec}\quad x+y\in
\tilde W_h \ ,\ A_i\in E_{-i}(h)\cap \goth g_0\ ,\ 
i=1,2$$  
et  $x_2+x_0+x_1\in  E'_2(h)\cap \goth g_1\oplus
E_0(h)\cap \goth g_1\oplus E_1(h)\cap \goth g_1,$ par l'hypoth\èse   :
$$\Psi(A,x+y)=x+[A_1,x]+(y+\frac{1 }{ 2}ad(A_1)^2(x)+[A_2,x])$$
\noi Ainsi il suffit de r\ésoudre les
\équations :
$$ x=x_2 \ ,\
[A_1,x]=x_1\ ,\
y=x_0-{1\over 2} ad(A_1)^2(x)-[A_2,x]\quad \hbox{et}\quad
[y,x]=0.$$
\noi Notons que:
$$[x_2,ad(A_1)^2(x_2)]=[x_2,[A_1,x_1]]=-[x_2,[x_1,A_1]]=
-[x_1,[x_2,A_1]]=[x_1,x_1]=0,$$
\noi 
car $ad(x_1)$ et  $ad(x_2)$ commutent 
$(E_{-3}(h)\cap \goth g_{-1}=\{0\}$  et est en bijection par
$ad(x_2)^3$ avec $E_3(h)\cap \goth g_2$).\\

\noi On obtient :
$$
 A_1=[ x_2^{-1},x_1]\ ,\ A_2={1\over
2}[x_2^{-1},x_0]\ ,x=x_2\ ,y=x_0-{1\over
2}(ad(x_1))^2(x_2^{-1})+ {1\over 2}ad(x_2^{-1})ad(x_2)(x_0). $$
 \fdem\\
 \end{enumerate}
 
 \begin{defi} 
 \begin{enumerate}
 \item Un \él\ément h, 1-simple, tel que $ n(2H_0)\subset p(h)  $ est appel\é 1-simple sp\écial.
 \item Un \él\ément h tel que h et $2H_0-h$ sont 1-simples sp\éciaux est appel\é tr\ès sp\écial.
 \end{enumerate}
 \end{defi}
 
 \bigskip
 
  \begin{rema} 
   \begin{enumerate}

  \item h est 1-simple sp\écial \ssi  le spectre de $ad(h)/\goth g_1$ est inclus dans $ \N.$
   \item h est 1-simple tr\ès sp\écial \ssi  $\{0,2\}\subset $ le spectre de $ad(h)/\goth g_1$   $\subset \{0,1,2\}.$
   \item Lorsque $\goth g_i=\{0\}$ pour $i\≥3,$  h est 1-simple sp\écial \ssi  le spectre de $ad(h)/\goth g_2$ est inclus dans $2 \Z.$ En particulier, dans les \PVs commutatifs (tr\ès r\éguliers) tous les \él\éments $1$-simples (distincts de  $2H_0$)
   sont  (tr\ès) sp\éciaux.   \end{enumerate}
 \end{rema}
 
 \bigskip

\noi  On montre que les   \él\éments $1$-simples tr\ès sp\'eciaux  n'apparaissent que dans des graduations courtes et induisent une ``conservation" des propri\ét\és initiales.  \\ \\

\subsection{\bf   Quelques propri\ét\és}

\bigskip

  \begin{lem} 
   Lorsque:
    \begin{enumerate}
  \item $(\goth g_0,\goth g_1,H_0)$ est un espace pr\éhomog\ène absolument irr\éductible et r\égulier,
  \item $H\in \goth a$ est 
   un \él\ément 1-simple  sp\écial distinct de $2H_0,$
   \end{enumerate}  \noi Le
pr\éhomog\ène $( E_0( H)_0,\ E_0(H)_1,H_0-\frac{H}{2})$ est \également  absolument irr\éductible et r\égulier et si $\chi_H$  d\ésigne
le caract\ère associ\é \à l'invariant relatif fondamental (de
degr\é $d_H$) on a pour
$u\in P(\check {\goth g})$ et commutant \à
H: $$\chi_H(h_u(t))=t^a\ \hbox{ avec}\ \
a= d_H. \frac{2 B(u,2H_0-H)}{B(2H_0-H,2H_0-H)}.$$
 \end{lem}
 
 \dem
 
 \bigskip
  
\noi Comme l'ensemble : $\{\alpha\in \Delta_0\ |\ \alpha(H)\geq 0\}$ est
une partie parabolique de $\∆_0,$ il existe un ordre pour
lequel $\∆_0^+\subset \{\alpha\in \∆_0\ |\ \alpha(H)\geq 0\}$ et  $\∆^+=\∆_0^+\cup_{i\geq 1}\∆_i.$
\noi Soit $\Sigma=\Sigma_0\cup \Sigma_1$ les racines simples correspondantes
   alors $\Sigma_1=\{\beta_0\}$ et $\beta_0(H)=0$  (en effet dans le cas contraire on a  
$E_0(H)\cap \goth g_1=\{0\}$  ce qui est absurde d'apr\ès  le
le lemme pr\éc\édent) d'o\ù l'\égalit\é  $\{\alpha\in \∆_{±1}\ |\
\alpha(H)=0\}=\∆_{±1}\cap (\oplus _{\{\alpha\in \Sigma\ |\  
\alpha(H)=0\}}\F\alpha)$  ainsi le syst\ème de racines
associ\é \à $ ( \goth U(\F H),\goth U(\F H)\cap \goth a)$
est donn\é par la composante connexe (not\ée $\∆_H$) contenant $\beta_0$ du
syst\ème de racines   $\∆\cap (\oplus _{\{\alpha\in \Sigma\ |\  
\alpha (H)=0\}}\F\alpha)$   ce qui montre
l'irr\éductibilit\é ainsi que l'absolue irr\éductibilit\é  puisque cette
d\émonstration convient \également sur $\overline{\F}.$ \\ 
 
\noi Pour le 2\ème point, il suffit de le v\érifier
dans le cas d\éploy\é c'est \à dire lorsque $\goth
a$ est une sous-alg\èbre de Cartan et pour $u\in \goth
a.$
Notons $\Sigma_0(H)$ les racines simples de
$\∆_H$ appartenant \à $\∆_0.$ On a:
$$\goth a=\goth Z\oplus \F (2H_0-H)\oplus_{\{\alpha\in \Sigma_0(H)\}}\F h_{\alpha},$$
$\goth Z$ \étant le centre de $E_0(H),$ d'o\ù la
d\écomposition de $u:$
$$u=z+b(H_0-{H\over 2})+\sum_{\alpha\in
\Sigma_0(H)}x_{\alpha}h_{\alpha}\ ,\ z\in \goth
Z\ ,\ \ b \in \F\ , \forall \alpha\in
\Sigma_0(H):x_{\alpha} \in \F.$$
Or $ \forall \alpha\in (\∆_H)_0$ on a $$\sum_{\mu \in
{(\∆_H)}_1}n(\mu,\alpha)=\sum_{\mu \in
{(\∆_H)}_1}n(s_\alpha(\mu),\alpha)=-\sum_{\mu \in
{(\∆_H)}_1}n(\mu,\alpha) \text{  donc }\sum_{\mu \in
{(\∆_H)}_1}n(\mu,\alpha)=0$$ d'o\ù:
$$\sum_{\mu \in
{(\∆_H)}_1}\mu(u)=b.\hbox{dim}(E_0(H)\cap \goth g_1)\ \ \hbox{avec}\ \  b=2 \frac{B(u,2H_0-H)}{B(2H_0-H,2H_0-H)}\ ,$$
car $B(z,H_0-\frac{H}{2})=0$ ($2H_0- H \in [E_0(H),E_0(H)]$ par le lemme 1.1.1);
or $\chi_H(g)^{\kappa}=det (g/_{E_0(H)\cap
\goth g_1})^2$ avec $\kappa=2\frac{\hbox{dim(}E_0(H)\cap
\goth g_1)}{d_H}\in \N^* $ d'o\ù le r\ésultat. On peut noter que $d_H.b$ est entier.
\fdem
\\
 \begin{lem} On suppose que le PV
   $(\goth g_0,\goth g_1,H_0)$ est
$1$-irr\éducible (i.e. il admet un unique invariant relatif
fondamental) et admet un \él\ément  1-simple tr\ès sp\écial alors $\goth g_p=\{0\}$  pour $|p|\geq 3.$
\end{lem}

\dem

\noi Soit  $(x,H,y)$ un $sl_2$-triplet $1$-adapt\é tel que $H$  soit  1-simple tr\ès sp\écial.  \\

   1) Soient $1\leq j<i, $ rappelons que les diff\érents sous-espaces $E_{ \pm i}(H)\cap
E_{±j}(2H_0-H)$ sont en bijection, or $H$ est $1$-simple sp\écial donc $E_{-i}(H)\cap
E_j(2H_0-H)=\{0\}$   d'o\ù $E_i(H)\cap
E_j(2H_0-H)=\{0\};$  comme $2H_0-H$ est \également $1$-simple sp\écial on a    $E_i(H)\cap
E_j(2H_0-H)=\{0\}$ pour $  i\not =j\geq 1.$\\ 
Soit $p\geq 3$   alors  $E_p(H)\cap
E_p(2H_0-H)=[x,E_{p-2}(H)\cap
E_p(2H_0-H)]=\{0\}$  d'o\ù $\goth g_p=E_{2p}(H)\cap \goth g_p\oplus
E_0(H)\cap \goth g_p$ lorsque $H$ et  $2H_0-H$  sont $1$-simples sp\éciaux, donc  $E_{0}(H)\cap \goth g_p$ et  $E_{2p}(H)\cap \goth g_p$ sont des $\goth g_0$-modules. 

\noi Notons que si $ E_{2p}(H)\cap \goth g_p\not =\{0\}$ (resp.$E_0(H)\cap \goth g_p\not =\{0\}$) le polynome  obtenu en prenant le d\éterminant de l'application $(ad(\pi_H(x))^{2p}:E_{-2p}(H)\cap \goth g_{-p}\mapsto E_{2p}(H)\cap \goth g_{p}$ identifi\é avec le dual de $E_{-2p}(H)\cap \goth g_{-p}$ par $B$ (idem  avec $ad(\pi_{2H_0-H}(x))^{2p}:E_{-2p}(2H_0-H)\cap \goth g_{-p}\mapsto E_{2p}(2H_0-H)\cap \goth g_{p}$) est un  polynome d\éfini sur $\goth g_1,$ relativement
invariant par $G$  et qui  ne s'annule pas sur
$E'_2(H)\cap
\goth g_1$ (resp.$E'_0(H)\cap
\goth g_1$).   \\
2) Rappelons que l'unique invariant relatif fondamental
s'annule exactement  sur le compl\émen\-taire de $\goth g'_1$
car
$2H_0$ est $1$-simple d'o\ù $\goth g_p=\{0\}$ pour $p\geq 3$  par 1) ci-dessus.
 \fdem
 
\begin{rema} \begin{enumerate}
\item Sous les hypoth\èses du lemme  pr\éc\édent, $\goth g$ se d\écompose tr\ès simplement:
$$\goth g_0=\oplus_{-2\leq i\leq 2}E_{i}(H)\cap \goth g_0\ \ ,\ \ \goth g_1=\oplus_{0\leq i\leq 2}E_{i}(H)\cap \goth g_1\ \ ,
\ \ \ \goth g_2=\oplus_{0\leq i\leq 2}E_{2i}(H)\cap \goth
g_2.$$
Lorsque $\goth g_1$ est un $\goth g_0$-module irr\éductible, le sous-espace $E_ 1(H)\cap \goth g_0\oplus E_2(H)\cap \goth g_0$ n'est pas r\éduit \à $\{0\}.$\\

\item Lorsque $\g_i=\{0\}$ pour $|i|\≥3,$ pour tout \elt $1$-simple, $h,$ on aura $E_{-i}(h)\cap\g_1=\{0\}$ pour $i\≥2.$
\end{enumerate}
\end{rema}

\bigskip
\begin{prop}
Soient $\goth g$  une alg\èbre absolument
simple  telle que
$(\overline{\goth g}_0,\overline{\goth g}_1)$ 
est 1-irr\éductible, $\goth s$  la sous-alg\èbre engendr\ée
par un $sl_2-$triplet 1-adapt\é  dont l'\él\ément
1-simple est tr\ès sp\écial alors le pr\éhomog\ène $( \overline{\goth U(\goth s)}_0, \overline{\goth U(\goth s)}_1)$ est \également  
1-irr\éductible \à l'exception d'un unique \elt $1$-simple tr\ès sp\écial (\à l'action de $G$ pr\ès) lorsque $(\overline{\goth g}_0,\overline{\goth g}_1)$ est de type orthogonal $(B_{3n},\alpha_{2n})$ ou $(D_{3n+2},\alpha_{2n+1})$ ($n\≥1$).
\end{prop}

\dem

\noi 0) Lorsque $\g_2=\{0\},$ le r\ésultat d\écoule du lemme 1.2.1 puisque $\goth U(\goth s)=\goth U (\F H)$, $H$ d\ésignant l'\elt $1$-simple du \Sl .\\

\noi 1) Lorsque $\g_2\not=\{0\},$ n'ayant pas de d\émonstration \à
priori, on se propose de d\éterminer, \à l'action de $G$ pr\ès, tous les  \él\éments $1$-simples tr\ès sp\éciaux pour chaque syst\ème de racines irr\éductible  gradu\é apparaissant dans cette d\émonstration,  c'est \à dire d\éterminer les diagrammes \à
poids  associ\és \à un \él\ément $H\in \goth a$
v\érifiant l'ensemble des propri\ét\és suivantes:
$$\begin{array}{r c l }
&i)& \ \forall \alpha\in \∑,\ \alpha(H)\geq 0;\\
&ii)& \ \forall
\alpha\in \∆_1,\ \alpha(H)\leq 2;\\
&iii) & \ \forall \alpha\in
\∆_2,\
\alpha(H)\in\{0,2,4\};\\
&iv)&  \ \{\alpha\in \∆_1,\
\alpha(H)=2\}\not=\emptyset;\\
&v)&\ ( E_0(H)_0, E_0(H)_1)
 \ \hbox{ est un \PV tr\ès r\égulier } \end{array}$$
(c'est \à dire que $2H_0-H$  est  $1$- simple),
 et d'indiquer dans chaque cas l'\él\ément $x_0\in \goth g_1$ d'un  $sl_2-$triplet $1-$adapt\é correspondant: $(x_0,H,x_0^{-1}).$\\

\noi  Le
diagramme de Dynkin de $\goth U(\F H)$   est alors donn\é par la
r\éunion des composantes connexes de
$\{\alpha\in ∑\ |\ \alpha(H)=0\}$ contenant  au moins
une racine de $\∑_1;$  on note $(\∆_H,\∑_H)$ le  diagramme de Dynkin gradu\é du  PV $( E_0(H)_0, E_0(H)_1).$ \\
\noi Comme il suffit de faire  la d\émonstration sur une
cl\ôture alg\ébrique de $\F$ et que $\goth g_3=\{0\}$ par le lemme pr\éc\édent,   la
classification des PV
$1$-irr\éductibles \établie dans \cite{rubthese} nous indique qu'ils sont
tous tr\ès r\éguliers et qu'ils sont
tous irr\éductibles \à l'exception du cas non commutatif de syst\ème de racines gradu\é de type:
$(A_{2p+q-1},\{\alpha_p,\alpha_{p+q}\})$ avec  $q>p\geq 1;$ lorsque $∑_1$ contient une seule racine, la liste des
diagrammes gradu\és associ\és est donn\ée par la table 1 de \cite{rubthese}, celle-ci permet \également de donner les degr\és des invariants relatifs fondamentaux  (utilis\és ult\érieurement). \\
 \noi  $\goth s$ \étant l'alg\èbre engendr\ée par $\{x_0,H,x_0^{-1}\},$ on 
 v\érifie dans chaque cas que le \PV:
$( \goth U(\goth s)_0, $ $\goth U(\goth s)_1)$ est 1-irr\éductible
soit :

$\bullet$ en d\éterminant le diagramme gradu\é associ\é en calculant
directement le sous-espace $ \goth U(\goth s)_1,$ 

$\bullet$ soit en
v\érifiant que le centre de
$ (\goth U(\F H)_0)_{\theta },$ avec $ \theta=\theta_H(-1),$ est de dimension $1$ ce qui implique que le centre de 
$ \goth U(\goth s)_0$ est \également de dimension $1$ donc le \PV  $( \goth U(\goth s)_0, \goth U(\goth s)_1)$ est absolument irr\éductible r\égulier.\\

\noi En effet on consid\ère l'alg\èbre r\éductive $ \goth U( \F H)_{\theta}=\{u\in  \goth U( \F H)\ |\ \theta(u)=u\},$ en utilisant la d\écomposition de $\goth g$ en $\goth s$-modules irr\éductibles (sachant que $\g_i=\{0\}$ pour $|i|\≥3$), on v\érifie que 
 $$ \goth U( \F H)_{\theta}=
 \goth U(\goth s)\oplus \goth
I\ \hbox{avec}\ \goth I=ad(y_0)^2(E_4(H)\cap
\goth g_2)\cap \goth U(\F H)\subset \goth g_0\
,$$\noi ainsi $\goth U(\goth s)_{\pm i}  
 =(\goth U( \F H)_{\theta})_{\pm i}$ sont des $(\goth U( \F H)_{\theta})_0-$modules pour $i\geq 1$ or  ils engendrent $\goth U(\goth s) $ d'o\ù
$\goth U(\goth s)$ est un id\éal de $\goth U( \F H)_{\theta}.$ Comme la restriction de $B$ \à $\goth U( \F H)_{\theta}\times \goth U( \F H)_{\theta}$ et \à $\goth U(\goth s)\times \goth U(\goth s)$ est non d\ég\én\ér\ée et que $\goth I$ est l'orthogonal de $\goth U(\goth s)$ pour $B,$
$\goth I$ est \également un id\éal de $\goth U( \F H)_{\theta}$  d'o\ù $[\goth U(\goth s),\goth  I]=0$ donc le centre de $ \goth U(\goth s)_0$ est inclus dans le centre de $ (\goth U(\F H)_{\theta})_0=(\goth U(\F H)_0)_{\theta}.$

\noi  Lorsque la d\écomposition de $\goth U(\F H)_0$ en
id\éaux simples est de la forme:
$$\goth U(\F H)_0=Z\oplus_{i\in I}(\goth L^{(i)}\oplus
\theta (\goth L^{(i)})),\ Z \ \hbox{\étant le centre de }\
\goth U(\F H)_0,$$ 
$\goth L^{(i)},\theta (\goth L^{(i)}),i\in I,$ \étant  les id\éaux
simples non commutatifs, on v\érifie ais\ément \à l'aide de la
commutativit\é que le centre de $(\goth U(\F H)_{\theta})_0$
est  donn\é par les points fixes de $\theta$ dans $Z.$ \\

2) Le cas $A_n$ avec comme diagramme
gradu\é: $(A_{2p+q-1},\{\alpha_p,
\alpha_{p+q}\})$ avec $q>p\geq 1.$

\noi  $\∆_H$ est soit connexe, soit se compose de deux composantes connexes. \\

$\alpha)$  Lorsque $\∆_H$ est  connexe on a 
 $\∑_1\subset \∑_H$ donc par $v)$ il existe $r $ avec $1\leq r\leq p-1,$ tel que 
 $\∑_H=\{\alpha_j, p-r+1\leq j\leq p+q+r-1\}$ et
 $\alpha_{p-r}(H)=\alpha_{p+q+r}(H)=2$ en raison
de $iv)$ et de $iii), $ de plus $\alpha_j(H)=0$ pour $j\not=p-r,p+q+r$ par $ii)$  d'o\ù:
$$x_0=\sum_{i=1}^{p-r}(X_{\beta_i}+X_{\gamma_i})\ ,\
H=2\sum_{i=1}^{p-r}H_{\beta_i +\gamma_i}\ ,\ \beta_i =
\sum_{j=i}^{p+i-1}\alpha_j\ \ \hbox{et}\ \ \gamma_i =
\sum_{j=p+i}^{2p+q-i}\alpha_j.$$
Le calcul de $\goth U(\goth s)_1$ donne:
$$ E_0(\goth s)_1=\oplus_{\alpha\in (\∆')_1}\goth
g^{\alpha}\ \hbox{avec}\ \∆'=\{\epsilon_i -\epsilon_j,\
 i,j\in \{p-r+1,...,p,2p-r+1,...,p+q+r,\ i\not=j\}\ \},$$
donc un syst\ème de racines simples du pr\éhomog\ène $( \goth U(\goth
s)_0,\goth U(\goth s)_1)$ est donn\é par:   $$\∑_{\goth s}=\{\beta\ ,\ \alpha_m\  \hbox{avec}\ \ p-r+1\leq m\leq p-1 \ \hbox{ ou bien }\  2p-r+1\leq m\leq p+q+r-1\}$$ avec $\beta=\sum_{p\leq j\leq 2p-r}\alpha_j$  et $\∑_{\goth s}\cap \∆_1=\{\beta,\alpha_{p+q}\}.$\\

\noi Ainsi le \PV  $( \goth U(\goth
s)_0,\goth U(\goth s)_1)$ est de type $(A_{2r+(q-p+r)-1},\{\alpha_r,\alpha_{q-p+2r}\}).$  \\

$\beta)$ Dans le second cas, le diagramme de Dynkin de
$(\goth U(
\F H)_0,\goth U( \F H)_1)$ a deux composantes connexes.

\noi Par $v),$ on a n\écessairement $\∑_H=\{\alpha_j$ avec $1\leq j\leq 2p-1$ ou bien $q+1\leq j\leq 2p+q-1\}$  donc $2p\≤q$ et:\\

$\alpha_j(H)=\left\{ \begin{array}{ll}    0\  \textrm{pour}\  j\not =2p\ ,\ j\not=q \\  
  2\  \ \textrm{ pour}\ j=2p=q\\
  1  \  \ \textrm{ pour}\ j=2p\not=q  \ \textrm{ou bien}\ j=q\not=2p,\end{array}\right. $\\
  
$$x_0=\sum_{i=1}^p(X_{\mu_i}+X_{\nu_i}),\
H=\sum_{i=1}^p(h_{\mu_i} +h_{\nu_i}),\ \mu_i =
\sum_{j=i}^{p+q-i}\alpha_j\ \ \hbox{et}\ \ \nu_i =
\sum_{j=p+i}^{p+q+i-1}\alpha_j$$
($2p$  racines orthogonales : $\theta=\prod_{1\≤i\≤p}(\theta_{\mu_i}\theta_{\nu_i})).$\\

\noi  Soit $\goth l^{(1)}$ (resp.$\goth l^{(2)}$)
l'alg\èbre (simple) engendr\ée par 
$X_{\pm\alpha_i},i=1,...,p-1$ (resp.$p+1,...,2p-1$), la
d\écomposition de $\goth U(\F H)_0$ en id\éaux simples est
donn\ée par:
$$ \goth U(\F H)_0= \goth l^{(1)}\oplus\goth
l^{(2)}\oplus\theta\goth l^{(1)}\oplus\theta\goth
l^{(2)}  \oplus Z\ ,\ Z= \F H_1\oplus\F\theta
H_1, \
H_1=\sum_{i=1}^ph_{\epsilon_i-\epsilon_{2p-i+1}} ,
$$ \noi $Z$ \étant de dimension $2,$ $Z_{\theta}$ est de dimension $1$ donc engendr\'e par $ 2H_0-H.$\\ 
  
3)  Les cas classiques restants: $ B_n,C_n,D_n.$\\ 
 
\noi Rappelons la description du syst\ème de racines gradu\é des
cas classiques $ B_n,BC_n,C_n,D_n$ avec
 $\Sigma_1=\{\alpha_k\}$  tel que  $1\≤k\≤n-1$ pour
$C_n$  (resp. $n-2$ pour  $D_n$):  
$$\begin{array}{ll}  \∆_0^+  =  \{\epsilon_i -\epsilon_j\ ,\ 1\≤i <j\≤k\ ,\ k+1\≤i <j\≤n\ ,\ \epsilon_l\ ,\ \epsilon_l+\epsilon_m\ ,\ k+1\≤l\≤m\≤n \}\cap \∆ \\
 \∆_1  =  \{\epsilon_i \pm \epsilon_j\ ,\ 1\≤i\≤k<j\≤n\ ,\ \epsilon_i\ ,\ 1\≤i\≤k\} \cap \∆ \\ \∆_2   = \{\epsilon_i +\epsilon_j\ ,\ 1\≤i\≤j\≤k  \}\cap \∆\end{array}$$
(notations des planches II,III,IV de \cite{bourbakigal6}).\\

\noi Dans le cas d\éploy\é et r\égulier, $k$ est pair pour $C_n,$
 $3k\≤2n$ dans les cas $C_n$ et $D_n$ et $3k\≤2n+1$  
 dans le cas $B_n$. \\   
 \noi L'invariant relatif fondamental, $F,$ est de degr\é $2k$ dans le cas orthogonal (i.e. $B_n$ ou $D_n$) et de degr\é $k$ dans le cas $C_n$, de plus:
 $$N=\frac{\text{dim}(\g_1)}{\text{degr\é de }F}=\begin{cases}n-k+\frac{1}{2}\text{ dans le cas }B_n,\\
  n-k \text{ dans le cas }D_n,\\
  2 (n-k)\text{ dan le cas }C_n.\end{cases}$$

\noi  Soient $p=\left\{ \begin{array}{ll}   n & \textrm{dans le cas}\ B_n,\\   n-1 & \textrm{dans le cas}\ 
 C_n,\\
   n-2 & \textrm{dans le cas}\ D_n
   \end{array}\right.$ et $\delta =\left\{ \begin{array}{ll} \sum_{j>p}\alpha_j  \textrm{ lorsque}\  p<n,\\
   0 \textrm{ sinon.}  \end{array}\right.$\\

\noi Pour
$j=1,...,p-1,$  
$\gamma_j=\sum_{
1\≤l\≤j}\alpha_l+2\sum_{j+1\≤l\≤p}\alpha_l+\delta=\epsilon_1+\epsilon_{j+1}$ est une
racine de $\∆_1$ (resp.$\∆_2$) lorsque $j\≥k$ (resp.$j<k$)  
ainsi  les seuls diagrammes \à poids possibles sont de deux
types.\\

 $\alpha)$ Il existe une  valeur $j$ comprise
entre $1$ et $k-1$ telle que  $\alpha_j(H)\not=0.$ On peut
toujours supposer que $\alpha_l(H)=0$ pour $1\≤l\≤j-1,$ en
prenant une valeur minimale de $j.$

\noi  En
consid\érant les valeurs $\gamma_l(H),$ on v\érifie
que  $\alpha_l(H)= 0$ pour $l\≥j+1$  
et que  
$\alpha_j(H)=2.$\\ 
\noi Dans ce cas le PV irr\éductible et r\égulier : $(\goth U( \F H)_0,\goth U( \F H)_1)$ a un diagramme gradu\é $(R_{n-j},\alpha_{k-j})$ analogue \à celui
d'origine (i.e. $R=B$ ou $D$ dans le cas orthogonal et $R=C$ dans le cas $C_n$) d'o\ù $j$ est pair dans le cas
$C_n$ ( \cite{rubthese}). \\

\noi Soit $j_0=[\frac{j}{2}],$   on a pour $j\≥2:$
$$ H=2(\sum_{1\≤i\≤j_0}h_{\beta_i+\delta_i})+H_1\ ,\
x_0=\sum_{1\≤i\≤j_0}(X_{\beta_i}+X_{\delta_i})+X_1$$
 avec   $\beta_i=\epsilon_{2i-1}+\epsilon_{n-i+1}$ et $\delta_i=\epsilon_{2 i}-\epsilon_{n-i+1}$ pour $1\≤i\≤j_0.$\\
 
\noi  $H_1=X_1=0 $ si $j$ est pair et sinon  
 $$H_1= \left\{ \begin{array}{ll}   h_{\epsilon_j}\  \textrm{dans le cas}\ B_n,\\  
 h_{\epsilon_j-\epsilon_{n-j_0}}+
h_{\epsilon_j+\epsilon_{n-j_0}} \ \textrm{dans le cas}\ 
 D_n, 
   \end{array}\right. \quad \textrm{et}\quad  X_1=\left\{ \begin{array}{ll}   X_{\epsilon_j}\  \textrm{dans le cas}\ B_n,\\  
 X_{\epsilon_j-\epsilon_{n-j_0}}+
X_{\epsilon_j+\epsilon_{n-j_0}} \ \textrm{dans le cas}\ 
 D_n, 
   \end{array}\right.$$

 \noi  Le cas $j=1$  se r\éduit \à : $x_0=X_1$ et $H=H_1.$\\
 
\noi Le calcul de $\goth U(\goth s)_1$ donne:
$$\goth U(\goth s)_1=\oplus_{\alpha\in \∆'_1}\goth
g^{\alpha}\ \oplus\ \goth B$$
avec
$$  \∆'_1= \left\{ \begin{array}{ll}   \{\pm\epsilon_i\pm \epsilon_q,\
j+1\≤i\not=q\≤n-j_0\ ,\ \pm\epsilon_i\ ,\
\pm 2\epsilon_i,\ j+1\≤i\≤n-j_0\}\cap \∆_1\  \textrm{lorsque}\ j\  \textrm{est pair},\\  
  \{\pm\epsilon_i\pm
\epsilon_q,\ j+1\≤i\not=q\≤n-j_0\}\cap \∆_1\ \textrm{lorsque}\ j\  \textrm{est impair}\ \textrm{et}\
 \∆=B_n, \\
  \{\pm\epsilon_i\pm
\epsilon_q,\ j+1\≤i\not=q\≤n-j_0-1\}\cap \∆_1\ \textrm{ lorsque}\ j\  \textrm{est impair}\ \textrm{et}\
 \∆=D_n
   \end{array}\right.$$
   et 
   $$\goth B = \left\{ \begin{array}{ll}    \{0\}\  \textrm{lorsque}\  \∆=B_n\  \textrm{ou}\ C_n\  \textrm{et lorsque }\ \∆=D_n\  \textrm{avec}\ j=2j_0 ,\\  
\oplus_{j+1\≤l\≤k} \F X_l\ \    \textrm{avec}\ \ X_l=[X_{\epsilon_j+\epsilon_{n-j_0}}-
X_{\epsilon_j-\epsilon_{n-j_0}},X_{\epsilon_l-\epsilon_j}]\  \textrm{lorsque}\  \∆=D_n \  \textrm{et}\
j=2j_0+1 . 
   \end{array}\right.$$

\noi Les r\ésultats sont analogues pour $\goth U(\goth s)_{-1}$ 
d'o\ù le syst\ème de racines du pr\éhomog\ène
$(\goth U(\goth s)_0,\goth U(\goth s)_1)$ est du  type
suivant :

$$\left\{ \begin{array}{ll}     (D_{n-j_0-j},\alpha_{k-j})\  \textrm{lorsque}\  \ \∆=B_n\  \textrm{et}\ \ j=2j_0+1\ \textrm{ou bien } \ \∆=D_n\   \textrm{et}\ j=2j_0\ ,\\  
  (B_{n-j_0-j},\alpha_{k-j})\ \     \  \textrm{lorsque}\  \∆=B_n\   \textrm{et}\ j=2j_0\,\\  
  (B_{n-j_0-j-1},\alpha_{k-j})  \  \textrm{lorsque}\  \∆=D_n\   \textrm{et}\ j=2j_0+1\,\\  
  (C_{n-j_0-j},\alpha_{k-j}) \  \textrm{lorsque}\  \∆=C_n,
 \end{array}\right.$$
A l'exception des $2$ cas:\\

$\bullet$ $(D_{3n+2},\alpha_{2n+1})$ et $j=2n$ donc $(E_0(H)_0,E_0(H)_1)$ est de type $(D_{n+2},\alpha_1),$\\

$\bullet$ $(B_{3n},\alpha_{2n})$ et $j=2n-1$ donc $(E_0(H)_0,E_0(H)_1)$ est de type $(B_{n+1},\alpha_1),$\\

\noi pour lesquels on a $A_1\oplus A_1.$\\

\noi Dans tous les cas, on a $N=N(H)=\ds\frac{\text{dim}(E_0(H)_1)}{\text{degr\é de l'\irf}}$ .

\bigskip  

$\beta)$ Pour $l=1,...,k-1$ on a $\alpha_l(H)=0,$
   on note
$j+1$ le premier indice pour lequel
$\alpha_{j+1}(H)\not=0 ,$  dans
ce cas $\∑_H=\{\alpha_m,1\≤m\≤j\}$ d'o\ù $(\∆_H,\∑_H)$  
 est de type 
$(A_j,\alpha_k),$ et en raison de la r\égularit\é (v)) on a $j=2k-1$ d'o\ù $2H_0-H=\sum_{1\≤i\≤k}h_{\mu_i}$ avec $\mu_i=\epsilon_i-\epsilon_{k+i}$ or $2H_0=\sum_{1\≤i\≤k}(h_{\mu_i}+h_{\nu_i})$ avec 
$\nu_i=\epsilon_i+\epsilon_{k+i}$ donc
 $$H=\sum_{1\≤i\≤k}h_{\nu_i}\ ,\  x_0=\sum_{1\≤i\≤k}X_{\nu_i}\quad (k\ \textrm{tds de type } \ A_1\ :\ \theta=\prod_{1\≤i\≤k}\theta_{\nu_i}).$$
On proc\ède exactement comme dans 2) $\beta).$
Soit $\goth l$ l'alg\èbre (simple) engendr\ée par 
$X_{\pm\alpha_i},i=1,...,k-1;$  la d\écomposition de $\goth U(\F H)_0$ en id\éaux simples est
donn\ée par:
$$ \goth U(\F H)_0= \goth l \oplus \theta\goth l\oplus 
  Z\ ,\ Z= \F (2H_0-H) 
$$ \noi d'o\ù le r\ésultat.\\ 

\noi On peut noter que le degr\é de l'\irf du \PV  commutatif: $(E_0(H)_0,E_0(H)_1)$ est de degr\é $k,$  il en est de m\ême pour le \PV $(\goth U(\goth s)_0,(\goth U(\goth s)_1)$ dans le cas orthogonal par orthogonalit\é des racines $\mu_i,i=1,...,k.$ Dans le cas $C_n,$ le calcul de $\goth U(\goth s)_1$ permet de v\érifier que le \PV  $(\goth U(\goth s)_0,(\goth U(\goth s)_1)$  est de type $(D_k,\alpha_k)$ donc l'\irf est de \dg  $\displaystyle\frac{k}{2}.$ Notons que $H$ et $2H_0-H$ sont dans la m\ême orbite de $G.$\\

c) Les cas exceptionnels: cf.la d\émonstration du lemme  suivant.\fdem\\

\subsection{\bf    Existence}

\bigskip

  \begin{lem} Un \PV  $(\goth g_0,\goth g_1)$ absolument
irr\éductible et r\égulier de type exceptionnel (i.e. $\∆$ est
de type exceptionnel) tel que $\goth g_3=\{0\}$ admet toujours un  \él\ément 1-simple
 tr\ès sp\écial sauf dans les cas suivants:
 \begin{enumerate}
 \item    $\∆=G_2$ 
 \item $(\∆,\∑_1)=(F_4,\{\alpha_4\})$ et $\goth g$ est
d\éploy\é.
 \end{enumerate}
  \end{lem}
 
 \dem

On la fait cas par cas et lorsque 
 $\∆_2$ est non vide donc le coefficient de l'unique racine de $\∑_1$ dans la
d\écomposition de $\tilde \alpha$ suivant $\∑$ prend la valeur
 $2$, ce qui donne la liste suivante des \PVs  concern\és par le lemme :\\

\noi $(G_2,\alpha_2);$ $(F_4,\alpha_i),$ $i=1$ ou $4; $
$(E_6,\alpha_i),$ $i\not=1,4,6;$ $ (E_7,\alpha_i),$ $i\in
\{1,2,6\};$ $(E_8,\alpha_i),$ $i=1$ ou $8.$

\noi (notations: planches de \cite{bourbakigal6}).\\
 
 \noi  Lorsque $\∆$ est de type $E_i,i\in \{6,7,8\},$ l'alg\èbre
$\goth g$ est d\éploy\ée (cf.tables de \cite{veisfeiler} et \cite{warner}) donc  pour 
$\∆=E_6$ on a $\∑_1=\{\alpha_2\}$ par r\égularit\é (cf.table 1 de \cite{rubthese}).\\

\noi  Les \él\éments tr\ès sp\éciaux sont  les  \él\éments  $H\not= 2H_0,$ $1-$simples, dont le spectre de   $ad(H) /\goth g_2$ ne comprend que des valeurs paires  c'est \à dire tels que $\alpha (H)\in 2\N$ pour tout  $\alpha \in \∆_2.$\\

1) Le cas $(E_7,\alpha_6).$\\

\noi Posons: $\beta_1=\alpha_6,$
$\beta_2=\alpha_2+\alpha_3+2(\alpha_4+\alpha_5)+\alpha_6,$
$\beta_3=\alpha_5+\alpha_6+\alpha_7.$

\noi Les repr\ésentants des orbites $1$-simples sont donn\és par:
$$\hbox{Pour}\ i=1,2,3:\ \sum_{j=1}^ih_{\beta_j},\ 2h_{\tilde
\alpha},\ 
 2h_{\tilde \alpha}+h_{\beta_1},\
h_{\beta_1}+h_{\beta_3},\  2h_{\tilde
\alpha}+h_{\beta_1}+h_{\beta_3},\  2H_0,$$
(par exemple prop.2.6 de \cite{mullerJA1})\\

\noi mais seul $H=2h_{\tilde \alpha}$ convient d'o\ù le diagramme
de Dynkin gradu\é du PV: $(\goth U(\F H)_0,\goth U(\F H)_1)$ est de type $(D_6,\alpha_2)$  et on v\érifie facilement que $2H_0-H$ et $H$ sont dans la m\ême orbite de $G.$

\noi On peut noter
que l'invariant relatif fondamental de ce \PV est de degr\é
$4$ comme celui de d\épart.\\

\noi On termine la d\émonstration de la proposition pr\éc\édente.\\
\noi Un $sl_2-$triplet $1$-adapt\é, $(x_0,H,y_0),$ est alors
donn\é par :
$$x_0=X_{\gamma_1}+X_{\gamma_2}\ \hbox{avec}\
\gamma_1=\sum_{1\leq i\leq 7}\alpha_i+\alpha_4+\alpha_5,\ 
\gamma_2=\sum_{1\leq i\leq 6}\alpha_i+\alpha_3+\alpha_4.$$
Le calcul donne:
$$\goth U(\goth s)_1=\oplus_{\alpha\in (\∆')_1}\goth
g^{\alpha}\ \hbox{avec}\ \∆'=(\oplus_{1\leq i\leq3}\Z\delta_i\oplus \Z\alpha_2\oplus \Z\alpha_4)\cap
\∆,\hbox{ et}$$
$$\delta_1=\alpha_5+\alpha_6,\ \delta_2=\alpha_6+\alpha_7,\
\delta_3=\alpha_3+\alpha_4+\alpha_5.$$
On v\érifie ais\ément que le diagramme de Dynkin gradu\é
associ\é au \PV: $(\goth U(\goth s)_0,\goth U(\goth
s)_1),$ est  donn\é par $(\∆',\∑_1=\{ \delta_1,\delta_2\})$ et que $\∑=\{\delta_1,\alpha_4,\alpha_2,\delta_3, \delta_2\}$ est un syst\ème de racines simples de $\∆'.$ Le \PV 
$(\goth U(\goth s)_0,\goth U(\goth
s)_1)$ est de type $(A_5,\{\alpha_1,\alpha_5\}),$ il n'est pas irr\éductible mais est $1$-irr\éductible et  son \irf est de degr\é $2.$\\

2) Dans tous les autres cas, il existe $p$ racines
orthogonales de $\∆_1$: $\lambda_1,...,\lambda_p,$ telles que
$\sum_{1\leq i\leq p}h_{\lambda_i}=2H_0,$ et pour tout \él\ément
$1$-simple, $h,$ il existe un sous-ensemble $I$ de
$\{1,...,p\}$ tel que $h$ et $\sum_{i\in I}h_{\lambda_i}$
soient dans la m\ême orbite de $G$ (\cite{mullerJA1},proposition 6.6 et
corollaire 4.6). Les racines $\lambda_1,...,\lambda_p,$ sont toujours fortement orthogonales et l'\irf est de degr\é $p$ \à l'exception de l'unique cas o\ù $\∆=F_4$ et $\∑_1=\{\alpha_4\}.$\\

\noi Les valeurs possibles de $p$ sont les suivantes:\\

 \noi a) $p=2$ alors $\∆=G_2$ et on a $3$ orbites $1$-simples de
repr\ésentants:
$2H_0,h_{\lambda_1},h_{\lambda_2}$ avec $\lambda_1=\alpha_1$ et $\lambda_2=2\alpha_1+\alpha_1$ mais $\tilde
\alpha(h_{\lambda_1})=1$ donc il n'existe aucun \él\ément
$1$-simple  tr\ès sp\écial.\\

\noi b) $p=4$ mais alors $\omega=\frac{1}{2}(\sum_{1\≤i\≤4}\lambda_i)$ est une racine de  $ \∆_2.$

\noi Les repr\ésentants des  orbites $1$-simples sont \à prendre
parmi: $H_i=\sum_{1\≤j\≤i}h_{\lambda_j},i=1,...,4;$ comme
$\omega(H_i)=i,$   $H_2$  est le seul   \él\ément
$1$-simple  tr\ès sp\écial (\à l'action de $G$ pr\ès) lorsque $(\∆,\∑_1)\not=(F_4,\alpha_4).$\\

\noi Dans le cas $(F_4,\alpha_4),$ il convient de distinguer le
cas o\ù $\goth g$ est d\éploy\ée (seuls $H_1$ et $H_4$ sont
$1$-simples) du cas o\ù $\goth g$ ne l'est pas.

\noi Dans ce
dernier cas la consultation des $\F$-formes des alg\èbres
absolument simples (\cite{veisfeiler},\cite{warner}) donne comme diagramme de Dynkin gradu\é
possibles pour $(\overline {\goth g}_0,\overline {\goth
g}_1)$ la liste suivante : \\

$\bullet$ $(E_8,\alpha_1)$ ($\F=\R$) donc  $ ({\goth g}_0,{\goth
g}_1)$ est un PV presque commutatif au sens de \cite{mullerJA1}  d'o\ù $H_2$ est $1$-simple
(derni\ère remarque de \cite{mullerJA1}) et c'est l'unique  \él\ément
$1$-simple  tr\ès sp\écial (\à l'action de $G$ pr\ès).\\

$\bullet$ $(E_7,\alpha_6)$ avec la $\F$-forme EVI (notation des tables de \cite{warner}).\\

\noi Il est facile de v\érifier que la plus grande racine $\tilde \omega$ de $\overline\∆$ est une  $\F$-racine d'o\ù $ 2h_{\tilde \omega}\in \goth g$ donc est $1$-simple
 tr\ès sp\écial pour le \PV: $(\goth
g_0,\goth g_1)$  et on a $H_2=2h_{\tilde  \omega}.$\\

\noi Dor\énavant $\{\lambda_1,...,\lambda_p\}$  est la chaine
canonique de racines (\cite{mullerJA1}, la liste explicite de 
$\{\lambda_1,...,\lambda_p\}$ et une description
compl\ète de $\∆_2$ \à l'aide des $(\lambda_i)_{1\≤i\≤p}$ sont donn\ées dans les tables de \cite{mullerJA2}).
 
\noi En regardant les diverses valeurs $\beta(H),\beta\in \∆_2$
et $H$ d\écrivant les repr\ésentants des orbites $1$-simples,
on obtient:\\

c)
$p=7$ alors $(\∆,\∑_1)=(E_7,\alpha_2),$ seuls $2$
repr\ésentants conviennent :
$H=h_{\lambda_1}+h_{\lambda_2}+h_{\lambda_3}$  et
$2H_0-H=h_{\lambda_4}+h_{\lambda_5}+h_{\lambda_6}+h_{\lambda_7}=2h_{\tilde \alpha}
.$\\

d) $p=8$  alors $(\∆,\∑_1)=(E_8,\alpha_1),$ seul
$H=\sum_{1\≤j\≤4}h_{\lambda_j}=2h_{\omega}$  est $1$-simple
 tr\ès sp\écial.\\
 
 \noi On termine la d\émonstration de la proposition pr\éc\édente.\\
\noi  Dans tous ces cas, $H$ est de la forme : $H=\sum_{i\in
I}h_{\lambda_i},$ un $sl_2$-triplet $1$-adapt\é, not\é
$(x_0,H,y_0),$ est alors donn\é par:
$$x_0=\sum_{i\in I}X_{\lambda_i},\ y_0=\sum_{i\in
I}X_{-\lambda_i}.$$
$(\goth U (\F H)_0,\goth U( \F H)_1)$ (resp.$(\goth U(\goth
s)_0,\goth U(\goth s)_1)$) est un PV quasi commutatif qui
admet $\{\lambda_i,i\notin I\}$ comme syst\ème orthogonal
maximal (d\émonstration du lemme 3.2 de \cite{mullerJA2}).

\noi Comme  $(\goth U (\F H)_0,\goth U( \F H)_1)$ est
irr\éductible (lemme 1.2.1) on a en appliquant la d\émonstration de la partie
2) de la prop.5.1.1 de \cite{mullerJA1} (v\érification imm\édiate): 
$$\forall i\not=j\notin I\ ,\ \exists k\not=l\ 
\hbox{tels que}\ :\ ( k\ \hbox{et}\ l\in I) \ \hbox{ou}\
(k\ \hbox{et}\ l\notin I)\
 \hbox{ et  }\  \omega={1\over
2}(\lambda_i+\lambda_j+\lambda_k+\lambda_l)\in \∆_2,$$
  donc
 le \PV 
$(\goth U(\goth s)_0,\goth U(\goth s)_1)$  est \également
irr\éductible (2) de la prop.5.1.1 de \cite{mullerJA1}:
  lorsque $k,l\in I$ on notera
que
$u_{\omega}=[X_{-\lambda_k},X_{\omega}]-
[X_{-\lambda_l},X_{\omega}]\in E_1(h_i)\cap E_1(h_j)\cap
\goth U(\goth s)_1$). \\

\noi Notons que l'  \irf du \PV  $(\goth U (\F)_0,\goth U( \F)_1)$ est de \dg $p$-cardinal$(I)$ comme celui du \PV  $(\goth U(\goth s)_0,\goth U(\goth s)_1).$\fdem

\bigskip

\noi Le lemme pr\éc\édent se compl\ète par:\\

\begin{prop}  Soit $\goth g$ soit une alg\èbre absolument simple
engendr\ée par $\goth g_0$ et $\goth g_{±1},$ qui n'est pas de
rang 1, ni de type $G_2$ et telle que

\begin{enumerate}

 \item $\goth g_p=\{0\}$ pour $|p|\≥3$

 \item $(\overline{\goth g}_0,\overline{\goth g}_1)$ a un unique
invariant relatif fondamental de degr\é plus grand que deux
\end{enumerate}
alors il existe des \elts 1-simples tr\ès sp\éciaux appartenant \à $\goth a$ de somme $2H_0$ et $2H_0$ est 1-simple.
\end{prop}

\dem

1) $2H_0$ est $1$-simple par consultation des tables de \cite{rubthese}.\\

\noi 2) Il suffit d'exhiber un \elt  $1$-simple tr\ès sp\écial lorsque   $\Delta_2$ est non vide.\\

\noi a)  Lorsque $\∆$ est de type exceptionnel, cela r\ésulte du
lemme pr\éc\édent  puisque l'invariant relatif fondamental est une
forme quadratique lorsque $\goth g$ est une alg\èbre d\éploy\ée
de type $F_4$ munie de la graduation induite par $\alpha_4,$ ce qui est exclu dans l'\énonc\é.\\
 
 \noi b) Dans les diff\érents cas classiques $B_n,BC_n,C_n$ (resp.
$D_n$ avec $n\≥4$) et $\Sigma_1=\{\alpha_k\}$ (resp.
$2\≤k\≤n-2$),   le r\ésultat est \évident pour les cas d\éploy\és par la \demo de la proposition 1.2.4 et sinon on applique le lemme 2.4.2.
  \\

\noi c) $\∆$ est de type $A_n$ donc $∑_1=\{\alpha_p,\alpha_q\}.$

\noi Comme  
$\Delta$   est simplement lac\é, $\tilde
\alpha$  est la somme de deux racines, $\alpha$ et $\beta$ de
$\Delta_1$ dont la diff\érence n'est pas une racine
  donc $2h_{\tilde
\alpha}$ est $1$-simple ce qui termine la d\émonstration
lorsque $H_0$ et $h_{\tilde
\alpha}$  ne sont pas proportionnels. Sinon
$\Sigma_1$ est compos\é des racines simples
reli\ées \à $\tilde
\alpha $ dans le  graphe de Dynkin  compl\ét\é donc $p=1$ et
$q=n.$

\noi Par consultation des tables de
\cite{veisfeiler}, \cite{warner}, on v\érifie que $\overline \∆$ est \également de
type $A_m$ puis par la classification de Rubenthaler
(prop.3.3.7 de \cite{rubthese}) on obtient que $n\≥3;$ l'\él\ément 
$H=h_{\alpha_1}+h_{\alpha_n}$ est $1$-simple puisque les
racines  $\alpha_1$ et $\alpha_n$  de $\Delta_1$ sont
fortement orthogonales et  $\forall \alpha\in \∆_1$   on a  $\alpha (H)\in  \{0,1,2\}$ et $H\not=2H_0.$ \fdem

\bigskip

\subsection{\bf Une classe de \PVs faiblement sph\ériques}

\bigskip

\noi Les r\ésultats des paragraphes pr\éc\édents nous am\ènent \à consid\érer les \PVs de type paraboliques ayant les propri\ét\és suivantes:

\begin{enumerate}
\item $\goth g$ est une alg\èbre absolument simple engendr\ée par $\goth g_{\pm 1}$,

\item $ ({\goth g}_0, {\goth g}_1,H_0)$ est un PV  absolument irr\éductible et r\égulier,  

\item il existe des \elts $1$-simples tr\ès sp\éciaux appartenant \à $\goth a$ de somme $2H_0.$

\end{enumerate}

\bigskip

\noi Ce sont  les \PV  $(\goth g_0,\goth g_1)$  absolument irr\éductibles, r\éguliers qui ne sont pas de rang $1,$ ni de type $G_2$ ou $(F_4,\alpha_4)$ d\éploy\é et pour lesquels $\goth g_p=\{0\}$ pour $|p|\≥3.$\\

\noi Notons  $H_1,...,H_p$ les \elts $1$-simples tr\ès sp\éciaux appartenant \à $\goth a$ de somme $2H_0,$ $\goth t$  la sous-alg\èbre qu'is engendrent, $\∆_R(\goth t)$ les restrictions, qui sont non nulles, des racines de $\∆$ \à $\goth t$ que l'on munit de  l'ordre
suivant:
$$ \lambda\succ 0 \ \Leftrightarrow\
\lambda(H_k)>0\quad \hbox{avec}\quad k=sup\{j\ ,\
\lambda(H_j)\not=0\} \ ,$$
$P_{\goth t}$ le  sous-groupe parabolique de $G$ associ\é \à $\goth
t$ avec cet ordre : $P_{\goth t}=G_{\goth t}.N _{\goth
t}$ avec $N_{\goth t}=exp(ad(n _{\goth t})), $  
 $n_{\goth t}=\oplus_{\lambda\succ 0,\lambda\in
\Delta_0}\goth g^{\lambda} $ et $G_{\goth t}$ est le centralisateur de $\goth t$ dans $G.$\\

\noi Soit $\goth p_{\goth t}=E_0(\goth t)\oplus n _{\goth t},$ on rappelle que $E_0(\goth t)=\goth t^0\oplus \goth m,$ $\goth t^0=\cap_{\lambda\in \∆_0|\lambda/\goth t=0}\text{Ker}\lambda$ est le centre de $E_0(\goth t)$ et $\goth m(=E'_0(\goth t))$ est l'orthogonal de $\goth t^0$ dans $E_0(\goth t)$ pour la forme de Killing, et que $\goth p_{\goth t}=\goth t^0\oplus \goth m \oplus n_{\goth t}$ est la d\écomposition de Langlands de la sous-alg\èbre parabolique  $\goth p_{\goth t}.$\\

\noi Dans cette situation, on a une d\écomposition particuli\èrement simple: $$\goth p_{\goth t}=E_0(\goth t)\oplus_{1\≤i<j\≤p}(E^{i,j}_{-1,1}\oplus E^{i,j}_{-2,2})$$avec
$E^{i,j}_{a,b}=\{x\in \goth g\ |\ [H_q,x]=0$ pour $q\not =i,j$ et $[H_i,x]=ax,[H_j,x]=bx\}$ et on notera $P_{\goth t}=P(H_1,...,H_p)$ lorsqu'on d\ésire mettre en \évidence l'ordre associ\é.

\begin{rema} 

\begin{enumerate}
\item $\goth t$ est unique (\à l'action de $G$
pr\ès) et est de dimension $2$ dans les cas exceptionnels non
commutatifs (d\émonstration du lemme  1.3.1).

\item  Lorsque $\goth t$ est inclus dans
l'alg\èbre de Lie engendr\ée par $\goth g_{\pm 2},$ 
   $G_{\goth t}$
contient le sous-groupe distingu\é $H,$ 
 centralisateur de   $\goth g_2$ dans $G$  ainsi 
$P_{\goth t}$ est un sous-groupe parabolique provenant
du groupe quotient $\Gamma = G / H.$

\item Puisque tout sous-groupe parabolique est conjugu\é \à un sous groupe parabolique standard,    \étant donn\é  l'ordre fix\é \à priori dans $\∆$ tel que $\∆_1\subset \∆^+,$ il existe des \él\éments $1$-simples tr\ès sp\éciaux, $H'_1,...,H'_p,$  tels que $ \lambda\succ 0 \ \Leftrightarrow\ \lambda >0$ et $\lambda/\goth t'\not=0,$  $\goth t'$ \étant la sous-alg\èbre engendr\ée par  $H'_1,...,H'_p.$

\noi Soit $\∑_{\goth t'}=\{\alpha\in \∑\ |\ \alpha/\goth t'=0\}$ et $<\∑_{\goth t'}>$ le syst\ème de racines engendr\é par $\∑_{\goth t'},$
alors $\∑_{\goth t'}\subset \∑_0$ et 
$$\{\lambda\in \∆_0\ |\ \lambda/\goth t'=0\}=<\∑_{\goth t'}> \quad ,\quad
\{\lambda\in \∆_0\ |\ \lambda \succ 0 \}=\∆_0^+-<\∑_{\goth t'}>,$$
$$\goth t'^0=\cap_{\lambda\in  \∑_{\goth t'}} \text{Ker}\lambda$$
\noi  et a pour dimension le nombre d'\elts de $\∑_1- \∑_{\goth t'}.$
 \end{enumerate}
\end{rema}
\bigskip

\begin{defi} Les sous-groupes paraboliques construits \à partir d'\elts 1-simples tr\ès sp\éciaux qui commutent et de somme $2H_0$ sont appel\és paraboliques tr\ès sp\éciaux.
\end{defi}

\bigskip

\noi Pour $k=1,...,p-1,$ soient $h_k=H_1+...+H_k,$ $\goth U^+(k)=\goth U(\F (2H_0-h_k))$ et $F_k$  
un invariant relatif fondamental du pr\éhomog\ène
$$(  \goth U^+(k)_0,\goth U^+( k)_1)=(E_0(h_k)\cap \goth g_0,E_2(h_k)\cap \goth g_1 )   $$
 gradu\é par $ad(\frac{h_k}{2})$,  soient  $\goth U^-(p-k)=\goth U(\F ( h_k))$ et  $F^*_{p-k}$ un invariant relatif fondamental du pr\éhomog\ène
$$ 
 (\goth U^-( p-k))_0,\goth U^-(p-k))_{-1})= (E_0(h_k)\cap \goth g_0,
E_0(h_k)\cap
\goth g_{-1} ) $$
\noi  gradu\é par $ ad(H_0-\frac{h_k}{2}), $  ces $2$ \PVs sont   absolument irr\éductibles par le
 lemme 1.2.1.\\
 
\noi On prolonge  
$F_k$  
  sur $\goth g_1$  et  $F^*_{p-k}$ sur $\goth
g_{-1}$ en conservant la m\ême notation, gr\âce au choix
naturel des suppl\émentaires : $\oplus_{i\not= 2}E_i(h_k)\cap
\goth g_1$ dans $\goth g_1$ et $\oplus_{i\not= 0}E_i(h_k)\cap
\goth g_{-1}$ dans  $\goth g_{-1}$. \\
\noi On note
$\chi_k$ (resp.$\chi^*_k$) le caract\ère associ\é \à $F_k$
(resp.$F^*_k$) que l'on \étend sur $  P_{\goth t}$ par
$\chi_k(gn)=\chi_k({\tilde g}),$
${\tilde g}$ d\ésignant la restriction de $g$ \à
$\goth U(\F (2H_0-h_k)).$

\noi On notera \également $F=F_p$  et $\chi=\chi_p$
(resp. $F^*=F^*_p $ et $\chi^*=\chi^*_p$), on pose :
$\chi_0=\chi^*_0=1,$ $h_0=0$ et $h_p=2H_0.$\\

\noi On a alors les  inclusions suivantes:
$$\goth U^+(1)_1\subset \goth U^+(2)_1 \subset ...\subset \goth U^+(p)_1=\goth g_1\quad  \hbox{et}\quad 
 \goth U^-(1)_{-1}\subset  \goth U^-(2)_{-1} \subset ...\subset \goth U^-(p)_{-1}=\goth g_{-1}.$$
 
 \noi Ces polynomes g\én\éralisent la notion de mineurs principaux
pour une matrice carr\ée et ont \ét\é largement \étudi\és notamment lorsque
le \PV est commutatif (cf. l'introduction ainsi que les travaux de  H.Rubenthaler et G.Schiffmann,
 et Y.Angely etc).\\

\begin{lem} Les polynomes $F_k,k=1,...,p$
(resp.$F^*_k$) sont relativements invariants par
$  P_{\goth t}$ de caract\ère
$\chi_k$ (resp.$\chi^*_k$).
\end{lem}

\dem Pour
$k=1,...,p,$ $F_k$  (resp. $F^*_k$)  est
relativement invariant par $G_{\goth t}$  de caract\ère
 $\chi_k$ (resp.$\chi^*_k$) car $G_{\goth t}$ centralise
$h_k.$  

\noi Il suffit de montrer l'invariance par $N_{\goth t}$ pour
$k=1,...,p-1,$ ce qui se fait simplement par r\écurrence sur
$p.$ V\érifions-le sur  $\goth g_1.$

\noi Soit $\goth 
B=\goth U(\F H_p);$ tout \él\ément
$x$ de
$\goth g_1$ se d\écompose suivant $ad(H_p):$
$$x=x_2+x_1+x_0,\ x_i\in E_i(H_p)\cap \goth g_1,\ \ 
\hbox{donc}\ \ F_k(x)=F_k^{\goth B}(x_0),\ k=1,...,p-1,$$
$F_k^{\goth B}$ d\ésignant l'extension de l'invariant
fondamental du pr\éhomog\ène $(\goth U(\F
(2H_0-H_k))_0,$ $\goth U(\F (2H_0-H_k))_1)$ \à
$\goth B_1.$

\noi Or $N_{\goth t}=N_1.N_2$ avec $N_i=exp(ad(\goth n_i)),i=1,2,$
$\goth n_1=\oplus_{\{\lambda\in \∆_0\ |\lambda\succ
0,\lambda(H_p)=0\}} \goth g^{\lambda}\subset \goth B_0$ et
$\goth
n_2=\oplus_{\{\lambda\in \∆_0\ | 
 \lambda(H_p)>0\}} \goth g^{\lambda}$ (\cite{seligman}) d'o\ù:
$$F_k(nx)=F_k^{\goth B}(n_1x_0)=F_k^{\goth
B}(x_0)=F_k(x),\ n=n_1n_2, $$ 
en utilisant l'hypoth\èse de r\écurrence.
\fdem

\bigskip

\noi Les \PV duaux,  obtenus
  \à l'aide de la forme de Killing,
correspondent aux sous-espaces vectoriels :
$${\goth U^+(k)_1}^*=E_{-2}(h_k)\cap \goth g_{-1}=E_0(H_p+...+H_{k+1})\cap \goth g_{-1}\ ,\ $$
$${\goth U^-(p-k)_{-1}}^*=E_0(h_k)\cap \goth g_1=E_2(H_p+...+H_{k+1})\cap \goth g_1\ ,$$
 ils sont donc associ\és \à l'indexation $H_p,...,H_1$ d'o\ù les invariants relatifs correspondants prolong\és sur $\goth g_{-1}$ et $\goth g_{1},$ que l'on note $P^*_k$ pour $\goth U^+(k)_1^*$ et $P_{p-k}$ pour $\goth U^-(k)_{-1}^*,$ sont relativement invariants sous l'action du parabolique oppos\é \à $P_{\goth t}:$ 
$ P^-_{\goth t}=G_{\goth t} N^-_{\goth
t} $ avec $N^-_{\goth
t}=exp(ad(\oplus _{-\lambda \succ 0,\lambda\in
\Delta_0}\goth g^{\lambda})).$   
  \\

\noi Le lemme  suivant est une  g\én\éralisation du lemme de
Gauss qui \énonce qu'une matrice sym\étrique dont tous les
mineurs principaux sont non nuls peut \être diagonalis\ée \à
l'aide d'une matrice triangulaire, ces derni\ères 
 \étant remplac\ées par $P_{\goth t}.$\\
 
\begin{lem}
Soit $W_{\goth t} 
=\{x_1+...+x_p\ ,\ x_i\in E'_2(H_i)\cap \goth g_1\ ,\
  [x_i,x_j]=0 \ ,i,j=1,...,p\}$ et  $\goth
g"_1=\{
x \in   
\goth g_1\ |$ tel que
$\prod_{1\≤k\≤p}F_k(x)\not=0\}$  alors $G_{\goth t}$ op\ère
dans $W_{\goth t}$ avec un nombre fini d'orbites et les
repr\ésentants des orbites de $G_{\goth t}$ dans $W_{\goth
t}$ sont \également des  repr\ésentants des orbites de
$P_{\goth t}$ dans $\goth
g"_1.$ 

\noi On a le m\ême   r\ésultat  avec $\goth
g"_{-1}$ et $W^*_{\goth t} 
=\{y_1+...+y_p\ ,\ y_i\in E'_{-2}(H_i)\cap \goth g_{-1}\ ,\
  [y_i,y_j]=0 \ ,i,j=1,...,p\}$
\end{lem} 
 
 \dem La premi\ère assertion est \évidente puisque $G_{\goth t}$  a
un nombre fini d'orbites dans $\oplus_{1\≤i\≤p}E_2(h_i)\cap\goth g_1$ (\cite{vinberg}) et op\ère sur  $W_{\goth t} .$

\noi Pour la seconde assertion, on v\érifie que $N_{\goth t}. 
  W_{\goth t}=\goth
g"_1$ par r\écurrence sur $p,$ le cas $p=2$ r\ésultant du
lemme 1.1.1. On suppose que $p\≥3$ et
on reprend les notations de la d\émonstration du lemme  pr\éc\édent.\\
\noi Par r\écurrence on a :
$$\goth B"_1=\{
x \in   
\goth B_1\ |\ \hbox{ tel que}\
\prod_{1\≤k\≤p-1}F^{\goth B}_k(x)\not=0\}=N_1.W^{\goth B}\quad
(1)$$
avec 
$W^{\goth B}=\{x_1+...+x_{p-1}\ ,\ x_i\in E'_2(H_i)\cap
\goth g_1\ ,\
  [x_i,x_j]=0 \ ,i,j=1,...,p-1\}.$ 

\noi Par le lemme 1.1.1 appliqu\é \à $h=H_1+...+H_{p-1},$ pour $x\in
\goth
g"_1,$ il existe $n\in N_2,y\in \goth B'_1$ et
 $x_p\in E'_2(H_p)\cap
\goth g_1$ tels que $[y,x_p]=0$ et $x=n(y+x_p).$

\noi Or pour $k=1,...,p-1$ on a $F_k(x)=F^{\goth B}_k(y)$ donc
$y\in \goth B"_1$ et par (1) $\exists n_1\in N_1$ et $z\in
W^{\goth B}$ tels que $y=n_1z$ d'o\ù
$x=n(n_1z+x_p)=nn_1(z+n_1^{-1}x_p)=nn_1(z+x_p)$ et $z+ x_p\in
W_{\goth t}.$

\noi Pour finir,  soit  $z\in W_{\goth t},$ il suffit de noter que
l'\équation $nz\in W_{\goth t},$ avec $n\in N_{\goth t}$
admet comme seules solutions le centralisateur de $z$ dans  
$N_{\goth t}.$\fdem

\bigskip
 
 \begin{prop}
  \begin {enumerate} 
 \item  L'action de 
$\bf P_{\goth t}$ dans
$\goth g_1$  (resp.$\goth g_{-1}$) est g\éom\étriquement
pr\éhomog\ène et
$F_1,...,F_p$ (resp.$F^*_1,...,F^*_p$) en sont les
invariants relatifs fondamentaux.

\item On suppose que $\F=\C .$

\noi Soit $M=\cap_{1\≤k\≤p}Ker \chi_k,$ l'anneau $\C[\goth g_1]^M$ des polyn\ômes d\éfinis sur $\goth g_1$  qui sont $M-$invariants est \égal \à l'anneau $\C[F_1,...,F_p].$

 \end {enumerate}

 \end{prop}
 
 \dem On proc\ède comme dans la proposition 7.4 de \cite{mullerJA2} .
  \fdem

\bigskip

\begin{rema} 

 \begin {enumerate} 
\item Pour $x\in \goth g"_1,$ l'orbite  $G.x$    est faiblement sph\érique au sens de \cite{sato4} et $(G,G_x,P_{\goth t})$ est un triplet sph\érique (i.e. $\overline {P_{\goth t}}$ a une orbite Zariski-ouverte dans $\overline G.x$).\\

\item Dans le cas r\éel (i.e. $\F=\R$), soit $\theta$  une involution de Cartan  de $\goth g$   telle que $\theta/\goth a=-Id.$\\
\noi Lorsqu'il existe $x\in W_{\goth t}$ tel que $(x,2H_0,\theta (x))$ soit  un sl$_2$-triplet, soit 
$\theta_x=\theta_{x,2H_0}(-1)$ l'involution associ\ée alors $\theta_x$ et $\theta$ commutent et la sous-alg\èbre  parabolique $\goth p_{\goth t}$    est  $\theta_x \theta$-stable minimale au sens de Van den Ban puisqu'on peut prendre un ordre sur $\∆_0$ pour lequel:
$$\{\lambda\in \∆_0\ |\ \lambda \succ 0 \}\subset \∆_0^+.$$
 
\noi Rappelons que :

\noi $\bullet$ Une condition suffisante pour que $W_{\goth t,\theta}=\{x\in W_{\goth t}\ |\ (x,2H_0,\theta (x)$ soit $1$  \Sl  $\}\not=\emptyset$ est  l'existence de racines fortement orthogonales : $\{\lambda_1,...,\lambda_n\}$ de $\∆_1$ et d'une partition $I_1,...,I_p$ de $\{1,...,n\}$ telles que $H_i=\sum_{j\in I_i}h_{\lambda_j}$ pour $i=1,...,p$ (cf.lemme 1.17,1 p.18 de \cite{boppruben}), cette condition est toujours v\érifi\ée dans le cas commutatif (cf. lemme 2.3.2).

\noi $\bullet$ Lorsque le \PV est commutatif on a
$G^{\theta_x}=G_x$ (prop.4.4.5 p.$146$ de \cite{rublivre}) .

\noi  Ce r\ésultat est en g\én\éral faux lorsque le \PV n'est pas commutatif.

\noi Par exemple, lorsque le sous-espace vectoriel $U=Ker(adx:\g_2\mapsto \g_3)$ n'est pas r\éduit \à $\{0\}$ (ce qui est toujours notre cas lorsque $\g_2\not=\{0\}$ puisque $U=\g_2$)  le sous-espace vectoriel $V=ad(x^{-1})^2(U)$  est inclus dans $(\g_0)_{\theta_x}=\{z\in \g_0\ |\ \theta_x(z)=z\}$ car les \elts non nuls de $U$ sont primitifs de poids $4$  mais l'application $ad(x^{-1})^4: \g_2\mapsto \g_{-2}$ est injective.\\

\item Lorsque $\goth t$ est une sous-alg\èbre de l'alg\èbre de Lie engendr\ée par $\g_{\pm 2},$ ce qui est toujours possible dans les cas ``classiques'' (i.e. $(\overline\g_0,\overline\g_1)$ de type $(B_n$ ou $C_n$ ou $D_n,\alpha_k)$ avec $2\≤k<n$), le sous-groupe parabolique $P'_\goth t,$ projection de $P_\goth t$ sur $\Gamma=G/H$ (cf.2. de la remarque 1.4.1), est un \sgp de $\Gamma$. Les \irfs $F_1,...,F_p,$ peuvent \être obtenus \à partir des \irfs associ\és \à un \PV de type commutatif (cf.les travaux de A.Mortajine dans \cite{mortajine} notamment le iii) du th\éor\ème 4.1.1), \également \à partir des repr\ésentations d'alg\èbres de Jordan (cf.\cite{clerc} et dans ce cas $H$ est compact).

\noi Dans le cas complexe par exemple, $\Gamma$ op\ère sur $\C[\g_1]^M$ avec multiplicit\é $1$ lorsque  $\goth t$  a pour dimension   le rang de $\Gamma$ (ce qui est toujours possible cf.prop.2.4.3), on retrouve le th\éor\ème 4.2.2 de \cite{mortajine} en appliquant la m\ême \demo que dans le th\éor\ème 5.3.1 de \cite{rublivre}.
 \end {enumerate}

\end{rema}

\bigskip

\noi Les diff\érents invariants relatifs fondamentaux et les
caract\ères correspondants v\érifient encore les relations
bien connues dans le cas des \PV commutatifs:  

\begin{lem}
Pour $k=1,...,p$ on a:
$\chi^*_{p-k} =\chi_k. \chi_p^{-1}$ en tant que caract\ères
de $G_{h_k},$ et pour
$x\in E_2(h_k)\cap \goth g_1$ et
$y\in
 E_0(h_k)\cap \goth g_1$  qui commutent,
  $F_p(x+y) $ est proportionnel \à
$F_k(x) P_{p-k}(y),$ \à l'exception 

\begin{enumerate}

\item du cas $(\∆,\∑_1)=(C_n,\alpha_k)$ avec $\goth t=\F(\sum_{1\≤i\≤k}h_{\epsilon_i-\epsilon_{k+i}})+\F(\sum_{1\≤i\≤k}h_{\epsilon_i+\epsilon_{k+i}}),$

\item des $\F-$formes du cas
exceptionnel $(E_7,\alpha_6)$ 
\end{enumerate}

\noi o\ù les relations deviennent :
$\chi^*_1 =\chi_2^{-2}. \chi_1$ et $F_2^2(x+y) $ est proportionnel \à
$F_1(x) P_1(y)$ pour $x\in
E_2(h_1)\cap \goth g_1$ et
$y\in
 E_0(h_1)\cap \goth g_1$  qui commutent.
 
 \bigskip
   
\end{lem}

\dem 1) Pour $k=1,...,p-1,$ soient $d_k$ le degr\é de $F_k,$ $d'_k$
celui de $P_{p-k}$  et $p_{i,j}(h_k)$  la dimension du
sous-espace : $E_{i,j}(h_k)=E_i(h_k)\cap E_j(2H_0-h_k).$

\noi Introduisons \également les  rationnels suivants:
$$ N=\frac{ \hbox {dim}(\goth g_1)} { d_p}
  ,\quad  m_k=\frac{1} {Nd_k} ( p_{2,0}(h_k) +\frac{
p_{1,1}(h_k)}{  2})\ ,\quad  m'_k=\frac{ 1} {Nd'_k}( p_{0,2}(h_k) +\frac{ p_{1,1}(h_k)}{2})\ ,$$
donc $m_kd_k+m'_kd'_k=d_p.$\\

\noi Par  unicit\é des
invariants relatifs fondamentaux , on a pour
$(i,j)\not =(0,0) $ et $g\in G_{h_k}$ :
 $$(det(g/_{ E_{i,j}(h_k)})^2=\chi_k(g)^{\frac{i p_{i,j}(h_k) }
{d_k}}
 \psi(g)^{\frac{j p_{i,j}(h_k) }{ d'_k}}\ , \quad\psi\ \hbox{
\étant le caract\ère associ\é \à}\ P_{p-k} 
 $$ 
\noi  donc   on obtient:
$$(det(g/_{ \goth
g_1})^2=\chi_k(g)^{ 2Nm_k}.
 \psi(g)^{ 2Nm'_k }
=\chi(g)^{ 2N}\quad (1)$$
\noi d'o\ù pour  $x\in  E_2(h_k)\cap  \goth g_1$ et $y\in
 E_0(h_k)\cap \goth g_1$
 tels que $[x,y]=0,$   
$F(x+y)^{2N}$ et $F_k(x)^{2Nm_k }P_{p-k}(y)^{2Nm'_k}$
sont proportionnels puisque
 dans le cas alg\ébriquement clos  l'ensemble
$\{x+y\ |\ x\in  E'_2(h_k)\cap  \goth g_1\ ,$ $\ y\in
 E'_0(h_k)\cap  \goth g_1\ |\ [x,y]=0\}$ est inclus dans
une seule orbite de
$G_{h_k}.$

\bigskip

2) Par dualit\é (cf. th\éor\ème 1.4.1 p.17 de \cite{rublivre}) on a : $\chi^*_p=\chi_p^{-1}$  et
$\psi=(\chi^*_{p-k})^{-1}.$\\

3) a) Lorsque $m_k=m'_k=1$ la \demo est termin\ée. Cela convient dans les cas classiques non cit\és dans le lemme (cf. \demo de la prop.1.2.4 puisque dans ces cas $p_{1,1}(h_k)=0$).\\

b) Dans les cas restants, lorsque $d_k+d'_k=d_p$ $(2)$ et que $d_k$ (resp.$d'_k$) est le degr\é de l'invariant
relatif fondamental du \PV absolument irr\éductible : $(\goth U(\goth s)_0,\goth
U(\goth s)_1),$  $\goth s$ \étant la tds engendr\ée par  le
$sl_2-$triplet $(y,2H_0-h_k,y^{-1})$ (resp.$(x,h_k,x^{-1})$)
  ainsi $F_k$ (resp. $P_{p-k}$) est
\également l'invariant relatif fondamental du \PV :
 $(\goth U(\goth s)_0,\goth U(\goth s)_1)$ par cons\équent,
 pour $(x,y)\in \{x+y\in E_2(h_k)\cap \goth g_1\oplus
E'_0(h_k)\cap \goth g_1\ |\ [x,y]=0\},$ il existe un entier
$n$ et une fonction, not\ée $c,$ tels que
$F(x+y)=c(y)F_k^n(x)$  d'o\ù $n=m_k$
donc
$m_k\in \N^*,$ de m\ême $m'_k\in \N^*,$ d'o\ù par la
relation (2) on a $m_k=m'_k=1$ et on applique a) (3).\\
 
c)  Il reste les cas  \énum\ér\és dans le lemme et pour lesquels on
a $d_2=d_1=d'_1,$  le calcul donne $m_1=m'_1=\frac{1}{2}$ (cf.d\émonstration de la prop. 1.2.4 et du lemme 1.3.1).\fdem

\bigskip

\begin{rema}

Les \PVs introduits dans ce paragraphe v\érifient une propri\ét\é de descente que l'on peut \énoncer ainsi.\\

\noi $\bullet$  $\forall k=1,...,p-1,$ l'alg\èbre $\goth U^+(k)$ 
a les m\êmes propri\ét\és que $\goth g:$ 

\begin{enumerate}
\item elle est absolument simple engendr\ée par $\goth U^+(k)_{\pm 1},$

\item $(\goth U^+(k)_0,\goth U^+(k)_1,\frac{h_k}{2})$ est un PV absolument irr\éductible et r\égulier,

\item lorsque $k\geq 2,$ $(H_1,...,H_k)$ sont des \elts $1$-simples tr\ès sp\éciaux de la sous-alg\èbre d\éploy\ée maximale, $\goth a\cap \goth U^+(k),$ et de somme $h_k.$

\end{enumerate}
\noi $F_1,...,F_k$ sont les invariants relatifs fondamentaux du PV: $(P^{\goth U^+(k)}(H_1,...,H_k),\goth U^+(k)_1)$, $P^{\goth U^+(k)}(H_1,...,H_k)$ \étant le sous-groupe parabolique de $G_{h_k}$ associ\é \à $H_1,...,H_k$ (avec cet ordre), il est donn\é par: 
$$P^{\goth U^+(k)}(H_1,...,H_k)=  \{
\begin{array}{c c}
G_{\goth t} &  \  \hbox{pour}\ k=1 \\
G_{\goth t} exp(ad(\oplus_{1\leq i<j\leq k}(E_{i,j}(-1,1)\oplus E_{i,j}(-2,2))\ )) &  \  \hbox{pour}\ k\≥2
\end{array}$$
\noi  et on a: $P^{\goth U^+(k)}(H_1,...,H_k)\subset P(H_1,...,H_p)_{h_k}.$\\

\noi $\bullet$ Soit $(x,H_1,x^{-1})$ un {\Sl} et $\goth s$ l'alg\èbre qu'il engendre alors :

\begin{enumerate}
\item $(\goth U(\goth s)_0,\goth U(\goth s)_1,H_0-\frac{H_1}{2})$ est $1$-irr\éductible et absolument irr\éductible r\égulier pour $p\≥3,$ 

\noi lorsqu'on supposera de plus dans le cas orthogonal :
$(\overline{\g}_0,\overline\g_1)$ de type $$\begin{cases} (D_n,\alpha_k)\text{ que }3k\≤2n-2,\\
(B_n,\alpha_k) \text{ que }3k\≤2n-1.\end{cases}$$
Notons que dans ce cas, lorsque $(\overline{\goth U(\goth s)}_0,\overline{\goth U(\goth s)}_1)$ est encore de type orthogonal, alors il v\érifie les m\êmes conditions (cf.\demo de la prop.1.2.4).

\item Pour $p\≥3,$ $H_2,...,H_p$ sont des \elts $1$-simples tr\ès sp\éciaux de somme $2H_0-H_1$ appartenant \à une sous-alg\èbre ab\élienne d\éploy\ée de $\goth U(\goth s).$
 \end{enumerate}
Le sous-groupe parabolique de $G_{\goth U(\goth s)}$ associ\é \à $H_2,...,H_p$, not\é $P^{\goth U(\goth s)}(H_2,...,H_p),$ v\érifie:
$$P^{\goth U(\goth s)}(H_2,...,H_p)= \begin{cases}
 (G_{\goth t})_x   \  \hbox{pour}\ p=2, \\
(G_{\goth t})_xexp(ad(\oplus_{2\leq i<j\leq p}(E_{i,j}(-1,1)\oplus E_{i,j}(-2,2))\ ))  \  \hbox{pour}\ p\≥3
\end{cases}$$
\noi  et on a: $P^{\goth U(\goth s)}(H_2,...,H_p)\subset P(H_1,...,H_p)_{ \goth s}.$\\

\noi $F_2,...,F_k$ sont des invariants relatifs du PV: $(P^{\goth U(\goth s)}(H_2,...,H_p),\goth U(\goth s)_1),$ fondamentaux lorsque $p\≥3.$
\end{rema}

\bigskip

\bigskip

Dans le but de montrer l'existence d'\équations fonctionnelles pour les fonctions Z\étas associ\ées \à  un\sg parabolique tr\ès sp\écial, $P,$ dans le cas $\goth p-$adique ceci \à l'aide du th\éor\ème $k_{\goth p}$ de F.Sato (p.477 de \cite{sato3}), on commence par \établir les  propri\ét\és n\écessaires portant sur les orbites singuli\ères de $P$.\\

\newpage

\section {\bf  Un r\ésultat sur les orbites des paraboliques tr\ès sp\éciaux}

\bigskip

\subsection{Le r\ésultat }

\bigskip
\noi Les notations sont celles du paragraphe 1.4.\\

\noi Soient $S_{P_{\goth t}}=   \{ x\in {\goth g_{1}} \
 | \ \prod_{1\≤i\≤p}F_i(x) =0\}$ 
(resp. $S^*_{P_{\goth t}}
  = \{ x\in {\goth g_{-1}} \
 | \ \prod_{1\≤i\≤p}F^*_i(x) =0\}),$  et $N=exp(ad(\oplus_{\alpha\in
\∆_0^+}\ \goth g^{\alpha})).$  \\
\begin{theo}
\begin{enumerate}
\item $P_{\goth t}$ a un nombre fini d'orbites dans
$\goth g_1$ et $\goth
g_{-1}.$

\item Pour tout
$x\in S_{P_{\goth t}}
  $  (resp. $S^*_{P_{\goth t}}$)
il existe $i\in \{1,...,p\}$  
tel que
$\chi_i /( (P_{\goth t})_x)^0\not=1$
 (resp.
$(\chi^*_i(/(P_{\goth t})_x)^0\not=1).$
\end{enumerate}
\end{theo}

\bigskip

\noi Il a d\éj\à \ét\é \établi  pour  quelques exemples de \PV commutatifs (cf.\cite{sato3} et \cite{hironaka2}).\\

\dem Par r\écurrence sur la dimension de $\goth t$ mais on remplace le point $2)$ par :

\begin{rien} Pour tout
$x\in S_{P_{\goth t}}
  $  (resp. $S^*_{P_{\goth t}}$)
 il existe  $z\in
 P_e.x$ et
$u\in P( \check{\goth g}), $ commutant \à $z$ et \à $\goth t,$ 
 tel que
$B(u, .)/_{\goth t}\not=0,$   $P_e$ \étant le groupe engendr\é par $N$ et $h_v(t)$ pour $v\in P(\check{\goth g})$ et $t\in \F^*.$ 

\end{rien}

\noi qui implique 2),   en effet
soit   
  $i_0=$inf$\{i\in
\{1,...,p\}\ |\ B(u,H_i)\not=0\},$ l'\él\ément
$h_u(t)\in (P_z)^0$ et v\érifie $\chi_{i_0}(h_u(t))=t^a$ avec
$a\not=0$ (cf.   lemme 1.2.1).  \\

\noi On note  $ P_{\goth t}=P(H_1,...,H_p).$

\noi  Lorsque $p=1,$ on a $ P_{\goth t}=G$ et le r\ésultat est
\évident (pour 2.1.2 cf.6. des notations).\\

\noi  Soit  
 $x$ non nul
appartenant \à
$\goth g_1$  que l'on  d\écompose
suivant
$ad(H_1):$  $x=x_2+x_1+x_0$ avec $x_i\in E_i(H_1)\cap
\goth g_1$ pour $i=0,1,2.$ \\

1) Lorsque $x_2$ admet
$H_1$ comme \él\ément $1$-simple (i.e. $F_1(x)\not=0$), il
existe $g\in N$ tel que  
$ x'=g(x)=x_2+y_0$ avec  
$y_0\in \goth U(\goth s)_1,$ $\goth s$ \étant l'alg\èbre
(r\éductive dans $\goth g$) engendr\ée par le $sl_2$-triplet
 construit avec  $x_2$ et $H_1$  (lemme 1.1.1).

\noi Le point 1) d\écoule  de l'hypoth\èse
de r\écurrence appliqu\ée au \PV: $(\goth U(\goth s)_0,\goth
U(\goth s)_1)$ avec le parabolique associ\é \à
$(H_2,...,H_p)$ (cf.remarque 1.4.8).

\noi Lorsque $F_i(x)=0,$   on a
$ F_1(x_2)F_{i-1}^{\goth U(\goth s)}(y_0)=0$ (lemme 1.4.7
appliqu\é au \PV: $(\goth U^+(i)_0,\goth U^+(i)_1)$ et le r\ésultat
d\écoule \à nouveau  de l'hypoth\èse de r\écurrence appliqu\ée au
\PV:
$(\goth U(\goth s)_0,\goth U(\goth s)_1)$ avec le
parabolique associ\é \à
$(H_2,...,H_p).$  

\noi Notons que $f.P(\check{\goth U(\goth s) })
\subset  P(\check{\goth g }),$ $f$ \étant l'indice de
connexion du syst\ème de racines  de $\goth U(\goth s).$\\
 
 2) Lorsque $F_1(x)=0,$ on montre le th\éor\ème cas
par cas suivant le type de
syst\ème de racines dans les paragraphes qui suivent
en omettant le cas $x_1=x_2=0$ puisque le point  1) d\écoule
de l'hypoth\èse de r\écurrence appliqu\ée au \PV: $(
\goth U(\F H_1)_0, (  \goth U(\F H_1)_1)$ avec le
parabolique associ\é \à
$H_2,...,H_p$ et $u=H_1$ convient pour le point 2).\\

\noi Notons qu'il suffit
de prouver 1) pour  tous les 
 paraboliques tr\ès sp\éciaux et minimaux  pour l'inclusion,  
 \à l'action de $G$ pr\ès, ceci dans le cas
alg\ébriquement clos car
$\F$ est un F-corps (th\éor\ème 5, chapitre III,paragraphe 4.4 de
\cite{serre}). \\
 
  \noi  Les r\ésultats pour
$\goth g_{-1}$ s'obtiennent de mani\ère analogue et la
d\émonstration est omise.\fdem \\ 

\noi La d\émonstration restante consiste \à construire des repr\ésentants plus simples des \él\éments de $S_{P_{\goth t}}$  (cf.2.1.2) et pour ceci on  regroupe  quelques r\ésultats techniques dans le paragraphe suivant avant d'effectuer la v\érification cas par cas.\\

\subsection{\bf  Lemmes techniques ($\bf \F$ est de caract\éristique $\bf 0$)}
  
  \noi Soit $x\in \goth g_1,$ on note $s(x)=\{\mu \
|\ X_{\mu}\not=0\},$    
$x=\sum_{ \mu\in \∆_1 } X_{\mu}$ d\ésignant la
d\écomposition de $x$  suivant les sous-espaces radiciels.\\

\begin{lem}
Soit $(\goth g_0, \goth g_1)$ un \PV \  pour lequel $\goth g_i=\{0\}$ pour $|i|\≥3.$

Soit  $x$ non nul dans
$\goth g_1,$ tel que s(x) contienne une racine $\mu_0$
de
$\∆_1,$ longue dans $\∆_1$ et de hauteur minimale parmi les
racines de s(x)
 alors il existe au plus une racine
$\mu_1\in \∆_1$ telle que 

  $\bullet\ n(\mu_0,\mu_1)=-1$

 $\bullet \ N.x\cap\ (\goth g^{\mu_0}\oplus
\goth g^{\mu_1}
  \oplus\ E_0(h_{\mu_0})\cap \goth g_1\ ) \not= \emptyset.$
 \end{lem}

\dem 

\noi Soit $\beta$ une racine  
longue dans
$\∆_1,$ on rappelle que $E_2(h_{\beta})\cap \goth g_1=\g^{\beta}$ et que pour $\alpha\in \∆_0$  on a $\alpha + \beta\in \∆\  \Leftrightarrow\
n(\alpha, \beta)<0.$\\

1) Prenons $g=exp(ad(A)),$ avec $A=\sum_{\mu\in
s(x)\ |\ n(\mu,\mu_0)=1}X_{\mu-\mu_0},X_{\mu-\mu_0}$
choisi tel que $[X_{\mu-\mu_0},X_{ \mu_0}]=-X_{\mu } ,$
on
v\érifie facilement que
 $s( y) \subset \{ \mu\in \∆_1\ | \quad n(\mu,\mu_0)
\≤0\}\cup
\{\mu_0\}$ avec $y=g(x).$ 

2) Soit $s_{-1}=\{\mu\in s( y)\ |\
n(\mu,\mu_0)=-1\}$ et $\mu\in s_{-1}.$ Comme $\∆_3=\emptyset,$
la consid\ération de la $\mu$-chaine d\éfinie par $\mu_0$
permet de montrer que toutes les racines de $s_{-1}$ sont
\également longues dans $\∆_1;$ de plus soit $\mu'\in s_{-1},$
 comme
$\mu_0+\mu+\mu'\notin \∆$  on a
$0\≤n(\mu_0+\mu,\mu')=n(\mu_0,\mu')+n(\mu,\mu')$ d'o\ù
 $n(\mu,\mu')\≥1.$

Lorsque 
$s_{-1}$ a au moins $2$ \él\éments,  soit
$\mu_1$ la racine de hauteur minimale de $s_{-1}$ et soit
$g'=exp(ad(B))$ avec
$B=\sum_{\mu\in s_{-1}-\mu_1}X_{\mu-\mu_1},X_{\mu-\mu_1}$
choisi tel que
$[X_{\mu-\mu_1},X_{ \mu_1}]=-X_{\mu } ,$  alors
  $g'(y)\in \goth g^{\mu_0}\oplus
\goth g^{\mu_1}
  \oplus\ E_0(h_{\mu_0})\cap \goth g_1.$ \fdem\\ \\

\noi Indiquons deux r\ésultats \él\émentaires lorsque  $(\goth g_0, \goth g_1)$ un \PV \  pour lequel $\goth g_i=\{0\}$ pour $|i|\≥3.$

\noi Dans les $2$ lemmes qui suivent, $H_1,...,H_p$ sont simplements des \él\éments de $\goth a$ et on note toujours $\goth t=\oplus_{1\≤i\≤p}\F H_j.$

\noi $N(H_1,...,H_p)$ est le sous-groupe engendr\é par $\oplus_{\{\alpha \in \∆_0|\alpha\succ 0\}}\goth g^{\alpha},$ avec $\{\alpha\succ 0$ \ssi \ $\alpha (H_{j_0})<0\}$ et $j_0=$inf$\{j|\alpha (H_j)\not=0\}$ 
et $P_e(H_1,...,H_p)$ d\ésigne le groupe engendr\é par $exp(ad(A)),$ $A$ \elt\ nilpotent de $\goth g_0$ commutant \à $\goth t$ ou bien $A\in \oplus_{\{\alpha \in \∆_0|\alpha\succ 0\}}\goth g^{\alpha}.$  \\

\begin{lem} Soit  
  $(x,h,x^{-1})$  un  $sl_2-$triplet tel que:

\begin{enumerate} 

\item  h et  $\goth t$ commutent,

 \item $x\in E_2(H_1)\cap \goth g_1$ 

\end{enumerate}

pour $i\in \Z$ soient $y_i\in E_i(h)\cap E_0(H_1)\cap \goth g_1$ alors $x+y_{-1}+y_0\in N(H_1,...,H_p)(x+\sum_{i\in \Z}y_i) $ et $ N(H_1,...,H_p)(x+\sum_{i\in \Z}y_i)\cap (x+y_{-1}+\goth U(\goth s)_1)\not=\emptyset,$ $\goth s$ \étant l'alg\èbre engendr\ée par $x,h,x^{-1}$ et $H_1.$

  \end{lem}

 \dem  Comme : $E_{-i}(h)\cap \goth g_1=\{0\}$ pour $i\≥2,$ on a   $ y=\sum_{i\in \Z}y_i=y_{-1}+y_0+\sum_{i\≥1}y_i\     
  ,$
 or pour $i\≥1$ il existe $A_i \in E_{i-2}(h)\cap 
E_{-2}(H_1)\cap
\goth g_0$ tel que $[A_i,x_2]=-y_i$  d'o\ù $exp(ad(\sum_{i\≥1}A_i))(x_2+y)=x_2+\sum_{i\≥1}[A_i,x_2]+y=x_2+y_{-1}+y_0.$

\noi Comme $[x,y_0]\in E_2(h)\cap E_0(H_1-h)\cap \goth g_2,$ il existe $A\in  E_{-2}(h)\cap E_0(H_1-h)\cap \goth g_0\subset E_{-2}(H_1)\cap \goth g_0(\subset \oplus_{\alpha\in \∆_0^+}\goth g^{\alpha}$) tel que $ad(x)^2(A)=[x,y_0]$ et on a $exp(ad(A))(x+y_{-1}+y_0)=x+y_{-1}+y'_0$ avec $y'_0=[A,x]+y_0$ d'o\ù $[x,y'_0]=0.$\fdem\\

\noi On continue \à simplifier les repr\ésentants des orbites:

 \begin{lem}
Soient
$(x,h,x^{-1})$ et $(y,h',y^{-1})$ 2 $sl_2-$triplets tels que

\begin{enumerate} 

\item  h,   h' et  $\goth t$ commutent,

 \item $x\in E_2(H_1)\cap \goth g_1$ et $y\in E_0(H_1)\cap
E_{-1}(h)\cap \goth g_1,$

\item $ \oplus_{i\≥2} E_{-i}(H_1-h)\cap \goth g_1=\{0\}.$

\end{enumerate}

\noi et soit $z\in E_0(h)\cap E_0(H_1)\cap \goth g_1$ tel que
$[[z,y^{-1}],\goth t]=0$ alors   
$P_e(H_1,...,H_p) (x+y+z+E_2(H_1)\cap E_3(h)\cap \goth g_1)$ rencontre  $x+y+ 
E_0(H_1)\cap E_0(h)\cap E_0(h')\cap \goth g_1+ E_2(H_1)\cap E_3(h)\cap \goth g_1.$
 
 \end{lem}
 
 \dem Comme $[y,E_{-1}(h')\cap E_0(h)\cap \goth g_1]\subset
E_{-1}(h)\cap \goth g_2=\{0\}$ et que pour $i\≥2$ on a
$ad(y^{-1})^i(E_i(h')\cap E_0(h)\cap \goth g_1)\subset
E_{-i}(h')\cap E_i(h)\cap \goth g_{1-i}=\{0\},$ $z$ admet la
d\écomposition suivante relativement \à $ad(h'):$ $z=z_0+z_1$
avec $z_i\in  E_i(h')\cap   E_0(H_1)\cap E_0(h)\cap
\goth g_1 $ pour $i=0,1.$ 

\noi Soit $v\in  E_2(H_1)\cap E_3(h)\cap \goth g_1$  et soit $A=[z_1,y^{-1}],$ alors $[A,\goth t]=0$ d'o\ù $exp(ad(A))\in P_e(H_1,...,H_p)$  et $[A,y]=-z_1$  donc  $exp(ad(A))(x+y+z+v)=x+y+z_0 +[A,z_0]+\frac{1}{2}[A,z_1]+([A, x]+v )$ par 3).
 
\noi On termine en appliquant la d\émonstration du   lemme  pr\éc\édent puisque 
$$[A,x]+v \in E_2(H_1)\cap E_3(h)\cap \goth g_1\ ,\ [A,z_0]+\frac{1}{2}[A,z_1] \in E_0(H_1)\cap E_1(h)\cap \goth g_1$$
\noi et que $[E_{-2}(H_1)\cap E_{-1}(h)\cap \goth g_0,E_0(H_1)\cap E_1(h)\cap \goth g_1\oplus  E_2(H_1)\cap E_3(h)\cap \goth g_1]=\{0\}$ par 3).\fdem\\

\noi{\it Remarque}: L'hypoth\èse $3)$ du lemme est v\érifi\ée lorsque $H_1-h$ est $1$-simple et  $[[z,y^{-1}],\goth t]=0$  pour $p=2.$\\

\noi Ce dernier lemme  utilise des descriptions faites ult\érieurement.\\

\begin{lem}  Soient  
$P(H_1,...,H_p) $ un sous-groupe parabolique standard tr\ès sp\écial , 
$S=\{\lambda_1,...,\lambda_q\}$ $q(\≥1)$ racines  
fortement orthogonales de $\∆_1,$   longues dans
$\∆_1,$ telles que
$\lambda_i(H_1)=2$ pour
$i=1,...,q$ et soit
$\∆_{-1}(h)=\{\mu\in \∆_1\ |\ \mu (H_1)=0$ et
$ \mu (h)=-1\}\not= \emptyset $ avec
$h=\sum_{1\≤i\≤q}h_{\lambda_i}\not=H_1.$

 Pour  $x_2$ de support $S$ et $y\in E_0(H_1)\cap
\goth g_1,$  
  il existe $r(\≤q)$ racines fortement
orthogonales, $\mu_1,...,\mu_r,$
    de   $\∆_{-1}(h)$ telles que pour $i=1,...,r$ on
a $n(\mu_i,\lambda_i)=-1$ (\à la num\érotation des
$\lambda_i$ pr\ès) et  
$P_e(H_1,...,H_p)(x_2+y )\cap (\oplus_{1\≤j\≤q}
\goth
g^{ \lambda_j}\oplus_{1\≤i\≤r}
\goth
g^{ \mu_i}\oplus E_0(H_1)\cap E_0( h)\cap
\goth g_1)\not=\emptyset$

\noi (pour $p\≥3 \ :\ P_e(H_1,...,H_p)(x_2+y )\cap ( x_2 
 \oplus_{1\≤i\≤r}
\goth
g^{ \mu_i}\oplus E_0(H_1)\cap E_0( h)\cap
\goth g_1)\not=\emptyset).$

\end{lem}

\dem   D'apr\ès le lemme  2.2.1, on peut supposer que $y=y_{-1}+y_0$ avec $y_i\in E_i(h)\cap E_0(H_1)\cap \goth g_1$ pour $i=-1,0$    et $y_{-1}\not=0.$  \\

\noi Comme $E_1(h_{\lambda_1})\cap
E_{-1}(h)\cap \goth g_1\subset
E_{-2}(h_{\lambda_2}+...+h_{\lambda_q})\cap \goth g_1=\{0\}$(cf.2) de la remarque 1.2.3)
on a: 
$$  \∆_{-1}(h)= \cup_{1\≤i\≤q}\∆_{-1}(h_{\lambda_i})  \ \hbox{avec}$$
$$\∆_{-1}(h_{\lambda_i}) =  \{\mu\in
\∆_1\ |\ \mu(H_1)=0\ ,\ n(\mu,\lambda_i)=-1\ ,\
n(\mu,\lambda_j)=0\ \hbox{pour}\ 1\≤j\not=i\≤q\},  $$
par cons\équent,  en proc\édant comme dans la d\émonstration du lemme  2.2.1 et en changeant \éventuellement l'indexation des $\lambda_i,i\≥2,$ on peut supposer qu'il existe
$\gamma_1\in \∆_{-1}(h_{\lambda_1})$ tel que:
$$ y_{-1}=X_{\gamma_1}+z_0\ \hbox{avec}\ z_0\in
E_0(h_{\lambda_1})\cap E_{-1}(h)\cap E_0(H_1)\cap \goth
g_1.$$
\noi Lorsque $q=1,$ la d\émonstration est termin\ée. \\

\noi Pour $q\≥2$ 
 et $z_0\not=0,$ on continue la r\éduction en
d\écomposant $z_0$ relativement \à $ad(h_{\lambda_2})$:
$$z_0=t_{-1}+t_0\ ,\ t_i\in 
E_0( h_{\lambda_1})\cap E_i( h_{\lambda_2})\cap
E_{-1}(h)\cap E_0(H_1)\cap
\goth g_1\ ,\ i=-1,0,$$
et on peut supposer que $t_{-1}\not=0$ en changeant \éventuellement l'indexation des $\lambda_i,i\≥2.$\\
 
\noi Pour $\gamma\in s(t_{-1})$ on a $n(\gamma,\gamma_1)\in \{0,1\}$ et soit $a=$Min$\{n(\gamma,\gamma_1),\gamma\in s(t_{-1})\},$ prenons $\gamma_2$ de hauteur minimale dans $A=\{\gamma\in s(t_{-1})\ |\ n(\gamma,\gamma_1)=a\}$. 

\noi  Pour $\gamma\in s(t_{-1})-\{\gamma_2\}$  et $\alpha\in \{\lambda_i,i=1,...,q,\gamma_1,\gamma'\in s(t_{-1})-\{\gamma_2\}\}$, $ \gamma-\gamma_2+\alpha \not\in \∆$ puisque $n( \gamma-\gamma_2,\alpha )\≥0$ et que $\alpha$ est longue dans $\∆_1;$ notons que $\goth g^{ \gamma-\gamma_2  }$ normalise les sous-espaces $V=E_0(h)\cap E_0(H_1)\cap
\goth g_1$ et $W=E_0( h_{\lambda_1})\cap E_0( h_{\lambda_2})\cap
E_{-1}(h)\cap E_0(H_1)\cap
\goth g_1$.

\noi Ainsi:\\

$\bullet$ Lorsque $p=2$, on fait op\érer  $exp(ad(\oplus_{\gamma\in
s(t_{-1})-\gamma_2}\goth g^{\gamma-\gamma_2}))(\subset G_{\goth t})$ pour se ramener \à $t_{-1}=X_{\gamma_2}.$  \\

$\bullet$ $\bullet$ Lorsque $p\≥3$,  $\∆$ est de type classique
$B_n,C_n,BC_n,D_n$ avec $\Sigma_1=\{\alpha_k\}.$  On a
$k\≤n-1$ dans les cas $B,C$ ou $BC$   et
$k\≤n-2$ dans le cas
$D$ car  
$\∆_{-1}(h)\not=\emptyset.$

\noi On reprend les notations du cas classique ($\S 2.4$).

\noi Par la
  proposition     2.4.3, il existe: $l_1,...,l_p$ tels que  
 $ l_{p+1}=0<l_p< ...<l_ 2<l_1=k$ tels que  
  $H_i=\sum_{l_{i+1}+1\≤j\≤l_i}h_{\epsilon_j}$  pour
$i=1,...,p.$ Ainsi, pour $j=1,...,q,$  il existe $p_j\in\{l_2+1,...,k\}$ et $q_j\in\{ k+1,...,n\}$ tels que $\lambda_j=\epsilon_{ p_j}-(\pm)\epsilon_{q_j},$ et
 $$\{\lambda\in \∆_1\ |\ \lambda (H_1)=0,\lambda (h_{\lambda_i})=-1\}=\{\epsilon_m+(\pm)\epsilon_{q_i},1\≤m\≤l_2\}.$$

\noi Si  $\∆_{-1}(h_{\lambda_i})\not=\emptyset,$ on a $q_j\not= q_i$ pour $1\≤j\not= i\≤q $ donc par orthogonalit\é des racines $\lambda_j$ on a $p_j\not= p_i$   (d'o\ù $\frac{\lambda_j-\lambda_i}{2}\notin \∆$) pour $1\≤j\not= i\≤q .$\\

\noi On fait op\érer  $exp(ad(\oplus_{\gamma\in
 A-\gamma_2}\goth g^{\gamma-\gamma_2}))(\subset N)$ pour se ramener \à $t_{-1}=X_{\gamma_2}+u.$
 
 \noi Lorsque $u\not=0,u\in \oplus_{\gamma\in
 B }\goth g^{\gamma},$ avec $B=\{\gamma\in s(t_{-1})\ |\ n (\gamma,\gamma_1)=1\}\not= \emptyset$ et $n(\gamma_2,\gamma_1)=0;$ il existe $m_1\≤l_2$ tel que $\gamma_1=\epsilon_{ m_1}+(\pm)\epsilon_{q_1}$ donc  $B=\{
\delta=\epsilon_{ m_1}+(\pm)\epsilon_{q_2}\}$ d'o\ù  $\delta-\gamma_1$ s'annule sur $\goth t.$

\noi Prenons $B\in \goth g^{\delta-\gamma_1}$ tel que $[B,X_{\gamma_1}]=-u$ et $C\in \goth g^{\mu}$ avec $\mu=\delta-\gamma_1+\lambda_2-\lambda_1$ tel que $[C,x_2]=-[B,x_2] $ et soit $g=exp(ad(C))exp(ad(B))(\in G_{\goth t})$  alors $$g(x_2+y)=x_2
+X_{\gamma_1}+X_{\gamma_2}+t_0+y'_0\ ,\ t_0\in W,y'_0\in V  $$  

\noi car $\goth g^{\delta-\gamma_1}$ commute \à $\oplus_{1\≤i\not=2\≤q}\goth g^{\lambda_i}\oplus \goth g^{\gamma_2}\oplus \goth g^{\delta}\oplus W\ $ et normalise $V,$ que $\goth g^{\mu}$ commute \à $\oplus_{2\≤i\≤q}\goth g^{\lambda_i}\oplus \goth g^{\mu+\lambda_1}\oplus\goth g^{\gamma_1} \oplus \goth g^{\gamma_2}\oplus W $ et normalise $V.$   

\noi On peut donc toujours supposer que $t_{-1}=X_{\gamma_2}.$ \fdem de $\bullet$ $\bullet$\\

\noi On revient au cas g\én\éral. Lorsque $q=2$ la d\émonstration est finie.\\

\noi Remarquons que, lorsque $n(\gamma_2,\gamma_1)=1,$ on peut supposer que $\gamma_2-\gamma_1\in \∆_0^+ $ (en changeant l'indexation des $\lambda_i$) et en proc\édant comme ci-dessus, il existe $g\in exp(ad(\goth g^{\nu}))exp(ad(\goth g^{\gamma_2-\gamma_1}))(\subset P_e(H_1,...,H_p)$ avec $\nu=\gamma_2-\gamma_1+\lambda_2-\lambda_1,$ tel que $g(x_2+y) =x'_2
+X_{\gamma_1}+ t_0+y'_0 ,$ avec $ t_0\in W,$ $y'_0\in V  $ et $ \{\lambda_i,i=1,3,...,q\}\subset s(x'_2)\subset \{\lambda_i,i=1,...,q\}$ car  $E_2(h_{\lambda_2})\cap \goth g_1=\goth g^{\lambda_2}.$

\noi  Lorsque $\displaystyle{\ \gamma_2-\gamma_1\not= \frac{\lambda_1-\lambda_2}{2}}$ on a $s(x'_2)=s(x_2).$ \\

\noi Lorsque $q\≥3:$ \\

\noi $\bullet$ Dans les cas classiques,  on continue la d\émonstration
par r\écurrence sur $q,$ exactement comme cela a d\éj\à \ét\é
fait, dans le cas
$P(H_1,...,H_p)$ cit\é ci-dessus et dans l'unique cas restant
$P'_0$ (cf.prop.2.4.3).\\
 
 \noi $\bullet$ $\bullet$ Dans les cas exceptionnels, dont la description est rappel\ée
dans le   $\S 2.5,$ on a toujours
$q\≤3.$ De plus,  lorsque
$q=3,$ $\∆$ est   simplement lac\é,  $\goth g$ est d\éploy\ée et
$H_1-\sum_{1\≤i\≤3}h_{\lambda_i}=h_{\lambda_4},$ 
  $\lambda_4$ \étant une racine de $\∆_1$ orthogonale aux
racines de
$S.$   

\noi Dans les
notations pr\éc\édentes, soit
  $\gamma\in s(t_{-1}),$  
  comme
$n(\gamma,\lambda_4)=1$  on a $\gamma+\lambda_2-\lambda_4\in
\∆_1$ or $\gamma_1+\lambda_1\in \∆_2$ donc
$n(\gamma_1+\lambda_1,\gamma+\lambda_2-\lambda_4)=
n(\gamma_1,\gamma)-n(\gamma_1,\lambda_4)=
n(\gamma_1,\gamma)-1\≥0$
d'o\ù $n(\gamma_1,\gamma)=1.$ Ainsi dans ce cas on se ram\ène
toujours \à
$t_{-1}=0$ c'est \à dire que $r=1.$\fdem\\

\subsection{\bf  Le cas commutatif ($\bf\goth
g_2=\{0\}$)}

\noi  Dans tout ce paragraphe   on suppose que $\goth
g_2=\{0\}.$\\
 
 \begin{prop} $(\F$ de caract\éristique $0)$ Pour tout $x\in \goth
g_1$ (resp. $\goth g_{-1}$), il existe un ensemble S de
racines fortement orthogonales de
$\∆_1$  (resp. $\∆_{-1}$)  telles que
$N.x\cap\ (\oplus_{ \mu\in S
 }\goth g^{\mu})\ \not= \emptyset.$ 
 \end{prop}
 
\dem Par r\écurrence sur le
rang de
$\∆,$ le cas de rang $1$  \étant \évident.\\

 \noi  Soit $\mu_0$ une racine de hauteur minimale de $s(x).$
Lorsque $\mu_0$ est longue dans $\∆_1,$ on applique le 
 lemme 2.2.1 puis l'hypoth\èse de r\écurrence  au \PV 
$((\goth U(\F h_{\mu_0})_0,(\goth U(\F h_{\mu_0})_1),$ lorsqu'on prend pour sous-alg\èbre d\éploy\ée maximale: $\goth a\cap \goth U(\F h_{\mu_0})$ et pour ordre sur $\∆^{\bot \mu_0}=\{\alpha\in \∆ |
(\alpha, \mu_0)=0\}$
 qui est le syst\ème de racines associ\é,   l'ordre de $\∆$ c'est \à dire que $(\∆^{\bot  \mu_0})^+=\{\alpha\in \∆^+ |(\alpha,\mu_0)=0\},$ 
 ce qui termine la d\émonstration du cas
simplement lac\é.

\noi Lorsqu'il existe une racine $\mu\in \∆_1$ telle que $s(x)\subset \mu^{\bot}$, il suffit \également d'appliquer l'hypoth\èse de r\écurrence  au \PV 
$((\goth U(\F h_{\mu})_0,(\goth U(\F h_{\mu})_1).$

\noi Dans les cas restants,    $\mu_0$ est  courte et par classification le \PV a   pour syst\ème
de racines gradu\é soit
 
 \begin{enumerate}

 \item  $(C_n,\alpha_n)$ alors
$\∆_1=\{ \epsilon_i+\epsilon_j,i,j=1,...,n\}$ d'o\ù
$2\epsilon_n\notin
s(x)$   et $s(x)\cap
\{\epsilon_i+\epsilon_n,i=1,...,n-1\}\not=\emptyset.$  
Soit
$i_0=sup\{i\ |\ \epsilon_i+\epsilon_n\in s(x)\},$ \à l'aide
d'un \él\ément convenable, $g,$  de
$exp(ad(\oplus_{1\≤i\≤i_0-1 }\goth
g^{\epsilon_i-\epsilon_{i_0}})),$ on peut se ramener \à
$s(g(x))\subset
\{ \epsilon_{i_0}+\epsilon_n, 
\epsilon_i+\epsilon_j, 1\≤i\≤j <n\};$ comme
$ad(X_{\epsilon_{i_0}+\epsilon_n})$ est une surjection de
$\oplus_{1\≤i\≤n-1}\goth g^{\epsilon_{i}-\epsilon_n}(\subset
\oplus_{i=-1,0}E_i(h_{\epsilon_{i_0}+\epsilon_n}))$ sur
$\oplus_{1\≤i\≤n-1}\goth g^{\epsilon_{i}+\epsilon_{i_0}}$, on
peut se ramener au cas o\ù
$s(g(x))\subset \{  \epsilon_{i_0}+\epsilon_n, 
\epsilon_i+\epsilon_j, i,j\not=i_0,n\}$ et on applique
l'hypoth\èse de r\écurrence \à l'\él\ément
$y=g(x)-X_{\epsilon_{i_0}+\epsilon_n}$  dans le PV:
$(\goth U(\goth s)_0,\goth U(\goth s)_1)$ avec  $\goth s=\F h_{ 2\epsilon_{i_0}}\oplus
\F h_{2\epsilon_n } .$\\

 \item $(B_n,\alpha_1)$ alors
$\∆_1=\{\epsilon_1,\epsilon_1\pm\epsilon_j, j=2,...,n\}$
d'o\ù $\mu_0 =\epsilon_1  $ et  $s(x) =
\{\epsilon_1,\epsilon_1+\epsilon_j,j=2,...,n\} .$ 

\noi Comme
$ad(X_{\epsilon_1 })$ est une surjection  de
 $E_0(h_{\epsilon_1})\cap \goth g_0$ sur
$E_2(h_{\epsilon_1})\cap \goth g_1,$  
\à l'aide d'un \él\ément convenable de $
exp(ad(\oplus_{2\≤j\≤n}\goth g^{\epsilon_j })),$ not\é $g,$ on a
$s(g(x))=\{\epsilon_1\}.$ \fdem
\end{enumerate} 
\vskip 5pt

\noi Pr\écisons la forme des sous-groupes paraboliques dans le cas commutatif:\\
  
 \begin{lem}
 \begin{enumerate}
\item  Il existe un unique ensemble maximal ordonn\é de racines fortement orthogonales (longues) de
$\∆_1,$  $\lambda_1,...,\lambda_n,$ telles que le sous-groupe parabolique 

\noi
 $P_0=P( h_{\lambda_1},...,h_{\lambda_n})$ soit standard.

\item  Soit
$P_{\goth t}$ un sous-groupe parabolique  standard
 alors $P_{\goth t}\supset P_0$ et
  il existe  p, $ 2\≤p\≤n,$ et des
entiers:
$ l_0=0<1\≤l_1< ...<l_p=n$ tels que $P_{\goth
t}=P(H_1,...,H_p)$ avec  
$H_i=\sum_{l_{i-1}+1\≤j\≤l_i}h_{\lambda_j}$ pour $i=1,...,p.$
\end{enumerate}
 \end{lem}
 
 \dem
 On rappelle qu'un sous-groupe parabolique $P(H_1,...,H_p),$ avec 
 
 \noi $\goth t=\oplus_{1\≤i\≤p}\F H_i\subset \goth a,$ est standard si $\{\lambda\in \∆_0\ |\lambda \succ 0\}=\{\lambda\in \∆_0^+| \lambda/_{\goth t}\not=0\}$ (cf.3) de la remarque 1.4.1).\\
 
 1) Il est bien connu (par exemple lemme 6.3 de \cite{mullerJA1}) qu'un ensemble maximal de racines fortement orthogonales  de
$\∆_1$ peut \être construit de mani\ère canonique par orthogonalisations successives  en prenant $\lambda_1$ l'unique racine simple de $\∑_1$ puis $\lambda_2$ l'unique racine simple de
 $\∆^{(1)}=\{\lambda\in \∆|\ n(\lambda,\lambda_1)=0\}$ muni de l'ordre $\∆^{(1)+}=\∆^{(1)}\cap \∆^+$ et qui appartient  \à  $\∆_1,$ 
 puis $\lambda_3$ l'unique racine simple de
 $\∆^{(2)}=\{\lambda\in \∆^{(1)}|\ n(\lambda,\lambda_2)=0\}$ muni de l'ordre $\∆^{(2)+}=\∆^{(2)}\cap \∆^+$ et qui appartient  \à  $\∆_1,$ etc. jusqu'\à \épuisement des racines. L'ensemble $\{\lambda_1,...,\lambda_n\}$ ainsi construit est un ensemble maximal  de racines fortement orthogonales  de
$\∆_1$ donc $\sum_{1\≤i\≤n}h_{\lambda_i}=2H_0.$

Soit $\alpha\in \∆_0^+$ ne s'annulant pas sur $\goth t_0=\oplus_{1\≤i\≤n}\F h_{\lambda_i}$ alors par la construction pr\éc\édente  il existe $(i,j)$ tels que $1\≤i<j\≤n$  et $n(\alpha,\lambda_j)=-1=-n(\alpha,\lambda_i)$ donc  
 $P_0=P( h_{\lambda_n},...,h_{\lambda_1})$  est un sous-groupe parabolique standard.
 
 Soit $\{\mu_1,...,\mu_n\}$ un ensemble maximal de racines fortement orthogonales  de
$\∆_1$ tel que le sous-groupe parabolique $P_1=P( h_{\mu_1},...,h_{\mu_n})$ soit standard, $P_0$ et $P_1$ sont conjugu\és par $Aut_0(\goth g,\goth a)_{H_0}$ puisqu'il existe $w\in W_0,$ qui est le groupe de Weyl de $\∆_0,$ tel que pour $i=1,...,n$ on ait $w(\mu_i)=\lambda_{n-i+1}$ (par exemple, prop.4.4,2) de \cite{mullerJA1} avec l'irr\éductibilit\é de\gog). Comme $P_0$ et $P_1$ sont conjugu\és et standards, ils sont \égaux d'o\ù l'\elt $g,$ ant\éc\édent de $w$ dans $Aut_0(\goth g,\goth a)_{H_0}$,  appartient \à $ P_1\cap Aut_0(\goth g,\goth a)_{H_0}$ donc pour $i=1,...,n$ il existe $u_i\in \oplus_{\{\alpha\in \∆_0^+| \alpha\succ 0\}}\goth g^{\alpha}$ tels que $g(h_{\mu_i})=h_{\mu_i}+u_i=h_{\lambda_{n-i+1}}$ d'o\ù $g(h_{\mu_i})=h_{\mu_i}=h_{\lambda_{n-i+1}}.$

2) Pour $i=1,...,p,$ le \PV $(E_0(H_i)\cap \goth g_0,E_2(H_i)\cap \goth g_1,\frac{H_i}{2})$ est commutatif et $H_i$ est $1$-simple; il admet $\∆'^{(i)}=\{\alpha\in \∆\ |\ \alpha (H_j)=0$ pour $1\≤j\not=i\≤p\}$ comme syst\ème de racines lorsqu'on prend pour sous-alg\èbre d\éploy\ée maximale $\goth a\cap \goth s,$ avec $\goth s=\oplus_{1\≤j\not=i\≤p}\F H_j,$ on munit $\∆'^{(i)}$ de l'ordre induit par $\∆,$ alors il existe un
ensemble maximal de racines fortement orthogonales  de
$\∆'^{(i)}_1:$ $\{\mu_1^{(i)},...,\mu_{n_i}^{(i)}\},$ tel que le sous-groupe parabolique correspondant $P( h_{\mu_1^{(i)}},...,h_{\mu_{n_i}^{(i)}})$ soit standard et on a $H_i=\sum_{1\≤j\≤n_i}h_{\mu_j^{(i)}}.$

$\{\mu_j^{(i)},1\≤i\≤p,1\≤j\≤n_i\}$ est alors un
ensemble maximal de racines fortement orthogonales  de $\∆_1$ (prop.4.4,1) de \cite{mullerJA1}).

Soit $l_0=0$ et pour $i=1,...,p,$ $l_i=n_1+...+n_i$ (donc $l_p=n$ par maximalit\é), pour $l_{i-1}+1\≤j\≤l_i$ on pose $\mu_j=\mu_{j-l_{i-1}}^{(i)}$ alors $\{\mu_1,...,\mu_n \}$ est   un
ensemble maximal ordonn\é de racines fortement orthogonales  de $\∆_1$ et il est facile de v\érifier que le sous-groupe parabolique 
$P( h_{\mu_1},...,h_{\mu_n})$  est standard pour $\∆_0^+$ car chaque $P( h_{\mu_1^{(i)}},...,h_{\mu_{n_i}^{(i)}})$ est standard pour $\∆'^{(i)+}_0$ et $P(H_1,...,H_p)$ est standard pour $\∆_0^+.$\fdem\\

\begin{rema} Soient $\goth t_0=\oplus_{1\≤i\≤n}\F h_{\lambda_j}$ et $
(\goth t_0)^0=\cap_{\lambda\in \∆_0^{\bot}} Ker \lambda$ avec $\∆_0^{\bot}=\cap_{1\≤i\≤n}\{\lambda\in \∆_0\ |\ n(\lambda,\lambda_i)=0\}$  alors $\goth t_0\subset (\goth t_0)^0\subset \goth a$ et  il est facile de v\érifier par des consid\érations cas par cas que que $\goth t_0=(\goth t_0)^0$ dans tous les cas sauf $(A_{2n-1},\alpha_n)$ pour lequel $(\goth t_0)^0=\goth a$ et  $\goth t_0$ est de dimension $n.$

\noi Le sous-groupe parabolique $P_0$ est un sous-groupe parabolique minimal uniquement dans les $3$ cas : $(A_{2n-1},\alpha_n),$   $(B_2,\alpha_1)$ et $(C_n,\alpha_n).$
\end{rema}

\bigskip

\bigskip

\noi D\émonstration de {\bf  2.1.2} dans le cas commutatif:\\

\bigskip
\noi Le point 1) d\écoule de
la proposition  pr\éc\édente   puisque les racines de $S$  sont
lin\éairement ind\épendantes.  \\
 
\noi  On suppose
que le sous-groupe parabolique 
$P_{\goth t}\ $ est standard et on reprend les notations du lemme
 pr\éc\édent.\\

\noi Terminons la d\émonstration  du point 2) lorsque la 
 d\écomposition de $x\in \goth g_1$  relativement \à
$ad(H_1)$ est: $x=x_2+x_1+x_0$ avec $F_1(x)=0.$ \\

\noi  Comme le PV  : $(\goth U(\F (2H_0-H_1))_0,\goth U(\F (2H_0-H_1))_1)$ est commutatif, on peut supposer que
$s(x_2)\subset S_1=\{\lambda_j, 1\≤j\≤l_1\}$ (prop. 5.2.2,
\cite{mullerJA1}, les \él\éments de $\goth U(\F (2H_0-H_1))_0$ utilis\és  commutent \à $\goth t$) et
lorsque
$x_2\not=0,$ soit $h=\sum_{\lambda\in s(x_2)}h_{\lambda}.$

\noi Lorsque $x_2\not=0,$ la d\écomposition de $x_1$ relativement
\à $ad(h)$ est de la forme: $x_1=y_0+y_1$ avec $y_i\in
E_i(h)\cap E_{1-i}(H_1-h)\cap \goth g_1,i=0,1,$ \à l'aide de
$exp(ad(
 E_{-1}(h)\cap
E_0(H_1-h)\cap \goth g_0))$ on peut toujours supposer que
$y_1=0$ ainsi $y_0+x_0\in \goth U(\goth s)_1$ avec $\goth
s=\oplus_{\lambda\in s(x_2)}\F h_{\lambda}.$

\noi La  proposition 2.3.1 appliqu\ée au \PV commutatif: $(\goth U(\goth
s)_0,\goth U(\goth s)_1)$  
 permet de supposer que 
$s(x)=\{\mu_1,...,\mu_k\},$ les
racines $\mu_1,...,\mu_k$ \étant fortement orthogonales et
telles que  $\{\mu\in s(x)\ |\ \mu(H_1)=2\}\subset 
S_1,$ rappelons que $\{\mu\in s(x)\ |\ \mu(H_1)=2\}\not=
S_1.$

\noi Ce r\ésultat est encore vrai lorsque $x_2=0.$\\
  
\noi Pour toute racine $\alpha$ de $\∆_1,$ la quantit\é $\sum_{1\≤i\≤k}n(\alpha,\mu_i)\in \{0,1,2\},$ puisqu'elle repr\ésente une valeur propre de $adh /\g_1,$ $h$ \étant l'\elt $1$-simple $h=\sum_{1\≤i\≤k}h_{\mu_i}.$ On rappelle \également que $$B(h_\alpha,H_1)=\frac{1}{2}\sum_{1\≤i\≤l_1}n(\lambda_j,\alpha)B(h_{\lambda_j},h_{\lambda_j}) = \frac{1}{2}B(h_{\lambda_1},h_{\lambda_1})(\sum_{1\≤i\≤l_1}n(\lambda_j,\alpha)).$$
Soient $\lambda\in S_1-s(x)$ et
$a=\sum_{1\≤i\≤k}n(\lambda,\mu_i),$ 
examinons les diff\érents cas:

\begin{enumerate}

\item  Lorsque $a=0,$ l'\él\ément $u=h_{\lambda}\in
\cap_{1\≤i\≤k}Ker\mu_i$ et v\érifie $B(u,H_1)\not=0.$
 
\item  Lorsque $a=1,$  $ n(\lambda,\mu_1)=1$ et $
n(\lambda,\mu_i)=0 $  pour $i\≥2$  (\à l'indexation pr\ès), donc $\mu_1(H_1)=1$ et comme la
racine $\lambda$ est longue on a
$n(\mu_1,\lambda)= 1$  
  d'o\ù $u=2h_{\lambda}- h_{ \mu_1}\in  
\cap_{1\≤i\≤k}Ker\mu_i$ et v\érifie $B(u,H_1)\not=0.$

\item Lorsque $a=2,$ on a soit

 $\bullet \ 
n(\lambda,\mu_1)=n(\lambda,\mu_2)=1$  et
$n(\lambda,\mu_i)=0 $  pour $i\≥3$  (\à l'indexation pr\ès) 
  donc, comme pr\éc\édemment, $u=2h_{\lambda}- h_{ \mu_1}-h_{
\mu_2}\in  
\cap_{1\≤i\≤k}Ker\mu_i$ et v\érifie $B(u,H_1)\not=0,$

\noi soit

$\bullet \ 
n(\lambda,\mu_1)=2$ et $n(\lambda,\mu_i)=0 $  pour $i\≥2$  (\à
l'indexation pr\ès)  d'o\ù $\mu_1(H_1)=1.$ Or
$\mu_1(2H_0)=\sum_{1\≤j\≤n}n(\mu_1,\lambda_j)=2$ donc il
existe un unique $j  >l_1$ tel que $
n(\mu_1,\lambda_j)=1$ d'o\ù $n( \lambda_j,\mu_1)=2$ car les
racines $\lambda_j$ et $\lambda$ ont la m\ême longueur par
cons\équent  $n(\lambda_j,\mu_i)=0 $  pour $i\≥2$ d'o\ù $u=h_{
 \lambda }-h_{\lambda_j }   \in  
\cap_{1\≤i\≤k}Ker\mu_i$ et   $B(u,H_1)=B(h_{
 \lambda },H_1)\not=0.$
 
 \noi \fdem
 \end{enumerate}

  \subsection {\bf Les cas classiques}
 
 \bigskip
 
 \noi Ce sont les cas $(R_n,\alpha_k)$ avec $R=B,BC,C$ ou $D.$ \\

\noi     Les descriptions de $\∆_0,\∆_1$ et
$\∆_2$ sont donn\ées dans la d\émonstration de la proposition
1.2.4, 3) b) (dont on reprend les notations) ceci lorsque $k\≤n-1$ dans le cas $C_n$ et
 $k\≤n-2$ dans le cas $D_n),$ en particulier:
$$\∆_1=\{\epsilon_i\pm \epsilon_j\ |\ 1\≤i\≤k<j\≤n\ ,\ \epsilon_i,1\≤i\≤k\}\cap \∆\ ,\ \∆_2=\{\epsilon_i + \epsilon_j\ |\ 1\≤i\≤j\≤k\}\cap \∆.$$

\noi Notons que dans le cas $\∆=D_n,$ la description ci-dessus pour $k=n-1$ correspond \à la graduation donn\ée  par $\∑_1=\{\alpha_{n-1},\alpha_n\}$.

\noi Les \PVs  commutatifs correspondent aux cas $(B_n,\alpha_1),(C_n,\alpha_n)$ et $(D_n,\alpha_k)$ avec $k\in \{ 1,n-1,n\}.$\\ \\

\begin{prop} Dans le cas alg\ébriquement clos, soit $(\goth g_0,\goth g_1)$ un \PV \ de type $(R_n,\∑_1 )$ avec:
\begin{enumerate}
\item  $\∑_1=\{\alpha_k\}$ lorsque $R=B_n,$
 \item $\∑_1=\{\alpha_k\}$ avec $k\≤n-1$ lorsque $R=C_n,$ 
\item  $\∑_1=\{\alpha_k\}$ avec  $k\≤n-2$ ou bien   $\∑_1=\{\alpha_{n-1},\alpha_n\}$ lorsque  $R=D_n.$ \end{enumerate}
\noi Pour tout $x\in
\goth g_1$ (resp. $\goth g_{-1}$), il existe un ensemble S de
racines lin\éairement ind\épendantes de $\∆_1$   
 (resp. $\∆_{-1}$)  telles que $P_e(h_{\epsilon_k},...,h_{\epsilon_1}).x\cap\ (\oplus_{\mu\in S 
 }\goth g^{\mu})\ \not= \emptyset$ et S contient au
plus une racine courte qui est orthogonale aux autres
racines de S.
\end{prop}

\dem Il suffit de la faire lorsque $\∑_1=\{\alpha_k\}$ avec $k\≤n-1$ car le cas $(D_n,\{\alpha_{n-1},\alpha_n\})$ provient du lemme 2.2.1 ($S$ a   au plus $3$ racines) et le cas
$(B_n, \alpha_n)$ \étant \évident, il est omis ($S$
 a   au plus une seule racine). 

\noi On proc\ède par r\écurrence sur le
rang de
$\∆,$ les cas $(B_n,\alpha_1) $ et $(D_n,\alpha_1)$
 d\écoulent de la prop.2.3.1,   les cas $(C_2,\alpha_1)$ et $(D_4,\alpha_2)$ de ce qui suit.

\noi On suppose donc  la proposition v\érifi\ée  pour les syst\èmes de racines
classiques de rang
$\≤n-1,$  et on notera $\∆^{(i)}=\{\alpha\in \∆\ |\
\alpha(h_{\epsilon_i})=0\}$ et plus g\én\éralement 
$\∆^{(i,j,...)}=\{\alpha\in \∆\ |\
\alpha(h_{\epsilon_i})=\alpha(h_{\epsilon_j})=...=0\},$ ce
sont des syst\èmes de racines gradu\és,
$\∆^{(i,..)}=\cup_{-2\≤j\≤2} \∆^{(i,..)}\cap \∆_j,$ et de m\ême type.\\
  
\noi Soit  
$u\in
\goth g_1.$
 Lorsque $s(u)\subset \∆^{(k)}_1$ (donc
$k\≥2$), $u$ est dans un PV donn\é par  un syst\ème de racines
gradu\é  de type
$(R_{n-1}, 
\∑_1= \{\alpha_{k-1}\}\ ),$ sauf dans le cas $(D_4,\alpha_2)$  o\ù $\∑_1= \{\alpha_2,\alpha_3\},$  le r\ésultat d\écoule de   l'hypoth\èse de r\écurrence et pour $(D_4,\alpha_2)$ du lemme 2.2.1. 
  
 \noi Lorsque  $s(u)$ n'est pas inclus dans $ \∆^{(k)}_1\ ,$  soit
$A=\{\epsilon_k±\epsilon_j,j\≥k+1\}\cap
 s(u).$
 
 \noi 1) Si $A$ est vide alors $\epsilon_k\in s(u)$ et $\∆=B_n$
 avec $k<n,$ il existe $v\in exp(ad(\oplus_{1\≤i\≤k-1}\goth g^{\epsilon_i-\epsilon_k}))u$ tel que $v=v_1+v_2$  avec $s(v_1)=\{ \epsilon_k\}$ et $s(v_2)\subset  \∆^{(k)}_1\cap D_n$ et on conclut en appliquant l'hypoth\èse de
r\écurrence \à $v_2$ dans $\∆^{(k)}_1\cap D_n.$

\noi  2) Lorsque $A$ est non vide et en utilisant le groupe de Weyl associ\é \à $\∆\cap (\oplus_{k+1}^n\Z \epsilon_j),$ on peut supposer que $\alpha_k\in A$ puis que $s(u) \subset \{\epsilon_k±\epsilon_{k+1},\epsilon_l+ \epsilon_{k+1}\}\cup \∆^{(k,k+1)}_1$ en appliquant le lemme 2.2.1, ainsi $u=u_1+u_2$ avec $s(u_1)\subset \{\epsilon_k±\epsilon_{k+1},\epsilon_l+ \epsilon_{k+1}\}$ et $s(u_2)\subset \∆^{(k,k+1)}_1,$   $\∆^{(k,k+1)}$ est un syst\ème de racines gradu\é
de type
$(R_{n-2},\alpha_{k-1})$ \à l'exception de l'unique
cas $(D_n,\alpha_{n-2})$ pour lequel il est de type
$(D_{n-2},\{\alpha_{n-2},\alpha_{n-3}\}).$  

\noi La d\émonstration est termin\ée lorsque $k=1$ ou $k=n-1$ car $\∆^{(k,k+1)}=\emptyset.$   \\

\noi Examinons les diff\érentes situations lorsque $2\≤k\≤n-2:$\\

a) $s(u_1)\subset \{\epsilon_k±\epsilon_{k+1}\},$ on conclut en appliquant  l'hypoth\èse de r\écurrence \à
 $u_2.$  \\
 
 b)  $s(u_1) = \{\epsilon_k-\epsilon_{k+1},\epsilon_l+\epsilon_{k+1}\},$ 
alors $u_1=x+y $ avec $s(x)=\epsilon_k-\epsilon_{k+1},$ $s(y)=\epsilon_l+\epsilon_{k+1},$  on applique le lemme 2.2.3 avec $(x,$ $h=h_{\epsilon_k-\epsilon_{k+1}},x^{-1}),$    $(y,$ $h'=h_{\epsilon_l+\epsilon_{k+1}},y^{-1})$   et $z=u_2$ ainsi toutes les hypoth\èses du lemme 2.2.3 sont v\érifi\ées et on peut se ramener au cas o\ù $s(z)\subset \∆^{(l,k,k+1)}_1,$ puis appliquer l'hypoth\èse de r\écurrence \à $z$ lorsque $\∆^{(l,k,k+1)}_1\not=\emptyset$ ($\∆^{(l,k,k+1)}_1 =\emptyset$ dans le cas $(D_4,\alpha_2)$).\\

c) $s(u_1) = \{\epsilon_k±\epsilon_{k+1},\epsilon_l+\epsilon_{k+1}\},$\\

\noi i) $R=C$

\noi Il existe $v\in exp(ad(\goth g^{2\epsilon_{k+1}}))u$ tel que $s(v)\subset \{\epsilon_k-\epsilon_{k+1},\epsilon_l+\epsilon_{k+1}\}\cup \∆^{(k,k+1)}_1$ et on est ramen\é \à b).\\

\noi ii) $R=B$

\noi Il existe $v\in exp(ad(\goth g^{\epsilon_{k+1}}))u$   tel que $ \epsilon_k \in s( v)\subset
\∆^{(k,k+1)}_1\cup\{  \epsilon_k -\epsilon_{k+1}, 
 \epsilon_k, \epsilon_i+ \epsilon_{k+1},\ i=1,...,k-1\} ,$    puis
$w\in exp(ad(\goth g^{-\epsilon_{k+1}}))v$ tel que $ \epsilon_k \in s(w)\subset \{
\epsilon_k ,\epsilon_i±\epsilon_{k+1},\ i=1,...,k-1\}\cap\∆^{(k,k+1)}_1$ et finalement  $z\in
\prod_{1\≤i\≤k-1}exp(ad(\goth g^{ \epsilon_i-\epsilon_k}))w$
tel que $s(z)\subset \{\epsilon_k\}\cup
\∆^{(k )}_1\cap D_n$ et on applique l'hypoth\èse de
r\écurrence \à $\∆^{(k )}\cap D_n.$\\

\noi ii) $R=D$

\noi Notons $u=X+y,$ avec
$s(X)=\{ \epsilon_k-\epsilon_{k+1}, 
 \epsilon_k+\epsilon_{k+1}\}$ et $y\in
E_0(h_{\epsilon_k})\cap \goth g_1.$

\noi Comme $(X,h_{\epsilon_k})$ se compl\ète en un $sl_2$-triplet
$1$-adapt\é:
$(X,h_{\epsilon_k},$ $X^{-1}=X_{-\epsilon_k+\epsilon_{k+1}}+X_{-\epsilon_k-\epsilon_{k+1}}),$ on peut supposer que  $y\in
\goth U(\goth s)_1,$  
 $\goth s$  \étant l'alg\èbre engendr\ée par $\{X,$ $h_{\epsilon_k},$ $X^{-1}\}$ (lemme 2.2.2).  

\noi  Or 
$$ \goth U(\goth s)_1=\oplus_{\mu\in
\∆^{(k,k+1)}_1}\goth g^{\mu}\oplus_{1\≤i\≤k-1}\F X_{\Lambda_i}\ 
\hbox{avec}$$
$$   X_{  \Lambda_i}=[ X_{-\epsilon_k+\epsilon_{k+1}}-X_{-\epsilon_k-\epsilon_{k+1}},X_{\epsilon_i+\epsilon_k}]\in
\goth g^{\epsilon_i-\epsilon_{k+1}}\oplus \goth
g^{\epsilon_i+\epsilon_{k+1}}, $$
\noi il est facile de v\érifier que $\goth h \cap 
 \goth U(\goth s)=\oplus_{1\≤i\≤n,i\not=k,k+1}\F h_{\epsilon_i}$ est une sous-alg\èbre de Cartan de $\goth
U(\goth s)$ et que $(\goth U(\goth s)_0,\goth U(\goth s)_1)$
est de type $(B_{n-2},\alpha_{k-1})$, les racines longues
appartiennent \à $\∆^{(k,k+1)}$ et on peut prendre   comme 
racines simples:  
$\alpha_1,...,\alpha_{k-2},
\beta=\alpha_{k-1}+\alpha_k+\alpha_{k+1},\alpha_{k+2},...,
\alpha_{n-1}$ et $  \Lambda_n,$ qui est 
 la restriction de $\epsilon_n$ \à $\goth h \cap 
 \goth U(\goth s).$

\noi Appliquons l'hypoth\èse de r\écurrence \à $z:$ il existe $g$ 
 tel que $gz=\sum_{\mu\in S} X_{\mu},$ $S$ ensemble de
racines lin\éairement ind\épendantes. Si toutes les racines
de $S$ sont longues, la d\émonstration est termin\ée, sinon il
existe une seule racine courte $  \Lambda_i,$ les autres \étant
longues, lin\éairement ind\épendantes et orthogonales \à
$  \Lambda_i$ d'o\ù 
$$ g(X+z)=X_1+\sum_{\mu\in S-  \Lambda_i} X_{\mu}\
\hbox{avec} \quad  
X_1=X_{\epsilon_k-\epsilon_{k+1}}+X_{\epsilon_k+\epsilon_{k+1}}
+ X_{  \Lambda_i }. $$
Il reste \à r\éduire l'\elt  $X_1$ \à l'aide de
$exp(ad(\goth g^{\epsilon_i-\epsilon_k})).$

\noi Ceci termine la d\émonstration du cas $(D_4,\alpha_2)$ et de la proposition.\fdem\\

\noi Avant de d\éterminer les \sgs paraboliques  standards tr\ès sp\éciaux, notons:\\
  
  \begin{lem} Soit \gog un \PV    \air de type classique $(R_n,\alpha_k)$ ($1\≤k\≤n-1$ lorsque R=C et $1\≤k\≤n-2$ lorsque R=D) alors 
 $h_{\epsilon_1}$ est 1-simple \à l'exception du cas
classique  d\éploy\é  R=C.
\end{lem}

\dem
 Comme
$h_{\epsilon_1}=h_{\epsilon_1-\epsilon_n}+
h_{\epsilon_1+\epsilon_n},$ $h_{\epsilon_1}$ est $1$-simple
dans tous les cas sauf peut-\être   $R=C.$ Dans ce
dernier cas on a $h_{\epsilon_1}=2h_{\tilde \alpha} $ et il
est ais\é de v\érifier que
$ 2h_{\tilde \alpha} $ n'est pas 1-simple dans le cas  
d\éploy\é (par exemple le \PV $(\goth U(\F h_{\epsilon_1})_0,\goth U(\F h_{\epsilon_1})_1)$ est de type $(C_{n-1},\alpha_{k-1})$ qui n'est pas r\égulier car $k-1$ est impair).

\noi   Lorsque $\goth g$ n'est pas d\éploy\ée, on note $\overline
h_{\alpha}$ la co-racine associ\ée \à la racine $\alpha\in
\overline \∆.$

\noi Par classification des diagrammes (\cite{veisfeiler},\cite{warner}), les
types possibles sont :
\vskip  5pt
 
i) $(\overline \∆,\alpha_0)=(C_{2n},\alpha_{2k})$ de
diagramme de Satake: 
\vskip 5pt

\hskip 10pt \hbox to 4,2 cm {\offinterlineskip \lower 2pt\hbox{$\bullet$} \hglue -3,2pt
{\vrule height  0,4pt depth 0pt width 0,5 cm}\lower 2pt\hbox{$\circ$}
\hglue -3,2pt
{\vrule height  0,4pt depth 0pt width 0,5 cm}\lower 2pt\hbox{$\bullet$}
  \dotfill \hbox to 1,3 cm {\offinterlineskip\lower 2pt\hbox{$\circ$}\kern -1pt   \hrulefill\kern -3,4pt\lower 2pt 
\hbox{  \offinterlineskip  {$\bullet$} \hglue -7pt\vbox{ {\hrule height 0,3pt width 0,6cm}\vskip 3,5pt{\hrule height 0,3pt width 0,6cm}
\vskip 0,3pt} \hglue -16pt $<$   \hglue -1,8pt{$\circ$}}}}

\bigskip 
ii) $(\overline \∆,\alpha_0)=(D_{2n},\alpha_{2k})$ :
 \bigskip

\hskip 10pt \hbox to 5 cm {\offinterlineskip \lower 2pt\hbox{$\bullet$}
\hglue -3pt{\vrule height .4pt depth 0pt width 0,5cm}\lower 2pt\hbox{$\circ$}
\hglue -7pt{ \vrule height .4pt depth 0pt width 0,5cm}\lower 2pt\hbox{$\bullet$}
 \dotfill \hbox to 2 cm {\lower 2pt\hbox{$\circ$}\hglue -1,3pt
\hrulefill\lower 2pt
\hbox{$\bullet$}\hrulefill \lower 2pt \vtop {\offinterlineskip \hbox{$\circ$} \hbox to 5pt{\hfill \vrule height 12pt width 0,3pt\hfill} \hbox{$\circ$}}\hrulefill\lower 2pt\hbox{$\bullet$} } }

 \bigskip
iii) $\overline \∆=E_7$ et $\overline \∑_1=\{\alpha_1\}$ ou
$\overline \∑_1=\{\alpha_6\},$  
$\F=\R$  et le diagramme de Satake est:
 \bigskip
 
\hskip 10pt \hbox to 3,5 cm {\lower 2pt\hbox{$\circ$}\hrulefill\lower 2pt
\hbox{$\bullet$}\hrulefill \lower 2pt \vtop {\offinterlineskip \hbox{$\bullet$} \hbox to 5pt{\hfill \vrule height 12pt width 0,3pt\hfill} \hbox{$\bullet$}}\hrulefill\lower 2pt\hbox{$\bullet$}\hrulefill\lower 2pt\hbox{$\circ$}\hrulefill\lower 2pt\hbox{$\circ$} }

\bigskip
 \noi Dans tous les cas, il existe une racine $\gamma \in
\overline \∆_2$  telle que
$2h_{\tilde
\alpha }=2\overline h_{\gamma},$ $\gamma$ est la somme 
 de deux racines,
$\alpha$ et
$\beta,$ de
$\overline {\Delta}_1$ dont la diff\érence n'est pas une
racine, ainsi
  $(X_{\alpha}+X_{\beta},2h_{\tilde
\alpha},2(X_{-\alpha}+X_{-\beta}))$ est   un
$sl_2-$triplet
$1$ adapt\é  (tds principale de type $A_2$) au sens du \PV: $(\overline{\goth
g}_0,\overline{\goth g}_1)$ d'o\ù $2h_{\tilde
\alpha}$ 
est \également $1$-simple  au sens du
PV: $( \goth g_0, \goth g_1).$ 

\noi   $\gamma$ est la plus grande racine de $\overline \∆$ pour ii) et iii) 
et
$\gamma=\alpha_1+2(\sum_{i=2}^{2n-1}\alpha_i )+
\alpha_{2n}$ dans le cas i). \fdem\\

\noi Dor\énavant, dans tout ce \S, le \PV \gog v\érifie les hypoth\èses du lemme 2.4.2, ce qui permet de d\éfinir:
 
$$p_0=\left\{ \begin{array}{ll}    \frac{k}{2} \  \textrm{ dans le cas}\  \ R=C\  \textrm{avec  }\ \goth g\ \textrm{ d\éploy\ée }   ,\\  
   k\    \textrm{ sinon}\    \end{array}\right.$$
\noi et pour $i=1,...,p_0:$
$$H_i=\left\{ \begin{array}{ll}   2(h_{2\epsilon_{k-2i+1}}+h_{2\epsilon_{k-2i+2}})    \  \textrm{ dans le cas}\  \ R=C\  \textrm{avec  }\ \goth g\ \textrm{ d\éploy\ée }   ,\\  
   h_{\epsilon_{k-i+1}}\    \textrm{ sinon.}\    \end{array}\right.$$
\noi Le \sg parabolique $P_0=P(H_1,...,H_{p_0})$ est alors un \sg 
 parabolique standard tr\ès sp\écial.\\
 
 \noi  Lorsque $2k\≤n,$ soient
$H'_1=\sum_{i=1} ^kh_{\epsilon_i-\epsilon_{2k-i+1}}$  et
$H'_2=\sum_{i=1}^kh_{\epsilon_i+\epsilon_{2k-i+1}},$ alors  $P_0' =P(H'_1,H'_2)$  est  \également un \sg 
 parabolique standard tr\ès sp\écial.\\
 
 \noi Soit $\goth t_0=\oplus_{i=1}^{p_0}\F H_i$ (resp $\goth t'_0=  \F H'_1\oplus \F H'_2$) alors $\∑_{\goth t_0}=\{\alpha_i, i\≥k+1\} $ dans tous les cas sauf le cas $C_n$ d\éploy\é  pour lequel $\∑_{\goth t_0}=\{\alpha_i, i\≥k+1,\alpha_{2j-1},1\≤j\≤\frac{k}{2}\} $ ($\∑-\∑_{\goth t'_0}=\{\alpha_k,\alpha_{2k}\}$) et $(\goth t_0)^0=\goth t_0$ (resp $(\goth t'_0)^0=\goth t'_0$).
 
 \noi $P_0$ est un \sg parabolique minimal dans les cas $(R_n,\alpha_n)$ avec $R=B$ ou $BC.$\\

\begin{prop} Soit \gog un \PV    \air de type classique $(R_n,\alpha_k),$ $R=B,C,BC$ ou D ($1\≤k\≤n-1$ lorsque R=C et $1\≤k\≤n-2$ lorsque R=D)et soit $P=P(H_1,...,H_p)$ un \sg  parabolique standard tr\ès sp\écial alors 

$\bullet$  soit  $P=P'_0$ et $2k\≤n,$

$\bullet$  soit $P\supset P_0,$   $k\≥2$  et  il existe des
entiers:
$ l_{p+1}=0<l_p< ...<l_2<l_1=k$ tels que  $H_i=\sum_{l_{i+1}+1\≤j\≤l_i}h_{\epsilon_j}$ pour $i=1,...,p.$
\end{prop}

\dem Bien que $\g$ ne soit pas d\éploy\ée en g\én\éral et de la m\ême mani\ère que dans le 3) de la proposition 1.2.4, la
consid\ération du diagramme \à poids associ\é \à $H_p$ donne
$2$ cas possibles:
\vskip 2mm
Premier cas: il existe   $j>k$  
tel que $\alpha_j(H_p)\not=0$ et $\alpha_l(H_p) =0$ pour
$1\≤l\≤j-1.$

\noi Comme $(\goth U(\F H_p)_0,\goth U(\F H_p)_1)$ est un
\PV de type $(A_{j-1},\alpha_k),$ il est commutatif et
gradu\é par $2H_0-H_p$ qui est $1$-simple dans ce \PV donc
$j=2k$ d'o\ù $H_p=H'_2.$ Lorsque $p\≥3,$ dans ce \PV le parabolique associ\é \à $H_1,...,H_{p-1}$ est standard donc par le lemme 2.3.2,
il existe $i\in
\{1,...,k-1\}$ tel que
$H_1=\sum_{1\≤j\≤i}h_{\epsilon_{k-j+1}-\epsilon_{k+j}}$ ce qui est absurde ($\epsilon_1+\epsilon_{k+1}(H_1)=-1$) d'o\ù $p=2$ et $H_1=H'_1.$ 

Deuxi\ème cas:
il existe une unique valeur $j$ , $1\≤j<k$
telle que $\alpha_j(H_p)=2$ donc $k\≥2$ et
$H_p=\sum_{1\≤i\≤j}h_{\epsilon_i}$ ce qui termine la d\émonstration dans le cas $p=2.$  Lorsque $p\≥3,$  $(\goth U(\F H_p)_0,\goth U(\F H_p)_1)$ est un
\PV de m\ême type : $(R_{n-j},\alpha_{k-j})$  et dans ce \PV le parabolique associ\é \à $H_1,...,H_{p-1}$ est standard mais il n'est pas de type $P'_0$ qui est associ\é \à $H'_1=h_{\epsilon_{j+1}-
\epsilon_{2k-j}}+...+h_{\epsilon_k-\epsilon_{k+1}}$ (resp.$H'_2=h_{\epsilon_{j+1}+
\epsilon_{2k-j}}+...+h_{\epsilon_k+\epsilon_{k+1}}$) car $H'_1$ et $H'_2$ ne sont pas $1$-simples sp\éciaux donc on est \à nouveau dans le 2\ème cas et on obtient le r\ésultat par r\écurrence sur  $n.$\fdem\\

\noi Fin de la d\émonstration du {\bf th\éor\ème 2.1.1.}\\

\noi 1) Le point 1 d\écoule du $\S 2.5.$ lorsque $\overline
\∆$ est de type exceptionnel et lorsque $\overline \∆$ est de
type classique cela r\ésulte de la proposition   2.4.1 pour
$P_0$ et des  r\ésultats de \cite{hillerohrle} pour $P'_0$
puisque le parabolique associ\é aux racines
$\{\alpha_i,i\not= k,2k\ ,\ 1\≤i\≤n\}$ est de longueur $\≤4$
donc op\ère avec un nombre fini d'orbites dans son radical
unipotent.\\

\noi 2) On omet la v\érification pour $P_0'$ et on termine la
d\émonstration de 2.1.2 dans les cas classiques pour 
  les \sgs paraboliques $P_{\goth t}\supset P_0,$ dont la description figure dans la proposition  2.4.3.
  
\noi   Soit $x\in \goth g_1,$  sa d\écomposition   relativement \à
$ad(H_1)$  est de la forme 
$x=x_2+ x_0,$ on suppose donc que   $x_2$ est non nul et que 
$F_1(x)=0.$

\noi  A l'action de
$P_{\goth t}$ pr\ès, et \à l'aide de la description de $\∆_1$ et de $\∆_2,$ on peut supposer que $(x_2,H,x_2^{-1})$
est un
$sl_2-$triplet avec comme formes possibles de $H$ (prop.2.5
de \cite{mullerJA1}):

\quad $(1) :$
$H=\sum_{l_2+1\≤i\≤l}h_{\epsilon_i}+
\sum_{1\≤j\≤q}h_{\gamma_j} $  avec $\gamma_j=\epsilon_{l+j}+ 
\epsilon_{k+j}$ ou bien $\gamma_j=\epsilon_{l+j}- 
\epsilon_{k+j}$ ce que l'on note $\gamma_j=\epsilon_{l+j}+(±
\epsilon_{k+j}),$   ou bien

\quad $(2) :$
$H= \sum_{1\≤j\≤q}h_{\gamma_j}$ avec $l=l_2$  et on a: $l+q\≤k,$   
ou bien

\quad $(3) : $
$H=\sum_{l_2+1\≤i\≤l}h_{\epsilon_i}$ avec  $ l< k.$

\noi On peut toujours supposer que
$x_0=y_{-1}+y_0$ avec $y_i\in E_i(H)\cap E_0(H_1)\cap \goth
g_1$ ( lemme 2.2.2).
\vskip 2mm
\noi i) Lorsque  $y_{-1} =0,$ $u=H_1-H$ commute \à $x$ et v\érifie
$B(u,H_1)\not=0.$ 
\vskip 2mm
\noi ii)  Lorsque  $y_{-1} \not= 0,$ donc $H$ est de la forme $(1)$ ou $(2),$ posons
$h=\sum_{1\≤j\≤q}h_{\gamma_j}, $ $h$ est
$1$-simple car les
$q$ racines $ \gamma_j$ sont fortement
orhogonales.  

\noi Or: $E_2(H)\cap \goth g_1=E_2(H-h)\cap \goth g_1\oplus
E_2(h)\cap \goth g_1$ avec:
$$E_2(H-h)\cap \goth g_1=(\oplus_{l_2+1\≤i\≤l,j\≥k+q+1}\goth g^{\epsilon_i ±\epsilon_j}\oplus \goth g^{\epsilon_i})\cap \goth g_1\ ,\ E_2(h)\cap \goth g_1=\oplus_{1\≤i,j\≤q}\goth g^{\epsilon_{l+i}+ 
(±\epsilon_{k+j})}$$
et de m\ême pour  $E_{-2}(H-h)\cap \goth g_{-1}$,  ainsi $E_2(H-h)\cap \goth g_1\oplus E_{-2}(H-h)\cap \goth g_{-1}$ commute  avec $E_2(h)\cap \goth g_1\oplus E_{-2}(h)\cap \goth g_{-1}.$ 

\noi Or le \PV: 
$(\goth U((2H_0-h)\F)_0,\goth U((2H_0-h)\F)_1=E_2(h)\cap \goth g_1),$ 
est commutatif, gradu\é par $\frac{h}{ 2}$
et r\égulier d'o\ù, \à l'action de $Aut_e(\goth U((2H_0-h)\F)_0)$ pr\ès, qui centralise $\goth t+\F h+\F H$ (prop.5.2.2,\cite{mullerJA1}), on peut supposer que 
$x_2=X_1+X_2$ avec
$ X_1\in E_2(H-h)\cap \goth g_1 $ et
$X_2=\sum_{1\≤j\≤q}X_{\gamma_j}.$ 

\noi Soit $\goth s=\oplus_{l_2+1\≤i\≤l}\F h_{\epsilon_i}$  dans le cas $(1)$ (resp.$\goth s=\{0\}$ 
 dans le cas $(2)$) on v\érifie que:
$$E_{-1}(H)\cap \goth g_1\cup   E_0(H)\cap \goth g_1 \subset
\goth U(\goth s)_1$$
donc   $X_2+y_{-1}+y_0\in  \goth U(\goth s)_1.$  

\noi Dans le cas $(1),$ le  \PV:
$(\goth U(\goth s)_0,\goth U(\goth
s)_1)$ est de type classique avec pour
syst\ème de racines gradu\é: $\∆'_i=\∆_i^{(l_2+1,...,l)}$ et le \sg parabolique
$P(H'_1,...,H_p)$ avec
$H'_1=\sum_{l+1\≤i\≤k}h_{\epsilon_i}$ est un \sg parabolique standard. \\

\noi On se place dans  le  \PV:
$( P(H'_1,...,H_p),\goth U(\goth
s)_1)$ (avec $H'_1=H_1$ dans le cas $(2)$) et soit $\goth t'=\F H'_1\oplus_{i=2}^p\F H_i.$

\noi En appliquant le lemme 2.2.4 \à $X_2+y_{-1}+y_0$   on peut
supposer que $s( y_{-1})=\{\mu_1,...,\mu_r\}$ avec $1\≤r\≤q$ et
 $\mu_i=
\epsilon_{p_i}-(\pm \epsilon_{k+i}),$ les $p_i$ \étant  tous
distincts et inf\érieurs ou \égaux \à $l_2,$ mais alors
$$[[E_0(H'_1)\cap E_0(h)\cap \goth U(\goth
s)_1,\oplus_{1\≤i\≤r}\goth
g^{-\mu_i}],\goth t']]=0$$ donc  par le
lemme 2.2.3 on peut supposer que 
  $s(y_0)\subset \{ \epsilon_i\pm
\epsilon_j,$ $\epsilon_i,$ $1\≤i\≤l_2\ ,\ i\not=p_1,...p_r,$
 $l+q+1\≤j\≤n\}\cap \∆_1.$\\

\noi L'\él\ément $u=h_{\epsilon_{p_1}}+(\pm
h_{\epsilon_{k+1}})-h_{\epsilon_{l+1}}$  commute \à  
 $x=X_1+\sum_{j=1}^qX_{\gamma_j}+\sum_{j=1}^rX_{\mu_j}+y_0$ et v\érifie $B(u,H_1)\not=0.$\fdem  

\bigskip

 \subsection {\bf  Les cas exceptionnels non commutatifs}
  
    \bigskip
    
    \subsubsection {\bf G\én\éralit\és}
    
     \bigskip
     
\noi    La liste des diagrammes de Dynkin gradu\és correspondants a \ét\é \établie dans la d\émonstration du lemme 1.3.1 et est donn\ée par:

\begin{enumerate}

\item    
$(F_4,\alpha_1),(E_6,\alpha_2),(E_7,\alpha_1),(E_8,\alpha_8)$   pour lesquels  $\goth g_2$ est de dimension $ 1,$  

\item  $(F_4,\alpha_4)$ et $\goth g$ est non d\éploy\ée (ce
qui correspond \à des $\F$-formes de 3), 4) et 5)),

\item $(E_7,\alpha_6),$

\item  $(E_7,\alpha_2),$ 

\item  $(E_8,\alpha_1).$ 
\end{enumerate}

\noi Rappelons que $\goth g$ est  d\éploy\ée dans les cas
$E_i,i=6,7,8$ et notons que $H_0$ est proportionnel \à $\sum_{\omega \in \∆_2}h_{\omega}$ (puisque celui-ci est non nul et est dans le centre de $\goth g_0$).
 \bigskip
 \begin{lem} Il y a un unique sous-groupe parabolique standard tr\ès sp\écial,   $P(H_1,H_2).$   Lorsque dim($\goth g_2)\≥2,$ $H_2=2h_{\tilde \alpha}$ et
  $\goth t$ est une sous-alg\èbre de l'alg\èbre de Lie engendr\ée par $\goth g_{\pm 2}.$
  
 \noi $P(H_1,H_2)$ est un \sg parabolique maximal de G \à l'exception du cas $(E_6,\alpha_2)$ pour lequel $\∑_0-\∑_{\goth t}=\{\alpha_1,\alpha_6\}.$
  \end{lem}
  
  \vskip 5pt
  \begin{rema}  Pr\écis\ément:
  
$$\∑-\∑_{\goth t}= \left \{
    \begin{array}{l} \{\alpha_1,\alpha_4\}\ \hbox{lorsque }\ \∆=F_4 ,\\ 
  \{\alpha_1,\alpha_6\}\ \hbox{lorsque }\ \∆=E_7\ \hbox{et}\ \∑_1\subset   \{\alpha_1,\alpha_6\},\\ 
   \{\alpha_1,\alpha_8\}\ \hbox{lorsque }\ \∆=E_8 ,\\ 
   \{\alpha_1,\alpha_2\}\ \hbox {dans le cas  }\ (E_7,  \alpha_2).\\  \end{array} \right.    
$$
  
  \end{rema}

  \vskip5pt
 
 \dem  L'unicit\é r\ésulte de la d\émonstration du lemme 1.3.1 puisqu'il y a toujours un unique \sg parabolique tr\ès sp\écial (\à conjugaison pr\ès)
 sauf pour $(E_7,\alpha_2)$ mais  alors l'un des $2$ \elts tr\ès sp\éciaux est donn\é  par $ 2h_{\tilde \alpha}.$
 
 \noi Comme $ 2h_{\tilde \alpha}$ est toujours $1$-simple sp\écial pour $E_7$ et $E_8$ lorsque  dim($\goth g_2)\≥2,$ il en est de m\ême pour leur  formes r\éelles car $\tilde \alpha$ est une racine r\éelle donc $H_2= 2h_{\tilde \alpha}$  lorsque dim($\goth g_2)\≥2.$

 \noi Comme  $P(H_1,H_2)$ est un \sg parabolique standard, $\∑_{\goth t}=\{\alpha\in \∑_0|\alpha (H_2)=0\}$ et
 lorsque $H_2=2h_{\tilde \alpha},$  on  d\étermine imm\édiatement  $\∑_{\goth t}$ \à partir  du diagramme de Dynkin compl\ét\é (cf.planches V \à VIII de [BO 1]).
 
 \noi  Lorsque dim$(\goth g_2)=1,$ soit $ \∑_1=\{\lambda_1\},$ par orthogonalisations successives avec conservation de l'ordre de $\∆,$ on construit $\{\lambda_2,\lambda_3,\lambda_4\},$ ainsi  on a $4$ racines orthogonales, longues dans $\∆_1$ et de somme $2H_0$ (cf.le cas commutatif ainsi que 6.4 \à 6.6 de [\cite{mullerJA2});  pour $\alpha\in \∑_0$ on a soit $n(\alpha,\lambda_i)=0$ pour $i=1,...,4$ soit $n(\alpha,\lambda_{i_0})<0, $ avec 
 $i_0=$inf$\{i|n(\alpha,\lambda_i)\not=0\}$ 
 donc le \sg parabolique $P(h_{\lambda_1}+h_{\lambda_2},h_{\lambda_3}+h_{\lambda_4})$ est standard;  par unicit\é $H_2=h_{\lambda_3}+h_{\lambda_4}$ d'o\ù $\∑_{\goth t}=\{\alpha\in \∑_0\ |\ 
 n(\alpha,\lambda_1)<0$ et $n(\alpha,\lambda_2>0\}\cup \{\alpha\in \∑_0\ |\
 n(\alpha,\lambda_1)=n(\alpha,\lambda_2=0\}$  et on le d\étermine  dans chaque cas \à l'aide du tableau II de \cite{mullerJA2}. Les r\ésultats sont report\és dans la remarque ci-dessus.  \fdem\\
 
 \noi Rappelons qu'il existe $\ell ,$  valeur rappel\ée dans chaque cas, racines orthogonales de $\∆_1$ et
longues dans $\∆_1$ (longues dans $\∆$ \à l'exception du
cas $(F_4,\alpha_4)$), not\ées
$\lambda_1,...,\lambda_{\ell },$ telles que
$H_1=\sum_{1\≤i\≤\ell }h_{\lambda_i}.$  
\noi Ainsi le \PV: $(E_2(H_1)\cap \goth g_0,E_2(H_1)\cap \goth
g_1)$ est dans tous les cas``presque-commutatif" (\cite{mullerJA1}
prop.6.6 et derni\ère remarque) d'o\ù tout \él\ément non nul
$x\in
  E_2(H_1)\cap \goth
g_1$ est dans la $G_{H_1}$-orbite
 d'un  \elt de la forme :
$$x=\sum_{1\≤i\≤q}X_{\lambda_i}\ ,\  1\≤q\≤\ell \ \hbox{ avec}\ \ 
  X_{\lambda_i}\in
\goth g^{\lambda_i}-0\ \ \hbox{ et 
}\ \ [ X_{\lambda_i},X_{\lambda_j}]=0 \ \hbox{pour}\ 1\≤i<j\≤q,  $$
$$\hbox{ avec }\ \ell=\begin{cases}2\text{ dans les cas }(F_4,\alpha_1),(E_6,\alpha_2),(E_7,\alpha_1),(E_8,\alpha_8),\\  
\text{ ainsi que } 
(F_4,\alpha_4) \text{ lorsque }\goth g\text{  est non d\éploy\ée},\\ 
3  \text{ dans le cas  }(E_7,\alpha_2),\\
4 \text{ dans les cas }
 (E_7,\alpha_6) \text{ et }
(E_8,\alpha_1).\end{cases}
$$

\subsubsection{\bf Les cas exceptionnels pour lesquels $\bf \goth g_2$ est de dimension $\bf 1$} 

\bigskip

\noi On reprend les notations du lemme 1.3.1 et on rappelle que:    
 $\goth g_1=E_2(H_1)\cap \goth g_1\oplus E_1(H_1)\cap \goth g_1\oplus E_0(H_1)\cap \goth g_1$  et $\∆_2=\{\tilde{\alpha}=\frac{1}{2}(\sum_{1\≤i\≤4} \lambda_i)\}$ avec $\∑_1=\{\lambda_1\}$ et que:

 \bigskip
 
 \begin{rappel}  Soient $\beta_i,1\≤i\≤4,$ $4$ racines orthogonales de
$\Delta_1,$  donc longues et de somme
$2\tilde{\alpha}.$ \\

1) Soient $x=X_{\beta_1} +X_{\tilde
\alpha-\beta_1},$ alors $F(x)\not=0$ et pour $(s,t)\in \F^2.$  

\noi Il existe $A\in 
  \goth g^{\tilde \alpha-\beta_2-\beta_4}$ tel que
$[A,X_{\tilde \alpha-\beta_1}]=tX_{\beta_3},\ B\in 
\goth g^{\tilde \alpha-\beta_2-\beta_3}$ tel que
$[B,X_{ \tilde
\alpha-\beta_1}]=sX_{\beta_4} $ et
$C\in \goth g^{\tilde
\alpha-\beta_1-\beta_2}$ tel que
$[C,X_{ \beta_1}]=-[B,tX_{\beta_3}] ,$ posons
   $g=exp(ad(C))exp(ad(B))exp(ad(A)),$ on a 
$g(x)=X_{\beta_1} +X_{\tilde
\alpha-\beta_1}+tX_{\beta_3}+sX_{\beta_4}.$ 

\noi $g\in P_e(H_1,H_2)$ lorsque $ \beta_2(H_1)=2.$\\
 
2)  Soit $x\in
\goth g_1$   tel que:
$s(x)= \{\tilde{\alpha}-\beta_1,\beta_1,\beta_2,\beta_3\}$ (resp.$s(x)= 
\{\tilde{\alpha}-\beta_1,\beta_i,1\≤i\≤4\}),$ il existe $g\in
exp(ad( \goth
g^{\tilde \alpha-\beta_2-\beta_3}))exp(ad( \goth
g^{\tilde \alpha-\beta_1-\beta_3}))exp(ad( \goth
g^{\tilde \alpha-\beta_1-\beta_2}))$ tel que: 
$$s(g(x))=\{\beta_i,i=1,...,4\} \quad (\hbox{resp.}\ s(g(x))\subset
\{\beta_i,i=1,2,3,4\}).$$ (\cite{mullerNAG}, d\émonstration du lemme
2.3.1). 
 
\noi $g\in P_e(H_1,H_2)$ lorsque $\beta_1(H_1)=2$ et $\beta_4(H_1)=0$ ou bien lorsque 
 $\beta_i(H_1)=1$ pour $i=1,...,4$.
\vskip 2mm

\end{rappel}
 
 \begin{rappel} $\forall \sigma\in S_4$ il existe $w\in W_0$ tel que $w(\lambda_i)=\lambda_{\sigma (i)}$ pour $i=1,...,4.$

 \end{rappel}
 \noi (cf.remarque du $\S 2.2$ \cite{mullerNAG} car   pour tout $ \sigma\in S_4$ on a   $\cap_{ i=1,2}E_1(h_{\sigma (\lambda_i)}) \cap_{ j=3,4} E_0(h_{\sigma (\lambda_j)}) \not=\{0\})$ \\
 
 \noi  Fin de la d\émonstration du th\éor\ème 2.1.1.\\
 
 \noi  Soit $x\in
\goth g_1.$\\

 \noi {\bf Premier cas:} $x$ a une composante non nulle suivant
$E_2(H_1)\cap \goth g_1.$\\ 

\noi A l'action de $G_{H_1}$-pr\ès on
peut supposer que $\lambda_1\in s(x)$ et que $x\in \goth
g^{\lambda_1}\oplus \goth g^{\tilde{\alpha}-\lambda_1}\oplus
E_0(h_{ \lambda_1})\cap \goth g_1$ (lemme 2.2.1); comme le  \PV $(E_0(h_{ \lambda_1})\cap \goth
g_0,E_0(h_{ \lambda_1})\cap \goth g_1)$ est {\it commutatif},
lorsque $s(x)/\{\tilde{\alpha}-\lambda_1, \lambda_1\}$ est non vide, on peut supposer que
$s(x)/\{\tilde{\alpha}-\lambda_1, \lambda_1\}\subset S,$ $S$ 
 \étant un ensemble  de racines orthogonales entre
elles et \à  $\lambda_1$ telles que $\sum_{\mu \in
  S}\mu+ \lambda_1=2\tilde \alpha$ (proposition 2.3.1) donc $\sum_{\mu \in
  S}\mu=\lambda_2+\lambda_3+\lambda_4.$ Comme:
  $$E_1(H_1)\cap \goth g_1=\oplus_{\sigma \in I } E_0(h_{\sigma (\lambda_1)}) \cap E_1(h_{\sigma (\lambda_2)}) \cap E_1(h_{\sigma (\lambda_3)}) \cap E_0(h_{\sigma (\lambda_4)})$$ 
  \noi avec $ I= \{1,(1,2),(3,4),(1,2)(3,4)\},$ l'ensemble  $S_1=\{\mu\in S\ | \
n(\mu,\lambda_2)=0\}$ est non vide, en effet si $S_1=\emptyset$ alors pour $\mu\in S$ on a $n(\mu,\lambda_2)=1$ donc $S=\{\mu_1,\mu_2\}$ d'o\ù l'une des $2$ valeurs : $n(\mu_1,\lambda_i)+n(\mu_2,\lambda_i),i=3,4,$ est diff\érente de $2$ ce qui est absurde.

\noi A l'aide de $G_{h_{\lambda_1}}\cap G_{h_{\lambda_2}}$ on peut se ramener \à 
$S_1  = \{\lambda_4 \}$ ou bien $S_1  = \{\lambda_3,\lambda_4\}.$\\

\noi $\bullet$ Lorsque $\lambda_2 \in s(x),$  donc $F_1(x)\not=0,$ on a $S=\{ \lambda_i,i=2,3,4\}$ et $s(x)\subset  \{ \tilde \alpha- \lambda_1, \lambda_i,i=1,...,4\}.$ A l'aide du rappel 2.5.3, on peut supposer que $s(x)\subset \{\lambda_1,\tilde \alpha- \lambda_1, \lambda_2\}$ ou que 
$s(x)= \{\lambda_i,i=1,...,4\}$ ou que $s(x)\subset \{\lambda_i,i=1,2,3\}.$  

\noi Lorsque   $s(x)\subset \{\lambda_i,i=1,2,3\},$  l'\elt $u=h_{\lambda_4}$ centralise $x$ et v\érifie $B(u,H_2)\not=0,$  dans les autres cas  $F(x)\not=0$ (cf.rappel 2.5.3,1)).\\

\noi  $\bullet$ Lorsque $\lambda_2 \notin s(x),$  donc $F_1(x)=0,$ il existe $\mu_2$ et $\mu_3\in S$ tels que
$n(\mu_2,\lambda_2)=n(\mu_3,\lambda_2)=1$ donc  $S_1=\{\lambda_4\}$  et $n(\mu_i,\lambda_j)= 1$ pour $i,j=2,3.$ 

\noi L'\elt    
$u=h_{\lambda_2} -h_{\lambda_3} $   centralise $x$ et v\érifie $B(u,H_1)\not=0.$ 

\noi De plus, si on note $  
\beta_1=\lambda_1,$
$\beta_2=\mu_2,$ 
$\beta_3=\mu_3,$ 
$\beta_4=\lambda_4,$ alors on peut supposer que 
$s(x)\subset \{ \beta_1,\tilde{\alpha}-\beta_1,\beta_j,\beta_k\}$  
 ou bien que 
$s(x)\subset\{\beta_j ,i=1,...,4\}$ en appliquant le  rappel 2.5.3\\

\noi {\bf Deuxi\ème cas :} $x$ n'a pas de composante suivant
$E_2(H_1)\cap \goth g_1.$  

\noi Lorsque $x\in E_1(H_1)\cap \goth g_1,$  $H_1-H_2$ centralise $x$ et v\érifie $B(H_1-H_2,H_1)\not=0.$\\ 
 
a) Le cas simplement lac\é

\noi En appliquant le lemme 2.2.1 on peut  supposer que $x\in \goth g^{\mu_1}\oplus \goth g^{\tilde \alpha -\mu_1}\oplus E_0(h_{\mu_1})\cap \goth g_1,$ or le \PV $(E_0(h_{\mu_1})\cap \goth g_0,E_0(h_{\mu_1})\cap \goth g_1)$ est commutatif , on peut donc supposer comme pr\éc\édemment que
$\mu_1\in s(x)\subset
\{\mu_1,\tilde{\alpha}-\mu_1,\mu_2,\mu_3,\mu_4\},$ les
racines $\mu_i,i=1,...,4,$ \étant orthogonales de somme
$2\tilde {\alpha},$ et que $\forall \beta\in s(x)$ on
a $\beta(H_1)\≤1.$
\vskip 2mm
\noi i) $s(x)$ a cinq racines

\noi Comme $2\tilde{\alpha}(H_1)=4=\sum_{i=1}^4\mu_i(H_1),$ on a   $\mu_i(H_1)=1$
pour $i=1,...,4$ d'o\ù $x\in E_1(H_1)\cap \goth g_1,$  et on applique le 2) du rappel 2.5.3 aux
racines
$\beta_i=\mu_i,i=1,...,4.$  
\vskip 2mm
\noi ii) $s(x)$ a au plus $4$ racines donc elles sont lin\éairement ind\épendantes.

\noi Seul le cas o\ù $\exists i$ tel que
$\mu_i(H_1)=2$ est \à consid\érer mais alors $\mu_i\notin s(x)$ et  lorsque $\{\mu_1,\tilde
{\alpha}-\mu_1\}\subset   s(x),$  par le 1) du rappel 2.5.3  appliqu\é
 aux racines 
$\beta_i=\mu_i,1\≤i\≤4,$  il existe $y\in Px$ tel
que
$s(y)=\{\tilde {\alpha}-\mu_1,\mu_1\}$  d'o\ù  $y\in E_1(H_1)\cap \goth g_1.$\\

b) Le cas restant: $(F_4,\alpha_1).$ \\

\noi Rappelons que $\∆_1=\{ \tilde{\alpha}-\lambda_i\
, 1\≤i \≤4,\ \frac{1} {2}(\lambda_i+ \lambda_j)\
,\ 1\≤i,j\≤4\}.$ \\ 

 \noi  A l'aide de $G_{H_1}$ on peut
supposer que
$\{\lambda\in s(x)\ |\ \lambda(H_1)=0\}\subset \{
\lambda_3,\lambda_4\}$ d'o\ù $s(x) \subset \{
\lambda_3,\lambda_4,\frac{1}{2}(\lambda_i+\lambda_j),i\in
\{1,2\}\ ,j\in\{3,4\}\ \}.$ 

\noi Lorsque $s(x)$ contient des racines courtes, on peut supposer que
$\frac{1}{2}(\lambda_1+\lambda_3 )\in s(x)$ (cf.rappel 2.5.4) et que 
$\frac{1}{2}(\lambda_2+\lambda_3)\notin s(x)$ (on utilise
$exp(ad(\goth g^{ \frac{1}{2}( \lambda_2-\lambda_1)}\ );$ de
m\ême on peut supposer que 
$\frac{1}{2}(\lambda_2+\lambda_4 )$ ou  $\frac{1}{2}(\lambda_1+\lambda_4 )$ n'est pas dans  $s(x)$ (sinon on fait agir  $exp(ad(\goth g^{\frac{1}{2}(
\lambda_1-\lambda_2 )})\ ).$ On a ainsi $2$ cas:
\vskip 2mm
\noi iii) ) $s(x)\subset
\{\lambda_3,\lambda_4,\frac{1}{2}(\lambda_1+\lambda_3),\frac{1}{2}(\lambda_1+\lambda_4)\}$ et  $ h_{\lambda_2}$ convient.

\noi Comme $x\in E_0(h_{\lambda_2})\cap \goth g_1$ et que le \PV $(E_0(h_{\lambda_2}\cap \goth g_0,E_0(h_{\lambda_2}\cap \goth g_1)$ est commutatif, on termine en appliquant la proposition 2.3.1.
 \vskip 2mm
\noi iv)  $\frac{1}{2}( \lambda_1+\lambda_3),
\frac{1}{2}(\lambda_2+\lambda_4 )
\}\subset s(x)\subset
\{\lambda_3,\lambda_4,\frac{1}{2}(\lambda_1+\lambda_3 ),
\frac{1}{2}(\lambda_2+\lambda_4)\}$

 \noi L'action de $exp(ad(\goth g^{\frac{1}{2}(\lambda_3-\lambda_1)}))$
permet de se ramener \à $\tilde \alpha-\lambda_1\in s(x)
\subset
\{\tilde \alpha-\lambda_1,
\lambda_4,\frac{1}{2}(\lambda_1+\lambda_3 ),
\frac{1}{2}(\lambda_2+\lambda_4)\}$ puis avec  un \elt 
appropri\é de $exp(ad(\goth g^{\tilde
\alpha-\lambda_2-\lambda_3 }))$ on arrive \à
$s(x)\subset \{\tilde
\alpha -\lambda_1,
 \frac{1}{2}(\lambda_1+\lambda_3),
\frac{1}{2}(\lambda_2+\lambda_4)\}$ ($3$ racines lin\éairement ind\épendantes) et  $
h_{\frac{1}{2}(\lambda_2-\lambda_4)} $  convient. \fdem\\

\noi On a \également \établi:\\
\begin{prop}
Dans les cas $(E_6,\alpha_2),(E_7,\alpha_1),(E_8,\alpha_8),$ soit $ x\in \goth g_1$  il existe $ y\in P_e(H_1,H_2).x$ et un ensemble de 4 racines fortement orthogonales $\beta_1,...,\beta_4,$ tels que $s(y)$ comprenne au plus 4 racines appartenant \à $\{\tilde \alpha-\beta_1,\beta_i,1=1,...,4\}$ et dans le cas $(F_4,\alpha_1)$ il faut y adjoindre l'ensemble $\{\tilde \alpha-\lambda_1,\frac{1}{2}(\lambda_1+\lambda_3),\frac{1}{2}(\lambda_2+\lambda_4)\}.$

\end{prop}
\bigskip

 \subsubsection{\bf D\émonstration du th\éor\ème 2.1.1 dans les  
cas exceptionnels restants}

\bigskip

\noi On rappelle que $\goth g_1=E_2(H_1)\cap \goth g_1\oplus E_0(H_1)\cap \goth g_1.$\\

\noi Il reste \à consid\érer les \elts du type $x=x_2+x_0,$
avec $ x_0$ non nul appartenant \à $E_0(H_1)\cap \goth g_1$ 
et
$$ x_2   =  \sum_{1\≤j\≤q}X_{\lambda_j} \hbox{ avec  }
   [X_{\lambda_i},X_{\lambda_j}]=0\ \hbox{pour } 1\≤i<j\≤q\ \hbox{  et } 1\≤q\≤\ell-1$$ 
\noi Ainsi $(x_2,h=\sum_{1\≤j\≤q}h_{\lambda_j})$ se compl\ète en un
$sl_2$-triplet not\é $(x_2,h,x_2^{-1})$ et on peut supposer
que la d\écomposition de $x_0$ relativement \à
$adh$ est de la forme:
$x_0=y+z$ avec $y\in E_{-1}(h)\cap E_0(H_1)\cap
E_2(H_2)$ et $z\in E_0(h)\cap E_0(H_1)\cap
E_2(H_2)$ (lemme 2.2.2).\\

$\bf 1)$ Lorsque $y=0$\\

\noi On peut supposer que  $z\in \goth U(\goth s)_1,$ $\goth s$ \étant
l'alg\èbre (r\éductive dans $\goth g$) engendr\ée par
$x_2,h,x_2^{-1}$ et $H_1$ (lemme 2.2.2) et  le \PV: $(\goth U(\goth
s)_0,\goth U(\goth s)_1)$ a un nombre fini d'orbites; l'\elt $u=H_1-h$
commute \à $x$ et v\érifie $B(u,H_1)\not=0$  ce qui termine la d\émonstration dans ce cas.\\

$\bf 2)$  Lorsque $y\not=0$\\

\noi $\bf a)$ Le cas $(F_4,\alpha_4)$ non d\éploy\é\\

\noi  Les racines de $\∆_1$ sont toutes courtes et il existe $4$
racines orthogonales, $\{\lambda_i\}_{1\≤i\≤4},$  telles que:
$$\sum_{1\≤i\≤4}h_{\lambda_i}=2H_0\ ,\ \
H_1=h_{\lambda_1}+h_{\lambda_2}\ ,\ \ \mu={1\over 2}( 
 \lambda_1+\lambda_2+\lambda_3+\lambda_4)\in \∆_2 $$
$$\∆_1= \{\lambda_i\ ,\mu-\lambda_i\ \ 1\≤i\≤4\ \}  \ ,
\  \∆_2= \{\lambda_i+\lambda_j\ ,\ 1\≤i<j\≤4\ ,\mu \ \}.$$ 
\noi   
 Nous avons donc:
$x_2=X_{\lambda_1},y=X_{\mu-\lambda_1}$ et on peut supposer $z=0$ en appliquant le lemme 2.2.3 puisque  $$E_3(h)\cap \goth g_1=\{0\}\ \ \hbox{et}\ \ E_0(h)\cap E_0(H_1)\cap
E_0(h_{\mu-\lambda_1})=\{0\}.$$  
\noi L'\elt $u=h_{\lambda_2}-h_{\lambda_3} $ centralise $x$ et v\érifie $B(u,H_1)\not=0.$\\

\noi   
 $\bf  b )$ Dans tous les autres cas, $\goth g$ est d\éploy\ée,
$\Delta$ est simplement lac\é,  les racines $\lambda_i,
1\≤i\≤\ell$ sont fortement orthogonales et  
$1\≤q\≤3$  donc $E_3(h)\cap E_2(H_1)\cap \goth g_1=\{0\}$ sauf
pour $q=3$ mais alors $\ell=4$ et $E_3(h)\cap E_2(H_1)\cap \goth
g_1=
\goth g^{\beta}$ avec
$\beta={1\over 2}(\sum_{1\≤i\≤3}\lambda_i-\lambda_4).$\\

\noi Par la d\émonstration  du lemme 2.2.4  on peut  supposer que soit:\\

i) $y=X_{\gamma_1 }$ avec $\gamma_1 \in \∆_1$ tel
que $\gamma_1 (H_1)=0 ,$  
$\gamma_1 (h)=n(\gamma_1 ,\lambda_1)=-1,$  \\

\noi soit: \\

 ii) $q=2$ et $y=X_{\gamma_1 }+X_{\gamma_2 }$ 
avec $\gamma_1$ et  $\gamma_2\in \∆_1$ tels
que $\gamma_1 (H_1)=\gamma_2 (H_1)=n(\gamma_1 ,\gamma_2)=n(\gamma_2
,\lambda_1)=n(\gamma_1
,\lambda_2)=0$ et  $n(\gamma_1 ,\lambda_1)=n(\gamma_2 ,\lambda_2)=-1.$ \\

\noi On compl\ète $y$ en $1$ \Sl : $(y,h',y^{-1})$ ($h'=h_{\gamma_1 }$ dans le cas i) et $h'=h_{\gamma_1 }+h_{\gamma_2 }$ dans le cas ii)) donc:\\

iii) On peut  supposer que $z=z_0+v$ avec  $z_0\in E_0(h)\cap E_0(h')\cap
E_0(H_1)\cap
\goth g_1$ et $v\in E_3(h)\cap E_2(H_1)  \cap \goth g_1$   puisque $H_1-h$ est $1$-simple et $p=2$ (lemme 2.2.3 et sa remarque).  \\
 
\noi $\bullet$ Lorsque $q=1$  ou bien $y =X_{\gamma_1}+X_{\gamma_2}$ (donc  $x_2=X_{\lambda_1}+X_{\lambda_2}$) on a $x_2\in E_{-1}(h')\cap \goth g_1$ donc on peut
supposer que $z_0\in \goth U(\goth s)_1,$ (vrai pour $z_0$ lorsque $q=1$ et pour $q=2$ cf. d\émonstration du lemme 2.2.2)  $\goth s$ \étant
l'alg\èbre (r\éductive dans $\goth g$) engendr\ée par $H_1$
et les deux
$sl_2$-triplets $(x_2,h,x_2^{-1})$ et
$(y_{-1},h',y_{-1}^{-1}),$ ce qui donne le nombre fini
d'orbites.

\noi Dans les deux cas l'\elt 
$u=-3H_1+4h+2h'$ commute \à $x$ et v\érifie   $B(u,H_1)\not=0.$
 \vskip 2mm

\noi Il reste donc les cas suivants:

\noi $\bullet$  Lorsque $q\≥2$ et $y_{-1}=X_{\gamma_1},$  
 on r\éduit encore
$z_0+v$ relativement \à
$ad(h_{\lambda_1}).$\\

\noi Soit $\delta\in s(z_0),$ comme
$n(\delta,\lambda_1+\gamma_1)=n(\delta,\lambda_1)\≥0$ car
$\∆_3=\emptyset, $ on a $z_0=w_0+w_1$ avec  
$w_i\in E_i(h_{\lambda_1})\cap E_0(h)\cap E_0(h')\cap
E_0(H_1)\cap \goth g_1$ pour $i=0,1$ donc il existe $A\in
\oplus_{\delta\in s(w_1)}\goth g^{\delta-\lambda_1}$ ($\subset$ $E_{-1}(h_{\lambda_1})\cap E_{-2}(h)\cap E_1(h')\cap
E_{-2}(H_1)\cap \goth g_0$) tel que
$[A, X_{\lambda_1}]=-w_1$ d'o\ù
$$
x'=exp(ad(A))(x)=x_2+X_{\gamma_1}+w_0+[A,v]+ [A,\sum_{2\≤i\≤q}X_{\lambda_i}]+v  .
$$
\noi On supprime l'\elt $[A,v]\in E_1(h)\cap E_0(H_1)\cap \goth
g_1$ \à l'aide d'un \elt convenable  de $ exp(ad(E_{-1}(h)\cap E_{-2}(H_1)\cap \goth
g_0)). $ Comme  $[A,\sum_{2\≤i\≤q}X_{\lambda_i}]\in E_1(h')\cap E_0(H_1)\cap E_0(h)\cap \goth g_1,$  \à l'aide du lemme 2.2.2 on  peut supposer que   $x'=x_2+X_{\gamma_1}+w_0+v'$ avec $v'\in E_2(H_1)\cap E_3(h)\cap \goth g_1
 .$
Lorsque $v'\not=0,$ donc $q=3,$  on supprime $v'$ \à l'aide de $exp(ad(\goth
g^{\beta-\lambda_1})$ dont la restriction \à $
E_0(h_{\lambda_1})\cap E_0(h_{\lambda_4}))\cap E_0(H_1)\cap
\goth g_1$  est l'identit\é, ce qui donne finalement:
$x"= x_2+X_{\gamma_1}+w_0.$\\

\noi  Or $w_0\in \goth U(\goth s')_1,$ $\goth s'$ \étant
l'alg\èbre (r\éductive dans $\goth g)$ engendr\ée par
$H_1-h,h_{\gamma_1},h_{\lambda_1}$ et 
 $h-h_{\lambda_1},$ donc on en d\éduit le nombre fini d'orbites de $P_\goth t$ dans $\g_1.$\\
 
 \noi L'\elt $u=3(H_1-h)-2h_{\gamma_1}- h_{\lambda_1}$
commute avec $x"$ et v\érifie
$B(u,H_1)=(3(l-q)-1)B(h_{\lambda_1},h_{\lambda_1})\not=0.$\fdem

 \newpage

\section{\bf  Fonctions Z\étas associ\ées}
\bigskip 

\noi A l'aide  des travaux de \cite{boppruben},  de \cite{sato3} et du th\éor\ème 2.1.1,  il est facile d'\établir  l'existence d'une \équation fonctionnelle abstraite v\érifi\ée par la fonction Z\éta locale associ\ée au \PV $(P_{\goth t},\goth g_1)$ d\éfinit dans le \S 1.4 et  c'est l'objet de ce paragraphe. 

\noi La plupart des d\émonstrations sont analogues \à celles faites dans le cas commutatif r\éel (\cite{boppruben}).

\noi  On reprend toutes les hypoth\èses du \S 1.4, c'est \à dire:

\begin{enumerate}
\item $\goth g$ est une alg\èbre absolument simple engendr\ée par $\goth g_{\pm 1}$,

\item $ ({\goth g}_0, {\goth g}_1,H_0)$ est un \PV  absolument irr\éductible et r\égulier,  

\item $P_{\goth t}=P(H_1,...,H_p)$ est un \sg parabolique tr\ès sp\écial lorsque $p\≥2$ et pour $p=1$ on pose $P(2H_0)=G.$

\end{enumerate}

\bigskip

\noi On commence par  pr\éciser quelques normalisations usuelles avant d'\établir l'existence de l'\équation fonctionnelle abstraite v\érifi\ée par les fonctions Z\étas.\\

\subsection{\bf  Normalisation de  la forme de Killing}

\bigskip

\noi On convient de poser $\widetilde{B}=cB$ (avec $c$ r\éel dans le cas complexe) et on choisira dans les applications $c=-\ds\frac{(\text{degr\é de }F_p)}{2B(H_0,H_0)},$ ce qui correspond \à une normalisation de la forme
de Killing de $\goth g$ ind\épendante du choix de $P_{\goth t}.$   \\

\noi Soit $\tau$ un caract\ère de $\F$ qui d\éfinit la
transformation de Fourier des fonctions de $ {\EuScript  S}(\F),$ espace de Schwartz de $\F$:
 
\begin{enumerate}
 \item$\tau$ est un caract\ère d'ordre $\nu$ lorsque $\F$ est un corps
$\goth P$-adique;

\item $\tau (x)=e^{2i\pi x}$ dans le cas r\éel 

\item$\tau (z)=e^{2i\pi(z+\overline{z})}$ dans le cas complexe. 
\end{enumerate}
\vskip 2mm
\noi La transformation de Fourier ${\EuScript  F}$ d'une fonction $f\in  {\EuScript  S}(\goth g_1)$ est alors la fonction ${\EuScript  F}(f)\in {\EuScript  S}(\goth g_{-1})$ d\éfinie par: 
 $${\EuScript  F}(f)(y)=\int_{\goth g_1}f(x)\tau(\widetilde B(x,y))dx\quad  y\in \goth g_{-1}.$$
\noi La  transformation de Fourier inverse d'une fonction  $g\in  {\EuScript  S}(\goth g_{-1})$ est donn\ée par:
$$ {\overline{\EuScript  F}}(g)(x)=\int_{\goth g_{-1}}  \overline{\tau(\widetilde B(x,y))}g(y)dy\quad  x\in \goth g_{1}.$$

\noi  Si
$\goth U$ est une alg\èbre de Lie semi-simple gradu\ée , on
note 
$B_{\goth U}$ la forme de Killing associ\ée; la
transformation de Fourier est alors d\éfinie par le caract\ère
:
$$\tau_{\goth U}=\tau\circ(\widetilde{B_{\goth U}})$$
\subsection {\bf Normalisation des mesures de $\goth g_1$ et
$\goth g_{-1}$ }

\bigskip

\noi Pour $x\in \F, |x|_{\F}$ (not\é \également $|x|$) d\ésigne la valeur absolue de $x$ c'est \à dire :
$$\hbox{lorsque}\ \F=\C\ :\ |x|_{\C}=x\overline x\ \ ,\ \ \hbox{lorsque}\ \F=\R\ :\ |x|_{\R}=\hbox{max}(x,-x).$$
Lorsque
$\F$ est un corps
$\goth P$-adique, on note $q$  sa caract\éristique r\ésiduelle,   $\goth O$ l'anneau des entiers de $\F,$ $\goth P$ l'unique id\éal maximal de $\goth O$ et $\pi$ un g\én\érateur de $\goth P$ alors $x$ se d\écompose: $x=\pi^nu,$ $u\in \goth O$ ($n=$ordre de $x$) et $ |x|_{\F}=q^{-n}.$\\

\noi On munit une fois pour toute $\F$ d'une mesure autoduale pour $\tau.$\\

 \noi On rappelle que lorsque $\F=\C$ , on consid\ère les structures d'espaces vectoriels r\éels de $\goth g_1$ et $\goth g_{-1}.$\\  

 On  
 note $\Phi$ l'application  d\éfinie sur $\goth
g'_1$ qui \à $x$ associe $x^{-1}$ l'unique  \elt de $\goth g'_{-1}$ tel
que  $(x,2H_0,x^{-1})$ soit un $sl_2$-triplet ($[x^{-1},x]=2H_0$), $\Phi$
commute \à l'action de $G$ et on  normalise   les polynomes
irr\éductibles fondamentaux des  \PVs:
$(G,\goth g_1)$ et $(G,\goth g_{-1})$ par: $$
F^*(\Phi(x))=\frac{1} { F(x)}\quad (N1).$$ 
\noi Pour $x\in \goth
g'_1,$ notons $\theta_x=\theta_x(-1)$ (cf.d\éfinition 7, introduction du $\S 1$) et  rappelons que $\theta_x$ est une involution de $\goth g$ qui v\érifie $\theta_x(x)=\Phi (x)=x^{-1}.$ 

\noi Comme
$\theta_{gx} =g\circ\theta_{x} \circ g^{-1},$ $det(\theta_x/ 
 \goth g_1)$
est une  application rationnelle relativement invariante par $G$
(de caract\ère associ\é $(det(g/\goth g_{-1})^2$)  il
existe une constante $c$ telle que $det(\theta_x/ 
 \goth g_1)=cF(x)^{-2N}$ avec $N=\displaystyle{\frac{\hbox{dim(}\goth g_1)}
 {\hbox{degr\é de\ } F}}.$
\vskip 2mm
\noi Notons que pour   $x_0\in \goth
g'_1,$ on a $G_{x_0}=G_{x_0^{-1}}$ et  $\theta_{x_0}$ \échange l'orbite   $G.x_0$ dans $\goth g_1$ avec l'orbite $G.x_0^{-1}$ dans $\goth g_{-1}$ car l'application $\Theta_{x_0}$ d\éfinie pour $g\in G$ par: $
\Theta_{x_0} (g)=\theta_{x_0} .g.\theta_{x_0} ^{-1},$ est un automorphisme de $G.$

\noi Pour $x\in \goth g_1,$ l'application polynomiale
$F^*(\theta_{x_0}(x))$ est   non identiquement nulle,
relativement invariante par $G$ donc par la normalisation (N1) on a
$F^*(\theta_{x_0}(x))=\displaystyle{\frac{F(x)} {F(x_0)^2}} $ (3).\\

\begin{lem} Il existe une unique mesure de Haar, $dx,$ sur $\goth g_1$ et une unique mesure de Haar, $dy,$ sur $\goth g_{-1}$ telles que:
 
(1) $\overline{\four}\circ\four=Id_{\EuScript  S(\goth g_1)}.$

(2) $\forall f\in L^1(\goth g_{-1})$ et $\forall x_0\in \goth g'_1$ on a
$\int_{\goth
g_1}f(\theta_{x_0}(x))dx= |F(x_0)|^{2N 
}.
\int_{\goth
g_{-1}}f(y)dy .$
 
\noi Pour $g\in L^1(\goth g'_{-1})$ et $\forall x_0\in \goth g'_1$ on a
$\int_{\goth
g_1}g(\theta_{x_0}(x))|F(x)|^{-N}dx=   
 \int_{\goth
g_{-1}}g(y)|F^*(y)|^{-N}dy .$
 
 \end{lem}

\dem Munissons $\goth g_1$ et  $\goth g_{-1}$ de
deux bases :  $\goth B$    et
  $\goth B^*,$  duales pour $\tilde{B},$
 et soient $d_{\goth B}x$ et   $d_{\goth B^*}y$ les mesures
de Haar correspondantes alors il existe une constante
positive $\lambda$ telle que $dx=\lambda d_{\goth B}x$ donc 
par la propri\ét\é de dualit\é demand\ée dans (1) on a $dy={1\over \lambda}
d_{\goth B^*}y$ puisque $\F$ est muni d'une mesure de Haar autoduale pour $\tau,$ et par changement de variable
  on a:
$$\int_{\goth
g_1}f(\theta_{x_0}(x))dx=\frac{\lambda^2}{  |c|}.|F(x_0)|^{2N 
}.
\int_{\goth
g_{-1}}f(y)dy\ \ \hbox{d'o\ù}\ \
\lambda=\sqrt{|c|}(=|det(\theta_{x_0}/\goth g_1)|^{1\over
2}|F(x_0)|^N)  $$
par (2) et (3). \fdem\\

\noi On munit $\goth g_1$ et $\goth g_{-1}$ des mesures uniques d\éfinies dans le lemme 3.2.1.\\

\noi Pour $x_0\in \goth g'_1,$ $G/G_{x_0}$ est muni de la mesure $G-$invariante qui v\érifie pour  $f\in \EuScript C_C(G.x_0)$ :
$$\int_{\goth
g_1}f(x)|F(x)|^{-N }dx=\int_{G /
G_{x_0}}f(gx_0)d{\dot g}\quad (N2)$$
\noi puisque  $|F(x)|^{-N }dx$ est une mesure
$G$-invariante.\\

\noi Notons que les  mesures $G-$invariantes $|F^*(y)|^{-N }dy$ et
$|F(x)|^{-N }dx$ sont  ind\épendantes du choix de
$\goth B$ et de la normalisation de $F$ mais la mesure $dx$
d\épend de l'invariant relatif  choisi dans l'\égalit\é  (2) du lemme pr\éc\édent,
on dira que la mesure est normalis\ée \à l'aide de $F$ et ceci
 sera pr\écis\é \à chaque fois.\\
 
 \noi Avec la normalisation des mesures impos\ées par les $4$ conditions: (N1), (N2) et (1) et (2) du lemme 3.2.1, on a les propri\ét\és suivantes:\\

\begin{lem}  Pour $  f\in {\EuScript S}( \goth g'_{-1})$ on a les relations suivantes:
\begin{enumerate}

 \item Si $\sigma$ est une involution de $\goth g$ telle que:
 
 \noi $ \sigma(H_0)=-H_0$ et $\{x\in \goth g_1|(x,2H_0, \sigma (x))$ est un \Sl $\}\not=\emptyset$ alors:
 \vskip 5pt  
\noi  $ \int_{\goth g_1}f\circ  \sigma(x)|F(x)|^{-N }dx=
 \int_{\goth g_{-1}}f(y)|F^*(y)|^{-N }dy .$
 
\item   $\int_{G /
G_{x_0^{-1}}}f(g{x_0}^{-1})d{\dot g}= \int_{G.{x_0}^{-1}
}f(y) 
 |F^*(y)|^{-N }dy.$
 
\item    
$ \int_{\goth g_1}f\circ \Phi(x)|F(x)|^{-N }dx=
 \int_{\goth g_{-1}}f(y)|F^*(y)|^{-N }dy .$
 
\end{enumerate}

 \end{lem}  
 
\dem   1) Soit $x_0$ tel que $(x_0,2H_0, \sigma (x_0))$ soit un \Sl alors:
$$ \int_{\goth g_1}f\circ \sigma (x)|F(x)|^{-N }dx=\int_{\goth g_1}f\circ ( \sigma\theta_{x_0})(\theta_{x_0}(x))|F(x)|^{-N }dx=\int_{\goth g_{-1}}f ( \sigma\theta_{x_0} (y))|F^*(y)|^{-N }dy,$$
par le lemme pr\éc\édent, d'o\ù le r\ésultat par changement de variable
 car $ \sigma\theta_{x_0}$ est une involution de $\goth g_{-1}$ puisque $\sigma$ et  $\theta_{x_0}$ commutent.\\

\noi 2) Soit:$$A=\int_{G /
G_{x_0^{-1}}}f(gx_0^{-1})d{\dot g}= \int_{G/
G_{x_0}}f(g\theta_{x_0}(x_0))d{\dot g}= \int_{G/
G_{x_0}}f\circ \theta_{x_0}(\Theta_{x_0}(g)x_0)d{\dot g}\ ,$$
\noi Comme
$\Theta_{x_0}^2$ est un automorphisme int\érieur de $G$ et que $G$ est
unimodulaire, on a:
$$A=  \int_{G/
G_{x_0}}f\circ \theta_{x_0}( gx_0)d{\dot g}=\int_{\goth g_1}{\bf 1}_{G.x_0^{-1}}(\theta_{x_0}(x)) f\circ \theta_{x_0}(x)|F(x)|^{-N}dx \ \hbox{par }(N2)$$
$$A=  \int_{G.{x_0^{-1}}}f( y)|F^*(y)|^{-N}dy\quad \hbox{par le lemme pr\éc\édent}.$$ 
\noi 3) Soit $\{x_1,...,x_m\}$ un ensemble de repr\ésentants des orbites de $G$ dans $\goth g'_1$ et soit:
$$B= \int_{\goth g_1}f\circ \Phi(x)|F(x)|^{-N }dx=\sum_{i=1}^m \int_{G.x_i}f\circ \Phi(x)|F(x)|^{-N }dx=\sum_{i=1}^m \int_{G/G_{x_i}}f(g^{-1}x_i^{-1}) d{\dot g},$$
\noi donc par le changement de variable: $g\rightarrow g^{-1}$ et en appliquant le r\ésultat pr\éc\édent, on obtient:
$$B=\sum_{i=1}^m \int_{G/G_{x_i^{-1}}}f(gx_i^{-1}) d{\dot g}= \sum_{i=1}^m \int_{G.x_i^{-1}}f (y)|F^*(y)|^{-N }dy =\int_{\goth g_{-1}}f (y)|F^*(y)|^{-N }dy .\hfill \Box $$\\

\subsection {\bf Normalisation des
invariants relatifs fondamentaux}
\bigskip

\noi  On rappelle les diff\érents \PVs qui interviennent:
\begin{enumerate}
\item   $(P_{\goth t},\goth g_1)$ avec les 
 invariants relatifs fondamentaux $F_1,...,F_p$    et le  \PV dual 
: $(P_{\goth t},\goth g_{-1})$  pour lequel les 
 invariants relatifs fondamentaux sont donn\és par $F^*_1,...,F^*_p$ 
 ($P_{\goth t}=P(H_1,...,H_p)$),
\item  $( {P_{\goth t}}^-,\goth g_1)$  avec les 
 invariants relatifs fondamentaux $P_1,...,P_p$ et le   \PV dual 
: $(   {P_{\goth t}}^-,\goth g_{-1})$   pour lequel les 
 invariants relatifs fondamentaux sont donn\és par $P^*_1,...,P^*_p$ ($   {P_{\goth t}}^- =P(H_p,...,H_1)$).
\end{enumerate}

\bigskip

\noi On normalise tous les invariants associ\és \à l'action du parabolique oppos\é \à partir des invariants associ\és \à $P_{\goth t}.$\\

\noi De mani\ère pr\écise, \étant donn\é 
$F_1,...,F_p,$ on choisira
 $P_1,...,P_k,$   normalis\és de
fa\çon que pour $k=1,...,p,$ $x\in E_2(h_k)\cap \goth g_1,$
$y\in E_0(h_k)\cap \goth g_1$ qui commutent on ait : 
$$(\ F_p  (x+y)\ )^{\frac{1}{m_k}}=F_k(x)P_{p-k}(y)\quad \hbox{(R1)}$$
avec $m_k=1$ \à l'exception des $\F$-formes de
$(E_7,\alpha_6)$ et du cas classique $(C_n,\alpha_k)$ avec le parabolique $P'_0$ pour lesquels $m_k=\frac{1}{2}$ (cf.lemme 1.4.7).\\

\noi 
Ceci d\étermine la normalisation des invariants relatifs
fondamentaux attach\és au pr\éhomog\ène  ``dual" en posant pour $z\in
 \goth g_{-1},$ pour
$k=1,...,p-1$ et pour $x\in E'_0(h_{p-k})\cap \goth g_{-1}$
(resp.$y\in E'_{-2}(h_{p-k})\cap \goth g_{-1}$):
$$F^*(z)=\frac{1} {
F(z^{-1})}\quad  ,\quad F^*_k(x)=\frac{1 }{
P_k(x^{-1})}\quad (\hbox{resp.}\quad  P^*_{p-k}(y)=\frac{1}{
F_{p-k}(y^{-1})}\ )\quad \hbox{(R2)}\ ,$$
qui correspond simplement \à la normalisation (N1) pour chaque \PV $(\goth U^+(k)_0,\goth U^+(k)_1)$ avec son ``dual" et pour chaque 
 \PV $(\goth U^-(k)_0,\goth U^-(k)_1)$ avec son ``dual".\\
 
\noi On aura  donc \également les relations :
  $$(\ F^*_p(x+y)\ )^{\frac{1}{m_k}}=F^*_k(x)P^*_{p-k}(y)\quad \hbox{(R3)}$$
lorsque $x$ et $y$ commutent.\\
 
\noi R\éciproquement, la donn\ée de
$F^*_1,...,F^*_p,$ d\éfinit \également tous les autres
invariants relatifs fondamentaux \à partir des relations impos\ées (R3) et R(2) qui impliqueront la relation (R1).\\

\noi Une fa\çon d'obtenir cette normalisation consiste \à choisir un  \elt
$X_0=\sum_{1\≤i\≤p}X_i\in W_{\goth t},$  puis de
normaliser tous les polynomes  en posant:  
$$k=1,...,p:\ F_k(X_0)=F^*_k(\Phi (X_0))=1\ ,\
k=1,...,p-1:\ P_k(X_0)=P_k^*(\Phi (X_0))=1.$$ 
 
  {\bf Dans le cas  complexe} on fixe un syst\ème de Chevalley $(X_{\alpha})_{{\alpha}\in \∆}$ alors $\goth g=\goth g_{\R}\otimes_{\R}\C,$ $\goth g_{\R}=\sum_{\alpha \in \∆}\R h_{\alpha}\oplus_{\alpha \in \∆}(\R X_{\alpha}\oplus \R X_{-\alpha})$ \étant l'alg\èbre simple r\éelle construite avec le syst\ème de Chevalley pr\éc\édent.  
    
\noi On choisira toujours les  \irfs, $F_1,...,F_p,$  r\éels sur $\goth g_{\R},$ ce qui est toujours possible puisque ceux-ci sont les extensions \à $\goth g_1$ des \irfs des \PVs irr\éductibles $(\goth U^+(k)_0, \goth U^+(k)_1)$ et que ces \PVs sont d\éfinis sur $\R$ puisque $h_k\in \goth g_{\R}$  pour $k=1,...,p$ (lemme 1.1 p.135 de \cite{sato-shintani}).

\noi Ainsi $W_{\goth t,\R}=\{x\in W_{\goth t}\ |\ F_k(x)\in\R^*$ pour $k=1,...,p\}\not=\emptyset.$ 

\noi  Soit $\overline x$ la conjugaison sur $\goth g$  associ\ée \à $\goth g_{\R},$ on a alors pour $x\in \goth g_1$  et $y\in \goth g_{-1},$ $F_k(\overline x)=\overline{F_k(x)}$ et  $F^*_k(\overline y)=\overline{F^*_k(y)}$ pour $k=1,...,p.$\\
 
 \subsection{\bf  Le cas  archim\édien } 
 
 \bigskip
 
 \noi  Dans les cas r\éels et complexes et pour $k=1,...,p,$ on d\éfinit les op\érateurs diff\érentiels
usuels (cf.par exemple \cite{sato-shintani}, \cite{farautkoranyi}):\\

\begin{enumerate}

\item $F_k(\partial)$
 est l'op\érateur diff\érentiel \à coefficients
constants  d\éfini  sur $\goth g_{-1}$   par
 $ F_k(\partial)e^{
\widetilde{B}(x,y)}=F_k(x)e^{
\widetilde{B}(x,y) }$  
et $F^*_k(\partial)$ celui d\éfini  sur $\goth g_1$ par
 $F^*_k(\partial)e^{
\widetilde{B}(x,y)}=F^*_k(y)e^{
\widetilde{B}(x,y) },$ avec $x\in \goth g_1$ et  $y\in \goth g_{-1}.$\\

\item {\bf Lorsque $\bf \F=\C$},  $F_k(\overline{\partial})$
 est l'op\érateur diff\érentiel \à coefficients
constants  d\éfini  sur $\goth g_{-1}$   par
 $ F_k(\overline{\partial})e^{\overline{
\widetilde{B}(x,y)}}=\overline{F_k(x)}e^{
\overline{\widetilde{B}(x,y)} }$  
et $F^*_k(\overline{\partial})$ celui d\éfini  sur $\goth g_1$ par
 $F^*_k(\overline{\partial})e^{
\overline{\widetilde{B}(x,y)}}=\overline{F^*_k(y})e^{\overline{
\widetilde{B}(x,y) }},$ avec $x\in \goth g_1$ et  $y\in \goth g_{-1}.$

\noi Notons $\alpha_1,...\alpha_{Nd_p}$ les racines de $∆_1$ et $F(\sum_{1\≤i\≤Nd_p}x_iX_{\alpha_i}\ )=\sum_{n=(n_1,...,n_{Nd_p})}a_n\prod_{1\≤i\≤Nd_p}x_i^{n_i}$ l'expression de $F_k$ dans un syst\ème de Chevalley de $\g$ alors:

\noi $ F_k(\partial)=\sum_{n=(n_1,...,n_{Nd_p})}a_n\prod_{1\≤i\≤Nd_p}c_i^{n_i}\prod_{1\≤i\≤Nd_p}(\ds\frac{\partial}{\partial x_i})^{n_i} $ et 

\noi $ F_k(\overline{\partial})=\sum_{n=(n_1,...,n_{Nd_p})}a_n\prod_{1\≤i\≤Nd_p}c_i^{n_i}\prod_{1\≤i\≤Nd_p}(\ds\frac{\partial}{ \partial \overline{x_i}})^{n_i} ,$

\noi $c_i$ \étant le nombre r\éel  $c_i=\widetilde B(X_{\alpha_i},X_{-\alpha_i}),$ en prenant la base duale dans $\g_{-1}$ relativement \à  $\widetilde B.$

\end{enumerate}

\bigskip

\noi $\ell $ d\ésignant la repr\ésentation r\éguli\ère gauche de $P_{\goth t}$ sur ${\EuScript C}^{\infty}(\goth g_{\pm 1})$:
$\forall p\in P_{\goth t}\quad l(p)f)(x)=f(p^{-1}x),$ on a pour $k=1,...,p:$
$$\forall p\in P_{\goth t} :\quad \ell(p)\circ F_k(\partial)=\chi_k(p)^{-1}F_k(\partial)\circ \ell(p)\quad \hbox{et}\quad \ell(p)\circ F^*_k(\partial)=\chi^*_k(p)^{-1}F^*_k(\partial)\circ \ell(p)\quad (*).$$

\noi On rappelle \également leur action sur la transformation de Fourier:\\ 
  
\begin{lem}
1.- Soit $f\in S(\goth g_1),$  pour $m>0$ et
$k=1,...,p,$ on a:
$$\four((F_k)^m.f)=(2i\pi)^{-md_k}F_k(\partial)^m({ \four 
f}) \ ,\  \four (\ F^*_k(\partial)^mf\
)=(-2i\pi)^{md'_{p-k}}({F^*_k})^m.{\four f},$$
2.- Lorsque $\F=\C$ on a \également:
$$\four((\overline{F_k})^m.f)=(2i\pi)^{-md_k}F_k(\overline{\partial})^m({ \four 
f}) \ ,\  \four (\ F^*_k(\overline{\partial})^mf\
)=(-2i\pi)^{md'_{p-k}}(\overline{{F^*_k}})^m.{\four f}.$$
$d_k$  \étant le \dg  de 
$F_k$ et $d'_{p-k}$ celui de
$F^*_k.$\\
 \end{lem}
 
 \bigskip
 
 \noi    Soit $s=(s_1,...,s_p)\in \C^p,$ on
pose
$F^s=\prod_{1\≤i\≤p}F_i^{s_i}\ ,\ 
 {F^*}^s=\prod_{1\≤i\≤p}{F_i^*}^{s_i}  ,$ 
  $|F| ^s=\prod_{1\≤i\≤p}|F_i|_{\F}^{s_i}$ et 
$ |F^*| ^s=\prod_{1\≤i\≤p}{|F_i^*|_{\F}}^{s_i}  .$\\ 

\noi $\omega_{-1}$ est le caract\ère signe d\éfini sur $\R^*$ par $\omega_{-1}(x)=\frac{x}{|x|},x\in \R^*.$\\
 
  \bigskip
\begin{lem}
1.- Pour $k=1,...,p$ il existe des polynomes $b_k$ et $b^*_k\in
\R[s]$ tels que 
$$F_k(\partial)F^{*s}=b_k(s)F^{*s-1_p+1_{p-k}}\ \ \hbox{et}\
\ F^*_k(\partial)F^s=b^*_k(s)F^{s-1_p+1_{p-k}},$$
 \à l'exception des formes r\éelles de
$(E_7,\alpha_6)$ et du cas classique $(C_n,\alpha_k)$ avec le parabolique $P'_0$ pour lesquels $p=2$ et
les relations pour
$k=1$ deviennent : 
$$F_1(\partial)F^{*s}=b_1(s)F_1^{*s_1+1}F_2^{*{s_2-2}}\quad \hbox{et}\quad 
 F^*_1(\partial)F^s=b^*_1(s)F_1^{s_1+1}F_2^{s_2-2}.$$
  2.- Lorsque $\F=\C$ on a \également pour $k=1,...,p:$  $$F_k(\overline{\partial})\overline{F^*}^s=b_k(s)\overline{F^*}^{s-1_p+1_{p-k}}\ \ \hbox{et}\
\ F^*_k(\overline{\partial})\overline{F}^s=b^*_k(s)\overline{F}^{s-1_p+1_{p-k}},$$
$$F_k( {\partial})F_k(\overline{\partial}) |F^*|^s=b_k(s)^2 |F^*|^{s-1_p+1_{p-k}}\ \ \hbox{et}\
\ F^*_k( {\partial})F^*_k(\overline{\partial}) |F|^s=b^*_k(s)^2 |F|^{s-1_p+1_{p-k}},$$
 \à l'exception des formes r\éelles de
$(E_7,\alpha_6)$ et du cas classique $(C_n,\alpha_k)$ avec le parabolique $P'_0$ pour lesquels $p=2$ et
les relations pour
$k=1$ deviennent : 
$$F_1(\overline{\partial})\overline{F^*}^s=b_1(s)\overline{F_1}^{*s_1+1}\overline{F_2}^{*{s_2-2}}\quad \hbox{et}\quad 
 F^*_1(\overline{\partial})\overline{F}^s=b^*_1(s)\overline{F_1}^{s_1+1}\overline{F_2}^{s_2-2},$$
 $$F_1( \partial)F_1(\overline{\partial}) |F^*|^s=b_1(s)^2 {|F_1|_{\C}}^{*s_1+1} {|F_2|_{\C}}^{*{s_2-2}}\quad \hbox{et}\quad 
 F^*_1( \partial)F^*_1(\overline{\partial}) |F|^s=b^*_1(s)^2 {|F_1|_{\C}}^{s_1+1}{ |F_2|_{\C}}^{s_2-2}.$$
 3.- Lorsque $\F=\R$ on a \également pour $k=1,...,p:$  $$F_k( {\partial}) |F^*|^s=b_k(s)\omega_{-1}(F^*_p.F^*_{p-k}) |F^*|^{s-1_p+1_{p-k}}\ \ \hbox{et}\
\ F^*_k( {\partial}) |F|^s=b^*_k(s)\omega_{-1}(F_p.F_{p-k}) |F|^{s-1_p+1_{p-k}},$$
 \à l'exception des formes r\éelles de
$(E_7,\alpha_6)$ et du cas classique $(C_n,\alpha_k)$ avec le parabolique $P'_0$ pour lesquels $p=2$ et
les relations pour
$k=1$ deviennent : 
$$F_1( {\partial}) |F^*|^s=b_1(s) \omega_{-1}(F^*_1){|F_1|_{\R}}^{*s_1+1} {|F_2|_{\R}}^{*{s_2-2}}\quad \hbox{et}\quad 
 F^*_1( {\partial}) |F|^s=b^*_1(s) \omega_{-1}(F_1){|F_1|_{\R}}^{s_1+1} {|F_2|_{\R}}^{s_2-2}.$$
 ($1_l$ d\ésignant l'\elt de
$\N^p$ dont toutes les composantes sont nulles \à
l'exception de la $l$-i\ème qui vaut
$1$ et $1_0$ ayant toutes les composantes nulles).\\
\end{lem}

\dem  1. Donnons-la dans le cas g\én\éral pour la premi\ère \égalit\é. Par la relation (*) et le lemme 1.4.7, $F_k(\partial)F^{*s}$  et $F^{*s-1_p+1_{p-k}}$ sont relativement invariants par $\ell$ de m\ême caract\ère donc pour $s\in \Z^p$ ils sont proportionnels d'apr\ès la proposition 1.4.5,  ainsi il existe $f_k:\Z^p\rightarrow \C$ telle que $\forall s\in \Z^p$ on ait: $F_k(\partial)F^{*s}=f_k(s)F^{*s-1_p+1_{p-k}}.$ Soit $\EuScript D$ le domaine de d\éfinition de $F^{*s};$ pour $x\in   \EuScript D$ et $s\in \C^p$ posons $b_{k,x}(s)= \displaystyle{\frac{(F_k(\partial)F^{*s})(x)}{F^{*s-1_p+1_{p-k}}(x)}},$ alors $b_{k,x}\in {\EuScript C}^{\infty}( \EuScript D)[s]$   et v\érifie $\forall n\in \Z^p:$  $b_{k,x}(n)=f_k(n)$ donc $b_{k,x}$ est ind\épendant de $x$, et \à coefficients r\éels puisque $F_k$ et $F^*_k$
sont \à coefficients r\éels pour $k=1,...,p.$

\noi 2 et 3. d\écoulent de 1.  \fdem\\

\begin{rema} Les polynomes $b_k$
et $b^*_k,$ $k=1,...,p,$ ne d\épendent pas de la
normalisation de
$F_1,...,F_p$ mais d\épendent de $\tilde B$ (qui les d\éfinit \à un facteur multiplicatif pr\ès).  

\noi On a $b_k(s_1,...,s_{p-1},0)=0$ puisque
$F_k(\partial)F^{*s}$ est un polynome pour $s\in \N^p.$
\end{rema}

\bigskip

\noi  Pour $s\in \C$ on note $b_{\goth g }(s)=b_p(0,...,s)$ le
polynome de Bernstein obtenu \à partir de de l'action de $F_p(\partial)$ sur ${F_p}^{*s}$ c'est \à dire  celui associ\é au 
\PV:
$( G ,\goth g_1)$  (et appel\é usuel)  et pour $s\in \C^p$ soit $b_{\goth g,P_{\goth t} }(s)=b_p(s_1,...,s_p)$ le
polynome de Bernstein  obtenu \à partir de l'action de $F_p(\partial)$ sur $\prod_{1\≤i\≤p}{F_i}^{*s_i}$ dans le \PV $(P_{\goth t},\goth g_1).$\\

\noi Dans le but de montrer que les polynomes 
$b_{\goth g,P_{\goth t}}$ sont des produits de polynomes de Bernstein  usuels associ\és \à  certains \PVs (provenant de centralisateurs d'alg\èbres de type  $sl_2$), on  relie les  polynomes  $b_k$  \à certains polynomes de type  $``b_{\goth g,P_{\goth t} }" $ dans la proposition qui suit.

\noi Pour ceci on introduit les notations suivantes:\\

\noi Soient:
 \begin{enumerate}
\item $X_0=\sum_{1\≤i\≤p}X_i\in W_{\goth t},$
$X' =\sum_{1\≤i\≤k}X_i$ la projection de
$X_0$  sur
$\goth U^+(k)_1=E_2(h_k)\cap
\goth g_1$ et  
$X =\sum_{k+1\≤i\≤p}X_i$ la
projection de
$X_0$   sur
  $\goth U^-(p-k)_1=E_0(h_k)\cap
\goth g_1,$ avec $1\≤k\≤p-1,$

\item $Y'$ la projection de $X_0^{-1}$   sur $\goth U^+(k)_{-1}=E_{-2}(h_k)\cap
\goth g_{-1}$ et 
  $Y$ la
projection de
$X_0^{-1}$   sur
 $\goth U^-(p-k)_{-1}=E_0(h_k)\cap
\goth g_{-1},$ 

\item $\goth s'_k$  l'alg\èbre engendr\ée par $X'$ et $Y'$ et $\goth s_{p-k} $
l'alg\èbre engendr\ée par  $X$
et
$Y$, 
 
\item $\goth U=\goth U(\goth s_{p-k})\subset \goth U^+(k)$ et 
$\goth U'=\goth U(\goth s'_k)\subset \goth U^-(p-k).$ 

\item  $A_k=\displaystyle{\frac{\widetilde B(h_k,h_k)}{\widetilde {B_{\goth U}}(h_k,h_k)}}$ et $B_k=\displaystyle{\frac{\widetilde B(2H_0-h_k,2H_0-h_k)}{\widetilde {B_{\goth U}}(2H_0-h_k,2H_0-h_k)}}$ alors $\widetilde B/\goth U=A_k\widetilde {B_{\goth U}}$ et $\widetilde B/\goth U'=B_k\widetilde {B_{\goth U'}}$ lorsque $\goth U$ et $\goth U'$ sont absolument irr\éductibles ce qui est toujours v\érifi\é lorsqu'on suppose  que dans le cas orthogonal $(\overline\g_0,\overline\g_1)$  on a la condition $3k\≤2n-2$ pour le type $(\overline \∆,\overline \lambda_0)=(D_n,\alpha_k)$ et $3k\≤2n-1$ pour le type $(\overline \∆,\overline \lambda_0)=(B_n,\alpha_k).$

\end{enumerate}

\noi On note le lieu ``non singulier" des \PVs associ\és \à l'action de \sgs paraboliques par $"$  par exemple $\g_1"=\g_1-S_{P_{\goth t}}$ et $\g_{-1}"=\g_{-1}-{S^*}_{P_{\goth t}}.$
\bigskip 

\noi  Alors on a:

 \vskip 5pt

\begin{prop}
 $$b_k(s_1,...,s_p)=A_k^{-d_k}.b_{\goth U 
 ,P(H_1,...,H_k)}  (s_{p-k+1},...,s_p)$$  
$$b^*_{p-k}(s_1,...,s_p)=B_k^{-d'_{k}}.b^*_{\goth U' 
 ,P(H_{k+1},...,H_p)}  (s_{k+1},...,s_p)$$ 
 \à l'exception des $\R-$formes du cas
exceptionnel $(E_7,\alpha_6)$ et du cas classique $(C_n,\alpha_k)$ avec le parabolique $P'_0$  pour lesquels les relations  
  deviennent : 
$$ b_1(s_1,s_2)=A_1^{- d_1}b_{\goth U 
 }(s_2)b_{\goth U 
 }(s_2-1)\quad ,\quad    b^*_1(s_1,s_2)=B_1^{-d'_1}b^*_{\goth U' 
 }(s_2)b^*_{\goth U' 
 }(s_2-1),$$ 
en ajoutant dans le cas orthogonal la condition:  $3k\≤2n-2$ pour le type $(\overline \∆,\overline \lambda_0)=(D_n,\alpha_k)$ et $3k\≤2n-1$ pour le type $(\overline \∆,\overline \lambda_0)=(B_n,\alpha_k).$

\end{prop}

\dem On  adapte une d\émonstration faite dans le cas commutatif (cf.lemme 5.6 dans  \cite{rubschiff}): soit
 $V$ l'espace affine suivant:
$$V=Y +\goth U^+(k)_{-1}=Y+E_{-2}(h_k)\cap \goth g_{-1}  \quad
(Y=\sum_{k+1\≤i\≤p}X_i^{-1}),$$ pour $f\in {\EuScript C}^{\infty}(\goth g_{-1}"),$ on d\éfinit $R_f\in {\EuScript C}^{\infty}( \goth U^+(k)_{-1}" )$ par $R_f=f(Y+.),$  il est facile de v\érifier que: 
$$(F_k(\partial)f)/_V=(F_k /_ {E_2
(h_k)\cap \goth
g_1})(\partial)(R_f)$$
La d\érivation \étant d\éfinie relativement \à la restriction de $\widetilde B$ \à $\goth U^+(k);$  ainsi
il suffit de calculer $F^{*s}/_V .$

\noi Or $$V=Y+\goth
U_{-1} +[Y,E_{-2}(h_k)\cap \goth
g_0]=exp(ad(E_{-2}(h_k)\cap
\goth g_0)(Y+\goth
U_{-1})\subset N_{\goth t}(V')  $$
avec
$V'=Y+\goth U_{-1}$ donc $F^{*s}/ _V=F^{*s}/ _{V'}$   d'o\ù, par orthogonalit\é des
diff\érents sous-espaces relativement \à $\widetilde B$, on a $$(F_k /_ {E_2(h_k)\cap \goth
g_1} )(\partial)(F^{*s}/ _V)=(F_k /_{\goth U_1 
 })(\partial)  (F^{*s}/ _{V'}),$$ 
 La d\érivation \étant d\éfinie relativement \à la restriction de $\widetilde B$ \à $\goth U.$
 
\noi En dehors des  formes r\éelles du
cas exceptionnel $(E_7,\alpha_6)$ et du cas classique $(C_n,\alpha_k)$ avec le parabolique $P'_0,$ les invariants relatifs
fondamentaux du \PV $(P(H_1,...,H_k),\goth U_1)$ sont donn\és pour $i=1,...,k$ par la restriction de $F_i$  \à $\goth U_1,$ que l'on note $G_i.$  

\noi Comme on a  pour
$u\in E_2(h_j)\cap \goth U_1$ et $v\in E_0(h_j)\cap \goth
U_1$ qui commutent (lemme 1.4.7) :
$$ F_p(X+u+v)=F_k(u+v)P_{p-k}(X)=
F_j(u)P_{p-j}(X+v)\ \ \
( X=\sum_{k+1\≤i\≤p}X_i),$$ 
on en d\éduit tous les invariants relatifs normalis\és ce qui donne  pour $j=1,...,k-1:$ 
$$P^{\goth U}_{k-j}(v)=\frac{P_{p-j}(X+v)}{P_{p-k}(X)}\  \hbox{
d'o\ù }\ 
G^*_j =\frac{F^*_{p-k+j}(Y+ .\ )}{F^*_{p-k}(Y)} \  \hbox{pour}\ j=1,...,k  \  \hbox{   donc}$$ $$\text{pour }y\in \goth U_{-1}\text{ on a }
F^*_j(Y+y)= \left \{
    \begin{array}{l}  F^*_j(Y)\  \hbox{pour }\ j=1...,p-k,\\ 
   G^*_{j-(p-k)}(y).F^*_{p-k}(Y)\ \hbox { pour  }\  j=p-k+1,...,p.\\  \end{array} \right.   $$ 
Ainsi il existe une
constante $C$ (explicite ) telle que: $ \prod_{1\≤i\≤p}F^{*s_i}(Y+y)=
C\prod_{1\≤i\≤k}G_i^{*s_{p-k+i}}(y)$ d'o\ù le r\ésultat en tenant compte du fait que $F_k(\partial)f)/\goth U=A_k^{-d_k}G_k(\partial)f.$  

\noi Pour les 2 cas restants on a $p=2$ et on
prend commes invariants relatifs fondamentaux $G
=F_2(X+.\ )$ donc
$G^*=F_2^*(Y+. \ )$  d'o\ù
$$ (F_1^{*s_1}F_2^{*s_2})(Y+.\
)= F_1^{*s_1}(Y)G^{*s_2}\
 ,  \ F_1= {G^2\over P_1(X)}=F_1^*(Y)G^2 \ 
 \hbox{
 donc}$$
 $$(F_1/_{\goth
U_1})(\partial)=A_1^{-d_1}F_1^*(Y)G(\partial)^2. $$
La d\émonstration est analogue pour $b^*_{p-k}.$\fdem

\bigskip

\subsection  {\bf   Fonctions Z\étas: d\éfinition et \équations
fonctionnelles abstraites}

\bigskip

Suivant F.Sato (\cite{sato3}) et J.I.Igusa (\cite{igusa5},\cite{igusa12}), au \PV $(P_{\goth t},\goth g_ 1)$ et \à son ``dual" sont associ\és des fonctions Z\étas locales dont on rappelle la d\éfinition. \\

\noi On note $\Omega(\F^{*p})$ le groupe des caract\ères continus de $\F^{*p}$ c'est \à dire que $\pi=(\omega ,s)\in \Omega(\F^{*p})$    si 
$\pi=(\pi_1,...,\pi_p),$\   
$ \pi_1,..., \pi_p$ \étant des caract\ères continus de
$\F^*,$ ils sont d\éfinis pour $i=1,...,p$ par
$\pi_i (x)=\omega_i(x) {|x|_{\F}}^{s_i},$
  $s_i\in
\C,$   $ \omega_i$ \étant un caract\ère unitaire de $\{x\in \F\ |\ |x|_{\F}=1\}$ dans le cas archim\édien et $\omega_i$
\étant un caract\ère unitaire de $\F^*$ dans le cas $\goth
p-$adique,
$\omega=(\omega_1,...,\omega_p)\in
\widehat{\F^{*p}}$ et $s=(s_1,...,s_p)\in \C^p.$ 

\noi Pour $x=(x_1,...,x_p)\in \F^{*p}$ on a $\pi (x)=\prod_{1\≤i\≤p}\pi_i(x_i),$ $\omega (x)=\prod_{1\≤i\≤p}\omega_i(x_i)$ et $|x|^s=\prod_{1\≤i\≤p} |x_i|_{\F}^{s_i}.$
\vskip 5pt
\noi On note $\Re(\pi )=\Re(s)= inf_{1\≤i\≤p}\Re(s_i).$\\

\noi Dans le cas r\éel (ou dans le cas $\goth p-$adique), pour $a\in
\F^*$ on note 
  $  \tilde\omega_a 
 $ le caract\ère (d'ordre $2$) de $\F^*$ (trivial sur
$\F^{*2}$) d\éfini par
$  \tilde\omega_a(x)=(x,a)$
 ($(.,.)$ symbole de Hilbert).\\
 
\begin{defi} Soit  $\pi=(\omega ,s)\in \Omega(\F^{*p}).$ Pour $f\in {\EuScript
S}(\goth g_1)$ et $g\in {\EuScript
S}(\goth g_{-1}),$ les fonctions Z\étas locales des \PVs $(P_{\goth t},\goth g_1)$ et $(P_{\goth t},\goth g_{- 1})$ sont d\éfinies pour toute orbite $O$   de
$P_{\goth t}$ dans
$\goth g"_1$  et toute orbite  $O^*$   de
$P_{\goth t}$ dans $\goth g"_{-1}$ par:
$$ 
Z_O(f; \pi)=Z_O(f;\omega, s)=\int_Of(x) 
\pi (F(x))\ dx  \ \ ,\ \
  Z^*_{O^*}(g;\pi)= Z^*_{O^*}(g;\omega, s)=\int_{O^*}g(y) 
 \pi(F^*(y))\ dy.$$
On pose  
 $$\ Z(f;\pi)=Z(f;\omega, s)=\int_{\goth g_1}f(x) 
\pi (F(x))dx\ \ ,\ \
 Z^*(g;\pi)=Z^*(g;\omega,
s)=\int_{\goth g_{-1}}g(y) 
\pi (F^*(y))dy\ .$$
Lorsque $\F=\R,$ on  d\éfinit de m\ême  les fonctions Z\étas associ\ées aux orbites de $(P_{\goth t})_{\R}$ dans  
$\goth g"_1$  et $\goth g"_{-1},$ avec $(P_{\goth t})_{\R}=(G_{\goth t})_{\R}.N_{\goth t},$ $(G_{\goth t})_{\R}$ \étant la composante connexe r\éelle de $G_{\goth t}.$
\end{defi}

\bigskip

\noi Notons par * l'involution de $ \Omega (\F^{*p})$
 et de $\C^p$ d\éfinie par :

$$\begin{array}{l}   \pi^*=(\pi_1,...,\pi_p)^*=(\pi_{p-1},...,
\pi_1, (\pi_1...\pi_p)^{-1}),\\ 
    s^*=(s_1,...,s_p)^*=(s_{p-1},...,s_1,
-(s_1+...+s_p).\\  \end{array}     $$
\à l'exception des
 $\F$-formes du cas exceptionnel $(E_7,\alpha_6)$ et du cas classique $(C_n,\alpha_k)$ avec le parabolique $P'_0$ pour lesquels   les
relations deviennent : $$
(\pi_1,\pi_2)^*=(\pi_1,(\pi_1^2\pi_2)^{-1}) \ \ ,\ \ s^*=(s_1,s_2)^*=(s_1,-(2s_1+s_2)).$$

\begin{theo} Equation fonctionnelle abstraite

\noi Pour $f\in S(\goth g_1)$ et $\omega\in  \widehat{\F^{*p}}$ les int\égrales $Z_O(f;\omega,s)$ et $Z^*_{O^*}({\EuScript F} (f);
\omega,s )$ convergent pour $Re( s)\≥0,$ se prolongent en
des fonctions m\éromorphes sur $\C^p$. 

\noi De plus il existe des fonctions m\éromorphes
$a_{O^*,O}(\omega,\ )$ et $a^*_{O,O^*}(\omega,\ )$ telles que
l'on ait pour toute orbite
$O^*$   de
$P_{\goth t}$ (resp.$(P_{\goth t})_{\R}$) dans $\goth g"_{-1}$  :
$$Z^*_{O^*}( {\EuScript F}(f)
;\omega,s) =\sum_{\hbox{orbites}\
O\ \hbox{de}\ P_{\goth t}\
\hbox{dans}\ \goth g"_1}a_{O^*,O}(\omega,s)Z_{O }( 
f;\omega^*,s^*-N1_p)  $$  
\noi et  pour toute orbite O  de
$P_{\goth t}$ (resp.$(P_{\goth t})_{\R}$) dans $\goth g"_{1}$:
$$Z_{O}( 
f;\omega,s) =\sum_{\hbox{orbites}\
O^*\ \hbox{de}\ P_{\goth t}\
\hbox{dans}\ \goth g"_{-1}}a^*_{O,O^*}(\omega,s)Z^*_{O^* }( 
  {\EuScript F}(f);\omega^*,s^*-N1_p)\quad   )\ .$$  

\end{theo}

\noi $(N=\frac{\hbox{ \tiny dim(} \goth g_1)}
 {\hbox{\tiny degr\é de } F})$\\

\dem Lorsque $p=1,$ on a simplement la fonction Z\éta du \PV:
$(\goth g_0,\goth g_1)$ et les r\ésultats sont dus \à \cite{sato-shintani} dans le cas archim\édien et \à J.I.Igusa \cite{igusa5} dans le cas $\goth p-$adique.\\

\noi Lorsque $p\≥2:$

\begin{enumerate}
\item Dans le cas $\goth p-$adique,on applique
 le th\éor\ème $k_{\goth p}$ de \cite{sato3}  et on utilise
notamment le th\éor\ème 2.1.1.

\item Dans le cas  complexe,  on proc\ède comme dans  \S 3 de 
\cite{boppruben1}; on utilise notamment les lemmes
3.4.1 et 3.4.2 qui impliquent pour $k=1,...,p:$
$$f\in S(\goth g_1)\ \ Z^*_{O^*}( \goth F(F_k f);
\omega,s )=(-2i\pi)^{d_k}b_k(s)Z^*_{O^*}(  {\EuScript F}(f);
\omega,s-1_p+1_{p-k} ) $$
avec l'am\énagement convenable lorsque $k=1$ pour les deux exceptions.
\item Dans le cas r\éel, on proc\ède comme dans le  $\S 5.2$ de  \cite{boppruben} soit avec le groupe $P_{\goth t}$ soit avec sa composante connexe, $(P_{\goth t})_{\R},$  et comme dans  le cas  complexe on utilise notamment les lemmes
3.4.1 et 3.4.2.

\noi On peut noter qu'il est inutile d'introduire  le caract\ère  $\tilde\omega_{-1}$ lorsqu'on consid\ère les fonctions z\étas associ\ées \à 
$(P_{\goth t})_{\R}.$  \end{enumerate}
\fdem\\

Les coefficients des \équations fonctionnelles sont ind\épendants des normalisations.  \\

 \noindent En effet soient $F'_1=a_1F_1,...,F'_p=a_pF_p,$ une autre famille d'invariants  relatifs fondamentaux du \PV $(P_\goth t,\g_1),$ on note les fonctions Z\étas associ\ées par $Z'$ et $Z'^*,$ $\four'$ la transformation de Fourier associ\ée et les coefficients de l'\équation fonctionnelle par $a'_{O^*,O}.$\\
 
 \noi Pour $u=(u_1,...,u_p)\in \bigl(\FF)^p$ on d\'efinit les ouverts
(\'eventuellement vides):\\

$O_u= \{x\in \goth
g_1\ |\ F_1(x) \F^{*2}=u_1\ ,\ F_2(x)  \F^{*2}=u_1u_2\ ,...,F_p(x)  \F^{*2}=u_1...u_p  \},$

$O^*_u= \{x\in
\goth g_{-1}\ |\ F^*_1(x) \F^{*2} = u_p\   ,\ F^*_2(x)  \F^{*2}= u_{p-1}u_p \ ,...,F^*_p(x)  \F^{*2}= u_p...u_1 \}.$ 
 \\
 
 \noi Ainsi que les analogues pour les invariants $F'_i,i=1,...,p,$ que l'on note $O'_u$ et $O'^*_u.$\\
 
 \noindent  On d\'efinit \'egalement  les fonctions Z\'etas associ\'ees, pour $f\in S(\goth
g_1)$ (resp.$h\in S(\goth g_{-1})$) :
$$  
Z_u(f;\omega)=Z(f{\bf 1}_{O_u};\omega) \quad (\hbox{resp.}
 Z^*_u(h;\omega)=Z^*(h{\bf 1}_{O^*_u};\omega)\  )\ .  $$
 $$  
Z'_u(f;\omega)=Z'(f{\bf 1}_{O'_u};\omega) \quad (\hbox{resp.}
 Z'^*_u(h;\omega)=Z'^*(h{\bf 1}_{O'^*_u};\omega)\  )\ .  $$

 \begin{lem} Soit $\pi\in \widehat{(\F^*)^p}$ alors:
 
 \begin{enumerate}
 \item $a'_{O^*,O}(\pi)=a_{O^*,O}(\pi).$
 
 \item On exclut les $\F-$formes du cas exceptionnel $(E_7,\alpha_6)$ ainsi que le cas classique $(C_n,\alpha_k)$ avec le parabolique $P'_0.$\\
i) On suppose que pour tout $u,v\in (\FF)^p$ il existe un coefficient $a_{v,u}(\pi)$ tel que pour tout $f\in S(\goth
g_1)$ on ait l'\équation fonctionnelle:
$$Z^*_v(\four (f); \pi)=\sum_{u\in  (\FF)^p}a_{v,u}(\pi)Z_u(f;\pi^*|\ |^{-N1_p})$$
alors on a  l'\équation fonctionnelle:
 $$Z'^*_v(\four' (f); \pi)=\sum_{u\in  (\FF)^p}a_{v',u'}(\pi)Z'_u(f;\pi^*|\ |^{-N1_p}),$$
$w=(w_1,...,w_p)\in \FF\rightarrow w'=(w'_1,...,w'_p)\in \FF$  d\éfini  par $w'_i=a_{i-1}a_iw_i .$\\

ii) On suppose que pour tout $\pi\in \widehat{(\F^*)^p}$ il existe un coefficient $A(\pi)$ tel que pour tout $f\in S(\goth
g_1)$ on ait l'\équation fonctionnelle:
$$Z^*(\four (f); \pi)= A(\pi)Z(f;\pi^*|\ |^{-N1_p})$$
alors on a les \équations fonctionnelles:\\

\noi a) $Z'^*(\four (f); \pi)= A(\pi)Z'(f;\pi^*|\ |^{-N1_p}).$\\

\noi b) $Z'^*_v(\four' (f); \pi)=\sum_{u\in  (\FF)^p}a'_{v,u}(\pi)Z'_u(f;\pi^*|\ |^{-N1_p})$ avec  $$a'_{v,u}(\pi)=\ds\frac{1 }{f^p}\sum_{(b_1,...,b_p)\in
\bigl(\FF\bigr)^p}\biggl(\prod_{1\≤i\≤p}  (b_i,u_pv_p...u_{p-i+1}v_{p-i+1})\biggr)\ A(\pi(\tilde\omega_{b_1},...,\tilde\omega_{b_p})).$$ 

 \end{enumerate}
 \end{lem}
 
 \dem  Soient $\goth B=(e_i)_{1\≤i\≤m}$ une base de $\g_1$ et $\goth B^*=(e_i^*)_{1\≤i\≤m}$ la base duale pour $\widetilde B$ alors pour $Re(\pi)>0$
 $$Z(f;\pi)=\lambda\int_{\F^m}f(\sum_{1\≤i\≤m}x_ie_i)\prod_{1\≤j\≤p}\pi_j(F_j(\sum_{1\≤i\≤m}x_ie_i))dx_1...dx_m$$
 avec $\lambda=|(det\ (\theta_{x_0}/\g_1)\ )_\goth B^{\goth B^*}|^{\frac{1}{2}} |F_p(x_0)|^N.$\\
 
\noi Soient $\goth B'=(e'_i)_{1\≤i\≤m}$ une autre base de $\g_1$ et $\goth B'^*=({e'}_i^*)_{1\≤i\≤m}$ la base duale pour $\widetilde B$ alors pour $Re(\pi)>0$
 $$Z'(f;\pi)=\lambda'\int_{\F^m}f(\sum_{1\≤i\≤m}x'_ie'_i)\prod_{1\≤j\≤p}\pi_j(F'_j(\sum_{1\≤i\≤m}x'_ie'_i))dx'_1...dx'_m$$
 avec $\lambda'=|a_p|^N |(det\ (\theta_{x_0}/\g_1)\ )_{\goth B'}^{\goth B'^*}|^{\frac{1}{2}} |F_p(x_0)|^N=|a_p|^N|det(P)|\lambda,$ $P$ \étant la matrice de passage de $\goth B$ \à  $\goth B' $ donc par prolongement m\éromorphe:
 $$Z'(f;\pi)=|a_p|^N \prod_{1\≤j\≤p}\pi_j(a_j)\ Z(f;\pi)\ ,\ Z'_O(f;\pi)=|a_p|^N \prod_{1\≤j\≤p}\pi_j(a_j)\ Z_O(f;\pi)\ ,$$
 $${Z'}_{u}(f;\pi)=|a_p|^N \prod_{1\≤j\≤p}\pi_j(a_j)\  {Z}_{u'}(f;\pi).$$
 Pour $f\in S(\goth
g_1)$ et $z\in \g_{-1}$ on a $\four'(f)(z)=|a_p|^N \four (f)(z).$\\

\noi Par le choix de la  normalisation des invariants relatifs fondamentaux  (cf.$\S 3.3$) on a pour $j=1,...,p,$ $F'^*_j=c_jF^*_j $ avec $c_j=\ds\frac{a_{p-j}}{a_p },$ donc par la normalisation des mesures du $\S 3.2$ on a pour $g\in S(\goth
g_{-1})$:
$$Z'^*(g;\pi)=|a_p|^{-N} \prod_{1\≤j\≤p}\pi_j(c_j)\ Z^* (g;\pi)\ ,\ Z'^*_{O^*} (g;\pi)=|a_p|^{-N} \prod_{1\≤j\≤p}\pi_j(c_j)\ Z^*_{O^*}(g;\pi)\ ,$$
$$\ Z'^*_{u}(g;\pi)=|a_p|^{-N} \prod_{1\≤j\≤p}\pi_j(c_j)\  Z^*_{u'}(g;\pi),$$
 d'o\ù le r\ésultat de 1. et 2 i),ii)a).\\
 
 \noi  Pour 2)ii),b) soit $f$ le cardinal de
$\bigl(\FF\bigr).$

Pour  $v\in \bigl(\FF\bigr)^p$,  
$h\in \EuScript S(\goth g_{-1})$ et
$Re( \pi)>0$ on a:
$$Z^*_v(h;\ \pi)=\ds\frac{1 }{f^p}\sum_{(b_1,...,b_p)\in
\bigl(\FF\bigr)^p}\biggl(\prod_{1\≤i\≤p}  (b_i,v_p...v_{p-i+1})\biggr)\
\ Z^*(h; \pi. (\tilde\omega_{b_1},...,\tilde\omega_{b_p}))\ .$$
 
\noi Par prolongement m\'eromorphe, cette \'egalit\'e est vraie pour
tout caract\`ere $ \pi$ donc, en appliquant l'\'equation
fonctionnelle abstraite   \`a $Z^*(\four( f); \pi. (\tilde\omega_{b_1},...,\tilde\omega_{b_p})),$  on obtient:
$$Z'^*_v(\four (f); \pi)=\sum_{u\in  (\FF)^p}a'_{v,u}(\pi)Z'_u(f;\pi^*|\ |^{-N1_p})$$
avec $a'_{v,u}(\pi)=\ds\frac{1 }{f^p}\sum_{(b_1,...,b_p)\in
\bigl(\FF\bigr)^p}\biggl(\prod_{1\≤i\≤p} (b_i,u_pv_p...u_{p-i+1}v_{p-i+1})\biggr)\ A(\pi(\tilde\omega_{b_1},...,\tilde\omega_{b_p})).$  \fdem
 
\begin{rema}
\end{rema}

\begin{enumerate}
 
\item Plus g\én\éralement, soit $\goth g'$ une alg\èbre de Lie isomorphe \à $\goth g,$ notons $T$ l'isomorphisme de $\goth g$  sur $\goth g',$ consid\érons  le \PV $(P(H'_1 ,...,H'_p),\goth g'_1)$ avec $ H'_i= T(H_i)$ pour $i=1,...,p)$ et $\tilde B_{\goth g}(H_0,H_0)=\tilde B_{\goth g'}(H'_0,H'_0),$ alors pour toute orbite $O$ de $P(H_1,...,H_p)$ dans $\goth g"_1$ et toute orbite $O^*$ de $P(H_1,...,H_p)$ dans $\goth g"_{-1}$ on a:
$$a^{\goth g'}_{T(O^*),T(O)}=a^{\goth g}_{O^*,O}.$$

\item Lorsque $p\≥2,$ et en excluant les 2 exceptions du lemme 3.5.7, l'\équation fonctionnelle admet une \écriture plus simple en utilisant les conventions usuelles (cf.\cite{farautkoranyi}, \cite{clerc}, \cite{hironaka1}) que l'on rappelle.

\noi On d\éfinit pour $\pi=( \pi_1,..., \pi_p)\in \Omega (\F^{*p}),$ $i( \pi)=( \pi_p,..., \pi_1)$ et $\tilde { \pi}=( \pi_1 \pi_2^{-1},..., \pi_{p-1} \pi_p^{-1}, \pi_p)$ ainsi que l'analogue sur $\C^p:$
$$\hbox{pour}\quad s=(s_1,...,s_p)\quad  \tilde s=(s_1-s_2,...,s_{p-1}-s_p,s_p)\quad ,\quad i(s)=(s_p,...,s_1),$$
alors $(\widetilde  \pi)^*=\widetilde{i( \pi)}^{-1},$  $(\widetilde  s)^*=-\widetilde{i( s)} $  et l'\équation fonctionnelle s'\écrit simplement lorsqu'on pose:
$${\widetilde Z_O}(\ ;  \pi)=Z_O(\ ,\widetilde  \pi)\ ,\  \widetilde Z^*(\ ;  \pi)=Z^*(\ ;\widetilde  \pi)\quad \hbox{ 
alors}$$
$$ {\widetilde Z}^*_{O^*}( 
 {\EuScript F}(f);  \pi ) =\sum_{\hbox{orbites}\
O\ \hbox{de}\ P_{\goth t}\
\hbox{dans}\ \goth g"_1}\widetilde a_{O^*,O}( \pi)\widetilde Z_{O }( 
f;  i(\pi)^{-1}|\ |^{-N1_p})  \ ,$$  
 $$\widetilde Z_{O}( 
f;\  \pi ) =\sum_{\hbox{orbites}\
O^*\ \hbox{de}\ P_{\goth t}\
\hbox{dans}\ \goth g"_{-1}}\widetilde a^*_{O,O^*}( \pi)\widetilde Z_{O^* }( 
  {\EuScript F}(f);  i(\pi)^{-1}|\ |^{-N1_p})   $$
$$\hbox{avec}\quad
\widetilde a_{O^*,O}(  \pi)=a_{O^*,O}( \widetilde \pi)\ \hbox{et}\ 
\widetilde a^*_{O,O^*} (\pi)= a^*_{O,O^*}(\widetilde  \pi).$$
De plus, par le choix de la normalisation des \irfs, on a  la formule habituelle:
 $$\forall x\in \goth g"_1\quad  \widetilde  \pi (F^*)(x^{-1})=\widetilde{ i( \pi)}^{-1}(F)(x).$$
\noi Dans le cas archim\édien, pour $s=(s_1,...,s_p)\in \C^p$ on pose comme dans \cite{farautkoranyi} :$$F_s=F^{\tilde s}=F_1^{s_1-s_2}...F_{p-1}^{s_{p-1}-s_p}.F_p^{s_p} \quad ,\quad (F^*)_s=(F^*)^{\tilde s}\quad \hbox{alors}\quad F^*_s(x^{-1})=F_{-i(s)}(x)$$
et pour  $m=(m_1,...,m_p)\in \N^p$ tel que $m_1\≥m_2\≥...\≥m_p\≥0=m_{p+1}$ soient $F_m(\partial)=(F_1(\partial))^{m_1-m_2}...(F_p(\partial))^{m_p}$ et $F^*_m(\partial)=(F^*_1(\partial))^{m_1-m_2}...(F^*_p(\partial))^{m_p},$ alors  on a  les  relations habituelles:
 $$ F^*_m(\partial)F_s=B_m(s)F_{s-i(m)}\  {et}\  F_m(\partial)F^*_s=B^*_m(s)F^*_{s-i(m)}$$
  $$\hbox{avec}\ B_m(s)=\prod_{\{k=1,...,p | m_k>m_{k+1}\}}[\prod_{j=m_{k+1}}^{m_k-1}{\tilde b}_k(s_1,...,s_p-j)\ ]\quad \hbox{ et}\quad  {\tilde b}_k(s )=b_k(\tilde s),$$
   $$\hbox{avec}\ B^*_m(s)=\prod_{\{k=1,...,p | m_k>m_{k+1}\}}[\prod_{j=m_{k+1}}^{m_k-1}{\tilde b^*}_k(s_1,...,s_p-j)\ ]\quad \hbox{ et}\quad  {\tilde b^*}_k(s )=b^*_k(\tilde s)$$
par application du lemme 3.4.2 et de la proposition 3.4.4. 
\end{enumerate}

\bigskip

\subsection{\bf  Deux exemples fondamentaux}

\bigskip
Le calcul des coefficients apparaissant dans l'\équation v\érifi\ée par les fonctions Z\étas associ\ées aux \PVs $(P_{\goth t},\goth g_1)$ (ou bien  $(G,\goth g_1)$) se fera par descente en se ramenant \à des \PVs de dimensions plus petites dont les \irfs ont des degr\és plus petits. \\

\noi Dans ce paragraphe on rappelle  2 exemples  utilis\és ult\érieurement,  ils correspondent aux cas o\ù l'\irf de $(G,\goth g_1)$  est de degr\é $1$ ou $2.$\\

\noi On rappelle \également que pour la fonction Z\éta de $(G,\g_1)$ (resp.$(G,\g_{-1})$  avec un seul \irf $F$ (resp. $F^*)$ on a simplement lorsque:\\

\noi$\omega \in \widehat \F$ et $s\in \C,$ $u\in \F^*/\F^{*2},$   $f\in {\EuScript S}(\g_1)$ et $g\in {\EuScript S}(\g_{-1})$ : 
 $$\begin{array}{lll}Z_u(f;\omega,s)&=Z(f1_{\{x\in \g_1|F(x)\F^{*2} =u };\omega,s)&=\ds\frac{1}{| \F^*/\F^{*2}|}\sum_{a\in  \F^*/\F^{*2}} \tilde\omega_a (u)Z(f; \tilde\omega_a .\omega,s)\ ,\\
 \\
 Z^*_u(g;\omega,s)&=Z(g1_{\{x\in \g_{-1}1|F^*(x) \F^{*2}=u\}};\omega,s)&=\ds\frac{1}{| \F^*/\F^{*2}|}\sum_{a\in  \F^*/\F^{*2}} \tilde\omega_a (u)Z^*(g; \tilde\omega_a .\omega,s)\ ,\end{array}$$
$| \F^*/\F^{*2}|=$ cardinal de $ \F^*/\F^{*2}.$  

  \bigskip
  
\subsubsection{\bf $\g_1$ est de dimension $1$}
 \bigskip

C'est l'exemple le plus simple de  \PV : celui associ\é \à l'action de $Gl_1(\F)$ sur $\F,$ c'est \à dire au cas o\ù $\goth g$ est une alg\èbre d\éploy\ée de rang $1:$
$$\goth g=\F X_{\alpha}\oplus \F h_{\alpha}\oplus \F X_{-\alpha}  $$
$(X_{\alpha},h_{\alpha},X_{-\alpha})$ \étant un \Sl, alors $H_0=\frac{1}{2}h_{\alpha}$ et, ainsi qui'il est dit dans le $\S 3.1,$ la normalisation choisie est :$$\tilde B_{\goth g}(X_{\alpha},X_{-\alpha})=1(=-2 \tilde B_{\goth g}(H_0,H_0))\ \hbox{afin que}\ \forall (x,y)\in \F\times \F\ \hbox{on ait}\ \tau( \tilde B_{\goth g}(xX_{\alpha},yX_{-\alpha}))=\tau(xy).$$
On rappelle que la mesure de Haar choisie sur $\F$ est autoduale pour $\EuScript F.$\\

\begin{theo} (Tate \cite{tate} p.319)
 Soit $f\in {\EuScript S}(\F),$ $Z(f;\omega,s)=\int_{\F}f(x) \omega(x)|x|^sdx$ admet un prolongement m\éromorphe \à $\C$ et satisfait \à l'\équation fonctionnelle:
$$ Z({\EuScript F (f)}; \omega,s)= \omega(-1) 
 \rho ( \omega,s+1)Z(f; \omega^{-1},-s-1). $$
 \end{theo}
 
 \bigskip
   
On rappelle que les  coefficients $\rho,$ plus pr\écis\ément $ \omega (-1)\rho
( \omega,s)= \omega (-1)\rho ( \omega|\ |^s)$ (appel\és \également $\Gamma_{\F}( \omega|\ |^s)$ (cf.\cite{sato3}) ou $b_{ \omega}(s)$ (cf.\cite{igusa3}) ou  $\Gamma (  \omega| \ |^s) $ dans le cas $\goth p$-adique (cf.\cite{sally}) sont donn\és par:\\

\begin{enumerate}
 \item Dans le cas complexe: pour 
  $m\in \Z$   soit $ \overline\omega_m(t)=(\frac{t}{\sqrt{|t|_{\C}}} )^m,$  on a :
 $$  \overline\omega_m(-1)\rho (   \overline\omega_m,s)=i^{|m|_{\R}}(2\pi)^{1-2s}\frac{\Gamma(s+\frac{|m|_{\R}}{2})} {
\Gamma(1-s+\frac {|m|_{\R}}{2})}.$$
 
\item Dans le cas r\éel on a :
$$\begin{array}{r c c}
\rho (  \tilde\omega_1,s)&=\pi^{-s+\frac{1}{2}}\ds\frac{\Gamma(\frac{s} {2})} {
\Gamma(\frac{1-s} {2})}  &=2.(2\pi)^{-s}\Gamma (s)\cos(\frac{\pi s} {
2}) \\
\\
   \tilde\omega_{-1}(-1)\rho(\tilde\omega_{-1},s)&= i\pi^{-s+\frac{1}{2}}\ds\frac{\Gamma(\frac{1+s} {2})}{ 
\Gamma(1-\frac{s}{ 2})}& =2 i.(2\pi)^{-s}
\Gamma(s)\sin(\frac{\pi s}{ 2})\ , \\
\end{array}$$
 ce que l'on peut \écrire en utilisant la convention $\sqrt 1=1$ et  $\sqrt {-1}=i:$
 $$  \tilde\omega_a(-1)\rho(   \tilde\omega_a,s)=2.(2\pi)^{-s}\Gamma (s)\sqrt{   \tilde\omega_a(-1)}\cos(\frac{\pi } {
2}(s+\frac{   \tilde\omega_a(-1)-1}{2})).$$
\item  Lorsque
$\F$ est un corps
$\goth P$-adique  de  caract\éristique r\ésiduelle $q,$  
on a 

\noi$\displaystyle{ \rho (  \tilde\omega_1,s)=q^{\nu (s-\frac{1}{2})}\frac{1-q^{s-1}}{1-q^{-s}}  }$
et
 $   
  \rho( \omega,s)=q^{\nu (s-\frac{1}{2})}C_{ \omega}q^{m( \omega)(s-1/2)}$ 
   si   $\omega$ est un caract\ère de 
$\goth O^*=\{x\in \F\ |\ |x|_{\F}=1\} $
  ramifi\é de degr\é  $ m( \omega)\≥1,$ prolong\é sur $\F^*$ en posant $\omega(\pi)=1;$
$C_{ \omega}$  est une somme de Gauss
$(C_{\omega^{-1}}.C_{\omega}=\omega (-1))$.\\
 
 \end{enumerate}
 
 \bigskip

 \begin{defi}\begin{enumerate}
 \item On pose $ \rho' (   \omega,s):=\omega(-1) \rho (   \omega,s).$
 
 \item 
  Pour $x\in  \F^*/\F^{*2},$  on d\éfinit: :
$$  \rho (   \omega,s;x)=\rho ( 
  \omega|\ |^s;x):=\frac{1}{  |\F^*/(\F^*)^2|}\sum_{y\in \F^*/(\F^*)^2}(x,y)\rho(  \omega \tilde\omega_y|\ |^s).$$
  \end{enumerate}
\end{defi}

\noi Dans les notations de F.Sato (\cite{sato3}) on a: $\rho' (   \omega,s)= \Gamma_{\F} (  \omega,s)$ et $ \rho (   \omega,s;x)=  \omega (-1) \Gamma_{\F} (  \omega,s;-x).$ \\

\noi Pour $b$ et $x\in  \F^*/\F^{*2},$  on a $\rho (   \omega \tilde\omega_b,s;x)=
 \tilde\omega_b(x)\rho (   \omega,s;x)$ et $\sum_{x\in \F^*/\F^{*2}}(x,\delta)\rho (  \omega,s;xb)=\rho ( \omega \tilde\omega_\delta,s)(b,\delta).$\\
 
\noi Dans le cas $\goth P-$adique, on a:
$ \rho (   \omega,s;x) =q^{m(\omega)(s-\frac{1}{2})} \rho (   \omega,\frac{1}{2};x)$ lorsque $m(\omega)>1.$\\

\noi  Alors du th\éor\ème  pr\éc\édent on d\éduit que:\\
 
 \begin{cor}
 Pour $v\in  \F^*/\F^{*2}$ et $f\in {\EuScript S}(\F),$   $Z_v(f;\omega,s)$ admet un prolongement m\éromorphe \à $\C$ et satisfait aux \équations fonctionnelles:
$$ Z_v({\EuScript F (f)};\omega,s)= \sum_{u\in  \F^*/\F^{*2}}
\omega (-1)  \rho (   \omega,s+1;-uv)Z_u(f;\omega^{-1},-s-1),$$
$$ Z_v( f;\omega,s)= \sum_{u\in  \F^*/\F^{*2}}
   \rho (   \omega,s+1;uv)Z_u(\four (f);\omega^{-1},-s-1).$$
Pour $v,w\in  \F^*/\F^{*2},$ on a l'\égalit\é:
$$ \sum_{u\in  \F^*/\F^{*2}}
   \rho (   \omega,s+1;u)\rho (   \omega^{-1},-s;-uvw)=\begin{cases} \omega (-1)\text {  si }v=w\\0\text{ sinon}.\end{cases}$$
 \end{cor}

 \subsubsection{\bf L'\irf est une forme quadratique}
 
 \bigskip
 
 \noi La fonction Z\éta associ\ée \à une forme quadratique non
d\ég\én\ér\ée a \ét\é \étudi\ée par Rallis et Schiffmann \cite{rallisschiffmann});  ils
ont, entre autre, \établi les \équations fonctionnelles. On
redonne celles associ\ées \à notre situation.  
\\
    
\noi  Rappelons que, suivant A.Weil, \à un caract\ère quadratique, $\tau\circ \phi,\phi$
\étant une forme quadratique non d\ég\én\ér\ée d\éfinie sur un espace
vectoriel de dimension $m,$ on associe un nombre complexe
$\gamma(\tau\circ \phi),$ de module $1.$ Avec le choix de $\tau$ on a:

\begin{enumerate}
\item $\F=\C:$ $\gamma(\tau\circ \phi)=1,$

\item $\F=\R:$ 
 $\gamma(\tau\circ \phi)=
 e^{i\frac{\pi}{4}(p-q)},$  $(p,q)$ \étant la signature de   $\phi,$
 
 \item $\F$ $\goth P$-adique:  $\gamma(\tau\circ \phi)=\alpha(1)^{m-1}\alpha(D)h,$
$D$ et $h$ \étant respectivement le discriminant et l'invariant de Hasse
de $\phi.$ 

\noi On peut noter que la formule  $\gamma(\tau\circ \phi)=\alpha(1)^{m-1}\alpha(D)h,$ est v\érifi\ée \également sur $\C$ et sur $\R.$
\end{enumerate}

\bigskip

\noi Pour $t\in \F^*,$ l'application
$\gamma_{\tau\circ \phi}(t)=\gamma( \tau\circ (t\phi))$ d\éfinit une application de 
$\F^* /\F^{*2}$ dans $\C^*.$ \\

\noi On note  $\alpha $ l'application
associ\ée \à $\phi(x)=x^2, x\in \F,$ et $\alpha_t=\ds\frac{1}{| \F^*/\F^{*2}|}\sum_{b\in \F^*/(\F^*)^2}\ \alpha(b)(t,b)$ pour $t\in \F^*/(\F^*)^2.$

\noi On rappelle que :
$$\alpha (1)\alpha (xy)=\alpha (x)\alpha (y)(x,y)\quad x,y\in \F^* /\F^{*2}\quad (\cite{rallisschiffmann}),$$
et on regroupe par commodit\é dans le lemme suivant quelques r\ésultats \él\émentaires utilis\és partiellement pour simplifier les coefficients de certaines\eqs.

\begin{lem}
A) Pour  $\pi \in \widehat{\F^*}$ on pose:
$$ h(\pi):=\sum_{t\in \F^*/(\F^*)^2}\ \alpha(t)\rho(\pi;t)\ (\text{resp. }h(s)=h(| \ |^s)\ ,\ s\in \C) $$
alors pour  $ u ,c,\delta\in  \F^*/(\F^*)^2 $ on a:
$$\ \sum_{t\in \F^*/(\F^*)^2}\ \alpha_t \  (t,u)\rho(\pi \tilde\omega_t)= \sum_{b\in \F^*/(\F^*)^2}\ \alpha (bu)\rho(\pi;b)= \ds\frac{\alpha (u)}{\alpha (1)} h(\pi \tilde\omega_u),$$
 $$\  \sum_{b\in \F^*/(\F^*)^2}\ \alpha (bu)(b,\delta)\rho(\pi;bc)= (c,\delta)\ds\frac{\alpha (uc)}{\alpha (1)} h(\pi \tilde\omega_{uc\delta}) .$$
 B)  Pour $a\in \{0,1\},$ $\pi_1,\pi_2\in \widehat{\F^*}$ et $u,v,\delta\in  \F^*/(\F^*)^2$ on pose:
$$ A_{  \pi_1,\pi_2}^a(u,v,\delta):=
\sum_{t\in \F^*/(\F^*)^2}\ (\delta,t)\ \alpha(t)^a\
\rho(  \pi_1;tu)  
\ \rho(  \pi_2;tv)(= A_{  \pi_2,\pi_1}^a(v,u,\delta))\ ,\ $$
$$\omega \in \widehat{\F^*}\quad A_{ \omega,s_1,s_2}^a =A_{\omega| \ |^{ s_1},\omega|\ |^{s_2}}^a(u,v,\delta)\ ,\ A_{  s_1,s_2}^a =A_{ | \ |^{ s_1}, |\ |^{s_2}}^a(u,v,\delta).$$
 \noi  On a les \égalit\és suivantes pour $b\in \F^*/(\F^*)^2:$  
 $$\begin{array}{lll}1)&A^a_{  \pi_1\tilde\omega_b, \pi_2\tilde\omega_b}(u,v,\delta)&=  (b ,uv) A^a_{  \pi_1
 ,\pi_2}(u,v,\delta)\ ,\\
 \\
2)& A^0_{  \pi_1,\pi_2}(bu,bv,\delta)&=   (b,\delta) A^0_{  \pi_1
  ,\pi_2}(u,v,\delta)\ ,\\
 \\
3)& \alpha (b\delta)A^1_{  \pi_1,\pi_2}(u,v,b\delta)&=\alpha   (\delta) A^1_{  \pi_1,\pi_2}(bu,bv,\delta)\ ,\\
\\
4)& \sum_{u\in \F^*/\F^{*2}}\tilde\omega_b(u)A^0_{  \pi_1,\pi_2}(u,v,\delta)&= (v,\delta b)  \rho (\pi_1\tilde\omega_b)\rho (\pi_2\tilde\omega_{\delta b})\ ,\\
\\
5)& \sum_{u\in \F^*/\F^{*2}}\tilde\omega_b(u)A^1_{  \pi_1,\pi_2}(u,v,\delta)&= (v,\delta b)\ds\frac{ \alpha (v)}{\alpha (1)} \rho (\pi_1\tilde\omega_b)h (\pi_2\tilde\omega_{\delta bv}) \\
6)&\ A^1_{   \pi_1,\pi_2}(u,v,\delta)&=(v,\delta )\ds\frac{ \alpha (v)}{\alpha (1)}\frac{1}{|\F^*/(\F^*)^2|} \sum_{b\in \F^*/\F^{*2}} (b,uv)  \rho(\pi_1\tilde\omega_b ) h( \pi_2\tilde\omega_{bv\delta})  .    \end{array}$$
 \end{lem}

\bigskip

\dem A) On utilise les d\éfinitions:$$\begin{array}{lll}A&=& \sum_{t\in \F^*/(\F^*)^2}\ \alpha_t\   (t,u)\rho(\pi\tilde\omega_t)=  \sum_{b\in \F^*/(\F^*)^2}\ \alpha (b)  \rho(\pi;ub) 
=  \sum_{c\in \F^*/(\F^*)^2}\ \alpha (uc)  \rho(\pi;c) \\
\\
&=& \frac{\alpha (u)}{\alpha (1)} \sum_{c\in \F^*/(\F^*)^2}\ \alpha (c)  (u,c)\rho(\pi;c) 
=\ds\frac{\alpha (u)}{\alpha (1)} \sum_{c\in \F^*/(\F^*)^2}\ \alpha (c) 
\rho(\pi\tilde\omega_u;c)\ 
=\ds\frac{\alpha (u)}{\alpha (1)}h(\pi\tilde\omega_u). \end{array}$$
$$\begin{array}{lll} B&=&\sum_{b\in \F^*/(\F^*)^2}\ \alpha (bu)(b,\delta)\rho(\pi;bc)\ =\sum_{b\in \F^*/(\F^*)^2}\ \alpha (bcu)(bc,\delta)\rho(\pi;b)   \\
&=&\ds\frac{\alpha(1)}{\alpha (\delta)}(u,\delta)\sum_{b\in \F^*/(\F^*)^2}\ \alpha (bcu\delta) \rho(\pi;b)=\ds\frac{\alpha(uc\delta)}{\alpha (\delta)}(u,\delta)h(\pi \tilde\omega_{uc\delta}) 
=(c,\delta)\ds\frac{\alpha (uc)}{\alpha (1)} h(\pi \tilde\omega_{uc\delta}) . \end{array}$$

\noi B) 1)2)3) sont imm\édiats.\\

\noi Pour 4) et 5) on rempace $A^a$ par sa valeur d'o\ù:
$${S_b}(a)=\sum_{u\in \F^*/\F^{*2}}\omega_b(u)A^a_{  \pi_1,\pi_2}(u,v,\delta)= \rho (\pi_1\tilde\omega_b)\sum_{t\in \F^*/\F^{*2}}(t,\delta b)\alpha (t)^a \rho (\pi_2;tv)$$
$$= \rho (\pi_1\tilde\omega_b)(v,\delta b)\sum_{t\in \F^*/\F^{*2}}\alpha (t)^a\rho (\pi_2\tilde\omega_{\delta b};tv)= \rho (\pi_1\tilde\omega_b)(v,\delta b)\sum_{t\in \F^*/\F^{*2}}\alpha (vt)^a\rho (\pi_2\tilde\omega_{\delta b};t)$$
or
$$\sum_{t\in \F^*/\F^{*2}}\alpha (vt)^a\rho (\pi_2\tilde\omega_{\delta b};t)=\begin{cases}\rho (\pi_2\tilde\omega_{\delta b})\text{ si }a=0\\
\\
\frac{\alpha(v)}{\alpha (1)}\sum_{t\in \F^*/\F^{*2}}\alpha (t)\rho (\pi_2\tilde\omega_{\delta vb};t)=\ds\frac{\alpha(v)}{\alpha (1)}h(\pi_2\tilde\omega_{\delta vb}) \text{ si }a=1.\end{cases}$$

\noi   6) se d\éduit de 5) puisque $\ds\frac{1}{|\F^*/(\F^*)^2|}S_b(1)$ est le coefficient de la fonction $f(\ )= A^1_{  \pi_1,\pi_2}(\ ,v,\delta)$ dans la base $\{\tilde\omega_b,b\in \F^*/(\F^*)^2\}$ des fonctions d\éfinies sur $\F^*/(\F^*)^2.$ \fdem

\bigskip

 \begin{theo}(Rallis-Schiffmann)
 On suppose que $\goth g$ est simple, que $(\goth
g_0,\goth g_1)$ est 1-irr\éductible, que $2H_0$ est 1-simple
et que l'invariant relatif fondamental, F, est une 
 forme quadratique de discriminant D. 

\noi Soient $N= \frac{1}{2}dim(\goth
g_1),$ 
$C= \tilde B(H_0,H_0),$  $\gamma_{\tau\circ F}= 
\sum_{a\in \F^* / \F^{*2}} \beta_a(F)  \tilde\omega_a$  et $\delta=(-1)^{[N]}D.$
$$1)\frac{ Z^*( {\EuScript F}(f)
;  \omega,s)} {  \omega(C)^{-2}|C|_{\F}^{-2s-N}}= \rho
(  \omega,s+1)\sum_{a\in \F^*/ \F^{*2}}\overline{\beta_a(F)}  (a ,-1) \rho
(  \omega \tilde\omega_a,s+N)Z(f;  \tilde\omega_a  \omega^{-1},-s-N).$$  
2) Lorsque
$N$ est un entier, on a :
 $$\frac{Z^*(  {\EuScript F}(f);  \omega,s)}{   \omega
(C)^{-2}|C|_{\F}^{-2s-N}}=  \rho (  \omega,s+1) \gamma(\tau\circ F)  \rho
(  \omega \tilde\omega_{\delta}
 ,s+N)Z(f;\  \omega^{-1}  \tilde\omega_{\delta  }
,-s-N).$$
3) Lorsque
$N$ n'est pas un entier, on a :
$$\frac{Z^*(  {\EuScript F}(f);  \omega,s)}{    \omega
(C)^{-2}|C|_{\F}^{-2s-N} }=\frac{ \alpha (-1)}{ \alpha (\delta)} {\gamma(\tau\circ F)} \rho(  \omega,s+1) \sum_{u\in    F(\g_1)^*/\F^{*2}} \alpha(\delta u)  h (  \omega \tilde\omega_{\delta u}|\ |^{
 s+N})Z_u(f;\  \omega^{-1}   
,-s-N).$$
4) Lorsque $\F=\R$ et F est  anisotrope,
 $Z^*( {\EuScript F}(f); s)=A_N( s)Z(f; -s-N),$ avec
$$A_N( s)= -2|C|_{\R}^{-2s-N}(2\pi)^{-2s-N-1}\Gamma(s+1)\Gamma(s+N)\sin(\pi s )  .$$
5) Pour      
$v\in   F(\g_1)^*/\F^{*2}$  on a:
$$ Z_v^*( {\EuScript F}(f)
;  \omega,s) =  \omega
(C)^{-2}|C|_{\F}^{-N-2s} \  \sum_{u\in    F(\g_1)^*/\F^{*2}}   a^{(F)}_{v,u}(  \omega,s) Z_u(f;\  \omega^{-1}  
,-s-N),$$
$$ Z_v(  f
;  \omega,s) =  \omega
(C)^{2}|C|_{\F}^{N+2s} \  \sum_{u\in    F(\g_1)^*/\F^{*2}}   a^{(F)}_{v,u}(  \omega,s) Z^*_u(\four(f);\  \omega^{-1}  
,-s-N)\ \text{avec}$$
 $$\begin{array}{ll}a^{(F)}_{v,u}(  \omega,s)&=\sum_{t\in  \F/\F^{*2}}\gamma _{\tau\circ F}(t)\rho(\omega,s+1;tv)\rho(\omega,s+N;tu)\\
 \\
 &=  \alpha(-1)^{a_N}  \gamma(\tau\circ F)  A^{a_N}_{  \omega,s+1,s+N}(v,u,\delta)  \
 ,\ a_N=\begin{cases} 0\ \text{si }N\in \N,\\
1\ \text{si }N\notin \N.\end{cases}
\end{array}$$

6) Pour      
$u,v\in   F(\g_1)^*/\F^{*2}$  on a:
$$\sum_{w\in    F(\g_1)^*/\F^{*2}} a^{(F)}_{v,w}(  \omega,s) a^{(F)}_{w,u}(  \omega^{-1},-s-N)=\delta_{u,v},$$
$$ \sum_{w\in    F(\g_1)^*/\F^{*2}} A^{a_N}_{  \omega,s+1,s+N}(v,w,\delta)A^{a_N}_{  \omega^{-1},-s-N+1,-s }(w,u,\delta)= \alpha(-1)^{2a_N}\gamma(\tau\circ F)^2\delta_{u,v}.$$

 \end{theo}

\dem 1) Elle consiste simplement \à exprimer les r\ésultats de \cite{rallisschiffmann}
dans nos notations. 
 
\noi Soit $E$ un espace vectoriel de dimension $n$, $F$ une forme quadratique d\éfinie sur $E$ non d\ég\én\ér\ée,
$$\beta (x,y)=F(x+y)-F(x)-F(y)\ ,\ x,y\in E$$
la forme bilin\éaire sym\étrique associ\ée, induisant un isomorphisme lin\éaire entre $E$ et $E^*$ d\éfini par:
$$f\in {\EuScript S }(E)\ :\ \hat f(x)=\int_Ef(y)\tau(\beta (x,y))d_Fy\ ,\ x\in E\ ,$$
$d_Fx$ \étant la mesure de Haar de $E$ autoduale correspondante:
$$\int_E\hat f(y)\tau(-\beta (x,y))d_Fy=f(x) ,\ x\in E\ .$$
Si on pose pour $  \omega|\ |^s$ caract\ère continu de $\F^*:$
$$f\in {\EuScript S }(E)\ :\ Z'_f(  \omega,s)=\int_Ef(x)  \omega(F(x))|F(x)|^{s-\frac{n}{2}}d_Fx$$
on a le th\éor\ème 2.13 p.521 de  \cite{rallisschiffmann}:
$$Z'_f(  \omega,s)=\rho(  \omega,s-\frac{n}{2}+1)\sum_{a\in \F^*/ \F^{*2}}\overline{\beta_a(F)}  \tilde\omega_a(-1)\rho
(  \omega \tilde\omega_a,s)Z'_{\hat f}(  \tilde\omega_a  \omega^{-1},-s+\frac{n}{2}).$$
Ici $E=\goth g_1,$ $dx$ \étant la mesure de Haar de $\goth g_1$ normalis\ée par $F$ au sens du \S3.2 il existe une constante $\mu$ telle que $d_{F}x=\mu dx.$

\noi Comme $2H_0$ est $1$-simple, $F$ est non d\ég\én\ér\ée  et on note $\Psi$ l'isomorphisme de $\goth g_1$ dans 
 $\goth g_{-1}$ 
 d\éfini par 
$$x,y\in \goth g_1\quad \tilde
B(x,\Psi(y))=\beta(x,y) \ ,$$
alors pour
$f\in {\EuScript S}(\goth g_1)$  et $x\in \goth g_1$ on a

$$  \hat f(x)=\mu( {\EuScript F}f)(\Psi(x))\ ,\   f(x)=\mu\int_{\goth g_1}({\EuScript F}f)(\Psi(y)) \tau (\tilde B(-x,\Psi(y)))d_Fy.$$

\noi On reprend les notations et les r\ésultats de la \demo du lemme 3.2.1.

\noi Soit $\goth B$ une base de  $\goth g_1,$ $\goth B^*$ la base duale dans $\goth g_{-1}$ pour $\tilde
B,$  comme:
$$ \int_{\goth g_1}g(\Psi(x))d_Fx =\mu \lambda^2 |det(\Psi)|_{\F}^{-1}
   \int_{\goth
g_{-1}}g(y)dy \quad \hbox{avec}\quad dx=\lambda d_{\goth B}x,$$
 on a:
$$f=\mu^2 \lambda^2|det(\Psi)|_{\F}^{-1}\overline{\EuScript F}({\EuScript F}f)=\mu^2 \lambda^2|det(\Psi)|_{\F}^{-1}f\ \hbox{ 
d'o\ù}\  \mu=\frac{\sqrt{|det( \Psi)|_{\F}}} { \lambda}$$  
et pour $g\in S(\goth g_{-1})$ on a
$ \int_{\goth
g_{-1}}g(y)dy=\mu\ \int_{\goth g_1}g(\Psi(x))d_Fx 
   .$
   
\noi Comme $\forall x\in \goth g_1$ et $g\in G$ on a $\Psi(gx)=\chi(g)\ g.\Psi(x),$ le polynome $F^*\circ\Psi$ est
relativement invariant par $G$ donc il existe une constante
non nulle, not\ée $\alpha,$ telle que
$F^*\circ\Psi=\alpha F.$ Il est facile de v\érifier que $\alpha$  est ind\épendant de la normalisation
de
$F$ et par cons\équent de la $\F$-forme $\goth g$ 
consid\ér\ée.  

\noi Soit $x_0\in \goth g'_1,$ on a pour tout $x\in \goth g_1,$
$\frac{1} 
{\alpha}F^*(\Psi(x))=F(x)=F^*(\theta_{x_0}(x)).F(x_0)^2  $ par (3) du \S3.2 donc les $2$ formes quadratiques 
$F(x_0)^2.F^*\circ\theta_{x_0}$ et  $\frac{1} 
{\alpha}F^*\circ\Psi$ sont \égales d'o\ù elles ont m\ême discriminant donc $$\alpha^{-2N}det(\Psi)^2=det(\theta_{x_0}/\goth g_1)^2F(x_0)^{4N} \quad \hbox{or}\quad   
 \lambda=|det(\theta_{x_0}/\goth g_1)|_{\F}^{\frac{1}{2}}|F(x_0)|_{\F}^N \quad \hbox{ d'o\ù }\quad \mu=|\alpha|_{\F}^{\frac{N} {2}}.$$

\noi Il reste \à d\éterminer $\alpha.$  Comme   $[ \Psi (x),x]$  
est dans le centre de $\goth g_0$ (prop.1.1.10  \cite{rubthese}), il est proportionnel \à
$H_0$ lorsque $\goth g_1$ est un $\goth g_0$ module
absolument simple. On v\érifie que ceci est \également vrai dans  
 l'unique cas restant puisqu'on peut supposer que 
l'alg\èbre $\goth g$ est d\éploy\ée et qui est, vu les
hypoth\èses, de type
$(A_n,\{\alpha_1,\alpha_n\})$ ( \cite{rubthese},propositions 3.3.7 et
3.3.8).

\noi Ainsi  $\Psi(x)=-\frac{1 }{ C}F(x)\Phi(x)$ pour $x\in
\goth g'_1$ donc $F^*(\Psi(x))={\displaystyle
\frac{1}{ C^2}F(x)}$  et $\alpha=\frac{1}{C^2}.$

\noi On applique le th\éor\ème 2.13 de  \cite{rallisschiffmann} en notant que pour 
 $f\in S(\goth g_1),$ on a:
  $$Z^*( {\EuScript F}(f);   \omega,s )= \omega (C)^{-2}|C|_{\F}^{-2s} Z'_{\hat f}(  \omega,s+ N)\quad
\hbox{et}\quad Z'_f(   \omega,s)= )| C|_{\F}^{ -N}\ Z( f;\  \omega,s-  N).  $$
2) De l'\égalit\é:
 $$(*)\quad \gamma_{\tau\circ F}(t)=
\gamma(\tau\circ F)  \tilde\omega_{\delta}(t)\ \Bigl(\alpha(t) 
\alpha(-1)\Bigr)^{a_N}\quad (\cite{rallisschiffmann}) ,\ $$
on d\éduit que:
$$\beta_a(F)=\begin{cases}0\text{ si }a\not=\delta\\
\gamma(\tau\circ F)\text{ si }a=\delta,\end{cases}\text{ lorsque }a_N=0$$
et pour $a\in \F^*:$
$$\beta_a(F)=\alpha(-1)\gamma(\tau\circ F)\frac{1}{  |\F^*/(\F^*)^2|}\sum_{y\in \F^*/(\F^*)^2}(y,a\delta)\alpha (y)=\alpha(-1)\gamma(\tau\circ F)\alpha_{a\delta}\text{ pour }a_N=1.$$
Comme $\overline{\gamma(\tau\circ F)}=\gamma_{\tau\circ F} (-1)=\gamma(\tau\circ F)(\delta,-1)(\alpha(-1))^{2a_N},$ on a:
$$\overline{\beta_a(F)}=\begin{cases}0\text{ si }a\not=\delta\\
\gamma(\tau\circ F)(\delta,-1)\text{ si }a=\delta,\end{cases}\text{ lorsque }a_N=0$$
d'o\ù 2) et lorsque $a_N=1:$
$$ \overline{\beta_a(F)}=\alpha(-1)\gamma(\tau\circ F)(a,-1)\frac{1}{  |\F^*/(\F^*)^2|}\sum_{y\in \F^*/(\F^*)^2}(y,a\delta)\alpha (y) $$
d'o\ù 3) en utilisant le A) du lemme 3.6.4.

5) On a: $$ a^{(F)}_{v,u}(  \omega,s)=
\ds\frac{1}{|\FF|}\sum_{a,b\in ( \F/\F^{*2})^2}(a,uv)(b,-u)\ \overline{\beta_b(F)}\ \rho(\omega\tilde\omega_a|\ |^{s+1}) \rho(\omega\tilde\omega_{ab}|\ |^{s+N})$$
or $\beta_b(F)=\ds\frac{1}{|\FF|}\sum_{t\in \FF}(b,t)\gamma_{\tau\circ F}(t)$ donc $ \overline{\beta_b(F)}=\ds\frac{(b,-1)}{|\FF|}\sum_{t\in \FF}(b,t)\gamma_{\tau\circ F}(t)$ d'o\ù la 1\ère \égalit\é et la seconde r\ésulte de la relation (*) appliqu\ée \à la 1\ère \égalit\é.\fdem
 \bigskip
\begin{rema}
1) Lorsque $\F=\R$ on a $b_{\goth g}(s)=b^*_{\goth g}(s)=\displaystyle{\frac{ 
 1} {C^2}}s(s+N-1).$
 \end{rema}
 
 \noi En effet il suffit de consid\érer une base $\goth B=\{e_i,1\≤i\≤2N\}$ de $\goth g_1$ dans laquelle
$$F(\sum_{1\≤i\≤2N}x_ie_i)=
\sum_{1\≤i\≤2N}a_ix_i^2;$$
dans la base de $\goth g_{-1},$
duale pour $\tilde B,$ et not\ée $(e^*_i)_{1\≤i\≤2N},$ on aura:
$$\Psi(\sum_{ 1\≤i\≤2N}x_ie_i)=2\sum_{ 1\≤i\≤2N}a_ix_ie^*_i  \quad
\hbox{d'o\ù}\quad F^*(\sum_{ 1\≤i\≤2N}y_ie^*_i)=\frac{1}{4C^2}
\sum_{ 1\≤i\≤2N} \frac{y_i^2}{ a_i}$$ par le calcul pr\éc\édent.\fdem\\
 
   \bigskip
 2) Rappelons que le coefficient  
 est ind\épendant de la normalisation de $F$ c'est \à dire que pour $u $
et $v$ \elts de $F(\goth g_1)^*/(\F^*)^2$ et $t\in \F^*$  on a: $a^{(tF)}_{tv,tu}=a^{(F)}_{v,u}.$

\bigskip

3) i) Dans le cas r\éel, on a pour $x=\pm 1:$ 
$\rho (  \tilde\omega_{\pm 1},s;x)= \tilde\omega_{\pm 1}(x)(2\pi)^{-s}\Gamma (s)e^{-ix\frac{\pi}{2}s}$
et pour $u,v=\pm 1:$ $m\in \N:$
$$A_{
s_1,s_2}^a(u,v,(-1)^m)=2.(2\pi)^{-(s_1+s_2)}\Gamma
(s_1)\Gamma (s_2)\ i^m\cos {\pi\over 2}(s_1u+s_2v+m-{a\over
2}).$$  ii) Dans le cas
$\goth p-$adique, lorsque la caract\éristique r\ésiduelle est
diff\érente de $2$, on note $\epsilon$ un \elt de $\goth O^*-(\goth O^*)^2,$ alors
$\{1,\epsilon,\pi,\epsilon\pi\}$ est un ensemble de
repr\ésentants dans $\F^*$ de $\F^*/(\F^*)^2.$
  
\noi Soit $C_0$ la constante
d\éfinie par:
$$C_0=\int_{|x|=q}\tau (x)\ \chi_0(x)|x|^{-\frac{1} { 2}}dx$$
$\chi_0$ \étant l'unique caract\ère non trivial de $\goth
O^*$ tel que $m(\chi_0)=1,$ on a $C_0^2=(\pi,-1)=(-\frac{1 }{ q})$ (symbole de Jacobi) et 
$$\Gamma(\chi_0|\ |^s)=C_0q^{s-\frac{1}{2}}=\chi_0(-1)\rho (\chi_0|\ |^s)=(\pi,-1)\rho (\chi_0|\ |^s) \ \ 
\hbox{( \cite{sally})}.$$ 
Soient $s_0\in \C$ d\éfini par $|\pi|^{s_0}=(\pi,-1)$ et $s_1\in \C$ d\éfini par $|\pi|^{s_1}= -1,$ on a:
$$\rho (  \tilde\omega_{\pi}|\ |^s)=\rho (\chi_0|\ |^{s+s_0})=C_0q^{s-\frac{1}{2}}\ ,\ \rho  ( \tilde\omega_{\epsilon\pi}|\ |^s)=\rho (  \tilde\omega_{\pi}|\ |^{s+s_1})=-C_0q^{s-\frac{1}{2}}\ ,$$
$$\rho (  \tilde\omega_{\epsilon} |\ |^s)=\rho ( |\ |^{s+s_1})=\frac{1+q^{s-1}}{1+q^{-s}}\  \hbox{ et}\  \alpha (1)=\alpha (\epsilon)=1,\alpha (\pi)=C_0=-\alpha (\epsilon\pi)\ \hbox{(prop.4-1 p.537 de \cite{rallisschiffmann}) }$$ 
d'o\ù $\alpha_1=\alpha_\epsilon=\frac{1}{2}$ et $\alpha_\pi= -\alpha_{\epsilon\pi}=\frac{1}{2}C_0(\pi,-1).$\\

\noi 

\noi Pour $s\in \C$ soit $f(s)=q^{-\frac{1}{2}}(q^s-q^{-s}),$  alors:\\
 \begin{lem}
 $\F$ est un corps $\goth p$-adique de
caract\éristique r\ésiduelle diff\érente de 2 et $\tau $ est
d'ordre $0$ alors:

A) Pour $u\in \F^*/\F^{*2}$ et $s\in \C$ on a:
$$ \alpha (u)  h(\tilde\omega_u|\ |^s )=\ds\frac{1+  \tilde\omega_{\epsilon}(u)}{ 2}\ [(u,\pi)q^{s-\ds\frac{1}{2}}+\frac{1-q^{-1}}{1-q^{-2s}}\ ]-\ds\frac{1-  \tilde\omega_{\epsilon}(u)}{ 2}.\ds\frac{q^{s-1}-q^{-s}}{1-q^{-2s}}.$$
 
B) Pour $a\in \{0,1\}$ et $u,v,\delta\in \F^*/\F^{*2}$ on a:
$$ 
2(1-q^{-2s_1})(1-q^{-2s_2})A_{s_1,s_2}^a(u,v,\delta)$$
$$= \sum_{\epsilon_1,\epsilon_2,
\epsilon_3\in
\{-1,1\}}P_{s_1,s_2}^{a,\epsilon_1,\epsilon_2,\epsilon_3} (u,v,\delta)
\frac{1+\epsilon_1\  \tilde\omega_{\epsilon}(u)}{ 2}.\frac{1+\epsilon_2
  \tilde\omega_{\epsilon}(v) }{ 2}.\frac{1+\epsilon_3  \tilde\omega_{\epsilon}(\delta) }{ 2} $$
avec:
\begin{enumerate}
  \item  $P_{s_1,s_2}^{0,1,1,1}(1,x,y)=(1-q^{-1})^2\ +\   \tilde\omega_{\pi}( 
 y) f(s_1-\frac{1}{2})f(s_2-\frac{1}{2}) \ + 
  \tilde\omega_{\pi}(- xy) f(s_1)f(s_2),$

  \item $P_{s_1,s_2}^{0,1,-1,1}(1,x,y)=-(1-q^{-1}).\Bigl(f(s_2-\frac{1}{2})+
   \tilde\omega_{\pi}( 
 y) f(s_1-\frac{1}{2})\Bigr),$ 
 
 \item  $P_{s_1,s_2}^{0,1,1,-1}(1,x,y) =\ -C_0   \tilde\omega_{\pi}( 
- y)\ . \Bigl(  f(s_1)f(s_2-\frac{1}{2})+
  \tilde\omega_{\pi}( 
 x)f(s_2)f(s_1-\frac{1}{2}) \Bigr) ,$
 
 \item $P_{s_1,s_2}^{0,1,-1,-1}(1,x,y) =\ C_0  (1-q^{-1}). \Bigl(f(s_2 )
   \tilde\omega_{\pi}( x) +  \tilde\omega_{\pi}(- y)f(s_1)\Bigr),$ 
 
  \item  $P_{s_1,s_2}^{1,1,1,1}(x,y,1)=(1-q^{-1})^2-  \Bigl(   \tilde\omega_{\pi}(x)f(s_1)f(s_2-\frac{1}{2})+
    \tilde\omega_{\pi}(y)f(s_1-\frac{1}{2})f(s_2)\Bigr),$

\item $P_{s_1,s_2}^{1,-1,-1,1}(x,y,1)=f(s_1-\frac{1}{2})f(s_2-\frac{1}{2}) +  \tilde\omega_{\pi}(- xy) f(s_1)f(s_2)   ,$ 

 \item $P_{s_1,s_2}^{1,1,-1,1}(x,y,1)=(1-q^{-1}).\Bigl(  \tilde\omega_{\pi}( 
 x) f(s_1)-f(s_2-\frac{1}{2}) \Bigr),$

\item $P_{s_1,s_2}^{1,-1,1,1}(x,y,1)=(1-q^{-1}).\Bigl(  \tilde\omega_{\pi}( 
 y) f(s_2)-f(s_1-\frac{1}{2}) \Bigr).$

 \end{enumerate}
 
  \end{lem}
 
  \dem
   
\noi  Elle r\ésulte d'un simple calcul et des r\ésultats de 3.6.4.\fdem

 \bigskip
 
  \noi Pour $s\in \C,$ posons $f_+(s)=q^{-\frac{1}{2}}(q^s+q^{-s})$ et 
  $$  {A^1}_{s_1,s_2}(x,y,1)=2(1-q^{-2s_1})(1-q^{-2s_2})P(x,y),
  $$Le calcul explicite donne donc :

  \begin{enumerate}
  
  \item Lorsque $x,y\in\{1,\epsilon\}:$
  
 $$ P(x,y)=\begin{cases}  (1-q^{- 1})^2-2q^{-\frac{1}{2}}f_+(s_1+s_2-\frac{1}{2})(\pi,x)+(\pi,x)(1+q^{-1})f_+(s_2-s_1)\text{ lorsque }x=y,\\
 \\
 (1-q^{- 1})^2 +(\pi,y)(1-q^{-1})f(s_2-s_1)\text{ lorsque }x\not=y.\end{cases}$$
 \item Lorsque $x\in \{\pi,\epsilon \pi\}$ et $y\in\{1,\epsilon\}:$
$$ P(x,y)=(1-q^{- 1})\biggl((\pi,y)f(s_2)-f(s_1-\frac{1}{2})\biggr).$$
\item Lorsque $x\in\{1,\epsilon\}$ et $y\in \{\pi,\epsilon \pi\}:$  
$$ P(x,y)=(1-q^{- 1})\biggl((\pi,x)f(s_1)-f(s_2-\frac{1}{2})\biggr).$$
\item Lorsque $x,y\in \{\pi,\epsilon \pi\}:$
$ P(x,y)=$
$$-q^{-\frac{1}{2}}f_+(s_2-s_1) \biggl(1+(\pi,-xy)\biggr)+ \biggl(1+(\pi,-xy)q^{-1}\biggr) q^{-\frac{1}{2}}\biggl((\pi,-xy)q^{s_1+s_2-\frac{1}{2}} +  q^{-(s_1+s_2-\frac{1}{2})}\biggr).$$

 \end{enumerate}

\bigskip

\noi Notons \également que :\\

\begin{lem} Dans le cas r\éel ou  $\goth p$-adique de caract\éristique r\ésiduelle diff\érente de $2$ avec $\tau $ d'ordre $0,$ on a:
\begin{enumerate}
\item
$A^1_{s,s+\frac{1}{2}}(1,x,1)=0 \text{  pour }x\not=1$
  et \\
  
  \noi $A^1_{s,s+\frac{1}{2}}(1,1,1)=|2|_\F^{-2s+\frac{1}{2}}\rho (|\ |^{2s})= h(\tilde\omega_u|\ |^s  )\rho (\tilde\omega_u|\ |^{s + \frac{1}{2}} ) \text{  pour }u\in \FF.$\\
  
\item $ h(\tilde\omega_u|\ |^s  )=|2|_\F^{-2s+\frac{1}{2}}\tilde\omega_u(-1)\rho (\tilde\omega_u|\ |^{-s  +\frac{1}{2}} ) \rho (|\ |^{2s})$ pour $u\in \FF.$
\end{enumerate}
\end{lem}

\dem 1)  R\ésulte d'un simple calcul \à l'aide de la formule explicite dans le cas r\éel (et on utilise la formule de duplication de Legendre) et dans le cas $\goth p$-adique de caract\éristique r\ésiduelle diff\érente de $2$ (relation B)5, du lemme 3.6.7). 

\noi  On v\érifie que la fonction d\éfinie par:  
$$x\in \F^*/\F^{*2} \quad f(x)=h(\tilde\omega_x |\ |^s )\rho(\tilde\omega_x |\ |^{s  +\frac{1}{2}} )    $$
est constante, soit en effectuant le calcul directement, soit par la relation 5)B du lemme 3.6.4 puisque $$\frac{1}{|\F^*/\F^{*2}|}\sum_{x\in \F^*/\F^{*2}}\tilde\omega_x(a)f(x)=A^1_{  s +\frac{1}{2} , s }(a,1,1)= A^1_{ s , s +\frac{1}{2}}(1,a,1)=0\ \text{pour }a\not=1 $$ donc $f(x)=f(1).$\\

\noi  La derni\ère \égalit\é s'obtient avec l'identit\é $\rho(\tilde\omega_u |\ |^{s +\frac{1}{2}} )\rho(\tilde\omega_u |\ |^{-s  +\frac{1}{2}} )=\tilde\omega_u(-1).$\fdem

  \vskip 4mm
  
\noi{\bf Remarques}:\\

1) Au vu du lemme 3.6.8, il semble inutile d'introduire la fonction $h(\pi)$ (ici quotient de fonctions $\rho$ de Tate) mais je ne sais pas si ce r\ésultat subsiste en caract\éristique r\ésiduelle diff\érente de $2$ lorsque $\F$ est un corps $\goth p-$adique.\\

2) Par le 5) du th\éor\ème 3.6.5, les valeurs prises par $A^1_{s,s+\frac{1}{2}}(\ ,\ ,1)$ sont les les coefficients de l'\eq de la fonction Z\êta associ\ée \à une forme quadratique isotrope de discriminant $-1$, cas trait\é dans \cite{dat}, malheureusement nos r\ésultats (B) du lemme 3.6.7) ne sont pas en accord avec les r\ésultats de la proposition 2.5 de \cite{dat}.\\ \\
 
 \subsection{\bf Application au cas  archim\édien}
 \bigskip
 
 L'\équation fonctionnelle abstraite permet de d\éterminer le
polynome $b_{\goth g,P_{\goth t}}$, et par cons\équent tous les
polynomes $b_k,k=1,...,p,$ par descente.
 
  $O$ \étant une
orbite de $P$ dans $\goth g"_1$  et  $O^*$ \étant une
orbite de $P_{\goth t}$ dans $\goth g"_{-1},$ les op\érateurs diff\érentiels introduits au \S 3.4 permettent
classiquement de d\éterminer la forme des coefficients
$a_{O,O^*}.$ Ce coefficient est not\é simplement $a(\overline\omega_m,s)$ dans le cas complexe puisque $\goth g"_1$  (resp. $\goth g"_{-1})$  est une seule orbite.

 \begin{lem} On suppose que $p\≥1.$ $\forall s\in
\C^p$ et $k=1,...,p:$
 \begin{enumerate}
\item Dans le cas complexe, soit $a(s)=a(\overline\omega_0,s):$
$$\begin{array}{ccl} a ( s)&=& 
( 2i\pi)^{-2d_k}\Bigl(b_k(s)\Bigr)^2a (s-\frac{1}{m_k}1_p+1_{p-k})\\
 &= &
(2i\pi)^{-2d'_{p-k}}\Bigl(b^*_k(s^*-1_{p-k} -(N-1)1_p)\Bigr)^2a (s-1_k).
\end{array}$$
\item Dans le cas r\éel: 
$$\begin{array}{ccl} a_{O^*,O}(\tilde\omega_{\pm 1};s)&=& 
(-2i\pi)^{-d_k}b_k(s)a_{O^*,O}( \tilde\omega_{\pm 1}\tilde\omega_{-1}^{1_p}\tilde\omega_{-1}^{1_{p-k}};s-\frac{1}{m_k}1_p+1_{p-k})\\
 &= &
(2i\pi)^{-d'_{p-k}}b^*_k(s^*-1_{p-k} -(N-1)1_p)a_{O^*,O}(  \tilde\omega_{\pm 1}\tilde\omega_{-1}^{1_k};s-1_k).
\end{array}$$
 \end{enumerate}
\noi $m_k=1$  \à l'exception  du cas classique $(C_n,\alpha_k)$ avec le parabolique $P'_0$  et  des formes de $(E_7,\alpha_6)$ pour lesquels $\frac{1}{m_1}=2$ (et $k=1$).
 \end{lem}
 
 \dem
 
Pour \établir ce r\ésultat "classique," dans le cas r\éel on calcule de 2
mani\ères  $Z^*_{O^*}( \goth F(F_k f);\omega ,s)$ et $Z^*_{O^*}(
\goth F(F^*_k(\partial ) f);\omega\tilde\omega_{-1}^{1_k} ,s-1_k).$

Pour $Z^*_{O^*}( \goth F(F_k f);\omega ,s),$ on applique  
d'une part  le
th\éor\ème 3.5.2, d'autre part on applique d'abord le lemme
3.4.1 \à 
$\goth F(F_k f)$ puis une
int\égration par partie  suivie du lemme 3.4.2 et pour finir
\à nouveau  le th\éor\ème 3.5.2.

Pour $Z^*_{O^*}(
\goth F(F^*_k(\partial ) f);\omega \tilde\omega_{-1}^{1_k} ,s-1_k), $ on applique
d'une part le lemme 3.4.1 \à
$\goth F(F^*_k(\partial ))$ puis le th\éor\ème 3.5.2, 
  d'autre part  on applique d'abord le th\éor\ème 3.5.2 puis  
 une int\égration par partie   suivie du lemme 3.4.2. 
 
 \noi Dans le cas complexe, on proc\ède de m\ême avec $Z^*(  \goth F(|F_k|_{\C} f); s)$ et $Z^* (
\goth F(F^*_k(\overline{\partial })F^*_k(\partial ) f); s-1_k).$\fdem\\
 
 \noi En prenant $k=p$ puis en appliquant \à $2$ reprises (avec $s$ ``convenable") le lemme 3.7.1, on en d\éduit:\\
 
 \begin{lem}  \begin{enumerate}
\item Dans le cas complexe:

 Pour $p\≥1$ on a :   
 $b_p(s)^2= b^*_p(s^* -(N-1)1_p)^2.$  
 
\noi Pour $p\≥2$ et $k=1,..,p-1$ on a:
$b_p(s)^2=b_k(s)^2b^*_{p-k}(s^*-(N-1)1_p)^2,$

\noi \à l'exception des formes de
$(E_7,\alpha_6)$ et  du cas classique $(C_n,\alpha_k)$ avec le parabolique $P'_0$ pour
lesquels on a:
$ \Bigl(b_2(s)b_2(s-1_2)\Bigr)^2=\Bigl(b_1(s)b^*_1(s^*-(N-2)1_2)\Bigr)^2.$ 
\item  Dans le cas r\éel: 

 Pour $p\≥1$ on a :   
 $b_p(s)=(-1)^db^*_p(s^* -(N-1)1_p).$ 
  
\noi Pour $p\≥2$ et $k=1,..,p-1$ on a:
$b_p(s)=b_k(s)(-1)^{d'_k}\ b^*_{p-k}(s^*-(N-1)1_p),$

\noi \à l'exception des formes de
$(E_7,\alpha_6)$ et  du cas classique $(C_n,\alpha_k)$ avec le parabolique $P'_0$ pour
lesquels on a:
$ b_2(s)b_2(s-1_2)=b_1(s)b^*_1(s^*-(N-2)1_2).$ 
\end{enumerate}
 
\end{lem}

\noi Avec les notations de la proposition
3.4.4, on  obtient:

\begin{prop}
On suppose que $p\≥2.$ Soient $b=b_{\goth
g,P(H_1,..,H_p)},$
$P_1$ et
$P_2$ les sous-groupes paraboliques associ\és
respectivement   \à $H_1,..,H_k$ dans $Aut(\goth U)$ et   \à
$H_{k+1},..,H_p$ dans $Aut(\goth U').$

\noi On a:
$$b(s_1,...,s_p)=A_k^{-d_k}B_k^{-d'_k}
b_{\goth U,
P_1 }(s_{p-k+1},...,s_p)b_{\goth U',
P_2 }(s_1,...,s_{p-k-1}, \sum_{p-k\≤i\≤p}s_i+r_k)$$
avec $r_k=\frac{p_1+2p_2} {2d'_k}$ et pour $i=1,2:$
$p_i=$dim$(E_i(h_k)\cap
\goth g_i),$

\noi \à l'exception du cas classique $(C_n,\alpha_k)$ avec le parabolique $P'_0$ et 
des formes de
$(E_7,\alpha_6)$ pour lesquelles on a:
$$b(s_1, s_2)=±\ C^4. \
s_2(s_2+3)(2s_1+s_2+4)(2s_1+s_2+7)\ \hbox{avec}\ C=\frac{2 }{
\widetilde B(H_0,H_0)}.$$
 \end{prop}

\dem 1)  Dans le cas r\éel, on applique les lemmes 3.7.2 et 3.4.4.

\noi Soit $N_{\goth U'}$ la constante associ\ée au \PV $(\goth U'_0,\goth U'_1),$ dans les notations du lemme 1.4.7 et avec cette  \demo on a:
$$N=\frac{2p_{0,2}+p_1}{2d'_km'_k}=(\frac{p_1}{2d'_k}+\frac{dim(\goth U')+p_2}{d'_k})\frac{1}{m'_k}=\frac{r_k}{m'_k}+   N_{\goth U'}\ \hbox{ 
donc }\ N-N_{\goth U'}=\frac{r_k}{m'_k}.$$

\noi Lorsque $d_k+d'_k=d_p,$  la relation (3) de la \demo  du lemme 1.4.7 nous donne $m'_k=1$ d'o\ù le r\ésultat.

\noi  Sinon on est dans l'un des $2$ cas cit\és dans la proposition et $k=1$ avec $m'_1=\frac{1}{2},$ $d_1=d'_1=d$ et on a:
$$b_2(s_1,s_2)b_2(s_1,s_2-1)=(A_1B_1)^{-d_1} b_{\goth U}(s_2)b_{\goth U}(s_2-1)b_{\goth U'}(2s_1+s_2+ 2r_1-1)
b_{\goth U'}(2s_1+s_2+ 2r_1).$$
\noi Dans le cas des formes r\éelles de
$(E_7,\alpha_6),$  les invariants relatifs fondamentaux des pr\éhomog\ènes $(\goth U_0,\goth U_1)$ et $(\goth U'_0,\goth U'_1)$ sont des formes quadratiques et on applique la remarque 3.6.3 en  notant que $N=4r_1=8.$

\noi 2) Dans le cas complexe, on applique le r\ésultat du cas r\éel $((\goth g_{\R})_0,(\goth g_{\R})_1)$ avec le \sg parabolique associ\é \à $H_1,...,H_p$  car  les $2$ \PVs ont les m\êmes  polynomes de Bernstein.

\fdem

 \begin{rema}: 
 
 \begin{enumerate}
 
 \item On ne donne pas le r\ésultat lorsque le \PV est de type $(C_n,\alpha_k)$ avec l'action du \sg parabolique $P'_0$ car ce r\ésultat ne sera pas utilis\é ult\érieurement; lorsque l'alg\èbre est d\éploy\ée, donc $k$ est pair, $b(s_1,s_2)$ est proportionnel \à:
$$\prod_{j=0}^{\frac{k}{2}-1}\Bigl((s_2+2j)(2s_1+s_2+2n-2k-1-2j)\Bigr).$$
Dans le cas $(E_7,\alpha_6),$ avec la normalisation $C=1,$ on \établira que $b(s_1,s_2)$ est unitaire (cf corollaire 8.1.5).\\

\item A l'exception des cas indiqu\és dans la proposition 3.7.3, on a imm\édiatement par r\écurrence l'existence, pour $\ell=1,...,p,$ de nombres rationnels positifs, $\lambda_{\ell,j},j=1,...,d_{p-\ell+1}-d_{p-\ell},$ avec $d_0=0,$ tels que $b(s_1,...,s_p)$ soit proportionnel \à:
$$ \prod_{\ell=1}^p\biggl(\prod_{j=1}^{d_{p-\ell+1}-d_{p-\ell}}(s_{\ell}+...+s_p+
\lambda_{\ell,j})\biggr).$$
Si on note $\lambda^{\goth g}_{\ell,j} $ ceux associ\és au \PV $(P_{\goth t},\goth g_1)$ alors on a la formule de r\écurrence suivante:
$$\lambda^{\goth g}_{\ell,j} =\Biggl\{\begin{array}{lrl}&\lambda^{\goth U'}_{\ell,j} +r_k\  &\ell=1,...,p-k\\
\\
&\lambda^{\goth U}_{\ell-p+k,j}  \  &\ell=p-k+1,...,p.\end{array}$$
$b(s_1,...,s_p)$ sera donn\é explicitement dans chaque cas lorsque le\sg parabolique est maximal parmi les \sgs paraboliques tr\ès sp\éciaux.
\end{enumerate}
\end{rema}
\bigskip
\newpage
\section{D\écomposition des mesures sur $\goth g_1$ et
$\goth g_{-1}$ et application}  
\bigskip 
Dans le cas archim\édien, on a montr\é que le polynome de Bernstein, $b_{\goth g,P_{\goth t}},$ associ\é au \PV $(P_{\goth t},\goth g_{\pm 1})$ est un produit de $2$ polynomes de Bernstein associ\és \à $2$ \PVs de rang plus petits (cf.prop.3.7.3). Il en est de m\ême, sous une forme plus \élabor\ée faisant intervenir les orbites, pour les coefficients de l'\équation fonctionnelle v\érifi\ée par la fonction Z\éta (cf.$\S$  5).

Pour obtenir ce type de r\ésultat, on se propose de d\écomposer les mesures de Haar sur $\goth g_{\pm 1}$ sous l'action du radical unipotent du sous-groupe parabolique associ\é \à $\goth t_k=\F h_k\oplus\F H_0,$ avec $1\≤k\≤ p-1,$
de donner la d\écomposition correspondante sur la transformation de Fourier ainsi que l'expression de la mesure relativement invariante par $G_{h_k}$ sur $W_{h_k}$ (cf.lemme 1.1.1).

Ces r\ésultats g\én\éralisent ceux obtenus ant\érieurement dans le cas commutatif (\cite{muller1}) et utilis\és dans \cite{boppruben} (th.4.28). \\

\subsection{Notations}

\noi Pour $ i,j\in \Z $  soit
 $E_{i,j}=E_{i,j}( h_k)=E_i( h_k)\cap E_j(2H_0- h_k) $ et
  $ p_{i,j} $  sa dimension, not\ée \également $p_i$ lorsque $i=j$ $(p_i:=p_{i,i}).$ \\
  
  \noi On note:
  $E'_{2,0}=\{x\in E_{2,0}\ |\ (x,h_k )$ se compl\ète en un \Sl $1$-adapt\é $\}$ et $E'_{0,2}=\{x\in E_{0,2}\ |\ (x,2H_0-h_k )$ se compl\ète en un \Sl $1$-adapt\é $\}.$\\ 
  
\noi   On rappelle que:
$$W_k:=W_{ h_k}=\{x+y\ |\ x\in E'_{2,0} ,y\in E'_{0,2}  ,[x,y]=0\}$$
et  on note: $$\goth n_k:=\goth n_{\goth t_k}=E_{-1,1}\oplus E_{-2,2}\quad , \quad N_k:=N_{\goth t_k}=exp(ad( \goth n_k)).$$ 
Lorsque $x+y\in W_k,$   soit
$\theta_x:=\theta_{x, h_k}(-1),$  
$\theta_y:=\theta_{y,2H_0- h_k}(-1)$ et $\theta_{x+y}:=\theta_{x+y,2H_0 }(-1)$
 alors  $\theta_{x+y}=\theta_x .\theta_y=\theta_y .\theta_x$ est une bijection de
$E_{i,j}$ sur
$E_{-i,-j}$ donn\ée par:
$$  \theta_x /_{E_{-1,±1}}= -ad(x)/_{E_{-1,±1}}\quad ,\quad
 \theta_x /_{E_{1,±1}}= -ad(x^{-1}) /_{E_{1,±1}}\quad ,\quad
\theta_x^2
 /_{E_{i,j}}=(-1)^i.Id$$
$$  \theta_x /_{E_{-2,±2}}=\frac {1}{2}ad(x)^2/_{E_{-2,±2}}\quad
,\quad
 \theta_x /_{E_{2,±2}}=\frac {1}{ 2}ad(x^{-1})^2 /_{E_{2,±2}} $$
avec des relations analogues pour $\theta_y .$\\

\begin{defi}  Pour  $x+y\in W_k,$ soit $F_{x,y}$ la restriction
de 
${\tilde B}( \theta_{x+y}(\ ),\ )$ \à $ E_{-1,1}\times E_{-1,1},$   $Q_{x,y}$ 
 la forme quadratique associ\ée, pour $  A,B\in E_{-1,1} $ on a:
 $$ 
F_{x,y}(A,B)=\tilde{B}([x,[y^{-1},A]],B)\  ,\ 
Q_{x,y}(A)=\frac{1 }{2}F_{x,y}(A,A)=\frac{1 }{2}\tilde B((ad(A)^2(x),y^{-1}),$$ 
et $\gamma_k(x,y)$ la constante
associ\ée au caract\ère quadratique $\tau\circ Q_{x,y}$ lorsque
$E_{1,1}\not=\{0\},$ et  $\gamma_k(x,y)=1$ sinon.

\noi$\gamma_k(\ ,\ )$ est constante sur les orbites de $G_{ h_k}$ dans $ W_{ h_k}.$

\end{defi}

 \begin{defi}  $V$ \étant un sous-espace vectoriel de $\goth g$ et $V^*$ son dual dans $\goth g$ pour $\tilde B,$ on note $\four_V $ (resp.$\overline{\four_V }$) la transformation de Fourier d\éfinie sur $V$ par la la restriction de $\tau\circ \widetilde B$  (resp.$\overline{\tau\circ \widetilde B}$) \à $V\times V^*$
\end{defi}

\subsection{Mesures sur les sous espaces $E_{i,j}$ et $E_{-i,-j}$ pour $(i,j)\not=(0,0)$}
\bigskip 

\begin{lem} Pour $(i,j)\not=(0,0)$ il existe  une unique mesure de Haar sur $E_{i,j}$ et  une unique mesure de Haar sur $E_{-i,-j},$ not\ées $dx,$ telles que :
\begin{enumerate}
\item $\overline{\four_{E_{-i,-j}}}\circ \four_{E_{i,j}}=\four_{E_{-i,-j}}\circ\overline{\four_{E_{i,j}}}=Id_{E_{i,j}}$

\item Pour
tout
$  x+y\in W_ k$ et $  f\in L^1(E_{-i,-j})$  on a:
$$  \int_{E_{i,j}}f(\theta_{x +y}(z))\
dz=|F_k(x)|^{\frac{ip_{i,j} }{d_k}}|P_{p-k}(y)|^{\frac{jp_{i,j} }{
d'_k}}\int_{E_{-i,-j}}f(u)\ du .$$
 \end{enumerate}
\end{lem}

\dem Omise car elle est analogue \à celle du lemme 3.2.1 en se rappelant  que $|det(\theta_{x+y}/E_{i,j})|^{\frac{1}{2}}
|F_k(x)|^{\frac{ip_{i,j} }{ 2d_k}}|P_{p-k}(y)|^{\frac{jp_{i,j} }{
2d'_k}}$ est une constante (cf.\demo du lemme 1.4.7); dans une base $\goth B$ de $E_{i,j}$ (donc $E_{-i,-j}$
est muni de la base   $\goth B^*$  duale pour $\tilde B$) elle a pour
expression $dz=\lambda_{\goth B}d_{\goth B}z$ avec:
$$ \lambda_{\goth B}=|det(\ \theta_{x+y}: E_{i,j}(\goth B)\rightarrow E_{-i,-j}(\goth B^*)\ )|^ {\frac{1}{2}}
|F_k(x)|^{\frac{ip_{i,j} }{ 2d_k}}|P_{p-k}(y)|^{\frac{jp_{i,j} }{
2d'_k}} \ $$
(et $\ds\frac{1}{\lambda_{\goth B}}d_{\goth B^*}z$ sur $E_{-i,-j}$).\fdem\\

\noi Dor\énavant chaque sous-espace $E_{i,j}$ pour $(i,j)\not=(0,0)$ est muni de cette unique mesure de Haar et avec cette normalisation on a les lemmes suivants: \\

\begin{lem} Soient $(i,j)\not=(0,0)$ et $  g\in L^1(E_{i,j}).$
\begin{enumerate}
\item Pour    $x\in E'_{2,0}$ et $y\in E'_{0,2} ,$ on a:
$$  \int_{E_{i,j}}g( u)\ du  =|F_k(x)|^{\frac{ip_{i,j} }{ d_k}} \int_{E_{-i,j}}g(\theta_x(z))\
dz  
  =|P_{p-k}(y)|^{\frac{jp_{i,j} }{
d'_k}}\int_{E_{i,-j}}g(\theta_y(z))\
dz  .$$
\item Soit $\sigma$ une involution de $\goth g$ telle que $\sigma ( h_k)=- h_k,$ $\sigma (H_0)=-H_0$ et on suppose que $W_{ k,\sigma}=\{z\in W_k \ |\  \sigma (x)=\Phi(x)\}\not=\emptyset,$ alors 
$$\int_{E_{-i,-j}}g( \sigma(v))\ dv=|F_k(x)|^{\frac{-ip_{i,j} }{ d_k}}
|P_{p-k}(y)|^{\frac{-jp_{i,j} }{
d'_k}}\int_{E_{i,j}}g( u)\ du\ \ \hbox{
avec}\  x+y\in W_{ k,\sigma}.$$
\end{enumerate}
\end{lem}

\dem 1) Soit $x\in E'_{2,0} ,$ on le compl\ète par un \elt $z\in  E'_{0,2} $ de telle mani\ère que $x+z\in  W_k $ et on munit:

\noi $\bullet$ $E_{i,j}$ d'une base $\goth B$ et $E_{-i,-j}$ de la  base $\goth B^*,$ duale de $\goth B$ pour $\tilde B,$

\noi $\bullet$ $E_{-i,j}$ de la base $\goth B'=\theta_x(\goth B)$ donc $  \theta_x(\goth B^*)$ est la base duale de $E_{i,-j}$  pour $\tilde B.$

\noi Par le choix des bases, on a :$$|det(\ \theta_{x+y}: E_{i,j}(\goth B)\rightarrow E_{-i,-j}(\goth B^*)\ )|=|det(\ \theta_{x+y}: E_{-i,j}(\theta_x(\goth B))\rightarrow E_{i,-j}(\theta_x(\goth B^*))\ )|\ ,$$
donc $\ds\frac{ \lambda_{\goth B}}{ \lambda_{\goth B'}}=|F_k(x)|^{\frac{ip_{i,j} }{ d_k}}$ d'o\ù la premi\ère \égalit\é.
On proc\ède de m\ême pour la seconde \égalit\é.\\

2) Soit $x+y\in W_{k,\sigma},$ alors $\sigma$ et $\theta_{x+y}$ commutent donc $\theta=\sigma \theta_{x+y}$ est une bijection de $E_{i,j}$ pour laquelle $\theta^2=(-1)^{i+j}Id_{E_{i,j}}$ donc:
$$\int_{E_{i,j}}g( \theta(v))\ dv=\int_{E_{i,j}}g( v)\ dv=\int_{E_{i,j}}g\circ \sigma(\theta_{x+y}(v))\ dv$$
d'o\ù le r\ésultat. Notons que $|F_k(x)|^{\frac{-ip_{i,j} }{ d_k}}
|P_{p-k}(y)|^{\frac{-jp_{i,j} }{
d'_k}}$ est constant sur $W_{k,\sigma}$ et sa valeur d\épend de la normalisation choisie pour les invariants relatifs $F_k$ et $F_p.$ \fdem

\bigskip

\begin{lem} On suppose que $E_{-1,1}\not=\{0\}.$ Soit $x+y\in W_k,$ pour $u\in L^1( E_{-1,1})$ et $\four_{E_{-1,1}}(u)\in L^1(E_{1,-1})$ on a:
$$\int_{E_{-1,1}}\tau\Bigl(Q_{x,y}(A)\Bigr)\four_{E_{-1,1}}(u)\biggl(\theta_{x+y}(A)\biggr)dA=C(x,y)\int_{E_{-1,1}}\overline{\tau\Bigl(Q_{x,y}(A)\Bigr)}u(A)dA$$avec
$$ C(x,y)= \gamma_k (x,y)|F_k(x)|^{-\frac{p_{1}}{2d_k}}|P_{p-k}(y)|^{\frac{p_{1}}{2d'_k}} .$$

\end{lem}

\dem   C'est le r\ésultat concernant la transformation de
Fourier d'un caract\ère quadratique, d\û \à A.Weil (\cite{weil}) traduit dans notre situation. Rappelons-le
bri\èvement:
$E$ \étant un espace vectoriel de dimension finie sur $\F$ et $Q$
une forme quadratique non d\ég\én\ér\ée sur $E$, on note
$F(A,B)=Q(A+B)-Q(A)-Q(B)$ la forme bilin\éaire associ\ée , alors
$\tau\circ Q$ est un caract\ère quadratique non d\ég\én\ér\é du groupe additif
$E$ et on l'\égalit\é :
$$\int_E{\overline \tau}(Q(A))u(A)dA=\gamma^{-1}(\tau\circ Q)\mid
\rho\mid^{1\over 2}\int_E\tau(Q(A))u^*(\rho(A))dA$$
lorsque $u$ et $u^*$ sont int\égrables , avec $$
u^*\Bigl(\rho(A)\Bigr)=\int_Eu(B)\tau\Bigl(F(A,B)\Bigr)dB$$
$\rho$ \étant l' isomorphisme sym\étrique associ\é \à $Q$ de $E$ sur $E^*$ et
$|\rho|$ est d\éfinie par  $$\int_{E^*}g(u)du=|\rho|\int_Eg(\rho(A))dA$$
les mesures de Haar  sur $E$ et $E^*$ \étant duales.

\noi Dans notre situation $E=E_{-1,1},$ $E^*$ est identifi\é \à $E_{1,-1}$ \à l'aide de $\tilde B,$ 
 $Q:=Q_{x,y}$ donc $F:=F_{x,y}$ et $\rho:=\theta_{x+y}/E_{-1,1}$ donc  $u^* =\four_{E_{-1,1}}(u) $ et par le choix des mesures sur $E_{-1,1}$ et $E_{1,-1}$ ainsi que le lemme 4.2.1 on a  $|\rho|^{-1}=|F_k(x)|^{\frac{-p_{1}}{d_k}}|P_{p-k}(y)|^{\frac{p_{1}}{d'_k}}.$\fdem\\

 \begin{lem} $$\begin{array}{llrll}  \int_{\goth g_1}&f(x)\ dx\
&=&\int\int\int_{E_{2,0}\times E_{0,2}\times
E_{1,1}} &f(x+y+z)\ dx\ dy\ dz\  ,\\
   \int_{\goth g_{-1}}&g(y)\ dy\
&=&\int\int\int_{E_{-2,0}\times E_{0,-2}\times
E_{-1,-1}}&g(x+y+z)\ dx\ dy\ dz\ .  \end{array}$$
 \end{lem}

\bigskip

\dem Pour $f\in L^1(\goth g_1)$ et $g\in L^1(\goth g_{-1})$ on pose
$$\begin{array}{llrll}  \int_{\goth g_1}&f(x)\ d'x\
&:=&\int\int\int_{E_{2,0}\times E_{0,2}\times
E_{1,1}} &f(x+y+z)\ dx\ dy\ dz\  ,\\
   \int_{\goth g_{-1}}&g(y)\ d'y\
&:=&\int\int\int_{E_{-2,0}\times E_{0,-2}\times
E_{-1,-1}}&g(x+y+z)\ dx\ dy\ dz\ .  \end{array}$$ Il suffit de v\érifier les $2$ hypoth\èses du lemme 3.2.1, ce qui se fait imm\édiatement :

\noi $\bullet$  pour la transformation de Fourier  avec $f(x+y+z)=f_1(x)f_2(y)f_3(z)$

\noi $\bullet$  en prenant $x+y\in W_k$ et en appliquant le lemme 4.2.1 ce qui donne le r\ésultat puisque $$|F_k(x)|^{2Nm_k}|P_{p-k}(y)|^{2Nm_k}=|F(x+y)|^{2N}\  \ \hbox{
avec}\
Nm_k=\frac{\frac{p_{1} }{ 2}+p_{2,0}} {d_k}=\frac 
{\frac{p_{1} }{ 2}+p_{0,2}} {d'_k}$$ (d\ém. du lemme
1.4.7).\fdem
 
 \bigskip
 
 \noi Indiquons le dernier choix de mesures n\écessaire pour \établir le th\éor\ème 4.3.3 :  \bigskip

\begin{defi} Soit $x\in E'_2(h_k)\cap \goth g_1$ (resp.$y\in E'_2(2H_0-h_k)\cap \goth g_1$) et $\goth s$ la sous-alg\èbre engendr\ée par le \Sl  $(x,h_k,x^{-1})$ (resp$.(y,2H_0-h_k,y^{-1})$, on munit $\goth U(\goth s)_{\pm 1}$ des uniques mesures v\érifiant:
\begin {enumerate}

\item $\overline{\four_{\goth U(\goth s)_{-1} }}\circ \four_{ \goth U(\goth s)_1}= Id_{ \goth U(\goth s)_1}$

\item Pour
tout
$  x_0\in  \goth U'(\goth s)_1$   et $  f\in L^1( \goth U(\goth s)_{-1})$  on a:
$$ \begin{array}{rlll} \int_{ \goth U(\goth s)_1}f(\theta_{x_0}(z))\
dz &= & |P_{p-k}(x_0)|^{ \frac{p_{0,2}-p_{2}} {
d'_k}}&\int_{\goth U(\goth s)_{-1} } f(u)\ du  \\ 
(\hbox{resp.} &=& |F_k(x_0)|^{\frac{p_{2,0}-p_{2} }{d_k}}  &\int_{ \goth U(\goth s)_{-1}}f(u)\ du ).\end{array}$$
\end{enumerate}
\noi c'est \à dire que $\widetilde {B_{ \goth U(\goth s)}}=\widetilde B/      \goth U(\goth s)$ et les mesures sont adapt\ées \à la restriction de $P_{p-k}$ (resp.$F_k$) \à $ \goth U(\goth s)_1.$  

\end{defi}

\noi Ce choix est coh\érent avec le lemme 4.2.1 puisque  $\goth U(\goth s)_{\pm 1}={(E_{0,\pm 2})}_x={(E_{0,\pm 2})}_{x^{-1}}=E_{0,\pm 2}$
 (resp.${(E_{\pm 2,0})}_{y}={(E_{\pm 2,0})}_{y^{-1}}=E_{0,\pm 2}$)  lorsque le \PV est  commutatif. \\

\noi Avec ce choix et dans les notations de cette d\éfinition on a:

\begin{lem} \begin{enumerate}
\item Soit $x\in E'_{2,0},$ pour $f\in L^1(E_{ 0,2})$ et $g\in L^1(E_{ 0,-2})$ on a:
$$\begin{array}{lrl}
\int_{E_{0, 2}}f(z)dz&=&|2|_{\F}^{\frac{p_2}{2}}|F_k(x)|^{\frac{p_2}{d_k}}\int\int_{(u,v)\in E_{-2,2}\times {(E_{0,2})}_x}f([x,u]+v)dudv\\
\int_{E_{0, -2}}g(z)dz&=&|2|_{\F}^{\frac{p_2}{2}}|P^*_{k}(x^{-1})|^{\frac{p_2}{d'_{p-k}}}\int\int_{(u,v)\in E_{2,-2}\times {(E_{0,-2})}_x}f([x^{-1},u]+v)dudv.
\end{array}$$
\item Soit $y\in E'_{0,2},$ pour $f\in L^1(E_{ 2,0})$ et $g\in L^1(E_{-2,0})$ on a:
$$\begin{array}{lrl}
\int_{E_{ 2,0}}f(z)dz&=&|2|_{\F}^{\frac{p_2}{2}}|P_{p-k}(y)|^{\frac{p_2}{d'_k}}\int\int_{(u,v)\in E_{2,-2}\times {(E_{2,0})}_y}f([y,u]+v)dudv\\
\int_{E_{ -2,0}}g(z)dz&=&|2|_{\F}^{\frac{p_2}{2}}|F^*_{p-k}(y^{-1})|^{\frac{p_2}{d'_k}}\int\int_{(u,v)\in E_{-2,2}\times {(E_{-2,0})}_y}f([y^{-1},u]+v)dudv.
\end{array}$$
\end{enumerate}
\end{lem}

\dem 1) On a: 
$$E_{0,2}=[ x,E_{-2,2}]\oplus \
(E_{0,2})_x\quad \hbox{et}\quad E_{0,-2}=[
x^{-1},E_{2,-2}]\oplus
\ (E_{0,-2})_x\ .$$
Soient $(f_i)_{1\≤i\≤p_2}$ et
$(f_i^*)_{1\≤i\≤p_2}$  2 bases   de $E_{-2,2}$ et
$E_{2,-2},$ duales pour $\tilde B$, alors 
$([x,f_i])_{1\≤i\≤p_2}$ et $(\frac{1 }{
2}[x^{-1},f_i^*])_{1\≤i\≤p_2}$ sont 2 bases   de
$[x,E_{-2,2}]$ et
$[x^{-1},E_{2,-2}],$ duales pour $\tilde B,$ on les compl\ète par $2$ bases de $(E_{0,\pm 2})_x$ duales pour  $\tilde B.$

Notons $\lambda ,\lambda_2$ et $\lambda_ x$ les constantes associ\ées aux mesures  d\éfinies sur  $E_{0,2},E_{-2,2}$ par le lemme 4.2.1 et sur 
$({E_{0,2}})_x$  par la d\éfinition 4.2.5, ceci avec les bases duales indiqu\ées , il suffit de calculer $C=\frac{\lambda }{\lambda_2  .\lambda_ x}$ c'est \à dire (cf. \demo du lemme 3.2.1 pour ${(E_{0,2})}_x$ et du lemme 4.2.1 pour $E_{0,2}$ et $E_{-2,2}$):
$$\begin{array}{rll}C&=|det(\theta_{x+y}/E_{0,2})|^{\frac{1 }{
2} }
 |det(\theta_{x+y}/ E_{-2,2})|^{- \frac{1 }{
2} } 
 |det(\theta_y/(E_{0,2})_x)|^{- \frac{1 }{
2}}|F_k(x)|^{\frac{p_2}{d_k}} \\
 & = 
 |det(\theta_{x+y}/[x,E_{-2,2}])|^{ \frac{1 }{
2}}  |det(\theta_{x+y}/ E_{-2,2})|^{- \frac{1 }{
2} } |F_k(x)|^{\frac{p_2}{d_k}} \\
  &=|2|_{\F}^{\frac{p_2 }{ 2} }|F_k(x)|^{\frac{p_2}{d_k}}  \end{array}$$
 d'o\ù le r\ésultat.
 
 \noi 2) Idem\fdem

\bigskip

 \subsection{D\écomposition des mesures sur $\goth g_1$ et
$\goth g_{-1}$}

\bigskip

\begin{defi}
\begin{enumerate}
\item  Pour  $f$  appartenant \à  $L^1(\goth g_1),$ $x$ dans
$E'_{2,0}$  et $y$  dans  $E_{0,2},$  on pose :
$$\begin{array}{rll} S_f(x+y)&=&|2|_{\F}^{\frac{p_2 }{2}}|F_k(x)|_{\F} ^{
md'_k}\int_{\goth n_k}f(e^{ad(A)}(x+y))\ dA\\
T_f(x+y)&=& \Bigl\lbrace
\begin{array}{llr}
&|F_k(x)|_{\F}^{\frac{  p_1 }{ 2d_k}}
\int_{E_{-1,1}}f(e^{ad(A)}(x+y))dA&\quad \hbox{ si }\quad
E_{-1,1}\not=\{0\}\\ 
 & f(x+y) & \quad \hbox{ si }\quad
E_{-1,1} =\{0\}. \end{array}
\end{array}$$
\noi  avec $m=\displaystyle{\frac{r_k}{d_k}= \frac{\frac{p_1 }{ 2 }+p_2}{
d_kd'_k}}.$
 
\item Pour $g$   appartenant \à  $L^1(\goth g_{-1}),$ $x'$  dans 
$E_{-2,0}$  et
$y'$    dans
  $E'_{0,-2},$  soit
$$ \begin{array}{rll} S^*_g(x'+y')&=& |2|_{\F}^{\frac{p_2}{2}}| F^*_{p-k} (y')|_{\F}^{ md_k}
\int_{\goth n_k }g(e^{ad(A)}(x'+y'))\ dA\\
T^*_g(x'+y')&= &\Bigl\lbrace
\begin{array}{llr}
 &|{F^*_{p-k}}(y')|_{\F}^{  \frac{ p_1}{ 2d'_k}}
\int_{E_{-1,1}}g(e^{ad(A)}(x'+y'))dA  &\quad \hbox{  si }\quad
E_{-1,1}\not=\{0\} \\
&g(x'+y') &\quad \hbox{ si }\quad
E_{-1,1} =\{0\}. \end{array}
    \end{array}  $$  
\noi  $\goth n_k$ \étant muni de la mesure produit.
\end{enumerate}
\end{defi}

\bigskip
\begin{rema} Les normalisations de $S_f,S^*_g,T_f$ et $T^*_g$ sont  d\ûes au lemme 4.2.3. \\

\noi Soient  $A_1\in E_{-1,1},$  $A_2\in E_{-2,2},$ $ x\in E'_{2,0},$ $y\in E_{0,2},$ $ x'\in E_{-2,0},$ $y'\in E'_{0,-2},$ on a  (cf. \demo du lemme 1.1.1):
$$exp(ad(A_1+A_2))(x+y)=x+[A_1+A_2,x]+ y+\frac{1 }{ 2}ad(A_1)^2(x)= exp(ad(A_1)(x+[A_2,x]+ y)\  , $$
$$exp(ad(A_1+A_2))(x'+y')=x'+[A_1+A_2,y']+ y'+\frac{1 }{ 2}ad(A_1)^2(y')=exp(ad(A_1)(x'+[A_2,y']+ y')\  , $$
donc:
$$\begin{array}{rlll}S_f(x+y)&=&|2|_{\F}^{\frac{p_2 }{2}}|F_k(x)|_{\F}^{
\frac{p_2}{d_k}}&\int_{E_{-2,2} }T_f( x+[A,x]+ y)\ dA\\
 S^*_g(x'+y')&=&|2|_{\F}^{\frac{p_2}{2}}|F^*_{p-k}(y')|_{\F} ^{
\frac{p_2}{d'_k}}&\int_{E_{-2,2} }T^*_g( x'+[A,y']+ y')\ dA\ .\end{array}$$

\end{rema}

\bigskip

\begin{theo} \begin{enumerate}
\item  Pour f  dans  $L^1(\goth g_1)$   on a :
 $$ \begin{array}{rll}\int_{\goth
g_1}f(x)dx&=& \int\int_{
E_{2,0}\times E_{0,2}}T_f(x+y) |F_k(x)|^{\frac{p_1 }{
2d_k } }dxdy\\  
&=&  
\int_ {E_{2,0}}|F_k(x) |^{\frac{ p_1 }{ 2d_k}}\ (\
\int_{(E_{0,2})_x}S_f(x+y)\ dy\ )\ dx\ .\end{array}$$
\item  Pour g dans $L^1(\goth g_{-1})$   :
$$\begin{array}{rll} \int_{\goth
g_{-1}}g(x)dx&=& \int\int_{ E_{-2,0}\times
E_{0,-2}}T^*_g(x'+y')\mid F^*_{p-k} (y')\mid ^{\frac{p_{1} }{
2d'_k}}dx'dy'\\
 &=&    \int_{
E_{0,-2}}|F^*_{p-k}(y') |^{ \frac{p_1 }{2d'_k}}\ (\
\int_{(E_{-2,0})_{y'}}\ S^*_g(x'+y')\ dx'\ )\ dy'.\end{array}$$
 \end{enumerate}
\end{theo}

\dem  Par le lemme 4.2.4 on a 
$$\int_{\goth g_1}f(x)\ dx\
=\int\int_{E_{2,0}\times E_{0,2}}\biggl(\int_{E_{1,1}} f(x+y+z)\ dz\biggr)\ dxdy=\int\int_{E'_{2,0}\times E_{0,2}}\biggl(\int_{E_{1,1}} f(x+y+z)\ dz\biggr)\ dxdy$$
Or pour $x\in E'_{2,0}$ on a par le 1. du lemme 4.2.2:
$$\begin{array}{rll}
\int_{E_{1,1}} f(x+y+z)\ dz&=&|F_k(x)|^{\frac{p_1}{d_k}}\int_{E_{-1,1}} f(x+y+\theta_x(A))\ dA\quad \hbox{donc} \\
 \int_{\goth g_1}f(x)\ dx&=&\int\int_{E'_{2,0}\times E_{0,2}}\biggl(|F_k(x)|^{\frac{p_1}{d_k}}\int_{E_{-1,1}} f(x+y+\theta_x(A))\ dA\biggr)\ dxdy\\
&=& \int_{E'_{2,0}}\biggl(\int _{E_{0,2}}\biggl(|F_k(x)|^{\frac{p_1}{d_k}}\int_{E_{-1,1}} f(x+y+[A,x] \ dA\biggr)\ dy\biggr)dx.\end{array}$$
 D'o\ù   la premi\ère \égalit\é du 1. est obtenue en effectuant la translation  $y\rightarrow y+ \frac{1 }{ 2}ad(A)^2(x)$ et pour l'autre \égalit\é de 1., on applique le lemme 4.2.6 et la remarque 4.3.2 \à la premi\ère \égalit\é. 
 
\noi On proc\ède de m\ême pour 2.
\fdem\\

\begin{rema} Soit $f$  dans  $L^1(\goth g_1)$ (resp.
$g\in  L^1(\goth g_{-1})$),  $T_f$  (resp.$T^*_g$) est d\éfini presque partout et  $|F_k(x) |^{\frac{ p_{1} }{ 2d_k}}T_f(x+y)$ (resp.${{F^*}_{p-k}(y')}^{\frac{p_1}{2d'_k}}T^*_g(x'+y'))$ est int\égrable sur $E_{2,0}\times E_{0,2}$ (resp.$E_{-2,0}\times E_{0,-2})$ muni de la mesure produit.
\end{rema}

\bigskip

\begin{theo} Soit $f\in  \EuScript S(\goth g_1)$ alors:\begin{enumerate} 

\item ${\four}_{E_{0,2}}\Bigl(T_{f}(x+.\ )\Bigr)(y')=
 {\four}_{ (E_{0,2})_{x }}\Bigl(S_f(x+.\ )\Bigr)(y')$ lorsque $[x,y']=0.$

\item Pour $x'\in E_{-2,0}$ et $y'\in E'_{0,-2}$   tels que
$[x',y']=0$ on a :$$ S^*_{\bf
{\four (f)}}(x'+y')=   {\four}_{ (E_{2,0})_{y' }}\biggl(\ \
\gamma_k( x\ ,y'^{-1}). {\four}_{  E_{0,2}}(T_{\bf f}(x+.\ ))(y')\
\biggr)(x').$$
\end{enumerate}
\end{theo}

\dem 1. R\ésulte du lemme 4.2.6 et de la remarque 4.3.2.\\

\noi 2. \noi  Soient $C\in \goth n_k:$ $C=C_1+C_2,$ avec $C_i\in
E_{-i,i},$ $ x'\in E_{-2,0},$ $y'\in 
E'_{0,-2} $ et 
$$\begin{array}{rll}(1) &=&\four
(f)(exp(adC(x'+y'))) \\
 &=&
\int\int\int_{x_i\in E_{i,2-i},i=0,1,2} f(x_2+x_1+x_0)   \tau( {\tilde B}(
 (exp(adC)(x'+y'),x_2+x_1+x_0 ))  
dx_2dx_1dx_0 \end{array}$$par le lemme
4.2.4, alors par orthogonalit\é et en utilisant le
th\éor\ème de Fubini, on a:
$$\begin{array}{rll}(1)&=&\int\int_{x_i\in E_{i,2-i},i=1,2} {\four}_{E_{0,2}}\Bigr(f(x_2+x_1+\bullet)\Bigl)(y') 
\tau({\tilde B}\Bigr(x_2,x'+{1\over 2}ad(C_1)^2(y')+[C_2,y']\Bigl)).\\
& &\tau({\tilde
B}\Bigr(x_1,[C_1,y']\Bigl))dx_2dx_1\\
&=& \int_{E_{1,1}} (2) \tau({\tilde
B}\Bigr(x_1,[C_1,y']\Bigl))\   dx_1\ \hbox{avec:}\\
(2)&=&  \int_{ E_{2,0}} {\four}_{E_{0,2}}\Bigr(f(x_2+x_1+\bullet)\Bigl)(y') 
\tau({\tilde B}\Bigr(x_2,x'+{1\over 2}ad(C_1)^2(y')+[C_2,y']\Bigl))dx_2.\end{array}$$
On applique le lemme 4.2.6 \à $E_{2,0}$ avec $y'^{-1}$ d'o\ù:
$$(2)=|2|^{\frac{p_2 }{2}}  |F_{p-k}^*(y')|^{-\frac{p_2}{
d'_k}}\int\int_{(u,v)\in E_{2,-2}\times (E_{2,0})_{y'}} \four_{E_{0,2}}(\ f([y'^{-1},u]+v+x_1+\bullet)\ )(y')(3)dudv$$
 avec
$$\begin{array}{rll} (3)&=& 
 \tau( {\tilde B}([y'^{-1},u]+v ,x'+ \frac{1 }{ 2}ad(C_1)^2(y') +[C_2,y'])\\
 (3)&=&\tau( {\tilde B}( v,x'+\frac{1 }{ 2}ad(C_1)^2(y')).
\tau( {\tilde B}( u,-2C_2+\frac{1 }{ 2}\ [ad(C_1)^2(y'),
y'^{-1}]))\ \hbox{lorsque}\ [x',y']=0 \end{array}$$
 par orthogonalit\é.\\
 
 \noi Donc, lorsque $[x',y']=0,$ on a en appliquant le th\éor\ème de Fubini:
 $$ \begin{array}{rll}(2)&=& |2|^{\frac{p_2  }{2}}  |F_{p-k}^*(y')|^{-\frac{p_2}{
d'_k}}\int_{(E_{2,0})_{y'}} \four_{E_{2,-2}}(g(\bullet,v,x_1))(-2C_2+\frac{1 }{ 2}\ [ad(C_1)^2(y'),y'^{-1}]).\\
& & \tau( {\tilde B}( v,x'+\frac{1 }{ 2}ad(C_1)^2(y'))dv.\end{array}$$
avec $  g( u,v,x_1)=\four_{E_{0,2}}(\ f([y'^{-1},u]+v+x_1+\bullet)\ )(y') .$\\

\noi Comme $f\in \EuScript S(\goth g_1)$ et que $y'$ est fix\é dans $E'_{0,-2},$ on a $  g\in   \EuScript S(E_{2,-2}\times (E_{2,0})_{y'}\times E_{1,1}).$\\

 \noi On a   \à \évaluer :
 $$ \begin{array}{rll} (4)&=& S^*_{ \four f}(x'+y')\\
&=&|2|^{\frac{p_2}{ 2}} |
F_{p-k}^*(y')|^{md_k}\int_{E_{-1,1}}\ (\ \int_{E_{-2,2}}\
(1)\ dC_2\ )\ dC_1\\
 &=&|2|^{\frac{p_2}{ 2}}  
|F_{p-k}^*(y')|^{ md_k}. 
 \int_{E_{-1,1}}\ (5)\ dC_1\quad \hbox{avec}\\ 
(5)&=& 
\int_{E_{-2,2}}\ (\ \int_{E_{1,1}}(2)\tau({\tilde
B}(x_1,[C_1,y']))\ dx_1\ ) \ dC_2\\
&=&|2|^{\frac{p_2 }{ 2}}  
|F_{p-k}^*(y')|^{-\frac{p_2 }{d'_k}}. \\
 
& &\int_{E_{-2,2}}\ 
\bigl(\ \int_{E_{1,1}}\
\biggl(\int_{ 
(E_{2,0})_{y'}}\four _{E_{2,-2}}(\ g(\bullet,v,x_1)\ )
(-2C_2+\frac{1 }{ 2}\ [ad(C_1)^2(y'),y'^{-1}])\\
& &\tau( {\tilde B}( v,x'+\frac{1 }{ 2}ad(C_1)^2(y'))\ dv\biggr)\ 
\tau({\tilde
B}(x_1,[C_1,y']))\ dx_1\ \bigr)\ dC_2\ ,\end{array}$$

\noi mais cette int\égrale converge absolument donc on peut
appliquer le th\éor\ème de Fubini d'o\ù: 
$$ \begin{array}{rll} (5)&=&  \int_{ 
(E_{2,0})_{y'}} \tau( {\tilde B}( v,x'+ \frac{1 }{ 2}ad(C_1)^2(y'))\ 
(\int_{E_{1,1}}(6)\tau({\tilde
B}(x_1,[C_1,y']))\ dx_1\ )\ dv\ \ \hbox{avec}\\
 (6)&= &|2|^{\frac{p_2}{ 2}}
 |F_{p-k}^*(y')|^{-\frac{p_2 }{
d'_k}}\int_{E_{-2,2}}\ 
\four_{E_{2,-2}}(\ g(\bullet,v,x_1)\ )
(-2C_2+\frac{1 }{ 2}\ [ad(C_1)^2(y'),y'^{-1}]) \ dC_2 \\
&= & |2|^{-\frac{p_2}{ 2}}|F_{p-k}^*(y')|^{-\frac{p_2 }{
d'_k}}. g(0,v,x_1)\\
&=&   |2|^{-\frac{p_2 }{ 2}}
|F_{p-k}^*(y')| ^{-\frac{p_2 }{
d'_k}}. \four_{E_{0,2}}
(f(v+x_1+\bullet)(y')\end{array}$$
\noi Ainsi l'int\égration suivant le sous-espace $E_{-2,2}$ est
"annul\ée" par "double transformation de Fourier". Il reste 
encore \à simplifier l'int\égration suivant le sous-espace
$E_{-1,1}$ \à l'aide du lemme 4.2.3.\\

\noi Nous avons:
$$ \begin{array}{rll}(4)&=&S^*_{ \four f}(x'+y')\\
&=& |2|^{\frac{p_2}{ 2}}  | F_{p-k}^*(y')|^{md_k}.\\
 & &\int_{E_{-1,1}}\ (\ 
\int_{ 
(E_{2,0})_{y'}} \tau( {\tilde B}( v,x'+ \frac{1 }{ 2}ad(C_1)^2(y'))\ h(C_1,v)
 \ dv\ )\ dC_1\ ,\\ 
 \hbox{avec}\ h(C_1,v)&=&\int_{E_{1,1}}(6)\tau({\tilde
B}(x_1,[C_1,y']))\ dx_1\\
&=& |2|^{-\frac{p_2 }{ 2}}  | F_{p-k}^*(y')|^{ -\frac{p_2 }{ d'_k} }.\four_{E_{1,1}}\biggl(\four_{E_{0,2}}(f(v+x_1+\ .))(y')\biggr)([C_1,y']) \ ,\end{array}$$
Comme $f\in \EuScript S(\goth g_1)$ et que $y'$ est fix\é dans $E'_{0,2},$ on a $h\in \EuScript S(E_{-1,1}\times (E_{2,0})_{y'})$ donc l'int\égrale double figurant dans (4) converge absolument et
 on utilise \à
nouveau le th\éor\ème de Fubini d'o\ù:
 $$ \begin{array}{rll}(4)&=&\four_{(E_{2,0})_{y'}}(h_1)(x')  \quad  
  \hbox{et
pour }\quad v\in (E'_{2,0})_{y'} \\
h_1(v)&=&|F_{p-k}^*(y')|^{\frac{p_1 }{
2d'_k}}\int_{E_{-1,1}}
\tau\circ Q_{v,y'^{-1}}(C_1)\ ( \int_{E_{1,1}}  g(0,v,x_1).\tau({\tilde
B}(x_1,[C_1,y']))\ dx_1\ )dC_1\\
&=&|F_{p-k}^*(y')|^{\frac{p_1 }{
2d'_k}} |F_k(v)|^{\frac{p_1 }{ d_k}}.\int_{E_{-1,1}}
\tau (Q_{v,y'^{-1}}(C_1)) (
\int_{E_{-1,1}}u(A) \tau( -F_{v,y'^{-1}}(C_1,A) )dA )
 dC_1\ ,\end{array}$$
avec $u(A)= g(0,v, [A,v]),$ en appliquant le 1. du lemme 4.2.2 \à 
$E_{1,1}$ avec $\theta_v.$\:

\noi D'o\ù:
$$ \begin{array}{rll}h_1(v) & = &|F_k(v)|^{\frac{p_1 }{ d_k}}|F_{p-k}^*(y')|^{\frac{p_1 }{
2d'_k}}.\int_{E_{-1,1}}
\tau (Q_{v,y'^{-1}}(C_1)) (
\int_{E_{-1,1}}u(A) \tau( F_{v,y'^{-1}}(C_1,A)  )dA)
 dC_1\\
 &=& |F_k(v)|^{\frac{p_1 }{ d_k}}|F_{p-k}^*(y')|^{\frac{p_1 }{
2d'_k}}.\int_{E_{-1,1}}
\tau ( Q_{v,y'^{-1}}(C_1))\four_{E_{-1,1}}(u)(\theta_{v+y'^{-1}}(C_1))dC_1\\
&=&\gamma_k(v,y'^{-1})\ k(v)\quad \hbox{avec:}\\
k(v)&=& |F_k(v)|^{\frac{p_1}{2d_k}}  \int_{E_{-1,1}} \tau (- Q_{v,y'^{-1}}(C_1)\ )u(C_1)dC_1\ ,\\
 \end{array}$$
par application du lemme 4.2.3.\\
\noi Or:
$$\begin{array}{rll}k(v)&= &|F_k(v)|^{\frac{p_1 }{
2d_k}}.\int_{E{-1,1}}\tau (-Q_{v,y'}(A)) { 
\four_{E_{0,2}}}(f(v+[A,v]+\bullet))(y')dA\\
 &=&  |F_k(v)|^{\frac{p_1}{
2d_k}}.\int_{E_{-1,1}}\biggl(\int_{E_{0,2}}f(v+[A,v]+y)\tau({\tilde
B}(y,y')) \tau(-Q_{v,y'^{-1}}(A))dy\biggr)dA\\
 &=&
|F_k(v)|^{\frac{p_1 }{ 2d_k}}.\int_{E_{-1,1}}\biggl(\int_{E_{0,2}}f(
exp(adA)(v+y))\tau({\tilde B}(y,y'))dy\biggr)
 dA\end{array}$$ et comme il y a convergence absolue, on
peut intervertir l'ordre des int\égrations d'o\ù :
$$ \begin{array}{rll}k(v)&=& |F_k(v)|^{\frac{p_1 }{ 2d_k}}.\int_{E_{0,2}}\biggl(\int_{E_{-1,1}}f(exp(adA)(v+y)) dA\ \biggr)\ \tau( {\tilde B}(y,y'))dy \\
 &=&\four_{E_{0,2}}(T_f(v+.))(y').\hfill\Box\end{array}  $$
 
\begin{rema}  \begin{enumerate}
\item  Lorsque $E_{2,2}=\{0\},$  on a la forme
simple  du cas commutatif.\\
 
 \item Le th\éor\ème est encore vrai sans  irr\éductibilit\é
mais avec quelques am\énagements.  \\

\item Soit $U$ un ouvert de $\{x\in \goth g_1\ |\ F_k(x)F(x)\not=0\}$ 
et $f\in \EuScript C_C^{\infty}(U)$ alors pour $y\in E'_{0,-2}$ fix\é 
la fonction $\gamma_k(\ ,y^{-1})\four_{E_{0,2}}(T_f(x+\ ))(y)\in \EuScript C_C^{\infty}({(E'_{2,0})}_y)$ (cf.lemme 1.1.1 et ${(E'_{2,0})}_y$ est une r\éunion finie de ${(G_{h_k})}_y$-orbites ouvertes).\\

\noi Plus g\én\éralemnt, si $V$ un ouvert de $\{x\in \goth g_1\ |\ \prod_{1\≤i\≤k}F_i(x)F(x)\not=0\}$ 
et $f\in \EuScript C_C^{\infty}(V)$ alors pour $y\in E'_{0,-2}$ fix\é 
la fonction $\gamma_k(\ ,y^{-1})\four_{E_{0,2}}(T_f(x+\ ))(y)\in \EuScript C_C^{\infty}({(E''_{2,0})}_y)$ (cf.lemme 1.1.1 et ${(E''_{2,0})}_y$ est une r\éunion finie de ${ P(H_1,...,H_k)}_y$-orbites ouvertes).\\
\end{enumerate}\end{rema}

\bigskip

\noi Les r\ésultats obtenus dans le th\éor\èmes 4.3.3 sont reli\és \à l'action de $G_{h_k}$ sur les sous-vari\ét\és:
$$\begin{array}{rlllrl} {\tilde W}_k&=&\{x+y\ |x\in E'_{2,0}\ ,\ y\in E_{0,2}\ ,\ [x,y]=0\}  &\hbox{et}\ &W_k=&\{z\in {\tilde W}_k\ |\ F_k(z)F_p(z)\not=0\}\\
 {\tilde W}^*_k&=&\{x'+y'\ |x'\in E_{-2,0},y'\in E'_{0,-2},[x',y']=0\}  &\hbox{et}\ &W^*_k=&\{z\in {\tilde W}^*_k\ |\ F^*_{p-k}(z)F^*_p(z)\not=0\}.\end{array}$$\\

\subsection{  Expression des  mesures $\bf G_{h_k}$
invariantes sur $W_k$ et $W^*_k$}
\bigskip  
\noi $W_k$ et $W^*_k$ sont des r\éunions finies de $G_{h_k}-$orbites et on a:
\vskip 3mm
 \begin{prop}
 Soit $x_1+y_1\in W_k,$ l'orbite $U=G_{h_k}.(x_1+y_1)$ est
 ouverte dans ${ \tilde W}_k.$
Elle est munie d'une mesure  $\lambda,$
 $G_{ 
h_k}$-invariante et $\forall f\in L^1(U,\lambda)$ on a  :
$$ \begin{array}{rll}\int_Uf\ d\lambda\
&=& \int_{G_{h_k}.x_1}|F_k(x)|^{-\frac{p_{2,0}} {d_k}}\ (\
\int_{G_{h_k,x}.(g.y_1)}\ f(x+y)\ |P_{p-k}(y)|^{-\frac
{p_{0,2}-p_{2}} { d'_k}}.
  \ dy\ )\ dx \\
 &=&  \int_{G_{h_k}.y_1
}|P_{p-k}(y)|^{-\frac{p_{0,2} }{d'_k}}
\ (\ \int_{G_{h_k,y}.(g'.x_1)}\ f(x+y)\ |F_k(x)|^{-\frac
{p_{2,0}-p_{2}} { d_k}} 
  \ dx\ )\ dy \end{array}$$
$g$ (resp. $g'$) \étant un  \elt de $G_{h_k}$ tel que
$x=g(x_1).$
 (resp. 
 $y=g'(y_1 ))\ .$  
\end {prop}
\vskip 3mm
Les mesures sont choisies suivant le $\S 4.2.$
\vskip 3mm
 \dem
 1) L'espace tangent en $(x,y)$ \à $ W_k$ est donn\é par :
$$T_{(x,y)}( W_k)=\{(u,v)\in 
E_{2,0}\times E_{0,2}\ ,\ \hbox{ tels que }\ \
[x,v]+[u,y]=0\}$$  Comme  $x_1+y_1\in  W_k\subset \goth g'_1 $   on a $\text{ad}(x_1+y_1)(\goth g_0)=\goth g_1$ et $E_{0,0}=$Ker(ad$x_1)\cap$Ker(ad$y_1)\cap E_{0,0}\oplus $ad$(x_1)$ad$(y_1)(E_{-2,-2})$ donc  dim($\goth g_0)-$ dim($\goth g_1)=$
dim(Ker($\text{ad}(x_1+y_1)/\goth g_0)=$

\noi dim(Ker($\text{ad}x_1)\cap $Ker($\text{ad}y_1)/E_{0,0})+p_1+p_2$ d'o\ù 
dim($E_{0,0})$-dim(Ker($\text{ad}x_1)\cap $Ker($\text{ad}y_1)\cap E_{0,0})= p_{2,0}+p_{0,2}-p_2$ ainsi
  l'application de $E_{0,0}$ dans
$T_{(x_1,y_1)}(W_k)$ donn\ée par $A\rightarrow [A,x_1+y_1]$ est
surjective, donc l'application de $G_{h_k}$ dans $ W_k$ donn\ée
par
$g\rightarrow g(x_1+y_1)$ est submersive donc ouverte d'o\ù l'orbite
 $U$ est ouverte dans $  W_k$ par cons\équent elle est localement ferm\ée dans l'espace vectoriel topologique $E_{2,0}+E_{0,2}$ et  
hom\éomorphe \à
$G_{h_1}/_{G_1},$ $G_1$ \étant le centralisateur de  
$x_1+y_1 $ dans
$G_{h_k}$ (proposition $6, \S 5, n^\circ 3,$ chapitre IX de  \cite{bourbakitop}.
\vskip 2mm
2) On utilise les r\ésultats sur les mesures quasi-invariantes sur
les espaces homog\ènes.

\noi $G_{h_k} ,G_{{h_k},x_1},G_{{h_k},y_1},G_1$ sont des groupes
r\éductifs car ce sont 
les centralisateurs dans $Aut_0(\goth g)$  des sous-alg\èbres
suivantes, toutes
r\éductives dans
$\goth g$ : $$\F.H_0\oplus\F.h_k\ \ ,\ \ {\goth s_1}=
\F.h_k\oplus \F.H_0 \oplus 
\F.x_1
 \oplus \F.x_1^{-1}\ \ ,\ \ {\goth s_2}= \F.H_0\oplus
\F.h_k 
\oplus 
\F.y_1
 \oplus \F.y_1^{-1}\ ,\ {\goth s}={\goth s_1} +{\goth
s_2}$$ ainsi ce sont des groupes unimodulaires \footnote{ $G$ \étant un sous-groupe alg\ébrique r\éductif de $Gl_p(\F),$  sa composante connexe alg\ébrique, $G^0$, est le produit d'un groupe commutatif inclus dans le centre et du groupe d\ériv\é, l'intersection des $2$ \étant un sous-groupe fini (\cite{boreltits},prop.2.2 p.63),  donc $G^0$ est unimodulaire ( \cite{bourbakiint},Int\égration,chap.$7,\S 2, n^\circ 9$ prop. $14$) d'o\ù $G$ l'est aussi car  le groupe quotient de $G$ par $G^0$ est fini (b) prop.$10,\S 2, n^\circ 7$,m\ême r\éf. et \cite{chevalley},chap.2,$\S 3,n^\circ 3,$ th.2).} d'o\ù sur
$G_{h_k}/_{G_1}$ on a une mesure $G_{h_k}$ invariante \à
gauche, unique \à une constante multiplicative pr\ès (corollaire 2, $\S 2$ n$^\circ 6$ Int, chap.VII de \cite{bourbakiint}), ce qui
d\émontre le premier point.

\noi L'expression des   deux d\écompositions qui suivent proviennent des inclusions:
$$G_1\subset G_{{h_k},x_1}\subset G_{h_k}\qquad G_1\subset
G_{h_k,y_1}\subset G_{h_k}$$ et des d\écompositions
correspondantes des mesures.\\

\noi Soit $f\in L^1(U)$ alors $k(\dot g)=f(g(x_1+y_1))\in L^1(G_{h_k}/G_1),$ $G_{h_k}/G_1$ \étant muni de la mesure invariante \à gauche $d\dot g, $ $k_1(\dot g')=\int_{G_{{h_k},x_1}/G_1}f( g'g"(x_1+y_1))d\dot g"\in L^1(G_{h_k}/G_{{h_k},x_1}),$ $G_{h_k}/G_{{h_k},x_1}$ \étant muni de la mesure invariante \à gauche $d\dot g', $ et on a:
$$\int_{G_{h_k}/G_1}k(\dot g)d\dot g=\int_{G_{h_k}/G_{{h_k},x_1}}d\dot g'\biggl (
 \int_{G_{{h_k},x_1}/G_1}f( g'g"(x_1+y_1))d\dot g"\biggr )$$
(a) du corollaire 1, n$^\circ8,$ \S 2, chap
VII, \cite{bourbakiint}). \\

\noi  Posons $x=g'x_1,$ alors:
$$\begin{array}{ll} \int_{G_{{h_k},x_1}/G_1}f( g'g"(x_1+y_1))d\dot g" &=  \int_{G_{{h_k},x_1}/G_1}f( x+g'g"g'^{-1}(g'y_1))d\dot g"\\
&=\int_{G_{{h_k},x}/g'G_1g'^{-1}}f( x+ g"(g'y_1))d\dot g"\\
&=\int_{G_{h_k,x}.(g'.y_1)}\ f(x+y)\ |P_{p-k}(y)|^{-\frac
{p_{0,2}-p_{2,2}} { d'_k}}
  \ dy\ ,\end{array}$$ 
  par unicit\é (\à une constante multiplicative pr\ès) de la mesure invariante sur $G_{{h_k},x_1}/G_1$, la derni\ère \égalit\é est obtenue par consid\ération du \PV $((E_{0,0})_x,E_{0,2})_x)$ (cf.(N2) $\S3.2$).\\
  
  \noi De m\ême:
  $$\int_{G_{h_k}/G_{{h_k},x_1}}k_1(\dot g')d\dot g'=\int_{G_{h_k}x_1} k_1(\pi^{-1}(\dot g'))|F_k(x)|^{-\frac{p_{2,0}} {d_k}}dx\ ,$$
$\pi$ \étant la projection canonique de $G_{h_k}$ sur   $G_{h_k}/G_{{h_k},x_1},$ par unicit\é (\à une constante multiplicative pr\ès) de la mesure invariante sur $G_{h_k}/G_{{h_k},x_1} $ et  en consid\érant du \PV 
$(E_{0,0},E_{2,2}).$ On obtient ainsi la premi\ère \égalit\é.\\
 
  \noi On obtient la deuxi\ème relation en consid\érant les sous-groupes  $G_1 $ et 
$G_{h_k,y_1}$ de $ G_{h_k}.$\\

 3) Il reste \à donner la valeur de la constante  apparaissant
dans la deuxi\ème formule, ce qui d\écoule d'un simple calcul d'int\égrales.\\

\noi L'orbite  $U^*=G_{h_k}(x_1^{-1}+y_1^{-1})$  est \également
munie d'une mesure $G_{h_k}$-invariante (\à une constante
multiplicative pr\ès)  que l'on peut d\éfinir  par 
$$\int_{U^*}fd\lambda^*=\int_{U}f\circ\Phi d\lambda^.$$
Lorsque $h\in \EuScript C_K(U)$ et
$g\in  \EuScript C_K(\goth n_k),$ on d\éfinit la fonction $E(h,g)\in
 \EuScript C_K( \{x\in \goth g_1|F_p(x)F_k(x)\not=0\})$ par: 
$$E(h,g)(\Psi (A,x,y))=E(h,g)(exp(ad(A)(x+y))=h(x+y).g(A).$$

\noi On a la m\ême application, not\ée  $E^*$ , lorsque $k\in
\EuScript C_K(U^*)$ et
$g\in \EuScript C_K(\goth n_k)$ 
  et $E^*(k,g)\in \EuScript C_K( \{x\in \goth g_{-1}|F^*_p(x)F^*_{p-k}(x)\not=0\})$ (cf. lemme 1.1.1) et  on a:
$$\begin{array}{lrl}E(k\circ \Phi,g)&=&E^*(k,g)\circ \Phi\ ,\\
 S_{E(h,g) }/U&=&|2|^{\frac {p_2}{2}}C(g)|F_k(x)|^{md'_k}h\ \hbox{avec}\ C(g)= \int_{\goth n_k}g(u)du.\end{array}$$
\noi Soient $r_1=p_{2,0}-p_{2},$  $r_2=p_{0,2}-p_{2}$ rappelons que:
$$m_k=\frac{\frac{p_{1}}{2}+p_{2,0}}{Nd_k}=\frac{\frac{p_{1}}{2}+p_{0,2}}{Nd'_k} \ ,\ m=\frac{\frac{p_1}{2}+p_2}{d_kd'_k}\ ,\ d:=d_p=m_k(d_k+d'_k)$$
et que  pour $ x+y\in W_k:|F(x+y)|=\biggl(|F_k(x)||P_{p-k}(y)|\biggr)^{m_k}.$

\noi  Par le 1) du th\éor\ème 4.3.3 et en raison du support de $h,$ on a :
$$\int_{\goth g_1}E(h,g)(x)|F(x)|^{-\frac{r_2 }{
m_kd'_k}}.| F_k(x)|^{- \frac{md}{ m_k }}dx= 
 |2|^{\frac{ p_{2} }{ 2}}C(g)\int_Uhd\lambda\  .$$
On suppose $C(g)$ non nul alors pour $k\in S(U^*)$ on  a:
$$\begin{array}{rll}\int_{U^*}kd\lambda^*&=&\int_{U}k\circ \Phi d\lambda\\
&=&  |2|^{ -\frac{ p_{2} }{ 2}}C(g)^{-1}\int_{\goth g_1}E(k\circ \Phi,g)(x)|F(x)|^{-\frac{r_2 }{
m_kd'_k}}.| F_k(x)|^{- \frac{md}{ m_k }}dx\\
&=&|2|^{ -\frac{ p_{2} }{ 2}}C(g)^{-1} \int_{\goth
g_{-1}}E^*(k,g)(y)|F^*(y)|^{-\frac{r_1 }{ m_kd_k}}| 
F^*_{p-k}|^{-\frac{md}{m_k}}dy\end{array}$$
en utilisant le 3) du lemme 3.2.2 ainsi que les relations 
 suivantes pour $x\in  \{x\in \goth g_{1}|F_p(x) \not=0\}$ et $y\in  \{x\in \goth g_{-1}|F^*_p(y)F^*_{p-k}(y)\not=0\}$ (cf.(R2) $\S 3.3$):
$$F^*\circ\Phi(x)=F(x)^{-1}\qquad | 
F_k\circ\Phi^{-1}(y)|^{m_k}=|F^*(y)|^{-1}.| 
F^*_{p-k}(y)|^{m_k}.$$ 
Par le 2) du th\éor\ème 4.3.3 et par le choix du support de $k,$  on obtient:
$$ \int_{U^*}kd\lambda^*=  
\int_{G_{h_k}.y_1^{-1}}|F^*_{p-k}(y')|^{-\frac{p_{0,2} }{
d'_k}}
\biggl(\int_{G_{h_k,y'}. (g.x_1^{-1})}k(x'+y')|P^*_k(x')|^{-\frac{r_1 }{
d_k}}dx'\biggr)dy'.  $$
D'o\ù
$$\begin{array}{rll}\int_{U}hd\lambda&=&\int_{U^*}h\circ\Phi^{-1}d\lambda^* \\
&=&\int_{G_{h_k}.y_1^{-1}}|F^*_{p-k}(y')|^{-\frac{p_{0,2} }{d'_k}}
(\int_{G_{h_k,y'}.(g.x_1^{-1})}h(x'^{-1}+y'^{-1})|P^*_k(x')|^{-\frac{r_1 }{
d_k}}dx')dy'\\
&= &\int_{G_{h_k}.y_1}|P_{p-k}(y')|^{-\frac{p_{0,2} }{
d'_k}}
(\int_{G_{h_k,y}.(g.x_1)}h(x+y)|F_k(x)|^{-\frac{r_1 }{
d_k}}dx)dy \end{array}$$
par le choix des mesures adapt\ées \à $F_k$ et $P_{p-k}$ dans la d\éfinition 4.2.5  (cf.3) du lemme 3.2.2 ainsi que (R2) du $\S 3.3$).

\fdem\\

\noi  On en d\éduit les r\ésultats suivants:\\

\begin{cor} $W_k$ et $W^*_k$ sont munis de mesures $G_{h_k}$-invariantes: $d\lambda$ et $d\lambda^*.$
\begin{enumerate}
\item Pour $ f\in L^1(W_k)$ on a  :
$$  \begin{array}{rll}\int_{W_k}f\ d\lambda\
&=& \int_{ E'_{2,0}}|F_k(x)|^{-\frac{p_{2,0}} {d_k}}\biggl(\
\int_{ (E'_{0,2})_x}\ f(x+y)\ |P_{p-k}(y)|^{-\frac
{p_{0,2}-p_{2}} { d'_k}}
  \ dy\ \biggr)\ dx \\
 &=&  \int_{ E'_{0,2}
}|P_{p-k}(y)|^{-\frac{p_{0,2} }{d'_k}}
\biggl(\ \int_{ (E'_{2,0})_y}\ f(x+y)\ |F_k(x)|^{-\frac
{p_{2,0}-p_{2}} { d_k}} 
  \ dx\biggr)\ dy \end{array}$$
\item  Pour $ g\in L^1(W^*_k)$ on a  :
$$  \begin{array}{rll}\int_{W^*_k}g\ d\lambda^*\
&=& \int_{ E'_{-2,0}}|P^*_k(x)|^{-\frac{p_{2,0}} {d_k}}\biggl(\
\int_{ (E'_{0,-2})_x}\ g(x+y)\ |F^*_{p-k}(y)|^{-\frac
{p_{0,2}-p_{2}} { d'_k}}
  \ dy\ \biggr)\ dx \\
 &=&  \int_{ E'_{0,-2}
}|F^*_{p-k}(y)|^{-\frac{p_{0,2} }{d'_k}}
\biggl(\ \int_{ (E'_{-2,0})_y}\ g(x+y)\ |P^*_k(x)|^{-\frac
{p_{2,0}-p_{2}} { d_k}} 
  \ dx\biggr)\ dy \end{array}$$
\item Les r\ésultats sont analogues sur :
$$W'_k=\{x+y\ |x\in E'_{2,0}\ ,\ y\in E'_{0,-2}\ ,\ [x,y]=0\}\ ,\ \{x+y\ |x\in E'_{-2,0}\ ,\ y\in E'_{0,2}\ ,\ [x,y]=0\}.$$  
  \end{enumerate}

\end{cor}

\bigskip
\noi Le th\éor\ème 4.3.3 devient:\\
\begin{cor}
\begin{enumerate}
\item Soit $f$  dans  $L^1(\goth g_1),$   $S_f$  est d\éfini presque partout sur $(W_k,d\lambda)$  et  $|F_k(\ )  |^{md_k}|F(\ )|^{\frac{ p_{0,2} -p_{2}}{ m_kd'_k}}S_f $ est int\égrable  sur  $(W_k,d\lambda)$  et on a:
$$\begin{array}{rll}\int_{\goth g_1}f(x)dx&=&\int_{W_k} |F_k(u)  |^{md_k}|F(u)|^{\frac{ p_{0,2} -p_{2}}{ m_kd'_k}}S_f (u)d\lambda (u)\ ,\\ Z(f;\pi)&=&\int_{W_k} |F_k(u)  |^{md_k}|F(u)|^{\frac{ p_{0,2} -p_{2}}{ m_kd'_k}}\pi (F(u))S_f (u)d\lambda (u)\ \hbox{pour}\ \Re(\pi)>0.\end{array}$$
Soit O une orbite ouverte de $P_{\goth t}$ dans $\goth g_1$ alors:
$$Z_O(f;\pi)=\int_{O\cap W_k} |F_k(u)  |^{md_k}|F(u)|^{\frac{ p_{0,2} -p_{2}}{ m_kd'_k}}\pi (F(u))S_f (u)d\lambda (u)\ \hbox{pour}\ \Re(\pi)>0.$$
\item Soit $g$  dans  $L^1(\goth g_{-1}),$   $S^*_g$  est d\éfini presque partout sur $(W^*_k,d\lambda^*)$  et  $|F^*_{p-k}(\ )  |^{md'_k}|F^*(\ )|^{\frac{ p_{2,0} -p_{2}}{ m_kd_k}}S^*_g$ est int\égrable  sur  $(W^*_k,d\lambda^*)$  et on a:
$$\begin{array}{rll}\int_{\goth g_{-1}}g(x)dx&=&\int_{W^*_k}  |F^*_{p-k}(u)  |^{md'_k}|F^*(u)|^{\frac{ p_{2,0} -p_{2}}{ m_kd_k}}S^*_g(u)d\lambda^* (u)\\
Z^*(g;\pi)&=&\int_{W^*_k}  |F^*_{p-k}(u)  |^{md'_k}|F^*(u)|^{\frac{ p_{2,0} -p_{2}}{ m_kd_k}}\pi(F^*(u))S^*_g(u)d\lambda^* (u)\ \hbox{pour}\ \Re(\pi)>0.\end{array}$$
Soit $O^*$ une orbite ouverte de $P_{\goth t}$ dans $\goth g_{-1}$ alors:
$$Z_{O^*}^*(g;\pi)=\int_{O^*\cap W^*_k}  |F^*_{p-k}(u)  |^{md'_k}|F^*(u)|^{\frac{ p_{2,0} -p_{2}}{ m_kd_k}}\pi(F^*(u))S^*_g(u)d\lambda^* (u)\ \hbox{pour}\ \Re(\pi)>0.$$
\end{enumerate}
\end{cor}
\bigskip

\noi Soient :
$$\begin{array}{rllll}N_{k,0}&=&exp\bigl(ad(\goth n_{k,0}) \bigr) \hbox{avec}\ \goth n_{k,0}&=&\oplus_{k+1\≤i<j\≤p}(E_{i,j}(-1,1)\oplus E_{i,j}(-2,2))\ ,\\
N_{k,2}&=&exp\bigl(ad(\goth n_{k,2} )\bigr)\ \hbox{avec}\ \goth n_{k,2}&=&\oplus_{1\≤i<j\≤k}(E_{i,j}(-1,1)\oplus E_{i,j}(-2,2)) \ ,\end{array}$$
comme $\goth n_{k,0}$ et $\goth n_{k,2}$ commutent,  il en est de m\ême pour   $N_{k,0}$ et $N_{k,2}.$  

\noi $N_{k,0}$ centralise $E_{\pm 2,0}$ et $N_{k,2}$ centralise $E_{0,\pm 2}.$ \\

\noi $x_i,t_i,i=0$ ou $2$ sont des \elts de $E_{i,2-i}$ et $y_i,z_i,i=0$ ou $-2$ sont des \elts de $E_{i,-2-i}.$ 
\begin{rema}

 \begin{enumerate}
\item Soit O une orbite ouverte de $P_{\goth t}$ dans $\goth g_1$ et $t_2+t_0\in O\cap W_k$ alors $O\cap W_k=G_{\goth t} N_{k,0}N_{k,2}.(t_2+t_0)$ et 
$$x_2+x_0\in O\cap W_k\Leftrightarrow\ x_2\in G_{\goth t}N_{k,2}t_2 \ \hbox{et}\ x_0\in \bigl(G_{\goth t} \bigr)_{x_2}N_{k,0}(gt_0)\ ,\ \hbox{avec}\ x_2=gt_2.$$

\item Soit $O^*$ une orbite ouverte de $P_{\goth t}$ dans $\goth g_{-1}$ et $z_{-2}+z_0\in O^*\cap W^*_k$ alors $ O^*\cap W^*_k=G_{\goth t} N_{k,0}N_{k,2}.(z_{-2}+z_0)$ et 

$$y_{-2}+y_0\in O^*\cap W^*_k\Leftrightarrow\ y_0\in G_{\goth t}N_{k,0}z_0 \ \hbox{et}\ y_{-2}\in \bigl(G_{\goth t} \bigr)_{y_0}N_{k,2}(g'z_{-2})\ ,\ \hbox{avec}\ y_0=g'z_0.$$

\item Soit $t_2+z_0\in W'_k$ alors $x_2+x_0\in G_{\goth t} N_{k,0}N_{k,2}.(t_{2}+z_0)$ \ssi $x_2\in G_{\goth t}N_{k,2}t_2$  et $ x_0\in \bigl(G_{\goth t} \bigr)_{x_2}N_{k,0}(gz_{0}) $ avec $ x_2=gt_2$ \ssi $x_0\in G_{\goth t}N_{k,0}z_0 $ et $ x_2\in \bigl(G_{\goth t} \bigr)_{x_0}N_{k,2}(g' t_2) $ avec $x_0=g'z_0.$
\end{enumerate}
\end{rema} 
\bigskip

\noi En effet  dans le \PV faiblement sph\érique $(G_{\goth t}N_{k,2},E_{2,0}),$ l'orbite contenant $x_2$ rencontre $W_{\goth t'}$ avec $\goth t'=\oplus_{1\≤i\≤k}\F H_i$ et soit $x'_2=gn.x_2\in   W_{\goth t'},$ avec $g\in G_{\goth t}$ et $n\in N_{k,2}.$ Dans le \PV faiblement sph\érique $((G_{\goth t})_{x'_2}N_{k,0},(E_{0,2})_{x'_2}),$ l'orbite contenant $gnx_0=gx_0$ rencontre $W_{\goth t"}$ avec $\goth t"=\oplus_{k+1\≤i\≤p}\F H_i$  donc il existe $g'\in (G_{\goth t})_{x'_2}N_{k,0}$ tel que $g'gx_0\in W_{\goth t"}$ d'o\ù $g'gn(x_2+x_0)\in W_{\goth t}.$ On proc\ède de m\ême pour $t_0+t_2:$ il existe $g"\in G_{\goth t} N_{k,0}N_{k,2}$ tel que $g"(t_0+t_2)\in W_{\goth t};$ comme $g"(t_0+t_2)$ et $g'gn(x_2+x_0)$ sont dans $W_{\goth t}$ et dans la m\ême orbite de $P_{\goth t},$ ils sont dans la m\ême orbite de $G_{\goth t}$ par le lemme 1.4.4 d'o\ù $1.$

\noi Il en est de m\ême pour $2.$
\newpage
 \section{ Application aux fonctions Z\étas}

\subsection{Les coefficients de l'\équation fonctionnelle}
\bigskip 

On utilise les d\écompositions \établies dans le $\S 4$ relativement \à $ad(h_k)$ pour obtenir les 
coefficients de l'\équation fonctionnelle v\érifi\ée par les fonctions Z\étas par application \à $2$ reprises des \équations fonctionnelles associ\ées \à des centralisateurs de \Sls .  

\noi Les notations sont celles du $\S 4.$\\

Soit $O^*$ une orbite de $P_{\goth t}$
 dans $\goth g"_{-1},$  prenons $z\in O^*\cap W^*_k$
et soit $z=z_{-2}+z_0$ sa d\écomposition suivant $ad(h_k),
z_i\in E_i(h_k)\cap \goth g_{-1}.$ \\

On d\ésigne par
$t_i,i=1,...,l$ un ensemble de repr\ésentants des orbites de
$(G_{\goth t})_{z_0}N_{k,2}$ dans ${(E"_{2,0})}_{z_0},$ et pour chaque $t_i:$ $\{t_{i,j},j=1,...,p_i\}$
est un ensemble de repr\ésentants des orbites de $ (G_{\goth t})_{t_i}N_{k,0}$
dans
${(E"_{0,2})}_{t_i}.$ \\

 Comme $t_i+t_{i,j}$ et $t_i+t_{i,k},$  
$j\not=k,$ ne sont pas dans la m\ême orbite de $P_{\goth t}$ dans
$\goth g_1$ (cf.remarque 4.4.4), pour chaque orbite $O$ de $P_{\goth t}$ dans
$\goth g"_1$ et chaque $t_i$ il existe au plus une
valeur $j$ telle que $P_{\goth t}(t_i+t_{i,j})=O,$ lorsqu'elle
existe on  notera
  $t_{i(O)}$ l'\elt correspondant et $I_{O,O^*}=\{i\ |\
\exists j$ tel que $t_i+t_{i,j}\in O\},$   sinon 
$I_{O,O^*}=\emptyset.$

Lorsque $I_{O,O^*}\not=\emptyset,$ $\{t_i,i\in
I_{O,O^*}\}$ est un ensemble de
repr\ésentants des orbites de
$(G_{\goth t})_{z_0}N_{k,2}$ dans l'ouvert
non vide $(\pi_{h_k}(O))_{z_0}$ (notation du lemme 1.1.1) et 
on a  le sch\éma suivant pour les diff\érents repr\ésentants:
$$\begin{array}{lrlll}
\hbox{Pour}\quad O^*&:&\ z_{-2}+z_0&\in &W^*_k\cap O^*,\\
\hbox{Pour}\quad O&:&\ t_i+t_{i(0)}&\in &W_k\cap O,\\
\hbox{avec}\quad &:&\ z_{-2}\  \hbox{et}\ t_i\ \hbox{commutant avec le}\ sl_2-\hbox{triplet}\ &:&(z_0^{-1},2H_0-h_k,z_0 ),\\
\hbox{ et}\quad &:&\ z_0\ \hbox{et}\ t_{i(0)}\ \hbox{commutant avec le}\  sl_2-\hbox{triplet}\ &:&(t_i, h_k,t_i^{-1})\\
\hbox{ainsi que}&:&t_i+z_0^{-1}\in W_k
.\end{array}$$
La transformation de Fourier
\étant d\éfinie \à partir de la restriction de $\tilde B,$  
les diff\érents coefficients de l'\équation fonctionnelle de
la Fonction Z\éta  sont
alors reli\és par une relation donn\ée la proposition 5.1.1.\\

\noi Notations pr\éliminaires:\\

$\bullet$ $a_{ z_{-2},t_i}^{(z_0)}$ d\ésigne le coefficient associ\é aux orbites $((G_{\goth t})_{z_0}N_{k,2}.z_{-2}$ et $(G_{\goth t})_{z_0}N_{k,2}.t_i$ dans le \PV: $( G_{\goth t})_{z_0}N_{k,2},(E_{2,0})_{z_0})$ et  le \PV  dual  $(( G_{\goth t})_{z_0}N_{k,2},(E_{-2,0})_{z_0}),$\\

$\bullet$ $a_{ z_{0},t_{i(O)}}^{(t_i)}$ d\ésigne le coefficient associ\é aux orbites $((G_{\goth t})_{t_i}N_{k,0}.z_{0}$ et $(G_{\goth t})_{t_i}N_{k,0}. t_{i(O)}$ dans le \PV: $(( G_{\goth t})_{t_i}N_{k,0},(E_{0,2})_{t_i})$ et  le \PV  dual  $(( G_{\goth t})_{t_i}N_{k,0},(E_{0,-2})_{t_i}).$

\begin{prop}
Lorsque $ I_{O,O^*}=\emptyset$ on a
$a_{O^*,O}(\pi)=0$ et sinon
 
$$a_{O^*,O}(\pi)= \sum_{ i\in I_{O,O^*}  
 } a_{z_{-2},t_i}^{(z_0)}(\pi_{p-k+1},...,\pi_p).\gamma_k
(t_i,z_0^{-1}).a_{z_0,t_{i(O)}}^{(t_i)}(
\pi_1,...,\pi_{p-k-1},\pi_{p-k}^{\frac{1}{m_k}}\prod_{p-k+1\≤j\≤p}\pi_j|\ |^{\frac{r_k}{m_k}}) $$  $\displaystyle{ (r_k =\frac{  p_2+\frac{p_1}{ 2}} { d'_k}}$ et $m_k=1$ sauf  lorsque   $(\∆,\lambda_0)$ est une   $\F$-forme
de
$(E_7,\alpha_6)$ ou bien  de type $(C_n,\alpha_k)$ avec $\goth t =\goth t'_0,$   dans ce cas $p=2,k=1$ et $m_1=\frac{1}{2})$

\end{prop}

\dem  Soit $f\in \EuScript C_c^{\infty}(O)$ et $\pi$ un caract\ère tel que
$Re(\pi)>0.$  

\noi Comme $Re(\pi)>0,$ 
 $ Z^*_{O^*}({ \four ( f)};\pi)$ est donn\é sous forme
int\égrale  donc par application  du 2. du  corollaire 4.4.3 on a :
$$\begin{array}{rll} Z^*_{O^*}( \four( f);\pi) 
 &=&\int_{O^*\cap W^*_k} |F^*_{p-k}(u)  |^{md'_k}|F_p^*(u)|^{\frac{ p_{2,0} -p_2}{ m_kd_k}}\pi(F^*(u))S^*_{\four (f)}(u)d\lambda^* (u)\\
&=&   \int_{O(z_0)} |F^*_{p-k}(y')) |^{md'_k-\frac{p_{0,2}}{d'_k}}I(y')dy'\ \hbox{avec}\\
I(y')&=&
 \int_ {O'}  
 S^*_{\four( f)}(x'+y') \pi (F^*(x'+y'))
  |P_k^*(x')|^{-\frac{p_{2,0}-p_2}{d_k}}|F^*_p(x'+y')|^{\frac{ p_{2,0} -p_2}{ m_kd_k}}dx' \end{array} $$
avec $m=\ds\frac{r_k}{d_k},$ en appliquant le 2 du corollaire 4.4.2,
avec $O(z_0)=G_{\goth t}N_{k,0}.z_0$ et  $O'=(G_{\goth t})_{y'}N_{k,2}(gz_{-2})\subset
(E_{-2,0})_{y'},$ $g$ \étant un  \elt de $G_{\goth t}.N_{k,0}$
 tel
que $gz_0=y'$ (cf.2.remarque 4.4.4).\\

1) On suppose que $(\∆,\lambda_0)$ n'est pas une $\F$-forme
de $(E_7,\alpha_6)$ et que $\goth t\not=\goth t'_0$  lorsque   $(\∆,\lambda_0)$ est de type $(C_n,\alpha_k)$  c'est \à dire qu'on suppose que  $m_k=1.$\\  

\noi Le \PV faiblement sph\érique $( {(G_{\goth t})}_{y'}.N_{k,2},{(E_{2,0})}_{y'})$ 
est muni des invariants relatifs fondamentaux $G_j ,$ restrictions de $F_j,$ pour $j=1,...,k,$ donc les  invariants relatifs fondamentaux du \PV $( {(G_{\goth t})}_{y'}.N_{k,2},{(E_{-2,0})}_{y'})$ sont donn\és par:
$$j=1,...,k\ ,\ x'\in {(E_{-2,0})}_{y'}\quad :\quad  G^*_j(x')=\frac{F^*_{p-k+j}(x'+y')}{F^*_{p-k}(y')}$$
(cf.\demo de la proposition 3.4.4) donc $P^*_k=G^*_k$ et:
$$\prod_{1\≤i\≤p}\pi_i(F^*_i(x'+y'))=\prod_{1\≤i\≤p-k-1}\pi_i(F^*_i(y'))
(\prod_{p-k\≤i\≤p}\pi_i)(F^*_{p-k}(y'))\prod_{1\≤j\≤k}\pi_{p-k+j}(G^*_j(x'))  ,$$
soient $C(y')=\prod_{i=1}^{i= p-k-1}\pi_i(F^*_i(y'))
(\prod_{i=p-k}^{i= p}\pi_i)(F^*_{p-k}(y'))$ et $\pi'=(\pi_{p-k+1},...,\pi_p)$  alors:
$$I(y')=  C(y')|F^*_{p-k}(y')|^{\frac{p_{2,0}-p_2}{d_k}}Z^*_{O'}(S^*_{\four (f)}(\ +y' );\pi') \ ,$$
$Z^*_{O'}$ \étant la fonction Z\éta associ\ée \à l'orbite $O'$ du \PV $( {(G_{\goth t})}_{y'}.N_{k,2},{(E_{-2,0})}_{y'}).$\\ 
 
Posons
$h_f(x,y')=\four_{E_{0,2}}(T_f(x+\ ))(y'),$ alors pour $y'$ fix\é $h_f(\ ,y')\in \EuScript C_C^{\infty}({(E"_{2,0})}_{y'})$ (cf.3.remarque 4.3.6)  et pour
$i=1,...,l$ notons
  $O_i= (G_{\goth
t})_{y'}N_{k,2}.g(t_i),$ alors   par le  
th\éor\ème 4.3.5 on a :
$$ \begin{array}{rll} S^*_{\four ( f)}(x'+y') &=& 
\four_{(E_{2,0})_{y'}}(\gamma_k(\ ,y'^{-1})h_f(\ ,y')\
)(x')\\
 &=& \sum_{1\≤i\≤l}\gamma_k(t_i,z_0^{-1})\four_{(E_{2,0})_{y'}}(1_{ O_i}.h_f(\ ,y')\
)(x')\quad \hbox{    
d'o\ù}\\
   |F^*_{p-k}(y')) |^{md'_k-\frac{p_{0,2}}{d'_k}} I(y')&=&\sum_{1\≤i\≤l}\gamma_k(t_i,z_0^{-1}) C(y')|F^*_{p-k}(y')|^{\frac{p_1}{2d'_k}}I_i(y')
   \ \ \hbox{avec}\\
     I_i(y')&=&Z^*_{ O'}(\
\four_{(E_{2,0})_{y'}}(1_{ O_i}.h_f(\ ,y'))\ ;\pi'),\end{array}$$
en utilisant l'\égalit\é: $\displaystyle{\frac{\frac{p_1}{2}+p_{2,0}}{d_k}=
\frac{\frac{p_1}{2}+p_{0,2}}{d'_k}}$ (d\ém. du lemme 1.4.7).
\\
\noi  On applique l'\équation fonctionnelle dans le \PV faiblement sph\érique $( {(G_{\goth t})}_{y'}.N_{k,2},$ ${(E_{2,0})}_{y'})$ ainsi que 1. de la remarque 3.5.4 puisque les alg\èbres de Lie $\goth U(\F y'\oplus \F(2H_0-h_k)\oplus \F y'^{-1})$ et  $\goth U(\F z_0\oplus \F(2H_0-h_k)\oplus \F z_0^{-1})$ sont isomorphes, l'isomorphisme \étant donn\é par un \elt de $G_{\goth t}.N_{k,0}$ qui centralise $\oplus_{1\≤i\≤k}\F H_i,$ ce qui donne:$$Z^*_{ O'}(\
\four_{(E_{2,0})_{y'}}(1_{ O_i}.h_f(\ ,y'))\ ;\pi')=a_{z_{-2},t_i}^{(z_0)}(\pi')Z_{O_i}(h_f(\ ,y');\pi'^*|\ |^{-\frac{p_{2,0}-p_2}{d_k}{\bf 1}_k}).$$
\noi Comme $h_f(\ ,y')\in \EuScript C_C^{\infty}({(E"_{2,0})}_{y'})$ on a:
$$Z_{O_i}(h_f(\ ,y');\pi'^{*}|\ |^{-\frac{p_{2,0}-p_2}{d_k}})=\int_{O_i}
h_f(x ,y')\pi'^*(G(x))|G_k(x) |^{-\frac{p_{2,0}-p_2}{d_k}}dx$$
d'o\ù:
$$ Z^*_{O^*}({ \four ( f)};\pi)=\sum_{1\≤i\≤l}\gamma_k(t_i,z_0^{-1})a_{z_{-2},t_i}^{(z_0)}(\pi')J_i$$
avec:
$$J_i=\int_{O(z_0)} |F^*_{p-k}(y')) |^{ \frac{p_1}{2d'_k}}\ C(y')\biggl(\int_{O_i}
h_f(x ,y') \pi'^*(G(x)) |F_k(x) |^{-\frac{p_{2,0}-p_2}{d_k}}dx\biggr)dy'\ ,$$
et on a:
$$\pi'^*(G(x))=\prod_{j=1}^{j=k-1}\pi_{p-j} (F_j(x))\biggl(\prod_{p-k+1}^p\pi_j\biggr)^{-1}(F_k(x)).$$

\noi Soit   $V=G_{\goth t}.N_{k,0}.N_{k,2}(t_i+z_0),$  $V$ est ouvert dans $W'_k$ et $x+y'\in V$ \ssi $y'\in O(z_0)$ et $x\in O_i$ \ssi $x\in 
G_{\goth t}.N_{k,2}.t_i$ (not\é $O(t_i)$) et $y'\in {(G_{\goth t})}_{x}.N_{k,0}(g'z_0)$ avec $x=g't_i$ (cf.remarque 4.4.4) et qui sera not\ée $O(x).$  Soit $g$ la fonction mesurable sur $W'_k$
d\éfinie par:
$$g(x+y')=h_f(x,y') \pi'^*(G(x)) C(y')  |F^*_{p-k}(y') |^{\frac{\frac{p_1}{2}+p_{0,2}}{d'_k}} 
 {\bf 1}_V(x+y'),$$ alors 
par le 3. du corollaire 4.4.2  et de la remarque 4.4.4 on a:
$$ \int_{W'_k} gd\lambda '=\int_{O(t_i)}
   \pi'^*(G(x)) |F_k(x)|^{-\frac{p_{2,0}}{d_k}}   \biggl( \int_{O(x)}h_f(x,y') C(y') |F^*_{p-k}(y') |^ {md_k}dy'\biggr)dx\ .$$

Comme $f\in \EuScript C_C^{\infty}(O),$ $T_f\in  \EuScript C_C^{\infty}(W_k)$
donc $T_f\in \EuScript S(W_k)$ et est \à support compact dans
$E"_{2,0} $ donc relativement \à la variable $x$,
comme $Re(\pi)>0$  l' int\égrale double ci-dessus converge
absolument d'o\ù $g\in L^1(W'_k)$ et $J_i=\int_{W'_k} gd\lambda '.$\\
 
 \noi  On consid\ère maintenant les \PVs faiblement sph\ériques : $(G_{\goth t})_x.N_{k,0},(E_{0,\pm 2})_x,H_0-\frac{h_k}{2} )$ ($x=g'.t_i$). $((G_{\goth t})_xN_{k,0},  (E_{0,-2})_x)$   est muni des invariants relatifs fondamentaux
 provenant du \PV: $(P_{\goth t},\goth g_1)$ donn\és par la restriction
de $F^*_1,..,F^*_{p-k}$ et not\és $H^*_1,...,H^*_{p-k}.$ 

\noi Posons $\pi"=(\pi_ 1,...,\pi_{p-k-1},\prod_{j=p-k }^{j=p}\pi_j|\ |^{md_k} )$ et soit:
$$ \begin{array}{rll}D(x)&=&\int_{O(x)}h_f(x,y') C(y') |F^*_{p-k}(y') |^ {md_k}dy'\\
 &=&\int_{O(x)}\four_{(E_{0,2})_x}(S_f(x+\ ))(y')\pi"(H^*(y')) dy'\quad(1\text{ du th.}4.3.5)\\
 &=&  Z^*_{O(x)}(\four_{(E_{0,2})_x}(S_f(x+\ ) );\pi"  )\end{array}$$
 
 \noi  On applique \à nouveau l'\équation fonctionnelle
  dans  ce \PV faiblement sph\érique en tenant compte du 1. de la remarque 3.5.4 (conjugaison par $g'\in G_{\goth t}N_{k,2}$ qui centralise $\oplus_{j=k+1}^p\F H_i$) donc:
 $$\begin{array}{rll}D(x)&=&   \sum_{j} a^{(t_i)}_{z_0, t_{i,j}} 
 (  \pi ")Z_{i,j}(x)\ \hbox{avec}\\
Z_{i,j}(x)&=&Z_{ O(x,t_{i_j}) }(S_f(x+\ );   \pi"^* |\ |^
 {-\frac{p_{0,2}-p_{2}} {d'_k}{\bf 1}_{p-k}}) \ \hbox{et}\\
 O(x,t_{i_j})&=&{(G_{\goth t})}_x.N_{k,0}(g't_{i,j}).  \end{array}$$
 Comme le \PV $(G_{\goth t})_x.N_{k,0},(E_{0,-2})_x,H_0-\frac{h_k}{2} )$ est muni des \irfs $H^*_j,$ restriction de $F^*_j,j=1,...,p-k,$  les \irfs associ\és dans le 
 \PV $(G_{\goth t})_x.N_{k,0},(E_{0,2})_x,H_0-\frac{h_k}{2} )$ sont donn\és par $H_j(y)=\ds\frac{F_{j+k}(x+y)}{F_k(x)}$ pour $j=1,...,p-k$ (cf.$\S 3.3$). 
 
 \noi Comme $f\in \EuScript C_C^{\infty}(  O),$   $S_f\in \EuScript
C_C^{\infty}(W_k)$ donc $Z_{i,j }(x)$ est donn\é par
l'int\égrale: 
 $$ Z_{i,j}(x)=\int_{ O(x,t_{i,j})}S_f(x+y)(\prod_{1\≤i\≤p}\pi_i)^{-1}(F_p(x+y))(\prod_{p-k\≤i\≤p}\pi_i) (F_k(x))\prod_{k+1\≤j\≤p-1}(\pi_{p-j} (F_j(x+y)).$$
 $$\biggl(\frac{|F_k(x)|}{|F_p(x+y)|}\biggr)^{md_k+\frac{p_{0,2}-p_2}{d'_k}} dy$$
 d'o\ù
 $$J_i=\sum_ja^{(t_i)}_{z_0,t_{ij}}(\pi")\int_{W_k\cap P_{\goth t}(t_i+t_{i,j})}S_f(u)\pi^*(F(u))|F_p(u)|^{-N}|F_k(u)|^{md_k}|F_p(u)|^{\frac{p_{0,2}-p_2}{m_kd'_k}}d\lambda(u)$$
 puisque $N=\displaystyle{\frac{\frac{p_1}{2}+p_{0,2}}{d'_k}}$ (d\ém.du lemme 1.4.7) donc:
 $$\begin{array}{rll}J_i&=&\sum_ja^{(t_i)}_{z_0,t_{i,j}}(\pi")Z_{P_{\goth t}(t_i+t_{i,j})}(f;\pi^*|\ |^{-N{\bf 1}_p})\quad (1.\ \hbox{corollaire}\ 4.4.3)\\
 \\
 &=& \biggl\{{\begin{array}{llllll}&a^{(t_i)}_{z_0,t_i( O)}(\pi")Z_O(f;\pi^*|\ |^{-N{\bf 1}_p})&
 \ &\hbox{lorsque}&\ i\in I(O,O^*)\\
& 0&\ &\hbox{sinon}&
 \end{array}}
 \end{array}$$ 
par le choix du support de $f,$ d'o\ù le r\ésultat d'abord pour $\Re (\pi)>0$ puis en g\én\éral par prolongement m\éromorphe.\\

2) Dans les $2$ cas restants on a $p=2$ donc $k=1$ et $m_1=m'_1=\frac{1}{2}.$

\noi Lorsque $y'\in E'_{0,-2},$ le \PV
$(G_{H_1})_{y'},  (E_{2,0})_{y'})$ est munit de l'invariant relatif
fondamental $G(x)=F(x+y'^{-1})$ donc $G^*(x')=F^*(x'+y');$ de
m\ême pour $x\in E'_{2,0},$ le \PV $(G_{H_1})_x, 
(E_{0,2})_x)$ est muni de l'invariant relatif
fondamental $H(y)=F(x+y )$ donc $H^*(y')=F^*(x^{-1}+y').$ 

\noi On reprend les calculs pr\éc\édents, en prenant garde que les
mesures donn\ées dans les th\éor\ème 4.3.3,4.3.5 et  dans le $\S
4.4$ sont donn\ées pour
$F_1$ restreint \à $(E_{2,0})_{y'}$ et $P_1$ restreint \à
$(E_{0,2})_x$ ce qui apporte \à chaque fois une constante
lorsqu'on remplace les fonctions Z\éta  relatives \à 
$(E_{2,0})_{y'}$ et \à $(E_{0,2})_x$ par leurs valeurs mais
ce qui n'est pas le cas pour celles du dual
(compensation avec la transformation de Fourier).\fdem\\
 
\bigskip

\begin{rema} De la m\ême mani\ère on d\étermine $a^*_{O,O^*}(\pi)$ en prenant $x=x_2+x_0\in W_k\cap O$ $(x_i\in E_i(h_k)\cap \goth g_1)$ puis un ensemble de repr\ésentants des orbites de ${(G_{\goth t})}_{x_2}$ dans ${(E"_{0,-2})}_{x_2},$ not\é $\{u_i\}$  et pour chaque $u_i$  un ensemble de repr\ésentants des orbites de ${(G_{\goth t})}_{u_i}$ dans ${(E"_{-2,0})}_{u_i},$ not\é $(u_{i,j})_{1\≤j\≤r_i}.$ 
\noi On pose $I^*_{O,O^*}=\{i\ |\
\exists j$ tel que $ u_{i,j}+u_i\in O^*\}$ et lorsque $ I^*_{O,O^*}\not=\emptyset$ soit $u_{i(O^*)}$ l'unique repr\ésentant tel que $u_{i(O^*)}+u_i\in W^*_k\cap O^*.$ 

\noi Alors:

\noi $\bullet$ lorsque $ I^*_{O,O^*}=\emptyset$ on a
$a^*_{O,O^*}(\pi)=0$ et sinon
 
$$\bullet \ a^*_{O,O^*}(\pi)= \sum_{ i\in I_{O,O^*}  
 }a^{*(x_2)}_{x_0,u_i }(\pi_{k+1},...,\pi_p)\gamma_k^{-1}
(x_2,u_i^{-1})a^{*(u_i)}_{x_2,u_{i(O^*)} }(
\pi_1,...,\pi_{k-1},\pi_{k}^{\frac{1}{m_k}}\prod_{k+1\≤j\≤p}\pi_j|\ |^{\frac{md'_k}{m_k}}) $$
\end{rema}

\bigskip
 \subsection{Le cas complexe}
\bigskip
Dans les notations du paragraphe pr\éc\édent,  lorsque $\F=\C$ les \PVs ont tous une seule orbite non singuli\ère et $\gamma_k\equiv 1$ donc la proposition 5.1.1 donne:
\bigskip
\begin{cor} $\F=\C$
$$a(\pi)=a^{\goth U}(\pi_{p-k+1},...,\pi_p)a^{\goth U'}(\pi_1,...,\pi_{p-k-1},\pi_{p-k}^{\frac{1}{m_k}} \prod_{p-k+1\≤j\≤p}\pi_j |\ |^{\frac{r_k}{m_k}})\ ,$$
avec $\goth U:=\goth U(\goth s),$ $\goth s$ \étant l'alg\èbre engendr\é par $z_0$ et $z_0^{-1},$  $\goth U':=\goth U(\goth s'),$ $\goth s'$ \étant l'alg\èbre engendr\é par $z_{-2}$ et $z_{-2}^{-1}$ et $a^{\goth U},$ $a^{\goth U'}$ \étant les coefficients correspondant pour les fonctions Z\étas associ\ées.
\end{cor}
\bigskip
\noi Par cons\équent on relie classiquement le polynome de Bernstein au coefficient de l'\équation fonctionnelle, r\ésultat d\éj\à \établi pour quelques cas (cf. \cite{boppruben}, \cite{farautkoranyi}, \cite{clerc},...).
 
\bigskip
\begin{theo} 
\begin{enumerate}
\item On suppose que $(\∆,\lambda_0)$ est  une $\F$-forme
de $(E_7,\alpha_6)$  ou bien que $(\∆,\lambda_0)$ est de type $(C_n,\alpha_k)$ avec  $\goth t=\goth t'_0,$  et soit 
$$b_{\goth g,P_{\goth t}}(s_1, s_2)=C \prod_{j=1}^{\frac{d_1}{2}} \biggl(\ (s_2+\lambda_{2,j}) (2s_1+ s_2+
\lambda_{1,j})\biggr)$$
le polynome de Bernstein du \PV $(P_{\goth t},\goth g_1),$ alors il existe une constante $D$ telle que:
$$a(\omega_{q_1}.|\ |^{s_1}, \omega_{q_2}.|\ |^{s_2})=D\prod_{ j=1}^{\frac{d_1}{2}} \biggl( \rho'(\omega_{ q_2} ;s_2+ 
\lambda_{2,j}+1)\rho'(\omega_{ q_1}^2\omega_{ q_2} ;2s_1+s_2+ 
\lambda_{1,j}+1)\biggr).$$
\item  Dans tous les autres cas, soit 
$$b_{\goth g,P_{\goth t}}(s_1,...,s_p)=C \prod_{\ell=1}^p\biggl(\prod_{j=1}^{d_{p-\ell+1}-d_{p-\ell}}(s_{\ell}+...+s_p+
\lambda_{\ell,j})\biggr)$$
le polynome de Bernstein du \PV $(P_{\goth t},\goth g_1),$ alors il existe une constante $D$ telle que:
$$a(\omega_{q_1}.|\ |^{s_1},...,\omega_{q_p}.|\ |^{s_p})=D\prod_{\ell=1}^p\biggl(\prod_{j=1}^{d_{p-\ell+1}-d_{p-\ell}}\rho'(\omega_{\ell}....\omega_p;s_{\ell}+...+s_p+
\lambda_{\ell,j}+1)\biggr).$$
\end{enumerate}
\end{theo}

\dem Par r\écurrence sur $p.$ 

\noi Pour $p=1$ cela d\écoule du th\éor\ème 2 de J.I.Igusa  \cite{igusa5} et pour $p\≥2$ on applique le corollaire 5.2.1 ainsi que la remarque 3.7.4  et la proposition 3.7.3.

\noi Pour \éviter d'utiliser le  th\éor\ème 2 de J.I.Igusa, on aurait pu commencer par montrer le  th\éor\ème pour les paraboliques  tr\ès sp\éciaux minimaux en utilisant les r\ésultats des $2$ exemples fondamentaux.\fdem\bigskip

\noi{\bf Remarque 5.2.3} A l'aide des r\ésultats des sections suivantes, on v\érifie qu'avec les normalisations choisies on a:
$$C=D=1\ (1).$$
c'est \à dire que dans la notation suivante:

\noi  A toute application polynomiale $b:\C^p\mapsto \C,$ $b(s_1,...,s_p)=\prod_{1\≤j\≤q}(a_{1,j}s_1+...+a_{p,j}s_p+a_j  ),$  les $
 a_{r,j},$ \étant entiers  lorsque $r=1,...,p$ et $j=1,...,q$   et les coefficients $a_j$ rationnels  pour $ j=1,...,q,$  
 
 \noi et tout caract\ère continu $\pi=(\pi_1,...,\pi_p)$ de $(\hat \F^*)^p,$ on associe la quantit\é:
 $$\rho'_b(\pi)=\prod_{1\≤j\≤q}\rho'(\pi_1^{a_{1,j}}...\pi_p^{a_{p,j}}|\ |^{a_j +1 }),$$
 $\rho'(\pi_1)=\pi_1(-1)\rho(\pi_1)$ \étant le coefficient de Tate (cf.3.6.1) alors on a:
 $$a(\pi)=\rho'_{b_{\g,P_{\goth t}}}(\pi).$$
 \bigskip

\dem 
 \noi En effet, \à l'exception de l'unique cas $(C_n,\alpha_k),k$ pair,
 avec le \sgp $P_{\goth t'_0},$ (1) sera v\érifi\é dans les sections qui suivent  pour tous les \PVs $(P_{\goth t},\g_{\pm 1})$ lorsque le \sgp standard tr\ès sp\écial est maximal parmi ceux-ci  donc (1) est v\érifi\é pour tous les \sgps standards tr\ès sp\éciaux (puisque cela correspond \à certaines valeurs nulles de $s_{i_1},...$) sauf dans le cas $(C_n,\alpha_k),k$ pair,
 avec le \sgp $P_{\goth t'_0}.$
 
 \noi Mais  dans ce dernier cas $2H_0=H_1+H_2,$ $H_1$ et $H_2$ \étant dans la m\ême orbite de $Aut(\g)$ (cf.3) $\beta$ de la \demo de la prop.1.2.4) donc $$\widetilde B(\frac{H_1}{2},\frac{H_1}{2})=\frac{1}{4}\widetilde B(H_1,H_1)=\frac{1}{2}\widetilde B(H_0,H_0)=-\frac{\text{degr\é}(F_2)}{4}.$$
Or les \PVs $(\goth U_0, \goth U_1)$ et $(\goth U'_0, \goth U'_1)$ sont des \PVs commutatifs de type $(D_k,\alpha_k)$ et les \irfs sont de degr\é $\ds\frac{\text{degr\é}(F_2)}{2}=\ds \frac{k}{2}$ (m\ême ref. qu'avant) ainsi les normalisations sont coh\érentes avec la descente et les coefficients de la prop.3.7.3 valent $1$ ($A_1=B_1=1$) donc $C=\pm 1$ et par la proposition 5.1.1 et le th.6.2.1 ii) on a $a(\pi_1,\pi_2)=
 \rho'_B(\pi_1,\pi_2)$ avec $b_{\g,P_{\goth t}}(s_1,s_2)=C.B(s_1,s_2).$\\
 
 \noi Il reste \à v\érifier que $C=1.$ Pour ceci on consid\ère la m\ême situation sur $\R$ c'est \à dire le \PV $((\g_\R)_0,(\g_\R)_{\pm 1}) $ pour lequel on a encore $a(\pi_1,\pi_2)=
 \rho'_B(\pi_1,\pi_2)$, le coefficient $\rho'$ \étant cette fois-ci le coefficient de Tate {\bf r\éel}, mais alors par le 2) du lemme 3.7.1 on a:
 $$b(s_1,s_2)=(-2\sqrt{-1}\pi)^k\frac{( a(|\ |^{s_1},|\ |^{s_2})}{( a(|\ |^{s_1},\tilde\omega_{-1}|\ |^{s_2-1})}=\frac{(-2\sqrt{-1}\pi)^k}{(2\sqrt{-1}\pi)^k}B(s_1,s_2)=B(s_1,s_2).$$
 \fdem
 
 \bigskip

 La d\étermination des orbites pouvant s'av\érer difficile notamment dans la proposition 5.1.1 (cf.par exemple $\S 8.2.3$ et $8.2.4$) et au vu des 
r\ésultats du corollaire 3.6.3 et du th\éor\ème 3.6.5, on termine ce paragraphe par une situation particuli\èrement simple. 
\bigskip

\bigskip
\subsection{Un cas particulier}
      \bigskip

\noindent  Dans ce paragraphe,   on suppose que $(\∆,\lambda_0)$ n'est pas une $\F$-forme
de
$(E_7,\alpha_6)$ et que $\goth t\not=\goth t'_0$  lorsque   $(\∆,\lambda_0)$ est de type $(C_n,\alpha_k)$  c'est \à dire qu'on suppose que  $m_k=1.$\\

\noi On note par  $\mathbb H$ un sous-groupe particulier de $\F^*$ contenant $\F^{*2}  :$\\

\quad $\bullet$ soit $ \mathbb H:= \F^{*2},$\\

\quad $\bullet$  soit  $\mathbb H$ contient chaque $\chi_i(G_{\goth t})$
pour $i=1,...,p,$ \\
 
 \noindent Pour $u=(u_1,...,u_p)\in \bigl(\F^*/\mathbb{H}\bigr)^p$ on d\'efinit les ouverts
(\'eventuellement vides):\\

$O_u= \{x\in \goth
g_1\ |\ F_1(x)\mathbb H =u_1\ ,\ F_2(x) \mathbb H =u_1u_2\
  ,...,F_p(x) \mathbb H =u_1...u_p\  \},$ 

\noi (not\é \également $O_{u_1,...,u_p}$)\\

$O^*_u= \{x\in
\goth g_{-1}\ |\ F^*_1(x) \mathbb H =u_p\  ,\ F^*_2(x)\mathbb H =
 u_{p-1}u_p\   ,...,F^*_p(x) \mathbb H = u_p...u_1\  \}$ 
 
 \noi (not\é \également $O^*_{u_1,...,u_p}$).\\

 \noi On a:
 $$x\in O_u\Leftrightarrow x^{-1}\in O^*_u.$$
Chaque $O_u$ (resp.$O^*_u$) est invariant par $(\cap_{1\≤i\≤p}Ker \chi_i)N_{\goth t}.$\\

\noindent  On d\'efinit \'egalement  les fonctions Z\'etas associ\'ees, pour $f\in S(\goth
g_1)$ (resp.$h\in S(\goth g_{-1})$) :
$$  
Z_u(f;\omega)=Z(f{\bf 1}_{O_u};\omega) \quad (\hbox{resp.}
 Z^*_u(h;\omega)=Z^*(h{\bf 1}_{O^*_u};\omega)\  )\ .  $$
 
 \noindent Soit $x\in  W_k$ et $x=x_2+x_0$ sa d\'ecomposition   relativement \`a $ad(h_k),$ on suppose que $\gamma_k (x_2,x_0)$ (qui ne d\'epend que de la $G_{h_k}-$orbite de $x$) est
  ind\'ependant du choix de $x\in O_{u_1,...,u_p}$ et on pose: $\tilde \gamma_k (u_1,u_2,...,u_p)=\tilde \gamma_k (u):=\gamma_k(x_2,x_0).$\\
  
  \begin{prop}  Soit $w\in \bigl(\F^*/\mathbb{H}\bigr)^p.$\\
  
  \noi Lorsque $W^*_k\cap O^*_w\not=\emptyset,$  soit $z=z_0+z_{-2}\in W^*_k\cap O^*_w,$ on note  $\goth s$    l'alg\`ebre
engendr\'ee par
$z_0$ et $z_0^{-1},$  $\goth s'$  celle 
 engendr\'ee par
$z_{-2}$ et $z_{-2}^{-1},$ $\goth U:=\goth U(\goth s)$ et $\goth U':=\goth U(\goth s').$\\

On
suppose que:
 
\begin{enumerate}
\item  Pour tout $\pi'\in \Omega ((\F^*)^k)$ et $\forall (u,v)\in \bigl(\F^*/\mathbb{H}\bigr)^k\times \bigl(\F^*/\mathbb{H}\bigr)^k$ il existe une constante
$b_{v,u}^ {(w_{k+1},...,w_p)}(
\pi')$  telle que dans le \PV faiblement sph\érique : $((G_{\goth t})_{z_0}N_{k,2},\goth U_{\pm 1})$  on a les \équations fonctionnelles suivantes pour tout $\pi'\in \Omega ((\F^*)^k)$ et $\forall v\in \bigl(\F^*/\mathbb{H}\bigr)^k:$
$$ Z^*_v(\four
(f);\pi')=\sum_{u\in \bigl(\F^*/\mathbb{H}\bigr)^k}b_{v,u}^{ (w_{k+1},...,w_p)}( \pi')Z_u ( 
f;\pi'^*|\ |^{-  \frac{p_{2,0}-p_2}{ d_k}{\bf 1}_k})\ ,\
\forall f\in \EuScript S(\goth U_1),$$
avec   $b_{v,u}^{(w_{k+1},...,w_p)}=0$  si
l'un des deux ouverts est vide.

 \noi ($u$ est associ\é aux valeurs des \irfs $F_1,...,F_k$)

\item Pour tout $\pi"\in \Omega((\F^*)^{p-k})$ et $\forall (u',v')\in \bigl(\F^*/\mathbb{H}\bigr)^{p-k}\times \bigl(\F^*/\mathbb{H}\bigr)^{p-k}$ il existe 
une constante $c_{v',u'}^{( w_1,...,w_k)}(
\pi")$ telle que 
dans le \PV faiblement sph\érique : $((G_{\goth t})_{z_{-2}}N_{k,0},\goth U'_{\pm 1})$  on a les \équations fonctionnelles suivantes pour tout $\pi"\in \Omega((\F^*)^{p-k})$ et  $\forall v'\in \bigl(\F^*/\mathbb{H}\bigr)^{p-k} $   : 
$$Z^*_{v'}(\four
(g);\pi")=\sum_{u'\in \bigl(\F^*/\mathbb{H}\bigr)^{p-k}}c_{v',u'}^{( w_1,...,w_k)}(\pi")Z_{u'} ( 
g;\pi"^* |\ |^{-   \frac{p_{0,2}-p_2}{ d'_k}{\bf 1}_{p-k}})\ ,\ \forall g\in
\EuScript S(\goth U'_1) \ ,
$$
avec $c_{v',u'}^{( w_1,...,w_k)}=0$ si
l'un des deux ouverts est vide.

 \noi ($v'$ est associ\é aux valeurs des \irfs $F^*_1,...,F^*_{p-k}$)
\end{enumerate}

\noi alors pour $\pi\in \Omega((\F^*)^p)$ et $\forall \bigl(\F^*/\mathbb{H}\bigr)\in D^p $ et $\forall
h\in \EuScript S(\goth
g_1)$ on a:
$$ Z^*_v(\four (
h);\pi)=\sum_{u\in \bigl(\F^*/\mathbb{H}\bigr)^p}d_{v,u}(
\pi)  Z_u(h;\pi^*| \ |^{-N{\bf 1}_p}) \quad \hbox{avec}
$$
  $$ d_{v,u}(
\pi)=\tilde \gamma_k (u_1,...,u_k,v_{k+1},...,v_p)b_{(v_1,...,v_k),(u_1,...,u_k) }^{(v_{k+1},...,v_p)}( \pi').  
 c_{(v_{k+1},...,v_p),(u_{k+1},...,u_p)}^{(u_1,...,u_k)}(\pi") $$ 
lorsque $u=(u_1,...,u_p)$ et $v=(v_1,...,v_p),$

\noi avec $\pi '=(\pi_{p-k+1},...,\pi_p)$ et $\pi"=(\pi_1,...,\pi_{p-k-1},\ \prod_{p-k\≤j\≤p}\pi_j\ |\
|^{ r_k }).$

\end{prop}
 \vspace{2mm}
 \dem
 
 1) On reprend  bri\`evement la  \demo de la proposition 5.1 mais relativement aux ouverts $O_u$ et $O^*_v.$\\  
   
 \noindent Prenons $u,v\in \bigl(\F^*/\mathbb{H}\bigr)^p$ tels que $O_u$ et $O^*_v$ soient non vide,  $f\in \EuScript C_C^{\infty}( O_u)$ et $Re(\pi)>0.$\\
 
 \noindent  Posons $u'=(u_{k+1},...,u_p),$ $u"=(u_1,...,u_k) $ et
$v'=(v_{k+1},...,v_p),$ $v"=(v_1,...,v_k)$, on a :\\
$$  Z^*_v(\four ( f);\pi)\ 
  =  \int_{ E_{0,-2}\cap O^*_{v'}} C(y')|F^*_{p-k}(y')|^{\frac{p_1}{2d'_k}} Z^*_{v"}( \four_{(E_{2,0})_{y'}}(\gamma_k(\ ,y'^{-1})h_f(\ ,y')\
;\pi')dy'\ ,$$
\noindent par  l'hypoth\èse faite sur $\gamma$  on a:
$$ \four_{(E_{2,0})_{y'}\ }(\ \gamma_k(\ ,y'^{-1})h_f(\ ,y') \ )=\sum_{w\in \bigl(\F^*/\mathbb{H}\bigr)^k}\tilde \gamma_k (w,v')\four_{(E_{2,0})_{y'}}(1_{O_w} h_f(\ ,y') \ )$$
 donc on a \'egalement:
 $$Z^*_{v"}( \four_{(E_{2,0})_{y'}}(\gamma_k(\ ,y'^{-1})h_f(\ ,y')\
;\pi')=\sum_{w\in \bigl(\F^*/\mathbb{H}\bigr)^k}\tilde \gamma_k (w,v')Z^*_{v"}( \four_{(E_{2,0})_{y'}}(1_{O_w}  h_f(\ ,y')\
;\pi')$$
et par l'hypoth\èse 1:
$$Z^*_{v"}( \four_{(E_{2,0})_{y'}}(\gamma_k(\ ,y'^{-1})h_f(\ ,y')\
;\pi')=\sum_{w\in \bigl(\F^*/\mathbb{H}\bigr)^k}\tilde \gamma_k (w,v')b_{v",w}^{(v')}(\pi')
Z_{w}(h_f(\ ,y');\pi'^*|\ |^{-\frac{p_{2,0}-p_2}{d_k}{\bf 1}_k}) $$
 donc:
\begin{eqnarray*}  Z^*_v({\four f};\pi)\ 
 & = &  \sum_{w\in \bigl(\F^*/\mathbb{H}\bigr)^k}  \tilde \gamma_k ( w,v')b_{v",w}^{(v')} (\pi')I_{w}\ avec\\
 I_w&=&\int_{
E_{0,-2}\cap O^*_{v'}}\ C(y')|F^*_{p-k}(y')|^{\frac{p_1}{2d'_k}}
  Z_{w}(h_f(\ ,y');\pi'^*|\ |^{-\frac{p_{2,0}-p_2}{d_k}{\bf 1}_k}) dy'\\
 & = &  \int_{
E_{2,0}\cap O_{w}}
 Z^*_{v'}  
 (\four_{(E_{0,2})_x} (S_f(x+\  )) ;\pi").|F_k(x)|^{-{p_{2,0}\over d_k}} 
 \pi'^*(G(x))dx\\
 & =&   \sum_ {w'\in \bigl(\F^*/\mathbb{H}\bigr)^{p-k} }
 c_{ v',w'}^{(w)}(\pi") \int_{
E_{2,0}\cap O_{w}}Z_{w'}(S_f(x+\
);\pi"^* |\ |^{-\frac{p_{0,2}-p_{2} }{ d'_k}{\bf  1}_{p-k}}) .    |F_k(x)|^{-{p_{2,0}\over d_k}}   \pi'^*(G(x))dx \end{eqnarray*}
en utilisant l'hypoth\èse 2. d'o\ù:
$$
   I_w=  \sum_ {w'\in \bigl(\F^*/\mathbb{H}\bigr)^{p-k} }
 c_{ v',w'}^{(w)}(\pi") Z_{(w,w')}(f;\pi|\ |^{-N{\bf 1}_p})=c_{ v',u'}^{(u")}(\pi") Z_ u(f;\pi|\ |^{-N{\bf 1}_p})$$
 par choix du support de $f,$ donc:
 $$Z^*_v({\four f};\pi)=
  \tilde \gamma_k ( u",v')  
 b_{v",u"}^{(v')}(\pi') c_{ v',u'}^{(u")}(\pi")Z_u(f;\pi|\ |^{-N{\bf 1}_p}).   $$

 \noindent Ainsi   $\forall f\in C_c^{\infty}( O_u)$ et
$Re(\pi)>0,$ on a:
$$ Z^*_v(\hat
f;\pi)= d_{v,u}(
\pi)  Z_u(f;\pi^* |\ |^{-N{\bf  1}_p})\ , 
$$ 
et $d_{v,u}=0$ lorsque $O_u$ ou $O^*_v$ sont vides (puisque  $ O_u= \emptyset\ \Leftrightarrow\ O_{u"}= \emptyset$ ou $O_{u'}=\emptyset$ et idem pour $O^*_v$).\\

 \noindent  2) Lorsque $\mathbb H $ contient chaque $\chi_i(P_{\goth t})$ pour $i=1,...,p,$ les ouverts non vides $O_u$ (resp. $O^*_u$) sont r\'eunion de $P_{\goth t}-$orbites dans $\goth g''_1$ (resp. $\goth g''_{-1}$) donc 
 pour $Re(\pi)>0$ et $O\subset O_u$ on a : $d_{v,u}
( \pi )=\sum_{\{O^* \ |\ O^*\subset O^*_v\}}a_{O^*,O}( \pi) ,$  par prolongement m\'eromorphe, cette \'egalit\'e est vraie pour
tout caract\`ere $ \pi$
 donc $\sum_{\{O^* \ |\ O^*\subset O^*_v\}}a_{O^*,O}( \pi)$ est
ind\'ependante de l'orbite 
$O\subset O_u$ d'o\`u le r\'esultat. \\

\noindent  3)  Lorsque $\mathbb H=  \F^{*2},$  
pour  $v\in \bigl(\FF\bigr)^p$,  
$h\in \EuScript S(\goth g_{-1})$ et
$Re( \pi)>0$ on a:
$$Z^*_v(h;\ \pi)=\frac{1 }{ |\FF|^p}\sum_{(a_1,...,a_p)\in
\bigl(\FF\bigr)^p}\biggl(\prod_{1\≤i\≤p}(a_i,v_p...v_{p-i+1})\biggr)\
\ Z^*(h; \pi. (\tilde\omega_{a_1},..., \tilde\omega_{a_p}))\ .$$
 
\noi Par prolongement m\'eromorphe, cette \'egalit\'e est vraie pour
tout caract\`ere $ \pi$ donc, en appliquant l'\'equation
fonctionnelle abstraite   \`a $Z^*(\four( f); \pi. (\tilde\omega_{a_1},...,\tilde\omega_{a_p})),$ on
montre qu'il existe une  fonction m\'eromorphe en $ \pi$
d\'ependant de l'orbite $O,$ de $u $ et de $v$, not\'ee
$
\alpha_{O, \pi}(u,v),$  telle que 
$\forall f\in \EuScript S(\goth g_{1})$:
$$Z^*_v(\four( f); \pi))=\sum_{\{u\in
(\FF)^p\ ,\ \hbox{orbites}\ O \  |\
O\cap O_u\not=\emptyset\}}\alpha_{O, \pi}(u,v)\ Z(f{\bf
1}_{O\cap O_u}; \pi^*.|\ |^{-N{\bf  1}_p})  .$$ Or, pour $f\in
\EuScript C_C^{\infty}(O\cap O_u)$ et $Re( \pi) >0,$ on a:
$$\alpha_{O, \pi}(u,v)=d_{v,u}( \pi)$$
donc par prolongement m\'eromorphe, cette \'egalit\'e est vraie pour
tout caract\`ere $ \pi$ donc $\alpha_{O, \pi}(u,v)$ est ind\'ependant de l'orbite $O$ rencontrant $O_u$ d'o\`u le 
r\'esultat.   \fdem\\
 
 \begin{rema}

\begin{enumerate}
\item Lorsque $O_u$ est inclus dans une seule orbite de $P_{\goth t}$ dans $\goth g"_1,$ pour $u$ et $v$ dans $D^p,$ on dit que $u\sim  v\ \Leftrightarrow$ $O_u$ et $O_v$ sont dans la m\^eme orbite de $P_{\goth t}$  (ce qui est \'equivalent \`a $O^*_u$ et $O^*_v$ sont dans la m\^eme orbite de $P_{\goth t}$) alors   $\forall
f\in \EuScript S(\goth
g_1),$   $ \pi\in \Omega({(\F^*)}^p)$ et  pour toute orbite $O^*$de $P_{\goth t}$ dans $\goth g"_{-1},$ on a:
$$ Z^*_{O^*}(\four (
f);\pi)=\sum_{\{O|O\ \hbox{orbites de } P_{\goth t}\ \hbox{dans}\ \goth g"_1\}}A_{O^*,O}( \pi) Z_O(f; \pi^*|\ |^{-N{\bf 1}_p}) \quad \hbox{avec}
$$
  $$A_{O^*,O}(\pi)= 
 \sum_{\{w\in \bigl(\F^*/\mathbb{H}\bigr)^p,w\sim v\}}$$
 $$\tilde \gamma_k (u_1,...,u_k,w_{k+1},...,w_p)b_{ (w_1,...,w_k),(u_1,...,u_k)}^{(w_{k+1},...,w_p)}( \pi')
c_{ (w_{k+1},...,w_p),(u_{k+1},...,u_p)}^{(u_1,...,u_k)}( \pi")$$  $u=(u_1,...,u_p),$  $v=(v_1,...,v_p)$ choisis tels que $O_u\subset O$ et  $O^*_v\subset O^*$ (mais quelconques).\\

\item Lorsque $ \mathbb H=
\F^*,$ on normalise les invariants $F_1,...,F_p$ de telle mani\`ere qu'ils repr\'esentent tous $1$ alors
  $\forall
f\in \EuScript S(\goth
g_1)$ on a:
$$ Z^*(\four(
f); \pi)=A( \pi)  Z(f;\pi^*|\ |^{-N1_p}) \quad \hbox{avec}
$$
  $$A(\pi)= \tilde \gamma_k (1,...,1)b_{(1,...,1),(1,...,1)}^{(1,...,1)}( \pi')
c_{(1,...,1),(1,...,1)}^{(1,...,1)}( \pi")\ .$$ 
 \end{enumerate} 
    \end{rema}
 \bigskip
        \noi  Terminons par les 2 exemples r\éels suivants pour lesquels  l'\irf $F$ est de degr\é $4$  donc $\widetilde B(H_0,H_0)=-2$:\\

$\bullet$  $(E_6,\alpha_2)$ de type III r\éel c'est \à dire que $\g$ est une alg\èbre de Lie de rang $2$ et de diagramme de Satake de type EIII,  $(\overline{\g_0},\overline{\g_1})$ est de type $(E_6,\alpha_2),$  donc $(\g_0,\g_1)$ est de type $(BC_2,\lambda_1):$ 
  
 
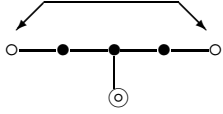
\begin{picture}(100,50)(-50,0)
\put(-25,5){\vector (-1,-1){10.5}}
\put(-25,5){\line (1,0){50}}
\put(25,5){\vector (1,-1){10.5}}
\end{picture}

  \hskip 10pt \hbox to 3 cm {\lower 2pt\hbox{$\circ$}\hrulefill\lower 2pt
\hbox{$\bullet$}\hrulefill \lower 2pt \vtop {\offinterlineskip \hbox{$\bullet$} \hbox to 5pt{\hfill \vrule height 12pt width 0,3pt\hfill} \hbox{$\circledcirc$}}\hglue -3,5pt\hrulefill\lower 2pt\hbox{$\bullet$}\hrulefill\lower 2pt\hbox{$\circ$} } 
 \\ \\

 ($\lambda_1$ est la restriction de $\alpha_2$ et $\lambda_2$ de $\alpha_1$ et \également de $\alpha_6$)\\

$\bullet$  $(E_7,\alpha_1)$ de type EVII r\éel c'est \à dire que $\g$ est une alg\èbre de Lie de rang $3$ et de diagramme de Satake de type EVII,  $(\overline{\g_0},\overline{\g_1})$ est de type $(E_7,\alpha_1),$  donc $(\g_0,\g_1)$ est de type $(C_3,\lambda_1):$   \\ \\
 
 \hskip 10pt \hbox to 3cm {\lower 2pt\hbox{$\circledcirc$}\hrulefill\lower 2pt
\hbox{$\bullet$}\hrulefill \lower 2pt \vtop {\offinterlineskip \hbox{$\bullet$} \hbox to 5pt{\hfill \vrule height 12pt width 0,3pt\hfill} \hbox{$\bullet$}}\hrulefill\lower 2pt\hbox{$\bullet$}\hrulefill\lower 2pt\hbox{$\circ$}\hrulefill\lower 2pt\hbox{$\circ$} }
\\  \\

($\lambda_1$ est la restriction de $\alpha_1,$   $\lambda_2$ de $\alpha_5$ et  $\lambda_3$ de $\alpha_6$)\\

 \noi   Dans les $2$ cas le  \sgp standard tr\ès sp\écial est donn\é par $P=P(H_1,H_2)$ avec $H_1=h_{\lambda_1}$ et $H_2=2H_0-H_1=h_\mu;$ $H_1$ et $H_2$ sont dans la m\ême orbite de $G.$\\
 
 \noi Indiquons \également les diff\érents sous-espaces qui interviennent ainsi que leurs  dimensions: $\g_2=E_{2,2}=\g^{\widetilde\lambda}$ est de dimension $1,$   $\widetilde\lambda$ \étant la plus grande racine; $E_{1,1}=\g^{\lambda_1+\lambda_2}$ est de dimension $4d,$  $\g^{\lambda_1}$ de dimension $d+2$ d'o\ù  $\g_1$ est de dimension   $8+6d.$\\\\

  \begin{tabular}{ |c||c|c|}

\hline

 Type r\éel
  
  & EIII&  EVII\\\hline

 $\mu$& $\lambda_1+2\lambda_2 $  &$ \lambda_1+2\lambda_2+\lambda_3$ 
 
 \\\hline
 &&\\
 
  $\widetilde\lambda$
& $ 2(\lambda_1+\lambda_2)$ &$2 (\lambda_1+\lambda_2)+\lambda_3$  
 
 \\\hline
  
 $d$
& $ 2$ & $ 4$\\\hline
 \end{tabular}\\ \\

La restriction de l'\irf  $F_1$  (resp. $P_1$) \à $\goth U(\R H_2)_1=\g^{\lambda_1}$ (resp. \à $\goth U(\R H_1)_1=\g^{\mu}$  est donc une forme quadratique anisotrope  (chaque $\goth U(\R H_i),$ $i=1,2,$ est de rang $1$).   \\

\begin{prop} \begin{enumerate} 
\item $b_1(s_1,s_2)=s_2(s_2+\ds\frac{1+d}{2})\ , \ b_2(s_1,s_2)=s_2(s_2+\ds\frac{1+d}{2})(s_1+s_2+\ds\frac{2d+1}{2})(s_1+s_2+\ds\frac{3d+2}{2}).$

\item  $P$ a une seule orbite dans $\g"_1.$

 \item Soient $(s_1,s_2)\in \widehat{(\C)^2},$ 
 $f\in \EuScript S(\g_1)$ alors:
 $$ Z^*(\four f;(s_1,s_2))=C_d(s_1,s_2)\sin(\pi s_2)\cos (\pi(s_1+s_2)) Z(f;(s_1,-s_1 -s_2 - \ds\frac{3}{2}d-2)\quad \text{avec} $$
 $$C_d(s_1,s_2) =4(2\pi)^{-2s_1-4s_2-3d-6}\Gamma(s_2+1) \Gamma(s_2+\ds\frac{d+3}{2}).$$
 $$\Gamma(s_1+s_2+\ds\frac{2d+3}{2})\Gamma(s_1+s_2+ \ds\frac{3}{2}d+2) . $$
 \end{enumerate}

\end{prop}

\dem 1) D\ûe \à la normalisation de $\widetilde B.$\\

2) Soit $x$ et $x'$  $2$ \elts de $W_{\goth t},$ d\écomposons $x=x_2+x_0$ et $x'=
x'_2+x'_0$ avec $x_i$ et $x'_i\in E_i(  H_1 )\cap \goth g_1,$ alors $x_2$ et $x'_2$ sont dans la  m\ême orbite de Int$(\goth U(\R H_1)_0)$ car l'alg\èbre $\goth U(\R H_1)$ est de rang $1$  et le \PV $(\goth U(\R H_1)_0,\goth U(\R H_1)_1)$ est commutatif (cf.par exemple d\émonstration de 6.1.7) donc on peut supposer que $x_2=x'_2$ ensuite on proc\ède de m\ême avec $x_0$ et $x'_0$ dans le \PV commutatif  $((\goth U_{x_2})_0,(\goth U_{x_2})_1)$ avec $\goth U_{x_2}:=\goth U(\R x_2\oplus \R H_2\oplus \R x_2^{-1}).$ \\

3) On applique le 4) du th\éor\ème 3.6.5 ainsi que la proposition 5.3.1 avec $\mathbb{H}=\R^*.$

\noi Soit $x\in W_{\goth t},$ $x=x_2+x_0$ avec $x_i\in E_i(  H_1 )\cap \goth g_1$ alors on a vu que   $-x_2+x_0$ et $x_2+x_0$ sont dans la m\ême orbite de $G_{H_1}, $ donc la forme quadratique d\éfinie sur $E_{-1,-1}$ par $\widetilde B(ad(\ )^2(x_2 ),x_0)$ est de type $(2d,2d)$ d'o\ù  $ \gamma_1(x_2,x_0)=  1.$\fdem\\

\noi Remarques: 1) Cette proposition est bien connue lorsque $s_1=0$ (\cite{muro1}).\\

2) La situation n'est plus aussi simple pour les autres formes qui seront trait\ées dans la section 8.2 (prop.8.2.6, cf.\également prop. 8.2.11).

 \newpage
\section{\bf Le cas  commutatif }

 Dans cette section, on traite les\eqs des \PVs de type commutatifs, bien que ces r\ésultats soient connus en grande partie, comme cela a \été rappel\é dans l'introduction, afin de les exprimer dans nos notations et normalisations pour une utilisation ult\érieure (cf. $\S 8.3$ qui traite du cas exceptionnel $(E_7,\alpha_2)$ et remarque 5.2.3).\\
 
 \noi {\bf Notation}: $\phi$ \étant une forme quadratique non d\ég\én\ér\ée, on note simplement $\gamma (\phi)$ la constante $\gamma(\tau\circ \phi)$ (cf.$\S 3.6.2$).\\

\subsection{ \bf Structure}

\smallskip

\subsubsection{ \bf Rappels}

\smallskip

\noi On consid\ère le\sg parabolique tr\ès sp\écial standard d\éfini dans le lemme 2.3.2, $P_0=P(h_{\lambda_1},...,h_{\lambda_n}),\{\lambda_n,...,\lambda_1\}$ \étant l'ensemble maximal canonique de racines orthogonales de $\∆_1$ (l'ordre est invers\é).\\

\noi On rappelle que le syst\ème de racines obtenu \à partir des restrictions non nulles de $\∆$ \à $\goth t_0=\oplus_{1\≤i\≤n}\F h_{\lambda_i}$ est de type $C_n,$
 qu'il existe une  alg\`ebre d\'eploy\'ee, not\'ee $\tilde {\goth g},$ 
 admettant $\goth t_0 $ 
  comme sous-alg\`ebre de Cartan (prop.2.2.1 \cite{mullerNAG}) et que: 
  $$\tilde B=\displaystyle{-\frac{d_nB}{2B(H_0,H_0)}}(=\displaystyle{-\frac{2d_1B}{B(h_{\lambda_i},h_{\lambda_i})}}),$$
($d_n=nd_1$ lemme 6.1.4) ainsi  pour toute racine longue $\alpha$ de $\∆$ on a $\tilde B(X_{\alpha}, X_{-\alpha})=d_1.$  \\

 \noindent Pour $i\not= j,$ soit $E^{i,j}_{u,v}=\{x\in \goth g\ |\ [h_{\lambda_i},x]=ux\ ,\  [h_{\lambda_j},x]=vx\  ,\ [h_{\lambda_k},x]=0 \ $ pour $ 1\leq k\not=i,j\leq n\}$ le sous espace associ\é \à la racine $\ds\frac{u\lambda_i+v\lambda_j}{2}$ de $R$ ($u,v=\pm 1$). \\     

\begin{prop} 
\begin{enumerate}
 
\item Pour $1\leq k\not= r\leq n,$ $E^{k,r}_{1,1}\not=\{0\}$ et
il existe un \elt  de $Aut_e(\g_0),$ $g_{k,r}$ dont la restriction \à 
$\oplus_{1\≤i\≤n}\g^{\lambda_i}$ soit une involution qui se r\éduise \à l'identit\é sur $\oplus_{1\≤i\not=k,r\≤n}\g^{\lambda_i}$ et $g_{k,r}(\g^{\lambda_k})=\g^{\lambda_r}.$

\item Il existe un syst\`eme de Chevalley , $(X_{\mu},h_{\mu},X_{-\mu})_{\mu\in R}$, de $(\tilde {\goth g},\goth t_0)$ tel que
toutes les formes quadratiques $f_{\frac{ \lambda_i- \lambda_j}{2}}$ d\'efinies sur $E^{i,j}_{-1,1}$ par 
$f_{\frac{ \lambda_i- \lambda_j}{2}}(A)=\frac{1}{2} \tilde B(ad(A)^2(X_{\lambda_i}),X_{-\lambda_j}) $
sont \'equivalentes et repr\ésentent $d_1.$

\item Dans le cas r\éel, on peut supposer de plus que $\Theta (X_{\lambda_i})=X_{-\lambda_i},i=1,...,n,$ $\Theta$ \étant une involution de Cartan telle que $\Theta/\goth a=-id.$
 
 \item Pour toute orbite non r\éduite \à $\{0\}$ de $G$ dans $\goth g_1,$ il existe $j\in \{1,...,n\}$ telle que cette orbite  rencontre $\oplus_{1\≤i\≤j}(\g^{\lambda_i}-\{0\}).$
\end{enumerate}
\end{prop}

\dem 1)-2) C'est la \demo  de la prop.4.1.1 \cite{mullerNAG})  qui convient \également lorsque dim$(\g^{\lambda_i})>1$ pour $i=1,...,n$ (cf.\également la \demo du 1) du lemme  6.1.7).

\noi 3) Lorsque $\F=\R$ on a $[\Theta(X_{\lambda_i}),X_{\lambda_i}]=c_ih_{\lambda_i}$ avec $c_i=\tilde B(\Theta(X_{\lambda_i}), X_{\lambda_i})>0 $ donc 
$$(\frac{1}{\sqrt c_i}X_{\lambda_i},h_{\lambda_i},\Theta (\frac{1}{\sqrt c_i}X_{\lambda_i}))$$ est encore un \Sl  d'o\ù le r\ésultat  (cf.lemme 1.1.7 de  \cite{boppruben}). 

\noi  4) est une cons\équence de la prop.5.2.2 de \cite{mullerJA1} et du 1) de cette proposition.\fdem

\bigskip

\noi  On normalise les invariants relatifs fondamentaux par:
$$i=1,...,n\ : \ F_i(\sum_{1\≤j\≤n}X_{\lambda_j})=1\quad \hbox{donc}\quad F^*_i(\sum_{1\≤j\≤n}X_{-\lambda_j})=1\ \text{pour}\ i=1,...,n.$$
On note $E:=E^{1,2}_{ - 1, 1},f:=f_{\frac{ \lambda_1- \lambda_2}{2}},$ $d$ est la dimension commune des sous-espaces $E^{i,j}_{\pm 1,\pm 1},$  et $\delta:=(-1)^{[\frac{d}{2}]}.$discriminant de $f.$  \\

\noi On rappelle que:
$$\begin{array}{lll}\g&=&\g_{-1}\oplus \g_0\oplus \g_1\\
\\
 \g_{-1}&=&\oplus_{1\≤i\≤n}\g^{-\lambda_i}\oplus_{1\≤i<j\≤n}E^{i,j}_{-1,-1}\\
 \\
 \g_{1}&=&\oplus_{1\≤i\≤n}\g^{\lambda_i}\oplus_{1\≤i<j\≤n}E^{i,j}_{1,1}\\
 \\
 \g_{0}&=&E(0) \oplus_{1\≤i<j\≤n}\bigl(E^{i,j}_{-1,1}\oplus E^{i,j}_{1,-1}\bigr)\\
 \\
 E(0)&=&\cap_{1\≤i\≤n}E_0(h_{\lambda_i})\ .\end{array}$$
 
 \noi En raison des relations de commutation:
 $$W_{\goth t_0}=\oplus_{1\≤i\≤n}\ (\goth g^{\lambda_i}-\{0\})\quad (\text{resp.}
\ W^*_{\goth t_0}=\oplus_{1\≤i\≤n}\ (\goth g^{-\lambda_i}-\{0\})\ )$$
et $F/W_{\goth t_0}=\prod_{1\≤i\≤n}G_i$ (resp.$F^*/W^*_{\goth t_0}=\prod_{1\≤i\≤n}G^*_i$), $G_i$ \étant l'\irf du \PV $(E_0(h_{\lambda_i})\cap \goth g_0,\goth g^{\lambda_i})$ normalis\é par $G_i(X_{\lambda_i})=1$ donc $G^*_i,$ \irf du \PV dual, v\érifie $G^*_i(X_{-\lambda_i})=1$ ($G_1=F_1$ et $F^*_1=G^*_n$).\\

\begin{defi} Le \PV \gog est dit presque d\éploy\é si les sous-espaces radiciels $\g^{\lambda_r},r=1,...,n$ sont de dimension $1.$

\end{defi}

\noi  Lorsque $\frac{\lambda_2-\lambda_1}{2}\in \∆,$ $f_0$ est la restriction de $f$ \à $\goth g^{\frac{\lambda_2-\lambda_1}{2}}$ et $e$ est la dimension de $\goth g^{\frac{\lambda_2-\lambda_1}{2}}$, on rappelle que:

 \begin{lem}  \begin{enumerate}
\item $f$ repr\'esente $d_1.$  
\item $\goth g$ est de rang $n$ $\Leftrightarrow$ $f$ est anisotrope.
 
\item Lorsque \gog est presque d\éploy\é:

i) $f \sim f_0\oplus f_1,$ $f_0$ est anisotrope et   $f_1 $ est hyperbolique.

ii) Lorsque $n\≥2,F_2\sim-f\oplus$ la forme hyperbolique \à deux variables.

 iii)  Lorsque $e=0$ ou bien $e=d,$
   $f$ et $af,$ $a\in \F^*,$ sont \'equivalentes $\Leftrightarrow$ $a$ est un \'el\'ement de $f(E)^*$ \à {\bf l'exception} du cas $d=3$ lorsque $\F$ est un corps $\goth p-$adique.
   
 iv)  $f(E )^*\subset \cap_{1\≤i\≤n-1} \chi_i(P_0).$
   
\end{enumerate}
\end{lem}

\dem 1. r\ésulte du 2. de la prop.6.1.1.

\noi 2. r\ésulte de la \demo du 2) du lemme 4.1.3 de \cite{mullerNAG}.

 \noi 3.i) Ce r\ésultat,  bien connu dans le cas r\éel (lemme 2-22 p.49, \cite{boppruben}), r\ésulte d'un simple calcul puisque pour $\mu\in \∆_{-1,1}=\{\mu\in \∆\ |n(\mu,\lambda_2 )=-n(\mu,\lambda_1 )=1\}$ on a $\mu'=-(\mu+\lambda_1-\lambda_2)\in \∆_{-1,1}$ or:
$$s_\mu(\lambda_1)(h_{\lambda_2})=(\lambda_1-n(\lambda_1,\mu)\mu)(h_{\lambda_2})=-n(\lambda_1,\mu)\≤2,$$
et l'\égalit\é : $n(\lambda_1,\mu)=-2\Rightarrow s_\mu(\lambda_1)=\lambda_2$ donc $\mu=\frac{\lambda_2-\lambda_1}{2}$ est une racine courte, ainsi
:
$$\∆_{-1,1}=\begin{cases} \{\mu_i\ , \mu_i':=-(\mu_i+\lambda_1-\lambda_2),i=1,...,l_1\}\ \text{si }e=0\\
 \{\frac{\lambda_2-\lambda_1}{2},\mu_i\ , \mu_i':=-(\mu_i+\lambda_1-\lambda_2),i=1,...,l_1\}\ \text{si }e>0,\end{cases}$$
 les racines $ \{\mu_i\ , \mu_i',i=1,...,l_1\}$ ayant m\ême longueur que $\lambda_1.$

\noi Par dualit\é de la forme de Killing, on a pour $A_0\in \goth g^{\frac{\lambda_2-\lambda_1}{2}}$ et $ x_i,y_i\in\F,i=1,...,l_1:$
 $$f(\sum_{i=1}^{l_1}(x_iX_{\mu_i}+y_iX_{\mu'_i})+A_0)=f(A_0)+\sum_{i=1}^{l_1}x_iy_i\tilde B([X_{\mu_i},X_{\lambda_1}],([X_{\mu'_i},X_{-\lambda_2}]) , $$
Tout $A\in \goth g^{\frac{\lambda_2-\lambda_1}{2}}-\{0\}$ se compl\ète en $1$ \Sl $(A,h_{\lambda_2}-h_{\lambda_1},B),B\in \goth g^{\frac{\lambda_1-\lambda_2}{2}}-\{0\},$ donc $ad(A)^2$ est une bijection de $ \goth g^{ \lambda_1}$ sur $ \goth g^{ \lambda_2  }$ qui sont de dimension $1$ d'o\ù $\exists x\in \F^*$ tel que $ad(A)^2(X_{\lambda_1})=xX_{\lambda_2}$ et $f(A)=\frac{x}{2}\not=0.$

\noi 3.ii) Soit $x=a_1X_{\lambda_1}+a_2X_{\lambda_2}+y,y\in E^{1,1}_{1,1},a_2\not=0,$ alors par le calcul habituel (cf.\demo du lemme 1.1.1):
$$x=exp(ad(A))(a_1X_{\lambda_1}+bX_{\lambda_2})\ \text{avec }b=a_2-\frac{B(ad(y)^2(X_{-\lambda_1},X_{-\lambda_2})}{2a_1B(X_{\lambda_2},X_{-\lambda_2})},$$
donc avec les normalisations choisies:
$$F_2(x)=a_1b=a_1a_2-f(\theta_{X_{\lambda_1}}(1)(y)) .$$

\noi 3)iii) Comme $f\sim f_1$ ou bien $f\sim f_0$,  3.iii) est  \évident dans le cas r\éel puisque $f_0$ est d\éfinie positive et dans le cas $\goth p-$adique  car alors $d\≤4$ et $f$ repr\ésente $1.$

\noi 3. iv)  r\ésulte de la \demo du 1) du lemme 4.1.3 de  \cite{mullerNAG}.\fdem

 \bigskip
 
 \begin{lem} Soit $\overline d=2\ds\frac{\text{dim}(\g_1)-d_n}{(d_n-1)d_n}$  alors:
 $$ d=(d_1)^2\overline d \  ,\ N=\frac{\text{dim}(\g_1)}{d_n}=\frac{\overline d }{2}(nd_1-1)+1 \ \text{et }\ d_k=k  d_1\ \text{pour }\ k=1,...,n.$$

 \end{lem}
 
 \dem 1) Comme tous les sous-espaces $\goth g^{\lambda_j},$ $j=1,...,n,$ sont  conjugu\és et qu'ils commutent,  la relation sur les degr\és des invariants s'\établit par r\écurrence sur $n$ (cf.lemme 1.4.7).\\

2) Lorsque \gog est quasi-d\éploy\é, $F_k$ est de degr\é $k$ pour $k=1,...,n$ (prop.3.2 de \cite{mullerNAG}) donc $d_1=1$ et $\overline d=d.$\\
 
 3) Lorsque dim$(\goth g^{\lambda_1})>1,$ sur une extension alg\ébrique convenable de dimension finie de $\F,$ not\ée $\mathbb E,$ 
$\overline {\goth g}:=\g\otimes_\F\mathbb E$ est d\éploy\ée, on a alors les \égalit\és suivantes:
$$\begin{array}{rll} d_n&=&nd_1\\
\text{dim}(\overline {\goth g^{\lambda_1}})&=&d_1+\ds\frac{(d_1-1)d_1}{2}\overline d\\
\text{dim}(\overline {\goth g_1})&=&d_n+\ds\frac{(d_n-1)d_n}{2}\overline d\\
\text{dim}( \goth g_1)&=& n\text{dim}( \goth g^{\lambda_1})+\ds\frac{(n-1)n}{2}d\quad \text{   d'o\ù le r\ésultat.}\quad \Box \end{array} $$
 
 \smallskip

\noi Pour des raisons de normalisations, on redonne le r\ésultat suivant d\éj\à connu (\cite{farautkoranyi},\cite{boppruben}):\\
 
\begin{lem} $\F=\R:$ le polyn\ôme de Bernstein associ\é \à $(P_0,\goth g_1)$ est donn\é par:
$$(s_1,...s_n)\in \C^n\ :\ b(s_1,...s_n)=\prod_{k=0}^{n-1}\ \ \prod_{j=0}^{d_1-1}\ \bigl(\sum_{n-k\≤i\≤n}s_i+\frac{\overline d}{2}(j+d_1k)\ \ \bigr).$$

\end{lem}

\dem On applique la proposition 3.7.3 en notant que par le choix de la normalisation de $B$ on a:
$$k=1,...,n\ \:\ \  \widetilde B(\sum_{1\≤i\≤k}h_{\lambda_i},\sum_{1\≤i\≤k}h_{\lambda_i})=-2d_k,$$
et on a $d_1=1$ ou $2$
donc pour $(s_1,...s_n)\in \C^n$ on a:
$$ b(s_1,...s_n)=(C^{-d_1})^n\prod_{k=0}^{n-1}\ \ \prod_{j=0}^{d_1-1}\ \bigl(\sum_{n-k\≤i\≤n}s_i+\frac{\overline d}{2}(j+d_1k)\ \ \bigr)$$ avec:\\

 $\bullet$ lorsque $d_1=1$, $C=\widetilde B( X_{\lambda_1},  X_{-\lambda_i})=1$\\
 
$\bullet$ lorsque $d_1= 2$ , $C=\widetilde B( \frac{h_{\lambda_1}}{2},    \frac{h_{\lambda_1}}{2})=-1$ (remarque 3.6.6).\fdem
 
\smallskip

\noi \begin{rema}
 \noi Rappelons la terminologie utilis\ée dans \cite{boppruben}), on dit que:\\

\noi $\bullet$ \gog est de type I lorsqu'il est presque d\éploy\é et $f$ est anisotrope \ssi rang($\g)= d_n.$

On a donc $e=d$ et $d\≤4$ dans le cas $\goth p-$adique.

\noi $\bullet$ \gog est de type II lorsqu'il est presque d\éploy\é et $f$ est  isotrope \ssi rang($\g)> d_n$.

\noi On a donc $\gamma (f)= \gamma (f_0),$  $e<d$  et $\goth g$ est d\éploy\ée \ssi $e\≤1$ (\demo du lemme 6.1.3).

\noi $e=0$ lorsque $n\≥3$ (cf. tableau 1).

\noi $\bullet$ \gog est de type III lorsqu'il n'est pas presque d\éploy\é  \ssi rang($\g)<d_n$.

 \end{rema}
\bigskip
\noi Pour $u\in \bigl( \F^*/\F^{*2}\bigr)^n,$ $O_u$ et $O^*_u$ ont \ét\é d\éfinis dans le $\S 3.5.$
\bigskip

\subsubsection{Le type III}

\bigskip

Par les tables de \cite{warner} et \cite{veisfeiler} (cf. \également le  tableau 1 et les tables de \cite{boppruben}p.222-224), les cas possibles, en dehors du cas o\ù $\goth g$ est
de rang $1,$ sont :
\vskip 3mm
\quad i) $(\overline \∆,\alpha_0)=(A_{2nm-1}, \alpha_{nm})$ et 
$(\∆,\lambda_0)=(A_{2n-1},\alpha_n)$ avec $m\≥2$

\quad ii) $(\overline \∆,\alpha_0)=(C_{2n},\alpha_{2n})$ et 
$(\∆,\lambda_0)=(C_ n,\alpha_n).$ 
\vskip 3mm
Les d\émonstrations sont faites cas par cas ainsi dans le cas ii) 
lorsque $\F$ est un corps $\goth P$-adique, on supposera la
caract\éristique r\ésiduelle diff\érente de $2,$ pour des raisons techniques.\\

\begin{lem}
Soit $x_0\in \goth g^{\lambda_1},$ on suppose que l'application de E dans $ \goth g^{\lambda_2} $ d\éfinie par $ad(A)^2(x_0)$ est surjective alors : 
\begin{enumerate}
 \item  Les orbites  de G et de Ker$(\chi_n)$ dans $\g_1-\g'_1$ sont les m\êmes et ont pour repr\ésentants: $0,\sum_{1\≤j\≤i}X_{\lambda_j}, i=1,..,n-1.$ \\
 
 \item Dans le \PV $(\goth U(\F h_{\lambda_n})_0,\goth U(\F h_{\lambda_n})_1),$ $(Aut_0(\goth U(\F h_{\lambda_n}))_{h_{n-1}}$ a n 
orbites de repr\ésentants: $0,\sum_{1\≤j\≤i}X_{\lambda_j}, i=1,..,n-1$ et $P_{\oplus_{1\≤j\≤n-1}\F h_{\lambda_j}}$ a une seule orbite dans $\goth U(\F h_{\lambda_n})_1".$
\end{enumerate}

\end{lem}

 \dem

 1) Tout \elt $A\in E_{-1,1}^{k,r}$  pour lequel il existe $x\in  \goth g^{\lambda_k}$ tel que $ad(A)^2(x)\not=0$  se compl\ète en un \Sl $(A,h_{\lambda_r} -h_{\lambda_k},B)$ avec $B\in E_{1,-1}^{k,r}$(cf.  \demo du lemme 2.2.2 de \cite{mullerNAG})  par cons\équent l'automorphisme \él\émentaire: 
 
 $$g_A=expad(B)expad(A)expad(B).g_{k,r}\in G_{\goth t_0}$$
se r\éduit \à l'identit\é sur  $\oplus_{1\≤j\≤n,j\not=k,r} \goth g^{\lambda_j},$ sur $\goth g^{\lambda_j}$ (resp. $\goth g^{\lambda_k}$) \à $\frac{1}{2}ad(A)^2\circ g_{k,r}$  (resp. 
 $\frac{1}{2}ad(B)^2\circ g_{k,r}$).\\
 
2)  Soit $x\in \g_1-\g'_1,$ on peut supposer que $x=\sum_{1\≤i\≤j}x_i,x_i\in \goth g^{\lambda_i}-\{0\}$ (4 prop.6.1.1). Soit $i\in \{1,...,j\},$   par 1) il existe $A_i\in E_{1,-1}^{i,n}$ et une involution $g_{A_i}\in G_{\goth t_0}\cap G_e$ r\éduite \à l'identit\é sur $\oplus_{1\≤k\≤n-1,k\not=i} \goth g^{\lambda_k}$ telle que $g_{A_i}(x)=\sum_{1\≤k\≤j,k\not=i}x_k+X_{\lambda_i}$ puisque l'application est surjective d'o\ù 1. et 2. puisque la restriction des $g_{A_i}$ \à $\goth U(\F h_{\lambda_n})$ est dans $Aut_0(\goth U(\F h_{\lambda_n})).$\fdem
 
 \begin{prop}

 \begin{enumerate}
 \item   Dans le cas r\éel ou bien dans le cas $\goth P$-adique  lorsque  $(\overline \∆, \alpha_0)=(A_{2nm-1},\alpha_{nm})$ et 
$(\∆,\lambda_0)=(A_{2n-1},\alpha_n)$ avec $m\≥2,$ $P_0$   a une seule orbite dans $ \g"_1$ et $G$  a $n+1$ orbites dans $\g_1$ de repr\ésentants: $0,\sum_{1\≤j\≤i}X_{\lambda_j}, i=1,..,n.$ Les orbites de Ker$(\chi_n)$ dans $\g_1-\g'_1$ sont les m\êmes que celles de $G.$\\
 
 \item
  Lorsque $\F$ est un corps $\goth P$-adique et que $(\overline
\∆,\alpha_0)=(C_{2n},\alpha_{2n})$ et  $(\∆,\lambda_0)=(C_n,\alpha_n),$  $G $ a
  trois orbites dans  $ \g'_1$ lorsque  $n=1$  et
 pour $n\≥2$ il  y a  quatre orbites  en bijection avec $\F^*/\F^{*2}$.
 
 \noi Les orbites de $P_0$  dans $ \g"_1$ (resp.$ \g"_{-1}$) sont donn\ées par $O_u$ (resp.$O^*_u$) avec $u$ d\écrivant $ (F_1(\goth g^{\lambda_1}-\{0\})/\F^{*2})^n.$   \end{enumerate}
 \end{prop}
 
 \dem
 
 1)  Lorsque $\goth g$ est de rang $1,$
on  montre que $Ker (\chi_1)$ agit transitivement  sur $\{x\in \goth g_1|F_1(x)=u\},$ $u\in \F^*$ \étant  repr\ésent\é par $F_1.$ \\

a) Lorsque $F_1$ est une forme quadratique (donc anisotrope), on sait que $SO(P)$
agit transitivement sur $\{x\in \goth g_1 / P(x)=t\}.$  Soit $U$ le
sous-groupe alg\ébrique connexe de $Aut(\goth g),$ d'alg\èbre de Lie
$\goth g'_0=[\goth g_0,\goth g_0],$   alors $U\subset (Ker\chi_1)^0.$

\noi On consid\ère l'application de restriction, not\ée $f,$ de $U$ dans 
$O(P),$ qui \à $g$ de $U$ associe la restriction de $g$ \à $\goth g_1.$
  $f$ \étant continue , $f(U)$ est un sous-groupe connexe  
de $O(P)$ donc un sous-groupe de $SO(P)$ , comme le noyau de $f$
est r\éduit \à l'identit\é , on a : $$ dim\ \goth L(U)=
dim\ \goth L(f(U))=dim\ \goth g'_0$$
Il suffit donc de v\érifier que $$(1)\qquad
dim\ \goth g'_0 =dim\ (\goth L(SO(P))) $$
pour avoir $\overline {f(U)}=\overline {SO(P)},$ ainsi tout
\él\ément de $SO(P)$ se prolonge en un  \elt de
$Aut_0(\overline{\goth g}) $ normalisant $\goth g_1$ , donc par
dualit\é , il normalise \également $\goth g_{-1},$ et par
engendrement \également $\goth g .$ Ainsi tout  \elt de $SO(P)$
se prolonge en un \él\ément de $Ker (\chi_1).$

\noi Notons que lorsque $\F=\R,$ $U$ est un groupe compact donc
$f(U)$ \également ainsi $f(U)$ est un sous-groupe analytique
compact de $SO(P),$ donc $f(U)$ est un sous-groupe alg\ébrique de
$SO(P)$ ayant la m\ême alg\èbre de Lie donc $f(U)=SO(P).$

\noi Or $$ \text{dim}\ (\goth L(SO(P)))=\frac
{\text{dim}(\goth g_1)(\text{dim} (\goth g_1)-1)}{  2}$$ La v\érification de
(1) est imm\édiate \à l'aide des tables de \cite{warner}, \cite{veisfeiler}
et du tableau
$1,$ en effet les diagrammes de Satake correspondants aux cas absolument
irr\éductibles, commutatifs, de $\F$-rang un, tr\ès r\éguliers,
ayant un invariant relatif fondamental de degr\é deux, sont les
suivants :\\

\hskip 10pt  \hbox to 1,5 cm {\offinterlineskip \lower 2pt\hbox{$\bullet$} \hglue -3,2pt
{\vrule height  0,4pt depth 0pt width 0,5 cm}\lower 2pt\hbox{$ \circledcirc$}
\hglue -3,2pt
{\vrule height  0,4pt depth 0pt width 0,5 cm}\lower 2pt\hbox{$\bullet$}}
 \hskip  4cm dim ${\goth g}'_0 =6$ , dim $\goth g_1 =4,$\\

\hskip 10pt \hbox to 4,2 cm {\offinterlineskip \lower 2pt\hbox{$ \circledcirc$} \hglue -3,2pt
{\vrule height  0,4pt depth 0pt width 0,5 cm}\lower 2pt\hbox{$ \bullet$}
\hglue -3,2pt
{\vrule height  0,4pt depth 0pt width 0,5 cm}\lower 2pt\hbox{$\bullet$}
  \dotfill \hbox to 1,3 cm {\offinterlineskip\lower 2pt\hbox{$ \bullet$}\kern -1pt   \hrulefill\kern -3,4pt\lower 2pt 
\hbox{  \offinterlineskip  {$\bullet$} \hglue -7pt\vbox{ {\hrule height 0,3pt width 0,6cm}\vskip 3,5pt{\hrule height 0,3pt width 0,6cm}
\vskip 0,3pt} \hglue -16pt $>$   \hglue -1,8pt{$ \bullet$}}}} \hskip  1,3cm dim $({\goth g}'_0
)=(n-1)(2n-1)$ , dim$ ({\goth g}_1 )=2n-1,$ \\

\hskip 10pt \hbox to 5 cm {\offinterlineskip \lower 2pt\hbox{$\circledcirc$}
\hglue -3pt{\vrule height .4pt depth 0pt width 0,5cm}\lower 2pt\hbox{$ \bullet$}
\hglue -7pt{ \vrule height .4pt depth 0pt width 0,5cm}\lower 2pt\hbox{$\bullet$}
 \dotfill \hbox to 2 cm {\lower 2pt\hbox{$ \bullet$}\hglue -1,3pt
\hrulefill\lower 2pt
\hbox{$\bullet$}\hrulefill \lower 2pt \vtop {\offinterlineskip \hbox{$ \bullet$} \hbox to 5pt{\hfill \vrule height 12pt width 0,3pt\hfill} \hbox{$ \bullet$}}\hrulefill\lower 2pt\hbox{$\bullet$} } }\hskip  0,6cm  dim $({\goth g}'_0)=(n-1)(2n-3)$ , dim$ ({\goth g}_1 )=2(n-1),$\\

\hskip 10pt \hbox to 1,2 cm { \offinterlineskip  {$\circledcirc$} \hglue -7pt\vbox{{\hrule height 0,3pt width 0,6cm}\vskip 3,5pt{\hrule height 0,3pt width 0,6cm}\vskip 0,3pt }\hglue -12pt $>$   \hglue -2pt $ \bullet $}
\hskip  4,3cm  dim $({\goth g}'_0)=$dim $({\goth
g}_1)=3.$\\

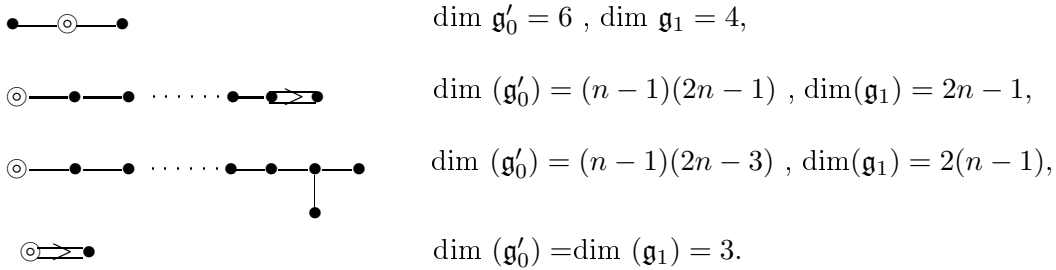
\qquad b) Lorsque $F_1$ n'est pas une forme quadratique , $\goth g$
est isomorphe \à $sl_2(\D),$ $\D$ \étant une alg\èbre \à division et
$\F$ est un corps p-adique.\\

\noi Donnons une br\ève description de $(\goth g_0,\goth g_1)$ lorsque
$\goth g$ s'identifie \à $sl_{2n}(\D).$ Ce cas est le m\ême que
$sl_{2n}(\F).$ On a
$$H_0=\left(\begin{array}{cc} I_n&O_n\\
0_n&-I_n\end{array}\right)$$
 $$\goth g_1=\{ \left(\begin{array}{cc}0_n&A_n\\
0_n&0_n \end{array}\right)\quad \hbox{avec}\ \ A_n\in M_{n\times n}(\D)\}\quad \goth
g_{-1}=\{ \left(\begin{array}{cc}0_n&0_n\cr
B_n&0_n \end{array}\right)\quad \hbox{avec}\ \ B_n\in M_{n\times n}(\D)\}$$
$\goth U(\F(2H_0-h_{\lambda_{n-i+1}})$ est l'ensemble des matrices ayant tous ses
\él\éments nuls sauf ceux dont les indices des lignes et des colonnes
appartiennent \à l'ensemble $\{n-i+1,n+i\}.$ $E$ 
est l'ensemble des matrices ayant tous ses
 \elts nuls sauf celui situ\é sur la $n-1$
ligne et $n$ colonne ou bien celui  situ\é sur la $n+1$ ligne et $n+2$
colonne.\\

\noi Dans le cas $n=1,$ $\goth g_1$ s'identifie \à $\D,$ l'invariant est alors
donn\é par la norme de $\D,$ not\ée $N,$ $G$ s'identifie \à $\D^*\times
\D^*$ op\érant par $x\rightarrow axb^{-1}.$ $U$ contient 
$\{x\in \D\ /\ N(x)=1\}$ qui op\ère transitivement sur l'ensemble 
$\{x\in \D\ /\ N(x)=t\}, t\in \F^*.$
 Notons qu'ici $G$ a deux orbites dans $\goth
g_1.$ \fdem du rang $1.$\\

\qquad 2) \'Evident dans le cas r\éel  puisque  $\chi_i(G_{\goth t_0 })=\R^{*+}$ et $F_i(x)\≥0$ pour $i=1,...,n,$ de plus $G_{\goth t_0}$ (resp. Ker$(\chi_n)$) op\ère transitivement sur $ \oplus_{1\≤i\≤j}{\goth g}^{\lambda_i}-0$ (resp. lorsque$j<n$) par 1) et toute $G-$ orbite   non r\éduite \à $\{0\}$ rencontre $ \oplus_{1\≤i\≤j}{\goth g}^{\lambda_i}-0.$\\

  \qquad 3) Pour le cas  $sl_{2n}(\D),$ $\D$
\étant une alg\èbre \à division. \\

\noi V\érifions que l'hypoth\èse du lemme 6.1.6 est v\érifi\ée et pour ceci on
peut supposer $n=2.$ 

\noi On a :
$$E=\{A= \left(\begin{array}{cc}{\left(\begin{array}{cc} 0&0\\
a&0\end{array}\right )}&0\\
0& {\left(\begin{array}{cc} 0&0\\
b&0\end{array}\right )} \end{array}\right ) \quad a,b\in \D\}$$
et pour ${\tilde x}= \left(\begin{array}{cc}0&{ \left(\begin{array}{cc}0&x\\
0&0  \end{array}\right )}\\ 0&0  \end{array}\right )$  on a
 $ad(A)^2({\tilde x})= \left(\begin{array}{cc}0& {\left(\begin{array}{cc}0& 0\\
2axb&0  \end{array}\right )}\\0&0 \end{array}\right )$
ainsi cette application est surjective sur $\goth g^{\lambda_2}$ lorsque $x$
est non nul. On peut appliquer au  \PV $(\goth g_0,\goth g_1)$
le 1) du lemme 6.1.6.\\

\noi Comme on  peut toujours plonger  $sl_{2n}(\D)$ dans  $sl_{2n+2}(\D)$ de
la mani\ère suivante :
$$A\rightarrow {\tilde A}=\left ( \begin{array}{ccc}0&...&0\\
.&A&.\\
0&...&0 \end{array}\right ).$$ 
On peut appliquer au  \PV $(\goth g_0,\goth g_1)$
le 2) du lemme 6.1.6.

\noi  Notons que dans le cas $\goth p$-adique, $\chi_i(G)=\F^*$ pour $i=1,...,n$ car la
norme est surjective sur $\F.$\\

\qquad 4) Le cas $(C_{2n},\alpha_{2n})$ non d\éploy\é

La  \demo se trouve dans le $\S 6.1.3.$  \\

\subsubsection{Le cas $\bf (C_{2n},\alpha_{2n})$  de type III  }
\bigskip

\noi Le cas trait\é correspond au diagramme de Satake suivant:
\hskip 10pt \hbox to 4,2 cm {\offinterlineskip \lower 2pt\hbox{$\bullet$} \hglue -3,2pt
{\vrule height  0,4pt depth 0pt width 0,5 cm}\lower 2pt\hbox{$\circ$}
\hglue -3,2pt
{\vrule height  0,4pt depth 0pt width 0,5 cm}\lower 2pt\hbox{$\bullet$}
  \dotfill \hbox to 1,3 cm   {\offinterlineskip\lower 2pt\hbox{$\circ$}\kern -1pt   \hrulefill\kern -3,4pt\lower 2pt 
\hbox{  \offinterlineskip  {$\bullet$} \hglue -7pt\vbox{ {\hrule height 0,3pt width 0,6cm}\vskip 3,5pt{\hrule height 0,3pt width 0,6cm}
\vskip 0,3pt} \hglue -16pt $<$   \hglue -1,8pt{$\circledcirc$}}}}\\

 \noi Lorsque $\F$ est un corps local non
archim\édien, on suppose  {\bf la caract\éristique r\ésiduelle  
diff\érente de $2.$}\\

 $\bullet$ {\bf  La description du cas $(C_l,\alpha_l)$ d\éploy\é} 
\vskip 3mm
\noi 
$\goth g$ s'identifie \à l'ensemble des matrices de la forme:
$
\left ( \begin{array}{cc}
A&B\\
C&D\end{array}\right)
$
o\ù $A,B,C,D$ sont des matrices carr\ées d'ordre $l$ telles que
 $D=-{^tA},B$ et $C$ sont sym\étriques.(cf.N.Jacobson,Lie Algebras,Interscience,1962) \\

\noi La forme de Killing est donn\ée par la relation suivante:$
B(U,V)=(2l+2)\text{trace}(U.V)
$
 
\noi $2H_0$ est alors donn\é par: $
2H_0= \left ( \begin{array}{cc}
I_l&0_l\\
0_l&-I_l\end{array}\right),
$
$O_l$ \étant la matrice carr\ée d'ordre $l$ dont tous les coefficients sont
nuls et $I_l$ la matrice carr\ée d'ordre $l$ dont tous les coefficients sont
nuls sauf ceux situ\és sur la premi\ère diagonale qui sont \égaux \à $1$ 

\noi Un calcul facile donne :
$$
\goth g_1=\lbrace  \left( \begin{array}{cc}0_l&U\\
                                0_l&0_l \end{array}\right)\quad U=^tU \in \goth M_{l,l}
 \rbrace \quad  
\goth g_{-1}=\lbrace  \left( \begin{array}{cc}0_l&0_l\\
                                U_l&0_l\end{array}\right)\quad U=^tU \in \goth M_{l,l}
  \rbrace\ 
$$
On identifie ainsi $\goth g_1$ et $\goth g_{-1}$  \à l'espace des
matrices sym\étriques ayant $l$ lignes.                         

\noi L'action de $Aut(\goth g)_{H_0}=Aut_0(\goth g)_{H_0}$ (\cite{bourbakigal8}, chapitre 8, n$^\circ 3$, corollaire $2$) (resp.$G_e$) s'identifie \à
l'action de l'ensemble des matrices :

$$
\lbrace g=  \left( \begin{array}{cc}A&0_l\\      
                    0_l&\mu{^tA^{-1}} \end{array}\right) \text{avec}
 \quad  \mu\in {\F^*}\quad A\quad \text{inversible}
\rbrace\quad
 \hbox{(resp: }\ \mu=1\  \hbox{et} \ det(A)=1)$$

\noi agissant par la conjugaison usuelle , ce qui donne avec les notations 
ci-dessus:
 $$
g(U)=µ^{-1}AU({^tA})
$$ 
on a ainsi, modulo l'action du centre, l'action usuelle de
$Gl_l(\F)$ sur les matrices  sym\étriques.\\

\noi$\goth a$ peut \être choisi comme l'ensemble des matrices
diagonales de $\goth g,$ et $P_0$ correspond aux \elts $g$ pour lesquels la matrice  $A$ est triangulaire sup\érieure.\\

\noi Pour $X\in \goth M_{l,l}, X_i$ d\ésigne la matrice tronqu\ée \à $i,$ lignes et $i$ colonnes obtenues en conservant les lignes et les colonnes comprises entre $l-i+1$ et $l.$ 

\noi Les invariants relatifs fondamentaux (non normalis\és) et les caract\ères associ\és sont alors
donn\és par:$$
F_i(U)=det(U_i)\qquad \chi(g)=µ^{-l}(det(A_i))^2.
$$\\

{\bf  $\bullet$ La description du cas $(C_{2n},\alpha_{2n})$  de type III}
 \bigskip
 
\noi Soit $\E$ l'extension galoisienne sur laquelle $\goth g$ se
d\éploie , on a $\E=\F(\sqrt \epsilon),$ $\epsilon$ \étant
une unit\é qui n'est pas un carr\é (\cite{veisfeiler}).

\noi On note par $\overline x$ la conjugaison de $x$ dans $\E$ et $(.,.)$
le symbole de Hilbert d\éfini sur $\F^*\times \F^*.$

\noi Si on identifie $\goth g\otimes_{\F}\E$ \à l'espace des matrices
donn\é pr\éc\édemment, $\goth g$ s'identifie alors \à $$
\{X\in {\goth g\otimes_{\F}\E}\ /\ T{\overline X}T^{-1}=X\}\quad \text{avec}$$
  $$T= \left (\begin{array}{cc}J&0_l\\
                     0_l&^tJ\end{array}\right ) \quad  ,\quad 
 J= \left (\begin{array}{ccc} 
  I_\beta&0&0\\
                 0& I_\beta&0\\
                ...&...&...\\
0&...0& I_\beta  \end{array}\right )\quad ,\quad 
I_\beta= \left(\begin{array}{cc}0&\beta\cr1&0\end{array}\right )\quad \text{avec}\quad ( \beta,\epsilon)=-1.$$
On peut remarquer que
$J^2= \beta Id.$\\

\noi   
En particulier $\goth g_1$ (resp.$\goth g_{-1}$) s'identifie \à l'espace des matrices 
sym\étriques de $\goth M_{2n,2n}(\E)$ v\érifiant la relation 
$$
J{\overline U}=U(^tJ)\quad (\hbox{resp.\ }(^tJ){\overline U}=UJ\ )
$$
On note avec un $\E$ en indice, tous les groupes pr\éc\édents
provenants de  $\goth g\otimes_{\F}\E.$\\
 
\noi Comme $\goth g_{\pm1}$ engendrent $\goth g,$ on a $$G_{\E}\cap
Aut(\goth g)=\{g\in G_{\E}\ /g(\goth g_1)\subset \goth
g_1\}$$ et son action s'identifie \à l'action des matrices $$
\{g=  \left (\begin{array}{cc}A&0_l\\     
                    0_l&\mu{^tA^{-1}}  \end{array}\right ) \ \text{avec}
 \quad  \mu\in {\F^*}\ ,\  J{\overline A}=\pm AJ\quad  A\
\text{inversible}\} $$
agissant sur $\goth g_1$ (resp.$\goth g_{-1}$) apr\ès identification matricielle par: $$
g(U)=\mu^{-1}AU(^tA)\quad (\hbox{resp.\ }g(U)=\mu ^tA^{-1}UA^{-1}\ ).
$$
Notons que  
$$\widetilde {G_{2n}(\E)}=\{ A\in G_{2n}(\E)\quad J{\overline
A}=\pm AJ\}=G_1\cup \sqrt\epsilon G_1 \quad
avec \quad G_1=\{ A\in G_{2n}(\E)\quad J{\overline A}=AJ\}.$$

\noi Donnons la description de chaque $h_{\lambda_{n-i+1}},1\≤i\≤n,$ ils sont
repr\ésent\és par les matrices 
$$h_{\lambda_{n-i+1}}=(h_{k,l}^{(i)})$$ 
avec 
$$ h_{k,l}^{(i)}=0\quad \text{si} \quad k\not=l\quad  \text{ou}\quad  \text{si}\quad 
k=l\notin\{2i-1,2i,2(n+i)-1,2(n+i)\}$$ 
$$ h_{2i-1,2i-1}^{(i)}=
h_{2i,2i}^{(i)}=1=-h_{2(n+i)-1,2(n+i)-1}^{(i)}=
-h_{2(n+i),2(n+i)}^{(i)}.$$
La sous-alg\èbre $\goth U(\F(2H_0-h_{\lambda_{n-i+1}})$ correspond aux 
matrices $U^{(i)}$ de $\goth g$ dont tous les \elts
$u_{k,l}$ sont nuls , sauf ceux tels que $k$ et $l$
appartiennent \à l'ensemble $\{2i-1,2i,2(n+i)-1,2(n+i)\}.$
Les sous-espaces $\goth g^{\lambda_{n-i+1}}$ ( resp.$\goth g^{-\lambda_{n-i+1}}$) correspondants sont les matrices $U^{(i)}$ de $\goth g_1$ dont tous les
 \elts $u_{k,l}$ sont nuls , sauf ceux tels que $k$ et $l$
appartiennent \à l'ensemble $\{2i-1,2i\}$ ( resp.
$\{2(n+i)-1,2(n+i)\}).$  En particulier, l'\écriture sous forme de
matrices $(2,2)$ donne  
$$ \hbox{pour}\quad \goth g^{{\lambda_{i}}}\quad :\quad
{\D}=\{X_{a,b}= \left (\begin{array}{cc}\beta \overline a&b\\      
                    b&a \end{array}\right )\quad  \hbox{avec}
 \quad  a\in \E\quad b\in \F\}$$
$$ \hbox{pour}\quad \goth g^{-{\lambda_{i}}}\quad :\quad
{\D}'=\{Y_{a,b}= \left (\begin{array}{cc}a&b\\     
                    b&\beta{\overline a} \end{array}\right )\quad  \hbox{avec}
 \quad  a\in \E\quad b\in \F\}$$
 On a: $$det(X_{a,b})={\overline a}a\beta-b^2=(x^2-\epsilon
y^2)\beta-b^2\qquad a=x+\sqrt{\epsilon}y,\ x,y\in \F$$ 
ce qui donne lorsque $\F=\R$  $(\epsilon=\beta=-1)$ :  
$det(X_{a,b})=-(x^2+y^2+b^2).$\\

\noi La normalisation des invariants relatifs fondamentaux   d\épend du choix du $sl_2$-triplet
de base : $(X_{\lambda_1},h_{\lambda_1},X_{-{\lambda_1}}).$ Dans cet exemple on convient de
prendre :
$$X_{{\lambda_1}} \ \ \hbox{correspond \à}\ \ X_{0,1}\ \  \hbox{alors}\ \ 
X_{-{\lambda_1}} \ \ \hbox{correspond \à}\ \ Y_{0,-1}$$
ainsi \à l'identification matricielle pr\ès: 
$$F_{\goth U(\F(2H_0-h_{\lambda_1}))}(X_{a,b})=-det(X_{a,b})\ \ \hbox{ et}\ \ F^*_{\goth
U(\F(2H_0-h_{\lambda_1}))}(Y_{a,b})=-det(Y_{a,b})$$
c'est \à dire que l'un des coefficients de la forme quadratique anisotrope \à $3$ variables et r\éduite sous forme de carr\és, est un carr\é alors $\epsilon=-$discriminant$(F_1).$

\bigskip

$G_{\goth t_0}$ correspond aux matrices dont les \elts $A$ sont de la forme:
$A=  \left (\begin{array}{ccc} 
   A_1&0&0\\
                 0& A_ 2&0\\
                ...&...&...\\
                 0&...0& A_n  \end{array}\right ),$
                 (resp. $P_0$ correspond aux matrices dont les \elts $A$ sont de la forme:
$A=  \left (\begin{array}{ccc} 
   A_1&*&*\\
                 0& A_ 2&*\\
                ...&...&...\\
0&...0& A_n  \end{array}\right )\quad )$
les matrices $A_i$,$i=1,...,n,$ ayant $2$ lignes et $2$ colonnes.  \\ \\
\noi Tout \elt de $\widetilde G_2\cap G_1$ s'\étend en un \elt de $G_{\goth t_0}\cap G_1$ centralisant $\oplus_{1\≤j\≤n,j\not=i}\goth g^{\lambda_j}$puisque toute matrice inversible \à $2$ lignes et $2$ colonnes v\érifiant l'\égalit\é $I_{\beta}\overline A=AI_\beta$ s'\étend en un \elt  $\tilde A$ par: $$A\rightarrow \tilde A= \left (\begin{array}{ccccc} 
   I_2& 0&...&...&0\\...&...&...&...&...\\
                 0&...& A&...&0\\
                 ...&...&...&...&...\\
                 ...&...&...&...&I_2  \end{array}\right ).$$

\noi  Description compl\ète du cas $\bf (C_4,\alpha_4):$\\
 
 $$\goth g=\{ \left (\begin{array}{cccc}
 A_1&A_3&B_1&B_3\\
             A_4&A_2&^tB_3&B_2\\
             C_1&C_3&-^tA_1&-^tA_4\\
             ^tC_3&C_2&-^tA_3&-^tA_2 \end{array}\right ),
\quad A_i,B_i,C_i,B,C\in {\goth M}_{2,2}(\E)\quad \hbox{avec} :  $$
$$ B_1\ ,\ B_2\ ,\ C_1\ ,\ C_2\  \ \hbox{sym\étriques} \ \hbox{et}:$$
$$1\≤i\≤4:\ \   I_\beta{\overline A_i}=A_iI_\beta\ ,\ 
 1\≤j\≤3: \ I_\beta{\overline B_j}=B_j(^tI_\beta) \ \hbox{et}\ 
(^tI_\beta){\overline C_j}=C_jI_\beta\ \}
$$

\noi $\goth g^{\lambda_1}$ ( resp.$\goth g^{-\lambda_1}$) correspond au
sous-ensemble pr\éc\édent dont tous les  \elts sont nuls sauf
peut-\être $B_2$ ( resp. $C_2$). $\goth U(\F(h_1)$ (resp.
$\goth U(\F(h_2)$ ) 
correspond au sous-ensemble pr\éc\édent dont tous les  \elts sont
 nuls sauf
peut-\être ceux qui portent un indice $1$ ( resp. $2$). $E=E_{-1,1}^{1,2}$
correspond \à l'ensemble des matrices repr\ésent\é par $A_3.$

\noi L'action de $E$ sur $\goth g^{\lambda_1}$ d\éfinie par 
 $\frac{1}{2}ad(\ )^2$  est repr\ésent\ée matriciellement par  $B_1:=A_3B_2({^tA_3})$   lorsque l'\elt de $E$ est repr\ésent\é par $A_3,$ celui de  $\goth g^{\lambda_1}$ par $B_2$ et celui de $\goth g^{\lambda_2}$ par $B_1.$\vskip 3mm

 \noi Pour $x\in \goth g^{\lambda_1},$ $y\in \goth g^{-\lambda_2},$  
 la forme quadratique 
$f_{x,y}$ est donn\ée dans la description
des matrices ci-dessus par : $$X\in \D\quad Y\in \D'\quad A\in
 \mathbb{H} \qquad f_{X,Y}(A)=cste .\text{Trace}(AX(^tA)Y)$$
$$\hbox{avec}\ \ \mathbb{H}=\{  \left (\begin{array}{cc}a&\beta{\overline b}\\
b&{\overline a}\end{array}\right ) \quad a,b\in \E\}=G_1\cup\{  \left (\begin{array}{cc}0&0\\
0&0 \end{array}\right ).$$ 
 \vskip 3mm

$\bullet$ On rappelle que  les orbites du cas de rang 1,
$(C_2,\alpha_2),$ sont en bijection avec $F_1(\goth g^{\lambda_1}-0)/{\F^{*2}}.$\\
 
\noi Lorsque $\F$ est un corps $\goth P$-adique,  
 les repr\ésentants de
$F_1(\goth g^{\lambda_1}-0)/{\F^{*2}}$ contiennent $ \{1,\beta,\epsilon\beta\},$  il y
a \égalit\é lorsque
la
caract\éristique r\ésiduelle est diff\érente de $2$  et $ X_{0,1},$ $X_{1,0},$ $X_{\alpha,0}$ (resp.$
 Y_{0,1},$ $Y_{1,0},$ $Y_{\alpha,0}),$ avec  $\alpha{\overline \alpha}=\epsilon,$ sont des repr\ésentants des
orbites non r\éduites \à $0$ dans $\D$ (resp. $ \D'$).\\
\noi De plus, comme le produit de deux  \elts de  
$F_1(\goth g^{\lambda_1}-0)/{\F^{*2}}$ d\écrit $\F^*/\F^{*2},$  $F_i$ d\écrit $\F^*$  pour $i\≥2.$  \\

\noi Lorsque deux matrices $X$ et $X'$ de ${\D}$ v\érifient
l'\égalit\é $\det(X)\equiv \det(X')$ (modulo $\F^{*2}$), il existe 
$A$ dans $\widetilde {G_2(\E)}\cap G_1$ tel que $gX(^tg)=X'$ puisque $X$ et $\epsilon X$ sont dans la m\ême $G_1-$orbite et que $X$ et $X'$ sont dans la m\ême $\widetilde {G_2(\E)}$ (cf.  le cas
 de rang $1$).
\\

\noi Soit $x\in O_u$ suppos\é non vide, alors il existe $g\in N_0$ et $x_i\in \g^{\lambda_i},i=1,...,n,$ tels que $x=g(\sum_{1\≤i\≤n}x_i)$ et $G_i(x_i)\equiv u_i$ (mod($\F^{*2}$)) pour $i=1,...,n,$ donc chaque $u_i\in G_i(\g^{\lambda_i})^*/\F^{*2}=F_1(\g^{\lambda_1})^*/\F^{*2}$ (toutes les formes quadratiques $G_i$ sont \équivalentes \à $F_1$ par construction) d'o\ù $O_u$ est une seule $P_0$-orbite par le r\ésultat de rang $1$ rappel\é ci-dessus et parce que
 chaque $\chi_i(P_0)=\F^{*2}$.\\

$\bullet$  { \bf   Fin de la  d\émonstration du 2. de la proposition 6.1.7}
\vskip 5mm
Elle se fait par
r\écurrence sur $n,$ le cas $n=1$ \étant rappel\é ci-dessus.Les orbites de $G$ et de $G_1$ sont donc les  m\êmes.
On peut toujours supposer que $x$ et $x'$ appartiennent \à
$\oplus_{1\≤i\≤n}( \goth g^{\lambda_i}-0 )$ et on note $X_i$
(resp.$X'_i$ ) les diff\érentes composantes. On a :$$
(-1)^nF_n(x)=\prod_{1\≤i\≤n}det(X_i)\equiv
(-1)^nF_n(x')=\prod_{1\≤i\≤n}det(X'_i)\quad (\F^{*2}) $$
 Notons $\{\epsilon_i\ ,i=1,2,3\}=Image(F_1)/{\F^{*2}}$ et soient $n_i$ (resp.
$n'_i$) le nombre de composantes $X_k$ ( resp. $X'_k$ ) de $x$
(resp. de $x'$ ) telles que $\det(X_k)\equiv -\epsilon_i$ $(\F^{*2})$
( res. $\det(X'_k)\equiv -\epsilon_i$ $(\F^{*2}) ).$

\noi Comme on a suppos\é que la caract\éristique r\ésiduelle est
diff\érente de 2, et que $n_1+n_2+n_3=n'_1+n'_2+n'_3=n ,$ on obtient
pour tout $i$ : $n_i\equiv n'_i$ (modulo 2) , ainsi on a deux cas :\\

a) il existe $i$ tel que $n_i\equiv n'_i\equiv 1,$ il suffit alors
d'appliquer le r\ésultat du cas $n-1,$ car on peut toujours
supposer que $X_1=X'_1$ ( cf. 1)prop.6.1.1).  \\

b) Pour tout $i,$ on a 
$n_i\equiv n'_i\equiv 0.$ Il reste alors \à montrer que deux
 \elts de la forme $$
  \left (\begin{array}{cc}X&0_2\\ 0_2&X\end{array}\right) \quad et  \  \left (\begin{array}{cc}Y&0_2\\
0_2&Y \end{array}\right) $$
sont dans la m\ême orbite lorsque $\det(X)\det(Y)$ n'est pas un
carr\é, ce qui termine la d\émonstration du cas $n=2.$ Il suffit
ensuite d'appliquer l'hypoth\èse de r\écurrence pour le cas $n-2.$
\vskip 3mm
\noi D\émonstration du cas $n=2$ :

i) Consid\érons le polynome du second degr\é, \à coefficient dans $\F,$
d\éfini par $$ l(u,v)=\det(uX-Y)-v^2\det(X)\quad u,v\in {\F}
$$
On peut l'\écrire $$
 \begin{array}{lll}l(u,v)&=&\det(X)(u^2-v^2)+2u\det(X)C(X,Y)+\det(Y)
\\
&=&\det(X)(u+C(X,Y)  )^2-v^2\det(X)+\det(Y)-\det(X)
C(X,Y) ^2\end{array}$$

\noi Comme la forme quadratique $u^2-v^2$ est isotrope , elle
repr\ésente $\F$ d'o\ù il existe un couple $(u,v)$ tel que
$l(u,v)=0.$ 

\noi Soit $(u,v)$ tel que $\det(uX-Y)=v^2\det(X),$ n\écessairement
$uv\not=0$ car $X$ et $Y$ ne sont pas dans la m\ême orbite ainsi il
existe un  \elt$B$ de $\widetilde{GL_2(\E)}\cap G_1$ tel que $$
Y=uX+B(uX)(^tB).$$

ii) Soit $A$ la matrice de $\goth M_{4,4}(\E)$ d\éfinie par :$$
A= \left (\begin{array}{cc}Id_2&B\\ -X(^tB)X^{-1}
&Id_2 \end{array}\right)\qquad  \hbox{alors:}$$
 $$ 
A \left (\begin{array}{cc}uX&0_2\\0_2&uX\end{array}\right){^tA}= \left (\begin{array}{cc}Y&0_2\cr
0_2&Y'\end{array}\right)\quad \hbox{avec}\quad Y'\not=0_2 $$
 donc $\det A\not=0$ et $\det(Y)\equiv \det (Y') \ (\F^{*2}),$ on v\érifie facilement que $A$ est un  \elt de
 $\widetilde{G_4(\E)}\cap G_1.$ \fdem
\vskip 5mm
\noi $\bullet$ Indiquons \également les caract\éristiques de $f_{x,y}:$\\

 \begin{lem}
 La classe d'\équivalence de $f_{X,Y},X\in\D-0,$ $Y\in
\D'-0,$ ne d\épend que de l'orbite de $X$ dans $\D$ et de l'orbite
de $Y$ dans $\D'$ et $f_{X,Y}$ a pour  discriminant un.

Lorsque $\F=\R,$  $f\sim f_{X,Y}$ a pour signature (2,2).

Lorsque $\F$ est un corps $\goth P$-adique de  caract\éristique
r\ésiduelle  diff\érente de deux, l'invariant de Hasse
 de $f_{X,Y}$ vaut $-(   \delta_1,\det(X).\det(Y)\ ).(-\det(X),-\det(Y) \ )=\gamma (f_{X,Y}),$ 
avec $\delta_1=-\text{disc}(F_1).$ \end{lem}
\vskip 3mm
 \dem
 
 \noi Il suffit de faire la d\émonstration pour la forme quadratique
$$f_{X,Y}(A)=c\text{Trace}(AX(^tA)Y)\quad c\not=0\quad A\in  \mathbb{H}$$
Un calcul imm\édiat donne $f_{BX^tB,^tCYC}(A)=f_{X,Y}(CAB)$ d'o\ù la
premi\ère assertion.

\noi Le discriminant (not\é "disc") de cette forme quadratique est un polynome en $X$ et $Y,$
relativement invariant donc \à $Y$ ( resp. $X$) fix\é, c'est un
invariant relatif de $\goth g^{\lambda_1}$ (resp. $\goth g^{-\lambda_2}$ ) d'o\ù 
$$\text{disc}(f_{X,Y})=cste.(\det(X))^p.(\det(Y))^q$$
Comme
$$\text{disc}(tf_{X,Y})=\text{disc}(f_{tX,Y})=\text{disc}(f_{X,tY})=t^4\text{disc}(f_{X,Y})$$
on a $p=q=2$ ainsi $f_{X,Y}$ a un discriminant constant qu'il reste \à
calculer.   \\

\noi Rappelons que par normalisation, la forme quadratique anisotrope  $F_1$ est \équivalente \à  $X^2+aY^2+bZ^2$ avec $\text{disc}(F_1)\equiv ab\equiv -\epsilon.$\\

Prenons $X=X_{0,1}$ et
$Y=Y_{0,1},$ le calcul donne: $$f_{X_{0,1},Y_{0,1}}(A)=2cste( {a\overline
a}+\beta b{\overline
b})\quad \hbox{avec}\quad A= \left (\begin{array}{cc}a&\beta{\overline b}\\ b&{\overline a}\end{array}\right)$$ d'o\ù le
discriminant de $f_{X,Y}$ vaut $1.$

1) Lorsque $\F=\R,$  $f\sim f_{X,Y}$ a pour signature $(2,2).$

2) Lorsque $\F$ est un corps $\goth P$-adique de  caract\éristique
r\ésiduelle  diff\érente de deux, on a $\gamma (f_{X,Y})=h_{X,Y},$  $h_{X,Y}$ d\ésignant l'invariant de Hasse de $f_{X,Y},$ il y a seulement deux classes
d'\équivalence de telles formes quadratiques selon qu'elle repr\ésente
$0$ (alors $h_{X,Y}=1 )$ ou non (alors 
  $h_{X,Y}=-1)$ (\cite{o'meara}).

\noi Soit $h'(X,Y)=-(\epsilon, \det(X)\det(Y)\ ).(- \det(X),- \det(Y) \ ),$  on v\érifie que $h_{X,Y}=h'( X,Y)$ dans chaque cas.  \\

\noi On a
$$h'(X_{d,0},Y_{d,0})= -(-\beta d{\overline
d},-\beta d{\overline d})= -(-\beta d{\overline
d},-1)= -(\beta d{\overline
d},-1) $$ 
$$\hbox{et}\quad f_{X_{d,0},Y_{d,0}}(
A)=\beta d{\overline d}f_1(a)+f_1(\beta {\overline d}b) \quad avec  \
f_1(a)=f_1(a_1+\sqrt\epsilon a_2)=a^2+{\overline
a}^2=2(a_1^2+\epsilon a_2^2)$$
 d'o\ù $f_{X_{d,0},Y_{d,0}}$ 
repr\ésente $0$ si et seulement si il existe $a$ dans
$E^*$ tel que $-\beta d{\overline d}f_1(a)$ soit repr\ésent\é par $f_1$ c'est \à
dire que $$h_{X_{d,0},Y_{d,0}}=(-\beta d{\overline d} 
,-\epsilon)=-(\beta d{\overline
d},-1) $$ 
 parce que la caract\éristique r\ésiduelle est
diff\érente de $2$ ( rappelons que dans ce cas $(-1,-1)=(\epsilon
,-1)=1$).

\noi On a:
$$f_{X_{\alpha,0},Y_{1,0}}(A)=\beta f_2({\overline a})+f_2(\beta {\overline  b})\quad (\hbox{resp.\ } f_{X_{1,0},Y_{\alpha,0}}(A)=\beta f_2( a)+f_2(\beta {\overline  b})\ ) $$
avec $ f_2(a)=\alpha a^2+{\overline \alpha}{\overline a}^2.$
Apr\ès r\éduction de $f_2,$ le calcul direct donne
$h_{X_{\alpha,0},Y_{1,0}}=h_{X_{1,0},Y_{\alpha,0}}=(\beta ,-1)=h'(X_{\alpha,0},Y_{1,0})$
et 
$$f_{X_{d,0},Y_{0,1}}(A)=2\beta (d{\overline {ab}}+{\overline d}ab)\quad (\hbox{resp.}\ f_{X_{0,1},Y_{d,1}}(A)=2\beta (ad{\overline {b}}+{\overline ad}b)\ )$$
  s'annule pour $a=d=1$ et $b=\sqrt\epsilon$ donc
$h_{X_{d,0},Y_{0,1}}=h_{X_{0,1},Y_{d,0}}=1=h'(X_{d,0},Y_{0,1}). $
 \fdem  \vskip 3mm
 \noi {\bf Remarque:} La constante $\gamma_{I(k),J(r)}$ donn\ée dans le ii) p.498 de la note aux CRAS Paris (I.Muller- D\écomposition orbitale des espaces \PVs r\éguliers de type parabolique commutatif et application, t.303, S\érie I, n$^\circ11,$ 1986, p.495-498) est erron\ée puisqu'elle est donn\ée par $\gamma (f_{X,Y})$ donc dans les notations de cette note on a $\gamma_{I(k),J(r)}=-(-D,a_{I(k)}a_{J(r)})(a_{I(k)},a_{J(r)}).$  
 
\subsubsection {Calcul de $\gamma_k$}
\vskip 3mm

\noi Pour  $x$ dans $\goth g^{\lambda_1}-\{0\},$  $y$  dans $\goth
g^{-\lambda_2}-\{0\},$ on consid\ère 
la forme quadratique
d\éfinie sur $E $ par
$$f_{x,y}(A)=\frac{1}{2}\tilde B(ad(A)^2(x), y), $$ 
et pour $u,v\in F_1( \goth g^{\lambda_1}-\{0\}),$ soit $\tilde \gamma(u,v):=\gamma (f_{x,y}),$ avec $F_1(x)=u$ et $G_2^*(y)=v,$ ce qui est bien d\éfini puisque $G_{\goth t_0}$ op\ère transitivement sur $\{x\in W_{\goth t_0}|F_i(x)\equiv u_1...u_i\ (\F^{*2})\ ,i=1,...,n\},(u_1,...,u_n)\in \F^{*n}$ (cf.\demo de la prop.6.1.7).\\

\begin{lem}
\begin{enumerate}
\item Pour le type I avec $d\equiv 0\ mod(4),$ ou bien le type II avec $e=0,$ ou bien le type III on a $\delta=1.$ \\

\item Dans le cas r\éel on a:
 $$\gamma (f) := \begin{cases}  1 \  \text{pour le type  III  et le type II avec  } e=0,\\
 (-1)^{\frac{d}{4}}  \  \text{pour le type  I avec }d\equiv 0\ mod(4).
 \end{cases}$$
et dans le cas $\goth p-$adique, $\gamma (f)=-1$ pour le type  III $(C_n)$ et pour le type  I avec  $d=4.$ 

\item $$\tilde \gamma(u,v)= \begin{cases} \gamma (f) \  \text{pour le type  III r\éel et le type III }(A_n),\\
-(  \delta_1 , uv)(u,v)\  \text{pour le type III } (C_n) \ \goth p-\text{adique}.
\end{cases}$$
Pour les types I et II:
$$\tilde \gamma(u,v)= \gamma (uvf) =\gamma (uvf_0)=\overline {\gamma(uvF_2)}  =
 \begin{cases}  1 \  \text{pour le type II avec }e=0 ,\\
 (\delta,uv)\gamma (f)\  \text{pour le type I avec }  d\text{ pair}.
 \end{cases}$$
\item Soient $k,1\≤k\≤n-1,$ $x=x_2+x_0\in  O_{(u_1,...,u_n)}$  avec $x_i\in E_{i}(h_k)\cap \goth g_{1} $ alors
$$ \gamma_k(x_2,x_0)= \begin{cases} \gamma (f)^{k(n-k)} \  \text{pour le type II avec }e=0,\text{le type III r\éel et le type III }(A_n),\\
\\
 \gamma (  u_1u_2f_0) \   \text{pour le type II avec }e>0,\\
\\
\prod_{1\≤i\≤k,k+1\≤j\≤n} \gamma (u_iu_jf) \  \text{pour le type I} , \\
\\
(\delta,F_k(x_2))^n) (\delta, F_k(x_2)F^*_{n-k}(x_0^{-1}))^k\gamma (f)^{k(n-k)}\  \text{pour le type I avec }  d\text{ pair}, \\
\\
(-1)^{k(n-k)}(  \delta_1,  F_k(x_2)^{n-k} F^*_{n-k}(x_0^{-1})^k)( F_k(x_2) , F^*_{n-k}(x_0^{-1}))\\
  \text{pour le type III } (C_n) \ \goth p-\text{adique}.

\end{cases}.$$
\end{enumerate}

\end{lem}

\dem 1) Pour le type III $(C_n)$ on applique le lemme 6.1.9.

2) Lorsque $tf_{x,y}\sim f_{x,y}$ $\forall t\in \F^*,$ on a $d$ pair et
 $\gamma (tf_{x,y})=\gamma (f_{x,y})=( (-1)^\frac{d}{2}\text{disc}(f_{x,y}),t)\gamma (f_{x,y})$ $\forall t\in \F^*,$ donc 
$(-1)^\frac{d}{2}\text{disc}(f_{x,y})\equiv 1\equiv \delta$ (mod($\F^{*2}$)).\\

\noi Cette situation se produit:\\

$\bullet$ Dans le type III ($A_n$) ou bien le type III r\éel (cf.\demo 2)3) de la prop.6.1.8 et le lemme 6.1.9), de plus $f_{x,y}\sim f.$\\

$\bullet$ Dans le type II avec $e=0$ (i)3) du lemme 6.1.3) alors $\gamma(f)=1$ (cf.remarque 6.1.6) et $f_{x,y}\sim f.$\\

$\bullet$ Dans le cas $\goth p-$adique de type I avec $d\equiv 0$ (mod($4$)) car alors $d=4$ donc $f$ repr\ésente $\F^*$ et $tf \sim f  $ $\forall t\in \F^*$ (3) du lemme 6.1.3) d'o\ù $\gamma (f)=-1$ \cite{o'meara}.\\

3) Les autres r\ésultats proviennent de  la dualit\é  pour la forme de Killing. \fdem\\

\subsection{\bf Equations fonctionnelles}

\bigskip

On applique la proposition 5.3.1 dont les hypoth\èses sont v\érifi\ées par le lemme 6.1.10 puisque $\g_1$ est commutative.\\

\subsubsection{Le cas transitif}
\bigskip

\begin{theo}
Pour les types:\\

i) I avec $d=4$ dans le cas $\goth P-$adique,\\

ii) II avec $e=0,$\\

iii) III r\éel,\\

iv) III $(A_n)$ $\goth p-$adique \\

\noi 1) Alors $P_0$ (resp.G) a une seule orbite dans $\goth g"_1$ et $\goth g"_{-1}$ (resp.$\goth g'_1$ et $\goth g'_{-1}$) et 

 \noi $\chi_1(P_0)=...=\chi_n(P_0)=F_1(\g^{\lambda_1})^*=\begin{cases} \R^{*+}\ \text{pour le type III r\éel},\\
\F^*\ \text{ dans tous les autres cas.}\end{cases}$\\

\noi 2)  $\forall f\in \EuScript S(\goth g_1)$ et $\pi=(\pi_1,...,\pi_n)\in (\widehat{\F^*})^n$ on a:
$$Z^*(\four (f);\pi)=a( \pi)Z(f;\pi^* |\ |^{-N1_n}) \quad
\text{avec}\quad  a( \pi)=\gamma(f)^{\frac{n(n-1)}{2}}\prod_{1\≤j\≤n}a^{(1)}(\pi_j...\pi_n|\ |^{\frac{d}{2d_1}(n-j) +1})$$ 
$$\text {et}\quad  
 a^{(1)}(\pi_1)=\begin{cases} \rho'(\pi_1)\ \text{pour  i) et ii)},\\
 \\
-2(2\pi)^{-2s-\frac{d}{8}-2}\Gamma (s+1)\Gamma(  s+\frac{d}{8}+1)\sin\pi s\ \text{pour le type III r\éel}\ \text{avec }\pi_1=|\ |^s,\\
\\
(-1)^{d_1-1}\prod_{0\≤j\≤d_1-1}\rho'(\pi_1 |\ |^{j} )\ \text{pour le type III}\ (A_n)\ \goth p-\text{adique}. 
\end{cases}$$
\end{theo}

\dem

1) C'est la proposition 6.1.8 dans les cas III r\éel, alors $F_i(\goth g_1)\≥0$ pour $ i=1,...,n,$ et III ($A_n$) $\goth p-$adique   pour lequel $F_i(\goth g_1)=\F$ pour $ i=1,...,n,$ puisque la norme d'une alg\èbre \à division est surjective sur $\F.$\\

\noi Dans les autres cas, les \PV  sont presque d\éploy\és et on a $f(E)^*=\F^*=\chi_i(G_{\goth t_0})$ pour $i=1,...,n.$

\noi En effet,   \gog apparait dans le tableau 3 de \cite{mullerNAG}, \à l'exception du cas  o\ù\gog est de type $(D_n,\alpha_1)$ avec  $\g$   d\éploy\ée (type II avec $e=0$ et $n=2$). 

\noi  Lorsque  \gog est de type $(D_n,\alpha_1)$ avec  $\g$   d\éploy\ée, soit $t\in \F^*$ et $\phi$ l'homomorphisme   de $\∆$ dans $\F^*$ d\éfini par $\phi( \alpha_1)=\phi(\alpha_n)^{-1}=t,$
$\phi( \alpha_i)=1$ pour $2\≤i\≤n-1,$ il se prolonge en l'\elt $g_{\phi}\in G_{\goth t_0}$ d\éfini par $g_{\phi}/\goth g^\mu=\phi(\mu)Id/\goth g^\mu$.

\noi Dans les autres cas il existe un  \PV  commutatif: $(\widetilde{\g_0},\widetilde{\g_1},\widetilde{ H_0},\widetilde{\goth t_0}=\oplus_{1\≤i\≤n+1}\F h_{\widetilde{ \lambda_i}}),$ 
tel que $\g=\goth U(\F h_{\widetilde\lambda_{n+1}})$ et $\ds{\widetilde{ H_0}=H_0+\frac{ h_{\widetilde \lambda_{n+1}}}{2}}$ (cf.prop 4.2.1 de \cite{mullerNAG}) et on applique le 2.du lemme 6.1.7 ainsi que le lemme 6.1.3,iv.\\

2) Par r\écurrence sur $n,$  le cas $n=1$ \étant connu puisqu'il est soit donn\é par l'\équation fonctionnelle de Tate dans les cas I et II (avec $\widetilde B(X_{\lambda_1},X_{-\lambda_{1}})=1$), dans le cas III r\éel c'est  l'\équation fonctionnelle v\érifi\ée par une forme quadratique d\éfinie positive (avec $\ds{\widetilde B(\frac{h_{\lambda_1}}{2},\frac{h_{\lambda_{1}}}{2})=-1}$) et dans le cas III $(A_n)$ $\goth p-$adique, elle est d\ùe \à Jaquet-Godement \cite{godementjacquet}; on applique ensuite la proposition 5.3.1 avec $\mathbb{H}=\F^*$ et $\gamma_k=\gamma(f)^{k(n-k)}$ (lemme 6.1.10).\fdem

\bigskip

\noi {\bf Remarque}: $\gamma (f)=\pm 1,$ $\gamma (f)=\begin{cases}-1 \text{ dans le cas i},\\
1 \text{ dans les cas ii) et iii)}\end{cases}$ (\demo du lemme 6.1.10), on n'a pas cherch\é sa valeur dans le cas iv).

\bigskip
\subsubsection{ Le cas non transitif  }
\bigskip

\noi Pour $u= (u_1,...,u_n)\in (\F^*/\F^{*2})^n$ on note:
$$p(u)_j:=\begin{cases}1\text{ pour }j=0\\
 u_1...u_j \text{ pour }j=1,...,n\end{cases}\ ,\ \ p(u):=p(u)_n\ \text{ et }h(u )=\prod_{1\≤i<j\≤n}(u_i,u_j).$$
 
\begin{theo}
Pour les types:\\

i) I r\éel et  I $\goth p-adique$ avec $d\≤3$ (alors $d_1=1$), \\

ii)  II avec $e>0$ (alors $d_1=1$ et $n=2$),\\

iii) III $(C_{2n},\alpha_{2n})$ $\goth p-$adique de caract\éristique r\ésiduelle diff\érente de $2$ (alors $d=2d_1=4$),\\
 
 \noi On a $\forall f\in \EuScript S(\goth g_1)$ et $\pi=(\pi_1,...,\pi_n)\in (\widehat{\F^*})^n$ ($\widetilde \pi=(\pi_1\pi_2^{-1},...,\pi_n)$): \\  
 
 \begin{enumerate}

\item $(v_1,...,v_n)\in   (F_1(\g^{\lambda_1})^*/\F^{*2})^n$:
$$Z_v^*(\four (f);\widetilde\pi)= \sum_{u=(u_1,...,u_n)\in (F_1(\g^{\lambda_1})^*/\F^{*2})^n} \widetilde a^{(n)}_{v,u}( \pi)\ Z_u(f;(\widetilde\pi)^* |\ |^{-N1_n}) \quad \text{avec}:$$

$$\begin{array}{lll}\widetilde a^{(n)}_{v,u}( \pi)&=&\Lambda_n(v,u)\prod_{1\≤i\≤n} C(\pi_{n-i+1}|\ |^{\frac{d}{2d_1 }(i-1)+1};v_i,u_i ) \text{ avec }\Lambda_1(v,u)=1\text{ et pour }n\≥2\\
\\
\Lambda_n(v,u)&=&\begin{cases}  \prod_{1\≤i<j\≤n}\gamma (u_iv_jf) \text{ dans les cas i) et ii),}\\

 \\
   (-1)^{\frac{n(n-1)}{2}}\ \gamma (F_1)^n\ \bigl(\prod_{1\≤i<j\≤n} (\delta_1,u_iv_j).(u_i,v_j)\ \bigr)\text{ dans le cas iii),}\end{cases} \\
\\
C(\pi;v,u)&=&\begin{cases} \pi(-1)\rho(\pi ;-uv) \text{ dans les cas i) et ii),}\\
\\
 A^1_{  \pi,  \pi|\ |^{\frac{1}{2}}}(v,u,\delta_1)\text{ dans le cas iii)}.\end{cases} \end{array}$$
 
\item Dans les cas i) et ii) avec  $d$ est pair:
$$Z^*(\four (f);\pi)=B(\pi)Z(f;(\pi_{n-1},...,\pi_1,\bigl(\prod_{1\≤i\≤n}\pi_i\bigr)^{-1}.\tilde\omega_{\delta}^{n-1}. |\ |^{-N }) \quad \text{avec}:$$
$$B(\pi)= \gamma (f)^{\frac{n(n-1)}{2}}\
\prod_{1\≤j\≤n}\rho'(\pi_{n-j+1}...\pi_n.\tilde\omega_{\delta}^{j-1}|\ |^{ \frac{d}{2 }(j-1)+1}).$$
\item Dans les cas i) avec  $d=1,$ on a:
$$Z^*(\four (f);\widetilde \pi)=C_n\prod_{1\≤i\≤n}\pi_i(-1)\sum_{u\in  ( \F^*/\F^{*2})^n} \epsilon(u)\ .{{\widetilde b}_u}^{(n)}(\pi)\ Z_u(f;(\widetilde\pi)^* |\ |^{-\frac{n+1}{2})1_n})\quad \text{avec}:$$
$$  {{\widetilde b}_u}^{(n)}(\pi)=\rho (\pi_n|\ |) \prod_{1\≤j\≤[\frac{n-1}{2}]} \rho(\pi_{n- 2j}|\ |^{j+1}\tilde\omega_{(-1)^j  p(  u)_{2j}})
\prod_{1\≤j\≤[\frac{n}{2}]}h(\pi_{n-2j+1}|\ |^{j+\frac{1}{2}}\tilde\omega_{ (-1)^jp(  u)_{2j}} ) $$
et : \\
 
$\bullet$ lorsque $\F$ est un corps $\goth p-$ adique,  soit $p_n=
\begin{cases} 1\text{ si }n\text{ est pair }, \\
0\text{ sinon ,}\end{cases}
$  

\noi alors:
$$ C_n= (- 1,-1)^{[\frac{n}{2}]}\alpha(1)^{[\frac{n}{2}]+\frac{ (n-1)n}{2}-(n-1)^2}\ ,\ 
\epsilon(u) =  (p(u),-1)^{[\frac{n}{2}]}h(u)^n \alpha(p(u))^{p_n}    $$
$\bullet$ lorsque $\F=\R,$ soit $q(u)$ le nombre de composantes n\égatives de $u,$ alors:
$$ C_n=e^{-i\frac{\pi}{4}(\frac{(n-1)n}{2}+[\frac{n}{2}])}\ ,\ 
\epsilon(u)=i^{(n-1)q(u)}(-1 )^{[\frac{q(u)}{2}]}.  $$

\item Dans le cas iii) avec $\pi_i=\tilde\omega_{a_i}|\ |^{s_i},i=1,...,n$ on a:
$$Z^*(\four (f);\widetilde \pi)=D_n\sum_{u\in  (F_1(\g^{\lambda_1})^*/\F^{*2})^n} (\delta_1,p(u))^{n}\ .\alpha(p(u))\ .{{\widetilde c}_u}^{(n)}(\pi)\ Z_u(f;(\widetilde\pi)^* |\ |^{-(n+\frac{1}{2})1_n})\quad \text{avec}:$$
$$\begin{array}{lll}D_n&=& \gamma (F_1)^n(-1)^{\frac{n(n-1)}{2}} \\
\\
 {{\widetilde c}_u}^{(n)}(\pi)&=& \prod_{1\≤i\≤n} \rho(\pi_{n-i+1}|\ |^i\tilde\omega_{  p(\delta_1 u)_{i-1}})h(\pi_{n-i+1}|\ |^{i+\frac{1}{2}}\tilde\omega_{ p(\delta_1 u)_i} )\ \text{et}\\
 \\
 p(\delta_1 u)_i&=&\begin{cases} 1\ \text{si }\ i=0\\
{\delta_1}^iu_1...u_i \ \text{si }\ i\≥1.\end{cases}\end{array}$$

\item Dans les cas i) et ii) avec $n=2$ et d impair, les orbites de $P_0$ dans $\g"_1$ (resp.$\g"_{-1}$) sont en bijection avec $\F^*/(\F^*)^2$ et donn\ées par:
$$w\in \F^*/(\F^*)^2\ :\ O'_w=\cup_{u\in \F^*/(\F^*)^2}O_{(u,wu)}\ (\text{resp. }O^{*'}_w=\cup_{u\in \F^*/(\F^*)^2}O^*_{(u,wu)})$$ 
et $\forall f\in \EuScript S(\goth g_1)$ et $\pi=(\pi_1, \pi_2)\in \widehat{\F^*}^2$ ($\widetilde \pi=(\pi_1\pi_2^{-1}, \pi_2)$): \\ 
$$Z_{O^{*'}_v}^*(\four (f);\widetilde\pi)= \sum_{u \in ( \F^*/\F^{*2})} \widetilde a' _{v,u}( \pi_1,\pi_2)\ Z_{O'_u}(f;(\widetilde\pi)^* |\ |^{-N1_n}) \quad \text{avec}:$$
$$\widetilde a' _{v,u}( \pi_1,\pi_2)= \alpha (-1) \gamma (F_2) (\pi_1.\pi_2)(-1)A^1_{ \pi_2|\ |,\pi_1|\ |^{\frac{d}{2}+1}}(v,u,\delta_2)\ ,\ \delta_2=(-1)^{[\frac{d}{2}+1]}\text{disc}(F_2).$$
\end{enumerate}
 \end{theo}
 
 \dem  1) On proc\ède par r\écurrence sur $n$ en appliquant la proposition 5.3.1 avec  $\mathbb{H}=\F^{*2}$ et $\gamma_k$  donn\é dans le lemme 6.1.10; le cas $n=1$ est rappel\é dans le corollaire 3.6.3 et le th\éor\ème 3.6.5.
 
 \noi 2) Lorsque $d$ est pair,  le r\ésultat  provient du calcul de $\sum_{ v\in (\F^*/\F^{*2})^n}Z_v^*(\four (f);\pi)$  en appliquant 1) ainsi que  la relation $\gamma (tf)= (\delta,t)\gamma (f).$
 
 \noi 3) On a pour $j\≥2:$
 $$ \prod_{1\≤i\≤j-1}\alpha (u_iv_j)=\gamma(v_j(\oplus_{ 1\≤i\≤j-1}u_iX_i^2))= (z_{j-1},v_j)(\alpha(v_j)\alpha(-1))^{p_j}\gamma_{j-1}(u)$$
avec:
$$\gamma_{j}(u)=\gamma( \oplus_{ 1\≤i\≤j}u_iX_i^2)=\prod_{ 1\≤i\≤j}\alpha(u_i)\ ,
  \ z_j=\begin{cases}1\text{ si }j=0\\
(-1)^{[\frac{j}{2}]}u_1...u_j\text{ pour }j=1,...,n\ ,\ \end{cases}  ,$$
$$ p_j=\begin{cases}1\text{ si }j\text{ est pair,}\\
0\text{ si }j\text{ est impair}\end{cases}$$
 donc $$\widetilde a^{(n)}_{v,u}( \pi)=\prod_{1\≤j\≤n}\pi_j(-1)\prod_{1\≤j\≤n-1}\gamma_{j}(u)\prod_{1\≤j\≤n}C_j(u,v_j)$$
avec 
 $$ C_j(u,v_j)=(z_{j-1},v_j)(\alpha(v_j)\alpha(-1))^{p_j}\rho(\pi_{n-j+1}|\ |^{\frac{j+1}{2}};-u_jv_j)$$
 or
 $$\sum_{v_j\in \F^*/\F^{*2} }C_j(u,v_j)= (z_{j-1},-u_j)(\alpha(-u_j)\alpha(-1)^2)^{p_j}.\begin{cases}\rho(\pi_{n-j+1}|\ |^{\frac{j+1}{2}}\tilde\omega_{z_{j-1}})\text{ si }j\text{ est impair,}\\
 h(\pi_{n-j+1}|\ |^{\frac{j+1}{2}}\tilde\omega_{z_{j}})\text{ si }j\text{ est  pair,}
\end{cases}$$
(cf.3.6.2 et 3.6.4)
d'o\ù  $$\begin{array}{lll}\sum_{v\in (\F^*/\F^{*2})^n}\widetilde a^{(n)}_{v,u}( \pi)&=&\prod_{1\≤j\≤n}\pi_j(-1)\prod_{1\≤j\≤n-1}\gamma_{j}(u)\prod_{1\≤j\≤n}
(\sum_{v_j\in \F^*/\F^{*2} }C_j(u,v_j))\\
\\
&=&\prod_{1\≤j\≤n}\pi_j(-1)\widetilde b_u^{(n)}(\pi)\prod_{2\≤j\≤n}\gamma_{j-1}(u)(z_{j-1},-u_j)
(\alpha(-u_j)\alpha(-1)^2)^{p_j}\\
\\
&=&\prod_{1\≤j\≤n}\pi_j(-1)\widetilde b_u^{(n)}(\pi)\prod_{2\≤j\≤n}\gamma(-u_j(\oplus_{1\≤i\≤j-1}u_iX_i^2)\alpha(-1)^{p_j}.
\end{array}$$
 Or:$$\gamma(u_j(\oplus_{1\≤i\≤j-1}u_iX_i^2)=\prod_{1\≤i\≤j-1}\alpha(u_iu_j)=\alpha(-1)^{j-1}\alpha(u_j)^{j-1}(u_j,p(u)_{j-1})\prod_{1\≤i\≤j-1}\alpha(u_i) ,$$
 
 donc
 $$\sum_{v\in (\F^*/\F^{*2})^n}\widetilde a^{(n)}_{v,u}( \pi)= \widetilde b_u^{(n)}(\pi)\alpha(1)^{\frac{n(n-1)}{2}-[\frac{n}{2}]}h(u)\overline{\gamma_n(u)} ^{n-1}\prod_{1\≤j\≤n}\pi_j(-1),$$
Il ne reste plus qu'\à exprimer $\gamma_n(u)$ suivant le corps $\F$ (cf.3.6.2).\\

   \noi 4) On a $F_1(\g^{\lambda_1})^*=\F^*-\delta_1\F^{*2}$ et 
 $$\alpha (\delta_1)A^1_{  \tilde\omega_a|\ |^s,\ \tilde\omega_a|\ |^{s+\frac{1}{2}}}(\delta_1,u,\delta_1)= (a,u\delta_1) A^1_{s,s+\frac{1}{2}}(1,u\delta_1,1)=0\ \text{ pour }\ u\delta_1\not=1$$
(lemme 3.6.4,B)1 et 3) et  1) du lemme 3.6.8) donc on applique 1) en sommant sur $\{v\in (\F^*/\F^{*2})^n\}$ comme pr\éc\édemment or:
$$\Lambda_n(v,u)=C(u) \prod_{1\≤j\≤n}(v_j,p(\delta_1u)_{j-1})\ ,\text{ avec }\ C(u)=D_n.(\delta_1,\prod_{1\≤i\≤n}u_i^{n-i}),$$
donc
$$\sum_{v\in (\F^*/\F^{*2})^n}\widetilde a^{(n)}_{v,u}( \pi)=C(u)\prod_{1\≤j\≤n}C_j,$$
avec
$$\begin{array}{lll}C_j&=&\sum_{v\in  \F^*/\F^{*2}}(v,p(\delta_1u)_{j-1})A^1_{\pi_{n-j+1}|\ |^{j },\pi_{n-j+1}|\ |^{j+\frac{1}{2}}}(v,u_j,\delta_1)\\
&=&(u_j,\delta_1p(\delta_1u)_{j-1})\alpha(u_j)\rho(\pi_{n-j+1}|\ |^{j }\tilde\omega_{p(\delta_1u)_{j-1}})h(\pi_{n-j+1}|\ |^{j+\frac{1}{2}}\tilde\omega_{p(\delta_1u)_j})
\end{array}$$
(B)5),lemme 3.6.4).

\noi On termine en notant que:
$$\prod_{1\≤i\≤n}\alpha(u_i)=\gamma(\oplus_{1\≤i\≤n}u_iX_i^2)=\alpha(1)^{n-1}\alpha(p(u))h(u)\quad  \cite{rallisschiffmann}.$$

 \noi 5) On a $\chi_2 (P_0)=\F^{*2}$ (prop.3.2,ii)\cite{mullerNAG}) et $\chi_1 (P_0)=\F^*$ d'o\ù les $P_0$-orbites dans $\g"_{\pm 1}.$
 
 \noi  Le reste d\écoule des \égalit\és suivantes:  $\delta_2=-\delta,\gamma (f)=\overline {\gamma (F_2)}$ (ii) du lemme 6.1.3) et $$\widetilde a' _{v,u}( \pi_1,\pi_2)=\sum_{v_1\in \F^*/\F^{*2}}\widetilde a^{(2)}_{(v_1,vv_1),(u_1,uu_1)}( \pi_1,\pi_2).\qquad  \Box $$ 
  
\bigskip
\noi {\bf Remarque}: \\

1) Ce th\éor\ème donne les coefficients des \équations fonctionnelles v\érifi\ées par la fonction Z\éta associ\ée \à l'action de $P_0$ puisque chaque $O_u$ (resp. $O^*_u$) avec $u\in (F_1(\g^{\lambda_1})^*/\F^{*2})^n,$ est une seule $P_0-$ orbite dans le cas iii).\\

\noi  Dans le cas non transitif  avec $n\≥3$ (donn\és dans i), les orbites d\épendent de la parit\é de $n.$

\noi   Pour \éviter cette distinction, on consid\ère le sous-groupe $P'_0=G'_0.N_0$ de $P_0$  avec
$ G'_0:=\{g\in G_{\goth t_0}|\chi_i(g)\in f(E)^*,i=1,...,n\}\ , $ alors, \à l'exception du cas $d=2$ lorsque $\F$ est un corps $\goth P-$adique, chaque $O_u$ (resp. $O^*_u$) avec $u\in ( \F^*/\F^{*2})^n,$ est une seule $P'_0-$orbite.\\

\noi Lorsque $\F$ est un corps $\goth P-$adique et $d=2,$ les $P'_0$-orbites sont donn\ées par :
$$  O'_{\epsilon}=\cup_{u\in (  \epsilon)}O_u\quad (\text{resp.}O^{*'}_{\epsilon}=\cup_{u\in ( \epsilon)}O^*_u)\quad \text{avec}$$
$$\epsilon=(\epsilon_1,...,\epsilon_n)\in\{-1,1\}^n\quad\text{et}\quad(\epsilon)=\{u\in ( \F^*/\F^{*2})^n, (\delta,u_i)=\epsilon_i,i=1,...,n\},$$
et $\forall f\in  \EuScript S(\goth g_1)$ et $\pi=(\pi_1,...,\pi_n)\in (\widehat{\F^*})^n$ ($\widetilde \pi=(\pi_1\pi_2^{-1},...,\pi_n)$) on a:
$$Z^*_{O^{*'}_{\epsilon'}} (\four (f);\widetilde\pi)= \sum_{\epsilon   \in \{-1,1\}^n} \widetilde b^{(n)}_{\epsilon',\epsilon}( \pi)\ Z_{O_{\epsilon}}(f;(\widetilde\pi)^* |\ |^{-N1_n}) \quad \text{avec}:$$
$$\widetilde b^{(n)}_{\epsilon',\epsilon}( \pi)= (\overline{\alpha(\delta)})^{\frac{n(n-1)}{2}} s(\epsilon',\epsilon)\prod_{i=1}^n C_{\delta}(\pi_{n-i+1}|\ |^i; \epsilon'_i\epsilon_i )  \ , $$
$$s(\epsilon',\epsilon)=\prod_{i=1}^n(\epsilon'_i)^{i-1}(\epsilon_i)^{n-i}\in\{-1,1\}\ \text{et}\ C_{\delta}(\pi ;t)= \frac{1}{2}(\rho'(\pi )+t\rho'(\pi \tilde\omega_{\delta})),t\in \C$$
$$\text{lorsque }\quad \epsilon  =(\epsilon_1,...,\epsilon_n)\quad ,\quad \epsilon'  =(\epsilon'_1,...,\epsilon'_n)$$
et la caract\éristique r\ésiduelle de $\F$ est diff\érente de $2.$\\

2) On retrouve \également les coefficients des \équations fonctionnelles v\érifi\ées par la fonction Z\éta associ\ée \à l'action de $G.$ \\

$\bullet$ Dans le cas $\goth p$-adique, lorsque $d=2$ et $f$ est anisotrope, $G$ a une seule orbite dans $\g'_1$ lorsque $n$ est impair et pour $n$ pair , les orbites sont s\épar\ées par les valeurs de $(\delta,F_n(x))$ pour $x\in \g'_1.$\\

\noi Pour $\epsilon=\pm 1$ et $f\in \ES (\g_1)$ (resp. $g\in \ES (\g_{-1}$) posons $Z_\epsilon(f;\pi)=Z(fI_{\{u\in \FF\ |\ (\delta,u)=\epsilon\}}(F_n(x));\pi)$ (resp. $Z^*_\epsilon(f;\pi)=Z*(gI_{\{u\in \FF\ |\ (\delta,u)=\epsilon\}}(F^*_n(x));\pi)$ alors il est imm\édiat que pour $ \epsilon'=\pm 1$ on a:
$$Z^*_{\epsilon'}(\four (f);\pi)=\sum_{ \epsilon=\pm 1}a_{\epsilon',\epsilon} (\pi)Z_\epsilon(f;\pi^{-1}|: |^-N)  \text{
avec  }a_{\epsilon',\epsilon} (\pi)=\frac{\epsilon^{n-1}}{2}(B(\pi)+\epsilon\epsilon'B(\pi\tilde\omega_\delta))$$
et $B(\pi)=\gamma(f)^{\frac{n(n-1)}{2}}\prod_{1\≤j\≤n}\rho'(\pi \tilde\omega_{\delta}^{j-1}|\ |^j)$ (2) th\éor\ème 6.2.2).\\

\noi On peut remarquer que $a_{\epsilon',\epsilon} (\tilde\omega_\delta\pi)=\epsilon .\epsilon' .a_{\epsilon',\epsilon} (\pi).$\\

\noi  Revenant au cas g\én\éral, soient $O$ et $O^*$  deux orbites de $G$ dans $\g'_1$ et $\g'_{-1}$ alors $O$ et $O^*$ contiennent des orbites $O_u$ et $O^*_u$ de $P_0$ dans $\g"_1$ et pour tout $h\in \EuScript S(\g_{-1})$ et $\pi\in \widehat {\F^*}$ on a $Z^*_{O^*}(h;\pi)=\sum_{v|O^*_v\subset O^*}Z^*_v(h;(id,...,\pi)).$

$\bullet$ Une premi\ère solution consiste \à calculer $$\sum_{v|O^*_v\subset O^*}\widetilde a^{(n)}_{v,u}(\pi,...,\pi)$$
et \à v\érifier que cette quantit\é est ind\épendante du choix de $u$ tel que $O_u\subset O,$ ce qui prouve \à nouveau l'existence de l'\équation fonctionnelle.

\noi Ceci s'av\ère difficile \à v\érifier et a \ét\é fait dans le cas r\éel (cf.lemme 5.40 de \cite{boppruben}).\\

$\bullet$ Une seconde solution consiste \à consid\érer  une orbite $O_u$ de $P_0$ dans $\g"_1$ telle que $O_u\subset O$ et \à prendre$f\in \EuScript S(O_u)$ telle que $Z(f;\omega |\ |^s)=Z_u(f;(id,...,\omega|\ |^s))$ ne soit pas identiquement nulle et de comparer les $2$ \équations fonctionnelles,  ce qui donne:
$$\begin{array}{rcl}Z^*_{O^*}(\four f;\pi)&=&a_{O^*,O}(\pi)Z(f;\pi^{-1}|\ |^{-N})\\
&=&\sum_{v|O^*_v\subset O^*}Z^*_v(\four f;(id,...,\pi))\\
&=&  (\sum_{v|O^*_v\subset O^*}\widetilde a^{(n)}_{v,u}(\pi,...,\pi)\ Z(f;\pi^{-1}|\ |^{-N}) \text{ d'o\ù} \\
 a_{O^*,O}(\pi)&=& \sum_{v|O^*_v\subset O^*}\widetilde a^{(n)}_{v,u} (\pi,...,\pi)\end{array}  $$
donc cette somme est ind\épendante du choix de $u$ tel que $O_u\subset O.$\\

\noi On indique uniquement les r\ésultats dans le cas $\goth p$-adique $(C_n\alpha_n)$ de type I (i.e.$\g$ est d\éploy\é) ou III (i.e.$(\g_0,\g_1)$ n'est pas presque d\éploy\é). 
\bigskip

  $\bullet$ Dans le cas  d\éploy\é, donc $d=1,$  avec $n\≥3,$ $n=2\ell$ ou $n=2\ell+1,$ soient:
  $$\begin{array}{ccl}X_1&=&\sum_{1\≤i\≤n-1}X_{\lambda_i}+(-1)^\ell X_{\lambda_n}\ ,\ 
  Y_1= \sum_{1\≤i\≤n-1}X_{-\lambda_i}+(-1)^\ell X_{-\lambda_n}, \\
  \\
   X_{-1}&=&\sum_{1\≤i\≤n-3}X_{\lambda_i}+wX_{\lambda_{n-2}}+ww'X_{\lambda_{n-1}}+ (-1)^\ell w'X_{\lambda_n} ,\\
   \\
  Y_{-1}&=&\sum_{1\≤i\≤n-3}X_{-\lambda_i}+\frac{1}{w}X_{-\lambda_{n-2}}+\frac{1}{ww'}X_{\lambda_{n-1}}+ (-1)^\ell \frac{1}{w'}X_{-\lambda_n}  ,\end{array}$$
  $w$ et $w'$ \étant $2$ \elts fix\és dans $\F^*$ tels que $(w',w(-1)^{\ell+1})=-(w,-1).$
  
\bigskip  

  \noi a) $n=2\ell$
  
\bigskip

\noi Pour $u\in \F^*/\F^{*2},$ soient 
$$O'_u=\{x\in \g_1|F(x) \in(-1)^\ell u\ \F^{*2}\}  \ ,\ O'^*_u=\{x\in \g_{-1}|F^*(x) \in(-1)^\ell u\ \F^{*2}\},$$
$$O'_{1,\epsilon}=G( X_\epsilon)    \ ,\  
 O'^*_{1,\epsilon}=G( Y_\epsilon)    \ ,\  \epsilon=\pm 1,$$
 alors $O'_{1,1}\cup O'_{1,-1}=O'_1$ et $O'^*_{1,1}\cup O'^*_{1,-1}=O'^*_1.$ \\

\noi Les orbites de $G$ dans $\g'_1$ (resp.$\g'_{-1}$) sont donn\ées par: $$O'_{1,1},O'_{1,-1},O_{u_0},u_0\in \F^*- \F^{*2}/ \F^{*2}\text{ (resp. }
O'^*_{1,1},O'^*_{1,-1},O'^*_{u_0},u_0\in \F^*- \F^{*2}/ \F^{*2})$$
 et 
 $$\cup_{u\in (\F^*/\F^{*2})^n| p(u)=(-1)^\ell u_0}O_{ u}\subset O'_{u_0} \text{ (resp. }\cup_{u \in (\F^*/\F^{*2})^n| p(u)=(-1)^\ell u_0}O^*_{ u}\subset O'^*_{u_0})\text{ lorsque }u_0\not=1,$$ 
$$\cup_{u \in (\F^*/\F^{*2})^n| p(u)=(-1)^\ell ,h(u )=\epsilon }O_u \subset O'_{ 1,\epsilon } \text{ (resp. }\cup_{u \in (\F^*/\F^{*2})^n| p(u)=(-1)^\ell ,h(u )=\epsilon}O^*_u \subset O'^*_{ 1,\epsilon }) \text{ avec }\epsilon =\pm 1.$$
 \bigskip

\noi   b) $n=2\ell+1$
  
\bigskip

\noi $G$ a $2$ orbites dans $\g'_1$ (resp.$\g'_{-1}$):
$O'_\epsilon =G(X_\epsilon)$ (resp.$O'^*_\epsilon=G(Y_\epsilon)$) avec $\epsilon=\pm 1$ et$$\cup_{u \in (\F^*/\F^{*2})^n| (p(u),(-1)^\ell )h(u )=\epsilon}O_u \subset O'_{\epsilon } \text{ (resp. }\cup_{u \in (\F^*/\F^{*2})^n| (p(u),(-1)^\ell )h(u )=\epsilon}O^*_u \subset O'^*_{ \epsilon }) \text{ avec }\epsilon =\pm 1.$$

\bigskip

\noi Dans les deux cas (pair et impair),  on a:
 $$\frac{a_{O'^*,O'}(\pi)}{\prod_{1\≤i\≤n-1}\alpha (u_i)^{n-i}}=\pi (-1)^{n}\alpha(-1)^{[\frac{n}{2}]}\sum_{v\in S(O'^*) }  \prod_{1\≤j\≤n}\alpha(v_{j})^{p_j}((-1)^{[\frac{j-1}{2}]} p(u)_{j-1},v_j)\rho (\pi|\ |^\frac{j+1}{2};-v_ju_j)   ,$$
 avec:
  $$S(O'^*_{\epsilon'})=\{ v \in (\F^*/\F^{*2})^n| (p(v),(-1)^\ell )h(v )=\epsilon'\}\ \text{lorsque }n\text{ est impair},$$
 et lorsque $n$ est pair:
  $$\begin{array}{ccl}S(O'^*_{1,\epsilon})&=&\{v \in (\F^*/\F^{*2})^n| p(v)=(-1)^\ell ,h(v )=\epsilon\}, \\
  \\
  S(O'^*_{v_0} )&=& \{v\in (\F^*/\F^{*2})^n| p(v)=(-1)^\ell  v_0\}\ \ v_0\not=1\ \ ,\end{array}$$
 \noi Dans les deux cas (pair et impair), $  O_u\subset O'$ quelconque.\\
 
 \noi Dans le cas impair, on peut choisir $u=(1,...,(-1)^\ell)$ pour $O'_1$ et $u= (1...,w,ww',w'(-1)^\ell)$ pour $O'_{-1}.$
 
  \noi Dans le cas pair, on peut choisir $u= (1...,w,ww',w'(-1)^\ell)$ pour $O'_{1,-1},$
$u=(1,...,u_0(-1)^\ell)$ pour $O'_{u_0}$ et $O'_{1,1}$ lorsque $u_0=1.$  \\
  
\noi On n'a pas cherch\é \à simplifier ces sommes.  \\

\noi  Lorsque la caract\éristique r\ésiduelle est diff\érente de $2,$ on donne  $a_{O'^*_1 ,O' }(s)+a_{O'^*_{ -1},O' }(s)$ dans le cas $n$ impair et  $a_{O'^*_v,O' }(s)$ lorsque  $v\not=1$  ainsi que $a_{O'^*_{ 1,-1},O' }(s)+a_{O'^*_{ 1,1},O' }(s)$ dans le cas $n$ pair.\\  

\noi Dans le cas $n$ impair, on introduit les ouverts de $\g_1$, $O'_{u_0,\epsilon}=O'_{\epsilon}\cap \{x\in \g_1|F(x)\in u_0\F^{*2}\},$  (resp.de $\g_{-1}$, $O'^*_{u_0,\epsilon}=O'^*_{\epsilon}\cap \{x\in \g_{-1}|F^*(x)\in u_0\F^{*2}\}$), $\epsilon=\pm 1$ ainsi que les fonctions Z\éta associ\ées: $Z_{u_0,\epsilon}(f;\ )=Z(f1_{O'_{u_0,\epsilon}};\ ).$

\noi On a $\cup_{\{u\in (\F^*/\F^{*2})^n|p(u)=u_0,h(u)=(u_0,-1)^\ell\epsilon\}} O_u \subset O'_{u_0,\epsilon}.$ \\

$\bullet$ Dans le cas  symplectique commutatif de type non d\éploy\é, i.e. $(\overline \g_0,\overline \g_1)$ est de type  $(C_{n},\alpha_{n})$ avec $n$ pair et $(\g_0,\g_1)$ est de type $(C_\frac{n}{2},\alpha_{\frac{n}{2}})$ avec $n\≥4,$ $\F$ \étant un corps $\goth p-$adique de caract\éristique r\ésiduelle diff\érente de $2.$   \\

\noi Les  orbites de $G$ dans $\g'_1$ (resp.$\g'_{-1}$) sont donn\ées par:
$$O'_{u}=\{x\in \g_1|F(x)\in u\F^{*2} \}\quad (\text{resp.}O'^*_{u}=\{x\in \g_{-1}|F^*(x)\in u\F^{*2} \}), \ u\in \F^*/\F^{*2}$$ et $$\cup_{u'\in (Im(F_1)^*/\F^{*2})^\frac{n}{2}|p(u')=u}O_{u'}\subset O'_{u} \quad (\text{resp.}\cup_{u'\in (Im(F_1)^*/\F^{*2})\frac{n}{2}n|p(u')=u}O^*_{u'}\subset O'^*_{u}).$$

\begin{defi} $\F$ est r\éel ou bien $\goth p-$adique.

\noi Pour $s\in \C$  et pour $n\in \N,$  soit : $$f_1(s)=1\ ,\ f_{n+1}(s)={|2|_\F}^{-2ns-n(n+1)-\frac{n}{2}}{\prod}_{\tiny{ 1\≤j\≤n}}\rho (|\ |^{2s+2j+1})\ ,\ n\≥1.$$

\noi Lorsque $\F$ est un corps $\goth p-$adique de caract\éristique r\ésiduelle diff\érente de $2$, on a:
$$ f_{n+1}(s)=(-1)^{n}q^{n(2s+n+1)} {\prod}_{\tiny{ 1\≤j\≤n}}\ds\frac{1-q^{-2(s+j)}}{1-q^{-2(s+j)-1}}\ ,\ n\≥1.$$
\end{defi}

 \begin{prop}
 
Cas symplectique commutatif  de rang $\≥2,$ sur un corps $\goth p-$adique de caract\éristique r\ésiduelle diff\érente de $2,$  i.e.  $(\overline\g_0,\overline\g_1)$ est de type $(C_n,\alpha_n)$ et $(\g_0,\g_1)$ est soit de type $(C_n,\alpha_n)$ (donc de type I) soit de type $(C_{\frac{n}{2}},\alpha_{\frac{n}{2}})$ (donc de type III  avec $n$ pair).  \\

\noi Soient $a\in \F^*/\F^{*2},s\in \C$ alors
pour $f\in \EuScript S(\g_1)$ on a:
\begin{enumerate}
 \item Lorsque $n$ est pair: 
 
 a)$$\frac{ Z^*(\four (f); \tilde\omega_a,s)}{ f_{\frac{n}{2}}(s)}=K_n\rho(\tilde\omega_a,s+1) 
 \sum_{u\in   \F^*/\F^{*2}}(u,\delta_1)^{\frac{n}{2}} \alpha(u)  h( \tilde\omega_{a \delta_1^{\frac{n}{2}}u}|\ |^{s+\frac{n+1}{2}})\    Z_u(f; \tilde\omega_a, -s -\frac{n+1}{2}) $$
 avec:
 $$\begin{array}{ccl}\delta_1&=&-\text{discriminant de }\begin{cases} F_2\text{ lorsque }\g\text{ est d\éploy\ée (i.e. de type I),}\\
 F_1\text{ lorsque }\g\text{ est  non d\éploy\ée (i.e. de type III)}  ,\end{cases}   \\
 \\
 \delta_1&=&-1 \text{  lorsque }\g\text{ est d\éploy\ée,}\end{array}   $$ $$K_n=\begin{cases} 1\text{ lorsque }\g\text{ est d\éploy\ée (type I),}\\
 \\
 \bigl(\  \gamma (F_1)^\frac{n}{2}(-1)^{\frac{n(n-2)}{8}} \text{ lorsque }\g\text{ est  non d\éploy\ée (de type III)} .\end{cases}
 $$
b) Pour $v\in \F^*/\F^{*2}$ on a:
$$Z_v^*(\four (f);  s)=K_nf_{\frac{n}{2}}(s) 
 \sum_{u\in  \F^*/\F^{*2}} A^1_{ s+1,s+\frac{n+1}{2}}(v,u,{\delta_1}^\frac{n}{2})  Z_u(f; -s-\frac{n+1}{2}).$$

 \item Lorsque $n$ est impair, donc $\g$ est d\éploy\ée (i.e. de type I):$$\frac{Z^*(\four (f); \tilde\omega_a,s)}{ f_{[\frac{n}{2}]}(s)}=\rho(\tilde\omega_a,s+1) (a,-1)    ( \ Z_{O'_1}(f;\tilde\omega_a, -s -\frac{n+1}{2})- Z_{O'_{-1}}(f; \tilde\omega_a, -s -\frac{n+1}{2}) ).$$
Pour $v\in \F^*/\F^{*2}$ on a:
$$Z_v^*(\four (f);  s)= f_{[\frac{n}{2}]}(s) 
 \sum_{u\in  \F^*/\F^{*2}, \epsilon=\pm 1}  \epsilon \rho(s+1;-uv)  Z_{u,\epsilon}(f; -s-\frac{n+1}{2}).$$
 
\end{enumerate}

\end{prop}

 \dem
$\bullet$ Dans le cas d\éploy\é,  on applique le 3) du th\éor\ème 6.2.2 avec $\pi=(\pi_1,...,\pi_1)$ et $\pi_1=\tilde\omega_a|\ |^s$ alors $C_n=1$ et pour $u\in ( \F^*/\F^{*2})^n$ on a $$\epsilon(u)=(p(u),-1)^{[\frac{n}{2}]}h(u)^n\alpha(p(u)^{p_n}$$ et:
$$\begin{array}{lll}{{\widetilde b}_u}^{(n)}(\pi)&=& \rho (\pi_1|\ |) \prod_{1\≤j\≤[\frac{n-1}{2}]} \rho(\pi_1 |\ |^{j+1}\tilde\omega_{ (-1)^j p(  u)_{2j}})
\prod_{1\≤j\≤[\frac{n}{2}]}h(\pi_1 |\ |^{j+\frac{1}{2}}\tilde\omega_{(-1)^j p(  u)_{2j}} )\\
\\
 
&=& \rho(\pi_1|\ |) f_{[\frac{n}{2}]}(s).\begin{cases}1\text{ lorsque }n\text{ est impair,}\\
h(\pi_1|\ |^{\frac{n+1}{2}}\tilde\omega_{  (-1)^\frac{n}{2} p( u)} )\text{ lorsque }n\text{ est  pair,} \end{cases}
\end{array}$$
(lemme 3.6.8).  \\

\noi Lorsque $n$ est pair, $\epsilon(u){{\widetilde b}_u}^{(n)}(\pi)$ ne d\épend que de $p(u),$  donc en sommant sur les $\{u\in (\F^*/\F^{*2})^n|p(u)= u_0\}$ avec $u_0 $ d\écrivant $ \F^*/\F^{*2}$, en obtient le premier r\ésultat.\\

\noi Comme on a $\forall v\in  \F^*/\F^{*2}:$  
$$Z_v^*( h; s)=\ds\frac{1}{|\FF|}\sum_{b\in  \F^*/\F^{*2}}(b,v)Z^*( h;\tilde\omega_b, s)$$ et que 
$\forall a\in  \F^*/\F^{*2}$:
$$Z^*(\four (f);\tilde\omega_a,s)=\sum_{u\in  \F^*/\F^{*2}}B_u(\tilde\omega_a,s)Z_u(f;\tilde\omega_a,-s-\frac{n+1}{2}),$$
 on a: $\forall v\in  \F^*/\F^{*2}$ on a:$$Z_v^*(\four (f); s)=\sum_{u\in  \F^*/\F^{*2}}A_{v,u}( s)Z_u(f; -s-\frac{n+1}{2}),$$
avec 
$$A_{v,u}( s)=\frac{1}{|\F^*/\F^{*2}|}\sum_{b\in \F^*/\F^{*2}}\tilde\omega_b(uv)B_u(\tilde\omega_b,s)$$
d'o\ù le second r\ésultat en appliquant 6,B du lemme 3.6.4..\\

\noi Lorsque $n$ est impair, $\epsilon(u){{\widetilde b}_u}^{(n)}(\pi)$ ne d\épend que de $ (p(u),-1)^{[\frac{n}{2}]}h(u),$  donc en sommant sur les $\{u\in (\F^*/\F^{*2})^n| (p(u),-1)^{[\frac{n}{2}]}h(u)=\epsilon\}$ avec $ \epsilon=\pm 1 ,$  en obtient le premier r\ésultat ce qui implique le second r\ésultat comme ci-dessus.\\

$\bullet$ Dans le cas non d\éploy\é, on applique le 4) du th\éor\ème 6.2.2  en rempla\çant $n$ par $\frac{n}{2},$  avec  $\pi=(\pi_1,...,\pi_1)$ et $\pi_1=\tilde\omega_a|\ |^s$ alors pour $u\in (F_1(\g^{\lambda_1})^*/\F^{*2})^n$ on a:
$$\begin{array}{lll}{{\widetilde c}_u}^{(\frac{n}{2})}(\pi)&=&\prod_{1\≤i\≤\frac{n}{2}} \rho(\pi_1|\ |^i\tilde\omega_{  p(\delta_1 u)_{i-1}})h(\pi_1|\ |^{i+\frac{1}{2}}\tilde\omega_{ p(\delta_1 u)_i} )\\
\\
&=& \rho(\pi_1|\ |)h(\pi_1|\ |^{\frac{n+1}{2}}\tilde\omega_{ \delta_1^\frac{n}{2}p( u)} )
\prod_{1\≤i\≤\frac{n}{2}-1}\rho(\pi_1|\ |^{i+1}\tilde\omega_{  p(\delta_1 u)_{i}})h(\pi_1|\ |^{i+\frac{1}{2}}\tilde\omega_{ p(\delta_1 u)_i} )\\
\\
&=& \rho(\pi_1|\ |)h(\pi_1|\ |^{\frac{n+1}{2}}\tilde\omega_{ \delta_1^\frac{n}{2}p( u)} )
f_{\frac{n}{2}}(s),\end{array}$$
(lemme 3.6.8)  ainsi ${{\widetilde c}_u}^{(\frac{n}{2})}(\pi)$ ne d\épend que de $p(u)$ d'o\ù le r\ésultat et on termine comme dans le cas d\éploy\é pair. \fdem\\

 \bigskip

 \bigskip
\begin{rema} Les coefficients obtenus dans le cas $n$ pair (1) de la proposition 6.2.4)  sont \à comparer \à ceux associ\és \à une forme quadratique d\éfinie sur un espace vectoriel de dimension
$n+1$ et de discriminant  $(-\delta_1)^\frac{n}{2}.$ \\

\noi Posons $n=2\ell,$ les valeurs explicites:
 $$\frac{K_{2\ell}\ }{2}\alpha (-\delta_1)^\ell(-1)^{\ell-1}q^{(\ell-1)(2s+\ell)}g_\ell(s)P^{1,\epsilon_1,
\epsilon_2,1}_{s+1,s+l+\frac{1}{2}}( \delta_1^\ell v,  \delta_1^\ell u,1),$$
 avec $g_1(s)=\ds\frac{1}{(1-q^{-2(s+1) })(1-q^{-2(s+ 3)})}$ et pour $\ell\≥2:$
$$g_\ell(s)=\frac{q_\ell(s)}{r_\ell(s)} \ ,\ r_\ell(s)=\prod_{j=1}^\ell(1-q^{-2(s+j)-1})\ ,\ q_\ell(s)=\begin{cases}1\text{ si }\ell=2\\
\prod_{j=2}^{\ell-1}(1-q^{-2(s+j)})\text{ si }\ell\≥3\end{cases}$$
sont donn\ées dans le lemme 3.6.7 (B).\\

\noi On rappelle que les racines du polynome de Bernstein associ\é au \PV de type $(C_{2\ell},\alpha_{\ell})$ sont donn\ées par: $\{-\frac{j}{2},j=0,...,2\ell-1\}$ donc les p\ôles de $g_\ell,$ qui sont simples et donn\és par $\{-(j+\frac{1}{2})+ik\frac{\pi}{\ln q},1\≤j\≤\ell,k\in \Z\},$ pour $\ell\≥2,$ ne font intervenir que les racines non enti\ères du polynome de Bernstein.
\end{rema}

\newpage
\section {  Les cas classiques}
\vskip 3mm
 
 Dans ce paragraphe, on d\étermine les polynômes de Bernstein ainsi que les coefficients de l'\eq associée au \PV $(P_0,\g_1),$ $P_0$ étant le \sg parabolique standard très spécial défini dans le $\S 2.4,$ (cf.$\S 7.1.2$) par la méthode de descente.\\

 Les résultats obtenus sont complets dans le cas symplectique, presque complets  dans les cas orthogonaux BI,DI ($\S 7.3$) et DIII réel (1) du th.7.4.4). Le cas DIII $\goth p-$adique est incomplet puisqu'on suppose certaines conditions restrictives (2) du th.7.4.4) et que l'on est incapable d'en déduire les résultats pour $(G,\g_1).$\\
 
 \noi {\bf Notation}: $\phi$ \étant une forme quadratique non d\ég\én\ér\ée, on note simplement $\gamma (\phi)$ la constante $\gamma (\tau\circ\phi)$ (cf.$\S 3.6.2$).\\

\subsection{  G\én\éralit\és}

\subsubsection{Description des cas classiques consid\ér\és}
\vskip 3mm
\noi Les cas classiques correspondent aux \PVs obtenus \à partir
d'une alg\èbre simple $\overline \g$ de type $B_n,C_n$
ou $D_n$ (not\ée $R_n$) munie d'une graduation de type
$(R_n,\alpha_k)$ avec:\\

\qquad   $\bullet$ $1\≤3k\≤ 2n-2$   et $k$    pair dans le cas $C_n,$
\vskip 2mm
\qquad   $\bullet$ $1\≤3k\≤ 2n-2$ dans le cas  $D_n,$
\vskip 2mm

\qquad    $\bullet$ $1\≤3k\≤ 2n-1$  dans le cas  $B_n,$ \\

\noi ces conditions plus restrictives que la condition pour que
$2H_0$  soit $1$- simple (\cite{rublivre},Table I p.137) sont dûes au choix du \sg parabolique standard très spécial. \\

\noi Dans les notations des planches de \cite{bourbakigal6}, le syst\ème de racines
associ\é \à $\overline \g$, not\é  
$\overline
\Delta ,$ est donn\é par: 
\vskip 2mm
\qquad $\bullet$ $R_n:=B_n\ \ : \ \ {\overline
\Delta}=\{{\overline \epsilon}_i\pm{\overline
\epsilon}_j\ ,\ 1\≤i\not =j\≤n\ ,\ {\overline \epsilon}_i\ ,\
 1\≤i\≤n\},$
\vskip 2mm
\qquad $\bullet$ $R_n:=C_n\ \ : \ \ {\overline
\Delta}=\{{\overline \epsilon}_i\pm{\overline
\epsilon}_j\ ,\ 1\≤i\not =j\≤n\ ,\ 2{\overline \epsilon}_i\ ,\
  1\≤i\≤n\},$
\vskip 2mm
\qquad $\bullet$ $R_n:=D_n:\ \ \ \ {\overline
\Delta}=\{{\overline \epsilon}_i\pm{\overline
\epsilon}_j\ ,\ 1\≤i\not =j\≤n\},$
 
 \noi et:
 
 \qquad $\bullet$  $ {\overline
\Delta}_1=\{{\overline \epsilon}_i\pm{\overline
\epsilon}_j\ , \ {\overline \epsilon}_i\ ,\
   \ 1\≤i \≤k<j\≤n\ \}\cap {\overline
\Delta},$
 \vskip 2mm
\qquad $\bullet$  $ {\overline
\Delta}_2=\{{\overline \epsilon}_i+{\overline
\epsilon}_j\ ,  
   \ 1\≤i \≤j\≤k\ \}\cap {\overline
\Delta}.$
\vskip 3mm
\noi $\goth g$ \étant une $\F$-forme de $\overline{\goth g}$
convenable (i.e. $H_0\in \goth g$), $(\goth g_0,\goth g_1)$
est \également de type classique : $(R_m,\lambda_p)$ avec
$R=B,C,BC$ ou $D.$ Plus pr\écis\ément:
\vskip 2mm

\begin{enumerate}

\item  $(\g_0,\g_1)$ est de type BI, ou BI$(n,k)$ ou $B_m(n,k)$ lorsque le diagramme de Satake de $\g$ est de la forme:\\

\hskip 10pt \hbox to 5cm {\offinterlineskip \lower 2pt\hbox{$\circ$} \hglue -3,2pt
{\vrule height  0,4pt depth 0pt width 0,5 cm}\lower 2pt\hbox{$\circ$}
\hglue -3,2pt
{\vrule height  0,4pt depth 0pt width 0,5 cm}\lower 2pt\hbox{$\circ$}
  \dotfill \hbox to 0,4 cm  
 { \offinterlineskip \lower 2pt\hbox{$\circ$} \hglue -3,2pt
{\vrule height  0,4pt depth 0pt width 0,5 cm}}\hglue 7pt{\lower 2pt\hbox{$\bullet$}}  
  \dotfill \hbox to 0,4 cm 
    { \offinterlineskip \lower 2pt \hbox{$\bullet$}\hglue -3,2pt
{\vrule height  0,4pt depth 0pt width 0,5 cm}\kern -1pt   \hrulefill     \kern -3,4pt\lower 2pt \hbox{ 
\offinterlineskip  {$\bullet$} \hglue -7pt\vbox{{\hrule height 0,3pt width 0,6cm}\vskip 3,5pt{\hrule height 0,3pt width 0,6cm}\vskip 0,3pt }\hglue -12pt $>$   \hglue -1,8pt $\bullet$}}}  \hskip 3cm (diagrammes BI de \cite{warner})

$\hskip 0,3cm \alpha_1.....\hskip 1,3cm.....\alpha_m.....\hskip 1,5cm.....\alpha_n$\\

\noi  i.e. $(\overline\g_0,\overline\g_)$ est de type $(B_n, \alpha_k)$  et $(\g_0,\g_1)$ est de type $(B_m, \lambda_k)$ avec $k\≤m\≤n,$ $\g$ est de rang $m$ et $\g$ est d\éploy\ée lorsque $m=n.$ 

\noi Lorsque $\overline R=B,$ on a toujours $(\g_0,\g_1)$ de type $BI(n,k) $ par classification. \\

\item  $(\g_0,\g_1)$ est de type DI, ou DI$(n,k)$ ou $D_m(n,k)$ lorsque le diagramme de Satake de $\g$ est de la forme:\\

\hskip 10pt \hbox to 6cm {\offinterlineskip \lower 2pt\hbox{$\circ$} \hglue -3,2pt
{\vrule height  0,4pt depth 0pt width 0,5 cm}\lower 2pt\hbox{$\circ$}
\hglue -3,2pt
{\vrule height  0,4pt depth 0pt width 0,5 cm}\lower 2pt\hbox{$\circ$}
  \dotfill \hbox to 0,4 cm  
 { \offinterlineskip \lower 2pt\hbox{$\circ$} \hglue -3,2pt
{\vrule height  0,4pt depth 0pt width 0,5 cm}}\hglue 7pt{\lower 2pt\hbox{$\bullet$}}  
 \dotfill \hbox to 2 cm {\lower 2pt
\hbox{$\bullet$}\hrulefill \lower 2pt \vtop {\offinterlineskip \hbox{$ \bullet$} \hbox to 5pt{\hfill \vrule height 12pt width 0,3pt\hfill} \hbox{$ \bullet$}}\hrulefill\lower 2pt\hbox{$\bullet$} } }   

$\hskip 0,3cm \alpha_1.....\hskip 0,9cm.....\alpha_m.....\hskip 0,8cm.....\alpha_n$\\

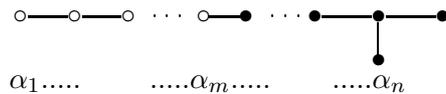
\noi ou bien toutes les racines du diagramme de Satake sont blanches (fl\éch\ées ou non) (diagrammes DI de \cite{warner}).\\

\noi  i.e. $(\overline\g_0,\overline\g_)$ est de type $(D_n,\alpha_k)$ et $(\g_0,\g_1)$ est de type $(D_n,\alpha_k)$ lorsque $\g$ est d\éploy\ée et sinon de type $(B_m,\alpha_k)$ avec $k\≤m<n,$ $m$ \étant le rang de $\g.$  \\

\noi Les cas BI ou DI sont appel\és "type I".\\

\item  $(\g_0,\g_1)$ est de type CII lorsque $\g$ n'est pas d\éploy\ée donc  le diagramme de Satake de $\g$ est de la forme:\\

\hskip 10pt \hbox to 5cm {\offinterlineskip \lower 2pt\hbox{$\bullet$} \hglue -3,2pt
{\vrule height  0,4pt depth 0pt width 0,5 cm}\lower 2pt\hbox{$\circ$}
\hglue -3,2pt
{\vrule height  0,4pt depth 0pt width 0,5 cm}\lower 2pt\hbox{$\bullet$}
  \dotfill \hbox to 0,4 cm  
  {\offinterlineskip\lower 2pt\hbox{$\circledcirc$}}  
  \dotfill \hbox to 0,4 cm 
    { \offinterlineskip \lower 2pt \hbox{$\circ$}\hglue -3,2pt
{\vrule height  0,4pt depth 0pt width 0,5 cm}\kern -1pt   \hrulefill     \kern -3,4pt\lower 2pt \hbox{ 
\offinterlineskip  {$\bullet$} \hglue -7pt\vbox{{\hrule height 0,3pt width 0,6cm}\vskip 3,5pt{\hrule height 0,3pt width 0,6cm}\vskip 0,3pt }\hglue -12pt $<$   \hglue -1,8pt $\circ$}}}\\

\noi ainsi que toutes les variantes ((diagrammes CII de \cite{warner}).\\

\noi  i.e. $(\overline\g_0,\overline\g_1)$ est de type $(C_n,\alpha_{2p}),$ i.e. $k=2p,$  et $(\g_0,\g_1)$ est de type $(C_m,\alpha_p)$  ou $(BC_m,\alpha_p).$\\

\item  $(\g_0,\g_1)$ est de type DIII ou DIII(n,p) lorsque le diagramme de Satake de $\g$ est de la forme:\\

\hskip 10pt \hbox to 5 cm {\offinterlineskip \lower 2pt\hbox{$\bullet$}
\hglue -3pt{\vrule height .4pt depth 0pt width 0,5cm}\lower 2pt\hbox{$\circ$}
\hglue -7pt{ \vrule height .4pt depth 0pt width 0,5cm}\lower 2pt\hbox{$\bullet$}
 \dotfill \hbox to 2 cm {\lower 2pt\hbox{$\circ$}\hglue -1,3pt
\hrulefill\lower 2pt
\hbox{$\bullet$}\hrulefill \lower 2pt \vtop {\offinterlineskip \hbox{$\circ$} \hbox to 5pt{\hfill \vrule height 12pt width 0,3pt\hfill} \hbox{$\circ$}}\hrulefill\lower 2pt\hbox{$\bullet$} } }\\

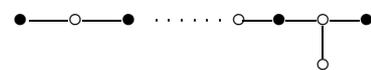
\noi ainsi que toutes les variantes (diagrammes DIII de \cite{warner}).\\

\noi  i.e. $(\overline\g_0,\overline\g_1)$ est de type $(D_n,\alpha_{2p})$ et $(\g_0,\g_1)$ est de type $(C_m,\alpha_p)$ (lorsque $n$ est pair) ou $(BC_m,\alpha_p)$ (lorsque $n$ est impair).\\

\noi  Dans le cas
$\goth p$-adique, on ne consid\èrera pas l'unique forme de $D_{2q+1}$  de diagramme de Satake:\\

\hskip 10pt \hbox to 5 cm {\offinterlineskip \lower 2pt\hbox{$\bullet$}
\hglue -3pt{\vrule height .4pt depth 0pt width 0,5cm}\lower 2pt\hbox{$\circ$}
\hglue -7pt{ \vrule height .4pt depth 0pt width 0,5cm}\lower 2pt\hbox{$\bullet$}
 \dotfill \hbox to 2 cm {\lower 2pt\hbox{$\bullet$}\hglue -1,3pt
\hrulefill\lower 2pt
\hbox{$\circ$}\hrulefill \lower 2pt \vtop {\offinterlineskip \hbox{$\bullet$} \hbox to 5pt{\hfill \vrule height 12pt width 0,3pt\hfill} \hbox{$\bullet$}}\hrulefill\lower 2pt\hbox{$\bullet$} } } \hskip 4cm (\cite{veisfeiler})\\

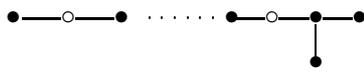
\noi et pour laquelle $rg(\goth g)={n-3\over 2}=q-1,$ ainsi on aura toujours $m=[\frac{n}{2}].$\\

\noi On supposera également que $n-3p\≥4$ dans le cas réel et que  $n-3p\≥2$ avec $n-p$ impair dans le cas $\goth p$-adique.
 
\end{enumerate}

 \vskip 2mm

\noi On note le syst\ème de racines associ\é \à $\g$ de la mani\ère suivante:
 \vskip 2mm
$  \Delta=\{ \epsilon_i\pm
\epsilon_j\ ,\ 1\≤i\not =j\≤m\ ,\ \epsilon_i\ ,\ 2\epsilon_i\
,\ 1\≤i\≤m\}\cap \Delta\supseteq \{ \epsilon_i\pm 
\epsilon_j\ ,\ 1\≤i\not =j\≤m\}.$ \\

\noi  Dans le cas orthogonal, i.e.$(\overline R_n,\alpha_k)$ avec $\overline R=B$ ou $D,$ on
a $\goth g_2=\{0\}$ lorsque $k=1$ (cas commutatif),  
$\goth g_2$ est de dimension $1$ lorsque $k=2$ et pour $k\≥3$
 les diff\érentes formes de $\goth g$ se
s\éparent simplement puisque:
\vskip 3mm
 \begin{lem}
 Soit $(\overline{\goth
g}_0,\overline{\goth
g}_1)$ un \PV  de type $(R_n,\alpha_k)$ avec 
\vskip 2mm
 \qquad i) Pour R=C: $2\≤k\≤n-1,$ 

\qquad ii) Pour R=D: $3\≤k\≤n-2,$\\

et soit $\tilde {\goth g}$ l'alg\èbre engendr\ée par $\goth
g_{\pm2},$ alors:
\begin{enumerate}
 \item $\goth g$ est de type I $\Leftrightarrow\ rg(\tilde
{\goth g})=k.$

  \item $\goth g$ est de type CII $\Leftrightarrow\ rg(\tilde
{\goth g})={k\over 2}.$

  \item $\goth g$ est de type DIII $\Leftrightarrow\ rg(\tilde
{\goth g})={k\over 2}$ et $r(\goth g)=[{n\over 2}].$
\end{enumerate}
\end{lem}

 \dem D'apr\ès la description de $\overline \∆_2,$ il est facile de
v\érifier que :
$\{\overline{\epsilon}_i-\overline{\epsilon}_{i+1}\ ,\
1\≤i\≤k-1,\mu\}$ avec :

\qquad $\bullet$ $\mu={\overline \epsilon}_{k-1}+{\overline
\epsilon}_k$ dans le cas orthogonal car $k\≥3,$

\qquad $\bullet$ $\mu=2{\overline \epsilon}_k$ dans le cas 
symplectique car $k\≥2,$

\noi est un ensemble de racines simples de $\overline{\tilde
{\goth g}},$ on en d\éduit le diagramme de Satake de 
$ \tilde{\goth g}$ \à partir du  diagramme de Satake de 
$  \goth g,$ d'o\ù le r\ésultat.  \fdem
 
 \subsubsection{ Sous-groupe parabolique standard tr\ès sp\écial et \PV associ\é}
 \bigskip
 
\noi On consid\ère le\sg parabolique $P(H_1,...,H_{p_0})=P_0$ donn\é dans le  $\S 2.4$ c'est \à dire que:\\
 
\qquad $\bullet$ Lorsque $\goth g$ est d\éploy\ée de type
$C_n:$ $H_{p_0-j+1}= 2 h_{ { 
\epsilon}_{2j-1}+
{ \epsilon}_{2j}}$ pour $j=1,...,p_0={k\over 2}\ ,$
 \vskip 2mm
 \qquad $\bullet$ dans tous les autres cas : $H_{p-j+1}=  
h_{{ \epsilon}_j}$ pour 
$j=1,...,p_0=p \ ,$\\
 \vskip 2mm
\noi ce qui donne dans les cas BI ou DI: $H_{p-j+1}=  h_{  
{\overline \epsilon}_j}$ et pour les cas CII et DIII:
$H_{p-j+1}= 2 h_{ { \overline
\epsilon}_{2j-1}+
{ \overline\epsilon}_{2j}}.$\\
\vskip 2mm
\noi On peut noter que $\forall \sigma\in \goth S_{p_0}, \exists
g_{\sigma}\in G$ tel que $\forall j=1,...,p_0$ on a
$g_{\sigma}(H_j)=H_{\sigma(j)}.$
\vskip 2mm
\noi Les descriptions de $\overline \∆_1$ et $\overline \∆_2$
  donnent
les d\écompositions suivantes particuli\èrement simples:\\
 
 $\begin{array}{lll}\qquad \bullet\  \goth
g_1&=&\oplus_{1\≤j\≤p_0}E_2(H_j)\cap \goth g_1\\
\\
 \qquad \bullet \ \goth g_2&=&\oplus_{1\≤j\≤p_0}E_4(H_j)\cap
\goth g_2
\oplus_{1\≤i<j\≤p_0}E_2(H_i)\cap E_2(H_j)\ . \end{array}(7.1)$\\

 \noi Soit:
$$d_0=\text{ dimension de } E_2(H_i)\cap E_2(H_j)\text{ pour }i\not =j =\begin{cases}
 1\text{ pour les types BI ou DI ,}\\
  4\text{  pour tous les autres cas (i.e C ou DIII),}\end{cases} $$
et on a:
$$\text{dim}(  E_4(H_j)\cap \g_2=\begin{cases}
 0\text{ pour les types BI ou DI ,}\\
  1\text{  pour   le type  DIII,}\\
  3\text{  pour   les types CI et CII.}\end{cases} $$

\noi  On rappelle que:\\
  
  \begin{lem} Pour $j=1,...,p_0$ on a degr\é de ($F_j)=d'.j$ avec $d'=\begin{cases}4 \text{ dans le cas }DIII,\\2\text{ sinon}.\end{cases}$
   \end{lem}
   
   \dem Il est bien connu par classification (\cite{rublivre},table I) que l'\irf de $(\overline \g_0,\overline \g_1)$ est de degr\é $2k$ dans les cas orthogonaux $(B_n$ ou $D_n,\alpha_k)$ et dans le cas symplectique $(C_n,\alpha_{2k}).$ Comme le \PV $(E_0(H_1)\cap \g_0,E_1(H_1)\cap g_1)$ est de type:
  
i) $(B$ ou $ D,\alpha_{k-1})$ lorsque $\overline R=B$ ou $D,$

ii)  $(C,\alpha_{2(k-1)})$ lorsque $\overline R=C,$
 
\noi (par exemple \demo de la proposition 1.2.4),
  on a le r\ésultat par r\écurrence sur $p_0.$\fdem\\
   
\vskip 4mm

\noi On rappelle \également que :\\

$\bullet\ N:=\ds\frac{\text{dim}(\g_1)}{d_{p_0}}=\begin{cases}
2(n-k)\text{ lorsque }(\overline \g_0,\overline \g_1)\text{ est de type } (C_n,\alpha_k),\\
n-k \text{ lorsque }(\overline \g_0,\overline \g_1)\text{ est de type } (D_n,\alpha_k),\\
n-k+\ds\frac{ 1}{2}  \text{ lorsque }(\overline \g_0,\overline \g_1)\text{ est de type } (B_n,\alpha_k).\end{cases}$\\
\bigskip

$\bullet\ \widetilde B=\ds\frac{-d_{p_0}B}{2B(H_0,H_0)}$ donc  $\widetilde B(H_j,H_j)=-2d'$ pour $j=1,...,p_0.$
\vskip 5mm

\noi Les coefficients des \équations fonctionnelles v\érifi\ées par les fonctions Z\étas sont  des sommes et produits de coefficients analogues associ\és aux \PVs obtenus \à partir de centralisateurs de tds,  ces derniers ont la m\ême  structure que le  \PV de d\épart,  en effet:  \\

\begin{prop} On suppose que $p_0\≥2.$

 Soit $X\in W_{\goth t_0},$ $\goth s$ l'alg\èbre
engendr\ée par les projections de $X$ sur $E_2(H_i)\cap \goth
g_1$ et de $X^{-1}$ sur $E_{-2}(H_i)\cap \goth
g_{-1},$ $\goth U=\goth U(\goth s)$ alors $(\goth U_0,\goth
U_1)$ est  de m\ême type que $(\goth g_0,\goth
g_1),$  plus pr\écis\ément:
\vskip 2mm
1) Lorsque $(\overline{\goth
g}_0,\overline{\goth
g}_1)$ est de type $(C_n,\alpha_{2p}),$
$(\overline{\goth
U}_0,\overline{\goth
U}_1)$ est de type $(C_{n-3},\alpha_{2(p-1)})$ et $\goth U$
est  \dpl (resp.de type CII) $\Leftrightarrow$ $\goth g$ 
est  \dpl (resp.de type CII).
\vskip 2mm
2)  
$(\goth U_0,\goth U_1)$
est de type BI(n-2,k-1)  (resp.DI(n-1,k-1))
lorsque $(\goth g_0,\goth g_1)$  est de type
 DI(n,k) (resp.BI(n,k) avec $n\≥5).$
\vskip 2mm
  3) 
$(\goth U_0,\goth U_1)$ est de type DIII(n-3,p-1) et rang $(\goth U)=[\frac{n-3}{2}]$ lorsque
$(\goth g_0,\goth g_1)$ est de type DIII(n,p) avec $p\≥3$  et rang$(\g)
=[\frac{n}{2}].$\end{prop}

 \dem 
Puisque tous les $H_i,i=1,...,p_0$ sont $G$-conjugu\és, il
suffit de le v\érifier pour $H_{p_0}.$ 

\noi Par la d\émonstration de la
proposition  1.2.4, $(\overline{\goth U}_0,\overline{\goth
U}_1)$ est de type :\\

\quad  i) $(D_{n-1},\alpha_{k-1})$ lorsque $(\goth
g_0,\goth g_1)$ est de type BI$(n,k),$

\quad  ii) $(B_{n-2},\alpha_{k-1})$
  lorsque $(\goth
g_0,\goth g_1)$ est de type  DI$(n,k),$

\quad  iii)
$(D_{n-3},\alpha_{2(p-1)})$ 
  dans le cas DIII$(n,p),$

\quad  iv) 
 $(C_{n-3},\alpha_{2(p-1)})$) dans les cas  symplectiques.

\noi ce qui termine la \demo  dans le cas ii) et dans le cas i) avec
$k$ pair.\\

\noi Il est facile de v\érifier que le diagramme de  Dynkin de
$\overline{\goth U}(\F H_{p_0})$ est donn\é par:
 
\hskip 2,9cm $\alpha_k$

\hskip 19pt \hbox to 5cm {\lower 2pt\hbox{$\circ$}\dotfill \hbox to 0,4 cm
{\offinterlineskip \lower 2pt\hbox{$\circ$} \hglue -3,2pt
{\vrule height  0,4pt depth 0pt width 0,5 cm}}\hglue 6,6pt

\lower 2pt\hbox{$\circledcirc$}
\hglue -4,2pt
{\vrule height  0,4pt depth 0pt width 0,5 cm}\hglue -1,2pt
\lower 2pt\hbox{$\circ$}
  \dotfill \hbox to 0,4 cm  
 { \offinterlineskip \lower 2pt\hbox{$\circ$} ...... }} \qquad avec 
 
\hskip 7mm $\alpha_\ell$ \hskip 0,9cm$\alpha_{k-1}$ \hskip 4mm  $\alpha_{k+1}$

$$ \ell= \begin{cases}
3 \hbox{ dans le cas  de type III ou le cas
symplectique}  \\
 2 \hbox{ sinon}  .\end{cases}$$ 
donc  $\goth U(\F H_{p_0})$ est de m\ême type que $\goth g$ et
$rg(\goth U(\F H_{p_0}))=\begin{cases}rg(\goth g)-2\text{ dans le cas CI},\\
                                                       rg(\goth g)-1\text{ dans les autres cas}\end{cases}.$
                                                       
                                                       \noi  Il en est de m\ême pour le \PV 
$(\goth U(\F H_{p_0})_0,\goth U(\F H_{p_0}) _1).
$

\noi Soit $\goth L$ l'alg\èbre de Lie
engendr\ée par $\goth U_2$ et
$\goth U_{-2},$ comme :
$$ \goth U_{\pm2}= \goth U(\F H_{p_0})_{\pm2}\ ,$$
On a:  
 
\quad $\bullet$ pour  $\goth g$ de type  BI:  $rg(\goth
L)=k-1$ lorsque $k-1\≥3$ (lemme 7.1.1 appliqu\é au \PV $(\goth U(\F H_{p_0})_0,\goth U(\F H_{p_0}) _1)$) donc $\goth U$ est de type I dans le cas i) avec $k\≥5$ et
impair 
 (lemme 7.1.1).\\

\noi   Lorsque $n\≥5$ et $k=3,$ $P(H_1,H_2)$ est un sous-groupe parabolique tr\ès sp\écial de $(\goth U_0,\goth U_1).$  Soit  $(x,H_1,x^{-1})$ un
$sl_2-$triplet du \PV
$(\goth U_0,\goth U_1)$, le \PV  $(\overline{\goth
U'}_0,\overline{\goth U'_1}),$ d\éfini par le
centralisateur dans $\goth U$ du $sl_2-$triplet
pr\éc\édent, est de type $(B_{n-3}, \alpha_1)$ puisque $(\overline{\goth
U}_0,\overline{\goth U_1})$  est de type $(D_{n-1},\alpha_2).$

\noi Or
les \PVs: $(\goth g'_0,\goth g'_1)$  de type
DIII$(n-1,1) $  admettent un seul sous-groupe parabolique tr\ès sp\écial standard non
trivial donn\é par
$P'_0$ (prop. 2.4.3) mais alors, soient $z\in
E'_2(H'_1)\cap \goth g'_1,$ $\goth s'$ l'alg\èbre engendr\ée
par   $z,$ 
$H'_1$ et $z^{-1},$ le \PV $(\overline{\goth U(\goth s')}_0
,\overline{\goth U(\goth s')}_1)$ est de type
$(A_3,\alpha_2)$  (cf.par exemple la \demo du lemme
7.4.1).

\noi Par cons\équent, $(\goth U_0,\goth U_1)$ est de type I.\\ 

\quad $\bullet$  pour  $\goth g$ de type CII (resp.DIII): 
$rg(\goth L)=p-1$ (lemme 7.1.1 appliqu\é au \PV $(\goth U(\F H_{p_0})_0,\goth U(\F H_{p_0}) _1)$) donc $\goth U$ est de type CII (resp.DIII car $p\≥3$) par le lemme 7.1.1.
 
\noi Il reste \à v\érifier que $rg(\goth U)\≥[ \frac{n-3}{2}]$ dans le
cas DIII(n,p) lorsque $n$ est  pair.   
 
\noi On peut supposer
que $X\in \goth g^{\epsilon_1-\epsilon_{p+1}}\oplus 
\goth g^{\epsilon_1+\epsilon_{p+1}}$ (cf.$\S 2.4$) ainsi
$\oplus_{2\≤i\≤p,\frac{n}{2}\≥j\≥p+2}\goth
g^{\epsilon_i\pm\epsilon_j}\subset \goth U_1$ (resp.
$\oplus_{2\≤i\≤p,\frac{n}{2}\≥j\≥p+2}\goth
g^{-\epsilon_i\pm\epsilon_j}\subset \goth U_{-1}$) donc
$rg(\goth U)\≥rg(\goth g)-2=\frac{n}{2}-2 $ d'o\ù le
r\ésultat.  \fdem

\vskip 5mm
\noi Pour terminer ce paragraphe de g\én\éralit\és, on redonne les
polynomes de Bernstein dans le cas  r\éel dans notre normalisation, bien qu'ils soient
d\éj\à connus (cf. travaux de \cite{clerc}, notamment le th\éor\ème
1).
\vskip 3mm
\noi Pour $x $ r\éel, on note $B_x$  le polynome: $
B_x(s)=s(s+x),$ alors:
\vskip 3mm
 \begin{prop}$\F=\R$
 
 Lorsque $({\overline \g}_0,{\overline \g}_1)$ est de type $(R_n,\alpha_k)$ et
que $F_1$ est de degr\é $2$ (i.e. \à l'exception du cas DIII), on a pour $j=1,...,p_0$  :
 $$ b_j(s_1,...,s_{p_0})=b^*_j(s_1,...,s_{p_0})=
\prod_{\ell=0}^{j-1}B_{N-(p_0-1)\frac{d_0 }{2}-1}\
(\sum_{p_0-\ell\≤i\≤p_0}s_i+\frac{d_0 }{ 2}\ell).$$
 \end{prop}
 
\dem 
 
\noi 
Par r\écurrence sur $p_0$ (donc sur $k$),
le cas $p_0=1$ r\ésultant
de la remarque  3.6.6, 1) puisque $F_1$ est une
forme quadratique et que $\widetilde B(\frac{H_1}{2},\frac{H_1}{2})=-1.$

\noi On suppose la proposition v\érifi\ée pour les \PVs:
$(P_0,\overline \g_1)$ de type $(R_m,\alpha_{q})$ avec
$1\≤q\≤k-1.$

\noi Lorsque le \PV:
$(P_0,\overline \g_1)$ est de type $(R_n,\alpha_k),$
soit $X\in W_{\goth t_0},I_r\subset \{1,...,p_0\}$ de cardinal
$p_0-r$ avec $1\≤r\≤p_0-1,$  
$\goth s_r$ l'alg\èbre engendr\ée par les projections de
$X$ et de $X^{-1}$ sur $E_{ 2}(\sum_{i\in I_r}H_i)\cap
\goth g_{  1}$ et $E_{ - 2}(\sum_{i\in I_r}H_i)\cap
\goth g_{ -1},$ on a:

$$\frac{dim(\goth U(\goth s_r)_1) }{ 2r}=
\ds\frac{\frac{dim(\goth g_1)}{
   p_0}.r-d_0r(p_0-r) }{ 2r}=
   \ds\frac{dim(\goth g_1)}{2p_0}-(p_0-r)\frac{d_0 }{ 2}\ ,$$
 donc les quantit\és qui suivent ont toutes la m\ême valeur ind\épendamment de $r:$ $$\frac{dim(\goth U(\goth s_r)_1) }{2r}-(r-1)\frac{d_0 }{ 2}-1=N-(p_0-1)\frac{d_0\ }{ 2}-1\quad (*)\ .$$

\noi Pour $j=1,...,p_0-1,$ on applique la  proposition 3.4.4 et
l'hypoth\èse de r\écurrence,  d'o\ù la formule donnant
$b_j$ en appliquant (*) ainsi que l'\égalit\é $b_j=b^*_j.$ 

\noi Pour le calcul de $b_{p_0}$ on applique 3.7.3 avec
$H_1,...,H_{p_0-1}$ et $H_{p_0}$, il est facile de v\érifier que la
constante
$r_{p_0-1}=(p_0-1)\frac{d_0 }{2}$  et on applique (*).

\noi Pour $b^*_{p_0}$ on
a:
$$ \begin{array}{llll}b^*_{p_0}(s)&=&b_{p_0}(s^*-(N-1)1_{p_0})\quad &\hbox{(lemme}\ 
3.7.2)\\ 
&=&b_1(s^*-(N-1)1_{p_0})\ b^*_{p_0-1}(s)\quad &\hbox{(lemme}\ 
3.7.2)\\
&=&b^*_1(s^*-(N-1)1_{p_0})\ b_{p_0-1}(s)\quad &\hbox{(v\érification
pr\éc\édente)}\\
 &=&b_{p_0}(s) \quad &\hbox{(lemme}\ 
3.7.2). \qquad\Box\end{array} $$

 \noi  Remarque :  Le calcul de    $N-(p_0-1)\frac{d_0 }{2}-1$
dans chaque cas donne : 
\vskip 2mm
\quad a) Dans le cas
symplectique, on a $k=2p_0$ donc  $N-(p_0-1) 
 \frac{d_0 }{2}-1=2(n-3p)+1.$  
\vskip 2mm
 \quad b) Dans le cas  orthogonal de type I, on a $k=p_0$ donc
$N-(p_0-1) \frac{d_0 }{2}-1=\ds\frac{2n-3k }{2}+\ds\frac{\delta-1 }{2}, $
 avec $\delta=0$
dans le cas $D$ et $1$ dans le cas $B.$   
\vskip 2mm
\noi Ainsi si on d\éfinit les constantes suivantes: 
$$p_C=q_C=1\quad ,\ \  p_D=p_B=\frac{1 }{2}=-q_D\quad ,\ \ 
q_B=0\ ,$$ 
not\ée $p_R$ et $q_R$
lorsque
$R=B,C$  on a:
$$ b_j(s_1,...,s_{p_0})=b^*_j(s_1,...,s_{p_0})=
\prod_{\ell=0}^{j-1}B_{p_R(2n-3k)+q_R}\
(\sum_{p_0-\ell\≤i\≤p_0}s_i+\frac{d_0 }{2}\ell).$$
  
 \begin{prop}
 Lorsque $(\goth g_1,\goth g_0)$ est de type
DIII(n,p), on a pour
$j=1,...,p$ :
$$ b_j(s_1,...,s_p)=b^*_j(s_1,...,s_p)=
\prod_{\ell=0}^{j-1}\ B_{n-3p-\frac{1 }{2}}\
(\sum_{p-\ell\≤i\≤p}s_i+\ell) B_{n-3p-\frac{1 }{2}}\
(\sum_{p-\ell\≤i\≤p}s_i+\ell+\frac{1 }{2}). $$
 \end{prop}
 
 \dem
 $(\overline{\goth g}_0,\overline{\goth g}_1)$ est de type
$(D_n,\alpha_{2p})$ et comme les polynomes de Bernstein ne
d\épendent pas de la forme r\éelle choisie mais des invariants
relatifs et de la normalisation choisie pour $B,$ on
consid\ère les polynomes obtenus dans la proposition 7.1.4
pour les formes de type I en notant que les normalisations sont les m\êmes.

Pour $j=1,...,p,$  soient :
$$b'_{2j}(s_1,...,s_{2p})=\prod_{\ell=0}^{2j-1}\ B_{n-3p-\frac{1 }{
2}}\ (\sum_{2p-\ell\≤i\≤2p}s_i+\frac{1 }{2}\ell)  \ ,$$
alors on a les \égalit\és suivantes :
$$ b_j(s_1,...,s_p)=b'_{2j}( 0,s_1,...,0,s_p)\ \ \hbox{et} \
\  {b_j}^*(s_1,...,s_p)={b'_{2j}}^*( 0,s_1,...,0,s_p)\ $$
d'o\ù le r\ésultat.\fdem

\bigskip

\subsection{ Le cas symplectique}

\bigskip

\noi Dans ce paragraphe, 
$\overline{\goth g}$ est une alg\èbre de Lie simple de type
$C_n$ avec $n\≥3,$  $(\overline{\goth
g}_0,\overline{\goth g}_1)$ est de type $(C_n,\alpha_{2p})$
avec $1\≤p\≤\ds\frac{n -2}{3}$ (cf.$n^\circ13$ de la classification de
\cite{satokimura}) et
 $ \goth g $ est soit d\éploy\ée (type CI) soit de type CII. \\

\noi Dans les cas r\éels et complexe, ils ont \ét\é \étudi\és par J.L
Clerc dans le cadre des repr\ésentations d'alg\èbres de Jordan
et dans le cas $\goth p$-adique, C.L.Pan a d\étermin\é
explicitement la fonction $Z_L$ associ\ée \à l'invariant
relatif fondamental.

\noi  Le proc\éd\é de descente s'appuie sur les r\ésultats du cas
$p=1$ donn\és dans le th\éor\ème 3.6.5 ce qui n\écessite le
r\ésultat suivant:\\

\begin{lem} Cas $p=1$

On  suppose que  $(\overline{\goth
g}_0,\overline{\goth g}_1)$ est de type $(C_n,\alpha_2)$
 $(n\≥3)$ alors F a pour discriminant 1 et $\gamma(F)=1$
dans le cas \dpl et $(-1)^n$ sinon.
\end{lem}

\bigskip

\dem 1) Dans le cas d\éploy\é, on a:
$$\goth g_1=\oplus_{3\≤j\≤n}\ (\ \goth
g^{\epsilon_1\pm\epsilon_j}\ \oplus\  \goth
g^{\epsilon_2\pm\epsilon_j}\ )\ ,$$
et comme $F(x)\not=0\Leftrightarrow\ (x,2H_0)$ se compl\ète
en un $sl_2$-triplet, il existe des constantes non nulles :
$a_i,b_i,i=1,...,n-2$ telles que:
$$F(\sum_{1\≤i\≤n-2}(x_iX_{\epsilon_1-\epsilon_{i+2}}+
y_iX_{\epsilon_2+\epsilon_{i+2}}+
z_iX_{\epsilon_1+\epsilon_{i+2}}+t_iX_{\epsilon_2-\epsilon_{i+2}})=
\sum_{1\≤i\≤n-2}\ (a_ix_iy_i\ +\ b_iz_it_i)$$
d'o\ù le r\ésultat.\\

2) On rappelle que $\overline \g=\g\otimes_{\F}\E$ est  d\éfinie sur l'extension $\E=\F[\sqrt
\epsilon];$ le diagramme de Satake de 
$\goth g,$ qui est de rang $q,$ est donn\é par:\\

\hskip 10pt \hbox to 5cm {\offinterlineskip \lower 2pt\hbox{$\bullet$} \hglue -3,2pt
{\vrule height  0,4pt depth 0pt width 0,5 cm}\lower 2pt\hbox{$\circ$}
\hglue -3,2pt
{\vrule height  0,4pt depth 0pt width 0,5 cm}\lower 2pt\hbox{$\bullet$}
  \dotfill \hbox to 0,3 cm  
  {\offinterlineskip\lower 2pt\hbox{$\circ$}}  
  \hglue -7,2pt
{\vrule height  0,4pt depth 0pt width 0,5 cm}\lower 2pt\hbox{$\bullet$}
  \dotfill \hbox to 0,4 cm 
    { \offinterlineskip \lower 2pt \hbox{$\bullet$}\hglue -3,2pt
{\vrule height  0,4pt depth 0pt width 0,5 cm}\kern -1pt   \hrulefill     \kern -3,4pt\lower 2pt \hbox{ 
\offinterlineskip  {$\bullet$} \hglue -7pt\vbox{{\hrule height 0,3pt width 0,6cm}\vskip 3,5pt{\hrule height 0,3pt width 0,6cm}\vskip 0,3pt }\hglue -12pt $>$   \hglue -1,8pt $\bullet$}}}

\hskip 2mm$\alpha_1$...\hskip 1,6cm....$\alpha_{2q}$\\

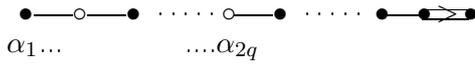
\noi Soit $\sigma$ la conjugaison associ\ée, un  \elt $x\in
\goth g_1$ s'\écrit:
$$ \begin{array}{lll}x&=&\sum_{1\≤i\≤2q}(x_i\overline
X_{\epsilon_1-\epsilon_{i+2}}+
\overline{x_i}\ \sigma (\overline
X_{\epsilon_1-\epsilon_{i+2}})+
y_i\overline
X_{\epsilon_1+\epsilon_{i+2}}+
\overline{y_i}\ \sigma(\overline
X_{\epsilon_1+\epsilon_{i+2}}))\\
&+&\sum_{1\≤i\≤n-2q}\
( z_i\overline X_{\epsilon_1-\epsilon_{i+2q}} +
 \overline{z_i}\ \sigma (\overline
X_{\epsilon_1-\epsilon_{i+2q}})+
 t_i\overline X_{\epsilon_1+\epsilon_{i+2q}} +
 \overline{t_i}\ \sigma (\overline
X_{\epsilon_1+\epsilon_{i+2q}})\ )\ ,\end{array}$$
\noi avec $x_i,y_i,z_i,t_i\in E.$

\noi En tenant compte des relations suivantes:\\

$\bullet$ $F(\sigma(x))=F(x),$

$\bullet$  $\sigma
(\epsilon_1\pm\epsilon_{2i+2})=\epsilon_2\pm\epsilon_{2i+1}$ et
$\sigma
(\epsilon_1\pm\epsilon_{2i+1})=\epsilon_2\pm\epsilon_{2i+2}$ pour $i=1,...,q-1,$

$\bullet$  $\sigma (\epsilon_1+\epsilon_j)=\epsilon_2-\epsilon_j$ et
$\sigma (\epsilon_1-\epsilon_j)=\epsilon_2+\epsilon_j$
pour $j=2q+1,...,n,$  \\

\noi il est
imm\édiat que $F$ est de la forme:
$$ \begin{array}{lll}F(x)&=&\sum_{1\≤i\≤q-1}\ (\
c_i\ x_{2i}\ \overline{y_{2i-1}}+ 
\overline{c_i\ x_{2i}}\
y_{2i-1}+d_i\ x_{2i-1}\ \overline{y_{2i}}+
\overline{d_i\ x_{2i-1}}\
y_{2i}\ )\\
&+&
\sum_{1\≤j\≤n-2q}\ (\ f_j\ z_j\overline{z_j}+ 
 g_j\ t_j\overline{ t_j}\ ) \quad ,\quad c_i,d_i\in \E^*\
,\ f_j,g_j\in \F^*.\end{array}$$ 

\noi Ainsi $F$ est la somme directe de
$q-1$ formes quadratiques
$Q_1$ de la forme $Q_1(x_1,x_2,y_1,y_2)=x\ \overline
y+\overline x\ y$ avec $x=x_1+\sqrt \epsilon\ x_2$ et 
$y=y_1+\sqrt \epsilon\ y_2,$ $Q_1$ a pour discriminant $1$ et
$\gamma (Q_1)=1,$ et de $n-2q$ formes quadratiques $Q_j$ de
la forme $Q_j(z_1,z_2,t_1,t_2)=\ f_j\ z\ \overline
z+\ g_j\ t\overline t $ avec $z=z_1+\sqrt \epsilon\ z_2$ et 
$t=t_1+\sqrt \epsilon\ t_2,$ qui sont de discriminant $1$
\également.\\

\noi Soit $j$ fix\é et soit $ \overline {\goth l}$ l'alg\èbre
engendr\ée par
$\overline {\goth g}^{\pm(\epsilon_1±\epsilon_j)}$ et 
$\overline {\goth g}^{\pm(\epsilon_2\pm\epsilon_j)},$ $ 
\overline {\goth l}$ est $\sigma-$stable et $\goth l= 
\overline {\goth l}\cap
\sigma( \overline {\goth l})$ est de rang $1,$ de  
diagramme de Satake: 

\hskip 10pt \hbox to 1,5cm{ \offinterlineskip \lower 2pt\hbox{$\bullet$}\hrulefill\lower 2pt\hbox{\offinterlineskip {$\circ$}\hglue -1,8pt\vbox{
 {\hrule height 0,3pt width 0,6cm}\vskip 3,5pt{\hrule height 0,3pt width 0,6cm}
\vskip 0,3pt} \hglue -16pt  $<$ \hglue -1,8pt  {$\bullet$}}} \\

\noi $Q_j$ est un invariant
relatif fondamental du \PV: $(\goth l\cap \goth g_0,\goth
l\cap\goth g_1)$ donc $Q_j$ est anisotrope d'o\ù
$\gamma(Q_j)=-1$ dans le cas r\éel et
$h(Q_j)=-(-1,-1)=-\alpha(1)^4$ dans le cas
$\goth p$-adique car le discriminant de $Q_j$ est $1$ donc 
$\gamma(Q_j)=-1.$

\noi Ainsi $F$ a pour discriminant $1$ et
$\gamma(F)=(-1)^{n-2q}=(-1)^n.$\fdem

\bigskip

On rappelle que:
$$b_{\g,P_0}(s_1,...,s_p)=\prod_{\ell=0}^{p-1}\biggl(\ (\sum_{p-\ell}^ps_i+2\ell) \ .
 (\sum_{p-\ell}^ps_i+2n-6p+2\ell+1)\ \biggr),$$
 et que pour $\pi=(\pi_1,...,\pi_p) \in (\widehat{\F^*})^p$ 
 $$\rho_{b_{\g,P_0}}(\pi)=\prod_{l=0}^{p-1}\rho \
(\pi_{p-l}...\pi_p|\ |^{2l+1}) \rho \
(\pi_{p-l}...\pi_p|\ |^{2(l+n-3p+1)}).$$
\noi Les facteurs $\rho$ ont \ét\é explicit\és dans
le $\S 3.6. $  \\

\begin{theo} 
 On suppose que $(\overline{\goth g}_0,\overline{\goth g}_1)$ est de
type $(C_n,\alpha_{2p})$ et soit $f\in S(\goth g_1),$ alors  pour $\forall f\in \EuScript S(\goth g_1)$ et $\pi \in (\widehat{\F^*})^p$ on a:
$$Z^*(\hat f;\pi)=\gamma^{np}\rho_{b_{\g,P_0}}(\pi) Z(f; 
\pi^*|\ |^{ -2(n-2p)1_p})$$  
avec $\gamma^n=\gamma(F_1),$ $\gamma=1$ (resp.$\gamma=-1$) lorsque
   $\goth g$  est  \dpl (resp.de type CII).
\end{theo}
 
 \dem On proc\ède par r\écurrence sur $p.$

\noi Le cas $p=1$ d\écoule du 2) du th\éor\ème 3.6.5, du lemme 7.2.1 et de
la normalisation choisie puisque dans ce cas $\g_1$ est de dimension $4(n-2)$ et $F_1$ de discriminant $1.$.

\noi On suppose le   th\éor\ème v\érifi\é lorsque  $(\overline{\goth
g}_0,\overline{\goth g}_1)$ est de type 
$(C_m,\alpha_{2(p-1)}).$ 

\noi Lorsque  $(\overline{\goth g}_0,\overline{\goth g}_1)$ est de
type $(C_n,\alpha_{2p}),$ on applique la proposition  
5.3.1 avec $H=\F^*$ et $k=1,$ puisque toutes les hypoth\èses de cette proposition sont v\érifi\ées. En effet, en reprenant les notations de la proposition  5.3.1  et $z\in W^*_{\goth t_0},$ soient $z_i$ (resp.$z_i^{-1})$ les
projections de $z$ sur $E_{-2}(H_i)\cap \goth g_{-1}$
(resp.sur
$E_2(H_i)\cap \goth g_1)$ et $\goth s_i$ l'alg\èbre engendr\ée
par $z_i$ et $z_i^{-1}$ alors $\goth U=\cap_{2\≤i\≤p}\ \goth
U(\goth s_i)$ or $\goth
U(\goth s_p)$ est de m\ême type que $\goth g$ (prop. 7.1.3) donc, par
 cette m\ême proposition, il en est de m\ême pour $\goth U,$ comme 
$(\overline{\goth U}_0,\overline{\goth U}_1)$ est de
type $(C_{n-3(p-1)},\alpha_2),$ on a pour $\pi_1\in \widehat{\F^*}:$
$$ b(\pi_1)=\gamma^{n-3p+3}\rho(\pi_1|\
|)\rho(\pi_1|\ |^{2(n-3p+1)}). $$
$(\overline{\goth U'}_0,
\overline{\goth U'}_1)$ est de type $(C_{n-3},\alpha_{2(p-1)})$
et $\goth U'$ est  \dple (resp. de type CII) lorsque
$\goth g$ l'est (prop.7.1.3), donc par r\écurrence pour $(\pi_1,...,\pi_{p-1})\in (\widehat{\F^*})^{p-1}:$
$$ c(\pi_1,...,\pi_{p-1})=\gamma^{(n-3)(p-1)}
\prod_{0\≤l\≤p-2}\rho \
(\pi_{p-1-l}...\pi_{p-1}|\ |^{2l+1}) \rho \
(\pi_{p-1-l}...\pi_{p-1}|\ |^{2(l+n-3p+1)}).$$
Notons que $r_1=2$ (cf. relation $7.1$)).\fdem\\

\bigskip

{\bf Remarques}

\begin{enumerate}
\item Dans le cas r\éel non \dpl, le\theor a \ét\é \établi dans un cadre plus g\én\éral  par J-L Clerc (th\éor\ème 2 de \cite{clerc}).

\item Dans le cas symplectique \dpl, on peut montrer que $\g"_1$ est une seule orbite pour $P_0$ \à l'aide du \theor de Witt.
 
\item On ne donne pas toutes les \équations fonctionnelles qui se d\éduisent du\theor 7.2.2 pour $G$ et $F_p$ en prenant $\pi_1=...=\pi_{p-1}=id$ (cf.$\S 6$).

\end{enumerate}

\bigskip

\subsection{Les cas
orthogonaux de type I  }
\vskip 3mm
\noi $(\g_0,\g_1)$ est de type $BI(n,k)$ ou $DI(n,k)$ avec l'hypoth\èse suppl\émentaire $3k\≤2(n-1)+\delta,$ et $P_0=P(H_1,...,H_k).$\\

\noi  Rappelons que $\delta=\begin{cases} 0\text{ dans le cas DI(n,k),}\\
1\text{ dans le cas BI(n,k).}\end{cases}$  \\

\noi On proc\ède comme dans le cas symplectique en ramenant le calcul des coefficients de l'\eq \à ceux associ\és aux formes quadratiques sous-jacentes, ce qui n\écessite dans cet exemple quelques choix suppl\émentaires.\\

\subsubsection{Normalisation des invariants relatifs}

Les
sous-espaces radiciels :
$\goth g^{\pm\epsilon_i\pm\epsilon_j}$ pour $1\≤i<j\≤k,$ sont de dimension $1$ et
engendrent une alg\èbre de Lie simple de type $D_k$ que l'on
munit d'un syst\ème de Chevalley, $(X_{\alpha}),$  tel que
les coefficients correspondants v\érifient 
les conditions:
$$\hbox{Soit}\quad
i<j\text{ et }\ell\not=i,j\quad \hbox{alors}\quad
 \ N_{\epsilon_i-\epsilon_j,\epsilon_j\pm\epsilon_{\ell}}=\begin{cases}-1\text{ si }
i<\ell<j \\
\ \ 1\text{ sinon.}\end{cases}.$$
Un tel syst\ème existe (\cite{bourbakigal8},chap.VIII,$\S13,n^\circ 4$),
on dit qu'il est $D-$adapt\é.Dans ces conditions on a:\\

\begin{lem} \begin{enumerate}
\item $N_{\epsilon_i-\epsilon_j, \pm\epsilon_{\ell}-\epsilon_i}=\begin{cases} 1\text{ si } i<\ell<j\text{ ou }j<i<\ell \text{ ou }\ell<j<i,\\
-1\text{ si } \ell<i<j \text{ ou }i<j<\ell \text{ ou }j<\ell<i.\end{cases}$\\

\item Pour $1\≤i<j\≤k,$  soit $\theta_{i,j}=
\theta_{X_{\epsilon_i-\epsilon_j},h_{\epsilon_i-\epsilon_j}}(-1),$ alors pour $\ell\not=i,j$ et $s,t=\pm 1,$ on a:
$$ \theta_{i,j}(X_{s\epsilon_i-t\epsilon_{ \ell}})=-N_{\epsilon_i-\epsilon_j, \epsilon_{\ell}-\epsilon_i}X_{s\epsilon_j -t\epsilon_{ \ell}}.$$
\end{enumerate}
\end{lem}

\dem 1) Lorsque $i<j$ on a:
$$N_{\epsilon_i-\epsilon_j, \pm\epsilon_{\ell}-\epsilon_i}=-N_{-\epsilon_i+\epsilon_j, \pm\epsilon_{\ell}-\epsilon_j}=-N_{\epsilon_i-\epsilon_j,\epsilon_j\mp\epsilon_{\ell}},$$
((3) du lemme 4 et proposition 7,$\S2,n^\circ4,$ chap.VIII,\cite{bourbakigal8})  
 et lorsque $i>j$ on a:
$$N_{\epsilon_i-\epsilon_j, \pm\epsilon_{\ell}-\epsilon_i}=N_{\epsilon_j-\epsilon_i,\epsilon_i \mp\epsilon_{\ell} } ,$$
 or le syst\ème est (D) adapt\é
d'o\ù le r\ésultat.

\noi 2) C'est un simple calcul, en remarquant que $N_{\epsilon_i-\epsilon_j,   \epsilon_{\ell}-\epsilon_i  }=N_{\epsilon_i-\epsilon_j,  \pm\epsilon_{\ell}-\epsilon_i  }.$\fdem

\bigskip

\noi On rappelle que $H_i=h_{\epsilon_{k-i+1}}$ pour $i=1,...,k$ et avec la normalisation choisie, on a pour toute racine  $\alpha =\pm\epsilon_i\pm\epsilon_j,$  avec $1\≤i<j\≤k,$  $\widetilde B(X_\alpha,X_{-\alpha})=1.$\\

\noi On d\éfinit les formes quadratiques suivantes pour $i=1,...,k$:
$$ x\in E_2(H_i)\cap \g_1\quad G_{i}(x)=b\begin{cases}\widetilde B(ad(x)^2(X_{-\epsilon_1-\epsilon_{k-i+1}}),X_{\epsilon_1-\epsilon_{k-i+1}})\text{ pour }1\≤i\≤k-1, \\
\\
\widetilde B(ad(x)^2(X_{-\epsilon_1-\epsilon_2}),
X_{-\epsilon_1+\epsilon_2})\text{ pour }i=k\end{cases}$$
avec $b\in\F^*.$\\

\noi Pour $i=1,...,k,$ $G_{i}$ repr\ésente l'invariant relatif
fondamental du \PV commutatif r\égulier:
 $$(\goth U(\oplus_{1\≤j\not=i\≤k}\F H_j)_0,\goth
U(\oplus_{1\≤j\not=i\≤k}\F H_j)_1=E_2(H_i)\cap \g_1),$$ c'est une forme
quadratique non d\ég\én\ér\ée et relativement invariante par
$\cap_{1\≤i\≤k}G_{H_i}.$

\noi Toutes ces formes quadratiques sont \équivalentes, en effet:

\vskip 3mm
 \begin{lem} 
 Pour $1\≤i<j\≤k$ on a
$G_{k-i+1}\circ \theta_{i,j}=G_{k-j+1}.$
\end{lem}
\vskip 3mm

\dem Pour $i\≥2$ ou bien $i=1$ et $j\≥3,$ le calcul donne:
$$ \theta_{i,j}(X_{\pm \epsilon_1-\epsilon_i})=- N_{\epsilon_i-\epsilon_j, \epsilon_1-\epsilon_i}X_{\pm \epsilon_1-\epsilon_j}=X_{\pm \epsilon_1-\epsilon_j} \quad \text{(lemme 7.3.1)} $$
et pour $i=1$ et $j=2$ on a $\theta_{1,2}(X_{ - \epsilon_1+\epsilon_2})=X_{  \epsilon_1-\epsilon_2}$ et $\theta_{1,2}(X_{  -\epsilon_1-\epsilon_2})=X_{ - \epsilon_1-\epsilon_2}.$  \fdem\\

{\bf Remarques:} 
\vskip 2mm
1) Le changement de syst\ème de Chevalley $D-$adapt\é a pour
cons\équence de multiplier toutes les formes quadratiques par
une m\ême constante  d\éfinie \à $\F^{*2}$ pr\ès. \\

\noi En effet, soit $X'_{\alpha}$ un autre syst\ème de Chevalley
 $D-$adapt\é alors il existe une
application $t$ du syst\ème de
racines 
$D_k$ dans
$\F^* $ telle que
  $X'_{\alpha}=t(\alpha )X_{\alpha}.$ Comme
les coefficients $N_{.,.}$ sont les m\êmes pour les deux
syst\èmes de Chevalley, $t$ est un morphisme du $\Z$-module associ\é \à $D_k$  
 et $t(\epsilon_1-\epsilon_
i)t(-\epsilon_1-\epsilon_i)t(\epsilon_1-\epsilon_2
)^{-1}t(-\epsilon_1-\epsilon_2)^{-1}\in \F^{*2}$ pour
$i=2,...,k.$
 \vskip 2mm
2) On rappelle que  dim$( \g^{\epsilon_1})=2(n-m)+\delta.$  \\

$\bullet$ Lorsque $\F=\R,$ on choisit la constante $b$
telle que $G_1$ soit de signature $(p,q)$ avec $p\≥q.$ \\  
 
$\bullet$ Dans le cas $BI,$   lorsque $\F$ est un corps $\goth p$-adique, on choisit la
constante
$b$ telle que discriminant$(G_1) =(-1)^{n-m}$ dans un souci de simplification des coefficients de l'\eq (cf.$\S 7.3.4$).\\

\noi Les invariants $F_1,F_2,...,F_k$ sont normalis\és par la
condition:
$$\hbox{pour}\quad x=\sum_{1\≤j\≤k}x_j\in W_{\goth t_0}\quad
F_i(x)=\prod_{1\≤j\≤i}G_j(x_j)\ ,$$
il est alors ais\é de v\érifier que:
$$\hbox{pour}\quad y=\sum_{1\≤j\≤k}y_j\in W^*_{\goth t_0}\quad
F^*_i(y)=\prod_{k-i+1\≤j\≤k}G^*_j(y_j)\ .$$\\

\noi Indiquons la nature de la forme quadratique $G_i$ et de sa restriction aux centralisateurs de tds.\\
 
 \begin{lem}    $G_1\sim G_a\oplus \widetilde G_1,$ $G_a=G_1/\goth g^{\epsilon_1}$ est une forme quadratique anisotrope, $\widetilde G_1$ \étant hyperbolique.\\

\noi  Cas r\éel:    $G_1$ a pour signature
$(p_0,q_0)$ avec $q_0=m-k$ et 
$p_0=2n-m-k+\delta $ ; on a : $
N=\displaystyle\frac{p_0+q_0 }{2}.$
\end{lem}

\dem 
 
\noi   1) Lorsque  $m=k,$  $\goth U(\oplus_{2\≤j \≤k}\F H_j)$ est de rang $1$ donc la forme
quadratique $G_1$ est anisotrope, donc d\éfinie positive par
normalisation dans le cas r\éel.

\noi  2) Lorsque  $m>k,$ prenons un syst\ème de Chevalley $D-$adapt\é de
l'alg\èbre simple de type $D_m$ d\éfinie par
$\{\pm\epsilon_i\pm\epsilon_j,1\≤i<j\≤m\}$ et  soit
$x\in E_2(H_1)\cap \goth g_1,$
$x=\sum_{k+1\≤\ell\≤m}(x_{\ell}X_{\epsilon_k-\epsilon_{\ell}}+
y_{\ell}X_{\epsilon_k+\epsilon_{\ell}})+Y,Y\in \goth
g^{ \epsilon_k },$ le calcul donne:
$$G_1(x)=2bb_1^2(\sum_{k+1\≤\ell\≤m}x_{\ell}y_{\ell})+G_1(Y)\ ,$$
$b_1^2$ \étant la constante provenant du syst\ème de Chevalley, or la restriction de $G_1$ \à $\goth
g^{ \epsilon_k }$ est anisotrope donc d\éfinie positive par choix de $b$ dans le cas r\éel d'o\ù
le r\ésultat.\fdem\\

\noi Lorsque $\Delta$ est de type $B_n,$ on note $Q_a$ la classe
d'\équivalence de la restriction d'une des formes quadratiques
$G_{k-i+1}$ \à $\goth g^{\epsilon_i};$ lorsque $\Delta$ est
de type $D_n,$   on pose $\forall t\in \F^*:$ 
$
\gamma_{\tau\circ
Q_a}(t)=1$ par extension. \\

\noi Soit $z\in
E'_2(H_1)\cap
\goth g_1=E'_2(h_{\epsilon_k})\cap
\goth g_1$ et
$\goth U^{(z)}=\goth U(\F z\oplus \F  H_1\oplus \F z^{-1}),$ alors pour
$i=2,...,k$ on a dim$(\goth U^{(z)}_1\cap E_2(H_i))=$
dim$(E_2(H_i)\cap \goth
g_1)-1$ et:\\

 \begin{lem} Pour   $i=2,...,k:$ $G_i\sim  G_i/_{\goth U^{(z)}_1\cap E_2(H_i)}\oplus G_1(z)X^2.$  
\end{lem}

 \dem  Pour   $i=2,...,k, $ on a:
$ E_2(H_i)\cap \g_1= \goth U^{(z)}_1\cap E_2(H_i)\oplus \F v,$ avec $v=[z,X_{-\epsilon_k+\epsilon_{k-i+1}}].$ 

1) Lorsque $k=2:$
$$[v,X_{ -\epsilon_1+\epsilon_2}]=[[z,X_{\epsilon_1-\epsilon_2}],X_{ -\epsilon_1+\epsilon_2}]=[z,[X_{\epsilon_1-\epsilon_2},X_{ -\epsilon_1+\epsilon_2}]]=-[z,h_{\epsilon_1-\epsilon_2}]=-z$$
car $[z,X_{ -\epsilon_1+\epsilon_2}]=0,$ donc 
 pour $y\in \goth U^{(z)}_1\cap E_2(H_2)$ on a:
$$[y,[v,X_{\ -\epsilon_1+\epsilon_2}]]=-[ y,z]=0$$
d'o\ù:
$$t\in \F\ :\ G_2(y+tv)=G_2(y)+t^2G_2(v)=G_2(y)+t^2G_1(z).$$

2) Lorsque $k\≥3.$

\noi Pour $i\≤k-1:$
$$[v,X_{\pm\epsilon_1-\epsilon_{k-i+1}}]= N_{-\epsilon_k+\epsilon_{k-i+1},\pm\epsilon_1-\epsilon_{k-i+1}}
[z,X_{\pm\epsilon_1-\epsilon_k }]=-[z,X_{\pm\epsilon_1-\epsilon_k }]$$
(1),lemme 7.3.1) car $[z,X_{\pm\epsilon_1-\epsilon_{k-i+1}}]=0$ donc $G_i(v)=G_1(z)$ et pour $y\in \goth U^{(z)}_1\cap E_2(H_i)$ on a:
$$[y,[v,X_{\pm\epsilon_1-\epsilon_{k-i+1}}]]= -[z,[y,X_{\pm\epsilon_1-\epsilon_{k}}]]=0,$$ d'o\ù:
$$t\in \F\ :\ G_i(y+tv)=G_i(y)+t^2G_i(v)=G_i(y)+t^2G_1(z).$$
De m\ême pour $i=k$ on a:
$$[v,X_{ -\epsilon_1\pm \epsilon_2}]= N_{\epsilon_1 -\epsilon_k, -\epsilon_1\pm\epsilon_ 2}
[z,X_{\pm\epsilon_2-\epsilon_k }]=[z,X_{\pm\epsilon_2-\epsilon_k }]\ \ (1),\text{lemme }7.3.1) $$donc pour $y\in \goth U^{(z)}_1\cap E_2(H_k)$ on a:
$$[y,[v,X_{ -\epsilon_1\pm\epsilon_k}]]= [z,[y,X_{ \pm\epsilon_2 -\epsilon_k}]]=0.$$ 
De plus:$$G_k(v)=b\widetilde B(ad(z)^2(X_{ \epsilon_2 -\epsilon_k}),X_{ -\epsilon_2 -\epsilon_k})=b\widetilde B(ad(z)^2(\theta_{1,2}(X_{ \epsilon_1 -\epsilon_k})),\theta_{1,2}(X_{- \epsilon_1 -\epsilon_k}))=G_1(z)$$
(2) du lemme 7.3.1) d'o\ù le r\ésultat.   \fdem

\bigskip

\begin{rema} disc$(G_i/_{\goth U^{(z)}_1\cap E_2(H_i)})=G_1(z).\text{disc}(G_1)$ et pour $t\in \F^*$  on a:
$$\gamma (tG_i/_{\goth U^{(z)}_1\cap E_2(H_i)})=\gamma
(tG_i)\alpha(-tG_1(z))=\gamma
(tQ_a)\alpha(-tG_1(z)).$$
\end{rema}  

\bigskip

\noi De la proposition 7.1.3, des $2$ lemmes  pr\éc\édents,  de leur \demo  et avec les normalisations choisies, on d\éduit
imm\édiatement le r\ésultat suivant:
 
 \begin{lem} Dans le cas r\éel, prenons $G_1(z)=\pm1$ alors  
  la restriction de $G_1$ \à $\goth U^{(z)}_1\cap E_2(H_i)$ a
pour signature $(p_0-\ds\frac{G_1(z)+1}{ 2},q_0-\ds\frac{1-G_1(z) }{2})$ et 
 $\goth U^{(z)}$ a pour rang :  $\hbox{rang}\ (\goth
g)-1+\ds\frac{G_1(z)-1 }{ 2}.$  
\end{lem}

\bigskip
 
\begin{rema} Lorsque rang$(\g)\≥2k$ on a rang$(\goth U^{(z)})\≥2(k-1).$    

\end{rema}
 En effet, comme $m>k,$  on prend un syst\ème de Chevalley de l'alg\èbre simple de type $D_m$ d\éfinie par
$\{\pm\epsilon_i\pm\epsilon_j,1\≤i<j\≤m\}$ et on peut supposer que
$z=cX_{\epsilon_k-\epsilon_{k+1}}+
dX_{\epsilon_k+\epsilon_{k+1}}$ (\à l'action des
automorphismes \él\émentaires pr\ès puisqu'on est dans un \PV
commutatif r\égulier) donc:
$$\goth U^{(z)}_1=\oplus_{1\≤j\≤k-1}\ \goth g^{\epsilon_j}
\ \oplus_{1\≤j\≤k-1,k+2\≤\ell\≤m}\ \goth
g^{\epsilon_j\pm\epsilon_{\ell}}
\ \oplus_{1\≤j\≤k-1}\ \F Y_j\ ,$$
$$\hbox{avec}\quad Y_j=cX_{\epsilon_j-\epsilon_{k+1}}-
dX_{\epsilon_j+\epsilon_{k+1}}\ .\qquad\Box$$

\bigskip

\subsubsection{Le r\ésultat ($\bf 3k\≤2(n-1)+\delta$)}
\vskip 3mm
\noi Pour $u=(u_1,...,u_k)\in (\F/\F^{*2})^k$ on  rappelle  que: 
\vskip 2mm 
$O_u=O_{u_1,...,u_k}=\{x\in \goth
g_1\ |\ F_1(x) \in u_1 \F^{*2}\ ,\ F_2(x) \in
u_1u_2\F^{*2}\ ,...,F_k(x) \in u_1...u_k\F^{*2} \},$\\

$O_u\cap W_{\goth t_0}=\{x=\sum_{1\≤i\≤k}x_i\in  W_{\goth t_0}|G_i(x_i)\in u_i\F^{*2} $ pour $i=1,...,k\},$\\

$O^*_u=O^*_{u_k,...,u_1} =\{x\in
\goth g_{-1}\ |\ F^*_1(x)\in u_k\F^{*2}\ ,\
F^*_2(x)\in
 u_{k-1}u_k\F^{*2}\ ,...,F^*_k(x)\in u_k...u_1\F^{*2}\},$\\
 
 $O^*_u\cap W^*_{\goth t_0}=\{x=\sum_{1\≤i\≤k}x_i\in  W^*_{\goth t_0}\ |\ G^*_i(x_i)\in u_i\F^{*2} $ pour $i=1,...,k\},$\\
\vskip 2mm 
\noi sont des ouverts (\éventuellement vides) et que 
  $Z_u(f;\pi)=Z (f{\bf 1}_{O_u};\pi)$  
 (resp.$Z^*_u(h;\pi)=Z^* (h{\bf 1}_{O^*_u};\pi)$) pour
$f\in S(\goth g_1)$ (resp.$h\in S(\goth g_{-1}).$ \\

\noi Dans le cas r\éel, on convient de prendre $u\in \{-1,1\}^k.$\\

\begin{lem} 
\begin{enumerate}

 \item  Dans le cas r\éel,
$O_u$ (resp.$O^*_u$) est non vide
$\Leftrightarrow$
Max $(-k,k-2q_0)\≤\sum_{1\≤i\≤k}u_i\≤$Min $(k,2p_0-k). $ \\

  \item  Dans le cas $\goth p$-adique, lorsque $2n-3k+\delta\≥3,$ $O_u$ (resp.$O^*_u$) est non vide pour tout $u\in (\F/\f)^k.$\\
  
\item Lorsque rang$(\goth g)\≥2k,$ $\forall u\in
(\F^*/\F^{*2})^k$ les ouverts $O_u$ et $O^*_u$ sont non vides. 
\end{enumerate}

\end{lem}

\dem  Dans les $3$ cas on proc\ède par r\écurrence sur $k,$ obtenue par descente.\\

\noi a) Lorsque $k=1,$ le r\ésultat est \évident dans le cas r\éel et dans le cas $\goth p$-adique puisque $n\≥3$ donc $\g_1$ est de dimension $\≥4$ d'o\ù $\F^*\subset G_1(\g_1)$ dans le cas $\goth p$-adique.

\noi Pour 3), on note que $G_1$ repr\ésente $0$ par
 le lemme 7.3.3 (cf sa d\émonstration).\\

\noi b) On suppose le lemme  v\érifi\é pour $k-1\≥1.$\\

\noi c) Lorsque $k\≥2,$ $Ou\not=\emptyset$ \ssi il existe $z\in E'_2(H_1)\cap \g_1$ tel que $G_1(z)=u_1$ et $O^{\goth U^{(z)}}_{u_2,...,u_k}\not=\emptyset$ dans le \PV $(\goth U^{(z)}_0,\goth U^{(z)}_1)$  utilis\é dans le lemme 7.3.4.

\noi $G_1$ est une forme quadratique non d\ég\én\ér\ée d\éfinie sur un espace vectoriel de dimension $2(n-k)+\delta$ et donn\ée dans le lemme 7.3.3, pour lequel on applique a).\\

\noi Le \PV $(\goth U^{(z)}_0,\goth U^{(z)}_1)$ est de type  
BI$(n-2,k-1)$ ou DI$(n-1,k-1)$ lorsque $n\≥5,$ par la proposition 7.1.3.

\noi Comme $2n-3k+\delta=2n'-3k'+\delta',$ avec $n'=n-2+\delta,$ $k'=k-1$ et $\delta'=1-\delta ,$ on conclut par r\écurrence pour 2).  
 
 \noi On proc\ède de m\ême dans le cas r\éel en appliquant \également le lemme 7.3.6.
 
 \noi Pour 3), on  applique la remarque 7.3.7 et l'hypoth\èse de r\écurrence  au \PV $(\goth U^{(z)}_0,\goth U^{(z)}_1)$ en notant que $G_1$ repr\ésente $0.$\\

\noi Dans le cas BI$(n,k)$ avec $n\≤4,$ on a $k=2$ puisque $k\≥2$ et que $3k\≤2(n-1)+\delta\≤7$ donc $n=4.$ 

\noi $(\overline {\goth U^{(z)}}_0,\overline {\goth U^{(z)}}_ 1)$ est de type $(A_3, \alpha_2)$ qui est encore un \PV dont  l'\irf est une forme quadratique non d\ég\én\ér\ée d\éfinie sur un espace vectoriel de dimension $4$ d'o\ù le r\ésultat dans le cas  $\goth p$-adique (2) et 3)).

\noi Dans le cas r\éel, on applique le lemme 7.3.6 pour 1).

\noi Lorsque $\g$ est de rang$\≥4,$ $\g$ est d\éploy\ée donc $G_2$ est de signature $(3,2)$ d'o\ù $G_2/\goth U^{(z)}_1$ repr\ésente $0$ par le  lemme 7.3.6.  
\fdem

\bigskip

{\bf Notations:}
\vskip 2mm 
\noi Pour $(u,v)\in (\F/\F^{*2})^k,$ $f$ application de $\F^*/\F^{*2}$ \à valeur complexe et $\pi_1,\pi_2\in
\widehat{\F^*},$  on note:
$$A_{v,u}^{(f)}( \pi_1;\pi_2)=\sum_{t\in \F^*/\F^{*2}}f(t)\rho(\pi_1|\ | ;tv)\rho(\pi_1\pi_2;tu).$$

On rappelle que $N=\ds\frac{\text{dim}(\goth g_1)}{  2k}=n-k+\frac{\delta}{2}.$\\
  
\begin{theo} 
 $ \goth g $ est de
type DI ou BI(n,k) avec $\bf 3k\≤2(n-1)+\delta.$\\

\noi Soit $v\in  (\F^*/\F^{*2})^k$ tel que $O^*_v$ soit
non vide,
  alors pour  $f\in
S(\goth g_1)$ et $\pi=(\pi_1,...,\pi_k)\in ( \widehat{\F^*})^k $ on a:
$$Z^*_v(\four f;  \pi)= \sum_{u\in
  (\F^*/\F^{*2})^k} a_{v,u}(\pi)Z_u(f; 
 \  \pi^* |\ |^{ - N1_k}),$$
avec $ a_{v,u}(\pi)=0$ lorsque $O^*_v$ est  
vide et sinon:
$$ a_{v,u}(\pi)= 
\prod_{1\≤\ell\≤k} 
A_{v_{\ell},u_{\ell}}^{(f_{\ell})} (
 \pi_{k-\ell+1}....\pi_k |\ |^{\frac{1 }{2}(\ell-1)},|\
|^{ \frac{1}{2}(2n-3k+\delta+1 )})\ ,\hbox{avec}$$
$$f_{\ell}(t)=\begin{cases}\gamma ( tQ_a)\text{ lorsque } k=1\\
\\
\gamma( tQ_a) 
\prod_{ 2\≤j\≤k}\alpha(-tv_j)\text{ lorsque }\ell=1\text{ et }k\≥2,\\
\\
\gamma ( tQ_a)\prod_{1\≤i\≤\ell-1}\alpha(-tu_i)
\prod_{\ell+1\≤j\≤k}\alpha(-tv_j)\text{ lorsque }2\≤\ell\≤k-1,\\
\\
\gamma ( tQ_a)\prod_{1\≤i\≤k-1}\alpha(-tu_i)
 \text{ lorsque }\ell=k .\end{cases}$$
  \end{theo}
 
 \bigskip
 
 \dem On proc\ède par r\écurrence sur $k,$
le cas $k=1$ est donn\é dans le th\éor\ème
3.6.5,5), on notera que $\widetilde B(\ds\frac{H_1}{2},\ds\frac{H_1}{2})=-1. $ 

\noi On suppose le r\ésultat vrai pour $k-1\≥1.$

\noi Lorsque $k\≥2,$ on applique la proposition 5.3.1, dont on reprend
les notations, avec l'\elt $1$-simple $H_1$ .

\noi Soient $w=(u_1,v_2,...,v_k)$ et $z=\sum_{1\≤i\≤k}y_i \in O^*_w\cap W^*_{\goth t_0},$ avec   $y_i\in E_{-2}(H_i)\cap \g_{-1},i=1,...,k,$ 
$\goth U=\goth U(\goth s)$ (resp.$\goth U'=\goth U(\goth
s')$), $\goth s$ (resp.$\goth s')$ \étant l'alg\èbre engendr\ée
par $z_0=\sum_{2\≤i\≤k}y_i$ et $z_0^{-1}=\sum_{2\≤i\≤k}y_i^{-1}$ (resp.$z_{-2}=y_1$
et $z_{-2}^{-1}$).\\

\noi $\goth U'$ est de type  
BI$(n-2,k-1)$ ou DI$(n-1,k-1)$ lorsque $n\≥5,$ par la proposition 7.1.3.
On a $2(n'-1)-3k'+\delta'=2(n-1)-3k+\delta\≥1,$ avec $n'=n-2+\delta,k'=k-1$ et $\delta'=1-\delta.$   

\noi Dans le cas BI$(n,k)$ avec $n\≤4,$ on a $k=2$ en raison de la condition impos\ée alors $(\overline {\goth U'}_0,\overline {\goth U'}_ 1)$ est de type $(A_3, \alpha_2)$ qui est encore un \PV dont l'invariant est une forme quadratique non d\ég\én\ér\ée.

\noi Dans les $2$ cas, les
coefficients sont  donn\és par la relation:
$$c_{(v_2,...,v_k),(u_2,...,u_k)}^{(z_{-2})}(\pi")=\prod_{1\≤\ell\≤k-1}A^{(f'_\ell)}_{v_{\ell+1},u_{\ell+1}}(\pi_{k-\ell}....\pi_k|\ |^{\frac{\ell}{2}};|\
|^{ \frac{1}{2}(2n-3k+\delta+1 )}),$$
ceci en appliquant l'hypoth\èse de r\écurrence lorsque $k\≥3$ et pour $k=2$ le th\éor\ème  3.6.5,5), 
avec $f'_\ell(t)=f_{\ell+1}(t)$ par la remarque 7.3.5 et par l'hypoth\èse de r\écurrence lorsque $k\≥3.$  \\

\noi $\goth U=\cap_{2\≤i\≤k}\ \goth
U(\F y_i\oplus \F H_i\oplus\F y_i^{-1})$ or $\goth
U(\F y_k\oplus \F H_k\oplus\F y_k^{-1})$ est de  
type   DI$(n-1,k-1)$ ou
BI$(n-2,k-1)$ par la proposition 7.1.3 lorsque $k\≥3,$ ou bien $k=2$ mais alors, comme ci-dessus, l'\irf est une forme quadratique non d\ég\én\ér\ée. 

\noi Dans tous les cas l'\irf du \PV $(\goth U_0,\goth U_1)$ est  la restriction de la forme quadratique non d\ég\én\ér\ée $G_1$ \à $\goth U_1$ de dimension $2n-3k+\delta+1(\≥3)$ et pour $t\in \F^*$ on a:
 $$\gamma (tG_1/\goth U_1)=\gamma
(tG_1)\prod_{2\≤i\≤k}\alpha (-tG^*_i(y_i))$$
en it\érant la remarque 7.3.5 donc 
$$b^{(z_0)}_{v_1,u_1}(\pi')=A^{(f_1)}_{v_1,u_1}(\pi_k;|\
|^{ \frac{1}{2}(2n-3k+\delta+1 )}),$$ 
 on termine en 
appliquant la proposition 5.3.1.\fdem

 \bigskip 
  
\noi Lorsque $\pi_\ell=\tilde\omega_{\epsilon_\ell}|\ |^{s_\ell},$ avec $\epsilon_\ell\in \F^*/\F^{*2},$ pour $\ell=1,...,k$  on a:

$$\begin{array}{cc}a_{v,u}(\pi_k\pi_{k-1}^{-1},...,\pi_2\pi_1^{-1},\pi_1)&= (\prod_{\ell=1}^k  (\epsilon_\ell,u_\ell v_\ell )\ )a_{v,u}( {s_k-s_{k-1}},..., {s_1})\text{ et} \\
\\
a_{v,u}( {s_k-s_{k-1}},... ,{s_1})&=\prod_{\ell=1}^k A_{v_\ell,u_\ell}^{(f_\ell)}(|\ |^{s_\ell+\frac{1 }{2}(\ell-1)};|\
|^{ \frac{1}{2}(2n-3k+\delta+1 )}).\end{array}$$\\

 \subsubsection{ Le cas r\éel $(p_0+q_0=2n-2k+\delta)$}
\vskip 3mm
On rappelle que $\goth g_1$ peut s'identifier \à
$M_{k,p_0+q_0}(\R)$ et que pour $\ell=1,...,k$ $F_\ell$ est alors proportionnel \à:
$$det_\ell(AI_{p_0,q_0}{^tA})\text{ avec }A\in M_{k,p_0+q_0}(\R) \text{ et }I_{p_0,q_0}=\left(\begin{array}{cc}1_{p_0}&0_{p_0}\\0_{q_0}&-1_{q_0}\end{array}\right)\ ,$$
$det_\ell(X)$ \étant le d\éterminant de la matrice tronqu\ée $X_\ell$ d\éfinie pour $\ell=1,...,k$ par :
$$  X={(x_{i,j})}_{1\≤i,j\≤k}\text{ alors }X_\ell={(x_{i,j})}_{k-\ell+1\≤i,j\≤k} \ \ .$$
 
Soient $ s_1,...,s_k\in \C$ on a:
$$ a_{v,u}( {s_k-s_{k-1}},... ,{s_1})=
C(s_1,...,s_k) 
{a'_{v,u}}^{(p_0,q_0;k)}(s_1,...,s_k)\ \ \hbox{avec :}$$

$$C(s_1,...,s_k)=2^k.\ (2\pi)^{-k(\frac{p_0+q_0+2 }{2})}(2\pi)^{-2(s_1+...+s_k)}\ .\prod_{\ell=1}^k\Gamma
(s_\ell+\frac{\ell+ 1 }{2})\Gamma  (s_{k-\ell+1}+\frac{p_0+q_0-\ell+1 }{ 2})$$

$$ {a'_{v,u}}^{(p_0,q_0;k)}(s_1,...,s_k)=
\prod_{\ell=1}^k\cos(\frac{\pi }{ 4}\ .\
\phi_\ell^{(p_0,q_0;k)}(v,u;s_\ell)
\ ,\
\hbox{avec}$$
$$ \phi_\ell^{(p_0,q_0;k)}(v,u;s_\ell)=2(u_\ell+v_\ell) \
(s_\ell+\frac{\ell }{2}) + u_\ell(p_0+q_0-k)+
(q_0-p_0)+\sum_{\ell\≤j\≤k}v_j +\begin{cases}0\text{ lorsque  }\ell=1,\\
 \sum_{1 \≤i\≤\ell-1}u_i  
\hbox{ lorsque }  \ell\≥2.\end{cases}$$
not\é \également: $a'_{v,u}$ et $\phi_\ell(v,u;s_\ell),$ lorsqu'il n'y a
pas d'ambiguit\é.\\

\begin{rema}
 \begin{enumerate}
 
 \item  Comme $ \phi_\ell^{(p_0,q_0;k)}(-v,-u;s_\ell)=- \phi_\ell^{(q_0,p_0;k)}(v,u;s_\ell),$ on a:
 $$ {a'_{-v,-u}}^{(p_0,q_0;k)}(s_1,...,s_k)= {a'_{v,u}}^{(q_0,p_0;k)}(s_1,...,s_k )\quad(R9).$$

\item Lorsque $p_0\≥k$ on a $a'_{(1,...,1),(1,...,1)}(s_1,...,s_k)= \prod_{\ell=1}^k\cos(\pi (s_\ell+\ds\frac{\ell+q_0}{2})).$\\

 \item Dans le cas particulier o\ù  $q_0=0,$ ce qui correspond \à $\g$ de rang $k,$ $\g"_1=O_{(1,...,1)}$ et $\g"_{-1}=O^*_{(1,...,1)}$ donc pour  $f\in
S(\goth g_1)$ et $s=(s_1,...,s_k)\in ( \C^*)^k $  on obtient:
$$Z^* (\four f; (s_k-s_{k-1},...,s_2-s_1,s_1))= C(s_1,...,s_k)\prod_{\ell=1}^k\cos(\pi s_\ell+\frac{\pi }{ 2}\ell)Z (f; s_2-s_1,...,s_k-s_{k-1},-s_k-\frac{p_0}{2}),$$
qui n'est qu'un cas tr\ès particulier du\theor 3 de \cite{clerc}.

\item Pour $u=(1,...,1,-1,...,-1)$ avec $u_1=...u_p=1$ et $u_{p+1}=...=u_k=-1$ on a:
$$a'_{u,u}(s_1,...,s_k)=\prod_{\ell=1}^p\cos\pi (s_\ell+\ds\frac{\ell+q_0-q}{2})\prod_{\ell=p+1}^k\cos\pi (s_\ell+\ds\frac{\ell+p_0-p}{2})\quad \text{avec}\ q=k-p.$$
\end{enumerate}
 \end{rema}
  
\vskip 3mm
La plupart
de ces coefficients sont nuls, en effet pour $u\in
\{-1,1\}^k,$ soit $p(u)$ (resp. $q(u)$) le nombre de
composantes positives (resp. n\égatives) de $u$
($p(u)-q(u)=\sum_{1\≤i\≤k}u_i,$  $p(u)+q(u)=k$) alors:\\

 \begin{lem}
 \begin{enumerate}
 \item Soient $u\not= v,$ on note $i_1,...,i_r,$ avec $ i_1<...<i_r,$ les indices
pour lesquels $u_{i_j}=-v_{i_j}.$ 

\quad i) Lorsque  $p_0+q_0+k$ est pair : $a_{v,u}=0$ pour
$r\≥2$ et pour $r=1$ on a $a_{v,u}\not=0\Leftrightarrow \
p(v)+p_0$ est pair.

\quad ii) Lorsque  $p_0+q_0+k$ est impair alors
$a_{v,u}\not=0\Leftrightarrow$ $v_{i_j}=(-1)^{p(v)+p_0+j-1}$
  pour $j=1,...,r.$ 
 \item
$a_{v,u}( s_1,...,s_k)=0$ lorsque
$|p(u)-p(v)|\≥2.$
 \end{enumerate}
 \end{lem}

\vskip 3mm
 \dem
Cela d\écoule simplement de  la formule:
$$j=1,...,r:\quad
\phi_{i_j}(v,u;s_{i_j})=(u_{i_j}-1)p_0+(q_0-k)(u_{i_j}+1)
 +2p(v)+2\alpha_j$$
avec $\alpha_1=0$ et pour
$j\≥2:\alpha_j=\sum_{1\≤\ell\≤j-1}u_{i_\ell}.$
\vskip 2mm
\noi Comme $p_0,q_0,p(v),\alpha_j$ sont entiers, $\frac{\pi }{
4}\phi_{i_j}\in \frac{\pi }{2}\Z;$
si   $u_{i_j}=u_{i_{j+1}}$ pour une valeur $j$ lorsque
$r\≥2,$ on a alors:
$$\phi_{i_{j+1}}=\phi_{i_j}+2u_{i_j}\ \ \hbox{ donc}\ \ 
 \cos(\frac{\pi\ }{ 4}\ \phi_{i_j}(v,u;s_{i_j})\ )  \cos(\frac{\pi }{
4}\
\phi_{i_{j+1}}(v,u;s_{i_j})\ )  =0.$$ 
Ainsi il reste \à regarder les cas pour lesquels  pour
$j=1,...,r-1$ on a $u_{i_{j+1}}=-u_{i_j}$ mais alors 
$\phi_{i_{j+2}}=\phi_{i_j}$ donc il suffit de regarder
$\phi_{i_1}$ et $\phi_{i_2}$ lorsque $r\≥2$ et $\phi_{i_1}$
lorsque $r=1$ d'o\ù le r\ésultat.

\noi 2) Se d\éduit simplement de 1).\fdem 
  
\vskip 3mm

{\bf Remarques:}
\vskip 3mm
1) De la d\émonstration du lemme pr\éc\édent, on d\éduit
imm\édiatement $a'_{v,u}$ lorsqu'il est non nul.
\vskip 2mm
\noi En effet, soit $I=\{i_1,...,i_r\}=\{i\ |\ u_i+v_i=0\ \},$ on
a :
$$a'_{v,u}(s_1,...,s_k)=B(I)\prod_{1\≤l\≤k,l\notin I}\cos(\
\frac{\pi
 }{2}\ [\ (2s_l+l)+(\frac{1-v_l }{2})p_0+(\frac{1+v_l }{2})(q_0-k)+v_lp(v)+v_lC_l(u)\ ]\ )$$
avec :
\vskip 2mm
$\bullet$ $C_l(u)=0$ si $I=\emptyset$ ou bien si  
  $l\≤i_1$ 
 sinon $C_l(u)=  \sum_{\{i\in I|1\≤i\≤l-1\}}u_i,$\vskip 2mm
$\bullet$ $B(I)=1$ si $I=\emptyset$ et sinon 
 $B(I)=(\ \cos({\pi \over 4}\phi_{i_1})\ )^{r-[{r\over
2}]}\ (\ \cos({\pi \over 4}\phi_{i_2})\ )^{[{r\over
2}]}.$
\vskip 4mm

2) Notons que l'on retrouve les
coefficients de l'\équation fonctionnelle associ\ée \à
$F=F_k,$ en effet il est clair que pour tout couple
d'entiers $(p,q)$ tels que :
$$p+q=k\ \ \hbox{et v\érifiant :}\ \ 
\hbox{Max}(-k,k-2q_0)\≤p-q\≤\hbox{Min}(k,2p_0-k)\ \ ,$$ 
l'ensemble : $O_{p,q}=\cup_{p(u)=p}O_u$
(resp.$O^*_{p,q}=\cup_{p(u)=p}O^*_u$) est inclus dans une seule orbite
sous l'action de $G.$ 

\noi  Pour $s\in \C,$ posons:
$${a'_{(p',q');(p,q)}}^{(p_0,q_0;k)}(s)=\sum_{v|p(v)=p'}{a'_{v,u}}^{(p_0,q_0;k)}(s,...,s)  \quad (R1)$$
avec $u$ quelconque tel que $p(u)=p$ (donc $q(u)=q)$ alors de mani\ère classique ${a'_{(p',q');(p,q)}}^{(p_0,q_0;k)}(s)$ ne d\épend pas du choix de $u$ (cf.2) de la remarque qui suit le  th\éor\ème 6.2.2) et on a:
$$a_{(p',q');(p,q)}( |\ |^s) =C(s,...,s)
 {a'_{(p',q');(p,q)}}^{(p_0,q_0;k)}(s).$$
Les coefficients $ {a'_{(p',q');(p,q)}}^{(p_0,q_0;k)}(s)$ sont bien connus, ils sont donn\és dans le  th\éor\ème A.1  de  \cite{sato4} puisqu'on a $ {a'_{(p',q');(p,q)}}^{(p_0,q_0;k)}(s)=C_{(p',q')}^{(p,q)}(I_{p_0,q_0};s)$.\\

\noi En effet, de mani\ère imm\édiate, on peut noter que:\\

i) $ {a'_{(p',q');(p,q)}}^{(p_0,q_0;k)}(s)=0$ si $|p-p'|\≥2$ (2) du lemme 7.3.11).\\

ii) $ {a'_{(q',p');(q,p)}}^{(q_0,p_0;k)}(s)= {a'_{(p',q');(p,q)}}^{(p_0,q_0;k)}(s)$ (1) remarque 7.3.10).\\

iii) Lorsque $p_0\≥k$ on a $$ {a'_{(k,0);(k,0)}}^{(p_0,q_0;k)}(s) =\prod_{\ell=1}^k\cos\pi (s+\ds\frac{\ell+q_0}{2}) \quad \text{(2) remarque }7.3.10).$$

iv) $$\begin{array}{ccc} {a'_{( p,q);(p,q)}}^{(p_0,q_0;k)}(s)
&=& \prod_
{\ell=1}^p\cos\pi (s +\ds\frac{\ell+q_0-q}{2})\prod_{\ell=p+1}^k\cos\pi (s +\ds\frac{\ell+p_0-p}{2})\\
&=&(-1)^{pq}\ \prod_{\ell=1}^q\cos\pi (s +\ds\frac{\ell+p_0}{2})
\prod_{\ell=q+1}^k\cos\pi (s +\ds\frac{\ell+q_0}{2}),\end{array}$$ lorsque $q_0-q\equiv p_0-p$ $(2)$ ou bien $p_0-p\equiv 0$ $(2)$ et $q_0-q\equiv 1$ $(2),$ par le 1) du lemme 7.3.11.\\

\noi Cependant le calcul g\én\éral est fastidieux car il fait apparaitre des sommes, il est alors plus simple de noter que:\\

v)$$\left(\begin{array}{cc}  {a'_{( 1,0);(1,0)}}^{(p_0,q_0;1)}(s)& 
{a'_{( 1,0);(0,1)}}^{(p_0,q_0;1)}(s)\\
 {a'_{( 0,1);(1,0)}}^{(p_0,q_0;1)}(s)& {a'_{( 0,1);(0,1)}}^{(p_0,q_0;1)}(s)\end{array}\right )=\left(\begin{array}{cc}  \cos \pi (s+\ds\frac{q_0+1}{2})&\sin\ds\frac{\pi p_0}{2}\\
\sin\ds\frac{\pi q_0}{2}&\cos \pi (s+\ds\frac{p_0+1}{2})
\end{array}\right ).$$

vi) Comme:
$$
\phi_k(v,u)^{(p_0,q_0;k)}(s)=\phi_1(v_k,u_k)^{(p'_0,q'_0;1)}
(s+\frac{k-1 }{2})\quad
(R2) $$ avec $p'_0=p_0-p(u)+\frac{1+u_k }{ 2}$ et 
$q'_0=q_0-q(u)+\frac{1-u_k }{ 2}$ (rappelons que $p(u)\≤p_0$ et
que $q(u)\≤q_0$)

\noi et que:
 $$ \prod_{l=1}^{k-1}\cos (\frac{\pi }{
4}\phi_l(v,u)^{(p_0,q_0;k)}(s)=
\prod_{l=1}^{k-1}\cos (\frac{\pi }{
4}\phi_l(v',u')^{(p"_0,q"_0;k-1)}(s)\quad
(R3)$$ avec
$u'=(u_1,...,u_{k-1}),$ $v'=(v_1,...,v_{k-1}),$
$p"_0=p_0-\frac{1+v_k }{2}$ et $q"_0=q_0-\frac{1-v_k }{ 2},$\\
 
\noi  les coefficients ${a'_{(p',q');(p,q)}}^{(p_0,q_0;k)}(s)$ v\érifient les relations
de r\écurrence du lemme A.3 de \cite{sato4}, relations de (A4) \à (A7) qui proviennent des relations $(R1),(R2),(R3).$ 

\noi Notons que (v) est la relation (A8), ii) est la relation (A9), iii) est le lemme A.4,  i)  le lemme A.6 et iv) une partie du lemme A.8.\\

\noi Ainsi  on a bien l'\égalit\é $ {a'_{(p',q');(p,q)}}^{(p_0,q_0;k)}(s)=C_{(p',q')}^{(p,q)}(I_{p_0,q_0};s),$ ces derniers coefficients \étant donn\és dans le th\éor\ème A.1 de \cite{sato4} et c'est l'\établissement des formules de r\écurrence
de cette partie de \cite{sato4} qui sont \à l'origine de ce
travail.

\noi Cependant la constante $C(s,...,s)$ diff\ère de la constante donn\ée par F.Sato qui n'a pas de facteur de la forme $2^{...}.$
\vskip 4mm

\subsubsection{Le cas $\goth p$-adique}
\vskip 4mm

\noi Posons $r=2n-3k+\delta$ et pour $(s_1,...,s_k)\in \C^k,$  $u,v\in (\F^*/\F^{*2})^k$ tels que $O_u$ et $O^*_v$  soient non vides:
$$\begin{array}{ccc}\widetilde{a}_{v,u}(s)&:=&a_{v,u}(s_k-s_{k-1},...,s_2-s_1,s_1)\\
&=&\prod_{\ell=1}^k A_{v_\ell,u_\ell}^{(f_\ell)}(|\ |^{s_\ell+\frac{1 }{2}(\ell-1)},|\
|^{ \frac{1}{2}( r+1 )}).\end{array}$$
Le calcul explicite donne:
$$\widetilde a_{v,u}(s)=\alpha (-1)^{ka_0}\ \gamma (Q_a)^k\ C(v,u)\prod_{1\≤\ell\≤k}A^{a_0}_{ s_\ell+\frac{1 }{
2}(\ell+1),s_\ell+ \frac{1 }{2}( r+\ell)}(v_\ell,u_\ell, w_\ell\delta'_0) $$ avec 
$\delta'_0=(-1)^{[\frac{k-\delta} {2}]}\delta_0$ et $\delta_0=(-1)^{n-m}.\text{discr}(Q_a),$ on rappelle que $\delta_0=1$ lorsque $\delta=1$ par normalisation ($\S 7.3.1$) et lorsque  $\delta=0$ pour $m=n$ ou $n-m=2,$
$$  a_0=\begin{cases}0\text{ lorsque }k+\delta \text{ ( i.e. }r )\text{ est impair ,}\\
1\text{ lorsque }k+\delta\text{ ( i.e. }r )\text{ est pair.}\end{cases}  \ ,
  \ C(v,u)=\prod_{1\≤i\≤k}\ \alpha(-u_i)^{k-i}\
\alpha(-v_i)^{i-1},$$
$$ w_\ell=\begin{cases}v_2.....v_k\text{ lorsque }\ell=1,\\
u_1.....u_{\ell-1}.v_{\ell+1}.....v_k \text{ lorsque }2\≤\ell\≤k-1,\\
u_1......u_{k-1}\text{ lorsque }\ell=k.\end{cases}$$
 
Notons que dans ce cas: 
$$b_k(s_k-s_{k-1},...,s_2-s_1,s_1)=\prod_{1\≤\ell\≤k}(s_\ell+\frac{\ell+1 }{
2}-1)(s_\ell+ \frac{\ell +r}{2}-1),$$
et que lorsque la caract\éristique r\ésiduelle de $\F$ est diff\érente de $2:$
$$(1-q^{-2(s_\ell+\frac{\ell+1 }{
2})})(1-q^{-2(s_\ell+\frac{\ell+r }{
2})})A^{a_0}_{ s_\ell+\frac{1 }{
2}(\ell+1),s_\ell+ \frac{1 }{2}( r+\ell)}(v_\ell,u_\ell,(-1)^{[\frac{ k-\delta}{2}]}w_\ell\delta_0)$$
est un polynome en $s_\ell$ dont les valeurs explicites sont donn\ées dans le lemme 3.6.7.\\

\noi Classiquement, les coefficients de l'\eq associ\ée \à $F_k$ s'obtiennent \à partir des coefficients $\widetilde a_{v,u}(s,...,s)$.\\

\noi Lorsque $2\≤2k\≤n,$ soit $H=\{s\in \F^*|sQ_a\sim Q_a\}$ ($H=\F^{*2}$ lorsque $\delta=1$) et    pour $u_0\in \FF$ et $\epsilon=\pm 1,$ soit $ S(u_0,\epsilon)=\{u\in (\FF)^k\ |$ il existe $t\in H$ v\érifiant $p(tu)=u_0$ et $h(tu)=\epsilon\}$ alors chaque $\cup_{u\in  S(u_0,\epsilon)}O_u$ est inclus dans une seule orbite de $G$ dans $\g'_1$ (cf.\cite{mullerJA2}). Comme dans le cas commutatif ($\S 6$), on n'effectue pas le calcul des sommes correspondantes qui sont peu maniables, on se contente de donner les sommes de coefficients : $\sum_{o^*}a_{O^*,O}$ lorsque la caract\éristique r\ésiduelle de $\F$ est diff\érente de $2.$
\\

\noi Lorsque $r\≥3,$ $O_u$ et $O^*_u$ sont non vides pour $u\in (\F^*/\F^{*2})^k$ (lemme 7.3.8), d\éterminons alors pour $a\in \FF$ et $s\in \C:$
$$\begin{array}{ll}\widetilde b_u(\tilde\omega_a|\ |^s)&=(a,p(u))\sum_{v\in (\F^*/\F^{*2})^k}\widetilde a_{v,u}(s,...,s)(a,p(v))\\
&=\alpha (-1)^{ka_0}\ \gamma (Q_a)^k \biggl(\prod_{1\≤i\≤k}\ \alpha(u_i)\biggr)^{-k}(a,p(u))b'_u(\tilde\omega_a|\ |^s)\end{array}$$
lorsque la caract\éristique r\ésiduelle de $\F$ est diff\érente de $2.$\\

\noi 
On rappelle les notations pour $z=(z_1,...,z_k)\in (\FF)^k:$ $$p(z)_i=\begin{cases}1\text{ lorsque }i=0\\
z_1...z_i\text{ pour }i=1,...,k,\end{cases} \text{ et }h(z)_i=\begin{cases}1\text{ lorsque }i=1\\
\prod_{j=1}^{k-1}(z_{j+1},z_1...z_j)\text{ pour }i=1,...,k\end{cases}$$
ainsi que $p(z):=p(z)_k$ et $h(z):=h(z)_k.$

\bigskip

\begin{prop} Lorsque la caract\éristique r\ésiduelle de $\F$ est diff\érente de $2,$ $k\≥2$ et $r\≥3$ on a:
$$\widetilde b_u(\tilde\omega_a|\ |^s)=  f_{[\frac{k+1}{2}]}(s)f_{[\frac{k+\delta+1}{2}]}(s+\frac{1}{2}r-\frac{1}{1+a_0})\rho (\tilde\omega_a|\ |^{s+1})\gamma (Q_a)^kB_u(\tilde\omega_a|\ |^s)$$
avec:
\begin{enumerate}
\item Lorsque r est impair (i.e. $a_0=0$):\\

\noi i)   Lorsque k est  pair  (donc  $\delta=1$ et $\delta_0=1$): 
 $$B_u(\tilde\omega_a|\ |^s)=h(u)\ \alpha(p(u))\ (p(u),-1) h(|\ |^{s+\frac{1}{2}(k+1)}\tilde\omega_{a(-1)^\frac{k}{2} p(u)}).$$
 ii) Lorsque k est impair (donc  $\delta=0$):
$$B_u(\tilde\omega_a|\ |^s)=(p(u),\delta_0) \rho(|\ |^{s+\frac{1}{2}(r+ k)}\tilde\omega_{  a\delta_0}).$$
 \item Lorsque r est pair (i.e. $a_0=1$):\\
 
i)   Lorsque k est  pair  (donc  $\delta=0$) :  $B_u(\tilde\omega_a|\ |^s)=    ((-1)^\frac{k}{2},a \delta_0)(a,p(u)) .$
$$\sum_{x\in\FF} 
 \alpha ( x)\   (x,a\delta_0) h(|\ |^{s+\frac{1}{2}(k+1)}\tilde\omega_{  ax} )\  \ A_{s+\frac{1}{2}(r+ 1),s+\frac{1}{2}(r+k)}^1(  x, (-1)^\frac{k}{2} p(u),1).$$
 
ii) Lorsque k est impair (donc  $\delta=1$ et $\delta_0=1$)
$$ B_u(\tilde\omega_a|\ |^s)= h(u)\alpha (p(u))  \ h(|\ |^{s+\frac{1}{2}(r+1)}\tilde\omega_{ a(-1)^\frac{k-1}{2}p(u)}).$$

\end{enumerate}
\bigskip

\end{prop}

\dem On remplace chaque $A^{a_0}_{ s_\ell+\frac{1 }{
2}(\ell+1),s_\ell+ \frac{1 }{2}( r+\ell)}(v_\ell,u_\ell, w_\ell\delta'_0)$ par sa valeur d'o\ù:
$$\prod_{\ell=1}^k A^{a_0}_{ s_\ell+\frac{1 }{
2}(\ell+1),s_\ell+ \frac{1 }{2}( r+\ell)}(v_\ell,u_\ell, w_\ell\delta'_0)= \prod_{\ell=1}^k (\sum_{t_\ell\in \FF}\alpha (t_\ell)^{a_0}\rho(s+\frac{1 }{2}( \ell+1);t_\ell v_\ell)\rho(s+\frac{1 }{2}( r+\ell);t_\ell u_\ell)(w_\ell \delta'_0,t_\ell)\ )$$
$$=\sum_{t=(t_1,...,t_k)\in (\FF)^k}\biggl(\prod_{\ell=1}^k f_\ell(t_\ell,u)\biggr) 
 \biggl(\prod_{\ell=1}^k \rho(s+\frac{1 }{2}( \ell+1);t_\ell v_\ell)\biggr) \biggl(\prod_{\ell=2}^k(v_\ell, p(t)_{\ell-1})\biggr) $$
avec
$$f_\ell(t_\ell,u)= \alpha (t_\ell) ^{a_0}  \rho(s+\frac{1 }{2}( r+\ell);t_\ell u_\ell) (t_\ell, p(u)_{\ell-1} \delta'_0\ ) $$
donc
$$b'_u(\tilde\omega_a|\ |^s)=\biggl(\prod_{1\≤i\≤k}\ \alpha(u_i)^i\biggr)\sum_{t=(t_1,...,t_k)\in (\F^*/\F^{*2})^k}\biggl(\prod_{\ell=1}^k f_\ell(t_\ell,u)\biggr) 
C(t,s)$$
avec:
$$\begin{array}{ll}C(t,s)&=\sum_{v\in (\F^*/\F^{*2})^k} \biggl(\prod_{\ell=1}^k \rho(s+\frac{1 }{2}( \ell+1);t_\ell v_\ell)\biggr) (\prod_{\ell=2}^k\alpha(-v_\ell)^{\ell-1}(v_\ell, p(t)_{\ell-1})\biggr)(a,p(v))\\
\\
&=\rho (\tilde\omega_a|\ |^{s+1})(t_1,a) \prod_{\ell=2}^kC_\ell(t)\text{ avec } \\
\\
 C_\ell(t)&=\sum_{v_\ell\in \FF}\alpha(v_\ell)^{1-\ell}(v_\ell, ap(t)_{\ell-1})\rho(s+\frac{1 }{2}( \ell+1);t_\ell v_\ell).\end{array}$$
 Le calcul donne (cf.3.6.2 et 3.6.4,A):
$$C_\ell(t)= \alpha(1)^{1-\ell} (t_\ell,a t'_{\ell-1}  )\begin{cases} 
 \rho(|\ |^{s+\frac{1 }{2}( \ell+1)}\tilde\omega_{ at'_{\ell-1} })\text{ lorsque }\ell\text{ est impair},\\
 \\
 \alpha(-t_\ell)h(|\ |^{s+\frac{1 }{2}( \ell+1)}\tilde\omega_{ at'_{\ell}})\text{ lorsque }\ell\text{ est pair}
\end{cases}$$
 en posant $t'_{\ell-1}=(-1)^{[\frac{\ell-1}{2}]} p(t)_{\ell-1},$ d'o\ù en r\éorganisant les produits: $C(t,s)\alpha(1)^{\frac{k(k-1)}{2}}=$
 $$  \rho(\tilde\omega_a|\ |^{s+1})h(t)\prod_{\ell=1}^{[\frac{k}{2}]}\alpha(-t_{2\ell})\prod_{\ell=1}^k(t_\ell,a(-1)^{[\frac{\ell-1}{2}]} ).
\prod_{\ell=1}^{[\frac{k}{2}]}h(|\ |^{s+\frac{1 }{2}( 2\ell+1)}\tilde\omega_{a t'_{2\ell}})\prod_{\ell=1}^{[\frac{k-1}{2}]}\rho(|\ |^{s+\frac{1 }{2}( 2\ell+2)}\tilde\omega_{a t'_{2\ell} })$$
donc par le 1) du lemme 3.6.8  et la d\éfinition 6.2.3:
$$C(t,s)\alpha(1)^{\frac{k(k-1)}{2}}= 
   \rho(\tilde\omega_a|\ |^{s+1})f_{[\frac{k+1}{2}]}(s)\ h(t)\ \prod_{\ell=1}^{[\frac{k}{2}]}\alpha(-t_{2\ell})\prod_{\ell=1}^k(t_\ell,a(-1)^{[\frac{\ell-1}{2}]} ).g_k(t_1...t_k)$$ avec:
$$ g_k(x)=\begin{cases}1\text{ lorsque }k\text{ est impair}\\
h(|\ |^{s+\frac{1 }{2}( k+1)}\tilde\omega_{ a(-1)^{\frac{k}{2}} x })\text{ lorsque }k\text{ est pair}.\end{cases}$$

\noi Dor\énavant on remplace $\alpha (1) $ par sa valeur: $\alpha (1)=1$ ainsi:
$$b'_u(\tilde\omega_a|\ |^s)=\rho(\tilde\omega_a|\ |^{s+1})f_{[\frac{k+1}{2}]}(s) c_u(s)$$avec:
$$c_u(s)=\biggl(\prod_{1\≤i\≤k}\ \alpha(u_i)^i\biggr)\sum_{t=(t_1,...,t_k)\in (\FF)^k}q(t_1,...,t_k)  \biggl(\prod_{\ell=1}^k     \rho(s+\frac{1 }{2}( r+\ell);t_\ell u_\ell) (t_\ell, U_{\ell-1})\biggr) g_k(t_1...t_k)$$
en posant $U_{\ell}=a(-1)^{[\frac{\ell}{2}]} p(u)_{\ell} \delta'_0 $ pour $\ell=0,...,k.$

\noi  Rempla\çant $h(t)$ par sa valeur, on a:$$q(t_1,...,t_k)=h(t)\ \prod_{\ell=1}^{[\frac{k}{2}]}\alpha(-t_{2\ell})\prod_{\ell=1}^k   \alpha (t_\ell) ^{a_0} =\begin{cases} \prod_{\ell=1}^{\frac{[k+1}{2}]}\alpha (t_{2\ell-1})\alpha(-t_1...t_k)\text{ lorsque }a_0=0\\
\\
 \prod_{\ell=1}^{\frac{[k}{2}]}\alpha (-t_{2\ell})\alpha(t_1...t_k)\text{ lorsque }a_0=1.\end{cases}$$
 
\noi Lorsque $k$ est impair, $k=2p+1:$  \\

$\bullet$ lorsque $a_0=0$ on a $ \alpha (t_{2p+1})\alpha(-t_1...t_k)= (t_{2p+1},t_1...t_{2p})\alpha(-t_1...t_{2p}),$ on somme d'abord sur $t_{2p+1}$  ce qui donne apr\ès simplification:
$$c_u(s)=\biggl(\prod_{1\≤i\≤2p}\ \alpha(u_i)^i\biggr)\sum_{t=(t_1,...,t_{2p})\in (\FF)^k} \prod_{\ell=1}^p\alpha (t_{2\ell-1})    \biggl(\prod_{\ell=1}^{2p}     \rho(s+\frac{1 }{2}( r+\ell);t_\ell u_\ell) (t_\ell, U_{\ell-1})\biggr) g(t_1...t_{2p})$$
avec: $g(x)=\alpha (-xu_{2p+1}) (u_{2p+1},(-1)^{p+1}U_{2p})\rho (|\ |^{s+\frac{1}{2}(r+2p+1)}\tilde\omega_{xU_{2p}});$\\

$\bullet$ lorsque $a_0=1,$ on commence par sommer sur $t_1$ d'o\ù:
$$c_u(s)=\biggl(\prod_{2\≤i\≤2p+1}\ \alpha(u_i)^i\biggr)\sum_{t=(t_2,...,t_{2p+1})\in (\FF)^k} \prod_{\ell=1}^p\alpha (-t_{2\ell})    \biggl(\prod_{\ell=2}^{2p+1}     \rho(s+\frac{1 }{2}( r+\ell);t_\ell u_\ell) (t_\ell, U_{\ell-1})\biggr) g(t_2...t_{2p+1})$$
avec: $g(x)=\alpha (u_1)\alpha (u_1x)(u_1, a\delta'_0)  h(|\ |^{s+\frac{1}{2}(r+ 1)}\tilde\omega_{xu_1a\delta'_0}).$\\

\noi  Lorsque $k$ est pair, $k=2p,$ on a \à \évaluer les quantit\és:\\

$\bullet$ lorsque $a_0=0:$
$$c_u(s)=\biggl(\prod_{1\≤i\≤2p}\ \alpha(u_i)^i\biggr)\sum_{t=(t_1,...,t_{2p})\in (\FF)^k} \prod_{\ell=1}^p\alpha (t_{2\ell-1})    \biggl(\prod_{\ell=1}^{2p}     \rho(s+\frac{1 }{2}( r+\ell);t_\ell u_\ell) (t_\ell, U_{\ell-1})\biggr) g(t_1...t_{2p})$$
avec: $g(x)= \alpha (-x)h(|\ |^{s+\frac{1}{2}(k+1)}\tilde\omega_{a(-1)^\frac{k}{2}x}).$ 
 
$\bullet$ lorsque $a_0=1:$ 
$$c_u(s)=\biggl(\prod_{1\≤i\≤2p}\ \alpha(u_i)^i\biggr)\sum_{t=(t_1,...,t_{2p})\in (\FF)^k} \prod_{\ell=1}^p\alpha (-t_{2\ell})    \biggl(\prod_{\ell=1}^{2p}     \rho(s+\frac{1 }{2}( r+\ell);t_\ell u_\ell) (t_\ell, U_{\ell-1})\biggr) g(t_1...t_{2p})$$
avec: $g(x)= \alpha (x)h(|\ |^{s+\frac{1}{2}(k+1)}\tilde\omega_{a(-1)^\frac{k}{2}x})=\alpha((-1)^px)((-1)^p,x)h(|\ |^{s+\frac{1}{2}(k+1)}\tilde\omega_{a(-1)^\frac{k}{2}x}).$\\

\noi Lorsque $a_0=0$, on obtient en appliquant le lemme 7.3.13, 3.6.8,1)  et la d\éfinition 6.2.3:
$$c_u(s)=f_{p+1}(s+\frac{1}{2}r-1)\ h(u)_{2p}\ (p(u)_{2p},-a\delta'_0)\ g(p(u)_{2p})$$ 
 avec :
$$\delta'_0=\begin{cases}(-1)^{p-1}\text{ lorsque }k=2p\text{ et }a_0=0,\\
(-1)^p\delta_0  \text{ lorsque }k=2p+1\text{ et }a_0=0.\end{cases}$$
Lorsque $a_0=1$, on obtient de m\ême:\\

$\bullet\ k=2p+1$ donc $\delta=1$ et $\delta'_0=(-1)^p:$
$$c_u(s)=f_{p+1}(s+\frac{1}{2}(r-1))\ h(u)_{2p+1}\ ( u_2...u_{2p+1},a\delta'_0)\ \prod_{i=2}^{2p+1}\alpha (u_i)\ g( u_2...u_{2p+1})$$

$\bullet\ k=2p$ donc $\delta=0$ et $\delta'_0=(-1)^p\delta_0:$
$$c_u(s)=f_{p}(s+\frac{1}{2}(r-1))\ h(u) \ (  p(u),a\delta'_0)\ \prod_{i=1}^{2p}\alpha (u_i)S(u) $$
avec:
$$S(u)=\sum_{x\in\FF}g((-1)^p  x)S_1(x)\text{ et}$$
$$\begin{array}{ll}S_1(x)&=\frac{1}{4}\sum_{z\in\FF}(a\delta'_0z , (-1)^pp(u)x)\  \rho (|\ |^{s+\frac{1}{2}(r+1)}\tilde\omega_{z })\ h(|\ |^{s+\frac{1}{2}(r+2p)}\tilde\omega_{z   (-1)^pp(u) })\\
\\
&=  ( a\delta'_0,(-1)^pp(u)x)\ \alpha (-(-1)^pp(u))A^1_{s+\frac{1}{2}(r+1),s+\frac{1}{2}(r+2p)}( x,  (-1)^pp(u),1)\end{array}$$
(lemme 3.6.4-B-6)) 
d'o\ù   en rempla\çant $g$ par sa valeur et en posant $$c_u(s)= f_{[\frac{k+\delta+1}{2}]}(s+\frac{1}{2}r-\frac{1}{1+a_0})c'_u(s),$$ on a: \\

\noi $\bullet$  Lorsque r est impair (i.e. $a_0=0$):
 $$c'_u( s)= h(u)\ \alpha (p(u))\ (p(u), a \delta'_0)\ .$$
$$\begin{cases}   
  h(|\ |^{s+\frac{1}{2}(k+1)}\tilde\omega_{a(-1)^\frac{k}{2} p(u)})\text{ lorsque }
 k \text{ est pair (donc }\delta=1 \text{  et }\delta_0=1)\\
 \\
   \rho(|\ |^{s+\frac{1}{2}(r+ k)}\tilde\omega_{  a\delta_0})\text{  lorsque } 
 k \text{ est impair (donc }\delta=0)\end{cases}$$
 $\bullet$ Lorsque r est pair (i.e. $a_0=1$):\\
 
i)   Lorsque $k$ est  pair  (donc  $\delta=0$) : $ c'_u( s)=  ( (-1)^\frac{k}{2},  a\delta_0p(u)). $
$$\sum_{x\in\FF} 
 \alpha ( x)\ (x,a\delta_0)   h(|\ |^{s+\frac{1}{2}(k+1)}\tilde\omega_{  a  x} )\  \ A_{s+\frac{1}{2}(r+ 1),s+\frac{1}{2}(r+k)}^1( x, (-1)^\frac{k}{2} p(u),1).$$
 
ii) Lorsque $k$ est impair (donc  $\delta=1$ et $\delta_0=1$)
$$ c'_u( s)=   (p(u),a(-1)^\frac{k+1}{2})\ h(|\ |^{s+\frac{1}{2}(r+1)}\tilde\omega_{ a(-1)^\frac{k-1}{2}p(u)}),$$
et on remplace $\prod_{1\≤i\≤k}\alpha (u_i)^{-k}=\begin{cases}(p(u),(-1)^\frac{k}{2})\text{ lorsque }k\text{ est pair,}\\
h(u)\alpha(p(u))(p(u),(-1)^\frac{k+1}{2})\text{ lorsque }k\text{ est impair.}\end{cases}$
 
\fdem

\bigskip
\noi {\bf Remarque:} Les calculs effectu\és dans cette proposition sont \également vrais sur $\R$ en conservant les facteurs $\alpha (1).$

\bigskip

\begin{lem}
Soient $p\≥1,$ $F:\FF\rightarrow \C$ et pour $i=1,...,p$ soient $c_i,u_i, v_i,U_i, V_i \in 
 \F^*,$ $s_i,s'_i\in \C $ et 
$$ S=\sum_{t_1,...,t_{p},y_1,...,y_p\in \FF} F(t_1...t_{p}y_1...y_p)\prod_{1\≤i\≤p}\alpha (c_iy_{i})(y_i,V_i)  (t_i,U_i)\rho (s_i;t_iu_i)\rho (s'_i;y_iv_i),$$
alors
$$ S=K\ \frac{1}{ f}\sum_{z,x\in \FF}\ (z,x)\  F(p(u)p(v)x)  \ \prod_{1\≤i\≤p}\ [\ \rho (|\ |^{s_i}\tilde\omega_{z U_i})\ h(|\ |^{s'_{i}}\tilde\omega_{zc_i v_{i}V_{i}})\ ]   $$avec
$$K=\alpha(1)^{-p}\ \prod_{1\≤i\≤p}\  \alpha(c_iv_i) \ \prod_{1\≤i\≤p}(u_i,U_i)\  \prod_{1\≤i\≤p}(v_i,V_i)$$
$ p(u)=\prod_{1\≤i\≤p}u_i,p(v)=\prod_{1\≤i\≤p}v_i,$ $f$ \étant le cardinal de $\FF .$

\end{lem}

\dem On change la sommation en posant pour $i=1,...,p:$ $x_i=t_1...t_{i}y_1...y_i$  donc $t_i=x_{i-1}x_iy_i$ d'o\ù:
$$A=\sum_{x_1,..x_p,y_1,...,y_p\in \FF} f(x_1,...,x_p) \prod_{1\≤i\≤p}\ [\alpha (c_iy_i)\rho (s_i;x_{i-1}x_iy_i u_i)
 \rho (s'_i;y_iv_{i})\ (y_i,U_iV_{i})\ ] $$
 avec 
 $$\begin{array}{lll}f(x_1,...,x_p)&=& F(x_p) \ \prod_{1\≤i\≤p}(x_{i-1}x_i, U_i)\quad \text{et}\quad x_0=1,\\
 \\
 &=& F(x_p) \ \prod_{1\≤i\≤p}(x_i, U_iU_{i+1})\quad \text{avec}\quad U_{p+1}=1\end{array}$$
donc:
$$A=\sum_{x_1,..x_p\in \FF} f(x_1,...,x_p)g(x)$$ avec:
$$\begin{array}{ll}g(x)&=\sum_{ y_1,...,y_p\in \FF} \prod_{1\≤i\≤p}\ \alpha (c_iy_i)\rho (s_i;x_{i-1}x_iy_i u_i)
 \rho (s'_{i};y_iv_{i})\ (y_i,U_iV_{i})\\
 \\
 &=\prod_{1\≤i\≤p}\ g_i(x)\end{array}$$
 avec:
 $$\begin{array}{lll}g_i(x)&=&\sum_{ y_i\in \FF}\alpha (c_iy_i)\rho (s_i;x_{i-1}x_iy_iu_i)
 \rho (s'_{i};y_iv_{i})\ (y_i,U_iV_{i})\   \\
 \\
 &=&(c_i,U_iV_i)A^1_{s_i,s'_i}(x_{i-1}x_ic_iu_i,c_iv_{i},U_iV_{i} ) \\
 \\
  &=&( v_{i},U_{i}V_{i})\ds\frac{\alpha (c_i v_{i})}{\alpha (1)}\ds\frac{1}{f}\sum_{z_i\in \FF}(z_i,x_{i-1}x_iu_{i}v_{i})\rho(|\ |^{s_{i}}\tilde\omega_{z_i})h(|\ |^{s'_{i}}\tilde\omega_{z_ic_iv_{i}U_iV_{i}})\end{array}$$
(B-6) du lemme 3.6.4)  donc:
  $$g(x)= K(\ds\frac{1}{f})^p\sum_{z_1,..z_p\in \FF} H(z_1,...,z_p)\prod_{1\≤i\≤p}\ (x_i,z_iz_{i+1})\quad \text{et}\quad z_{p+1}=1\ ,$$
 avec $$K=\prod_{1\≤i\≤p}[\  ( v_{i},U_{i}V_{i})\ds\frac{\alpha (c_i v_i)}{\alpha (1)}\ ]$$
 $$H(z_1,...,z_p)=\prod_{1\≤i\≤p}\rho(|\ |^{s_i}\tilde\omega_{z_i})h(|\ |^{s'_i}\tilde\omega_{z_ic_iv_iU_iV_i})(z_i,u_{i}v_{i})$$
 donc
 $$A= K\sum_{z_1,..z_p\in \FF} H(z_1,...,z_p)(\ds\frac{1}{f})^p\sum_{x_1,..x_p\in F} F(x_p) (x_p,z_pU_p)\ \prod_{1\≤i\≤p-1}(x_i,z_iz_{i+1}U_iU_{i+1})$$
   
$$= K\sum_{z_1,..z_p\in \FF} H(z_1,...,z_p)\prod_{1\≤i\≤p-1}[\ \ds\frac{1}{f}\sum_{ x_i\in F}(x_i,z_iz_{i+1}U_iU_{i+1})\ ](\ds\frac{1}{f}\sum_{ x_p\in \FF} F(x_p) (x_p,z_pU_{p})\ )$$
donc pour $i=1,...,p$ on a: $z_i=z_1 U_1U_i,$ on pose alors $z=z_1U_1$ d'o\ù $z_i=zU_i$ pour $i=1,...,p$ et on pose $ x=p(u)p(v)x_p.$  \fdem

 \bigskip

 \begin{rema}

 \begin{enumerate} 
 
 \item La proposition 7.3.12 est \également v\érifi\ée pour $k=1$ (cf.th.3.6.5-2) et 3), ici $N=\ds\frac{r+1}{2}$).\\
 
 \item Dans le cas $DI(n,k)$ (i.e. $\delta=0$)  avec $2n-3k\≥3$,les coefficients $B_u(\tilde\omega_a|\ |^s)$ ne d\épendent que de $p(u).$\\
 
 $\bullet$  Lorsque $k$ est impair, on a une \équation globale donn\ée par:$$Z^*(\four f;\tilde\omega_a|\ |^s)=B(\tilde\omega_a|\ |^s)Z(f;\tilde\omega_{a\delta_0}|\ |^{-s-n+k})\quad \text{avec}$$
 $$B(\tilde\omega_a|\ |^s)= \gamma (Q_a)^k f_{ \frac{k+1}{2}}(s) f_{\frac{k+1}{2}}(s+n-\frac{3}{2}k-1)\ \rho (\tilde\omega_a|\ |^{s+1})\rho(|\ |^{s+  n-k}\tilde\omega_{  a\delta_0}).$$
 
$\bullet$ Lorsque $k$ est  pair, on rappelle que pour $u\in \FF$  $Z_u(f;\pi)=Z(f1_{\{x\in \g'_1\ |F(x) \in u\FF\}};\pi)$ (resp.$Z^*_u(g;\pi)=Z(g1_{\{y\in \g'_{-1}\ |F^*(x) \in u\FF\}};\pi)$) alors classiquement (cf.$\S 6$ prop.6.2.4 et sa d\émonstration) on a pour tout $v\in \FF:$
$$Z_v^*(\four f;  s)=\gamma (Q_a)^k f_{ \frac{k}{2}}(s) f_{\frac{k}{2}}(s+n-\frac{3}{2}k-1)\sum_{u\in \FF}A_{v,u}( s)Z_u(f;  -s-n+k)$$
 avec: $A_{v,u}(s)=$
$$\begin{array}{lll} &=&\sum_{x\in \FF}A^1_{s+1,s+\frac{1}{2}(k+1)}((-1)^\frac{k}{2}v,x,1)((-1)^\frac{k}{2}x,\delta_0)A^1_{s+\frac{1}{2}(r+1),s+\frac{1}{2}(r+k)}(x,(-1)^\frac{k}{2}u,1),\\
&=& \sum_{x\in \FF}A^1_{s+1,s+\frac{1}{2}(k+1)}( v,x,(-1)^\frac{k}{2})( x,\delta_0)A^1_{s+\frac{1}{2}(r+1),s+\frac{1}{2}(r+k)}(x,u,(-1)^\frac{k}{2}).\end{array}$$
($r=2n-3k)$)\\

\noi Si on d\ésigne, pour $s\in \C,$ $A(s):={(A_{v,u}(s))}_{v,u\in \FF}$ la matrice "normalis\ée" des coefficients et
 $A_s(x,m):=(A^{a_{\frac{m}{2}}}_{s+1,s+\frac{m}{2}}(v,u,(-1)^{[\frac{m}{2}]}x))_{v,u\in \FF}$ 
la matrice "normalis\ée" des coefficients associ\és \à l'\eq pour une forme quadratique de discriminant $x$ sur un espace vectoriel de dimension $m,$ on a:
$$A(s)=A_s(1,k+1)DA_{s+\frac{1}{2}(r-1)}(1,k+1),$$
$D$ \étant la matrice diagonale \à $4$ lignes et $4$ colonnes de coefficients diagonaux $(x,\delta_0)_{x\in \FF}.$\\

\item Dans le cas $BI(n,k)$ (i.e. $\delta=1$)  avec $2n-3k\≥3$  et $k\≥3,$ les orbites de $G$ dans $\g'_1$ sont index\ées par $I=\{(u,\epsilon)|u\in \FF$ et $\epsilon=\pm 1\}$, notons $Z_{u,\epsilon}=Z_{O_{u,\epsilon}}$ (resp. $Z^*_{(u,\epsilon)}$) alors $$Z^*(\four f;  s)= C(s)\rho(|\ |^{s+1})\gamma (Q_a)^k\sum_{(u,\epsilon)\in S} A_{(u,\epsilon)}( s)Z_{(u,\epsilon)}(f;  -s-n+k-1)\quad \text{avec}$$
$$ C(s)= f_{ [\frac{k+1}{2}]}(s) f_{[\frac{k+2}{2}]}(s+n-\frac{3}{2}(k-1)-\frac{1}{1+a_0}) \ ,$$
$$A_{(u,\epsilon)}(s)=\epsilon\alpha (u) (u,(-1)^{k-1})h(|\ |^{s+\frac{1}{2}(b_0+1)}\tilde\omega_{(-1)^{[\frac{k}{2}]}u})$$
et $b_0=\begin{cases}k\text{ lorsque }k\text{ est pair,}\\
r\text{ lorsque }k\text{ est impair.}\end{cases}$\\

\noi Il est facile d'\établir comme dans 2) que pour $v\in \FF$ on a:
$$Z_{(v, -1)}^*(\four f;  s)+Z_{(v, 1)}^*(\four f;  s)=\gamma (Q_a)^k C(s). $$
$$\sum_{(u,\epsilon)\in S} \epsilon\  ((-1)^{[\frac{k-1}{2}]},u)A^1_{s+1,s+\frac{1}{2}(b_0+1)}(v,u,(-1)^{\frac{k}{2}} )   .Z_{(u,\epsilon)}(f;  -s-n+k-1) .$$
 \end{enumerate} 

 \end{rema}
\bigskip

 \subsection{Le cas DIII}
 
 \bigskip
 
\noi On rappelle que    
que  $(\overline{\goth g}_0,\overline{\goth g}_1)$ est de
type $(D_n,\alpha_{2p})$ avec $n\≥4$ et $3p\≤n,$ que
$\goth g$ est de type DIII donc en
particulier on a rang($\goth g)=[\ds\frac{n}{2}]$ (not\é
$m$) et $(\g_0,\g_1)$ est de type $(R_m,\lambda_p)$ avec $R=C$ lorsque $n$ est pair et $BC$ lorsque $n$ est impair (cf.tableau 2).\\

On donne les r\ésultats pour   $n-3p\≥2$ lorsque $\F$
est un corps $\goth p$-adique et $n-3p\≥4$ dans
le cas r\éel, hypoth\èses techniques permettant la descente avec des \PVs de type  DIII (cf. d\émonstration du
lemme 7.4.3),  plus restrictives  que la condition
pour que
$2H_0$ soit
$1$-simple  donn\ée par $n-3p\≥0;$ on suppose
\également 
{\bf dans le cas $\goth p$-adique   
 que l'entier $\bf n-p$ est impair},
hypoth\èse technique simplifiant les r\ésultats, en effet dans le cas $n-p$ pair il apparait un caract\ère difficile \à \évaluer de mani\ère intrins\èque
(cf.remarque 7.4.2,1)). \\

\noi Pour $k=1,..,p,$  $E_4(H_k) =\g^{2\epsilon_k}$ est de dimension $1$ et on note $X_k$ et $Y_k$ des g\én\érateurs
de $ E_4(H_k)$ et $ E_{-4}(H_k)$ tels que
$(X_k,\ds\frac{1  }{2}H_k,Y_k)$ soient des
$sl_2-$triplets.\\

\subsubsection{Le cas $p=1$}

\noi On rappelle que l'invariant relatif fondamental  du \PV $(\g_0,\goth
g_1)$ est donn\é par:
$$F(x)=\frac{2B(adx^4(Y_1),Y_1) }{B(H_1,H_1)}$$
et que $\chi (G)\subset \F^{*2}$ (lemme 2.3,\cite{mullerJA2}).
\vskip 2mm
\noi Comme $F$ est de degr\é $4$ et que l'on souhaite utiliser
les r\ésultats connus sur les formes quadratiques (th\éor\ème
3.6.5), l'\étude du \PV: $(\goth g_0,\goth g_1)$ se fait \à
travers l'\étude du  \PV: $(P',\goth g_1),$ $P'$ \étant le
sous-groupe parabolique d\éfini par $H'_1=h_{\lambda}$ et 
$H'_2=h_{\mu},$ avec $\lambda=\epsilon_1-\epsilon_2$ et
$\mu=\epsilon_1+\epsilon_2$  (dim$(\goth g^{\lambda})=4),$ alors $H'_1+H'_2=H_1.$\\

\noi On note $F'_1$ et $F'_2=F$ les invariants relatifs
fondamentaux associ\és \à $P'.$\\

\begin{lem}
1)  $F(\goth
g_1)= \F-\F^{*2}.$  

2) Les orbites de G dans $\goth g'_1$ (resp.$\goth
g'_{-1})$ sont donn\ées par
$O_u=\{x\in  \goth g_1\ |\ F(x) \F^{*2}=u\}$
(resp.$O^*_u=\{x\in  \goth g_{-1}\ |\ F^*(x) 
\F^{*2}=u\}$) avec $u\in \FF-1.$

3) Les orbites de P' dans $\goth g"_1$ (resp.$\goth
g"_{-1})$ sont donn\ées par
$O'_u=\{x\in  \goth g_1\ |\ F_1(x)\not=0,F(x) 
 \F^{*2}=u\}$ (resp.$O^*_u=\{x\in  \goth g_{-1}\ |\ \
F^*_1(x)\not=0\ ,\  F^*(x) \F^{*2}=u\}$) avec
$u\in \FF-1.$

4) Soient u et v dans $\FF-1, $ 

\quad i) $\F=\R$
$$ a_{-1,-1}
( s_1,s_2)=A_{\frac{3 }{ 2}}(s_2)A_{\frac{3 }{
2}}(s_1+s_2+n-\frac{7}{2}) 
=
4(2\pi)^{-2s_1-4s_2-2n+2}.$$
$$\Gamma (s_2+1)\Gamma
(s_2+\frac{3 }{ 2})\Gamma
(s_1+s_2+n-2).   \Gamma
(s_1+s_2+n-\frac{5 }{2})\ (-1)^n\ \sin ( \pi s_2)\cos \pi  (s_1+s_2)  
.$$
\quad ii) $\F$ est un corps $\goth p$-adique, on suppose
que n est pair alors $a_{v,u}(\omega_1,\omega_2;s_1,s_2)=$
$$     \sum_{w\in \FF-1}  
 A^1_{\omega_2,s_2+1, s_2+\frac{3 }{
2}}(v,w,1)A^1_{\omega_1\omega_2,s_1+s_2+n- \frac{5}{ 2},s_1+s_2+n-
2}  (w,u,1).$$
\end{lem}

\dem a) Soit $x$ non nul dans $\goth g_1,$ si
$s(x)=\{\epsilon_1\},$ \à l'aide de $exp(ad(\goth
g^{-\epsilon_2}))$ on peut se ramener \à $s(x)\cap
\{\epsilon_1±\epsilon_j,j\≥2\}\not =
\emptyset$ puis que $\epsilon_1-\epsilon_2\in s(x)$ donc
finalement que $\epsilon_1-\epsilon_2\in s(x)\subset
\{\epsilon_1\pm\epsilon_2\}$ par le lemme 2.2.1.

\noi Soit $x\in \goth g'_1$ (resp.$\goth g"_1)$, il existe ainsi
$g\in G$ (resp.$g\in P'$) tel que
$x=x_{\lambda}+x_{\mu}$ avec $x_{\lambda}\in \goth
g^{\lambda}-\{0\} $ et
$x_{\mu}\in \goth
g^{\mu}-\{0\}$  qui commutent d'o\ù 2).\\

\noi Le \PV commutatif associ\é \à $\goth
g^{\lambda}$ est de rang $1$ donc son invariant relatif
fondamental est une forme quadratique anisotrope et le
diagramme de  Satake associ\é  est de type :
 
\hskip 10pt \hbox to 1,5 cm {\offinterlineskip \lower 2pt\hbox{$\bullet$} \hglue -3,2pt
{\vrule height  0,4pt depth 0pt width 0,5 cm}\lower 2pt\hbox{$ \circledcirc$}
\hglue -3,2pt
{\vrule height  0,4pt depth 0pt width 0,5 cm}\lower 2pt\hbox{$\bullet$}}\\

\noi $\{x\in g^{\lambda}\ |\ F_1(x)\not=0\}$ est une seule orbite pour
$G_{H'_1}$ ( cf.$\S\ 6$) ainsi on peut
supposer que $F_1(x_{\lambda})=1$ et on a
$F(x_{\lambda}+x_{\mu})=P_1(x_{\mu}).$ Le \PV associ\é \à
$\goth U(\goth s),$ avec $\goth s=\F x_{\lambda}+\F h_{\lambda}+\F x_{\lambda}^{-1},$ est \également
commutatif de rang $1$ et dim$(\goth U(\goth s)_1)=3$ donc 
l'invariant relatif fondamental est une forme quadratique
anisotrope d\éfinie sur un espace de dimension $3$ d'o\ù  
 $F(\goth g_1)= \F-u\F^{*2},$ $-u$ \étant le discriminant de la restriction de $P_1$ \à $\goth U(\goth s)_1$ (\cite{o'meara}) et le diagramme de  Satake
associ\é  \à $(\goth U(\goth s)_0,\goth U(\goth s)_1)$ est de type :
\bigskip

\hskip 10pt \hbox to 1,2 cm { \offinterlineskip  {$\circledcirc$} \hglue -9pt\vbox{{\hrule height 0,3pt width 0,6cm}\vskip 3,5pt{\hrule height 0,3pt width 0,6cm}\vskip 0,3pt }\hglue -12pt $<$   \hglue -6pt $\bullet$}\\

\noi d'o\ù 3) (cf.$\S 6$ et 2).\\

b)   Montrons  que le discriminant de  $P_1 /\goth U(\goth s)_1$ vaut $-1$  d'o\ù $1).$

\noi Soit $E=\F[\sqrt \epsilon]$ une extension de dimension
$2$ de $\F$ sur laquelle $\goth g$ se d\éploie (\cite{veisfeiler}),
$( X_{\alpha})_{\alpha\in \overline ∆}$ une base de
Chevalley de $\goth g_E=\goth g\otimes_E\F$ et $\sigma$
l'antiautomorphisme de $\goth g_E$ ayant $\goth g$ comme
points fixes.\\

\noi On peut
supposer que :
$$\begin{array}{rl}x_{\lambda} &= 
X_{\overline\epsilon_2-\overline\epsilon_3}+
 \sigma (  X_{\overline\epsilon_2-\overline\epsilon_3}) 
=  X_{\overline\epsilon_2-\overline\epsilon_3}+p
   X_{\overline\epsilon_1-\overline\epsilon_4}\\
 x_{\mu}&=x\overline X_{{\overline \epsilon}_2+{\overline
\epsilon}_3}+
\overline{x}\sigma ( 
X_{\overline\epsilon_2+\overline\epsilon_3})+ y 
X_{\overline\epsilon_1+\overline\epsilon_3}
+\overline{y}\sigma ( 
X_{\overline\epsilon_1+\overline\epsilon_3})\\ &=x 
X_{\overline\epsilon_2+\overline\epsilon_3}+
\overline{x} q X_{\overline\epsilon_1+\overline\epsilon_4}+
y  X_{\overline\epsilon_1+\overline\epsilon_3}
+\overline{y} r 
X_{\overline\epsilon_2+\overline\epsilon_4}\ ,\ x,y\in
E.\end{array}$$  Prenons (ce qui ne change rien pour $F$) $Y_1=a 
X_{-\overline\epsilon_1-\overline\epsilon_2}\in \goth g_2,$
$a$ \étant un \él\ément non nul de $E$ v\érifiant la relation:
$$\frac{\overline a }{a}=pr\frac{N_1 }{N_2}\quad
\hbox{avec}\quad
N_1=N_{\overline\epsilon_1-\overline\epsilon_4,
\overline\epsilon_2+\overline\epsilon_4}\
,\ 
N_2=N_{\overline\epsilon_2-\overline\epsilon_3,
\overline\epsilon_1+\overline\epsilon_3}.$$
La relation $[x_{\lambda},x_{\mu}]=0$ se traduit par $ay\in 
\sqrt \epsilon\F.$

\noi En utilisant les propri\ét\és des coefficients
$N_{\alpha,\beta}$ (\cite{bourbakigal8}), chap.VIII,$\S 2,$n$^\circ4$, lemme 4 et
prop.7 ainsi que l'ex. 4), on obtient:
$$\begin{array}{rl}ad(x_{\lambda})^2(Y_1)&=-2apN_1N_3 
X_{-\overline\epsilon_3-\overline\epsilon_4}\\
ad(x_{\mu})^2(Y_1)&=2arN_4N_5\ (x\overline x\ \frac{q }{r}\ N\ +\ y\overline{y})\  
X_{\overline\epsilon_3+\overline\epsilon_4}\ \
\hbox{avec  } N=\pm 1,\\
N_3&=N_{\overline\epsilon_2-\overline\epsilon_3,-
\overline\epsilon_2-\overline\epsilon_4},\\
N_4&=N_{\overline\epsilon_1+\overline\epsilon_3,-
\overline\epsilon_1+\overline\epsilon_4},\\
N_5&=N_{\overline\epsilon_2+\overline\epsilon_4,-
\overline\epsilon_1-\overline\epsilon_2}\
,\ \hbox{d'o\ù}\\
F(x_{\lambda}+x_{\mu})&={a\over \overline a
}\ pr\ N_1N_3N_4N_5\ (x\overline x\ \frac{q\ }{
r}\ N\ a{\overline a}+\ ay\ \overline{ay})\\
&=N_2N_3N_4N_5\ (x\overline x\ \frac{q }{
r}\ N\ a{\overline a}+\ ay\ \overline{ay})\ ,\end{array}$$
le discriminant de cette forme quadratique vaut
$N_2N_3N_4N_5=-(N_2N_5)^2=-1$ car $N_i=±1, 2\≤i\≤5$ et on
applique l'ex.4 aux racines:
$$\alpha=\overline\epsilon_2-\overline\epsilon_3\ ,\
\beta=-\overline\epsilon_2-\overline\epsilon_4\ ,\
\gamma=\overline\epsilon_1+\overline\epsilon_3\ ,\
\delta=-\overline\epsilon_1+\overline\epsilon_4$$
ce qui donne $N_3N_4+N_{\gamma,\alpha}.N_{\beta,\delta}=0$
 or $N_{\gamma,\alpha}=-N_2$ et
$N_{\beta,\delta}=-N_5.$
\vskip 2mm

c) Il y a une seule classe de formes quadratiques anisotropes d\éfinies sur un espace
vectoriel de dimension 3  qui ne repr\ésentent pas
$1,$  elles ont pour discriminant $-1$ et pour invariant
de Hasse $-1$ (dans le cas r\éel et dans le cas $\goth
p$-adique), soit $G$ un repr\ésentant et $a_{v,u}^{(G)}$  le coefficient associ\é   aux orbites $O_u=\{x|G(x)\F^{*2}=u\}$ et 
 $O^*_v=\{x|G^*(x)\F^{*2}=v\}$ donn\é dans le th\éor\ème 3.6.5 (dans le cas r\éel et dans le cas $\goth
p$-adique):
 $$a_{v,u}^{(G)}(\omega,s)=-A^1_{\omega;s+1,s+\frac{3}{2}}(v,u,1).$$
 
 \noi Pour calculer les coefficients de l'\eq associ\ée au \PV $(P',\g_1),$  $a_{v,u}(\pi_1,\pi_2)$  avec $u,v\in \FF-1,$ on applique 
 la proposition 5.1.1 dont on reprend les notations.\\
 
 \noi Notons qu'ici
$\tilde B(\frac{H'_1 }{2},\frac{H'_1 }{2})=\tilde B(\frac{H'_2 }{2},\frac{H'_2 }{2})=-1.$\\

\noi Soient $O'_u$ et $O_v'^*$ deux orbites de $P'$ resp. dans $\g"_1$ et $\g"_{-1},$  $z=z_0+z_{-2}$ un
repr\ésentant dans
${{O_v'}^*}$ tel que ${{F'}_1}^*(z_0)=1$   et
${P'}_1^*(z_{-2})=v$ alors:\\

\noi $(E_{2,0})_{z_0}$ a pour \irf la forme quadratique anisotrope $F(...+z_0^{-1})/ (E_{2,0})_{z_0}$ qui ne  repr\ésente pas
$1,$ donc les coefficients correspondants sont donn\és par 
$a_{v,w}^{(
 G)}(\pi_2) $ avec $w\in \FF-1,$ soit $t_w\in (E_{2,0})_{z_0}$ tel que 
$F(t_w+z_0^{-1})=w$ alors 
$(E_{0,2})_{t_w}$ a pour \irf la forme quadratique anisotrope $F(t_w+... )/ (E_{0,2})_{t_w}$ qui ne  repr\ésente pas
$1,$ donc le coefficient correspondant est donn\é par 
$a_{w,u}^{(
 G)}(\pi_1\pi_2|\ |^{n-{7\over 2}})$ car:\\
 
   $\bullet \ p_{2,2}=1$ et 
 $ p_{1,1}=\hbox{dim}(\goth g_1)-8$ 
 
   $\bullet \ F(t_w+z_0^{-1})=w$ et $F(t_w+t_w(O'_u)=u,$\\
 
 \noi d'o\ù:
$$a_{v,u}(\pi_1,\pi_2)=
\sum_{w\in
\FF-1}\gamma(t_w,z_0^{-1})a_{v,w}^{(
 G)}(\pi_2) a_{w,u}^{   
(G)}(\pi_1.\pi_2|\ |^{n-{7\over 2}}) $$

d) Il reste \à calculer $\gamma(t_w,z_0^{-1})$ c'est \à dire \à
consid\érer la forme quadratique $Q(A)=\frac{1 }{ 2}{\tilde
B}([t_w,[z_0,A]],A)$ d\éfine sur
$E_{-1,1}=\oplus_{3\≤j\≤m}\goth
g^{\epsilon_2±\epsilon_j}\oplus E$ avec $E=\{0\}$ lorsque
$n$ est pair et $E=\goth g^{\epsilon_2}$ lorsque $n$ est
impair.
 
\noi Soit $A=\sum_{3\≤j\≤m}(A_j+B_j)+C$ avec $A_j\in \goth
g^{\epsilon_2-\epsilon_j},$ $B_j\in \goth
g^{\epsilon_2+\epsilon_j}$ et $C\in  E.$ 

\noi En raison des relations d'orthogonalit\é de $\tilde B$ on a:
$$Q(A)=\sum_{3\≤j\≤m} {\tilde
B}([t_w,[z_0,A_j]],B_j)+Q(C)\ .$$
Comme $ {\tilde
B}([t_w,[z_0,\goth
g^{\epsilon_2-\epsilon_j}]],\goth
g^{\epsilon_2+\epsilon_j})$ est une forme bilin\éaire non
d\ég\én\ér\ée, $\gamma (Q)=1$ lorsque $n$ est pair et sinon
$\gamma (Q)=\gamma (Q_1),$ $Q_1$ \étant la restriction de
$Q$ \à $E,$ espace vectoriel de dimension $4.$

\noi Soit $t\in \F^*,$ comme $tt_w+z_0^{-1}$ et $t_w+z_0^{-1}$
sont dans la m\ême $G_{H'_1}$ orbite par 3), $tQ_1$ et $Q_1$
sont \équivalentes. Ainsi $Q_1$ est de type $(2,2)$ dans le
cas r\éel et $\gamma (Q_1)=1;$ dans le cas $\goth p$-adique
$Q_1$ a pour discriminant $1$ et  $\gamma (Q_1)=1$ si $Q_1$
repr\ésente $0$ et $-1$ sinon.\fdem

\begin{rema}

\begin{enumerate}

\item Lorsque
$\F$ est un corps
$\goth p$-adique, on reprend les notations du b) de la
d\émonstration pr\éc\édente avec la description de $\goth
g_E$ donn\ée dans \cite{bourbakigal8} (chap.VIII,$\S\ 13,$ n$^\circ4$) et 
$\sigma$  donn\é par $\sigma (X)=T{\overline X}T^{-1}$  avec 

$$T=\left (\begin{array}{cccc}I_{\beta}&0&0... 0&0\\
0&I_{\beta}&0...0&0\\
.&...&...&.\\
0&0&...0&I_{\beta}\end{array}\right)\quad ,\quad
I_{\beta}=\left(\begin{array}{cc}0&\beta\\
1&0\end{array}\right )\quad ,\quad \beta\in \F^*\ \ \hbox{tel que}\ \
(\epsilon,\beta)=-1$$
alors $ a_{v,u}(\omega_1,\omega_2;s_1,s_2)  = $
$$\sum_{w\in \FF-1}\  
 \ (\epsilon, w)^n\  A^1_{\omega_2,s_2+1, s_2+\frac{3 }{2}}(v,w,1)A^1_{\omega_1\omega_2,s_1+s_2+n-\frac{5 }{  2},s_1+s_2+n- 
2}(w,u,1),$$ lorsque $n-1$ est impair,  la matrice des coefficients
est simplement le produit des 2 matrices des coefficients  des \irfs des centralisateurs:$$(a_{v,u}(\omega_1,\omega_2;s_1,s_2)_{v,u\in \FF-1}=A^{(G)}(\omega_2|\ |^{s_2}).A^{(G)}(\omega_1\omega_2|\ |^{s_1+s_2+n-\frac{7}{2}})$$ 
avec $A^{(G)}(\omega|\ |^{s})=(a^{(G)}_{v,u}(\omega|\ |^{s}))_{v,u\in \FF-1}.$\\

  \item Dans le cas r\éel, les polynomes de Bernstein associ\és au \PV $(P',\g_1)$ s'obtiennent \à partir de la prop.  3.4.4 et de 
 la  remarque 3.6.6,1) ce qui donne en raison de la
normalisation choisie:
$$b_1(s_1,s_2)=s_2(s_2+{1\over 2})\quad \hbox{et}\quad   
b_2(s_1,s_2)=s_2(s_2+{1\over 2})(s_1+s_2+n-{7\over 2})
(s_1+s_2+n-3).$$

\item Les \équations fonctionnelles associ\ées au \PV $(G,\g_{\pm 1})$ sont donn\ées par:\\

$\bullet$ Dans le cas r\éel:
$$Z^*(\four (f);s)=a^{(n)}(s)Z(f;-s-n+2)$$
avec
$$a^{(n)}(s) = 2(2\pi)^{ -4s-2n+2}\Gamma (s+1)\Gamma
(s+\frac{3 }{2})\Gamma
(s+n-2) 
 \Gamma
( s+n-\frac{5 }{2})(-1)^n\sin( 2\pi  s)   
.$$

$\bullet$ Dans le cas $\goth p$-adique , lorsque $n-1$ est impair, soit $v\in \FF-1$ alors:
$$Z_v^*(\four (f);s)=\sum_{u\in \FF-1} a^{(n)}_{v,u}(s)Z_u(f;-s-n+2)$$
avec
$$a^{(n)}_{v,u}(s ) = \sum_{w\in \FF-1} 
 A^1_{ s+1,s+\frac{3 }{ 2}}(v,w,1) A^1_{
 s +n- \frac{5 }{2},s+n-2}
   (w,u,1) .$$

\end{enumerate}
 \end{rema}

\subsubsection{Le cas $p\≥2$}

Pour $k=1,...,p$ soit: $$G_k(x)=\frac{2B(adx^4(Y_k),Y_k) }{
B(H_k,H_k)}\quad ,\quad x\in E_2(H_k)\cap \goth g_1\ .$$
Les invariants $F_1,F_2,...,F_p$ sont normalis\és par la
condition:
$$\hbox{pour}\quad x=\sum_{1\≤i\≤p}x_i\in W_{\goth t_0}\quad
F_k(x)=\prod_{1\≤i\≤k}G_i(x_i)\ ,$$
Il est alors ais\é de v\érifier que:
$$\hbox{pour}\quad y=\sum_{1\≤i\≤p}y_i\in W^*_{\goth t_0}\quad
F^*_k(y)=\prod_{p-k+1\≤i\≤p}G^*_i(y_i)\ .$$
On a $\chi_k(P_0)\subset \F^{*2}$  pour $k=1,...,p.$\\

  \begin{lem} On suppose que $n-3p\≥2$ lorsque $\F$
est un corps $\goth p$-adique et que $n-3p\≥4$ dans
le cas r\éel.

\noi 
 Soit $X\in W_{\goth t_0},$ $\goth s$ l'alg\èbre
engendr\ée par les projections de $X$ sur $E_2(H_i)\cap \goth
g_1$ et de $X^{-1}$ sur $E_{-2}(H_i)\cap \goth
g_{-1},$ $\goth U=\goth U(\goth s)$ alors $(\overline{\goth
U}_0,\overline{\goth
U}_1)$ est de type $(D_{n-3},\alpha_{2(p-1)})$ et $\goth U$
est de type DIII.
\end{lem}

  \dem
Il suffit de v\érifier que $\goth U$ est de type DIII (cf.\demo du 1) du lemme 7.1.3) avec $p=2$ puisque pour
$p\≥3$  on applique le lemme 7.1.3,3) et pour $i=1$ car $H_1$ et $H_2$ sont dans a m\ême orbite de $G.$

\noi Comme 
$\goth U(\F H_2)$ est de type DIII (cf.\demo du lemme 7.1.3) on a
$G_1(E_2(H_1)\cap \goth g_1)\subset
\F-\F^{*2}$ (lemme   7.4.1) donc $G_1({\goth
U}_1 )\subset
\F-\F^{*2}.$\\
 
 \noi Dans le cas $\goth p$-adique,  
on a $n\≥8 $ par hypoth\èse donc si $\goth U$  est de type DI alors $G_1(\goth U_1)=\F$ par le  lemme 7.3.8, 2)  ce qui est absurde donc $\goth U$  est de type DIII.\\

\noi Dans le cas r\éel, par la \demo du lemme   7.4.1, on peut supposer que $X\in \g^{\epsilon_1-\epsilon_2}\oplus \g^{\epsilon_1+\epsilon_2}$ donc $\oplus_{3\≤i\≤[\frac{n}{2}]}\g^{\epsilon_i }\subset \goth
U_0$ d'o\ù $rg(\goth U)\≥[\frac{n}{2}]-2$ or par hypoth\èse $n\≥10$ dans le cas r\éel d'o\ù $rg(\goth U)\≥3.$

\noi Supposons que  $\goth U$ soit de type DI et v\érifions alors que $G_1(\goth U_1)=\F$ ce qui est absurde .

\noi Il existe une constante $c$ telle que $cG_1/{\goth U_1}$ soit de type $(p_0,q_0)$ avec $q_0\≥1$ et $p_0\≥5$ (lemme 7.3.3) donc $O_{u_1,u_2}\not=\emptyset $ pour $0\≤u_1+u_2\≤2$ (lemme 7.3.8,1) d'o\ù $cG_1(\goth U_1)=\F.$
 \fdem
 \bigskip

\noi Soit  $u=(u_1,...,u_p)\in  (\FF-1)^p,$ on  rappelle que:\\

 $O_u=\{x\in \goth
g_1\ |\ F_1(x)\F^{*2} = u_1 \ ,\ F_2(x) \F^{*2} =u_1u_2\  \ ,...,F_p(x)\F^{*2} = u_1...u_p\  \}$ \\

$O^*_u=\{x\in
\goth g_{-1}\ |\ F^*_1(x) \F^{*2} =u_p\  \ ,\
F^*_2(x) \F^{*2} =
 u_{p-1}u_p\  \ ,...,F^*_p(x) \F^{*2} = u_p...u_1
\  \}$,\\

\noi  Les ouverts $O_u $
(resp.$O^*_u$) sont non vides et forment une partition de
$\goth g"_1$ (resp.$\goth g"_{-1}$).
 \vskip 3mm
 \begin{theo}
 On suppose que $ \goth g $ est de
type DIII et que $(\overline{\goth g}_0,\overline{\goth
g}_1)$ est de type $(D_n,\alpha_{2p}).$ 
 
  \begin{enumerate} 
\item Dans le cas r\éel, on suppose que $n-3p\≥4,$   soit $f\in
S(\goth g_1),$ alors:
$$Z^*(\hat f;s_1,...,s_p)=\bigl( \prod_{0\≤l\≤p-1} a^{(n-3p+3)}(
 s_{p-l}+...+ s_p+ l) \bigr)\  Z(f; 
 s^*  -(n-2p)1_p).$$
 \item Dans le cas $\goth p$-adique, on suppose que $n-3p\≥2$ et
que $n-p$ est impair.

 Soit $v\in (\FF-1)^p,$
  alors pour  $f\in
S(\goth g_1)$ on a:
$$Z^*_v(\hat f; \pi)= \sum_{u\in  (\FF-1)^p}a_{v,u}( s)Z_u(f; 
 s^*  -(n-2p)1_p),$$
avec $a_{v,u}(s)=\prod_{1\≤\ell\≤p}
a_{v_{\ell},u_{\ell}}^{(n-3p+3)}(
 s_{p-\ell+1}+...  +s_p +{\ell-1}).$
\end{enumerate}
\end{theo}
 \vskip 3mm
\dem On procède par récurrence sur $p,$
le cas $p=1$ ayant été fait dans le  lemme 7.4.1. On suppose le résultat vérifié lorsque  $(\overline{\goth
g}_0,\overline{\goth g}_1)$ est de type 
$(D_l,\alpha_{2(p-1)})$ avec $l-3(p-1)≥4$ dans le cas réel
(resp. $l-3(p-1)≥2$  et $l-(p-1)$ impair dans le cas $\goth
p$-adique)  avec 
$
\goth g
$  de type DIII.

Lorsque  $(\overline{\goth g}_0,\overline{\goth g}_1)$ est de
type $(D_n,\alpha_{2p}) ,$  avec $n-3p\≥4$ dans le cas réel
(resp. $n-3\≥2$  et $n-p$ impair dans le cas $\goth
p$-adique)  avec $ \goth g$  de
type DIII, on applique la proposition 5.3.1, dont on reprend
les notations, avec $k=1.$  

\noi  
Soient $z= z_0+z_{-2}\in O^*_v\cap W^*_{\goth t_0},$
$\goth U=\goth U(\goth s)$ (resp.$\goth U'=\goth U(\goth
s')$), $\goth s$ (resp.$\goth s')$ étant l'algèbre engendrée
par $ z_0$ et $ z_0^{-1}$ (resp.$z_2$
et $z_2^{-1}$).

\noi $(\overline {\goth U'}_0,\overline
{\goth U'}_1)$ est de type $(D_{n-3},\alpha_{2(p-1)})$ et 
$\goth U'$ est de type DIII par le lemme  7.4.3, 
$(\overline {\goth U}_0,\overline
{\goth U}_1)$ est de type $(D_{n-3p+3},\alpha_2 )$ et 
$\goth U=\cap_{2\≤i\≤p}\ \goth
U(\F y_i\oplus \F H_i\oplus\F y_i^{-1}),$  en décomposant  $z_0=\sum_{2\≤i\≤p}y_i$ et $y_i\in E_2(H_i)\cap \g_1,$ or $\goth
U(\F y_p\oplus \F H_p\oplus\F y_p^{-1})$ est de  
type DIII par le lemme 7.4.3 donc par récurrence
 $\goth U$ est également de type DIII.

\noi On applique la proposition 5.3.1 avec $H=\R^*$ dans le cas
réel  et $H=\F^{*2}$  dans le cas $\goth
p$-adique.\fdem
\bigskip

\noi{\bf Remarque} : Les sommes apparaissant dans le cas $\goth p-$adique semblent peu simplifiables aussi on ne donne aucune autre relation concernant les coefficients de l'\eq associ\ée \à $G$ ( cas $s_1=...=s_{p-1}=0$ et $s_p=s$).

         \newpage

\section{\bf Les cas exceptionnels}

\bigskip

Dans ce paragraphe, on d\étermine les polyn\ômes de Bernstein ainsi que les coefficients de l'\eq associ\ée aux \PVs $(\g_0,\g_1)$ pour lesquels $({\overline\g}_0,{\overline\g}_1)$ est de type exceptionnel, $P(H_1,H_2)$ \étant alors l'unique \sgp standard tr\ès sp\écial associ\é \à cette situation (cf.$\S 2.5$ et tableau 3).

\bigskip

\subsection{\bf Le cas $(E_7,\alpha_6)$}
\bigskip

\subsubsection{\bf G\én\éralit\és}

\bigskip

Les \PVs consid\ér\és dans cette section sont des $\F-$formes de $(E_7,\alpha_6)$. \\

Ils sont particuli\èrement simples \à \étudier parce que les \irfs des centralisateurs qui interviennent dans le calcul des coefficients de l'\eq
sont des formes quadratiques sur un espace de dimension paire ($8$) et dont les caract\éristiques (discriminant et coefficient $\gamma$) ne d\épendent que du rang de $\g$ (cf.prop.8.1.3).\\ \\

Par classification (\cite{warner},\cite{veisfeiler}), ils sont de trois types au plus de diagrammes de Satake  suivants:

\bigskip
1) Le cas d\éploy\é:\\

 \hskip 10pt \hbox to 3,5 cm {\lower 2pt\hbox{$\circ$}\hrulefill\lower 2pt
\hbox{$\circ$}\hrulefill \lower 2pt \vtop {\offinterlineskip \hbox{$\circ$} \hbox to 5pt{\hfill \vrule height 12pt width 0,3pt\hfill} \hbox{$\circ$}}\hrulefill\lower 2pt\hbox{$\circ$}\hrulefill\lower 2pt\hbox{$ \circledcirc$}\hrulefill\lower 2pt\hbox{$\circ$} }

\bigskip

2) Le cas EVI :

\bigskip

 \hskip 10pt \hbox to 3cm {\lower 2pt\hbox{$\circ$}\hrulefill\lower 2pt
\hbox{$ \circ$}\hrulefill \lower 2pt \vtop {\offinterlineskip \hbox{$ \circ$} \hbox to 5pt{\hfill \vrule height 12pt width 0,3pt\hfill} \hbox{$\bullet$}}\hrulefill\lower 2pt\hbox{$\bullet$}\hrulefill\lower 2pt\hbox{$\circledcirc$}\hrulefill\lower 2pt\hbox{$ \bullet$} }

\bigskip

Le diagramme de Dynkin du \PV $(\g_0,\g_1)$ est alors de type $(F_4,\lambda_4)$:
   \hskip 10pt \hbox to 2,5 cm {\offinterlineskip \lower 2pt\hbox{$\circ$} \hglue -5,2pt
{\vrule height  0,4pt depth 0pt width 0,6 cm}\hglue -0,8pt

    {\offinterlineskip\lower 2pt\hbox{ }\kern -3pt   \hrulefill\kern -3,4pt\lower 2pt \hbox{  \offinterlineskip  {$\circ$}
  \hglue -7pt\vbox{ {\hrule height 0,3pt width 0,6cm}\vskip 3,5pt{\hrule height 0,3pt width 0,6cm}
\vskip 0,3pt} \hglue -16pt $>$   \hglue -1,4pt{$\circ$}}}\hglue -0,8pt   {\vrule height  0,4pt depth 0pt width 0,6 cm}\hglue -0,8pt {\offinterlineskip\lower 2pt\hbox{$\circledcirc$}} }
\bigskip

3) Le cas EVII (alors $\F=\R$):

\bigskip

 \hskip 10pt \hbox to 3cm {\lower 2pt\hbox{$\circ$}\hrulefill\lower 2pt
\hbox{$\bullet$}\hrulefill \lower 2pt \vtop {\offinterlineskip \hbox{$\bullet$} \hbox to 5pt{\hfill \vrule height 12pt width 0,3pt\hfill} \hbox{$\bullet$}}\hrulefill\lower 2pt\hbox{$\bullet$}\hrulefill\lower 2pt\hbox{$\circledcirc$}\hrulefill\lower 2pt\hbox{$\circ$} }

\bigskip

 Le diagramme de Dynkin du \PV $(\g_0,\g_1)$ est alors de type $(C_3,\lambda_2)$: \hskip 10pt \hbox to 1,5cm{ \offinterlineskip \lower 2pt\hbox{$\circ$}\hrulefill\lower 2pt\hbox{\offinterlineskip {$\circledcirc$}\hglue -1,8pt\vbox{
 {\hrule height 0,3pt width 0,6cm}\vskip 3,5pt{\hrule height 0,3pt width 0,6cm}
\vskip 0,3pt} \hglue -16pt  $<$ \hglue -1,8pt  {$\circ$}}}

 \bigskip
Dans les $3$ cas, le \sg parabolique standard est donn\é par $P(H_1,H_2)$ avec $H_2=2h_{\widetilde\alpha}$ (cf.\demo du lemme 1.3.1) et $H_1=2(H_0-h{\widetilde\alpha})=2h_{\beta},$ $\widetilde\alpha$ \étant la plus grande racine de $\Delta$ et \\

$\bullet$ $\beta= 012221$ dans le cas d\éploy\é,
 
   \hskip 1,5cm   $1$  \\

$\bullet$ $\beta= 0122  $ dans le cas EVI (c'est la restriction de la racine pr\éc\édente),\\

$\bullet$ $\beta= 021  $ dans le cas EVII (c'est la restriction de la racine du cas d\éploy\é).\\

On note $(\∆,\lambda_0)$ le diagramme de Dynkin de $(\g_0,\g_1)$ et $\Sigma$ un ensemble de racines simples de $\∆:$   \\

$\bullet$ dans le cas d\éploy\é, $\Sigma=\{\alpha_i,1\≤i\≤7\}$ (notation de la planche VI de \cite{bourbakigal6}),\\

$\bullet$ dans le cas EVI, $\Sigma=\{\lambda_i,1\≤i\≤4\},$ $\lambda_i$ est la restriction de $\alpha_i,i=1,2,$  $\lambda_3$ est la restriction de $\alpha_4$ et $\lambda_4$ est la restriction de $\alpha_6,$\\
 
 $\bullet$ dans le cas EVII, $\Sigma=\{\lambda_i,1\≤i\≤3\},$ $\lambda_1$ est la restriction de $\alpha_1, $  $\lambda_2$ est la restriction de $\alpha_6$ et $\lambda_3$ est la restriction de $\alpha_7.$\\
 
$H_1$ et $H_2$ sont dans la m\ême orbite de $G.$\\

L'alg\èbre simple $  E=\goth U(\F H_2) $  (resp.$ F=\goth U(\F H_1)$) est gradu\ée par $\ds\frac{1}{2}H_1$ (resp.$\ds\frac{1}{2}H_2$) et  \\

$\bullet$ le \PV $(  E_0,  E_1)$ (resp.$(  F_0,  F_1)$) est de type
 $(D_6,\alpha_2)$ dans le cas d\éploy\é , \\

$\bullet$  $ E$ (resp.$F$) est de type DIII dans le cas EVI,\\

$\bullet$ $  E$ (resp.$F$) est de type DI dans le cas EVII,\\

$  E$  (resp.$F$)
admet comme sous-alg\èbre ab\élienne d\éploy\ée maximale $\goth h_0=\oplus_{\lambda\in \Sigma | n(\widetilde\alpha,\lambda)=0}\F h_\lambda$  (resp. $\goth h'_0=\oplus_{\lambda\in \Sigma | n( \beta,\lambda)=0}\F h_\lambda$), de plus:
$$\g_1=E_1\oplus F_1\quad \g_2=\g^{\widetilde \alpha}\oplus \g^\beta\oplus E_2(H_1)\cap E_2(H_2)$$
$\g^{\widetilde \alpha}$ et $\g^\beta$ sont de dimension $1,$    $\text{dim}(E_1)=\text{dim}(F_1)=2p_{2,2}=16.$  \\ \\
\subsubsection{\bf Pr\éliminaires}

\begin{lem} Soit $x\in E'_1$  (resp.$x\in F'_1$) et $\goth U_x=\goth U(\F x\oplus \F H_1\oplus \F x^{-1})$  (resp.$\goth U_x=\goth U(\F x\oplus \F H_2\oplus \F x^{-1})$) alors:

\begin{enumerate}
\item Dans le cas d\éploy\é,  $\text{rang}(\goth U_x) \≥3.$

\item Dans les cas EVI ou VII, $\text{rang}(\goth U_x) =\text{  rang}(\g)-2.$
\end{enumerate}

\end{lem}

\dem 1) On rappelle que rang$(\goth U_x)\≤$rang($\g)-2$ (proposition,appendice 1) ce qui termine le cas  EVII et dans le cas EVI on a l'in\égalit\é rang$(\goth U_x)\≤2.$  \\

2) Dans le cas d\éploy\é, on effectue le calcul de $(\goth U_x).$

\noi Comme $(E_0,E_1)$ est de type $(D_6,\alpha_3)$ on peut supposer que $x=\sum_{1\≤i\≤4}a_iX_{\mu_i},$ avec $$a_1a_2a_3a_4\not=0 \quad\text{et }$$
$$\mu_1=\alpha_6\ ,\ \mu_2=\alpha_5+\alpha_6+\alpha_7\ ,\
\mu_3=\alpha_6+2\alpha_5+2\alpha_4+\alpha_2+\alpha_3 \ ,\ \mu_4=\alpha_7+\alpha_6+\alpha_5+2\alpha_4+\alpha_2+\alpha_3$$
donc $x^{-1}=\sum_{1\≤i\≤4}a_i^{-1}X_{-\mu_i}$ (cf. par exemple prop.6.6 de \cite{mullerJA1} et  th.1 et tableau II de \cite{mullerJA2}).\\

\noi Soient :
$$\delta_1=\alpha_1+\beta\ ,\ \delta_2=\delta_1+\alpha_3+\alpha_4+\alpha_5\ ,$$
  $$X_1=a_2[X_{-\mu_1},X_{\delta_1}]-a_1[X_{-\mu_2},X_{\delta_1}]\ ,\ h_1=h_{\delta_1-\mu_1}+h_{\delta_1-\mu_2}\ ,\ Y_1=a_2^{-1}[X_{\mu_1},X_{-\delta_1}]-a_1^{-1}[X_{\mu_2},X_{-\delta_1}]\ ,$$
   $$X_2=a_2[X_{-\mu_3},X_{\delta_2}]-a_3[X_{-\mu_2},X_{\delta_2}]\ ,\ h_2=h_{\delta_2-\mu_2}+h_{\delta_2-\mu_3}\ ,\ Y_2=a_2^{-1}[X_{\mu_3},X_{-\delta_2}]-a_3^{-1}[X_{\mu_2},X_{-\delta_2}]\ ,$$
comme les racines  $\delta_1-\mu_1$ et $\delta_1-\mu_2$ (resp. $\delta_2-\mu_2$ et $\delta_2-\mu_3$) sont orthogonales, on v\érifie facilement que $(X_i,h_i,Y_i),i=1,2,$ sont $2$ \Sls de l'alg\èbre $\goth U_x$ donc $\F h_{\widetilde \alpha}\oplus \F h_1\oplus \F h_2=\F h_{\widetilde \alpha}\oplus \F h_{\alpha_2}\oplus \F h_{\alpha_3}$ est une sous-alg\èbre ab\élienne d\éploy\ée de $\goth U_x.$\\

3) Dans le cas EVI, dans une extension quadratique convenable de $\F,$ $\F'=\F\sqrt a,$ $\g'=\g\otimes_{\F}\F'$ est une alg\èbre d\éploy\ée, notons $\sigma$ la conjugaison ``associ\ée".\\

\noi  Comme $\sigma (\mu_1)=\mu_2$ et $\sigma (\mu_3)=\mu_4$,  $$\g^{\lambda_4}=\{xX_{\mu_1}+\overline xk_1X_{\mu_2}\ |\ x\in \F'\}\ ,\ \g^{\lambda_2+2\lambda_3+\lambda_4}=\{yX_{\mu_3}+\overline yk_3X_{\mu_4}\ |\ y\in \F'\}$$
avec $\sigma(X_{\mu_1})=k_1X_{\mu_2}$ et $\sigma(X_{\mu_3})=k_3X_{\mu_4}.$ 

\noi Notons \également que  $\sigma (\delta_1)=\delta_1$ donc $\sigma(X_{ \delta_1})=kX_{ \delta_1}$ avec $k\overline k=1.$\\

\noi Comme le \PV $(E_0,E_1)$ est de type $(C_3,\alpha_1),$ on peut supposer que $x\in E'_1 \cap (\g^{\lambda_4}\oplus\g^{\lambda_2+2\lambda_3+\lambda_4})$ donc s'\écrit sous la forme:
$$x=xX_{\mu_1}+\overline xk_1 X_{\mu_2}+yX_{\mu_3}+\overline yk_3X_{\mu_4}\ ,x,y\in \F'-\{0\}.$$
Soit $t\in \F'-\{0\}$ tel que $\ds\frac{t}{\overline t}=-\ds\frac{\overline k_1}{k_1}k$ et:
$$ X_1=t\overline xk_1[X_{-\mu_1},X_{\delta_1}]-xt[X_{-\mu_2},X_{\delta_1}]\ , \ Y_1=\ds\frac{1}{t\overline xk_1}[X_{\mu_1},X_{-\delta_1}]- \ds\frac{1}{tx}[X_{\mu_2},X_{-\delta_1}]\ ,   $$
alors $\sigma(X_1)=X_1,\sigma(Y_1)=Y_1$ et comme dans 2), $(X_1,h_1=h_{\delta_1-\mu_1}+h_{\delta_1-\mu_2},Y_1)$ est un \Sl-triplet qui commute \à $(x,H_1,x^{-1})$ donc  $\F h_{\widetilde \alpha}\oplus \F h_1 $ est une sous-alg\èbre ab\élienne d\éploy\ée de $\goth U_x$ de dimension $2$  or  $\goth U_x$ est de rang inf\érieur ou \égal \à  $2$ d'o\ù $\goth U_x$ est de rang $2.$\fdem\\

\begin{prop} Soit $x=x_2+x_0\in W_{\goth t}$ avec $x_i\in E_i(H_1)\cap \g_1,$ alors les \PVs $((\goth U_{x_i})_0,(\goth U_{x_i})_1),$ $i=0$ ou $2,$ sont des $\F$-formes de $(A_5,\{\alpha_1,\alpha_5\}).$ 

Dans les 2 cas l'\irf est une forme quadratique de discriminant 1 et d'invariant $\gamma=1$ sauf dans le cas EVI o\ù il vaut $-1.$

\end{prop}

\dem Les \PVs $((\goth U_{x_i})_0,(\goth U_{x_i})_1),$ $i=0$ ou $2,$ sont des $\F$-formes de $(A_5,\{\alpha_1,\alpha_5\})$ d'apr\ès la \demo du lemme 1.3.1 (2)). Mais alors le rang de $\goth U_{x_i}$  peut prendre ces valeurs dans $\{1,2,3,5\}$ par la  proposition  8.1.3 donc par le lemme 8.1.1 :\\

$\bullet$ le rang de $\goth U_{x_i}$ vaut $3$ ou $5$ dans le cas d\éploy\é d'o\ù l'invariant $\gamma(F(\ +x_i)/(\goth U_{x_i})_1))=1$ par la  proposition  8.1.3,\\

$\bullet$ le rang de $\goth U_{x_i}$ vaut $2$ dans le cas EVI donc  $\gamma(F(\ +x_i)/(\goth U_{x_i})_1))=-1$ par la  proposition 8.1.3,\\

$\bullet$ le rang de $\goth U_{x_i}$ vaut $1$ dans le cas EVII r\éel donc  $\gamma(F(\ +x_i)/(\goth U_{x_i})_1))=1$ par la  proposition 8.1.3.\fdem\\ \\

\begin{prop} Soit $(\g_0,\g_1)$ un \PV tel que $(\overline{\g_0},\overline{\g_1})$ soit de type $(A_{2n+1},\{\alpha_1,\alpha_{2n+1}\})$ et soit $p$ le rang de $\g$   alors $p=2n+1$ ou bien $1\≤p\≤n+1$ et l'\irf est une forme quadratique $F$ de discriminant 1 et d'invariant $\gamma(F)=\begin{cases} 1\text{ lorsque }p=2n+1,\\
(-1)^{n+1-p}\text{ sinon .}\end{cases}$

\end{prop}  

\dem On effectue le calcul explicite de $F$ dans chaque cas.\\

1) Le cas d\éploy\é\\

\noi Dans les notations de la planche I de\cite{bourbakigal6} et en prenant un syst\ème de Chevalley $(X_{\alpha})_{\alpha\in \Delta},$ il suffit de calculer $\widetilde B((adx)^2(X_{-\widetilde \alpha},(adx)^2(X_{-\widetilde \alpha})$ pour $x\in \g_1,$ $\widetilde \alpha$ \étant la plus grande racine de $\Delta.$ Notons que dans ce syst\ème de racines $N_{\alpha,\beta}\in \{0,\pm 1\}$ et rappelons que $\widetilde B(X_{\alpha},X_{-\alpha})=1.$\\

\noi Ecrivons $x=\sum_{i=2}^{2n+1}x_iX_{\epsilon_1-\epsilon_i}+y_i
X_{\epsilon_i-\epsilon_{2n+2}}$ alors:\\

\noi $adx(X_{-\widetilde \alpha})=\sum_{i=2}^{2n+1}x_iN_iX_{\epsilon_{2n+2}-\epsilon_i}+y_iN'_i
X_{\epsilon_i-\epsilon_{1}}$ avec $N_i=N_{\epsilon_1-\epsilon_i,-\widetilde \alpha}$ et $N'_i=N_{\epsilon_i-\epsilon_{2n+2},-\widetilde \alpha},$\\

\noi $(adx)^2(X_{-\widetilde \alpha})=\sum_{2\≤i\not=j\≤2n+1}x_iy_jA_{i,j}X_{\epsilon_{ j}-\epsilon_i}+ H$ avec:\\

\noi $A_{i,j}=N_{\epsilon_1-\epsilon_i, \epsilon_j-\epsilon_1}N'_j+N_iN_{\epsilon_j-\epsilon_{2n+2},\epsilon_{2n+2}-\epsilon_i}$ et $H=-(\ \sum_{i=2}^{2n+1}x_i y_i(N'_ih_{\epsilon_1-\epsilon_i}+
N_ih_{\epsilon_i-\epsilon_{2n+2}})\ ),$ \\

\noi or:  \\

\noi $N_{\epsilon_1-\epsilon_i, \epsilon_i-\epsilon_{2n+2}}=-N_{\epsilon_i-\epsilon_1, \widetilde \alpha}=-N_{\epsilon_1-\epsilon_i,-\widetilde \alpha}=-N_i$ et $N'_i=-N_{\widetilde \alpha,\epsilon_{2n+2}-\epsilon_i}=N_{-\widetilde \alpha,\epsilon_1-\epsilon_i}=-N_i$ (lemme $4(3),$ chap.VIII,$\S 2,n^\circ 4$ \cite{bourbakigal8}) donc\\

\noi $H= \ \sum_{i=2}^{2n+1}x_i y_i 
N_i( h_{\epsilon_1-\epsilon_i}-h_{\epsilon_i-\epsilon_{2n+2}}\ )$ et $A_{i,j}=-N_{\epsilon_1-\epsilon_i, \epsilon_j-\epsilon_1}N_j+N_iN_{\epsilon_j-\epsilon_{2n+2},\epsilon_{2n+2}-\epsilon_i} .$\\

\noi Pour $i\not=j,$ posons: $\alpha=\epsilon_1-\epsilon_i,\beta=\epsilon_j-\epsilon_1,\gamma= \epsilon_{2n+2}-\epsilon_j ,\delta=\epsilon_i-\epsilon_{2n+2},$ alors $\alpha+\beta+\gamma+\delta=0$ d'o\ù:\\

\noi $N_{\epsilon_1-\epsilon_i, \epsilon_j-\epsilon_1}N_{\epsilon_{2n+2}-\epsilon_j , \epsilon_i-\epsilon_{2n+2}}+N_{\epsilon_j-\epsilon_1,  \epsilon_{2n+2}-\epsilon_j}N_{\epsilon_1-\epsilon_i,\epsilon_i-\epsilon_{2n+2}}=0 $ (exercice $4$ p. $221$, chap.VIII,$\S 2$ \cite{bourbakigal8}) donc $N_{\epsilon_{2n+2}-\epsilon_j , \epsilon_i-\epsilon_{2n+2}}=-N_{\epsilon_1-\epsilon_i, \epsilon_j-\epsilon_1}N_jN_i.$\\

\noi Ainsi $A_{i,j}=-2N_j N_{\epsilon_1-\epsilon_i, \epsilon_j-\epsilon_1}$ et  
 $$ \widetilde B((adx)^2(X_{-\widetilde \alpha}),(adx)^2(X_{-\widetilde \alpha}))=\sum_{2\≤i\not=j\≤2n+1}x_iy_jx_jy_iA_{i,j}A_{j,i}+\widetilde B(H,H)$$
$$=-4\sum_{2\≤i\not=j\≤2n+1}x_iy_jx_jy_iN_iN_j+\sum_{i=2}^{2n+1}x_i y_iN_i((\epsilon_i-\epsilon_{2n+2})(H)-(\epsilon_1-\epsilon_i)(H)  \ )$$
$$=-6\bigl(\ \sum_{i=2}^{2n+1}x_i^2 y_i^2 +\sum_{2\≤i\not=j\≤2n+1}x_iy_jx_jy_iN_iN_j\bigr)=-6(\sum_{i=2}^{2n+1}x_iy_iN_i)^2\ ,$$
et finalement l'\irf est une forme quadratique $F(x)=c\sum_{i=2}^{2n+1}x_iy_iN_i,$ $c\in \F^*,$ par cons\équent $\gamma(F)=1$ et  $\text{discr}(F)=1.$\\

2) Le cas de rang $p$ avec $\g$ non d\éploy\éee.\\

Par classification, on ne peut avoir que le type AIII (tables de \cite{warner} et de \cite{veisfeiler}) par cons\équent $p\≤n+1$ et
sur une extension quadratique convenable $\F'=\F\sqrt a ,$ $\g$ est d\éploy\ée, notons $\sigma$ la "conjugaison associ\ée" alors pour $i=1,...,p-1$ on a $\sigma (\alpha_i)=\alpha_{2n+2-i},$ $\sigma (\alpha_i)=-\alpha_ i$ pour $i=p+1,...,2n+1-p$ lorsque $p\≤n$ et $\sigma (\alpha_p)=\sum_{p+1\≤j\≤2n+2-p}\alpha_j$ d'o\ù:\\

\noi  $\sigma (X_{ \epsilon_1-\epsilon_i})=k_i X_{ \epsilon_{2n+3-i}-\epsilon_{2n+2}}$ pour $i=1,...,p,2n+3-p,...,2n+1$ et $\sigma (X_{ \epsilon_1-\epsilon_i})=k_iX_{ \epsilon_i-\epsilon_{2n+2}}$ pour $i=p+1,...,$ $2n+2-p.$ \\

\noi Comme $\sigma (\widetilde \alpha)=\widetilde \alpha$ on a 
$\sigma (X_{\widetilde \alpha})=kX_{\widetilde \alpha}$ avec $k\overline k=1$ donc pour  $i=2,...,p:$ $$\sigma([X_{ \epsilon_1-\epsilon_i},X_{ \epsilon_i-\epsilon_{2n+2}}])=\ds\frac{k_i}{\overline {k_{2n+3-i}}}N_{2n+3-i}X_{\widetilde \alpha}=-N_ikX_{\widetilde \alpha}\Rightarrow  k_i=-\overline {k_{2n+3-i}}kN_{2n+3-i}N_i.$$
Soit $x\in\g_1$  d\écompos\é comme dans 1) et v\érifiant $\sigma (x)=x,$   alors   $y_{2n+3-i}=k_i\overline{x_i}$ pour $i=2,...,p,2n+3-p,...,2n+1$ et $y_i=k_i\overline{x_i} $ pour $i=p+1,...,$ $2n+2-p$ donc il existe $c\in \F'$ tel que $\ds\frac{\overline c}{c}=-k$ et:

$$\begin{array}{rl}F(x)&=c(\sum_{2\≤i\≤p,2n+3-p\≤i\≤2n+1}(\overline{x_i }{x_{2n+3-i}}k_iN_{2n+3-i}))+F_a(x')\\
\\
&=\sum_{i=2}^{p}(c\overline{x_i }{x_{2n+3-i}}k_iN_{2n+3-i}+c {x_i}\overline{x_{2n+3-i}}k_{2n+3-i}N_i)+F_a(x')\\
\\
 
&=\sum_{i=2}^{p}(x_i \overline{z_i}+z_i \overline{x_i})+F_a(x')\quad\text{avec }z_i= {c  k_i}N_{2n+3-i}x_{2n+3-i} \\
\\
F_a(x')&=0 \text{ lorsque }p=n+1 \text{ et  lorsque  }   p\≤n \text{ on a:}\\
\\
 F_a(x')&=\sum_{ p+1\≤i\≤2n+2-p}x_i\overline{x_i}ck_iN_i \text{ et }
x'=\sum_{ p+1\≤i\≤2n+2-p}x_iX_{\epsilon_i-\epsilon_{2n+1}}+\overline{x_i}\sigma(X_{\epsilon_i-\epsilon_{2n+1}}) , \end{array}$$

\noi donc $ \text{disc}(F)= 1$ et $\gamma (F)=1$  lorsque $p=n+1$  (pas de racines compactes) et sinon  $\gamma (F)=\gamma (F_a).$\\ 

\noi Soit $\goth b'$ la $\F'$-alg\èbre engendr\ée par $X_{\pm(\epsilon_1-\epsilon_{p+1})},$ $X_{\pm(\epsilon_{2n+1-p}-\epsilon_{2n+1})},$
$X_{\pm(\epsilon_{i}-\epsilon_{i+1})}$ pour $i=p+1,2n+1-p.$ Elle est gradu\ée par $h_{\tilde \alpha}$ qui appartient \à $\goth b'$ et le \PV $(\goth b'_0,\goth b'_1)$ est de type 
 $(A_{2(n+1-p)+1},\{\alpha_1,\alpha_{2(n+1-p)+1}\}).$
 
 \noi Elle est $\sigma$-stable et la $\F$-alg\èbre correspondante, $\goth b=\{x\in \goth b'|\sigma(x)=x\},$ est alors de $\F-$rang $1$ et $F_a$ est l'\irf du \PV $(\goth b_0,\goth b_1)$ donc $F_a$ est une forme quadratique anisotrope  sur un espace vectoriel de dimension $4(n+1-p)$ d'o\ù   $\gamma (F_a)=(-1)^{n+1-p}$ dans le cas r\éel et dans le cas $\goth p-$adique on a n\écessairement $p=n$ et $\gamma (F_a)=-1$ par classification des formes quadratiques.\fdem\\

\bigskip
\subsubsection{\bf Le r\ésultat}
 \bigskip

\begin{theo} Soient $(\pi_1,\pi_2)\in \widehat{(\F^*)^2},$ 
 $f\in \EuScript S(\g_1)$ alors:
 $$\begin{array}{rcl}Z^*(\four f;(\pi_1,\pi_2))&=&A(\pi_1,\pi_2) Z(f;(\pi_1,\pi_1^{-2}\pi_2^{-1}|\ |^{-8})\quad \text{avec} \\
  A(\pi_1,\pi_2)&=&\rho (\pi_2|\ |)\rho (\pi_2|\ |^4)\rho (\pi_1^2\pi_2|\ |^5)\rho (\pi_1^2\pi_2|\ |^8).\end{array}$$

\end{theo}

 \dem On reprend le calcul de la proposition 5.1.1 avec les m\êmes hypoth\èses et notations.\\
 
 On a:
 $$Z^*(\four f;(\pi_1,\pi_2))=\int_{E_{0,-2}}\pi_1(F^*_1(y))Z^*(\four_{(E_{2,0})_{y'}}(h_f(\ ,y'));\pi_2)dy'\ .$$
La d\écomposition de la mesure est donn\ée relativement \à la restriction de $P^*_1$ \à $(E_{-2,0})_{y'}$ et \à $F^*_1$ (cf.th.4.3.3 et 4.3.5) cependant dans le \PV 
$((\goth U_{y'^{-1}})_0,(\goth U_{y'^{-1}})_1)$
 il est plus commode de prendre comme \irfs :
$$G(x)=F(x+y'^{-1})\ \text{alors}\ G^*(x)=\ds\frac{1}{G(x'^{-1})}=\ds\frac{1}{F(x+y'^{-1})}=F^*(x'+y')\ ,$$
rappelons que $G$ est une forme quadratique sur un espace vectoriel de dimension $8,$ de discriminant $1$ et d'invariant $\gamma=\pm 1$ suivant le type de $\g$ (ind\épendant de $y'$) (proposition 8.1.2).\\

\noi Notons $Z',Z'^*$ les fonctions Z\étas correspondantes, $\hat f$ la transformation de Fourier associ\ée alors:
$$ \begin{array}{ll}Z^*(\four_{(E_{2,0})_{y'}}(h_f(\ ,y'));\pi_2)&=Z'^*(\widehat{h_f(\ ,y')};\pi_2)\\
&=\gamma\rho (\pi_2|\ |)\rho (\pi_2|\ |^4)Z'(h_f(\ ,y');\pi_2^{-1}|\ |^{-4})\quad (2)\text{th\éor\ème }3.6.5)\\
 &= c|F^*_1(y')|^{-2}Z(h_f(\ ,y');\pi_2^{-1}|\ |^{-4})\ \ \text{en posant  }c=\gamma\rho (\pi_2|\ |)\rho (\pi_2|\ |^4)\end{array}$$
 d'o\ù en permutant dans l'int\égrale::
 $$ Z^*(\four f;(\pi_1,\pi_2))=c\int_{E_{2,0}} |F_1(x)|^{-2}\int_{ (E_{0,-2})_{x}}h_f(x ,y')) \pi_1(F^*_1(y'))\pi'_2(F(x+y'^{-1}))dy'\ $$
 avec $\pi'_2={\pi'_2}^{-1}|\ |^{-4},$ la d\écomposition de la mesure est donn\ée relativement \à $P^*_1.$\\

 \noi Comme avant, dans le \PV $((\goth U_{x})_0,(\goth U_{x})_1)$
 il est plus commode de prendre comme \irfs :
$$H(y)=F(x+y)\ \text{alors}\ H^*(y')=\ds\frac{1}{H(y'^{-1})}=\ds\frac{1}{F(x+y'^{-1})}=F^*(x^{-1}+y')\ ,$$
$H$ est une forme quadratique sur un espace vectoriel de dimension $8,$ de discriminant $1$ et d'invariant $\gamma$ le m\ême que le pr\éc\édent  (proposition 8.1.2) d'o\ù avec les m\êmes notations qu'avant:
$$\int_{ (E_{0,-2})_{x}}h_f(x ,y')) \pi_1(F^*_1(y'))\pi'_2(F(x+y'^{-1}))dy'=\pi_1(F_1(x))Z'^*(\widehat{S_f(x+\ )};{\pi_1}^2{\pi'_2}^{-1})$$
$$\begin{array}{ll}&=\gamma 
\rho (\pi_1^2\pi_2|\ |^5)\rho (\pi_1^2\pi_2|\ |^8)\pi_1(F_1(x))Z'(S_f(x+\ );
\pi_1^{-2}\pi_2^{-1}|\ |^{-8})\quad (2)\text{th\éor\ème }3.6.5)\\
\\
&=d  
 \pi_1(F_1(x))|F_1(x)|^2Z(S_f(x+\ );
\pi_1^{-2}\pi_2^{-1}|\ |^{-8})\ \text{en posant  }d=\gamma \rho (\pi_1^2\pi_2|\ |^5)\rho (\pi_1^2\pi_2|\ |^8)\ ,\end{array}$$
la d\écomposition dans $Z$ est donn\ée relativement \à la restriction de $P_1$ \à $(\goth U_{x})_1$.\\

\noi Ceci nous donne:
$$Z^*(\four f;(\pi_1,\pi_2))=cd\int_{E_{2,0}} \pi_1(F_1(x))( \int_{(E_{2,0})_x}S_f(x+y)  \pi'(F(x+y) )dy)dx=cd{\ }Z(f;(\pi_1,\pi')) $$
(th.4.3.3) avec $\pi'=(\pi_1^2\pi_{2})^{-1}|\ |^{-8}.$\fdem\\ \\

 \begin{cor} Dans le cas r\éel, les polynomes de Bernstein sont donn\és par:
 
 $$b_1(s_1,s_2)=s_2(s_2+3)\ ,\ b_2(s_1,s_2)= s_2(s_2+3)(2s_1+s_2+4)(2s_1+s_2+7).$$
 \end{cor}
 \bigskip

 \dem Comme l'\irf du \PV $(\g_0,\g_1)$ est de degr\é $4,$ par normalisation   $\widetilde B(H_0,H_0)=-2$ donc $\widetilde B(\ds\frac{H_i}{2},\ds\frac{H_i}{2})=-1$ pour $i=1,2$ d'o\ù $b_1$ (1) remarque  3.6.6) et $b_2=\pm s_2(s_2+3)(2s_1+s_2+4)(2s_1+s_2+7)$ (cf.prop 3.7.3), il reste \à d\éterminer le signe ce qui se fait \à l'aide de la relation $$b_2(s_1,s_2)=(2\pi)^4\ds\frac{A((\omega_1,s_1),(\omega_1,s_2))}{A((\omega_{1},s_1),(\omega_{-1},s_2-1))}$$ (2) du lemme 3.7.1 pour $k=2$).\fdem
 \bigskip

\subsection{\bf $\goth g_2$ de dimension $1$}
 
\smallskip
\subsubsection{ Structure}

\smallskip

\noindent Leur diagramme de Dynkin sont donn\'es par : $(F_4,\alpha_1),$ $(E_6,\alpha_2),$ $(E_7,\alpha_1),$ et $(E_8,\alpha_8).$ \\
Ils ont tous une structure commune que l'on rappelle.\\

\noindent $\Delta_2=\{\tilde \alpha\}, $ soient $\lambda_1,\lambda_2,\lambda_3,\lambda_4$ $4$ racines orthogonales de $\Delta_1$ 
 alors $P=P(H_1,H_2)$ avec $H_1=h_{\lambda_1}+h_{\lambda_2}$. Lorsque $\{ \lambda_1,\lambda_2,\lambda_3,\lambda_4\}$ est l'ensemble canonique de racines orthogonales de $\Delta_1$ obtenues par orthogonalisations successives, $P(H_1,H_2)$ est un \sg parabolique standard lorsque $H_1=h_{\lambda_1}+h_{\lambda_2}$.

\noi Rappelons que $2H_0=\sum_{1\leq i\leq 4}h_{\lambda_i}$ et 
 qu'il existe une alg\`ebre d\'eploy\'ee, not\'ee $\tilde {\goth g},$ 
 admettant $\goth a=\oplus_{1\leq i\leq 4}\F h_{\lambda_i}$ 
  comme sous-alg\`ebre de Cartan et  le syst\`eme de racines correspondant, qui  est donn\'e par:
$$R=\{\frac{\pm\lambda_i\pm\lambda_j}{2}\ ,\ 1\leq i\ \not= j\leq 4\ \ ,\ \frac{\pm\lambda_1\pm\lambda_2\pm\lambda_3\pm\lambda_4}{2}\ , \ \pm\lambda_i\ ,\ 1\leq i\leq 4\ \}$$

\noi est de type $F_4.$\\

\smallskip
 \noi On rappelle que $\tilde B=\displaystyle{-\frac{2B}{B(H_0,H_0)}}(=\displaystyle{-\frac{2B}{B(h_{\lambda_i},h_{\lambda_i})}}),$ ainsi $\tilde B(\frac{H_1}{2},\frac{H_1}{2})=-1$ et pour toute racine longue $\alpha$ de $\∆$ on a $\tilde B(X_{\alpha}, X_{-\alpha})=1.$\\

\noindent Pour $i\not= j$ soit $E^{i,j}_{u,v}=\{x\in \goth g\ |\ [h_{\lambda_i},x]=ux\ ,\  [h_{\lambda_j},x]=vx\  ,\ [h_{\lambda_k},x]=0 \ $ pour $ 1\leq k\not=i,j\leq 4\};$ la dimension commune des sous-espaces $E^{i,j}_{\pm 1,\pm 1}$ est not\ée $d,$  on a $d\in\{1,2,4,8\},$   dim$(\goth g_1)=8+6d$ et $r_1=d+\frac{1}{2}.$\\

\noindent On rappelle que (\cite{mullerNAG} , prop.4.1.1):\\

\begin{prop} Il existe un syst\`eme de Chevalley , $(X_{\mu},h_{\mu},X_{-\mu})_{\mu\in R}$, de $\tilde {\goth g}$ tel que:

\begin{enumerate}
\item Pour $1\leq i\ \not= j\leq 4,$ toutes les formes quadratiques $f_{\frac{ \lambda_i- \lambda_j}{2}}$ d\'efinies sur $E^{i,j}_{-1,1}$ par 
$f_{\frac{ \lambda_i- \lambda_j}{2}}(A)=\frac{1}{2} \tilde B(ad(A)^2(X_{\lambda_i}),X_{-\lambda_j}) $
sont \'equivalentes et repr\ésentent $1.$
\item $ [X_{-\lambda_1},[X_{-\lambda_2},[X_{-\lambda_3},[X_{-\lambda_4},X_{\tilde \alpha},]]]]= X_{- \tilde\alpha} .$
\end{enumerate}
\end{prop}

  \bigskip

\noi On note $f=f_{\frac{ \lambda_1- \lambda_2}{2}}$  alors: \\

\begin{lem}  \begin{enumerate}
\item $f$ repr\'esente $1.$
\item Soit $a\in \F^*,$ $f$ et $af$ sont \'equivalentes $\Leftrightarrow$ $a$ est un \'el\'ement de $f(E^{1,1}_{-1,1})^*.$
\item $\goth g$ est de rang $4$ $\Leftrightarrow$ $f$ est anisotrope.
\item $f(E^{1,1}_{-1,1})^*\subset  \chi_1(P).$
\item $\chi (G)=\F^{*2}.$
\end{enumerate}
\end{lem}
 \bigskip
 
 \dem Pour 1), 2), 3) cf. le lemme 4.1.3 de \cite{mullerNAG} et pour 5) la proposition 3.2 de \cite{mullerNAG}.\\ 
 
\noi  Pour 4) on rappelle que $\forall x\in f(E^{1,1}_{-1,1})^*$ et $1\leq i<j \leq 4,$ il existe une involution $g_{i,j}(x)\in G_e$ telle que pour $1\leq k\leq 4,$ on a:
 $$g_{i,j}(x)(h_k)=h_{\tau (k)}\ ,\ g_{i,j}(x)(X_{\lambda_k})=f_{i,j,k}(x)X_{\lambda_{\tau (k)}}\ ,$$
avec  $f_{i,j,i}(x)=x=f_{i,j,j}(x)^{-1}$ et $f_{i,j,k}(x)=1$ sinon, $\tau$ d\'esignant la transposition $(i,j),$   (\demo du lemme 2.2.2 de \cite{mullerNAG}) d'o\`u $g'_{i,j}(x):= g_{i,j}(x) g_{i,j}(1)\in  \cap_{1\≤i\≤4}G_{ h_{\lambda_i}}$ v\'erifie $g'_{i,j}(x)(X_{\lambda_k})=f_{i,j,k}(x)X_{\lambda_k}$  et  $\chi_1(g'_{2,3}(x)) =x$ (cf. la  proposition suivante pour une formule explicite de $F_1$). \fdem\\ 

\begin{rema} $\text{disc}(f)\in \F^{*2}$ lorsque $d=1,4$ ou $8$ (cf.1) et 2) du lemme 8.2.2)\end{rema} 
 \vskip 2mm
\noindent Pour $x\in\goth g_1,$ soient: $ F(x)=\displaystyle\frac{1}{4!} \tilde B(ad(x)^4(X_{-\tilde \alpha}),X_{-\tilde \alpha})$ et $F_1(x)=\frac{1}{2}\tilde B(ad(x)^2(X_{-\tilde \alpha}),X'_{\tilde \alpha-\lambda_1-\lambda_2}),$ avec $X'_{\tilde \alpha-\lambda_i-\lambda_j}=[X_{ -\lambda_i},[X_{ - \lambda_j},X_{\tilde \alpha }]]$ avec $1\leq i\not= j\leq 4.$\\
  
\begin{prop}  
\begin{enumerate}
\item $F$ et $F_1$ sont les invariants relatifs fondamentaux du \PV:  $(P ,\goth g_1)$  v\'erifiant: $F(\sum_{1\leq i\leq 4}X_{ \lambda_i})=F_1(\sum_{1\leq i\leq 4}X_{ \lambda_i})=1.$ 
\item $\chi_1(P)=f(E^{1,1}_{-1,1})^*.$
\item Soit $x\in W_{\goth t},$ $x=x_2+x_0$ avec $x_i\in E_i(  H_1 )\cap \goth g_1$ alors:\\
 
$ \bullet$ $F_1/_{ (E_2(  H_1 )\cap \goth g_1)_{x_0}}\sim -P_1(x_0)X^2\oplus -f\oplus  Y^2-Z^2 ,$\\

 $ \bullet$   $P_1/_{ (E_0(  H_1 )\cap \goth g_1)_{x_2}}\sim -F_1(x_2)X^2\oplus -f\oplus  Y^2-Z^2,$\\

 $ \bullet$  $ \gamma_1(x_2,x_0)=  
 (-F(x),-F_1(x))^d.$     \end{enumerate}
\end{prop}

 \dem

\noi Pour tout $x\in W_{\goth t},$ il existe $g\in G_{H_1}$ et $(a_1,...,a_4)\in (\F^*)^4$ tels que $gx=\sum_{1\leq i\leq 4}a_iX_{  \lambda_i}$ (lemme 7.3 de \cite{mullerJA2}) donc on d\'emontre la proposition pour $x=\sum_{1\leq i\leq 4}a_iX_{  \lambda_i}.$\\

\begin{enumerate}

 \item   En raison des diff\'erentes relations de commutation, on a:
$$ad(x)^4(X_{-\tilde \alpha})=4!\ (\prod_{1\leq i\leq 4}a_i)\ X_{\tilde \alpha}\ \  \text{et}\ \ 
ad(x)^2(X'_{\tilde \alpha-\lambda_1-\lambda_2})=2a_1a_2 \ X_{\tilde \alpha}$$
d'o\`u 1.  

 \item De m\^eme on a:
$$ad(x)^2(X_{\tilde \alpha-\lambda_3-\lambda_4})=2a_3a_4 \ X_{\tilde \alpha}$$ donc $P_1(x)=\frac{1}{2}\tilde B(ad(x)^2(X_{-\tilde \alpha}),X'_{\tilde \alpha-\lambda_3-\lambda_4})$   

 \item  On a $E_2(H_1)\cap \goth g_1=\F X_{\lambda_1}\oplus \F X_{\lambda_2} \oplus \F X'_{\tilde \alpha-\lambda_3 }  \oplus \F X'_{\tilde \alpha -\lambda_4}  \oplus  E^{1,2}_{1,1},$ avec $ X'_{\tilde \alpha -\lambda_i} =[X_{-\lambda_i}, X_{\tilde \alpha}]$ pour $i=1,...,4,$ donc:\\
 
 \noi   $F_1(X)=F_1(A)+t_1t_2+z_1z_2$ lorsque 
 $X=t_1  X_{\lambda_1} +t_2X_{\lambda_2}   +z_1 X'_{\tilde \alpha-\lambda_3 }   +z_2X'_{\tilde \alpha -\lambda_4}+A.$   \\
 
 Soit $\theta=\theta_{X_{ \lambda_1}, h_{\lambda_1} }(-1)$ et $u\in E^{1,2}_{-1,1},$ alors :
$$F_1(\theta^{-1}(u))=\frac{1}{ 2}\tilde B(ad(u)^2(\theta X_{-\tilde \alpha}),\theta X'_{\tilde \alpha-\lambda_1-\lambda_2)   }    )  
 =- \frac{1}{ 2}\tilde B( [u,\theta X_{-\tilde \alpha}],v)\ \ \hbox{avec} $$
  $$v=[u,[ X_{\lambda_1},[X_{-\lambda_2},\theta X_{\tilde {\alpha} }]]  ] = [[[u,X_{\lambda_1}  ],X_{-\lambda_2}]  ,\theta (X_{\tilde \alpha}) ]$$   
   en raison des relations de commutation donc :
   \begin{eqnarray*} F_1(\theta ^{-1}(u))& = & \frac{1}{ 2}\tilde B([[u, \theta (X_{-\tilde \alpha})],\theta (X_{\tilde \alpha})],[[u,X_{\lambda_1}],X_{-\lambda_2}])\\
 & =  & -\frac{1}{ 2}\tilde B(u,[[u,X_{\lambda_1}],X_{-\lambda_2}])\\
& = & -f(u).\\
\end{eqnarray*}
\noi Comme $\chi_1(P)F_1\sim F_1$ et que $F_1\sim -f\oplus X^2-Y^2+Z^2-T^2,$ on v\'erifie, dans le cas r\'eel et le cas $\goth p$-adique, que $\chi_1(P)\subset f(E^{1,1}_{-1,1})^*$ lorsque $f(E^{1,1}_{-1,1})^*\not=\F^*$ d'o\`u $\chi_1(P) = f(E^{1,1}_{-1,1})^*$ en appliquant le lemme pr\'ec\'edent.\\ \\

\item On peut supposer que  $x_0=a_3X_{\lambda_3}+a_4X_{\lambda_4}$ et que $x_2=a_1X_{\lambda_1}+a_2X_{\lambda_2},$ un calcul imm\'ediat donne : 

$$(E_2(H_1)\cap \goth g_1)_{x_0}=\F X_{\lambda_1}\oplus \F X_{\lambda_2}   \oplus \F (a_3X'_{\tilde \alpha-\lambda_4}  -a_4X'_{\tilde \alpha-\lambda_3 }  )\oplus  E^{1,2}_{1,1},$$

 \noi de m\^eme pour $x_2$ d'o\`u le r\'esultat du 3).\\
 
  \item On a :
 $$E_{-1}(H_1)\cap E_1(H_2)=E^{1,3}_{-1,1}\oplus E^{1,4}_{-1,1}\oplus  E^{2,3}_{-1,1}\oplus E^{2,3}_{-1,1}\ .$$
 
 \noi Prenons
 comme pr\'ec\'edemment  $x_2=a_1X_{\lambda_1}+a_2X_{\lambda_2}$ et $x_0=a_3X_{\lambda_3}+a_4X_{\lambda_4},$ alors :

  $$\begin{array}{rll}Q_{x_2,x_0}(\sum_{i=1,2;j=3,4}A_{i,j })& = & \frac{1}{2}\tilde B([x_2,[x_0^{-1},\sum_{i=1,2;j=3,4}A_{i,j}]],\sum_{i=1,2;j=3,4}A_{i,j})\\
 &=& \sum_{i=1,2;j=3,4}\frac{a_i}{a_j} f_{\frac{\lambda_i-\lambda_j}{2}}(A_{i,j})\ ,\ \hbox{avec}\  A_{i,j}\in E^{i,j}_{-1,1}\ ,\end{array}$$
d'o\`u $ \gamma_1(x_2,x_0)= \begin{cases}
 \gamma (f)^4\  \text{si }\    d \ \text{est pair}, \\
(-F(x),-F_1(x))\  \text{si }\ d=1.
\end{cases}$ 

\noi Pour $d=4$ ou $8$ on a $\gamma (f)=\pm 1;$  pour $d=2$ on a $\gamma (f)=\alpha(disc(f))\alpha (1)$ d'o\ù $\gamma (f)^2=(-1,-disc(f))$ (cf.2) lemme 8.2.2)\fdem\\ \\

  \end{enumerate}
Lorsque $\F$ est archim\'edien, les $2$ polynomes de Bernstein s'obtiennent imm\'ediatement en appliquant les  propositions  3.4.4 et  3.7.3 puisque
  $F_1$ et $P_1$ sont des formes quadratiques (cf.remarque 3.6.6) d'o\`u:
  
  \begin{prop}  $$b_1(s_1,s_2)=s_2(s_2+\frac{d+1}{2})\ \hbox{et}\  b_2(s_1,s_2)=s_2(s_2+\frac{d+1}{2})(s_1+s_2+d+\frac{1}{2})(s_1+s_2+ \frac{3d}{2}+1).$$
  \end{prop}
  
  \bigskip
 \subsubsection{Une premi\ère \équation fonctionnelle} 
\bigskip 
\noi   On rappelle que pour $u=(u_1,u_2)\in (\FF)^2$ on d\'efinit les ouverts (tous non vides):\\

$O_u=O_{(u_1,u_2)}=\{x\in \goth
g_1\ |\ F_1(x)\F^{*2}=  u_1\  \ ,\ F_2(x)\F^{*2} = u_1u_2  
  \}$\\
 
$O^*_u=O^*_{(u_1, u_2)}=\{x\in
\goth g_{-1}\ |\ F^*_1(x) \F^{*2}=u_2 \ ,\ F^*_2(x) \F^{*2}=
 u_1u_2   \}$\\
 
\noi ainsi que les fonctions Z\'etas correspondantes:\\
  
$Z_u(f;\omega)=Z(f{\bf 1}_{O_u};\omega)$ pour $f\in S(\goth
g_1)$  
et
$ Z^*_u(h;\omega)=Z^*(h{\bf 1}_{O^*_u};\omega)$  pour $h\in S(\goth g_{-1})$.  \\ \\

\noi   Soit  $(-1)^{[\frac{d+3}{2}]}c(f) $ le discriminant de la forme quadratique $-f\oplus X^2-Y^2-Z^2,$ alors $c(f)=-1$ sauf lorsque $d=2$ avec $f$ anisotrope.  

\noi Notons que pour $d$ pair on a $\gamma (f)^2=(-1,-c(f)).$

\noi (cf.lemme 8.2.2 et \demo de la proposition 8.2.4).\\

\noi  Alors:\\

\begin{prop} 
\begin{enumerate}

\item Lorsque $\F=\C$ pour $f\in S(\goth
g_1),$  on a:
$$ Z^*( \four (
f);\omega,s)=a_d(\omega,s)Z\biggl(f;( \omega_1,( \omega_1\omega_2)^{-1});s_1,-(s_1+s_2+\frac{3}{2}d+2 )\biggr)$$
avec:
$$a_d(\omega,s)=\rho'(\omega_2;s_2+1)\rho'(\omega_2;s_2+\frac{d+3}{2})\rho'(\omega_1\omega_2;s_1+s_2+d+\frac{3}{2})\rho'(\omega_1\omega_2;s_1+s_2+\frac{3d}{2}+2).$$
\item Dans les autres cas, pour tout  $v=(v_1, v_2)$ dans $(\FF)^2$ et $f\in S(\goth
g_1),$  on a:
$$ Z^*_v( \four(
f);\ \omega,s)=\sum_{u\in (\FF)^2}d_{v,u}(
 \omega,s)  Z_u(f;( \omega_1,( \omega_1\omega_2)^{-1});s_1,-(s_1+s_2+\frac{3}{2}d+2 )\ ) \quad \hbox{avec}$$

$\bullet$  lorsque $d$ est pair et  $u=(u_1, u_2)$ :  $d_{v,u}(
\omega,s)= (-1, c(f)) .(u_1v_2,-c(f))$
$$\displaystyle{  A^1_{ \omega_2,s_2+1,s_2+\frac{d+3}{2}}(v_1v_2,u_1v_2,c(f))A^1_{ \omega_1 \omega_2,s_1+s_2+d+\frac{3}{2},s_1+s_2+\frac{3d}{2}+ 2} (u_1v_2,u_1u_2,c(f))}.$$
 
 $\bullet$  Pour $d=1:$  $d_{v,u}(
\chi,s)=$
 $$  \alpha (-1)\ \alpha (u_1v_2)  A^0_{\omega_2,s_2+1,s_2+2}(v_1,u_1,-v_2)A^0_{ \omega_1 \omega_2,s_1+s_2+\frac{5}{2},s_1+s_2+ \frac{7}{2}} (v_2,u_2,-u_1).$$ 
 \end{enumerate}

\end{prop}

\dem

\noi On applique la proposition 5.3.1 dont toutes les hypoth\`eses sont v\'erifi\'ees par le 3) de la proposition 8.2.4 ainsi que le  th\éor\ème 3.6.5,5).  \fdem\\

\begin{rema}  Dans le cas r\éel, soient:
 
$$\begin{array}{ll}C_d(s)&= 2.(2\pi )^{-(2s+\frac{d+5}{2})}\Gamma (s+1)\Gamma (s+\frac{d+3}{2})\quad , \quad s\in \C\ , \\
\\
 C_d(s_1,s_2)&=C_d(s_2).C_d(s_1+s_2+d+\frac{1}{2})\quad , \quad s_1,s_2\in \C\ , \\
 \\
&=4(2\pi)^{-(2s_1+4s_2+3d+6)}\Gamma(s_2+1)\Gamma(s_2+\frac{d+3}{2})\Gamma(s_1+s_2+d+\frac{3}{2})\Gamma(s_1+s_2+\frac{3d}{2}+2)\ ,\end{array}$$
et pour $a,b,c$ r\éels:$$\phi_s^{(f)}(a,b,c)=\frac{\pi}{2}\bigl( (a+b)s+\frac{d+1}{2}a-c\frac{c(f)}{2}\bigr)$$
alors pour $u,v\in \{-1,1\}^2$ et $\omega=Id$ on a:
$$\frac{d_{v,u}(s_1,s_2)}{C_d(s_1,s_2)}=\begin{cases} (-1)^{\frac{1+u_1}{2}\frac{1+v_2}{2}}\cos \phi_{s_2+1}^{(f)}(u_1,v_1,1+v_2).
\cos \phi_{s_1+s_2+\frac{5}{2}}^{(f)}(u_2,v_2,1+u_1) \text{ pour }\ d=1,\\
\\
(u_1v_2,-c(f ))\cos \phi_{s_2+1}^{(f)}(u_1,v_1,v_2).
\cos \phi_{s_1+s_2+d+\frac{3}{2}}^{(f)}(u_2,v_2,u_1)\  \ \ \text{ pour }\ d \ \hbox{pair}.\end{cases}$$ Indexons $\{-1,1\}\times \{-1,1\}:=\{a_1=(-1,-1),a_2=(-1,1),a_3=(1,-1),a_4=(1,1)\}$ et soit $M (d,c(f))$ la matrice \à $4$ lignes et $4$ colonnes et \à coefficients r\éels dont le $(i,j)-$\ème coefficient, not\é  $M_{i,j}(d,c(f)),$ est donn\é par:
$$   M_{i,j}(d,c(f))=\frac{d_{a_i,a_j}(s_1,s_2)}{C_d(s_1,s_2)}\ ,$$alors  on a:\\

 $\bullet$ $\bf M(1,-1)=$

$$\displaystyle\left(\begin{array}{cccc}-\sin\pi s_2\cos \pi(s_1+s_2)&0&0&0\\
0&\cos \pi s_2\cos \pi(s_1+s_2)&1&\sin \pi(s_1+s_2)\\
0&0&-\sin \pi s_2\sin \pi(s_1+s_2)&-\sin \pi s_2\\
0&-\cos \pi (s_1+s_2)&-\cos \pi s_2&-\cos \pi s_2\sin \pi(s_1+s_2)\end{array}\right)$$\\ 

$\bullet$   Pour $d=4$ ou $8,$ $\bf M(d,-1)=$

$$\displaystyle\left(\begin{array}{cccc}\sin\pi s_2\cos \pi(s_1+s_2)&\sin \pi s_2&\sin \pi (s_1+s_2)&0\\
0& -\cos\pi s_2\sin \pi(s_1+s_2)&0& 0\\
0&0& -\cos \pi s_2\sin \pi(s_1+s_2)& 0\\
0&\sin \pi(s_1+s_2)&  \sin \pi s_2&\sin \pi s_2 \cos \pi(s_1+s_2)\end{array}\right)$$\\

$\bullet$ $\bf M(2,-1)=$  

$$ \left(\begin{array}{cccc}-\cos\pi s_2\sin \pi(s_1+s_2)& 0& 0&0\\
 -\sin\pi s_2&\sin\pi s_2\cos \pi(s_1+s_2)&0&  \sin \pi(s_1+s_2)\\
  \sin \pi(s_1+s_2)&0&\sin \pi s_2\cos \pi(s_1+s_2)& -\sin\pi s_2\\
0& 0& 0&- \cos \pi s_2 \sin \pi(s_1+s_2)\end{array}\right)$$\\ \\

$\bullet$ $\bf M(2,1)=$  

$$ \left(\begin{array}{cccc}\sin\pi s_2\cos\pi (s_1+s_2)& -\sin\pi s_2& -\sin\pi(s_1+s_2) &0\\

  0&\cos \pi s_2\sin\pi(s_1+s_2)&0& 0 \\

 0 &0& \cos \pi s_2\sin \pi (s_1+s_2)& 0\\

0& -\sin\pi(s_1+s_2)&-\sin \pi s_2 & \sin \pi s_2 \cos \pi(s_1+s_2)\end{array}\right).$$
\end{rema}

\bigskip

\bigskip

\noi Dans les cas r\éels et $\goth p-$adique, indexons les \elts de $\FF$: $\FF=\{u_1,...,u_\ell\},$  et soit $A(\omega,s)$ la matrice \à $\ell$ lignes et $\ell$ colonnes des coefficients de l'\équation fonctionnelle "normalis\ée" relative \à la forme quadratique $Q_f=-f\oplus X^2-Y^2-Z^2$ lorsque $d$ est {\bf pair} c'est \à dire:
$$\begin{array}{ll}A(\omega,s)&=\biggl(a^{(Q_f)}_{u_i,u_j}(\omega,s)\biggr)_{1\≤i,j\≤\ell}\\
&=\epsilon_1 \biggl(A^1_{\omega,s+1,s+\frac{d+3}{2}}(u_i,u_j,c(f))\biggr)_{1\≤i,j\≤\ell}\\
&=\epsilon_2 \biggl(A^1_{\omega,s+1,s+\frac{d+3}{2}}(c(f)u_i,c(f)u_j, 1)\biggr)_{1\≤i,j\≤\ell}\end{array}$$
avec $\epsilon_1:=\alpha(1)^2(c(f),-1)\gamma (f)$ et $\epsilon_2=\bigg\{ \begin{array}{ll} (-1,-c(f)) \ &\hbox{ lorsque}\ d=2\\  \gamma (f)  \ &\hbox{ lorsque}\ d\in \{4,8\}\ ,\end{array}$

\noi (\demo de la proposition 8.2.4),  on a $\epsilon_1^2=(-1,c(f))$ et $\epsilon_2^2=1,$

\noi   alors:\\

\vskip 5mm
\begin{prop} On suppose que $d$ est pair.\\

Pour $u\in \FF,$ et $\pi\in \hat {(\F^*)^2}$ tel que $\Re (\pi)>0,$ on pose:

$\bullet$ $O_u=\{x\in \g_1$ tel que $F_2(x) (\F^*)^2=u\}$  et   pour $f\in \EuScript S(\g_1):$  $Z_u(f;\pi):=Z(1_{O_u}f;\pi),$

$\bullet$ $O^*_u=\{x\in \g_{-1}$ tel que $F^*_2(x)  (\F^*)^2=v\}$ et pour $g\in \EuScript S(\g_{-1}):$   $Z^*_u(g;\pi):=Z^*(1_{O^*_u}g;\pi),$

alors  $Z_u(f;\pi)$ et $Z^*_u(g;\pi)$ admettent un prolongement m\éromorphe sur $ \hat {(\F^*)^2}$ et  on a:
$$
\left(\begin{array}{ll}Z^*_{u_1}(\four (f);\omega,s)\\
...\ ...\\
Z^*_{u_\ell}\ (\four (f);\omega,s)\end{array}\right)=   B(\omega,s)\left(\begin{array}{ll}Z_{ u_1}(f; \omega^*|\ |^{-(\frac{3d}{2}+2)1_2})\\
...\ ...\\
Z_{u_\ell}\ ( f; \omega^*|\ |^{-(\frac{3d}{2}+2)1_2})\end{array}\right)$$
avec
$$B(\omega,s)=A(\omega_2,s_2).D.A(\omega_1\omega_2,s_1+s_2+d+\frac{1}{2}),$$
$D$ \étant la matrice diagonale : $D=((u_i,-c(f)).\delta_{i,j})_{1\≤i,j\≤\ell}.$ 
 \end{prop}

\dem Pour  $f\in \EuScript S(\g_1),$ $g\in \EuScript S(\g_{-1})$ et 
$\pi\in \hat {(\F^*)^2}$ tel que $\Re (\pi)>0,$ on a les \égalit\és:
$$Z_u(f;\pi)=\sum_{w\in \FF}Z_{w(1,u)}(f;\pi)\ ,\ Z^*_u(g;\pi)=\sum_{w\in \FF}Z^*_{w(1,u)}(g;\pi)\ ,$$
ce qui donne un prolongement m\éromorphe de $Z_u$ et $Z^*_u$ sur $ \hat {(\F^*)^2}.$

\noi  On a clairement la relation:
$$u,v\in (\FF)^2\ ,\ w\in \FF : d_{v,w.u}=d_{w.v,u}\ \hbox{
donc}:$$
$$\sum_{w\in \FF}d_{w(1,v),u_1(1,u)}=\sum_{w\in \FF}d_{(1,v),wu_1(1,u)}=\sum_{w\in \FF}d_{(1,v),w(1,u)}=\sum_{w\in \FF}d_{w(1,v),(1,u)},$$
d'o\ù  
$$ Z^*_v(\four (f);\pi)=\sum_{v\in \FF}a_{v,u}(\pi)Z_u(f;\pi^*|\ |^{-(\frac{3d}{2}+2)1_2}) \quad 
\hbox{avec}  $$
$$a_{v,u}= (-1,c(f))\sum_{w\in \FF} A^1_{ \omega_2,s_2+1,s_2+\frac{d+3}{2}}(v ,w,c(f))(w,-c(f))A^1_{ \omega_1 \omega_2,s_1+s_2+d+\frac{3}{2},s_1+s_2+\frac{3d}{2}+ 2} ( w,u,c(f)),$$
d'o\ù le r\ésultat.\\

\noi Notons que l'on a \également:
$$\begin{array}{ll}a_{v,u}&= \sum_{w\in \FF} A^1_{ \omega_2,s_2+1,s_2+\frac{d+3}{2}}(c(f)v ,c(f)w, 1)(w,-c(f))A^1_{ \omega_1 \omega_2,s_1+s_2+d+\frac{3}{2},s_1+s_2+\frac{3d}{2}+ 2} ( c(f)w,c(f)u, 1)\\
&=\sum_{w\in \FF} A^1_{ \omega_2,s_2+1,s_2+\frac{d+3}{2}}(c(f)v ,w, 1)(w,-c(f))A^1_{ \omega_1 \omega_2,s_1+s_2+d+\frac{3}{2},s_1+s_2+\frac{3d}{2}+ 2} (w,c(f)u, 1).\end{array}$$\fdem\\

\begin{rema}: \end{rema}

Soit $\tilde A(\omega,s)$ (resp.$\tilde A(s)$) la matrice \à $\ell$ lignes et $\ell$ colonnes ayant comme coefficient 

\noi $A^1_{ \omega,s+1,s+\frac{d+3}{2}}(u_i,u_j,1)$ (resp.$A^1_{  Id,s+1,s+\frac{d+3}{2}}(u_i,u_j,1)$) sur la $i-$\ème ligne et  $j-$ \ème colonne, alors : \\

$$
\left(\begin{array}{ll}Z^*_{c(f)u_1}(\four (f);\omega,s)\\
...\ ...\\
Z^*_{c(f)u_\ell}\ (\four (f);\omega,s)\end{array}\right)=   \tilde B(\omega,s)\left(\begin{array}{ll}Z_{ c(f)u_1}(f; \omega^*|\ |^{-(\frac{3d}{2}+2)1_2})\\
...\ ...\\
Z_{c(f)u_\ell}\ ( f; \omega^*|\ |^{-(\frac{3d}{2}+2)1_2})\end{array}\right).$$
 avec 
 $$\tilde B(\omega,s)=\tilde A(\omega_2,s_2).D.\tilde A(\omega_1\omega_2,s_1+s_2+d+\frac{1}{2}).$$\\

$\bullet$ Dans le cas r\éel :

\noi On prend $\frac{\R^*}{\R^{*+}}=\{u_1=-1,u_2=1\},$  pour $\omega=Id$ on a:

$$\tilde A(s)= C_d(s)\left(\begin{array}{ll}\sin (\pi s+\pi \frac{d}{4})&\cos \pi \frac{d}{4}\\
-\sin \pi \frac{d}{4}&-\cos (\pi s+\pi \frac{d}{4})\end{array}\right).$$

\noi On pose $B(s_1,s_2):=B( Id,(s_1,s_2)),$  alors:\\

  Lorsque $d=2$ et  $c(f)=-1$ (c'est \à dire $f$  isotrope):
$$B(s_1,s_2)=C_2(s_1,s_2)\left(\begin{array}{ccc}
  \sin \pi s_2\cos\pi (s_1+s_2)& \sin\pi (s_1+s_2)-\sin \pi s_2\\
   0& -\cos \pi s_2\sin\pi (s_1+s_2)\end{array}\right),$$ 
   et dans tous les autres cas ($d=2$ avec $f$ anisotrope  ou $d=4$ ou $d=8$):
   $$B(s_1,s_2)=C_d(s_1,s_2)\left(\begin{array}{ccc}
(-1)^{\frac{d}{2}+1}   \cos \pi s_2\sin\pi (s_1+s_2)&  0\\
  (-1)^{\frac{d}{2}}  (\sin \pi s_2+\sin\pi (s_1+s_2))&  \sin \pi s_2\cos\pi (s_1+s_2)\end{array}\right).$$
  
  $\bullet$ Dans le cas $\goth p$-adique, de caract\éristique r\ésiduelle diff\érente de $2,$ la matrice $\tilde A$ est donn\ée dans le lemme 3.6.7,B. \\ \\
  
Contrairement au cas $d$ pair, lorsque $d=1$ et $\F=\R,$ on constate que
$ \sum_{w\in \FF}d_{w(1,v),u_1(1,u)}(s_1,s_2)$ d\épend \également de $u_1$ lorsque $(s_1,s_2)$ d\écrit  $(\C)^2,$ cependant lorsque $s_1=0$ on a:\\

\begin{lem} Cas $d=1.$

\noi Dans le cas r\éel ou $\goth p$-adique de caract\éristique r\ésiduelle diff\érente de $2,$ on a pour $s\in \C,u$ et $v$ dans $\FF:$
$$ \sum_{w\in \FF}d_{w(1,v),u_1(1,u)}(0,s)=\alpha(-1)^2 {|2|_\F}^{-2s-\frac{7}{2}}\rho (|\ |^{2s+4})A^1_{s+1,s+\frac{7}{2}}(v,u,-1).$$
 \end{lem} 
 
  \dem On exprime 
  $$\begin{array}{ll}f(u_1)&= \alpha(-1)\sum_{w\in \FF}d_{w(1,v),u_1(1,u)}(0,s)\\
  &=\alpha(-1)\ \sum_{w,t,y\in \FF}\alpha(u_1vw)\ (-vw,t)\ (-u_1,y).\\
  &\rho (s+1;tw)\
 \rho (s+2;t u_1)\ \rho (s+\frac{5}{2};yvw)\ \rho (s+\frac{7}{2}; yuu_1)\end{array}$$
 On change les sommations:  $y$ en $yu_1,$  $t$ en $tu_1$  et $w$ en $wu_1v$ ce qui donne:
$$f(u_1)= \alpha(-1)\ \sum_{w,t,y\in \FF}\alpha(w)(-u_1w,tu_1)(-u_1,yu_1)\rho (s+1;tvw)\
 \rho (s+2;t  )\rho (s+\frac{5}{2};yw)\ \rho (s+\frac{7}{2}; yu )$$
 ensuite on somme sur $x=tyw,y,w$ d'o\ù: $\alpha(1)\ f(u_1)=$
 $$   \sum_{x\in \FF} (-u_1, x)\biggl(\sum_{y\in \FF}  \rho (s+1;xyv)\rho (s+\frac{7}{2}; yu ) 
 \biggl(\sum_{w\in \FF}
 \alpha(w)(w,xy)\rho (s+2;xyw  )\rho (s+\frac{5}{2};yw)\bigg)\biggr)$$
 or $$\sum_{w\in \FF}\alpha(w)(w,xy)\rho (s+2;xyw  )\rho (s+\frac{5}{2};yw)=A^1_{s+2,s+\frac{5}{2}}(xy,y ,xy)$$
 $$=\alpha(-xy)\alpha(1)A^1_{s+2,s+\frac{5}{2}}(1,x ,1)=\begin{cases} 0\text{ si }x\not=1,\\
 \alpha(-y)\alpha(1)A^1_{s+2,s+\frac{5}{2}}(1,1 ,1)\text{ si }x=1\end{cases}$$
(lemme 3.6.8,C))  d'o\ù:
 $$
 f(u_1)=  \ A^1_{s+2,s+\frac{5}{2}}(1,1 ,1) \sum_{y\in \FF} \alpha(-y)\rho (s+1;yv)\rho (s+\frac{7}{2}; yu ) $$
  $$=A^1_{s+2,s+\frac{5}{2}}(1,1 ,1)A^1_{s+1,s+\frac{7}{2}}(-v,-u ,1)
 =\alpha(-1)^2\   {|2|_\F}^{-2s-\frac{7}{2}}\rho (|\ |^{2s+4})A^1_{s+1,s+\frac{7}{2}}(v,u ,-1).$$
 (lemmes 3.6.4,B,3),3.6.8,1))\fdem \\
 
 \noi Rappelons que pour $u\in \FF:$  
 
 $\bullet$ $O'_u=\{x\in \g_1|F(x)\F^{*2}=u\}$ et $Z_u(f;\pi)=Z(f1_{O'_u};\pi),$\\
 
  $\bullet$ $O'*_u=\{x\in \g_{-1}|F^*(x)\F^{*2}=u\},$ et $Z^*_u(f;\pi)=Z(f1_{O'^*_u};\pi),$\\
 
 \noi alors la proposition 8.2.8 et le lemme 8.2.10 impliquent:\\
 
 \begin{prop}  On suppose que $\F=\R$  ou bien que $\F$ est un corps $\goth p$-adique de caract\éristique r\ésiduelle diff\érente de $2,$ alors \\
 
 \noi Pour tout $f\in \ES(\g_1),$ $ s\in \C$ et $v\in \FF$ on a:
 $$Z_v^*(\four(f); s)=\sum_{u\in \FF}a_{v,u}( s)Z_u(f;  -s-2-\frac{3}{2}d)$$
 avec 
$$\begin{array}{lll}a_{v,u}( s)&=&\alpha(-1)^2   {|2|_\F}^{-2s-\frac{7}{2}}\rho (|\ |^{2s+4})A^1_{s+1,s+\frac{7}{2}}(v,u,-1) \ \text{ lorsque }\ d=1,\\ 
\\  
 a_{v,u}(s)&=& \sum_{w\in \FF} A^1_{ s+1, s +\frac{d+3}{2}}(c(f)v ,w, 1)(w,-c(f))A^1_{ s +d+\frac{3}{2},s +\frac{3d}{2}+ 2} ( w,c(f)u, 1) \end{array}$$
lorsque $d$  est pair.
 \end{prop}
 
 \noi {\bf Remarque:} Les coefficients $a_{v,u}(s),v$ et $u\in \FF,$ de la proposition sont les  coefficients des \eqs   associ\ées \à un \irf de la forme : $$F_X(x)=\tilde B(adx^4(X),X),\ \text{ avec }\ X\in \g_{-2}\setminus\{0\}.$$
Lorsque $d$ est pair, la constante $c(f)$ associ\ée vaut toujours $-1$ sauf lorsque $d=2$  mais alors celle-ci est \égale au discriminant de l'une (quelconque) des formes quadratiques $A\in E_{-1,1}^{i,j}\rightarrow \tilde B(ad(A)^2(x),y),$ avec $x\in E_2(h_i)\cap \g_1\setminus \{0\},y\in E_{-2}(h_j)\cap \g_1\setminus \{0\},$ ou bien $B\in E_{1,1}^{i,j}\rightarrow \tilde B(ad(B)^2(z),y),$ avec $z\in E_{-2}(h_i)\cap \g_{-1}\setminus \{0\},y\in E_{-2}(h_j)\cap \g_{-1}\setminus \{0\};$ $c(f)\not= -1$ \ssi $\g$ est de rang $4.$\\
\noi Les r\ésultats explicites se d\éduisent de la remarque 8.2.9 en prenant $(s_1,s_2)=(0,s).$\\

 \subsubsection{Orbites}
 \bigskip
Dans ce paragraphe, on note $D$ un ensemble de repr\ésentants dans $\F^*$ de $\FF$, $D=\{-1,1\}$ dans le cas r\éel. On confond $D$ et $\FF$ c'est \à dire qu'un \elt $u\in D$ repr\ésente soit $u$ dans $\F^*$ soit $u\F^{*2}$ dans $\FF$ suivant le contexte mais il est  not\é $u$ dans les $2$ cas.  \\
\begin{prop} 
\begin{enumerate}
\item Lorsque $f$ est isotrope, les orbites de $G$ (resp.$P$) dans $\goth g_1'$ (resp.$\goth g_1"$) sont en bijection avec $\F^*/\F^{*2}$; un ensemble de repr\'esentants est donn\'es dans les $2$ cas par $X_u= \sum_{1\leq i\leq 3}X_{\lambda_i}+uX_{\lambda_4}$ pour $u\in D$ et $P.X_u=\cup_{(a_1,a_2)\in  (\FF)^2|a_1a_2 =u\F^{^2}}O_{a_1,a_2}.$
\item Lorsque $\F=\R$ et $f$ est anisotrope:

i) Il y a $3$ orbites de $G$ dans $\goth g_1'$ de repr\'esentants:
 $$X_0=\sum_{1\leq i\leq 4}X_{\lambda_i}\ ,\ X_1=-X_{\lambda_1}+\sum_{2\leq i\leq 4}X_{\lambda_i}\ ,\ X_2=-X_{\lambda_1}-X_{\lambda_2}+X_{\lambda_3}+X_{\lambda_4}.$$ 

ii)  Il y a $5$ orbites de $P$ dans $\goth g_1''$ de repr\'esentants:
$$X_0\ ,\ X_1\ ,\ X_2\ , \ X_3 =X_{\lambda_1}+ X_{\lambda_2}-X_{\lambda_3}+X_{\lambda_4}\ ,\ X_4=-X_{\lambda_1}+ X_{\lambda_2}-X_{\lambda_3}+X_{\lambda_4 }.$$
\item  Dans le cas $\goth p$-adique et $f$ est anisotrope:

i)  Les orbites de $G$ dans $\goth g_1'$ sont en bijection avec $\F^*/\F^{*2}$; un ensemble de repr\'esentants est donn\'es par $X_u= \sum_{1\leq i\leq 3}X_{\lambda_i}+uX_{\lambda_4}$ pour $u\in D.$

ii)  $d=2$

\noi  Il y a $2||\F^*/\F^{*2}||$ P-orbites
donn\'ees pour $u\in D$ par:
$$P.X_{u,\pm 1}=\{x\in \goth g''_{1}\ |\ F(x)\equiv u \ (\F^{*2})\ ,(F_1(x),-c(f))=\pm 1\}$$
 avec $X_{u,1}=X_u $ et $X_{u,-1}=v_0X_{\lambda_1}+   X_{\lambda_2}+ X_{\lambda_3}+\frac{u}{v_0}X_{\lambda_4},$ $v_0\in D$ v\'erifiant $(v_0,-c(f))=-1.$
 
 \item Les r\ésultats sont analogues dans $\g'_{-1}$ (resp. $\g"_{-1}).$
\end{enumerate}
\end{prop}

\dem

\noi On utilise les r\'esultats de \cite{mullerNAG} o\ù figurent les $G$-orbites dans $\goth g_1'$ (corollaire 4.3.3, corollaire 4.3.4 i)) d'o\ù 3)i) et 1) avec $f$ isotrope;   dans le cas r\éel on rappelle que $X= \sum_{1\leq i\leq 4}a_iX_{\lambda_i}\in \g'_1 $ et $Y= \sum_{1\leq i\leq 4}b_iX_{\lambda_i}\in \g'_1 $ sont dans la m\ême $G-$orbite \ssi il existe $t\in \R^*$ tel que les formes quadratiques $t(\oplus_{1\≤i\≤4}a_if)$ et $\oplus_{1\≤i\≤4}b_if$ soient \équivalentes d'o\ù 2)i).

\noi On rappelle \également que $\chi(G)=\F^{*2}$  et que $\chi_1(P)=f(E)^*.$

\noi De plus, pour $X=\sum_{1\leq i\leq 4}a_iX_{\lambda_i}$ tel que $\prod_{1\leq i\leq 4}a_i\not=0$   et 
pour tout $t$ non nul repr\'esent\'e par la forme quadratique $a_1f\oplus a_2 f$ (resp.$a_3f\oplus a_4 f$) il existe  un  automorphisme  \'el\'ementaire $g_t$ (resp.$g'_t$) qui centralise $H_1$ tel que:
$$\begin{array}{rll}g_t(\sum_{1\leq i\leq 4}a_iX_{\lambda_i} )&=&tX_{\lambda_1} +\ds\frac{a_1a_2}{t}X_{\lambda_2}+  a_3X_{\lambda_3} +a_4X_{\lambda_4}\ \\
 (\hbox{resp.}\ g'_t(\sum_{1\leq i\leq 4}a_iX_{\lambda_i} )&=& a_1X_{\lambda_1} +a_2X_{\lambda_2} +tX_{\lambda_3} +\ds\frac{a_3a_4}{t}X_{\lambda_4}\ )\end{array}$$ 
(cf.d\'emonstration 0) de la proposition 4.1.5 de \cite{mullerNAG}).

\noi Lorsque $a_1f\oplus a_2 f$ (resp.$a_3f\oplus a_4 f$) repr\ésente $1,$ on peut donc toujours supposer que le coefficient de $X_{\lambda_1}$ ou bien $X_{\lambda_2}$ vaut $1$ en utilisant $g_{1,2}(1)$ (resp. le coefficient de $X_{\lambda_3}$ ou bien $X_{\lambda_4} $  vaut $1$ en utilisant $g_{3,4}(1)$) (*).

\begin{enumerate}
\item Lorsque $f$ est isotrope, on  v\'erifie que $\cup_{(a_1,a_2)\in (\FF)^2|\ a_1a_2 = u\F^{*2}}O_{a_1,a_2}=P.X_u$ en utilisant les \'el\'ements $g'_{i,j}(x)$ d\'efinis dans la d\'emonstration du 4) du lemme 8.2.2  ainsi que $h_{\lambda_i}(t),i=1,...,4$ et $t\in \F^*,$ ce qui termine la \demo de 1).

\item Comme  $\chi_1(P)=\R^{*+}$ dans le cas r\éel, les \elts \énum\ér\és dans 2)ii) ne sont pas dans la m\ême $P-$orbite et tout \elt de $\oplus_{1\≤i\≤4}\R^*X_{\lambda_i}$ est dans la $P-$orbite de l'un des \elt \énum\ér\és dans 2)ii) par le r\ésultat (*) et la multiplication par $\R^*.$

\item Dans le cas  $\goth p$-adique anisotrope avec  $d=2.$

Comme $\chi_1(P)=f(E)^*=\{x\in \F^*|(x,-c(f))=1\}$ et $\chi (G)=\F^{*2},$ $2$ \elts quelconques de $\{X_{u,1},X_{u,-1},u\in D\}$ ne sont pas dans la m\ême $P-$orbite. 

\noi Pour $xy\in \F^*$ la forme quadratique $xf\oplus yf$ repr\ésente $\F^*$ donc par (*) toute $P-$orbite dans $\g"_1$ a un repr\ésentant de la forme $X=aX_{\lambda_1}+X_{\lambda_2}+X_{\lambda_3}+bX_{\lambda_4},$ soit $x\in\{1,v_0\}$ tel que $(xa,-c(f))=1$ alors $g'_{1,4}(\ds\frac{x}{a})X=xX_{\lambda_1}+X_{\lambda_2}+X_{\lambda_3}+\ds\frac{ab}{x}X_{\lambda_4}$ d'o\ù 3iii).

  \fdem
\end{enumerate}

\noi Remarques : \\

1) Lorsque $f$ est isotrope (donc  $d$ est pair et $c(f)=-1$), les coefficients de l'\équation fonctionnelle de la fonction Z\éta associ\ée \à l'action de $P$ sont ceux de la proposition 8.2.8:
 $$
\left(\begin{array}{ll}Z^*_{u_1}(\four (f);\omega,s)\\
...\ ...\\
Z^*_{u_\ell}\ (\four (f);\omega,s)\end{array}\right)=   B(\omega,s)\left(\begin{array}{ll}Z_{ u_1}(f; \omega^*|\ |^{-(\frac{3d}{2}+2)1_2})\\
...\ ...\\
Z_{u_\ell}\ ( f; \omega^*|\ |^{-(\frac{3d}{2}+2)1_2})\end{array}\right)$$
et  $B\bigg((id,\omega),(0,s)\bigg)$ est la matrice des coefficients de l'\équation fonctionnelle de la fonction Z\éta associ\ée \à l'action de $G,$ ce qui donne dans le cas r\éel:

$$B(0,s)= \frac{C_d(0,s)}{2}\left(\begin{array}{lll}(-1)^{\frac{d}{2}+1}\sin 2\pi s&0\\
  2(1+i^d)\sin\pi s&(-1)^\frac{d}{2}\sin 2\pi s\end{array}\right)\ ,\ (c(f)=-1)\ ,$$
  ce qui n'est pas en accord avec le r\ésultat obtenu par M.Muro pour  le coefficient de la $2$\ème ligne et $2$\ème colonne (\cite{muro1}: M.Muro obtient la moiti\é de notre coefficient).\\

2) Lorsque $d=1$ dans le cas $\goth p$-adique, $\chi_1(P)=\chi (G)=\F^{*2},$ on peut montrer que 
 chaque $O_u,u\in (\FF)^2,$ est exactement une $P$-orbite dans $\goth g_1''$ lorsque $u=(u_1,u_2)$ avec $u_1\not =u_2$ ou bien  $u_1=u_2=-1,$ en utilisant la multiplication par $\F^*$ et la propri\ét\é (*) de la \demo  de 8.2.12 puisque Ker$\tilde\omega_{-u_1}\cup$Ker$\tilde\omega_{-u_2}=\F^*$. De la m\ême mani\ère, on montre que  pour $u \in D$ tel que $-u\not\in \F^{*2},$   on a:
 $$O_{(u,u)} =P( X_{\lambda_1} +uX_{\lambda_2}+ X_{\lambda_3}+uX_{\lambda_4})\cup P( vX_{\lambda_1} +uvX_{\lambda_2}+ X_{\lambda_3}+uX_{\lambda_4}) $$  
 avec  $v\in D $ fix\é  tel que $(v,-u)=-1,$ mais je ne sais pas si $O_{(u,u)}$ est une seule $P-$orbite ou bien la r\éunion de $2$ $P-$orbites.\\

3) Dans le cas $\goth p$-adique, lorsque $d=2$ et $f$ est anisotrope, on a pour $u\in \FF$ 
$$P.X_{u,\pm 1}=\cup_{\{u_1\in \FF\ |\ (u_1,-c(f))=\pm 1\} }
O_{\{(u_1,u_1u)} \ ,\ P.Y_{u,\pm 1}=\cup_{\{u_2\in \FF\ |\ (u_2,-c(f))=\pm1\} }O^*_{(u_2u,u_2)}$$ avec $Y_{u, 1}=u X_{-\lambda_1}+\sum_{1\≤i\≤4}X_{-\lambda_i}$ et $Y_{u, -1}=\frac{u}{v_0} X_{-\lambda_1}+\sum_{ i=2,3}X_{-\lambda_i}+ v_0X_{-\lambda_4}.$ \\

\noi Il r\ésulte de la proposition 8.2.6 que les coefficients de l'\équation fonctionnelle de la fonction Z\éta associ\ée \à l'action de $P$ sont donn\és par:
$$a_{P.Y_{v,\epsilon'},P.X_{u,\epsilon}}=\sum_{\{w\in \FF|(w,-c(f))=\epsilon \epsilon'\}}d_{(v,1),(w,wu)}.$$

4) Dans le cas $\goth p$-adique,  pour $u$ et $v$ dans $\FF$ on a donc:
$$ a_{v,u}(\pi):=a_{G.X_v^{-1},G.X_u }( \pi)=\sum_{\{w\in \FF\}}d_{(w,wv),(1,u)}(Id,\pi).$$
 Pour $d$ pair, $B(id,\pi)$ est la matrice des coefficients de l'\équation fonctionnelle de la fonction Z\éta associ\ée \à l'action de $G,$ les r\ésultats obtenus ne semblent pas simples.\\ 
 
 Pour $d=1,$ en caract\éristique r\ésiduelle diff\érente de $2,$ on a simplement pour $a\in \FF$ et $s\in \C:$
 $$a_{v,u}(\tilde\omega_a|\ |^s)= (a,uv)\  \Frac{1-q^{2s+3}}{1-q^{-2s-4}}\ A^1_{s+1,s+\frac{7}{2}}(-v,-u,1),$$
 (lemme 8.2.10) les coefficients  sont proportionnels \à ceux associ\és \à une forme quadratique de discriminant $1.$\\

\bigskip 
 \subsubsection{Le cas r\éel de rang $4$ $(\∆=F_4$ et $f$ est anisotrope)}
 \bigskip 
 
 \noi  On applique la proposition 5.1.1 dans laquelle apparait les coefficients associ\és aux fonctions Z\éta d'une forme quadratique  qui est soit de signature $(2,d+1),$ soit de signature $(1,d+2)$ (3) de la \propo 8.2.4), cette forme quadratique, not\ée $Q$ dans ce rappel,  \étant l'\irf d'un \PV \air de type commutatif. \\
 
 \noi  On rappelle les r\ésultats bien connus suivants (cf.par exemple \cite{boppruben}),  $G_\R$ \étant la composante connexe r\éelle de $G$:  \\
  
$\bullet$ Lorsque la signature est $(2,d+1),$ on fait op\érer $G=\R^*\times G_\R,$ il y a donc $2$ orbites de $G$ dans $\goth g'_{\pm 1}:$ $\Omega_1=\{x\in \goth g_1|Q(x)>0\}$ et $\Omega_{-1}=\{x\in \goth g_1|Q(x)<0\}$ (idem ${\Omega^*}_{\pm 1}$), les coefficients correspondants sont not\és  $C(s)a_{\pm 1,\pm 1}(s)$, le premier indice correspond \à $\Omega_{\pm 1}^*$  et le second \à  $\Omega_{\pm 1},$ avec:
 
$$\biggl(\begin{array}{cc}a_{ 1, 1}(s)&a_{ 1, -1}(s)\\
a_{ -1, 1}(s)&a_{ -1, -1}(s)\end{array}\biggr)=
\biggl(\begin{array}{llc} -&\cos \pi (s+\frac{d}{2})&0  \\
 & \cos \pi \frac{d}{2} & \sin \pi s\end{array}\biggr)\ .$$
 (th.3.6.5,5)\\

$\bullet$ Lorsque la signature est $(1,d+2),$ on fait op\érer $ G_\R,$ il y a donc $3$ orbites de $G_\R$ dans $\goth g'_1$ not\ées $\Omega_i$  (idem dans $\goth g'_{-1}:$ ${\Omega_i}^* )$ avec $i=0,1,2$ et:

\noi  $\{x\in \goth g_1|Q(x)>0\}=\Omega_0\cup \Omega_2$ et $\Omega_1=\{x\in \goth g_1|Q(x)<0\}$ (idem pour ${\Omega_i}^*  $), 

\noi les coefficients correspondants sont not\és  $C(s)\epsilon_{ i,j}(s)$, le premier indice correspond \à ${\Omega_i}^*$   et le second \à $\Omega_j,$  avec:
 $$\Biggl(\begin{array}{lll} \epsilon_{ 0,0}(s)& \epsilon_{ 0,1}(s)&\epsilon_{ 0,2}(s)\\
 \\
 \epsilon_{ 1,0}(s)& \epsilon_{ 1,1}(s)&\epsilon_{ 1,2}(s)\\
 \\
  \epsilon_{ 2,0}(s)& \epsilon_{ 2,1}(s)&\epsilon_{ 2,2}(s)\end{array}\Biggr)=
\Biggl(\begin{array}{cccc} & \epsilon(s+\frac{d+3}{2})& \frac{1}{2}& {\epsilon(-s-\frac{d+3}{2})}\\
\\
 & -\sin \pi  \frac{d}{2} &-\cos \pi  s &  -\sin \pi  \frac{d}{2} \\
 \\
 &  \epsilon(-s-\frac{d+3}{2}) & \frac{1}{2}&\epsilon(s+\frac{d+3}{2})    \end{array}\Biggr)\ ,\  
 \epsilon(s)=\frac{1}{2}  e^{i\pi s }\ .$$
 (cf.par exemple la remarque 5.42 p.134 avec $k+1=2$ de \cite{boppruben})\\
 
Motiv\ées par la description des orbites ci-dessus et le mode de calcul des coefficients donn\é dans la proposition 5.1.1, on note  les  orbites de $P$ dans $\goth g"_1$ et $\goth g"_{-1}$  de la mani\ère suivante:
 
\noi  $\Omega_{i,j}$ et $\Omega^*_{i,j}$ avec $j\in \{0,1,2\}$ pour $i=0 $ et $j\in\{0,1\}$ pour $i=1.$ Plus pr\écis\ément dans $\goth g"_1$:
  $$\Omega_{0,1}=P.X_3=O_{(1,-1)}\ , \  \Omega_{1,0}=P.X_1=O_{(-1,1)}\ , \ \Omega_{1,1}=P.X_4=O_{(-1,-1)}.$$
$$\Omega_{0,0}=P.X_0 \ ,\ \Omega_{0,2}=P.X_2\ , \ O_{(1,1)}=\Omega_{0,0}\cup \Omega_{0,2}.$$
\noi Pour $\goth g"_{-1},$ soient:\\

\noi  $Y_0=\sum_{1\leq i\leq 4}X_{-\lambda_i}\ ,\ Y_1=-X_{-\lambda_1}+\sum_{2\leq i\leq 4}X_{-\lambda_i}\ ,\ Y_2=-X_{-\lambda_1}-X_{-\lambda_2}+X_{-\lambda_3}+X_{-\lambda_4},$
$Y_3=X_{-\lambda_1}+X_{-\lambda_2}-X_{-\lambda_3}+ X_{-\lambda_4}$ et $Y_4=-X_{-\lambda_1}+X_{-\lambda_2}-X_{-\lambda_3}+ X_{-\lambda_4},$ \\

\noi alors:

$$ \Omega^*_{0,1}=P.Y_1=O^*_{(-1,1)}\ \  , \ \ \Omega^*_{1,0}=P.Y_3=O^*_{(1,-1)}\ \ , \ \ \Omega^*_{1,1}=P.Y_4=O^*_{(-1,-1)}\ ,$$
$$\Omega^*_{0,0}=P.Y_0\ \ , \ \ \Omega^*_{0,2}=P.Y_2\ \ , \ \ O^*_{(1,1)}=\Omega^*_{0,0}\cup \Omega^*_{0,2},$$
c'est \à dire que  $\Omega_{i,j}\subset O_{((-1)^i,(-1)^j)}$ et $\Omega^*_{i,j}\subset O^*_{((-1)^j,(-1)^i)}.$\\ \\

 \begin{prop}
$$a_{\Omega^*_{j,k},\Omega_{q,\ell}}(s_1,s_2)=C_d(s_1,s_2)A_{\Omega_{q,\ell}}^
{\Omega^*_{j,k}} (s_1,s_2)\ ,$$   $A_{\Omega_{q,\ell}}^
{\Omega^*_{j,k}} (s_1,s_2)$ ayant les valeurs suivantes:\\

$$ A_{\Omega_{1,1}}^
{\Omega^*_{j,k}} (s_1,s_2)\ =\biggl \{ \begin{array}{ll} (-1)^d\sin\pi s_2\cos\pi (s_1+s_2)&\hbox{pour}\ k=j=1\\
0&\hbox{ sinon} \end{array}   $$

$$ A_{\Omega_{q,\ell}}^
{\Omega^*_{1,1}} (s_1,s_2)=\ \biggl \{ \begin{array}{ll}  \cos\pi\frac{d}{2} \sin\pi (s_1+s_2) &\hbox{pour}\  (q,l)=(0,1)\\
 \cos\pi\frac{d}{2} \sin\pi s_2&\hbox{pour}\  (q,l)=(1,0)\\
0&\hbox{pour} \  (q,l)=(0,0)\ \hbox{ou}\  (q,l)=(0,2).\end{array} $$ 

\noi Pour les autres valeurs:\\

\noi\begin{tabular}{|c|c|c|c|c|}
\hline
&&&&\\
$ A_{\Omega_{i,j}}^
{\Omega^*_{k,l}}$  &$\Omega^*_{0,0}$& $\Omega^*_{0,2}$  & $\Omega^*_{0,1}$  &   $\Omega^*_{1,0}$  \\
&&&&\\
\hline
&&&&\\
 $\Omega_{0,0}$ &$\tiny{\frac{(-1)^d}{2}\sin\pi(s_1+2s_2)}$ &$\tiny{-\frac{ 1}{2}\sin\pi s_1}$  & $\tiny{\frac{1-(-1)^d}{2}\sin\pi(s_1+s_2)}$ & $\tiny{\frac{(-1)^d-1}{2}\sin\pi s_2}$\\
&&&&\\
\hline
&&&&\\
$\Omega_{0,2}$&  $\tiny{-\frac{1}{2}\sin\pi s_1}$  &$\tiny{\frac{ (-1)^d}{2}\sin\pi (s_1+2s_2)}$   & $\tiny{\frac{1-(-1)^d}{2}\sin\pi(s_1+s_2)}$ &  $\tiny{\frac{(-1)^d-1}{2}\sin\pi s_2}$\\
&&&&\\
\hline
&&&&\\
$\Omega_{0,1}$ &$\tiny{\frac{(-1)^d}{2}\sin\pi(s_2+\frac{d}{2 })}$ &$\tiny{\frac{(- 1)^d}{2}\sin\pi (s_2+\frac{d}{2 })}$  & $\tiny{\frac{1-(-1)^d}{2}}$   &$\tiny{(-1)^{d+1} \cos\pi (s_2+\frac{d}{2 }) }$\\
&&&&$.\ \tiny{\sin\pi (s_1+s_2)}$\\

&&&&\\
\hline
&&&&\\
 $\Omega_{1,0}$ &  $\tiny{ \frac{(-1)^d }{2} } $ & $\tiny{\frac{ (-1)^d}{2} \sin\pi(s_1+s_2+\frac{d}{2 })}$ &$\tiny{ (-1)^{d +1} \cos\pi s_2   }$ &$0$\\
 &$.\ \tiny{\sin\pi(s_1+s_2+\frac{d}{2 })}$&&$ .\ \tiny{\sin\pi(s_1+s_2+\frac{d}{2 })}$&\\
&&&&\\
\hline
\end{tabular}

\end{prop} 
\bigskip
\dem On reprend les notations de la proposition 5.1.1 pour calculer le coefficient $a_{\Omega_{q,\ell}, \Omega^*_{j,k}}$ avec $\ell,k\in\{0,1,2\}$ lorsque $q,j=0$  et   $\ell,k\in\{0,1\}$ lorsque $q,j=1.$ \\

\noi Notons:
$$t_0=X_{\lambda_1} +X_{\lambda_2}\ ,\ t_1=-X_{-\lambda_1} +X_{\lambda_2}\ ,\ t_2=-X_{\lambda_1} -X_{\lambda_2}\ ,\ $$
 $$u_0=X_{\lambda_3} +X_{\lambda_4}\ ,\ u_1=-X_{\lambda_3} +X_{\lambda_4}\ ,\ u_2=-X_{\lambda_3} -X_{\lambda_4}.$$\\

\noi  Soit  $z=z_{-2}+z_0\in \Omega^*_{j,k} \cap \{Y_i,i=0,...,3\}$, alors $z_{-2}=t_k^{-1}$ et $z_0=u_j^{-1}.$ On a 2 cas:\\

$\bullet$   $j=1$ donc $F_1^*(z_0)<0,$ alors $\gamma_1(t_i,u_1)=1$ et $(G_{H_1})_{u_1}=\R^*\times {((G_{H_1})_{u_1})}_\R$  (puisque $c(-1)g_{3,4}(1)\in (G_{H_1})_{u_1})$ donc: $$a_{z_{-2},t_i}^{(u_1)}(s_2)=C_d(s_2)a_{(-1)^k,(-1)^q}(s_2)\ .$$
 $I_{O,O^*}=\{q\}$ et $t_q(\Omega_{q,\ell})=
u_\ell.$ \\

$\bullet$ Soit $j=0$ donc $F_1^*(z_0)>0,$  $z_0=u_0^{-1}$ et $(G_{H_1})_{z_0}= {((G_{H_1})_{z_0})}_\R$ (puisque $t_0+u_0$ et $t_2+u_0$ ne sont pas dans la m\ême orbite de $G$) d'o\ù 
$$a_{z_{-2},t_i}^{(u_0)}(s_2)=C_d(s_2)\epsilon_{k ,i}(s_2).$$  Lorsque $q=1,$ on a $I_{O,O^*}=\{1\}$, $t_1(\Omega_{1,\ell})=u_\ell$ et 
$\gamma_1(t_1,z_0)= 1.$

\noi  Lorsque $q=0$ on a $\gamma_1(t_0,z_0)= \gamma_1(t_2,z_0)=(- 1)^d$ et $I_{O,O^*}=\{0,2\};$ $t_0(\Omega_{0,\ell})=u_\ell $ et $$t_2(\Omega_{0,\ell})=  \bigg\{ \begin{array}{ll}&u_2\ \hbox{lorsque}\ \ell=0\\
&u_0\ \hbox{lorsque}\ \ell=2\\
&u_1\ \hbox{lorsque}\ \ell=1\end{array}\quad =\quad u_{2-\ell}.$$
Les r\ésultats sont analogues pour le calcul du coefficient $a_{z_0,t_i(\Omega_{q,\ell})}^{( t_i)}(s_1+s_2+d+\frac{1}{2})$ suivant que:\\
 
$\bullet$  $q=1$   alors  il vaut $ C_d(s_1+s_2+d+\frac{1}{2})a_{(-1)^j,(-1)^\ell}(s_1+s_2+d+\frac{1}{2})$\\

\noi ou que\\ 

$\bullet$   $q=0$  alors le coefficient  vaut $C_d(s_1+s_2+d+\frac{1}{2})\epsilon_{j, \ell' }(s_1+s_2+d+\frac{1}{2}),$ avec $\ell'=\ell$ si $i=0,1$ et $\ell'=2-\ell$ si $i=2,$\\

\noi D'o\ù
$a_{\Omega^*_{j,k},\Omega_{q,\ell}}(s_1,s_2)=C_d(s_1,s_2)A_{\Omega_{q,\ell}}^
{\Omega^*_{j,k}} (s_1,s_2)$ avec:\\

\noi $A_{\Omega_{1,\ell}}^
{\Omega^*_{1,k}} (s_1,s_2)= a_{ (-1)^k,-1}(s_2)a_{-1,(-1)^\ell} (s_1+s_2+d+\frac{1}{2}) ,$\\

\noi $A_{\Omega_{0,\ell}}^
{\Omega^*_{1,k}} (s_1,s_2)=a_{(-1)^k,1}(s_2)\epsilon_{1, \ell}(s_1+s_2+d+\frac{1}{2})   ,$\\

\noi $A_{\Omega_{1,\ell}}^
{\Omega^*_{0,k}} (s_1,s_2)= \epsilon_{k,1}(s_2)a_{1,(-1)^\ell}(s_1+s_2+d+\frac{1}{2}) ,$\\

\noi $A_{\Omega_{0,\ell}}^
{\Omega^*_{0,k}} (s_1,s_2)=(-1)^d\biggl(\epsilon_{k,0}(s_2) \epsilon_{ 0,\ell}(s_1+s_2+d+\frac{1}{2})  +\epsilon_{k,2}(s_2)  \epsilon_{ 0,2-\ell}(s_1+s_2+d+\frac{1}{2}) \biggr).$
\fdem
 
\vskip 5mm

\noi Remarque:  Notons: $O_i=G.X_i$ et $O^*_i=G.Y_i$ pour $i=0,1,2,$ les orbites de $G$ dans $\g'_1$ et $\g'_{-1}$ alors:\\

\noi $\Omega_{0,0}=O_0\cap \g"_1\ ,\Omega^*_{0,0}=O^*_0\cap \g"_{-1},$

\noi $O_2\cap \g"_1=\Omega_{0,2} \cup \Omega_{1,1}\ ,\ O^*_2\cap \g"_{-1}=\Omega^*_{0,2} \cup \Omega^*_{1,1},$

\noi $O_1\cap \g"_1=\Omega_{0,1} \cup \Omega_{1,0}\ ,\ O^*_1\cap \g"_{-1}=\Omega^*_{0,1} \cup \Omega^*_{1,0},$\\

\noi donc $a_{O^*_i,O_j}(s)=C_d(0,s)A_{i,j}(s)\ ,\ i,j=0,1,2,$ avec:\\

\noi $A_{0,i}(s)=  A_{\Omega_{0,i}}^{\Omega^*_{0,0}}(0,s)$   $(=
A_{\Omega_{1,i-1}}^{\Omega^*_{0,0}}(0,s)\hbox{ si }\ i\≥1)$\\

 \noi $A_{1,i}(s)=  A_{\Omega_{0,i}}^{\Omega^*_{1,0}}(0,s) +A_{\Omega_{0,i}}^{\Omega^*_{0,1}}(0,s) $  $(=
 A_{\Omega_{1,i-1}}^{\Omega^*_{1,0}}(0,s) +A_{\Omega_{1,i-1}}^{\Omega^*_{0,1}}(0,s)\ \hbox{ si }\ i\≥1)$\\

\noi $ A_{2,i}(s)=
   A_{\Omega_{0,i}}^{\Omega^*_{0,2}}(0,s) +A_{\Omega_{0,i}}^{\Omega^*_{1,1}}(0,s) $ $(=
 A_{\Omega_{1,i-1}}^{\Omega^*_{0,2}}(0,s) +A_{\Omega_{1,i-1}}^{\Omega^*_{1,1}}(0,s)$   si $\ i\≥1)$\\

\noi d'o\ù:
 $$\left(\begin{array}{llll}
  A_{2,2}( s)&A_{1,2}(s)&A_{0,2}(s)\\
   A_{2,1}( s)&A_{1,1}(s)&A_{0,1}(s)\\
   A_{2,0}( s)&A_{1,0}(s)&A_{0,0}(s)\end{array}\right)=$$
$$   \left(\begin{array}{llll}\frac{(-1)^d}{2}\sin2\pi s&  0 &0\\
   \frac{(-1)^{d}}{2} \sin\pi (s +\frac{d}{2})+\cos \frac{\pi d}{2}\sin\pi s&(-1)^{d+1} \cos\pi s \sin\pi (s +\frac{d}{2})&\frac{(-1)^d}{2}\sin\pi (s +\frac{d}{2})\\
   0& 0&\frac{(-1)^d}{2}\sin2\pi s\end{array}\right)$$
ce qui n'est pas en accord avec les r\ésultats obtenu par M.Muro pour $A_{2,1}(s)$ et $A_{0,1}(s)$ (de signe oppos\é dans (\cite{muro1}).\\

\noi Comme $A_{1,0}(s)=A_{1,2}(s)=0$ et que $A_{0,0}(s)+A_{2,0}(s)=A_{0,2}(s)+A_{2,2}(s),$   les fonctions Z\étas:\\

 $\bullet \ Z_u(f;\pi)=Z(f1_{\{x\in \g_1|uF(x) >0\}};\pi)$ pour $f\in \EuScript S(\g_1),$\\
 
$\bullet \ Z^*_u(g;\pi)=Z^*(g1_{\{x\in \g_{-1}|uF^*(x) >0\}};\pi)$ pour $g\in \EuScript S(\g_{-1}),$\\

\noi  avec $u=\pm 1, $ v\érifient l'\équation fonctionnelle suivante pour $f\in \EuScript S(\g_1):$
$$\begin{array}{ll}
\left(\begin{array}{ll}Z^*_{-1}(\four (f);s)\\
Z^*_{1}\ (\four (f);s)\end{array}\right)&= 
C_d(0,s)\left(\begin{array}{lll}
(-1)^{d+1} \cos\pi s \sin\pi (s +\frac{d}{2})&0\\
  (-1)^d\sin\pi (s +\frac{d}{2})+\cos \frac{\pi d}{2}\sin\pi s &\frac{(-1)^d}{2}\sin2\pi s\end{array}\right).\\
  \\
  &.\left(\begin{array}{ll}Z_{-1}(f;-s-2-\frac{3d}{2})\\
Z_{1}\ ( f;-s-2-\frac{3d}{2})\end{array}\right)\end{array}  
 $$
ce qui est en accord dans le cas $d$ pair avec la proposition 8.2.8 lorsque $s_1=0$ et $s_2=s$ (cf.remarque 8.2.9) et avec le lemme 8.2.10 pour $d=1.$\\ 
  
  \bigskip
  Dans les cas exceptionnels restants, on se contentera de donner les r\ésultats minimaux sans entrer dans le d\étail des $P$-orbites.
   \bigskip
   
\subsection {$\bf (E_7,\alpha_2)$}
\bigskip
On note $\{\lambda_1,...,\lambda_7\}$ les $7$ racines orthogonales obtenues canoniquement par orthogonalisations successives \à partir de $\lambda_1=\alpha_2$ (cf.\cite{mullerJA2}, tableau II pour une liste explicite et tableau III pour une description de $\Delta_2$),   soient $H_1:=\sum_{1\≤i\≤3}h_{\lambda_i}$ et $H_2=2H_0-H_1,$ on a:
$$\sum_{1\≤i\≤7}h_{\lambda_i}=2H_0\ , 
  H_2=2h_{\widetilde\alpha}\ ,\  \widetilde\alpha\ \text{ \étant la plus grande racine de }\Delta.$$
$P=P(H_1,H_2)$   est alors l'unique \sg parabolique standard tr\ès sp\écial du \PV  et il est associ\é \à la partie $\∑_0-\{\alpha_1\}$ de $\∆_0$ qui est de type $A_5.$\\

On fixe dor\énavant un syst\ème de Chevalley:   $(X_{\alpha},h_{\alpha},X_{-\alpha})_{\alpha\in \∆}$, de $ \goth g$ tel que pour toute racine $\omega\in \∆_2$ on ait:
$$ [X_{-\lambda_i},[X_{-\lambda_j},[X_{-\lambda_k},[X_{-\lambda_l},X_\omega]]]]= X_{-\omega}\ , \
\hbox{avec}\  \omega=\frac{1}{2}(\lambda_i+\lambda_j+\lambda_k+\lambda_l)\ ,$$
dont l'existence est assur\ée par le corollaire 4.3 de \cite{mullerJA2}.\\

$F$ est l'invariant relatif fondamental du \PV normalis\é par la condition : $$F(\sum_{1\leq i\leq 7}X_{\lambda_i})=1 \ \hbox{et on a pour}\ a_1,...,a_7\in \F^*\ :\quad F(\sum_{1\leq i\leq 7}a_iX_{\lambda_i})=\prod_{1\leq i\leq 8}a_i.$$

On normalise $B$ en posant  $\tilde B=\displaystyle{-\frac{7 B}{2B(H_0,H_0)}}(=\displaystyle{-\frac{2B}{B(h_{\lambda_i},h_{\lambda_i})}}),$ ainsi $\tilde B(\frac{H_1}{2},\frac{H_1}{2})=-\ds\frac{3}{2},$ $\tilde B(\frac{H_2}{2},\frac{H_2}{2})=- 2$ et pour toute racine longue $\alpha$ de $\∆$ on a $\tilde B(X_{\alpha}, X_{-\alpha})=1.$\\ 

\noi Rappelons que $E_1(H_1)\cap \g_1=\{0\},$ que $\g_2$ est de dimension $7$ et que $E_2(H_1)\cap \g_2$ est de dimension $6$ (cf.tableau III de \cite{mullerJA2}).\\

\begin{prop}
\begin{enumerate}
\item Le \PV $(E_0(H_1)\cap \goth g_0,E_2(H_1)\cap \goth g_1)$
est commutatif de type $(D_6,\alpha_6)$ et d'\irf $F_1$ normalis\é par la condition $F_1(\sum_{1\≤i\≤3}a_iX_{\lambda_i})=a_1a_2a_3.$ 

\item Pour $x\in \goth g_1$ on a alors $P_1(x)=\frac{1}{4!}\tilde B(adx^4(X_{-\omega},X_{-\omega})$ avec $\omega=\ds\frac{1}{2}
\sum_{4\≤i\≤7} \lambda_i.$

\item Soit $x\in W_{\goth t},$ $x=x_2+x_0$ avec $x_i\in E_i(  H_1 )\cap \goth g_1$ alors:\\
 
 a) Le \PV $((E_0(H_1)\cap \goth g_0)_{x_0},(E_2(H_1)\cap \goth g_1)_{x_0})$ est commutatif presque d\éploy\é de type:

$$(\∆,\lambda_0)=\bigg\{\begin{array}{l} (A_5,\alpha_3)\ \hbox{lorsque}\ -P_1(x_0)\in \F^{*2}\\
(C_3,\alpha_3)\ \hbox{sinon}\ ,\end{array}$$
et $f\sim X^2\oplus P_1(x_0)Y^2.$\\

 b) Le \PV $((E_0(H_2)\cap \goth g_0)_{x_2},(E_2(H_2)\cap \goth g_1)_{x_2})$ est de type $(F_4,\alpha_1).$
 \end{enumerate}
\end{prop}

\dem 1) Cela   r\ésulte du diagramme de Dynkin compl\ét\é de $E_7.$\

2) R\ésulte des valeurs:
 $$F(\sum_{1\≤i\≤7}a_iX_{\lambda_i})=\prod_{1\leq i\leq 8}a_i,F_1(\sum_{1\≤i\≤3}a_iX_{\lambda_i})=\prod_{1\leq i\leq 3}a_i,P_1(\sum_{4\≤i\≤7}a_iX_{\lambda_i})=\prod_{4\leq i\leq 8}a_i.$$

3) Pour tout $x\in W_{\goth t},$ il existe $g\in G_{H_1}$ et $(a_1,...,a_7)\in \F^*$ tels que $gx=\sum_{1\≤i\≤7}a_iX_{\lambda_i}$ (lemme 7.3 de \cite{mullerJA2}) donc on fait la \demo pour $x_2=\sum_{1\≤i\≤3}a_iX_{\lambda_i}$ et $x_0=\sum_{4\≤i\≤7}a_iX_{\lambda_i}.$\\

a) Le \PV $((E_0(H_1)\cap \goth g_0)_{x_0},(E_2(H_1)\cap \goth g_1)_{x_0})$ est un \PV \air  commutatif.  

\noi Posons: 
$$\omega_{1,2}=\ds\frac{1}{2}(\lambda_1 +\lambda_2 +{\lambda_4}+ {\lambda_7})\ ,\ \omega'_{1,2}=\ds\frac{1}{2}( {\lambda_1}+ {\lambda_2}+ {\lambda_5}+ {\lambda_6}),$$
$$\omega_{1,3}=\ds\frac{1}{2}( {\lambda_1}+ {\lambda_3}+ {\lambda_4}+ {\lambda_6})\ ,\ \omega'_{1,3}=\ds\frac{1}{2}( {\lambda_1}+ {\lambda_3}+ {\lambda_5}+ {\lambda_7}),$$
$$\omega_{2,3}=\ds\frac{1}{2}( {\lambda_2}+ {\lambda_3}+ {\lambda_4}+{\lambda_5})\ ,\omega'_{2,3}=\ds\frac{1}{2}( {\lambda_2}+ {\lambda_3}+ {\lambda_6}+ {\lambda_7}).$$
Soient:
$$ X_{1,2}=a_{7}[X_{-\lambda_{4}},X_{\omega_{1,2}}] -a_{4}[X_{-\lambda_{7}},X_{\omega_{1,2}}]\ ,\   
 X'_{1,2}=a_{6}[X_{-\lambda_{5}},X_{\omega'_{1,2}}] -a_{5}[X_{-\lambda_{6}},X_{\omega'_{1,2}}]\ ,  $$
$$X_{1,3}=a_{6}[X_{-\lambda_{4}},X_{\omega_{1,3}}] -a_{4}[X_{-\lambda_{6}},X_{\omega_{1,3}}]\ ,\  X'_{1,3}=a_{7}[X_{-\lambda_{5}},X_{\omega'_{1,3}}] -a_{5}[X_{-\lambda_{7}},X_{\omega'_{1,3}}]\ ,$$ 
$$X_{2,3}=a_{5}[X_{-\lambda_{4}},X_{\omega_{2,3}}] -a_{4}[X_{-\lambda_{5}},X_{\omega_{2,3}}]\ ,\ X'_{2,3}=a_{7}[X_{-\lambda_{6}},X_{\omega'_{2,3}}] -a_{6}[X_{-\lambda_{7}},X_{\omega'_{2,3}}]\ , $$  alors $(E_2(H_1)\cap \goth g_1)_{x_0}=\oplus_{1\≤i\≤3}\goth g^{\lambda_i}\oplus_{1\≤i<j\≤3}(\F X_{i,j}\oplus \F  X'_{i,j})  
$
donc le \PV est presque d\éploy\é (i.e. $d_1=1$) avec $d=2$ d'o\ù par classification, ce \PV est de type $(A_5,\alpha_3)$ lorsque la forme quadratique $f'$ ci-dessous est isotrope et sinon il est de type 
$(C_3,\alpha_3)$ (cf.tableau 1).

\noi Soit $f'$ la forme quadratique:
  $$\begin{array}{rl}f'(u,v)&=\frac{1}{2}\tilde B(ad(uX_{1,2}+vX'_{1,2})^2(X_{-\lambda_1}),X_{-\lambda_2})\\
  &=-\frac{1}{2}\tilde B( [uX_{1,2}+vX'_{1,2},X_{-\lambda_1}],[uX_{1,2}+vX'_{1,2},X_{-\lambda_2}])\\
  &=u^2a_4a_7+v^2a_5a_6\ ,\end{array}$$
$f'$ est isotrope \ssi $P_1(x_0)\equiv -1$ (mod $\F^{*2}$).\\

b) Comme dans a) on obtient: $$(E_4(H_2)\cap \goth g_2)_{x_2}=\g^{\widetilde\alpha}\ ,\ (E_2(H_2)\cap \goth g_1)_{x_2}=\oplus_{4\≤i\≤7}\goth g^{\lambda_i}\oplus_{1\≤i<j\≤3}(\F Y_{i,j}\oplus \F  Y'_{i,j}) \oplus_{4\≤i\≤7}\F [X_{-\lambda_i},X_{\widetilde\alpha}] 
$$ avec 
$$  Y_{i,j}=a_{j}[X_{-\lambda_{i}},X_{\omega_{i,j}}] -a_{i}[X_{-\lambda_{j}},X_{\omega_{i,j}}]\ ,\
 Y'_{i,j}=a_{j}[X_{-\lambda_{i}},X_{\omega'_{i,j}}] -a_{i}[X_{-\lambda_{j}},X_{\omega'_{i,j}}]\ .$$ 
En raison des dimensions de $(E_4(H_2)\cap \goth g_2)_{x_2}$ et de $(E_2(H_2)\cap \goth g_1)_{x_2},$ le \PV $((E_0(H_2)\cap \goth g_0)_{x_2},(E_2(H_2)\cap \goth g_1)_{x_2})$ est de type $(F_4,\alpha_1)$ (cf. $\S 7.1$ et $8.2$) et l'invariant associ\é est la restriction de $P_1$ \à $E_2(H_2)\cap \goth g_1)_{x_2}.$\fdem

\bigskip

Lorsque $\F$ est archim\'edien, les $2$ polynomes de Bernstein s'obtiennent imm\'ediatement en appliquant les  propositions  3.4.4 et  3.7.3 (ainsi que le lemme 6.1.5 et la proposition 8.2.5)   d'o\`u:
  
  \begin{prop}  Les polynomes de Bernstein du \PV de type $(E_7,\alpha_2)$  sont donn\és par: $$\begin{array}{ll}b_1(s_1,s_2)&=  s_2(s_2+1)(s_2+2)\ \hbox{et}\\ b_2(s_1,s_2)&=s_2(s_2+1)(s_2+2) (s_1+s_2+\frac{3}{2}) (s_1+s_2+\frac{5}{2}) (s_1+s_2+ 3)(s_1+s_2+4).\end{array}$$
  \end{prop}

  Remarque: $b_2(0,s)$ est le polynome de Bernstein usuel.
  \bigskip

\noi  On rappelle que pour $u=(u_1,u_2)\in (\FF)^2$ on d\'efinit les ouverts (tous non vides):\\

$O_u=O_{(u_1,u_2)}=\{x\in \goth
g_1\ |\ F_1(x) \F^{*2}= u_1\ \ ,\ F_2(x) \F^{*2}= u_1u_2 )
  \}$\\
 
$O^*_u=O^*_{(u_1, u_2)}=\{x\in
\goth g_{-1}\ |\ F^*_1(x) \F^{*2}= u_2 \   ,\ F^*_2(x) \F^{*2}=
 u_1u_2  \}$\\
 
\noi ainsi que les fonctions Z\'etas correspondantes:\\
  
$Z_u(f;\omega)=Z(f{\bf 1}_{O_u};\omega)$ pour $f\in \ES(\goth
g_1)$  
et
$ Z^*_u(h;\omega)=Z^*(h{\bf 1}_{O^*_u};\omega)$  pour $h\in \ES(\goth g_{-1})$.  \\  
 
\begin{prop} 
\begin{enumerate}

\item Lorsque $\F=\C$ pour $f\in \ES(\goth
g_1),$  on a:
$$ Z^*( \four (
f);(\omega_1,\omega_2),(s_1,s_2))=A( \omega_2,s_2) B(\omega_1\omega_2,s_1+s_2+\frac{3}{2})Z\biggl(f;( \omega_1,( \omega_1\omega_2)^{-1});s_1,-(s_1+s_2+ 5)\biggr)$$
avec:
$$A(\omega,s)=\rho'(\omega;s+1)\rho'(\omega;s+ 2)\rho'(\omega ;s + 3) $$
$$B(\omega,s)=\rho'(\omega;s+1)\rho'(\omega;s+ 2)\rho'(\omega ;s + \frac{5}{2}) \rho'(\omega ;s + \frac{7}{2}) .
$$
\item Dans le  cas r\éel ou bien lorsque $\F$ est un corps $\goth p-$adique de caract\éristique r\ésiduelle diff\érente de $2,$ on a pour tout  $v=(v_1, v_2)$ dans $(\FF)^2$ et $f\in \ES(\goth
g_1):$  
$$ Z^*_v( \four(
f); (s_1,s_2))=  |2|_\F^{-2(s_1+s_2)-\frac{13}{2}}\rho (|\ |^{2(s_1+s_2)+7})\sum_{u\in (\FF)^2}d_{v,u}(
  s)  Z_u(f; (  s_1,-(s_1+s_2+ 5)\ ) $$
avec $u=(u_1, u_2)$ et $v=(v_1,v_2)$ et:
 
\noi $d_{v,u}(
 s)= A^1_{s_1+s_2+ \frac{5}{2},s_1+s_2+5}(1,v_2u_2, -v_2)\frac{1}{f}\sum_{a\in \FF}(a,-u_1v_1)\prod_{\ell=1}^3\rho (\tilde\omega_a(\tilde\omega_{-v_2})^{\ell-1}|\ |^{s_2+\ell}) .$\end{enumerate}
\end{prop}
  
  \dem

\noi On applique la proposition 5.3.1 dont toutes les hypoth\`eses sont v\'erifi\'ees. En effet et dans les notations de la proposition 5.3.1, pour $z=z_0+z_{-2}\in O^*_{v_1,v_2}$ on a $$b_{v_1,u_1}^{(v_2)}( s)=\biggl(\alpha(1)\alpha(v_2) \biggr)^3(-v_2,-1)\frac{1}{f}\sum_{a\in \FF}(a,-u_1v_1)\prod_{\ell=1}^3\rho ( \tilde\omega_a(\tilde\omega_{-v_2})^{\ell-1}|\ |^{s +\ell})$$ en appliquant la proposition 8.3.1, le 2) des th\éor\èmes 6.2.2 et 6.1.2.1 ainsi que le lemme 3.5.3,ii)b) et $$c_{v_2,u_2}^{(w_1)}( s) =\alpha(-1)^2 |2|_\F^{-2s-\frac{7}{2}}\rho (|\ |^{2s+4})
 A^1_{s+1 ,s +\frac{7}{2}}(v_2,u_2, -1)$$ (prop.8.2.11 avec $d=1$).

\noi On termine en notant que: $\alpha(1)\alpha(v_2)^3(-v_2,-1)=\alpha(-1)\alpha(v_2)$ puis que 
$$ \alpha(-1)\alpha(v_2) (v_2,u_2, -1)=A^1_{s+1 ,s +\frac{7}{2}}(1,v_2u_2, -v_2) $$
(relation 3), lemme 3.6.4,B)).\fdem

\bigskip

\begin{prop}  Dans le cas r\éel ou le cas $\goth p$-adique de caract\éristique r\ésiduelle diff\érente de $2$ on a:

\begin{enumerate}
\item Soit $v\in \FF,$ alors:
$$ Z^*_{F_1^*F_2^*\F^{*2}=v}( \four(
f); (s_1,s_2))= K(s_1,s_2)\sum_{u\in \FF}D_{v,u}(
  s_1,s_2)  Z_{F_1F_2\F^{*2}=u}(f; (  s_1,-(s_1+s_2+ 5)\ ) $$
avec:

 $\begin{array}{lll}K(s_1,s_2)   & 
=&|2|_\F^{-2(s_1+s_2)-\frac{13}{2}}\rho(|\ |^{2(s_1+s_2)+7)}\rho (  |\ |^{s_2+ 1})\rho ( |\ |^{s_2+ 3})\ ,\end{array}$\\

$D_{v,u}(s_1,s_2)=   A^1_{s_1+s_2+ \frac{5}{2},s_1+s_2+5}(1,vu, -v)  \rho ( \tilde\omega_{-v} |\ |^{s_2+2}) .$\\

\item Soit $b\in \FF$ et $s\in \C,$ dans le \PV $(G,\g_{\pm 1})$ on a l'\eq:
$$ Z^*( \four(
f); (\tilde\omega_b,s))=C(s)B(\tilde\omega_b,s)Z (f; (\tilde \omega_b, -s - 5)\ ) $$
avec:
$$ C(s)= |2|_\F^{-4s-10}\rho ( |\ |^{2s+4})\rho ( |\ |^{2s+7})\ , $$
$$B(\tilde\omega_b,s)=  \rho' (\tilde \omega_b |\ |^{s+ 1})\rho' (\tilde\omega_b |\ |^{s+ 3})\rho' (\tilde\omega_b |\ |^{s+ 5}).$$
\end{enumerate}
\end{prop}

\dem 1) R\ésulte d'un simple calcul.\\

\noi 2) Calculons:
$$B_{u_1,u_2}(s)=(b,u_1u_2)\sum_{(v_1,v_2)\in (\FF)^2}d_{(v_1,v_2),(u_1,u_2)}(0,s)\ (b,v_1v_2).$$
En sommant selon $v_1$ on obtient:
$$B_{u_1,u_2}(s)=  \rho ( \tilde\omega_b |\ |^{s+ 1}) \rho (\tilde\omega_b |\ |^{s+ 3})B'_{u_1,u_2}(s)$$
avec:
$$B'_{u_1,u_2}(s)=(b, -u_2)\sum_{v_2\in \FF}A^1_{s +\frac{5}{2},s+5} (1, u_2v_2,-v_2)\rho ( \tilde\omega_{-bv_2} |\ |^{s+ 2})(b,v_2)$$
d'o\ù:
$$\begin{array}{ll}B'_{u_1,u_2}(s)&= \sum_{v_2,t\in \FF}(b,-u_2v_2)\alpha (t)  
\rho (  |\ |^{s+ \frac{5}{2}};t )
\rho (  |\ |^{s+ 5};tu_2v_2)(t,-v_2)\rho (\tilde\omega_{-bv_2} |\ |^{s+ 2}) \\
\\
&=\frac{1}{f}\sum_{z,v_2,t\in \FF}(z,tu_2v_2)\rho (\tilde\omega_z  |\ |^{s+ 5})(b,-u_2v_2)\alpha (t) (t,-v_2)
\rho (  |\ |^{s+ \frac{5}{2}};t )
\rho (\tilde\omega_{-bv_2} |\ |^{s+ 2})  \end{array}$$
 or:
 $$ \sum_{t\in \FF}\alpha (t) (t,-zv_2)\rho (  |\ |^{s+ \frac{5}{2}};t )=
 h(  |\ |^{s+ \frac{5}{2}}\tilde\omega_{-zv_2})$$
 (lemme 3.6.4,A)
puis:
$$\begin{array}{lll}f(z)&=&\frac{1}{f}\sum_{v_2\in \FF} h(  |\ |^{s+ \frac{5}{2}}\tilde\omega_{-zv_2})(v_2,bz)\rho (\tilde\omega_{-bv_2} |\ |^{s+ 2})\\
\\
&=&\frac{1}{f}\sum_{v_2\in \FF} h(  |\ |^{s+ \frac{5}{2}}\tilde\omega_{bzv_2})(-bv_2,bz)\rho (\tilde\omega_{v_2} |\ |^{s+ 2})\\
\\
&=&(-b,z)\ds\frac{\alpha (1)}{\alpha (zb)}A^1_{s + 2,s+\frac{5}{2}} (1, zb,1)
\end{array}$$
 (lemme 3.6.4,B)6) donc $f(z)=0$ pour $z\not=b$ (lemme 3.6.7,C) et finalement:
 $$B'_{u_1,u_2}(s)=A^1_{s + 2,s+\frac{5}{2}} (1, 1,1)\rho (\tilde\omega_b  |\ |^{s+ 5})(b,-1)
$$
On termine \à l'aide du lemme 3.6.8.\fdem

\begin{rema}
\end{rema}
1) Dans le cas r\éel, en n'utilisant pas la formule de duplication de Legendre,on a  \également:
$C(s)B(\tilde\omega_b,s)=D(s).D(b,s)$ avec $$D(s)=8.2\pi^{-(7s+21)}\Gamma (s+1)
\Gamma (s+2)\Gamma (s+\frac{5}{2})\Gamma (s+3)\Gamma (s+\frac{7}{2})\Gamma (s+4)\Gamma (s+5)\sin(2\pi s)$$
$$D(b,s)=\begin{cases}\biggl(\sin(\ds\frac{\pi s}{2})\biggr)^3\text{ lorsque }b=1,\\
\\
\sqrt{-1}\biggl(\cos(\ds\frac{\pi s}{2})\biggr)^3\text{ lorsque }b=-1,\end{cases}$$
ainsi on retrouve bien (\à une constante et une puissance de $2\pi$ pr\ès) un des r\ésultats de M.Muro (\cite{muro}), on rappelle que 
les calculs explicites des coefficients associ\és \à l'\eq  du \PV $(G,\g_{\pm 1})$ ont \ét\é faits par M.Muro (calcul micro-local).\\

Notons que les ouverts $O_{\pm 1,\pm 1}$ ne correspondent pas du tout \à la description de $\g''_1$ en $P-$orbites et pas plus \à la description de $\g'_1$ en $G$-orbites, par exemple $\sum_{1\≤i\≤7}X_{\lambda_i}$ n'appartient pas \à $G(\sum_{1\≤i\≤6}X_{\lambda_i}-X_{\lambda_7})=G(\sum_{1\≤i\≤3}X_{\lambda_i}- X_{\lambda_4}-X_{\lambda_5}+X_{\lambda_6}+X_{\lambda_7})$ (\cite{mullerJA2}).\\

2) Dans le cas  $\goth p$-adique, $ G$ a une seule orbite dans $\g'_1.$  

\noi Lorsque la caract\éristique r\ésiduelle est diff\érente de $2$, ce cas a \ét\é \étudi\é par I.J.Igusa, il a notamment calcul\é explicitement la fonction Z\éta associ\ée \à l'indicatrice d'un r\éseau $L$ ``convenable" de $\g_1$ (\cite{igusa7}) alors on v\érifie  que $\ds\frac{Z(1_L;s)}{Z(1_L;-s-5)}=C(s)B(Id,s).$

\bigskip
\subsection {$\bf (E_8,\alpha_1)$ }
\bigskip
\subsubsection {\bf Le cas d\éploy\é}
\bigskip
Ce cas est tout \à fait semblable au cas pr\éc\édent : $(E_7,\alpha_2)$.\\

On note $\{\lambda_1,...,\lambda_8\}$ les $8$ racines orthogonales obtenues canoniquement par orthogonalisations successives \à partir de $\lambda_1=\alpha_1$ (cf.\cite{mullerJA2} pour une liste explicite) alors:
$$\sum_{1\≤i\≤8}h_{\lambda_i}=2H_0\ ,$$
les racines de $\∆_2$ sont de la forme $\frac{1}{2}(\lambda_i+\lambda_j+\lambda_k+\lambda_l)\ ,$ avec $1\≤i<j\≤4$ et   $k=i+4,l=j+4,$ not\ée $\omega_{i,j},$ ou bien  $k=i'+4,l=j'+4,$ et  not\ée $\omega'_{i,j},$ avec $\{i,i',j,j'\}=\{1,2,3,4\}$ ainsi que  $\omega=\frac{1}{2}\sum_{1\leq i\leq 4}\lambda_i$ et 
$\widetilde\alpha=\frac{1}{2}\sum_{5\leq i\leq 8}\lambda_i$ qui est la plus grande racine de $\Delta.$\\
 
\noi $P=P(H_1,H_2),$ avec $H_2=2h_{ \widetilde\alpha}=\sum_{5\≤i\≤8}h_{\lambda_i}$ et $H_1=2h_{\omega}=\sum_{1\≤i\≤4}h_{\lambda_i}$  est   l'unique \sg parabolique standard tr\ès sp\écial du \PV , ce \sg parabolique est associ\é \à la partie $\∑_0-\{\alpha_8\}$ de $\∆_0$ qui est de type $D_7.$\\

\noi Dans cette situation on a simplement  $E_1(H_1)\cap \goth g_1=\{0\}.$
 \\
 
On fixe dor\énavant un syst\ème de Chevalley:   $(X_{\alpha},h_{\alpha},X_{-\alpha})_{\alpha\in \∆}$, de $ \goth g$ tel que pour toute racine $\omega\in \∆_2$ on ait:
$$ [X_{-\lambda_i},[X_{-\lambda_j},[X_{-\lambda_k},[X_{-\lambda_l},X_\omega]]]]= X_{-\omega}\ , \
\hbox{avec}\  \omega=\frac{1}{2}(\lambda_i+\lambda_j+\lambda_k+\lambda_l)\ ,$$
dont l'existence est assur\ée par le corollaire 4.3 de \cite{mullerJA2}.\\

On  rappelle que  l'\irf , $F$, est de degr\é $8,$  on le normalise par la  condition : $$F(\sum_{1\leq i\leq 8}X_{\lambda_i})=1 \ \hbox{et on a pour}\ a_1,...,a_8\in \F^*\ :\quad F(\sum_{1\leq i\leq 8}a_iX_{\lambda_i})=\prod_{1\leq i\leq 8}a_i\ ,$$
et que 
 $\tilde B=\displaystyle{-\frac{4 B}{B(H_0,H_0)}}(=\displaystyle{-\frac{2B}{B(h_{\lambda_i},h_{\lambda_i})}}),$ ainsi $\tilde B(\frac{H_1}{2},\frac{H_1}{2})=-2$ et pour toute racine  $\alpha$ de $\∆$ on a $\tilde B(X_{\alpha}, X_{-\alpha})=1.$\\ \\

\begin{prop}
\begin{enumerate}
\item Le \PV $(E_0(H_1)\cap \goth g_0,E_2(H_1)\cap \goth g_1)$
est de type $(E_7,\alpha_1)$ et d'\irf $F_1(x)=\frac{1}{4!}\tilde B(adx^4(X_{-\omega},X_{-\omega}).$ 

\item Pour $x\in \goth g_1$ on a $P_1(x)=\frac{1}{4!}\tilde B(adx^4(X_{-\widetilde\alpha},X_{-\widetilde\alpha}).$ 

\item Soit $x\in W_{\goth t},$ $x=x_2+x_0$ avec $x_i\in E_i(  H_1 )\cap \goth g_1$ alors:\\
 
 a) Le \PV $((E_0(H_1)\cap \goth g_0)_{x_0},(E_2(H_1)\cap \goth g_1)_{x_0})$ est de type:

$$(\∆,\lambda_0)=\bigg\{\begin{array}{l} (E_6,\alpha_2)\ \hbox{lorsque}\ -P_1(x_0)\in \F^{*2}\\
(F_4,\alpha_1)\ \hbox{sinon}\ ,\end{array}$$
et $c(f)\equiv  P_1(x_0)\quad (\text{mod. }\ \F^{*2}).$

 b) Le \PV $((E_0(H_2)\cap \goth g_0)_{x_2},(E_2(H_2)\cap \goth g_1)_{x_2})$ est de type:

$$(\∆,\lambda_0)=\bigg\{\begin{array}{l} (E_6,\alpha_2)\ \hbox{lorsque}\ -F_1(x_2)\in \F^{*2}\\
(F_4,\alpha_1)\ \hbox{sinon}\ ,\end{array}$$
et $c(f)\equiv  F_1(x_2)\quad (\text{mod.  }\ \F^{*2}).$
 
\end{enumerate}
\end{prop}

\dem

1) Cela r\ésulte du diagramme de Dynkin compl\ét\é de $E_8.$  \\

2)  $H_1$ et $H_2$ sont dans la m\ême orbite de $G,$ il suffit d'appliquer le 1) de la proposition 8.2.4.\\

3) Pour tout $x\in W_{\goth t},$ il existe $g\in G_{H_1}$ et $(a_1,...,a_8)\in \F^*$ tels que $gx=\sum_{1\≤i\≤8}a_iX_{\lambda_i}$ (lemme 7.3 de \cite{mullerJA2}) donc on fait la \demo pour $x_2=\sum_{1\≤i\≤4}a_iX_{\lambda_i}$ et $x_0=\sum_{5\≤i\≤8}a_iX_{\lambda_i}.$\\

a) Le \PV $((E_0(H_1)\cap \goth g_0)_{x_0},(E_2(H_1)\cap \goth g_1)_{x_0})$ est un \PV \air quasi-commutatif presque d\éploy\é pour lequel $(E_4(H_1)\cap \goth g_2)_{x_0}$ est de dimension $1$  or on a:
$$\begin{array}{rl}(E_2(H_1)\cap \goth g_1)_{x_0}&=\oplus_{1\≤i\≤4}\goth g^{\lambda_i}\oplus_{1\≤i<j\≤4}(\F X_{i,j}\oplus \F  X'_{i,j})\oplus_{1\≤i\≤4}[X_{-\lambda_i},X_\omega]\ \hbox{avec}\\ 
\\
X_{i,j}&=a_{4+j}[X_{-\lambda_{4+i}},X_{\omega_{i,j}}] -a_{4+i}[X_{-\lambda_{4+j}},X_{\omega_{i,j}}]\ ,\\
\\
X'_{i,j}&=a_{4+j'}[X_{-\lambda_{4+i'}},X_{\omega'_{i,j}}] -a_{4+i'}[X_{-\lambda_{4+j'}},X_{\omega'_{i,j}}]\ ,\  
\{i,i',j,j'\}=\{1,...,4\},\end{array}$$
donc, par classification, ce \PV est de type $(R,\lambda_0)=(F_4,\alpha_1)$ avec $d=2$ donc de type $(\∆,\lambda_0) =
 (E_6,\alpha_2)$ lorsque $f$ est isotrope et $(F_4,\alpha_1)$  sinon (cf.3),lemme 8.2.2) ainsi il suffit de consid\érer la forme quadratique:
  $$\begin{array}{rl}f'(u,v)&=\frac{1}{2}\tilde B(ad(uX_{1,2}+vX'_{1,2})^2(X_{-\lambda_1}),X_{-\lambda_2})\\
  &=-\frac{1}{2}\tilde B( [uX_{1,2}+vX'_{1,2},X_{-\lambda_1}],[uX_{1,2}+vX'_{1,2},X_{-\lambda_2}])\\
  &=c(u^2a_5a_6+v^2a_7a_8)\ ,\end{array}$$
$c$ \étant une constante non nulle, $f'$ est isotrope \ssi $P_1(x_0)\equiv -1$ (mod $\F^{*2}$).\\

b) Idem que pour a)\fdem \\ \\

Lorsque $\F$ est archim\'edien, les $2$ polynomes de Bernstein s'obtiennent imm\'ediatement en appliquant les  propositions  3.4.4 et  3.7.3 (ainsi que 8.2.5)   d'o\`u:
  
  \begin{prop} Soit  $B(s)= s(s+\frac{3}{2})(s+\frac{5}{2})(s+4)$ le polynome de Bernstein du \PV de type $(E_6,\alpha_2)$ alors: $$b_1(s_1,s_2)= B(s_2)\ \hbox{et}\ \ b_2(s_1,s_2)=B(s_2)B(s_1+s_2+3).$$
  \end{prop}
  
  \bigskip

\noi  On rappelle que pour $u=(u_1,u_2)\in (\FF)^2$ on d\'efinit les ouverts (tous non vides):\\

$O_u=O_{(u_1,u_2)}=\{x\in \goth
g_1\ |\ F_1(x) \F^{*2}= u_1\ \ ,\ F_2(x) \F^{*2}= u_1u_2 )
  \}$\\
 
$O^*_u=O^*_{(u_1, u_2)}=\{x\in
\goth g_{-1}\ |\ F^*_1(x) \F^{*2}= u_2 \   ,\ F^*_2(x) \F^{*2}=
 u_1u_2  \}$\\
 
\noi ainsi que les fonctions Z\'etas correspondantes:\\
  
$Z_u(f;\omega)=Z(f{\bf 1}_{O_u};\omega)$ pour $f\in \ES(\goth
g_1)$  
et
$ Z^*_u(h;\omega)=Z^*(h{\bf 1}_{O^*_u};\omega)$  pour $h\in \ES(\goth g_{-1})$.  \\ \\

\noi  Pour $u,v,w\in \FF,$  $\omega\in \hat \F^*$ et $s\in \C$ posons:
$$B_{v,u,w}(\omega,s)=\sum_{t\in \FF}A^1_{\omega,s+1,s+\frac{5}{2}}(wv,t,1)(t,-w)A^1_{\omega,s+\frac{7}{2},s+ 5}(t,wu,1)\ ,$$
(cf. 8.2.8 et sa d\émonstration ainsi que la prop.8.4.1).
Alors:

\begin{prop} 
\begin{enumerate}

\item Lorsque $\F=\C$ pour $f\in \ES(\goth
g_1),$  on a:
$$ Z^*( \four (
f);(\omega_1,\omega_2),(s_1,s_2))=a_2(\omega_2,s_2)a_2(\omega_1\omega_2,s_1+s_2+3)Z\biggl(f;( \omega_1,( \omega_1\omega_2)^{-1});s_1,-(s_1+s_2+ 8 )\biggr)$$
avec:
$$a_2(\omega,s)=\rho'(\omega;s+1)\rho'(\omega;s+\frac{ 5}{2})\rho'(\omega ;s +\frac{7}{2})\rho'(\omega ;s +5).$$
\item  Lorsque $\F=\R$ ou bien $\F$ est un corps $\goth p$-adique de caract\éristique r\ésiduelle diff\érente de 2, pour tout  $v=(v_1, v_2)$ dans $(\FF)^2$ et $f\in \ES(\goth
g_1),$  on a:
$$ Z^*_v( \four(
f);\  (\omega_1,\omega_2),(s_1,s_2))=\sum_{u\in (\FF)^2}d_{v,u}(
 \omega,s)  Z_u(f;( \omega_1,( \omega_1\omega_2)^{-1});s_1,-(s_1+s_2 +8 )\ ) \quad \hbox{avec}$$

   $u=(u_1, u_2)$ :  $d_{v,u}(
\omega,s)= B_{v_1,u_1,v_2}(\omega_2,s_2)B_{v_2,u_2,u_1}(\omega_1.\omega_2,s_1+s_2+3).$\end{enumerate}
\end{prop}
  
  \dem

\noi On applique la proposition 5.3.1 dont toutes les hypoth\`eses sont v\'erifi\'ees par les propositions 8.2.11 et  8.4.1, puisque $c(f)\equiv P_1(x_0)$ (resp. $c(f)\equiv F_1(x_2)$) dans le \PV $((E_0(H_1)\cap \goth g_0)_{x_0},(E_2(H_1)\cap \goth g_1)_{x_0})$ (resp.$((E_0(H_1)\cap \goth g_0)_{x_0},(E_2(H_1)\cap \goth g_1)_{x_0})$) et que $r_1=3.$ 
On notera que la restriction de $F_1$ (resp. $P_1$) \à $(E_2(H_1)\cap \goth g_1)_{x_0})$ (resp. $(E_2(H_2)\cap \goth g_1)_{x_2})$ est 
 l'\irf utilis\é dans le paragraphe 8.2 (cf. 1) de la proposition 8.2.4) \à un coefficient multiplicatif pr\ès appartenant \à $\F^{*2}$ ce qui est sans  incidence (cf.lemme 3.5.3).\fdem\\

\begin{rema}
\end{rema} 

1) Dans le cas r\éel,  posons:

$u_1=-1$ et $u_2=1$ alors   pour $u=\pm 1,$
$B_{u_i,u_j,u}(s)= B_{i,j}(u,s),$  $B_{i,j}(u,s),$  \étant le $i,j-$\ème coefficient de la matrice:
$$\frac{1}{2}C_2(0,s)\displaystyle\left(\begin{array}{cc} \sin2\pi s &0 \\
 -2(1+u) \sin \pi s & 
  u  \sin 2\pi s \end{array}\right)\ (\text{prop.}8.2.11)$$

Soient $a_1=(-1,-1),a_2=(-1,1),a_3=(1,-1),a_4=(1,1)$ et soit $D(s_1,s_2)=\biggl(D_{i,j}(s_1,s_2)\biggr)_
{1\≤i\≤j\≤4}$ la matrice:

$$\displaystyle\left(\begin{array}{cccc} \sin2\pi s_2\sin 2\pi(s_1+s_2)&0&0&0\\
0&- \sin 2\pi s_2\sin 2\pi(s_1+s_2)&0& 0\\
0&0&- \sin 2\pi s_2\sin 2\pi(s_1+s_2)& 0\\
0&4\sin \pi s_2\sin 2\pi (s_1+s_2)& 4\sin 2\pi s_2\sin \pi (s_1+s_2)& \sin 2\pi s_2\sin 2\pi(s_1+s_2)\end{array}\right)$$\\

\noi alors :
$$d_{a_i,a_j}(s_1,s_2)=\frac{1}{4}C_2(0,s_2)C_2(0,s_1+s_2+3)D_{i,j}(s_1,s_2).$$
 
\noi On a : \\

\noi $d_{a_1,a_1}(s_1,s_2)+d_{a_4,a_1}(s_1,s_2)=d_{a_1,a_4}(s_1,s_2)+d_{a_4,a_4}(s_1,s_2)$\\

\noi  $d_{a_2,a_1}(s_1,s_2)+d_{a_3,a_1}(s_1,s_2)=d_{a_2,a_4}(s_1,s_2)+d_{a_3,a_4}(s_1,s_2)=0$\\

\noi  $d_{a_2,a_2}(s_1,s_2)+d_{a_3,a_2}(s_1,s_2)=d_{a_2,a_3}(s_1,s_2)+d_{a_3,a_3}(s_1,s_2)$\\

\noi mais   $d_{a_1,a_2}(s_1,s_2)+d_{a_4,a_2}(s_1,s_2)\not =d_{a_1,a_3}(s_1,s_2)+d_{a_4,a_3}(s_1,s_2)$ pour $s_1\not=0$ donc on n'a pas l'analogue de la proposition 8.2.8 pour $(s_1,s_2)\in \C^2.$ \\

\noi Lorsque $\F$ est un corps $\goth p$-adique, les expressions analogues semblent compliqu\ées.\\

\noi Cependant lorsque $s_1=0,$ on a encore les \équations fonctionnelles habituelles.  Pour $u\in \FF,$ et $ s\in \C$ tel que $\Re (s)>0,$ on pose:\\

$\bullet$ $O_u=\{x\in \g_1$ tel que $F_2(x)  \F^{*2}=u\}$  et   pour $f\in \EuScript S(\g_1):$  $Z_u(f;s):=Z(1_{O_u}f;s),$\\

$\bullet$ $O^*_u=\{x\in \g_{-1}$ tel que $F^*_2(x)  \F^{*2}=u  \}$ et pour $g\in \EuScript S(\g_{-1}):$   $Z^*_u(g;s):=Z^*(1_{O^*_u}g;s),$\\

\noi alors  $Z_u(f;s)$ et $Z^*_u(g;s)$ admettent un prolongement m\éromorphe sur $ \C$   et  si $\FF=\{u_1,...,u_\ell\},$ on a:

\begin{prop} Dans le cas r\éel ou bien dans le cas $\goth p$-adique  avec caract\éristique r\ésiduelle diff\érente de $2, $ pour tout $f\in \EuScript S(\g_1)$ et $s\in \C$ on a:

$$\biggl(\begin{array}{ll}
&Z^*_{u_1}(\four( f);s)\\
&....\\
&Z^*_{u_\ell}(\four (f);s)\end{array}\biggr)=  f(s)A'(s).A'(s+\frac{11}{2})
\biggl(\begin{array}{ll}
&Z_{u_1}( f;-s-8)\\
&....\\
&Z_{u_\ell}(  f;-s-8)\end{array}\biggr).$$
avec $$f(s)={ |2|_\F}^{-4s-16}\rho (|\ |^{2s+7})\rho (|\ |^{2s+10}),$$
et $$A'(s)=\biggl(A^1_{s+1,s+\frac{5}{2}}(u_i,u_j,1)\biggr)_{1\≤i,j\≤\ell}\quad s\in \C$$
est la matrice des coefficients associ\és \à une forme quadratique de discriminant $1$ sur un espace de dimension $5.$

\end{prop}

\dem
 Les coefficients, $a_{v,u}(s),$ pour $u,v\in \FF,$ sont donn\és par (prop.8.4.3,2)):
$$ \begin{array}{ll}a_{v,u}(s)&=\sum_{v_1\in\FF}d_{(v_1,vv_1),(u_1,u_1u)}(0,s)\quad \text{avec }\ u_1\in \FF,\\
\\
&= \sum_{v_1\in\FF}B_{v_1,u_1,vv_1}( s)B_{vv_1,uu_1,u_1}(s+3)\\ 
\\
 &=\sum_{t_1,t_2\in \FF}A(t_1,t_2) \end{array}$$
avec $A(t_1,t_2)= \sum_{v_1\in\FF}$
$$ A^1_{ s+1,s+\frac{5}{2}}(v, t_1,1)(t_1,-vv_1)A^1_{ s+\frac{7}{2},s+ 5}(t_1,v v_1u_1,1)A^1_{ s+4,s+\frac{11}{2}}(u_1vv_1, t_2,1)(t_2,-u_1)A^1_{ s+\frac{13}{2},s+ 8}(t_2,u, 1)$$
$$= (t_1,-v) (t_2,-u_1)A^1_{ s+1,s+\frac{5}{2}}(v, t_1,1)A^1_{ s+\frac{13}{2},s+ 8}(t_2,u, 1).A'(t_1,t_2)$$
$$ \begin{array}{rcl}\text{avec }\ \ A'(t_1,t_2)&=&\sum_{w_1,w_2\in\FF}A(t_1,t_2,w_1,w_2)\alpha(w_1)\alpha(w_2) \rho(s+\frac{7}{2};t_1w_1)\rho(s+\frac{11}{2};t_2w_2) \quad \text{et}\\ 
\\
 A(t_1,t_2,w_1,w_2)&=& \sum_{v_1\in\FF}(v_1,t_1)\rho(s+ 5;w_1vv_1u_1)\rho(s+4;w_2u_1v_1v) \\
 \\
 &=& \sum_{v_1\in\FF}(v_1w_1w_2u_1v,t_1)\rho(s+ 5;w_2v_1 )\rho(s+4;w_1 v_1)\end{array}$$
donc:
$$ A'(t_1,t_2 )=( u_1v,t_1)\sum_{v_1\in \FF}\ (v_1 ,t_1)\ ( \sum_{w_1 \in\FF}  \alpha(w_1)  ( w_1 ,t_1)\rho(s+\frac{7}{2};t_1w_1) \rho(s+4;w_1 v_1) .$$
$$( \sum_{ w_2\in\FF}\alpha(w_2)  ( w_2 ,t_1)\rho(s+ 5;w_2v_1 )\rho(s+\frac{11}{2};t_2w_2)$$
$$=( u_1v,t_1)\sum_{v_1\in \FF}\ (v_1 ,t_1)\ A^1_{s+\frac{7}{2},s+4}(t_1,v_1,t_1)A^1_{s+5,s+\frac{11}{2} }( v_1,t_2,t_1)$$
$$=( u_1v,t_1)\sum_{v_1\in \FF}\ (v_1 ,t_1)\alpha(1)^2\alpha(-t_1)^2 A^1_{s+\frac{7}{2},s+4}(1,t_1v_1,1)A^1_{s+5,s+\frac{11}{2} }( t_1v_1,t_1t_2, 1)$$
$$A'(t_1,t_2 )=\begin{cases}( u_1v,t_1)  A^1_{s+\frac{7}{2},s+4}(1, 1,1)A^1_{s+5,s+\frac{11}{2} }(  1, 1, 1)\ \text{pour }t_1=t_2\\
0\quad   \text{pour }t_1\not=t_2\end{cases}$$
(lemme 3.6.8,1)),  d'o\ù: 
$$a_{v,u}(s)= \sum_{t\in \FF}A(t,t)=f(s)\sum_{t\in \FF}A^1_{ s+1,s+\frac{5}{2}}(v, t,1)A^1_{ s+\frac{13}{2},s+ 8}(t,u, 1) $$
avec:
$$f(s)=  A^1_{s+\frac{7}{2},s+4}(1, 1,1)A^1_{s+5,s+\frac{11}{2} }(  1, 1, 1)={ |2|_\F}^{-4s-16}\rho (|\ |^{2s+7})\rho (|\ |^{2s+10}).\qquad\ \Box$$

\begin{rema} 1) Les orbites de $G$ dans $\g_{\pm 1}$ sont donn\ées dans le th\éor\ème 2 de \cite{mullerJA2}.\\

2) Dans le cas $\goth p$-adique.\\

\noi  Les orbites de $G$ dans $\g'_1$ (resp.$\g'_{-1}$) sont donn\ées par $O'_u $ (resp.
$O'*_u $) pour $u\in \FF.$\\

\noi Lorsque la caract\éristique r\ésiduelle est diff\érente de $2, $ $\ell=4$ et la proposition 8.3.5 donne les  coefficients de l'\eq \ associ\ée au \PV $(G,\g_{\pm 1})$.\\
\noi On a alors: $$f(s)= \frac{(1-q^{2s+6})(1-q^{2s+9})}{(1-q^{-2s-7})(1-q^{-2s-10})}.$$

3) Dans le cas r\éel, la proposition pr\éc\édente  donne:
 $$
\left(\begin{array}{ll}Z^*_{ -1}(\four (f); s)\\
 
Z^*_{ 1}\ (\four (f); s)\end{array}\right)=  \frac{1}{4}C_2(0,s)C_2(0,s+3) E( s)\left(\begin{array}{ll}Z_{  -1}(f; -s-8)\\
 Z_{ 1}\ ( f; -s-8)\end{array}\right)$$
avec
$$E(s)=\displaystyle\left(\begin{array}{ll} -\sin^22\pi s &0 \\
4\sin \pi s \sin 2\pi s & \sin^22\pi s \end{array}\right)$$
ce qui n'est pas en accord avec les r\ésultats obtenus par Muro (\cite{muro}).\\

\noi (Pour m\émoire: $ \frac{1}{4}C_2(0,s)C_2(0,s+3)=$

\noi$4(2\pi)^{-8(s+\frac{9}{2})}\Gamma (s+1)  \Gamma (s+\ds\frac{5}{2})\Gamma (s+\ds\frac{7}{2})\Gamma (s+5)$ $\Gamma (s+4)  \Gamma (s+\ds\frac{11}{2})\Gamma (s+\ds\frac{13}{2})\Gamma (s+8)$)\\

\noi $G$ a $3$ orbites dans $\g'_1$  (resp.$\g'_{-1}$) :
$$O_1=G.(\sum_{1\≤i\≤8}X_{\lambda_i})\cup  G.(\sum_{1\≤i\≤6}X_{\lambda_i}-X_{\lambda_7}-X_{\lambda_8})\ ,\ O_{-1}=G.(\sum_{1\≤i\≤7}X_{\lambda_i}-X_{\lambda_7}).$$
Si l'on veut d\éterminer les coefficients de l'\eq associ\ée   aux orbites, il faut proc\éder de mani\ère analogue au paragraphe 8.2.4 ce que ne sera pas fait ici.\end{rema}

\subsubsection{\bf Le cas r\éel non d\éploy\é}

\bigskip

Dans ce paragraphe, on consid\ère la forme EIX de $E_8$ plus pr\écis\ément $\overline g$ est de type $E_8$ et $\g$ est de de type $F_4$, les \PVs correspondants sont les suivants:\\

$\bullet$  $({\overline\g}_0,{\overline\g}_1)$ est de type $(E_8,\alpha_1),$\\

$\bullet$ $(\g_0,\g_1)$  est de type $(F_4, \lambda_4),$ \\

et les diagrammes de Satake et de Dynkin sont:\\ \\

 \hskip 10pt \hbox to 4cm {\lower 2pt\hbox{$\circledcirc$}\hrulefill\lower 2pt
\hbox{$\bullet$}\hrulefill \lower 2pt \vtop {\offinterlineskip \hbox{$\bullet$} \hbox to 5pt{\hfill \vrule height 12pt width 0,3pt\hfill} \hbox{$\bullet$}}\hrulefill\lower 2pt\hbox{$\bullet$}\hrulefill\lower 2pt\hbox{$\circ$}\hrulefill\lower 2pt\hbox{$\circ$}\hrulefill\lower 2pt\hbox{$\circ$} }
\qquad \qquad \hbox to 2,5cm{ \offinterlineskip \lower 2pt\hbox{$\circ$}\hrulefill\lower 2pt\hbox{\offinterlineskip {$\circ$}\hglue -1,8pt\vbox{
 {\hrule height 0,3pt width 0,6cm}\vskip 3,5pt{\hrule height 0,3pt width 0,6cm}
\vskip 0,3pt} \hglue -16pt  $>$ \hglue -1,8pt  {$\circ$}}\hrulefill\lower 2pt\hbox{$\circledcirc$} }\\ \\

($\lambda_1$ est la restriction de $\alpha_8$, $\lambda_2$ est la restriction de $\alpha_7$, $\lambda_3$ est la restriction de $\alpha_6$, $\lambda_4$ est la restriction de $\alpha_1$)\\ 

Comme $\omega$ et $\widetilde\alpha$ sont des racines r\éelles, le \sgp standard est le m\ême que dans le cas d\éploy\é,  il est donn\é par $P(H_1,H_2)$ avec $H_1=2h_{\omega}$ et $H_2=2h_{\widetilde\alpha}$ qui sont dans la m\ême orbite de $G.$\\

L'alg\èbre simple $  E=\goth U(\R H_2) $  (resp. $ F=\goth U(\R H_1)$) est gradu\ée par $\ds\frac{1}{2}H_1$ (resp.$\ds\frac{1}{2}H_2$) et  son diagramme de Satake est de type EVII, le diagramme de Dynkin du \PV $(E_0,E_1)$ (resp. $(F_0,F_1)$) est de type $(C_3,\lambda_1),$ il est obtenu en prenant simplement $\Delta\cap (\oplus_{i=2,3,4}\Z \lambda_i).$ \\

\noi Ceci donne une situation tr\ès simple car:
 
\begin{lem} Soit $x\in E'_1$  (resp.$x\in F'_1$) et $\goth U_x=\goth U(\F x\oplus \F H_2\oplus \F x^{-1})$  (resp.$\goth U_x=\goth U(\F x\oplus \F H_1\oplus \F x^{-1})$) alors $\goth U_x$ est de rang 2 et de type EIII et le \PV  $({(\goth U_x)}_0,({\goth U_x)}_1)$ est de type $(BC_2,\lambda_1).$
  \end{lem}
 
 \dem  Par la proposition 8.4.1, le complexifi\é de  $((\goth U_x)_0, (\goth U_x)_1)$ est de type $(E_6,\alpha_2)$ et par la proposition de l'appendice 1, $\goth U_x$ est de rang inf\érieur ou \égal \à $2$ donc  $\goth U_x$ est de type EIII par classification (tables de \cite{warner}) d'o\ù le \PV  $({(\goth U_x)}_0,({\goth U_x)}_1)$ est de type $(BC_2,\lambda_1).$
 
 \fdem \\

\noi On en d\éduit:
  
  \begin{prop}
  
   \begin{enumerate} 
\item  $P$ a une seule orbite dans $\g"_1.$

 \item Soient $(s_1,s_2)\in \widehat{(\C)^2},$ 
 $f\in \EuScript S(\g_1)$ alors:
 $$ Z^*(\four f;(s_1,s_2))=C(s_1,s_2)\sin(2\pi s_2)\sin (2\pi(s_1+s_2)) Z(f;(s_1,-s_1 -s_2 - 8)\quad \text{avec} $$
 $$\begin{array}{ll}C(s_1,s_2)&=4(2\pi)^{-(4s_1+8s_2 +36)}\Gamma(s_2+1) \Gamma(s_2+\ds\frac{5}{2}). \Gamma(s_2+\ds\frac{7}{2})\Gamma(s_2+5)\Gamma(s_1+s_2+4)\\
 &.\Gamma(s_1+s_2+\ds\frac{11}{2})\Gamma(s_1+s_2+ \ds\frac{13}{2}) \Gamma(s_1+s_2+8) .\end{array} $$
 \end{enumerate}
 \end{prop}

\dem  On proc\ède comme dans la d\émonstration de la proposition 5.3.3 que l'on applique avec $s_1=0$ (cas $d=2$) ainsi que la proposition 5.3.1 avec $\mathbb{H}=\R^*.$\fdem\\

\noi Remarque: Cette proposition est bien connue lorsque $s_1=0$ (\cite{muro}).

   \newpage
 
\centerline{\bf \large Appendice 1}

\bigskip

On indique un r\ésultat sur le rang du centralisateur associ\é \à un \Sl  dont l'\él\ément simple est donn\é par $2h_{  \alpha},$ $h_{  \alpha}$ \étant la co-racine de la racine longue $\alpha$  dans un syst\ème de racines  irr\éductible et r\éduit.

\noi Ce r\ésultat permet de simplifier les calculs en "petit rang".\\ \\

On commence par un premier lemme:\\

\noi {\bf Lemme 1} {\it Soit $\Delta$ un syst\ème de racines irr\éductible et r\éduit , on note:\\

$\bullet$
$\Sigma$ une base de $\Delta,$\\

$\bullet$ $\widetilde \alpha$ la plus grande racine de $\Delta$ pour l'ordre induit par $\Sigma,$ \\

$\bullet$ $\Sigma'=\{\alpha\in \Sigma\ |\ n(\alpha,\widetilde \alpha)=0\}.$\\

Soit $\gamma$ une racine de $\Delta$  v\érifiant $n(\gamma,\widetilde \alpha)=1$  alors  il existe au moins une racine $\alpha$ de $\Sigma'$  (qui est non vide) telle que la quantit\é $n(\gamma,\alpha)$ soit non nulle.}\\

\dem Comme $\widetilde \alpha-\gamma$ est une racine, $\gamma$ est positive.\\

1) Dans le cas $A_n$ et dans les notations de la planche 1 de \cite{bourbakigal6}:\\

$\widetilde \alpha=\epsilon_1-\epsilon_{n+1},$   $\Sigma'=\{\alpha_i=\epsilon_i-\epsilon_{i+1},i=2,..,n-1\}$ et:\\

$\gamma=\epsilon_1-\epsilon_{i+1},$  avec $i=1,...,n-1$ ou bien 
$\gamma=\epsilon_i-\epsilon_{n+1},$ avec $i=2,...,n$ et dans les $2$ cas $n(\gamma,\alpha_i)\not=0$ pour $1=2,...,n-1;$ $n(\epsilon_1-\epsilon_2,\alpha_2)\not=0$ et $  n(\epsilon_n-\epsilon_{n+1},\alpha_{n-1})\not=0.$\\

2) Dans les cas restants, $\Sigma-\Sigma'=\{\alpha_0\}$ et $\gamma-\alpha_0$ est une combinaison lin\éaire (\à coefficients positifs ou nuls) d'\elts de $\Sigma'$.\\

\noi Proc\édons par l'absurde et supposons que $\gamma$ soit orthogonale aux racines de  $\Sigma'$  alors $(\gamma,\gamma)=(\gamma,\alpha_0)$ donc $n(\alpha_0,\gamma)=2$ d'o\ù $\gamma$ est une racine courte positive orthogonale \à $\Sigma'$ et $\Delta$ n'est pas simplement lac\é.\\

\noi  Il est facile de  v\érifier cas par cas (dans les notations des planches de \cite{bourbakigal6}) que c'est absurde, en effet:\\

a) $\Delta=G_2:$ alors $\gamma=\epsilon_i-\epsilon_3$ avec $i=1$ ou $2$ et $\Sigma'=\{\alpha_1\}=\{\epsilon_1-\epsilon_2\}.$\\

b)  $\Delta=F_4:$ alors $\gamma=\epsilon_i,i=1,2,3,4,$ ou bien $\gamma=\ds\frac{1}{2}(\pm\epsilon_1\pm\epsilon_2\pm\epsilon_3\pm\epsilon_4)$  et 

\noi$\Sigma'=\{\epsilon_3-\epsilon_4,\epsilon_4,$ $\ds\frac{1}{2}(\epsilon_1-\epsilon_2-\epsilon_3-\epsilon_4)\}.$\\

c) $\Delta=B_n:$ alors $\widetilde \alpha=\epsilon_1+\epsilon_2,$   $\gamma=\epsilon_i,i=1,2$ et $\Sigma'=\{\alpha_1=\epsilon_1-\epsilon_2,\alpha_i, 3\≤i\≤n\}.$\\

d) $\Delta=C_n:$ alors $\widetilde \alpha=2\epsilon_1,$   $\gamma=\epsilon_1\pm\epsilon_i,i=2,...,n$ et 

\noi$\Sigma'=\{\alpha_i=\epsilon_i-\epsilon_{i+1},  2\≤i\≤n-1,\alpha_n=2\epsilon_n\}.$\fdem\\ 

\noi {\bf Lemme 2} {\it Soit $\Delta$ un syst\ème de racines irr\éductible et r\éduit, $\alpha$ et $\gamma$  $2$ racines  de $\Delta$ telles que:

i) $n(\gamma,\alpha)=1,$

ii) $\alpha$ est une racine longue

alors il existe au moins une racine $\delta$ telle que $n(\gamma, \delta)\not=0$ et $n( \delta,\alpha)=0.$}\\

\dem Soit $\Sigma$ une base de $\Delta$ et $\widetilde \alpha$ la plus grande racine de $\Delta$ pour l'ordre induit par $\Sigma,$ alors $\alpha$ et $\widetilde \alpha$ ont m\ême longueur donc sont dans la m\ême orbite du groupe de Weyl  $W$ et il suffit d'appliquer le lemme pr\éc\édent.\fdem\\

\noi {\bf Proposition} {\it  Soit $\g$ une alg\èbre de Lie simple, $\goth a$ une sous-alg\èbre ab\élienne d\éploy\ée maximale de $\g$, $ \alpha$  une racine longue  de $\Delta,$ le syst\ème de racines de $(\g,\goth a)$ que l'on suppose r\éduit.\\

  $(x,2h_{  \alpha},y)$ \étant un \Sl , on appelle $\goth U_x$ l'alg\èbre d\ériv\ée du centralisateur dans $\g$ de l'alg\èbre  $\F x\oplus\F h_{  \alpha}\oplus\F y,$ alors:
  $$\text{rang de} (\goth U_x)\≤ \text{rang de }(\g)-2\ .$$}

\dem  Soit $E$ le centralisateur de $h_\alpha$ dans $\g$ alors:
$$E=c\oplus E'\ ,$$
$c$ \étant le centre de $E$ et contenant $h_\alpha$ et $E'$ \étant l'alg\èbre d\ériv\ée de $E.$

\noi Soit $\Delta'=\{\beta\in \Delta\ |\ n(\alpha, \beta)=0\}$ alors $\goth a'=\sum_{\beta\in \Delta'}\F h_\beta$ est une sous-alg\èbre ab\élienne d\éploy\ée maximale de $E',$  comme $ \goth U_x\subset E'$  on a  rang$(\goth U_x)\≤$rang$(E')=$dim$(\goth a')\≤$rang$(\g)-1$.

\noi Si  rang$(\goth U_x)\≤$rang$(\g)-2,$ la d\émonstration est termin\ée.\\

\noi Dans le cas contraire on a rang$(\goth U_x)$ $=$rang$(\g)-1=$dim$(\goth a').$

\noi Alors, soit $\goth a_x$ une sous-alg\èbre ab\élienne d\éploy\ée maximale de  $\goth U_x,$  $\goth a_x\oplus \F h_\alpha$ et $\goth a$ sont $2$ sous-alg\èbres ab\éliennes d\éploy\ées maximales de $\g$ 
 contenant $h_\alpha$ donc il existe un \elt $g$ dans $Aut(\g)_{h_\alpha}$ tel que $g(\goth a_x\oplus \F h_\alpha)=\goth a$ d'o\ù $(g(x),2h_{  \alpha},g(y))$ 
 est encore  un \Sl  (d'o\ù $g(x)\in \oplus_{\{\gamma|n(\gamma,\alpha)=1\}}\g^\gamma$) et l'alg\èbre $g(\goth U_x)=\goth U_{g(x)},$ admet
  $g(\goth a_x)$ comme sous-alg\èbre ab\élienne d\éploy\ée maximale
mais $g(\goth a_x)\subset \goth a\cap E'=\goth a'$ donc  $g(\goth a_x) =\goth a'$ d'o\ù $[g(x),\goth a']=0$ ce qui est absurde par le lemme pr\éc\édent \fdem

\newpage
\centerline{\bf \large Appendice 2}

\bigskip

 \centerline{\bf Remarques sur les mesures relativement invariantes sur les orbites singuli\ères}
\bigskip

\noi{\bf1. Conditions
suffisantes}

\bigskip

Dans ce paragraphe, on donne  des conditions
suffisantes d'existence des mesures  relativement invariantes sur le lieu singulier et temp\ér\ées.   

\noi Les hypoth\èses de ce  paragraphe sont les suivantes: $\F$ est un corps local de caract\éritique $0$ et on suppose que $(\goth g_0,\goth g_1)$ est un  \PV tel que:\\

\qquad (1)  $\goth g_1$ est un $\goth g_0$-module absolument simple,\\

\qquad (2) $2H_0$ est $1$-simple.\\

On consid\ère un $sl_2$-triplet $1$-adapt\é : $(x_0,h,y_0)$  non g\én\érique, c'est \à dire que $h\not=2H_0,$ et on suppose que les propri\ét\és suppl\émentaires suivantes
sont v\érifi\ées :\\

\qquad (3) le \PV  $(\goth U(\F h')_0, \goth  U(\F h')_1,\ds\frac{h}{2})$ admet un unique \irf  avec $h'=2H_0-h,$\\

\qquad (4) $h'$ est $1$-simple et le \PV $(\goth  U(\F h)_0,\goth  U(\F h)_1)$ admet un unique \irf .
 \\

Par les hypoth\èses (1),(2),(3),(4)  les  \PVs 
 $(\goth g_0,\goth g_1),$ $(\goth  U(\F h')_0,\goth  U(\F h')_1)$
et $(\goth  U(\F h)_0,\goth  U(\F h)_1)$ admettent des invariants
relatifs fondamentaux  non constants que l'on note respectivement :
$F,F_h,F_{h'}$ et $\chi,\chi_h,\chi_{h'}$ seront les caract\ères associ\és; on
note $m,m_h,m_{h'}$ les degr\és respectifs des polyn\ômes. \vskip 3mm
\noi Pour $(i,j)\in \Z,$ on d\ésigne par $E_{i,j}$ le sous-espace vectoriel:
$$E_{i,j}=E_i(h)\cap E_j(h')=\{x\in \goth g\ \ |\ \ [h,x]=ix\
\hbox{et}\  [h',x]=jx\ \}$$
et $p_{i,j}$ est sa dimension (notons que : $p_{\pm i,\pm j}=p_{i,j}$).

Soit $p_0=sup\{i\≥0\ \ /\ \ \goth g_i\not=\{0\}\}.$  
\vskip 3mm
$H$ est la composante connexe du noyau de $\chi,$ $H_h$ (resp.$H_{x_0}$)
est le centralisateur de $h$ (resp.$x_0$) dans $H,$ $N_h$ est le sous
groupe $exp(ad(\oplus_{i\≥1}E_i(h)\cap \goth g_0)),$ $P_h$ est le
sous-groupe parabolique $H_h.N_h$ et il est bien connu que
  $$H_{x_0}=(H_h)_{x_0}.N_{x_0}\subset P_h\subset
H.$$ 
 $\∆_{P_h}$ d\ésigne le module du sous-goupe $P_h$ lorsqu'il
est muni d'une mesure de Haar invariante \à gauche et not\ée $dp.$ 
\vskip 3mm
\noi {\bf Lemme 1} {\it $\∆_{P_h}^{-1}= |\chi_h|^r$ avec 
$r=\displaystyle\frac {\sum_{i=0}^{p_0}ip_{i,i}}{
m_h}. \frac{p\text{dim}(\goth g_p)}{\text{ trace}(ad(h')/\goth g_p)}$ avec $1\≤p\≤p_0.$

2) $H_{x_0}$ est unimodulaire.}\\

 \dem 1) Montrons que pour tout couple $(i,j)\in \Z\times \Z,
(i,j)\not=(0,0)$ et $g\in G_h$ on a: 
$$(a)\quad det(g/_{E_{i,j}})^2=\chi_h(g)^{{ip_{i,j}\over m_h}}
\chi_{h'}(g)^{\frac{jp_{i,j}}{ m_{h'}}}.$$
La d\émonstration est standard (cf.lemme 1.4.7) et repose sur les hypoth\èses (3) et (4).

\noi Supposons $i\≥1,$ pour $x\in E_{2,0}$ l'application d\éfinie par:
$$(adx)^i\quad E_{-i,j}\rightarrow E_{i,j}$$
est  bijective en $x_0,$  \à l'invariant relatif det($ad(x)^i/_{E_{-i,j}})$ est associ\é le caract\ère: $det(g/_{E_{i,j}}). det(g^{-1}/_{E_{-i,j}})$ donc par l'hypoth\èse
(3), celui-ci est une puissance de $\chi_h(g).$ Le coefficient se calcule
simplement \à l'aide de l'\él\ément $h_h(t)$ qui op\ère sur chaque $E_{i,j}$
par l'homoth\étie de rapport $t^i$ ce qui donne :
$$( b)\quad det(g/_{E_{i,j}}).
det(g^{-1}/_{E_{-i,j}})=\chi_h(g)^{{ip_{i,j}\over m_h}}.$$
Si $i\≤-1,$  par dualit\é donn\ée par $B$ on a :
$$det(g/_{E_{i,j}})=det(g^{-1}/_{E_{-i,-j}})$$
ainsi la relation (b) est v\érifi\ée pour $i\not=0.$

\noi On proc\ède de m\ême avec l'hypoth\èse (4) d'o\ù  
$$( c)\quad \hbox{pour}\quad j\not=0\quad det(g/_{E_{i,j}}).
det(g^{-1}/_{E_{i,-j}})=\chi_{h'}(g)^{\frac{jp_{i,j}}{ m_{h'}}},$$
ainsi de (b) et (c) on d\éduit  la relation (a).

\noi Exactement de mani\ère analogue, en raison des hypoth\èses (1) et (2), on
a pour $1\≤p\≤p_0$ :
$$(d)\quad \hbox{pour}\quad  g\in G\quad  (det(g/_{\goth
g_p})^2=\chi(g)^P\quad \hbox{avec}\quad P=\frac{2p\text{dim}(\goth g_p)}{ m}$$
D'apr\ès (a) on a pour $g\in G_h$ :
$$  \chi_h(g)^{M_h}\chi_{h'}(g)^{M_{h'}}=\chi(g)^P
\quad \hbox{avec}\quad M_h=\frac{\text{ trace}(ad(h)/\goth g_p)}{ m_h}
\quad M_{h'}=\frac{ \text{trace}(ad(h')/\goth g_p)}{m_{h'}}$$
d'o\ù $M_h\not=0,M_{h'}\not=0$ et la relation :
 $$( e)\quad \hbox{pour}\quad i\not=0\quad
\hbox{et}\quad g\in G_h\quad (det(g/_{E_{i,j}}))^{2M_{h'}}=\chi_h(g)^{M_{i,j}}\chi(g)^{\frac{jp_{i,j}}{m_{h'}}P} $$
 $$\hbox{avec}\quad M_{i,j}=\frac{p_{i,j}}{
m_hm_{h'}}(iM_{h'}m_{h'}-jM_hm_h).$$
 2) On a $\∆_{P_h}(N_h)=1$ et pour $g\in H_h$ :
$$ \begin{array}{ccl}{\∆_{P_h}}^{-1}(g)&=&\mid
\prod_{i=1}^{p_0}det(g/_{E_{i,-i}})\mid\\
&=&\mid \chi_h(g)\mid^r \quad \hbox{en utilisant la relation (e) et
le fait que }\ \chi(H)=1.\end{array}$$
3)  V\érifions que $H_{x_0}$ est unimodulaire. Pour ceci il suffit (comme
pour $P_h$) de calculer pour $g\in H_{h,x_0}$ la quantit\é:
$$\prod_{i=1}^{p_0}\mid det(g/_{(E_{i,-i})_{x_0}}\mid.$$
Comme   $g\in H_{h,x_0}$ on a $\chi_h(g)=1$ donc pour tout $i\not=0$
on a
$$\mid det(g/_{E_{i,j}})\mid=1.$$
Or pour $i\≥1$ on a $[x_0,E_{i,-i}]=E_{2+i,-i}$ et $gx_0=x_0$ donc 
$$1=\mid det(g/_{E_{i,-i}})\mid=\mid det(g/_{(E_{i,-i})_{x_0}})\mid
.\mid det(g/_{E_{2+i,-i}})\mid=\mid
det(g/_{(E_{i,-i})_{x_0}})\mid .\qquad \Box $$ 

\noi Ainsi il existe sur l'espace
homog\ène $P_h /H_{x_0}$ une mesure non nulle, relativement  
invariante \à gauche par $P_h,$ de multiplicateur $\∆_{P_h}^{-1}$
et unique \à une constante multiplicative pr\ès ( \cite{bourbakiint}, chap.7, th\éor\ème 3,$\S 2$,n$^\circ 6$), not\ée $d\dot p.$

 \vskip 3mm

\noi {\bf Lemme 2} {\it On suppose que  :
$$ l'(h)= \sum_{i=1}^{p_0}ip_{i,i}-
 (\sum_{i\≥2}i p_{i,2-i})+ (\sum_{i\≥2} p_{i,2-i}).\frac{\text{trace}(ad(h)/_{ \goth g_{p_0}})}{p_0 \text{dim}(\goth g_{p_0})}\ \≥0 ,$$
alors pour $f\in S(\goth g_1)$ et $g\in G$ 
$$\psi_f(g)=\int_{P_h /H_{x_0}}f(gg'x_0)\∆_{P_h}^{-1}(g')d\dot{g'}$$
converge absolument et v\érifie pour tout $p\in P_h$ et $g\in G$ 
$$\psi_f(gp)=\∆_{P_h}(p)\psi_f(g).$$}
 
\vskip 5mm
 \dem  Il suffit de montrer la convergence absolue, la relation
d'invariance \étant \évidente.

\noi On consid\ère l'action du groupe $P_h$ sur l'espace vectoriel:
$$E=E_2(h)\cap \goth g_1\oplus
E'\quad \hbox{avec}\quad E'=\oplus_{i\≥3}E_i(h)\cap \goth g_1$$ Un
calcul facile montre que:  $$P_h.x_0=H_h.x_0\oplus E'$$
Ainsi pour la topologie induite par $\F$ sur $E,$ $P_h.x_0$ est
une orbite ouverte de $E,$ de m\ême $P_h$ est un sous-groupe ferm\é
du groupe des automorphismes d'espace vectoriel de $\goth g,$ muni de la
topologie naturelle induite par $\F.$ Pour cette topologie $P_h$
op\ère continuement, est un groupe localement compact, d\énombrable \à
l'infini donc $P_h /H_{x_0}$ est hom\éomorphe \à l'orbite $ 
P_h.x_0$ (\cite{bourbakitop}, chapitre IX, prop.6,$\S 5,n^\circ3$).

\noi On consid\ère sur $E$ la mesure $P_h$-invariante donn\ée par 
$$\int_Ef(x)\mid F_h(projection(x)/_{E_2(h)\cap \goth
g_1}\mid^{-\delta}dx$$
$dx$ \étant la mesure sur $E$ invariante par translation, ``projection$(x)/_{E_2(h)\cap \goth
g_1}$" repr\ésentant la composante de $x$ suivant $E_2(h)\cap \goth
g_1.$  En effet le
changement de variable fait apparaitre la quantit\é
$$\prod_{i\≥2}\mid det(g)/_{E_{i,2-i}}\mid^{-1}=\mid \chi_h(g)\mid^{-\delta}$$
avec pour $1\≤p\≤p_0$ :
$$\delta=\frac{pdim(\goth g_p)}{ m_h(trace(ad(h')/_{\goth g_p})}.[
trace(ad(h)/_E)-\frac{dim(E)}{ pdim(\goth g_p)}trace(ad(h)/_{\goth
g_p})]$$ (relation (e) de la d\émonstration du lemme 1).

\noi En raison de l'hom\éomorphisme et de l'unicit\é (\à une constante
multiplicative pr\ès) de la mesure $ P_h$ invariante sur
 $P_h/_{{ P_h}_{x_0}}$ on a :
$$\psi_f(g)=\int_{{P_h}.{x_0}}f(x)\mid F_h(projection(x)/_{E_2(h)\cap \goth
g_1})\mid^{r-\delta}dx$$
d'o\ù cette derni\ère int\égrale converge pour $f\in S(\goth g_1)$ lorsque
$r-\delta\≥0$ d'o\ù la condition puisque trace$(ad(h')/\g_{p_0}>0.$\fdem

\vskip 3mm
On proc\ède comme dans \cite{rallisschiffmann}. $H$ est un groupe semi-simple, connexe
pour la topologie de Zariski donc $H$ est unimodulaire. Sur l'espace
${\it C}_K(H)$ des fonctions \à support compact, on consid\ère la mesure de
Haar $\mu; $ par la proposition 3 de \cite{bourbakiint} (chap.VII,$\S 2,n^\circ1$) il existe
une forme lin\éaire positive, non nulle, $H$-invariante \à gauche,
relativement born\ée, not\ée $\nu,$ d\éfinie sur l'ensemble $\it C$ des
fonctions continues d\éfinies sur $H$ \à valeur dans $\C$ telles que
$$\forall g\in H\ ,\  \forall p\in P_h\quad f(gp)=\∆_{
P_h}(p)f(g)$$
(rappelons que $\displaystyle H/_{P_h}$ est compact pour la
topologie induite par $\F$)
\vskip 3mm
 \noi {\bf Lemme 3} {\it On suppose que $l'(h)\≥0$ et pour $f\in S(\goth g_1)$ on pose:
$$\nu_{x_0}(f)=\nu(\psi_f)\quad \hbox{alors }$$
1) $\nu_{x_0}$ est une mesure positive, temp\ér\ée sur $\goth g_1,$
invariante par H et de support $ H.x_0$ que l'on peut \écrire sous la forme
$$\nu_{x_0}(f)=\int\int_{V_h\times P_h.x_0}f(v.w)\mid
F_h(projection(w)/_{E_2(h)\cap \goth g_1})\mid ^{l(h)}dvdw$$
avec $V_h=exp(\oplus_{i\≥1}E_{-i}(h)\cap \goth g_0,$ $V_h$ est    muni d'une
mesure de Haar dv et:$$l(h)=\frac{p_0\text{dim}(\g_{p_0})}{m_h\text{trace}(adh'/_{\g_{p_0}})}.\ l'(h).$$

2) Lorsque $H.x_0=G'.x_0,$ $G'$ \étant un sous-groupe ouvert de $G$
contenant $Ker\chi,$ $\nu_{x_0}$ est relativement invariante par $G'$ de caract\ère
$\mid \chi\mid^{-\alpha}$ avec $$\alpha=\frac {p_0dim(\goth g_{p_0})}{ m
(trace(ad(h')/_{\goth
g_{p_0}})}.[\sum_{i\≥1}ip_{i,i}-\sum_{i\≥1}ip_{2+i,-i}   ]\quad
(1\≤p\≤p_0)\ .$$ } \vskip 3mm
  
  \dem
1) Comme:
$$[\goth g_0,\goth g_0]=(\oplus_{i\≥1}E_{-i}(h)\cap \goth g_0)\oplus 
([\goth g_0,\goth g_0]\cap E_0(h)) \oplus (\oplus_{i\≥1}E_i(h)\cap \goth
g_0)$$ et que $\oplus_{i\≥1}E_{-i}(h)\cap \goth g_0$ est l'alg\èbre de Lie
du groupe alg\ébrique connexe $V_h$ et 
$[\goth g_0,\goth g_0]\cap E_0(h) \oplus_{i\≥1}E_i(h)\cap \goth
g_0$ celle du groupe alg\ébrique connexe $P_h,$ par la
 proposition 4 du chapitre VI, $\S 2$ p.373
de \cite{chevalley}
on a
$H=V_h.P_h$ \à un ensemble de mesure nulle pr\ès  et on peut donc
\écrire pour toute fonction $f:H\mapsto \C$ continue \à support compact: $$\int_Hf(h)dh=\int\int_{V\times P_h}f(v.p)\∆_{
P_h}^{-1}(p)dvdp=\nu(\int_{ P_h}f(g.p)\∆_{P_h}^{-1}(p)dp)$$
(\cite{bourbakitop}, $\S 2, n^\circ 9$,prop.3, $n^\circ 2$, prop.3)
d'o\ù la forme donn\ée pour $\nu_{x_0}.$
\vskip 2mm
2) Soit $g_0$ appartenant \à $G',$ on pose $\nu'(f)=\nu_{x_0}(f(g_0.))$ alors $\nu'$
est encore une mesure positive $H$ invariante de m\ême support que $\nu_{x_0},$
elle lui est donc proportionnelle. En raison des hypoth\èses (1) et (2),
on peut supposer que $\goth g_1$ et $\goth g_{-1}$ engendrent $\goth
g,$ ainsi le centre de $\goth g_0$ est de dimension un. Il suffit donc
de calculer le caract\ère sur les homoth\éties $c(t),t\in \F^*.$\fdem
\vskip 5mm
\noi Lorsque $h$ est $1$-simple tr\ès sp\écial, (3) et (4) sont alors v\érifi\és, cependant $l'(h)<0$ en g\én\éral sauf dans des cas orthogonaux tr\ès  particuliers ou bien: \\ \\

\noi{\bf Lemme 4} {\it Soit $h$ un \elt $1$-simple tr\ès sp\écial alors $l'(h)\≥0$ lorsque:

$\bullet$ $\g_2=\{0\}$,

$\bullet$ Lorsque l'\irf est de degr\é $4$ et $\overline \g $ est de type exceptionnel.}

\dem  Comme la condition porte sur des dimensions, il suffit de faire la v\érification lorsque $\g$ est d\éploy\ée.\\

1) Lorsque $\g_2=\{0\}$, on reprend les notations du $\S 6.1.1.$ On peut supposer que $h=\sum_{1\≤i\≤p}h_{\lambda_i}$ avec $1\≤p\≤n-1$ d'o\ù:
$$d_1=1\ ,\ p_{2,2}=0\ ,\ p_{1,1}=p(n-p)\overline d\ ,\ p_{2,0}=p+p(p-1)\frac{\overline d}{2}\ ,\ \text{dim}(\g_1)=n+n(n-1)\frac{\overline d}{2}$$
donc $l'(h)=\ds\frac{2p(n-p)}{n}((n-p+1) \frac{\overline d}{2}-1)\≥\frac{2p(n-p)}{n}(\overline d-1)\≥0.$

2) Lorsque l'\irf est de degr\é $4$ et $\overline \g $ est de type exceptionnel (donc distinct de $G_2$ par le choix de $h$),  par classification on a soit $\g_2$ de dimension $1$ et ce cas sera fait dans la proposition suivante, soit $(\overline \g_0,\overline \g_1)$ est de type $(E_7,\alpha_6)$ (cf.tableau 3) mais alors:
\vskip 2mm

\noi $p_{1,1}=0\ ,  \ p_{2,0}=\ds\frac{\text{dim}(\g_1)}{2}=2p_{2,2}=16\ ,\ p_{4,0}=p_{0,4}=1$  (cf.$\S 8.1.1)$) d'o\ù $l'(h)=0.$\fdem\\

\bigskip

\noi {\bf Proposition 1}  {\it   Soit $(\goth g_0,\goth g_1)$ un pr\éhomog\ène  tel que:

1) $\goth g_1$ est un  $\goth g_0-$module absolument simple,

2)  $2H_0$ est 1-simple et 

soit

5) $\g_2=\{0\},$

soit

6) $\overline g$ est de type exceptionnel  distinct de $G_2$ et
$\overline g_2$ est de dimension $1,$ \\

 \noi alors pour tout $sl_2$-triplet 1-adapt\é,  $(x_0,h,y_0)$  avec $h\not=2H_0,$ il existe une mesure positive,
temp\ér\ée, invariante par H et de support $H.x_0$ (resp. 
$H.y_0 ).$ 

\noi Pour tout $f\in S(\goth g_1)$ (resp.  $f\in S(\goth g_{-1})),$
on peut l'\écrire : $$\nu_{x_0}(f)=\int_{V_h}[\int_{H_h.x_0}f(v.z)\mid 
F_h(z)\mid^{l(h)} dz]dv$$ $$ (\quad resp. \ \
\nu^*_{y_0}(f)=\int_{N_h}[\int_{H_h.y}f(n.z)\mid  F^*_h(z)\mid^{l(h)}
dz]dn\quad )$$ 

\noi avec 
$V_h=exp(\oplus_{i\≥1}E_{-i}(h)\cap \goth g_0)$  (resp.
$N_h=exp(\oplus_{i\≥1}E_i(h)\cap \goth g_0))$  muni d'une  mesure de
Haar dv  (resp. dn) 

\noi $H_h$ \étant le centralisateur de h dans H et $F_h$
 (resp. $F^*_h$)
\étant le polynome d\éfini par $F_h$(projection de x sur
 $E_2(h)\cap \goth g_1)$  (resp. $F^*_h$(projection de y sur
 $E_{-2}(h)\cap \goth g_{-1})$ et:\\

$ l(h)=(1-\widetilde B(2H_0-h,H_0)) {\overline
d\over 2}-1$ dans le cas commutatif,

$ l(h)=(1-\widetilde B(2H_0-h,H_0))\frac { 
d} {2}+\begin{cases}1\text{ lorsque }\widetilde B(h,H_0)=-1,\\ -1\text{ lorsque }\widetilde B(h,H_0)=-2,\\
0 \text{ lorsque }\widetilde B(h,H_0)=-3\end{cases}$ lorsque dim$(\g_2 )=1.$}

\dem   Il suffit de faire le calcul dans le cas o\ù $\g$ est d\éploy\ée.  

\noi 1) Dans le cas commutatif (i.e. 1)2)5)), tout $h$ 1-simple distinct de $2H_0$ est $1$-simple tr\ès sp\écial (remarque 1.1.3) et on applique le lemme pr\éc\édent en notant que pour $h=\sum_{1\≤i\≤p}h_{\lambda_i}$ on a $l(h)=(n-p+1)\frac{\overline d}{2}-1$ et $\widetilde B(h,H_0)=-p.$\\

2) Dans les cas v\érifiant 1)2)6), ce qui correspond \à la description donn\ée au $\S 8.2.1$ que l'on reprend ici, on peut supposer que $h=\sum_{1\≤i\≤p}h_{\lambda_i}$ avec $1\≤p\≤3.$  

\noi On rappelle que $2H_0=\sum_{1\≤i\≤4}h_{\lambda_i}$ et que pour tout $i, j$ distincts et compris entre $1$ et $4,$ il existe un \elt du groupe de Weyl de $\∆_0$ qui permute les racines $\lambda_i$ et $\lambda_j$ donc les propri\ét\és 3) et 4) sont automatiquement v\érifi\ées puisque les \PVs sont absolument irr\éductibles et on a $p=-\widetilde B(h,H_0).$ 
 
 \noi Il reste \à v\érifier  que la quantit\é $l(h)$  est
positive ce qui se fait \à l'aide des d\écompositions donn\ées en 8.2.1\fdem\\

\bigskip

\noi {\bf Remarques:}  1) Lorsque le \PV est de type $(G_2,\alpha_2)$ d\éploy\é, on a $l(h)<0$ pour tout $h$ $1$-simple distinct de $2H_0$ donc la mesure relativement invariante dont le support est inclus dans $\g_1-\g'_1$ ne prend pas cette forme simple (cf. la distribution $\Sigma_1$ de \cite{formescubiques}). \\

2) Lorsque les hypoth\èses de la proposition 1 sont satisfaites, on peut v\érifier que l'application :$\begin{array}{cccl}&\g_1&\mapsto &\Q\\ &x&\rightarrow &\begin{cases} B(h,H_0), \ (x,h,y)\text{ \étant un  \Sl 1-adapt\é} ,\\
0\text{ lorsque }x=0,\end{cases}\end{array}$ 

\noi est semi-continue inf\érieurement.\\

\noi De plus on peut montrer que pour tout  $x_0\in \g_1$ on a:

 $$\overline{G.x_0}=\{0\}\cup G.x_0\cup_{u\in \goth I}G.u \ ,\ 
 \overline{H.x_0}=\{0\}\cup H.x_0\cup_{u\in \goth I'}H.u\ ,$$
\noi avec $\goth I=\{u\in \g_1\ |\ $ pour lesquels il existe $2$ \Sls 1-adapt\és qui commutent $(u,h,v)$ et $(u',h',v')$ tels que $u+u'\in G.x_0\}$ et pour $\goth I'$ on a la m\ême d\éfinition mais avec $u+u'$ et $x_0$ dans la m\ême orbite de $H.$

\noi Dans le cas $\goth p$-adique, on peut v\érifier classiquement (cf.th\éor\ème 2.3 de \cite{rubsln}) que les mesures $\nu_{x}, $ $x$ d\écrivant un ensemble de repr\ésentants des orbites de $H$ dans $\g_1-\g'_1$ est une base des distributions $H$- invariantes de support inclus dans $\g_1-\g'_1$ lorsqu'on pose $\nu_0(f)=f(0).$\\

3) Soit $G'=\begin{cases}\{g\in G\ |\ \chi (g)\in f(E)^*\}\text{ dans le cas commutatif de type I }\\
\qquad \qquad \qquad \qquad\text{( cf. remarque $\S 6.1.6)$ pour la d\éfinition), }\\
G\quad \text{    dans tous las autres cas,}\end{cases}$\\

\noi alors $H.x_0=G'.x_0$ (th.4.3.2 et remarque 2) du th.4.2.3 de \cite{mullerNAG} et $\S 6.1.2$ pour le type III commutatif) donc les mesures $\nu_{x_0}$ et $\nu^*_{x_0}$ sont relativement invariantes par $G'$ de caract\ère associ\é respectivement $|\chi|^\alpha $ et $|\chi|^{-\alpha} $ avec $\alpha=-\widetilde B(h,H_0)\frac{\overline d}{2}$ dans le cas commutatif et 

\noi $\alpha=-\widetilde B(h,H_0)\frac{  d}{2}+\begin{cases}1\text{ lorsque }\widetilde B(h,H_0)=-3,\\
\frac{1}{2} \text{ sinon}\end{cases}.$
\bigskip

\noi {\bf 2  Une  application aux mesures singuli\ères dans le  cas commutatif 
avec $\overline d$ pair} \vskip 5mm 

\noi On reprend toutes les notations du $\S 6$ et on donne la d\éfinition suivante de rang:\\

\noi{\bf D\éfinition}   {\it Soit $(x,h,y)$ un $sl_2$-triplet $1$-adapt\é,  on appelle ``rang de $x$ ", ``rang de $h$ "et ``rang de $y$ " la quantit\é $r(h)=-\ds\frac{\widetilde B(h,H_0)}{d_1}.$

\noi On pose $r(0)=0.$}\\

\noi Soit $(x_0,h,y_0)$ un $sl_2$-triplet $1$-adapt\é  associ\é \à une orbite de rang $p$ avec $1\≤p\≤n-1,$ alors on  a (\à l'action de $G$ pr\ès):  $$2H_0=\sum_{1\≤i\≤n}h_{\lambda_i}\ ,\  2H_0-h=\sum_{i=1}^{n-p}h_i\ \ ,\ h=\sum_{i=n-p+1}^nh_i \ ,  $$ la proposition 1 permet d'\écrire les mesures sur
les orbites singuli\ères de la mani\ère suivante : $$\forall f\in S(\goth
g_1)\quad \nu_{x_0}(f)=Z^{(p)}_O(T^h_f(\ . +0),{\overline d\over 2}-1)$$
$$\forall f\in S(\goth g_{-1})\quad
\nu^*_{y_0}(f)=Z^{*(n-p)}_{O^*}(T^{*(2H_0-h)}_f(0+\  .),{\overline d\over
2}-1)$$ o\ù $T^h_f$ d\ésigne la transformation introduite dans le $\S 4.3$ sur
$S(\goth g_1)$ pour le couple d'\él\éments $1$-simples : $(h,2H_0-h),$ et 
$T^{*(2H_0-h)}$  est la transformation introduite dans le  $\S 4.3$ sur
$S(\goth g_{-1})$ pour le couple d'\él\éments $1$-simples : $(2H_0-h,h)$

\noi $Z^{(p)}$ (resp. $Z^{*(n-p)}$ ) \étant les fonctions Z\éta du  \PV
$(H_h,E_2(h)\cap \goth g_1)$ (resp.$(H_h,E_{-2}(h)\cap \goth g_{-1})$
), $O=H_h.x_0$ et $O^*=H_h.y_0,$ les normalisations des mesures \étant choisies en concordance avec le  $\S 4.$

\noi  Le th\éor\ème 4.3.5 et les th\éor\èmes 6.2.1 et 6.2.2 permettent
d'\établir une g\én\éralisation du lemme 4.2 de \cite{sato-shintani} lorsque $\frac{\overline d}{2}\≥0.$\\
   \vskip 5mm
\noi{\bf Proposition 2 }
 {\it On suppose $\overline d$ pair.  
 
\noi Soit f une fonction de Schwartz
de support inclus dans $ \goth g'_1=\{x\in \goth g_1\ \ F(x)\not=0\}$
alors :
 $$ \nu^*_{y_0}(\hat
f)=K.C_p\tilde\omega_{\delta^{n-p+1}}(F^*_h(y_0)) \int_{\goth
g_1}f(x)\tilde\omega_{\delta^p}(F(x)) | F(x)|^{-{\overline d\over 2}pd_1}dx\quad
\hbox{avec }$$ 
$$C_p= \gamma(f)^{np+{p(p-1)\over 2}}
\prod_{j=0}^{p-1}b(\tilde\omega_{\delta^{j+1}}|\ |^{\frac{\overline d}{2}(1 +j  d_1)})$$
$ b(\pi)$ \étant d\éfini par 

\quad a)  Lorsque  le rang s\épare
les orbites de G dans $\goth g_1$ (alors $\delta=1)$ ou bien lorsque $\F$ est un corps
$\goth P$-adique ou bien lorsque $n=2$: $$b(\omega|\ |^s)=a^{(1)}( \omega|\ |^s)\quad \text{(def.th.6.2.1)}$$
\quad b) Lorsque $\F=\R$ et f est anisotrope  
 ($\delta=(-1)^{d/2},$ $\gamma(f)=(\sqrt{-1})^{d/2}$):
$$b(\tilde\omega_{\pm 1}|\ |s)=\rho (|\ |^s;-1)=({\sqrt{-1}\over 2\pi})^s. \ \Gamma(s)$$
K \étant la constante d\éfinie par 

\quad i)  Lorsque  le rang s\épare
les orbites de G dans $\goth g_1$ : K=1
\vskip 2mm
\quad ii)  Lorsque $\F=\R$ et f est anisotrope :
$$K=2^{\frac{1+(-1)^n}{2}}C_p^r $$ r d\ésignant le nombre de
 composantes positives d'un repr\ésentant de $H_h.y_0\cap
\oplus_{n-p+1\≤j\≤n}\goth g^{-\lambda_j}.$
\vskip 2mm
\quad iii)  Lorsque $\F$ est un corps $p$-adique et f est anisotrope
 (cas $d=\overline d=2$) on a $$K=\begin{cases} 1\text{ lorsque }n\text{ est pair et }p\text{ impair}\\
 \frac{1}{2}\text{ sinon.}\end{cases}$$}

  \vskip 3mm
 \dem
 1)  Soit $f$ une fonction de Schwartz
de support inclus dans $\goth g'_1=\{x\in \goth g_1\ \
F(x)\not=0\},$ montrons que:
 $$ \nu^*_{y_0}(\hat f)=\sum_{\{O\ orbites\ de\ G'\ dans\ 
{\goth g'_1}\}}a_O\int_Of(x)
\mid F(x)\mid^{-{\overline d\over
2}pk}dx\quad \hbox{avec }$$
$G'=\{g\in G\ |\ \chi(g)\in f(E)^*\}$ et
$$a_O=a_{H_h.{y_0},O'}(\tilde\omega_{\delta^{n-p}};{\overline d\over
2}-1)\gamma(f)^{p(n-p)}(\delta^p,F(x)).(\delta^n,F_h(x)),$$
  $x$ \étant un \él\ément de $O\cap\oplus_{1\≤i\≤n}\goth g^{\lambda_i},$ $O'$
\étant une orbite de $H_h$ dans $E_2(h)\cap \goth g_1$ contenant la
projection de $x$ sur  $E_2(h)\cap \goth g_1.$ \\

\noi 
 Soit $f_g$ la fonction d\éfinie par $f_g(x)=f(gx),$ comme: 
 
$$\nu^*_{y_0}(\hat {f_g})=\nu^*_{y_0}({\hat f}_g)\mid
\chi(g)\mid^{ -N} 
=\nu^*_{y_0}(\hat f)\mid
\chi(g)\mid^{ -N+\frac{\overline d}{2}pd_1} \text{ avec  }N=\frac{\text{dim}(\g_1)}{\text{degr\é}(F)}=\frac{\overline d}{2}(nd_1-1)+1,$$
si $f$ est une fonction \à support dans l'orbite ouverte $O$ de $G'$ dans
$\goth g_1,$ on a pour une constante convenable :
$$\nu^*_{y_0}(\hat f)=a_O\int_Of(x)\mid F(x)\mid^{-{\overline d\over
2}pd_1}dx$$
Il suffit de calculer la constante $a_O.$ Or par d\éfinition de 
$\nu^*_{y_0}$ et par le 2)
du th\éor\ème 4.3.5 on a pour $z\in E_{-2}(h)\cap \goth g_{-1}:$  
$$\nu^*_{y_0}(\hat f)=\int_{H_h.{y_0}}\mid    F^*_h(z)\mid ^{
{\overline d\over 2}-1}[\int_{E_0(h)\cap \goth g_1}\gamma(x,z)\goth
F_2(T_f^{(2H_0-h)}(x+\ ))(z)dx]\ \ dz\ $$
et:  $$\gamma(x,z)= \gamma(f)^{p(n-p)}(\delta^p,
F_{2H_0-h}(x)) (\delta^{n-p},  F^*_h(z))$$ 
(lemme 6.1.10 car $\overline d$ est pair).

\noi
Pour simplifier les calculs, prenons  $f$ \à support dans un
ouvert $O_u,$ $u=(u_1,...,u_n)\in (\FF)^n,$ tel que $O_u\subset O;$ comme on a suppos\é $\overline
d\≥2,$ on peut appliquer le th\éor\ème de Fubini ainsi:
$$\nu^*_{y_0}(\hat f)=\int_{E_0(h)\cap \goth g_1}\ [\
\int_{H_h.{y_0}}\ \gamma(x,z) \mid   F^*_h(z)\mid ^{
{\overline d\over 2}-1}\goth F_2(T_f^{(2H_0-h)}(x+\ ))(z)dz]\ \ dx$$
$$=\gamma(f)^{p(n-p)} .\int_{E_0(h)\cap \goth g_1 } 
\tilde\omega_{\delta^p}(F_{2H_0-h}(x)). \biggl(
\int_{H_h.{y_0} } \mid   F^*_h(z)\mid ^{ {\overline d\over
2}-1} \tilde\omega_{\delta^{n-p}}(  F^*_h(z))\goth F_2(T_f^{(2H_0-h)}((x+\ ))(z)dz\biggr) 
dx.$$  
On applique ensuite l'\équation fonctionnelle au \PV $(H_h,E_{-2}(h)\cap \goth g_{-1})$ avec les
$sl_2$-triplets $h_{p+1},...,h_n$ ce qui donne:
$$\nu^*_{y_0}(\hat f)= \gamma(f)^{p(n-p)}\sum_{O'}a_{H_h.{y_0},O'}(
\tilde\omega_{\delta^{n-p}};{\overline d\over 2}-1).$$
$$\int_{E_0(h)\cap \goth g_1 }\int_{O'}T_f^{(2H_0-h)}(x,y)
\tilde\omega_{\delta^p}(F_{2H_0-h}(x)) \tilde\omega_{\delta^{n-p}}(F_h(y))
\mid F_h(y)\mid^{-{\overline d\over 2}pd_1} dxdy.$$
En raison du support choisi pour $f,$ cette somme se simplifie et donne :
$$\nu^*_{y_0}(\hat f)=a_O  \int_{O_u}f(x)\mid F(x)\mid^{-{\overline d\over
2}pd_1}dx \text{  
avec  }$$
$$a_O=\gamma(f)^{p(n-p)}.\ \tilde\omega_{\delta^p}(\prod_{1\≤i\≤n}u_i).\ 
\tilde\omega_{\delta^{n}}(u_{n-p+1}...u_n). a_{H_h.{y_0},O'}(
\tilde\omega_{\delta^{n-p}};{\overline d\over
2}-1) ,$$
$O'$ \étant l'orbite de $H_h$ dans $E_2(h)\cap \goth g_1$ contenant
$O_{u"},$ avec $ u" =(u_{n-p+1},...,u_n),$ d'o\ù le r\ésultat.

2) Il ne reste plus qu'\à \évaluer le coefficient $a_{H_h.{y_0},O'}(
\tilde\omega_{\delta^{n-p}};{\overline d\over
2}-1)$
 ce qui  conduit \à des r\ésultats simples  en raison  de la valeur tr\ès particuli\ère de ${\overline d\over
2}-1.$

 \noi a) Dans le cas transitif pour $G,$ alors le rang s\épare les orbites de $G$ dans $\g_1,$ ce qui correspond aux types II avec $e=0,$ III r\éel, III $A_n$ en $\goth p$-adique et I ($d=4$ en $\goth p$-adique), $G=G'$  
et $H_h$ a \également une seule orbite ouverte dans 
$E_2(h)\cap \goth
g_1,$ (resp. dans $E_{-2}(h)\cap \goth
g_{-1}$ ) on applique alors le th\éor\ème 6.2.1. Notons que dans ce
cas $\gamma=\pm 1$ et $\delta=1$  (1) du lemme 6.1.10).\\

\noi b) Il reste le cas de type I et II avec $e>0$ (alors $d_1=1$ et $n=2$).\\

\noi Notons que pour:
 
\quad i) $n$  pair, on a $G=G'$ et les orbites de $H_h$
dans $E_2(h)\cap \goth
g_1$ (resp.dans $E_{-2}(h)\cap \goth
g_{-1})$  sont celles du  \PV $(\goth U(\F (2H_0-h)_0,
\goth U(\F (2H_0-h)_1),$  ce qui termine le cas II car alors $H_h$ a  une seule orbite ouverte dans 
$E_2(h)\cap \goth
g_1,$ (resp. dans $E_{-2}(h)\cap \goth
g_{-1}$ ) donc: $$ a_{H_h.{y_0},O'}(
\tilde\omega_{\delta^{n-p}};{\overline d\over
2}-1)=\rho'(\tilde\omega_{\delta }|\ |^{\overline d\over
2}) \text{ d'o\ù  }a_O=\gamma(f).\ \tilde\omega_{\delta}(u_1u_2).\rho'(\tilde\omega_{\delta }|\ |^{\overline d\over
2}).$$
 \vskip 2mm
\quad ii) $n$   impair, on a   $G'\not=G$, $\sum_{1\≤i\≤n}a_ix_i$ et 
$\sum_{1\≤i\≤n}b_ix_i$ sont dans la m\ême orbite de $G'$ si et seulement
si les formes quadratiques $f^a=\oplus_{1\≤i\≤n}a_if$ et $f^b=\oplus_{1\≤i\≤n}b_if$ sont \équivalentes ainsi que l'analogue pour le \PV $(\goth U(\F (2H_0-h)_0,
\goth U(\F (2H_0-h)_1)$.\\

\noi  Il reste deux situations possibles:\\

$\alpha$) Le cas $\goth p$-adique restant qui correspond \à $d=\overline d=2$ et $\delta\not=1$ ($ f$ est anisotrope).\\

\noi Lorsque $n$ est pair et $p$ est impair, il y a \également une seule orbite ouverte de $H_h$ dans $\goth U(\F h')_{\pm1}$ et $a_{H_h.y_0,O'}(\tilde\omega_{\delta^{n-p}};0)=\gamma(f)^{p(n-p)+\frac{p(p-1)}{2}}\prod_{1\≤j\≤p}\rho' (\tilde\omega_{\delta^{j}}|\ |^j)$ (2) du th.6.2.2).

\noi Dans les cas restants: 
$$\begin{array}{lll}a_{\epsilon',\epsilon}(\tilde\omega_{\delta^{n-p}})&=&(\epsilon'.\epsilon)^{n-p}a_{\epsilon',\epsilon}(Id)\\
&=&\frac{1}{2}\gamma(f)^{\frac{p(p-1)}{2}}(\epsilon'.\epsilon)^{n-p}\ .\ \epsilon^{p-1}\biggl(\prod_{1\≤j\≤p}\rho' (\tilde\omega_{\delta^{j-1}}|\ |^j)+\epsilon'.\epsilon
\prod_{1\≤j\≤p}\rho' (\tilde\omega_{\delta^{j}}|\ |^j)\biggr)\\
&=&\frac{1}{2}\gamma(f)^{\frac{p(p-1)}{2}}(\epsilon')^{n-p+1}.\epsilon^{n}\prod_{1\≤j\≤p}\rho' (\tilde\omega_{\delta^{j}}|\ |^j),\end{array}$$
puisque $\rho'(|\ |)=0$ (cf. $\S 6.2.2$, remarque),
d'o\ù en rempla\çant   $\epsilon=\tilde\omega_\delta(u_{n-p+1}...u_n)$ et $\epsilon'=(\delta,F^*_h(y_0))$ on a:$$ a_O=\frac{1}{2}\gamma(f)^{p(n-p)+\frac{p(p-1)}{2}}\tilde\omega_{\delta^p}(F(x))  \tilde\omega_{\delta^{n-p+1}}(F^*_h(y_0))\prod_{1\≤j\≤p}\rho' (\tilde\omega_{\delta^{j}}|\ |^j) \text{ pour }x\in O_u. $$

$\beta)$
  Le cas r\éel  alors $f$ est d\éfinie positive, $\delta=(-1)^{d/2}$ et
$\gamma(f)= (\sqrt{-1})^{d/2}.$ 

\noi Par le 1) du th\éor\ème 6.2.2, on doit \évaluer la quantit\é $a_O=\gamma(f)^{p(n-p)+\frac{p(p-1)}{2}}\tilde\omega_{\delta^p}(\prod_{1\≤i\≤n}u_i).b$ avec:
$$ b=\sum_{\{w\in \{-1,1\}^p\ |\ O^*_w\subset H_h.y_0\}}   (\delta,\prod_{1\≤j\≤p}u_{n-p+j}^{n-p+j}\ .\ w_{n-p+j}^{j-1}) \ \prod_{1\≤j\≤p}\rho'(\tilde\omega_{\delta^{n-p}}|\ |^{j\frac{d}{2}};-u_{n-p+j}w_{n-p+j}) .$$
Or pour $s\in \C,$ $\rho'(\tilde\omega_{\pm 1}|\ |^s;x)= \tilde\omega_{\pm 1}(-x)\rho(|\ |^s;x)$ (cf.formule donn\ée en 3) de la remarque 3.6.6) donc:
$$ b=(\delta^{n-p+1},F^*_h(y_0))\sum_{\{w\in \{-1,1\}^p\ |\ O^*_w\subset H_h.y_0\}}   \prod_{1\≤j\≤p}
( (-1)^{j\frac{d}{2}}, u_{n-p+j} .\ w_{n-p+j} )\  \rho( |\   |^{j\frac{d}{2}};-u_{n-p+j}w_{n-p+j}) ,$$
 De plus pour $n\in \N^*,$ il est imm\édiat que $\rho (|\ |^n;x)=((-1)^n,-x)\rho(|\ |^n;-1)$ (v\érification sur la formule explicite) d'o\ù (on rappelle que $d$ est pair):
$$ b=(\delta^{n-p+1},F^*_h(y_0))  \prod_{1\≤j\≤p}
 \rho( |\ |^{j\frac{d}{2}};-1)  |{\{w\in \{-1,1\}^p\ |\ O^*_w\subset H_h.y_0\}}|  .$$
 Soit $(p-r,r), 0\≤r\≤p,$ la signature de la forme quadratique associ\ée \à un repr\ésentant de l'orbite de $y_0$ alors:
 
 $\bullet$ si $n$ est impair, l'orbite de $H_h.y_0$  est caract\éris\ée par la signature donc
  
\noi $ \{w\in \{-1,1\}^p\ |\ O^*_w\subset H_h.y_0\}= \{w\in \{-1,1\}^p\ $ ayant $r$ coordonn\ées $-1$ et $p-r$ coordonn\ées $1\}$ donc  $ |\{w\in \{-1,1\}^p\ |\ O^*_w\subset H_h.y_0\}|=C_p^r.$
 
  $\bullet$ si $n$ est pair, $H_h$ contient l'homoth\étie de rapport $-1$ donc l'orbite de $H_h.y_0$  est caract\éris\ée par la signature $(p-r,r)$ ou bien $(r,p-r)$ donc $ |{\{w\in \{-1,1\}^p\ |\ O^*_w\subset H_h.y_0\}}|=2C_p^r.$\fdem
  
  \bigskip
\noi  {\bf Remarques:}

\noi Lorsque le rang s\épare les orbites de $G$ dans $\g_{\pm 1}, $ on convient de poser pour $1\≤p\≤n-1:$ $\nu_p:=\nu_{x}$ pour $x\in \g_1$  de rang $p$ et $\nu^*_p:=\nu^*_{y}$ pour $y\in \g_{-1}$  de rang $p.$\\

1) Lorsque $\overline d=d$ et $f(E)^*=\F^*$ on a
$C_p= \prod_{j=1}^p \rho( |\ |^{jd\over 2}) $
ce qui donne :
$$C_p=\begin{cases} 0\text{ lorsque }d=2,\\
\\
\ds\prod_{j=1}^p({1-q^{{d\over 2}j-1}\over 1-q^{-{d\over 2}j}})\text{ dans le cas }\goth p-\text{adique avec }d=4\text{ ou }8 ,\\
\\
 2^p(-1)^{{d\over 8}p(p+1)}(2\pi)^{-{d\over 4}p(p+1)}.\prod_{j=1}^p\
\Gamma(j.{d\over 2})\  \text{ dans le cas r\éel avec }d=4\text{ ou }8.\end{cases}$$
Lorsque $d=2,$ H.Rubenthaler a montr\é que  $\widehat {\nu^*_p}$ est proportionnelle \à $\nu_{n-p}$ (\cite{rubsln}).\\

2) Lorsque $\overline d=d$ et $f(E)^*\not=\F^*$ on a:

 \quad i) dans le cas $\goth P-$adique: 
 $$KC_p= K.(\delta,(-1)^{p(p+1)/2}). \gamma(f)^{np+\frac{p(p-1)}{ 2}}.  \prod_{j=1}^{[\frac{p}{ 2}]}\frac{1-q^{2j-1}}{1-q^{-2j}} \  \prod_{j=1}^{[\frac{p+1}{ 2}]}\rho (\tilde\omega_{\delta} |\ |^{2 j-1})\ .$$
\quad ii) dans le cas r\éel lorsque $f$ est d\éfinie positive: 
$$KC_p=2^{1+(-1)^n}. 
(2\pi)^{-{d\over 4}p(p+1)}.(\sqrt{-1})^{\frac{d}{ 2}np}.(\prod_{j=1}^p\
\Gamma(j.{d\over 2})\ ).C_p^r$$
lorsque $d=2,$ au coefficient $2^{1+(-1)^n}.2^{np-\frac{p(p+1)}{2}} $ pr\ès, on
retrouve le r\ésultat du lemme 4.2 de \cite{sato-shintani}. \\

3) Lorsque $\F$ est un corps $\goth p$-adique, on peut montrer que $\widehat {\nu^*_p} $ est proportionnelle \à $Z(\ ;|\ |^{-\frac{d}{2}p})$ pour $p=1,...,n-1$ lorsque $d=4$ ou $d=8.$  \\

4) Lorsque $\overline d=2$ et $d_1\≥2,$ on a encore $C_p=0$ dans le cas r\éel (il apparait le coefficient $\sin (\pi \overline d)=0$) et dans le cas $\goth p$-adique puisqu'il apparait le facteur
$a^{(1)}(|\ |)=(-1)^{d_1-1}\prod_{1\≤j\≤d_1}\rho (|\ |^j)=0.$

\noi La proportionalit\é entre $\widehat {\nu^*_p}$ et $\nu_{n-p}$ est encore vraie dans le cas r\éel (corollaire \à la proposition IV-15 de \cite{rais}).

     \newpage
 \centerline{\bf Tableau 1}
\vskip 5mm
        
\noi  Liste des \PVs de type parabolique absolument irr\éductibles,   r\éguliers et commutatifs \à laquelle il convient d'ajouter les cas BI(n,1) et DI(n,1) du tableau 2.\\

    \noi Le num\éro indiqu\é entre parenth\èse est celui qui lui correspond dans la classification de Sato-Kimura (\cite{satokimura},$\S7,$ table I p.144).\\ \\

\begin{tabular}{| c | c | c | c | c | c | c | }
\hline
 & & & & &  &\\
$(\Delta,\lambda_0)$ & Diagramme de Satake de $\g\otimes_\F\overline F$& ${\overline d}$ & $d_1$ & $n$ & $e$ &Type \\
& & & & & &
\\
\hline
& &  & & && \\
$(C_n,\alpha_n)$ & 
\hbox{\unitlength=0,6pt
\begin{picture}(300,30)(-50,0)
\put(10,40){\circle{10}}
\put(15,40){\line (1,0){30}}
\put(50,40){\circle{10}}
\put(60,40){\circle*{1}}
\put(65,40){\circle*{1}}
\put(70,40){\circle*{1}}
\put(75,40){\circle*{1}}
\put(80,40){\circle*{1}}
\put(85,40){\circle*{1}}
\put(90,40){\circle*{1}}
\put(95,40){\circle*{1}}
\put(100,40){\circle*{1}}
\put(110,40){\circle{10}}
\put(115,40){\line (1,-1){21}}
\put(140,15){\circle{10}}
\put(140,15){\circle{18}}
\put(115,-10){\line (1,1){21.5}}
\put(110,-10){\circle{10}}
\put(110,8){\vector(0,1){23}}
\put(110,22){\vector(0,-1){23}}
\put(100,-10){\circle*{1}}
\put(95,-10){\circle*{1}}
\put(90,-10){\circle*{1}}
\put(85,-10){\circle*{1}}
\put(80,-10){\circle*{1}}
\put(75,-10){\circle*{1}}
\put(70,-10){\circle*{1}}
\put(65,-10){\circle*{1}}
\put(60,-10){\circle*{1}}
\put(50,-10){\circle{10}}
\put(50,8){\vector(0,1){23}}
\put(50,22){\vector(0,-1){23}}
\put(15,-10){\line (1,0){30}}
\put(10,-10){\circle{10}}
\put(10,8){\vector(0,1){23}}
\put(10,22){\vector(0,-1){23}}
\end{picture}
}

& $2$ &$1$ &$n$ & $2$&I  \\
& &  & & & &\\
\hline
& & & & & & \\
$(C_n,\alpha_n)$ & 
\hskip 10pt \hbox to 5 cm {\offinterlineskip \lower 2pt\hbox{$\bullet$}
\hglue -3pt{\vrule height .4pt depth 0pt width 0,5cm}\lower 2pt\hbox{$\circ$}
\hglue -7pt{ \vrule height .4pt depth 0pt width 0,5cm}\lower 2pt\hbox{$\bullet$}
 \dotfill \hbox to 2 cm {\lower 2pt\hbox{$\circ$}\hglue -1,3pt
\hrulefill\lower 2pt
\hbox{$\bullet$}\hrulefill \lower 2pt \vtop {\offinterlineskip \hbox{$\circ$} \hbox to 5pt{\hfill \vrule height 12pt width 0,3pt\hfill} \hbox{$\circledcirc$}}\hglue -3pt\hrulefill\lower 2pt\hbox{$\bullet$} } }

& $4$ &$1$ &$n$ &$4$ & I   \\
& & & & & &\\
\hline
& & & & && \\
$(C_n,\alpha_n)$ &
\hskip 10pt \hbox to 5 cm {\offinterlineskip \lower 2pt\hbox{$\circ$} \hglue -3,2pt
{\vrule height  0,4pt depth 0pt width 0,5 cm}\lower 2pt\hbox{$\circ$}
  \dotfill \hbox to 1,3 cm {\offinterlineskip\lower 2pt\hbox{$\circ$}\kern -3pt   \hrulefill\kern -3pt\lower 2pt 
\hbox{  \offinterlineskip  {$\circ$} \hglue -8pt\vbox{ {\hrule height 0,3pt width 0,6cm}\vskip 3pt{\hrule height 0,3pt width 0,6cm}
\vskip 0,6pt} \hglue -16pt $<$   \hglue -1,8pt{$\circledcirc$}}}}
& $1$ &$1$ &$n$ & $1$&I\\
(2) & & & & &&\\
\hline
& & & & & &\\
$(C_3,\alpha_3)$ &

\hskip 10pt \hbox to 5cm {\lower 2pt\hbox{$\circ$}\hrulefill\lower 2pt
\hbox{$\bullet$}\hrulefill \lower 2pt \vtop {\offinterlineskip \hbox{$\bullet$} \hbox to 5pt{\hfill \vrule height 12pt width 0,3pt\hfill} \hbox{$\bullet$}}\hrulefill\lower 2pt\hbox{$\bullet$}\hrulefill\lower 2pt\hbox{$\circ$}\hrulefill\lower 2pt\hbox{$\circledcirc$} }

& $8$ &$1$ &$3$ &$ 8$&I   \\
&  cas r\éel uniquement& & && & \\
\hline
& & & && &\\
$(A_{2n-1},\alpha_n)$ &
\hbox to 4,5 cm {\offinterlineskip \lower 2pt\hbox{$\circ$} \hglue -3,2pt
{\vrule height  0,4pt depth 0pt width 0,5 cm}\lower 2pt\hbox{$\circ$}
\hglue -3,2pt
 \dotfill \hbox to 1,2 cm  
{\offinterlineskip \lower 2pt\hbox{$\circ$} \hglue -3,2pt
{\vrule height  0,4pt depth 0pt width 0,5 cm}\lower 2pt\hbox{$\circledcirc$}
\hglue -3,2pt
{\vrule height  0,4pt depth 0pt width 0,5 cm}\lower 2pt\hbox{$\circ$}}
 \hskip 5pt   \dotfill \hbox to 0,4 cm 
  {\offinterlineskip \lower 2pt\hbox{$\circ$}  
 \hglue -3,2pt
{\vrule height  0,4pt depth 0pt width 0,5 cm}\lower 2pt\hbox{$\circ$}}}
& $2$ &$1$ &$n$ &$ 0$&II \\
$\subset$  (1) &  & & & &&\\
\hline
& &  & & & &\\
$(D_{2n},\alpha_{2n})$ &  
\hskip 10pt \hbox to 5,3 cm {\offinterlineskip \lower 2pt\hbox{$\circ$}
\hglue -3pt{\vrule height .4pt depth 0pt width 0,5cm}\lower 2pt\hbox{$\circ$}
\hglue -7pt{ \vrule height .4pt depth 0pt width 0,5cm}\lower 2pt\hbox{$\circ$}
 \dotfill \hbox to 2 cm {\lower 2pt\hbox{$\circ$}\hglue -1,3pt
\hrulefill\lower 2pt
\hbox{$\circ$}\hrulefill \lower 2pt \vtop {\offinterlineskip \hbox{$\circ$} \hbox to 5pt{\hfill \vrule height 12pt width 0,3pt\hfill} \hbox{$\circledcirc$}}\hglue -3pt\hrulefill\lower 2pt\hbox{$\circ$} } }

& $4$ &$1$ &$n$ &$0$& II \\
(3) & & & && &\\
\hline
& & & & & &\\
$(E_7,\alpha_7)$ &

\hskip 10pt \hbox to 5cm {\lower 2pt\hbox{$\circ$}\hrulefill\lower 2pt
\hbox{$\circ$}\hrulefill \lower 2pt \vtop {\offinterlineskip \hbox{$\circ$} \hbox to 5pt{\hfill \vrule height 12pt width 0,3pt\hfill} \hbox{$\circ$}}\hrulefill\lower 2pt\hbox{$\circ$}\hrulefill\lower 2pt\hbox{$\circ$}\hrulefill\lower 2pt\hbox{$\circledcirc$} }

& $8$ &$1$ &$3$ &$0$&II\\
 (27)&   & && & & \\
\hline
& & & & &&\\
$(A_{2n-1},\alpha_{n})$ & \hskip -80pt\hbox to 5 cm {
\offinterlineskip \lower 2pt\hbox{$\bullet$}...\hglue -0,1pt
{\vrule height  0,4pt depth 0pt width 0,5 cm}\hglue-1pt\lower 2pt \hbox{$\circ$}\hglue -1,2pt
{\vrule height  0,4pt depth 0pt width 0,5 cm}\lower 2pt\hbox{$\bullet$}}\hglue -3,2pt...........\hglue -0,1pt
 {\offinterlineskip \lower 2pt\hbox{$\bullet$} \hglue -3,2pt
{\vrule height  0,4pt depth 0pt width 0,5 cm}\lower 2pt\hbox{$\circledcirc$}
\hglue -3,2pt
{\vrule height  0,4pt depth 0pt width 0,5 cm}\lower 2pt\hbox{$\bullet$}}........{\hglue -0,2pt
{\offinterlineskip \lower 2pt\hbox{$\circ$}}{\vrule height  0,4pt depth 0pt width 0,5 cm}...  \hglue -4,2pt{ \offinterlineskip \lower 2pt \hbox{$\bullet$} }
}

& $2$ &$m$ &$n$ &$0$&III   \\
& & & & &&\\
&$(\overline \Delta,\alpha_0)=(A_{2nm-1},\alpha_{nm})$ et $m=2$ dans le cas r\éel& & & & &  \\
 
\hline
& &  & & & &\\
$(C_{n},\alpha_{n})$ &
\hskip 10pt \hbox to 5 cm {\offinterlineskip \lower 2pt\hbox{$\bullet$} \hglue -3,2pt
{\vrule height  0,4pt depth 0pt width 0,5 cm}\lower 2pt\hbox{$\circ$}
  \dotfill \hbox to 1,3 cm {\offinterlineskip\lower 2pt\hbox{$\circ$}\kern -3pt   \hrulefill\kern -3pt\lower 2pt 
\hbox{  \offinterlineskip  {$\bullet$} \hglue -8pt\vbox{ {\hrule height 0,3pt width 0,6cm}\vskip 3pt{\hrule height 0,3pt width 0,6cm}
\vskip 0,6pt} \hglue -16pt $<$   \hglue -1,8pt{$\circledcirc$}}}}
& $1$ &$2$ &$n$ &$4$& III\\
& & & & & &\\
\hline

\end{tabular}    
\vskip 1mm
$d=\overline d.d_1^2$ et $d_n=nd_1$

  \newpage

 \centerline{\bf Tableau 2}
\vskip 5mm
        
\noi  Liste des \PVs de type parabolique absolument irr\éductibles,   r\éguliers de type classique  de la section $7.$

 \noi Le num\éro indiqu\é entre parenth\èse est celui qui lui correspond dans la classification de Sato-Kimura (\cite{satokimura},$\S7,$ table I p.144).\\ \\

\noi\begin{tabular}{| c | c | c | c | c | }
\hline
 & & & &   \\
Type &$(\overline {\Delta},\alpha_0)$& $(\Delta,\lambda_0)$ & Diagramme de Satake de $\g\otimes_\F\overline F$&  Conditions   \\
& & & &  
\\
\hline
& & & &  \\
BI(n,k)&$(B_n,\alpha_k)$ & $(B_m,\lambda_k)$ & 
  \hbox to 4cm {\offinterlineskip \lower 2pt\hbox{$\circ$} \hglue -3,2pt
{\vrule height  0,4pt depth 0pt width 0,5 cm}\lower 2pt\hbox{$\circ$}
\hglue -3,2pt
{\vrule height  0,4pt depth 0pt width 0,5 cm}\lower 2pt\hbox{$\circ$}
  \dotfill \hbox to 0,4 cm  
 { \offinterlineskip \lower 2pt\hbox{$\circ$} \hglue -3,2pt
{\vrule height  0,4pt depth 0pt width 0,5 cm}}\hglue 7pt{\lower 2pt\hbox{$\bullet$}}  
  \dotfill \hbox to 0,4 cm 
    { \offinterlineskip \lower 2pt \hbox{$\bullet$}\hglue -3,2pt
{\vrule height  0,4pt depth 0pt width 0,5 cm}\kern -1pt   \hrulefill     \kern -3,4pt\lower 2pt \hbox{ 
\offinterlineskip  {$\bullet$} \hglue -7pt\vbox{{\hrule height 0,3pt width 0,6cm}\vskip 3,5pt{\hrule height 0,3pt width 0,6cm}\vskip 0,3pt }\hglue -12pt $>$   \hglue -1,8pt $\bullet$}}}

& $ k\≤m\≤n $   \\
& (15) & &
$\hskip 1,3cm \alpha_1.....\hskip 0,6cm.....\alpha_m.....\hskip 1cm.....\alpha_n$  &$1\≤3k\≤2n-1$  \\
 & & & &   \\
\hline
& & & &  \\
DI(n,k)&$(D_n,\alpha_k)$ & $(D_n,\alpha_k)$&  \hbox{\unitlength=0,5pt
\begin{picture}(300,30)(-50,35)
\put(10,40){\circle{10}}
\put(15,40){\line (1,0){30}}
\put(50,40){\circle{10}}
\put(60,40){\circle*{1}}
\put(65,40){\circle*{1}}
\put(70,40){\circle*{1}}
\put(75,40){\circle*{1}}
\put(80,40){\circle*{1}}
\put(85,40){\circle*{1}}
\put(90,40){\circle*{1}}
\put(95,40){\circle*{1}}
\put(100,40){\circle*{1}}
\put(110,40){\circle{10}}
\put(115,40){\line (1,-1){21}}
\put(140,15){\circle{10}}
\put(115,40){\line (1,1){21}}
\put(140,65){\circle{10}}
 \end{picture}}
&$1\≤3k\≤2n-2$ \\
&(15)&&& \\
&&$(B_m,\lambda_k)$&

\hbox to 6cm {\offinterlineskip \lower 2pt\hbox{$\circ$} \hglue -3,2pt
{\vrule height  0,4pt depth 0pt width 0,5 cm}\lower 2pt\hbox{$\circ$}
\hglue -3,2pt
{\vrule height  0,4pt depth 0pt width 0,5 cm}\lower 2pt\hbox{$\circ$}
  \dotfill \hbox to 0,4 cm  
 { \offinterlineskip \lower 2pt\hbox{$\circ$} \hglue -3,2pt
{\vrule height  0,4pt depth 0pt width 0,5 cm}}\hglue 7pt{\lower 2pt\hbox{$\bullet$}}  
 \dotfill \hbox to 2 cm {\lower 2pt
\hbox{$\bullet$}\hrulefill \lower 2pt \vtop {\offinterlineskip \hbox{$ \bullet$} \hbox to 5pt{\hfill \vrule height 12pt width 0,3pt\hfill} \hbox{$ \bullet$}}\hrulefill\lower 2pt\hbox{$\bullet$} } }

& $k\≤m\≤n$    \\
& & &\hskip -0,8cm$\alpha_1$.....$\hskip 0,9cm.....\alpha_m.....\hskip 0,7cm.....\alpha_n$ &   \\
 & & &  
 \hbox{\unitlength=0,5pt
\begin{picture}(300,30)(-50,60)
\put(10,40){\circle{10}}
\put(15,40){\line (1,0){30}}
\put(50,40){\circle{10}}
\put(60,40){\circle*{1}}
\put(65,40){\circle*{1}}
\put(70,40){\circle*{1}}
\put(75,40){\circle*{1}}
\put(80,40){\circle*{1}}
\put(85,40){\circle*{1}}
\put(90,40){\circle*{1}}
\put(95,40){\circle*{1}}
\put(100,40){\circle*{1}}
\put(110,40){\circle{10}}
\put(115,40){\line (1,-1){21}}
\put(140,15){\circle{10}}
\put(115,40){\line (1,1){21}}
\put(140,65){\circle{10}}
\put(140,45) {\vector(0,-1){23}}
\put(140,35) {\vector(0,1){23}}
\end{picture}}&       \\
&&$(B_{n-1},\lambda_k)$&& $k\≤n-2$\\
  &&&&\\
  &&&&\\
  \hline
  &&&&\\
  CI&$(C_n,\alpha_{2p})$&$(C_n,\alpha_{2p})$ &
\hskip -25pt \hbox to 5cm {\offinterlineskip \lower 2pt\hbox{$\circ$} \hglue -3,2pt
{\vrule height  0,4pt depth 0pt width 0,5 cm}\lower 2pt\hbox{$\circ$}
\hglue -3,2pt
{\vrule height  0,4pt depth 0pt width 0,5 cm}\lower 2pt\hbox{$\circ$}
  \dotfill \hbox to 0,4 cm  
  {\offinterlineskip\lower 2pt\hbox{$\circledcirc$}}  
  \dotfill \hbox to 0,4 cm 
    { \offinterlineskip \lower 2pt \hbox{$\circ$}\hglue -3,2pt
{\vrule height  0,4pt depth 0pt width 0,5 cm}\kern -1pt   \hrulefill     \kern -3,4pt\lower 2pt \hbox{ 
\offinterlineskip  {$\circ$} \hglue -7pt\vbox{{\hrule height 0,3pt width 0,6cm}\vskip 3,5pt{\hrule height 0,3pt width 0,6cm}\vskip 0,3pt }\hglue -12pt $<$   \hglue -1,8pt $\circ$}}}
& $1\≤3p\≤n-1 $   \\
&(13) &  & &  \\

\hline
&&&&\\
  CII&$(C_n,\alpha_{2p})$&$(C_{\frac{n}{2}},\alpha_{p})$ &
 \hskip -20pt \hbox to 5cm {\offinterlineskip \lower 2pt\hbox{$\bullet$} \hglue -3,2pt
{\vrule height  0,4pt depth 0pt width 0,5 cm}\lower 2pt\hbox{$\circ$}
\hglue -3,2pt
{\vrule height  0,4pt depth 0pt width 0,5 cm}\lower 2pt\hbox{$\bullet$}
  \dotfill \hbox to 0,4 cm  
  {\offinterlineskip\lower 2pt\hbox{$\circledcirc$}}  
  \dotfill \hbox to 0,4 cm 
    { \offinterlineskip \lower 2pt \hbox{$\circ$}\hglue -3,2pt
{\vrule height  0,4pt depth 0pt width 0,5 cm}\kern -1pt   \hrulefill     \kern -3,4pt\lower 2pt \hbox{ 
\offinterlineskip  {$\bullet$} \hglue -7pt\vbox{{\hrule height 0,3pt width 0,6cm}\vskip 3,5pt{\hrule height 0,3pt width 0,6cm}\vskip 0,3pt }\hglue -12pt $<$   \hglue -1,8pt $\circ$}}}& $1\≤3p\≤n-1 $\\
& & & &  \\
&  & 
 $(BC_m,\lambda_{p})$& \hskip -20pt \hbox to 5cm {\offinterlineskip \lower 2pt\hbox{$\bullet$} \hglue -3,2pt
{\vrule height  0,4pt depth 0pt width 0,5 cm}\lower 2pt\hbox{$\circ$}
\hglue -3,2pt
{\vrule height  0,4pt depth 0pt width 0,5 cm}\lower 2pt\hbox{$\bullet$}
  \dotfill \hbox to 0,4 cm  
  {\offinterlineskip\lower 2pt\hbox{$\circledcirc$}}  
  \dotfill \hbox to 0,4 cm 
    { \offinterlineskip \lower 2pt \hbox{$\bullet$}\hglue -3,2pt
{\vrule height  0,4pt depth 0pt width 0,5 cm}\kern -1pt   \hrulefill     \kern -3,4pt\lower 2pt \hbox{ 
\offinterlineskip  {$\bullet$} \hglue -7pt\vbox{{\hrule height 0,3pt width 0,6cm}\vskip 3,5pt{\hrule height 0,3pt width 0,6cm}\vskip 0,3pt }\hglue -12pt $<$   \hglue -1,8pt $\bullet$}}}& $p\≤m$ \\
&&&&cas r\éel  \\
 & & & &   \\
\hline
& & & &   \\

DIII&$(D_n,\alpha_{2p})$ & $(C_m,\lambda_{p})$&

  \hskip -0,8cm\hbox to 5 cm {\offinterlineskip \lower 2pt\hbox{$\bullet$}
\hglue -3pt{\vrule height .4pt depth 0pt width 0,5cm}\lower 2pt\hbox{$\circ$}
\hglue -7pt{ \vrule height .4pt depth 0pt width 0,5cm}\lower 2pt\hbox{$\bullet$}
 \dotfill \hbox to 2 cm {\lower 2pt\hbox{$\circ$}\hglue -1,3pt
\hrulefill\lower 2pt
\hbox{$\bullet$}\hrulefill \lower 2pt \vtop {\offinterlineskip \hbox{$\circ$} \hbox to 5pt{\hfill \vrule height 12pt width 0,3pt\hfill} \hbox{$\circ$}}\hrulefill\lower 2pt\hbox{$\bullet$} } }

&  $n$ pair    \\
&  & & &  \\
 &&$(BC_m,\lambda_p)$
& \hskip 15pt \hbox to 5 cm {\offinterlineskip \lower 2pt\hbox{$\bullet$}
\hglue -3pt{\vrule height .4pt depth 0pt width 0,5cm}\lower 2pt\hbox{$\circ$}
\hglue -7pt{ \vrule height .4pt depth 0pt width 0,5cm}\lower 2pt\hbox{$\bullet$}
 \dotfill \hbox to 2 cm {\lower 2pt\hbox{$\bullet$}\hglue -1,3pt
\hrulefill\lower 2pt
\hbox{$\circ$}\hrulefill \lower 2pt \vtop {\offinterlineskip \hbox{$\bullet$} \hbox to 5pt{\hfill \vrule height 12pt width 0,3pt\hfill} \hbox{$\circ$}}\hrulefill\lower 2pt\hbox{$\circ$} } } \hbox{\unitlength=0,5pt
\begin{picture}(100,-20)(-50,0)
\put(-100,-30){\vector (1,1){20.5}}
\put(-100,-30){\vector (-1,-1){2.5}}
\end{picture}}
  
& $n$ impair   \\  & &  & &    \\ 

& & $m=[\ds\frac{n}{2}]$&$\R:n-3p\≥4$ &  \\
& & & $\goth p$-adique : $n-3p\≥2$ et $n-p$ impair& \\
& & & &    \\
\hline

\end{tabular}    
\vskip 1mm
 
 \newpage

 \centerline{\bf Tableau 3}
\vskip 5mm
        
\noi  Liste des \PVs de type parabolique absolument irr\éductibles,   r\éguliers de type  exceptionnels, i.e.$\overline \g$ est de type exceptionnel, de la section 8 et les 2 exemples r\éels du $\S 5.3.2.$

    \noi Le num\éro indiqu\é entre parenth\èse est celui qui lui correspond dans la classification de Sato-Kimura (\cite{satokimura},$\S7,$ table I p.144).\\ 
    
   \noi $\∑_P=\{\mu\in \∑_0\ |\ \mu/{\goth t}\not=0\},$ le ``degr\é" d\ésigne le degr\é de l'\irf.\\

\noi\begin{tabular}{| c | c | c | c | c | c | c |c|}
\hline
 & & & & & && \\
$(\overline {\Delta},\alpha_0)$& $(\Delta,\lambda_0)$ & Diagramme de Satake &  Conditions & type de $\g$ & dim($\g_2$)   & degr\é&
$\Sigma_P$\\
& &$\circledcirc\leftrightarrow\alpha_0\ ,\ \boxcircle\leftrightarrow \∑_P$ & & & &&
\\
& & & & && &\\
\hline
\hline
& & & & && &\\
$(E_7,\alpha_6)$& $(E_7,\alpha_6)$ &  \hskip 10pt \hbox to 3,5 cm {\lower 2pt\hbox{$\boxcircle$}\hrulefill\lower 2pt
\hbox{$\circ$}\hrulefill \lower 2pt \vtop {\offinterlineskip \hbox{$\circ$} \hbox to 5pt{\hfill \vrule height 12pt width 0,3pt\hfill} \hbox{$\circ$}}\hrulefill\lower 2pt\hbox{$\circ$}\hrulefill\lower 2pt\hbox{$ \circledcirc$}\hrulefill\lower 2pt\hbox{$\circ$} }  & &
  &10  &4&  $\alpha_1$  \\
& (20)& &
 & &  && \\
 & & & && & & \\
\hline
& & & & & &&\\
 $(E_7,\alpha_6)$&$(F_4,\lambda_4)$ & \hskip 10pt \hbox to 3cm {\lower 2pt\hbox{$\boxcircle$}\hrulefill\lower 2pt
\hbox{$ \circ$}\hrulefill \lower 2pt \vtop {\offinterlineskip \hbox{$ \circ$} \hbox to 5pt{\hfill \vrule height 12pt width 0,3pt\hfill} \hbox{$\bullet$}}\hrulefill\lower 2pt\hbox{$\bullet$}\hrulefill\lower 2pt\hbox{$\circledcirc$}\hrulefill\lower 2pt\hbox{$ \bullet$} } & &EVI & 10&4&$\lambda_1$\\
&&&& &&&\\
&&  & & & & &  \\
 
&  & &&&&&\\
  \hline
  &&&&&&&\\
 $(E_7,\alpha_6)$ &$(C_3,\lambda_{2})$&  \hskip 10pt \hbox to 3cm {\lower 2pt\hbox{$\boxcircle$}\hrulefill\lower 2pt
\hbox{$\bullet$}\hrulefill \lower 2pt \vtop {\offinterlineskip \hbox{$\bullet$} \hbox to 5pt{\hfill \vrule height 12pt width 0,3pt\hfill} \hbox{$\bullet$}}\hrulefill\lower 2pt\hbox{$\bullet$}\hrulefill\lower 2pt\hbox{$\circledcirc$}\hrulefill\lower 2pt\hbox{$\circ$} }
 &  $\F=\R$&EVII
 & 10&4&$\lambda_1$   \\
& &  & & & &&\\
 & & &  & && & \\
  
\hline
 \hline

& & &&&&&\\
 $(E_7,\alpha_2)$ &  $(E_7,\alpha_2)$&  \hskip 10pt \hbox to 3cm {\lower 2pt\hbox{$\boxcircle$}\hrulefill\lower 2pt
\hbox{$\circ$}\hrulefill \lower 2pt \vtop {\offinterlineskip \hbox{$\circ$} \hbox to 5pt{\hfill \vrule height 12pt width 0,3pt\hfill} \hbox{$\circledcirc$}}\hglue-3pt{\hrulefill\lower 2pt\hbox{$\circ$}}\hrulefill\lower 2pt\hbox{$\circ$}\hrulefill\lower 2pt\hbox{$\circ$} }
 &  
 &  & 7&7&$\alpha_1$ \\
&(6)&  & && && \\
 & && && & &  \\
  \hline
\hline
& & & & & & & \\

& & & & & &&  \\
  
$(E_8,\alpha_1)$&$(E_8,\alpha_1)$&  \hskip 10pt \hbox to 4cm {\lower 2pt\hbox{$\circledcirc$}\hrulefill\lower 2pt
\hbox{$\circ$}\hrulefill \lower 2pt \vtop {\offinterlineskip \hbox{$\circ$} \hbox to 5pt{\hfill \vrule height 12pt width 0,3pt\hfill} \hbox{$\circ$}}\hrulefill\lower 2pt\hbox{$\circ$}\hrulefill\lower 2pt\hbox{$\circ$}\hrulefill\lower 2pt\hbox{$\circ$}\hrulefill\lower 2pt\hbox{$\boxcircle$} }& &&14&8&$\alpha_8$\\

 &(24) &  & & & & & \\

& & & &  & && \\
\hline
 
 & & & && & & \\
  
$(E_8,\alpha_1)$&$(F_4,\lambda_4)$&  \hskip 10pt \hbox to 4cm {\lower 2pt\hbox{$\circledcirc$}\hrulefill\lower 2pt
\hbox{$\bullet$}\hrulefill \lower 2pt \vtop {\offinterlineskip \hbox{$\bullet$} \hbox to 5pt{\hfill \vrule height 12pt width 0,3pt\hfill} \hbox{$\bullet$}}\hrulefill\lower 2pt\hbox{$\bullet
$}\hrulefill\lower 2pt\hbox{$\circ$}\hrulefill\lower 2pt\hbox{$\circ$}\hrulefill\lower 2pt\hbox{$\boxcircle$} }&$\F=\R$ & EIX &14&8&$\lambda_1$\\

 & &  & &&  & & \\
 & & & && &&  \\

\hline
 \end{tabular}    
 \newpage
 
\noi\begin{tabular}{| c | c | c | c | c | c |c|c|}
\hline
& &  & & & & & \\ 
  $(\overline {\Delta},\alpha_0)$& $(\Delta,\lambda_0)$ & Diagramme de Satake &  Conditions & type de $\g$& dim($\g_2$) &degr\é&$\Sigma_P$\\
& &$\circledcirc\leftrightarrow\alpha_0\ ,\ \boxcircle\leftrightarrow \∑_P$ & && & &\\
& &  & &  && & \\

  \hline
  \hline
 & &  & & & & &  \\ 
  $(F_4,\alpha_1)$&$(F_4,\alpha_1)$ &  \hbox to 2,5cm{ \offinterlineskip \lower 2pt\hbox{$\circledcirc$}\hrulefill\lower 2pt\hbox{\offinterlineskip {$\circ$}\hglue -1,8pt\vbox{
 {\hrule height 0,3pt width 0,6cm}\vskip 3,5pt{\hrule height 0,3pt width 0,6cm}
\vskip 0,3pt} \hglue -16pt  $>$ \hglue -1,8pt  {$\circ$}}\hrulefill\lower 2pt\hbox{$\boxcircle$} }
&&&  1 (d=1)&4&$\alpha_4$ \\
 
& (14)& & & & & & \\

\hline
\hline
& & & & && & \\

 $(E_6,\alpha_2)$& $(E_6,\alpha_2)$&  \hskip 10pt \hbox to 3 cm {\lower 2pt\hbox{$\boxcircle$}\hrulefill\lower 2pt
\hbox{$\circ$}\hrulefill \lower 2pt \vtop {\offinterlineskip \hbox{$\circ$} \hbox to 5pt{\hfill \vrule height 12pt width 0,3pt\hfill} \hbox{$\circledcirc$}}\hglue -3,5pt\hrulefill\lower 2pt\hbox{$\circ$}\hrulefill\lower 2pt\hbox{$\boxcircle$} } &
&   &   1 (d=2)&4&$\alpha_1,\alpha_6$ \\
&(5)  & & && && \\
\hline
& &   \begin{picture}(100,1)(-50,-15)
\put(-25,-22){\vector (-1,-1){10.5}}
\put(-25,-22){\line (1,0){50}}
\put(25,-22){\vector (1,-1){10.5}}
\end{picture}
&&&&&\\
&&
  \begin{picture}(100,10)(-50,1)
\put(-6,3){\vector (-1,-1){10}}
\put(-6,3){\line (1,0){16}}
\put(10,3){\vector (1,-1){10}}
\end{picture}&&&&&\\
$(E_6,\alpha_2)$&$(F_4,\lambda_1)$&

  \hskip 10pt \hbox to 3 cm {\lower 2pt\hbox{$\boxcircle$}\hrulefill\lower 2pt
\hbox{$\circ$}\hrulefill \lower 2pt \vtop {\offinterlineskip \hbox{$\circ$} \hbox to 5pt{\hfill \vrule height 12pt width 0,3pt\hfill} \hbox{$\circledcirc$}}\hglue -3,5pt\hrulefill\lower 2pt\hbox{$\circ$}\hrulefill\lower 2pt\hbox{$\boxcircle$} } &&EII&1 (d=2)&4&$\lambda_4$\\
&  & & && && \\

\hline
    &  &

 \begin{picture}(100,30)(-50,3)
\put(-25,7){\vector (-1,-1){10.5}}
\put(-25,7){\line (1,0){52}}
\put(27,7){\vector (1,-1){10.5}}
\end{picture}&&&&&\\
 
  $(E_6,\alpha_2)$&$(BC_2,\lambda_1)$&
  \hskip 10pt \hbox to 3 cm {\lower 2pt\hbox{$\boxcircle$}\hrulefill\lower 2pt
\hbox{$\bullet$}\hrulefill \lower 2pt \vtop {\offinterlineskip \hbox{$\bullet$} \hbox to 5pt{\hfill \vrule height 12pt width 0,3pt\hfill} \hbox{$\circledcirc$}}\hglue -3,5pt\hrulefill\lower 2pt\hbox{$\bullet$}\hrulefill\lower 2pt\hbox{$\boxcircle$} } & $\F=\R$&EIII &1&4&$\lambda_2$  \\
  &  &  &  &  & ($\S 5.3)$& &\\ 
  
& &  & &  & && \\ 
\hline
\hline
&  & & &&& & \\
$(E_7,\alpha_1)$ &$(E_7,\alpha_1)$& \hskip 10pt \hbox to 3,5 cm {\lower 2pt\hbox{$\circledcirc$}\hrulefill\lower 2pt
\hbox{$\circ$}\hrulefill \lower 2pt \vtop {\offinterlineskip \hbox{$\circ$} \hbox to 5pt{\hfill \vrule height 12pt width 0,3pt\hfill} \hbox{$\circ$}}\hrulefill\lower 2pt\hbox{$\circ$}\hrulefill\lower 2pt\hbox{$\boxcircle$}\hrulefill\lower 2pt\hbox{$\circ$} }

&  
  
&  & 1 (d=4)&4&$\alpha_6$  \\  &(23) &  & &  & && \\

& & & &  & && \\
\hline
 &  && && &&  \\
$(E_7,\alpha_1)$ &$(F_4,\lambda_1)$& \hskip 10pt \hbox to 3,5 cm {\lower 2pt\hbox{$\circledcirc$}\hrulefill\lower 2pt
\hbox{$\circ$}\hrulefill \lower 2pt \vtop {\offinterlineskip \hbox{$\circ$} \hbox to 5pt{\hfill \vrule height 12pt width 0,3pt\hfill} \hbox{$\bullet$}}\hrulefill\lower 2pt\hbox{$\bullet$}\hrulefill\lower 2pt\hbox{$\boxcircle$}\hrulefill\lower 2pt\hbox{$\bullet$} }

&  
  
& EVI  & 1(d=4)&4&$\lambda_4$  \\  & & & & & & &  \\

\hline
 &  & && && & \\
$(E_7,\alpha_1)$ &$(C_3,\lambda_1)$& \hskip 10pt \hbox to 3,5 cm {\lower 2pt\hbox{$\circledcirc$}\hrulefill\lower 2pt
\hbox{$\bullet$}\hrulefill \lower 2pt \vtop {\offinterlineskip \hbox{$\bullet$} \hbox to 5pt{\hfill \vrule height 12pt width 0,3pt\hfill} \hbox{$\bullet$}}\hrulefill\lower 2pt\hbox{$\bullet$}\hrulefill\lower 2pt\hbox{$\boxcircle$}\hrulefill\lower 2pt\hbox{$\circ$} }

&  $\F=\R$&EVII&
  
1  &4&$\lambda_2$    \\  
 
 & & & & & ($\S 5.3$) & & \\ 
 & &  & & & & & \\ 
\hline
\hline
& &  & &  &&&  \\ 

$(E_8,\alpha_8)$&$(E_8,\alpha_8)$& \hskip 10pt \hbox to 4cm {\lower 2pt\hbox{$\boxcircle$}\hrulefill\lower 2pt
\hbox{$\circ$}\hrulefill \lower 2pt \vtop {\offinterlineskip \hbox{$\circ$} \hbox to 5pt{\hfill \vrule height 12pt width 0,3pt\hfill} \hbox{$\circ$}}\hrulefill\lower 2pt\hbox{$\circ$}\hrulefill\lower 2pt\hbox{$\circ$}\hrulefill\lower 2pt\hbox{$\circ$}\hrulefill\lower 2pt\hbox{$\circledcirc$} }&&&1 (d=8)&4&$\alpha_1$\\
 
 &(29) &&  & &  & & \\

& & && & &  & \\
\hline
& &  & & & &&  \\

$(E_8,\alpha_8)$&$(F_4,\lambda_1)$&  \hskip 10pt \hbox to 4cm {\lower 2pt\hbox{$\boxcircle$}\hrulefill\lower 2pt
\hbox{$\bullet$}\hrulefill \lower 2pt \vtop {\offinterlineskip \hbox{$\bullet$} \hbox to 5pt{\hfill \vrule height 12pt width 0,3pt\hfill} \hbox{$\bullet$}}\hrulefill\lower 2pt\hbox{$\bullet$}\hrulefill\lower 2pt\hbox{$\circ$}\hrulefill\lower 2pt\hbox{$\circ$}\hrulefill\lower 2pt\hbox{$\circledcirc$} }&$\F=\R$& EIX&1 (d=8)&4&$\lambda_4$\\
 & & & & &  && \\
\hline

\end{tabular}

\end{document}